\documentclass[twoside,12pt]{article}
\usepackage[sumlimits,intlimits,namelimits]{amsmath}
\usepackage{amsthm,amssymb,amstext,amsfonts}
\usepackage[mathscr]{eucal}
\usepackage{makeidx}
\usepackage{exscale}
\newcommand{\Em}[1]{\emph{#1}} 
\newcommand{\De}[1]{\emph{#1}} 
\newcommand{\unity}{{1\!\!\!\:\mathrm{l}}}
\newcommand{\var}{\delta}
\DeclareMathOperator{\vol}{vol}
\DeclareMathOperator{\tr}{tr}
\DeclareMathOperator*{\res}{Res}
\newcommand{\xx}[1]{\Hat{#1}} 
\newcommand{\yy}[1]{\Check{#1}} 
\newcommand{\fourier}[1]{\Hat{#1}} 
\newcommand{\triv}[1]{\Tilde{#1}} 
\newcommand{\comp}{} 
\newcommand{\Op}[1]{\mathsf{#1}} 
\DeclareMathAlphabet{\Mat}{U}{eur}{m}{n} 
\DeclareMathAlphabet{\Set}{U}{eus}{m}{n} 
\newcommand{\Spa}[1]{\mathrm{#1}} 
\newcommand{\Sh}[1]{\mathcal{#1}} 
\newcommand{\func}[1]{#1} 
\newcommand{\parameter}[1]{#1} 
\newcommand{\Func}[1]{\mathcal{#1}} 
\newcommand{\lattice}{\Lambda}
\newcommand{\dual}{^{\!\ast}}
\newcommand{\torus}{\mathbb{R}^2/\lattice}
\newcommand{\willmore}{\mathscr{W}}
\newcommand{\bloch}{\mathcal{B}}
\newcommand{\fermi}{\mathcal{F}}
\newcommand{\moduli}{\mathcal{M}}
\newcommand{\sobolev}[1]{W^{#1}}
\newcommand{\banach}[1]{L^{\!#1}}
\newcommand{\parity}{\mathscr{P}}
\renewcommand{\qed}{{\bf\hspace{\fill}q.e.d.}}
\newtheorem{Lemma}{Lemma}[section]
\newtheorem{Theorem}[Lemma]{Theorem}
\newtheorem{Corollary}[Lemma]{Corollary}
\newtheorem{Proposition}[Lemma]{Proposition}

\newtheorem{Remark}[Lemma]{Remark}
\begin{document}
\date{}

\title{A proof of the Willmore conjecture}

\author{Martin U. Schmidt\thanks{Supported by DFG,
SFB~288 ``Differentialgeometrie und Quantenphysik''}\\
Institut f\"ur Theoretische Physik\\
Freie Universit\"at Berlin\\
Arnimallee 14\\
D-14195 Berlin\\
email: {martin.schmidt@physik.fu-berlin.de}}
\begin{titlepage}
\maketitle
\end{titlepage}
\tableofcontents 
\newpage

\section{Introduction}\label{section introduction}
\subsection{The main result}\label{subsection main result}

The Willmore functional
\index{Willmore!functional}
associates to each immersion of the
two--dimensional torus $\mathbb{T}^2$ into the three--dimensional
Euclidean space $\mathbb{R}^3$ the integral
$\int\limits_{\mathbb{T}^2} H^2 d\mu$ of the square
of the mean curvature of the induced metric with respect to the
corresponding measure over the two--dimensional torus. This functional
has an obvious extension to immersions of two--dimensional manifolds into
higher--dimensional Riemannian manifolds. In this article we shall
mainly consider immersions of $\mathbb{T}^2$ into $\mathbb{R}^3$ and
sometimes into $\mathbb{R}^4$. This functional is invariant under
conformal transformations of $\mathbb{R}^3$. The
Willmore conjecture states that up to conformal transformations
\index{Willmore!conjecture}
this functional has a unique absolute minimum.
The corresponding immersion is the \Em{Clifford torus},
\index{Clifford torus}
which is the surface of a rotated circle of
radius $r$ around some axis in the plane spanned by the
circle with distance $\sqrt{2}r$ from the center of the circle.
The corresponding Willmore functional is equal to $2\pi^2$.

All smooth metrics on $\mathbb{T}^2$ are conformally equivalent to
a flat metric and each flat metric corresponds to exactly one complex
\index{conformal!class!of an immersion}
structure on $\mathbb{T}^2$. Hence, the conformal classes of these
immersions are equal to the real two--dimensional manifold $\moduli_1$
\index{moduli space!$\moduli_1$ of conformal classes}
\index{conformal!class!moduli space $\moduli_1$ of $\sim$es}
of moduli of compact Riemann surfaces of genus one. Consequently, the
Willmore functional may be restricted to all immersions of some given
conformal class, and the infima of these restrictions fit together to
some function on this moduli space $\moduli_1$. The main result of
this article is the calculation of another function on $\moduli_1$,
which yields a lower bound of the former function. Furthermore,
this latter function has a unique minimum at the conformal class
corresponding to the flat torus $\mathbb{R}^2/\mathbb{Z}^2$.
Finally, for all rectangular conformal classes
this lower bound is sharp and the corresponding immersions
are unique up to conformal transformations of $\mathbb{R}^3$.
The minimum is realized by the \Em{Clifford torus}.
We conjecture that the latter function describes
the minimum of the restriction of the Willmore functional
to conformal immersions of two--dimensional tori $\mathbb{T}^2$
into $\mathbb{R}^4$
in dependence of the conformal class of $\mathbb{T}^2$.

We now outline the strategy of the proof. It contains

\subsection{A transformation of the variational problem}
\label{subsection transformation}

Usually the determination of the minimum splits into two problems:
\begin{description}
\item[Proof of the existence of a minimizer.]
\index{minimizer!existence of a $\sim$}
\index{existence!of a minimizer}
\item[Determination of the lowest critical point.]
\end{description}
As part of the \Em{proof of the existence} we shall specify the
class of objects (immersions) we consider. In view of the
first problem this class should be chosen large enough in order
to contain enough relative minimizers.
But in view of the second problem this class should be chosen
small enough in order to simplify the classification.
For the Willmore functional both problems seem to be difficult to solve.

The first problem was solved by L.\ Simon
\cite{Si1,Si2} using geometric measure theory. With his methods,
however, it seems hopeless to prove or disprove the existence
of a minimizer for the restrictions to arbitrary conformal classes.

The critical points of the Willmore functional are called
Willmore tori\index{Willmore!torus},
and  the critical points of the restrictions to some
conformal classes are called \De{Constrained Willmore tori}
(the precise definition is given in 
\De{Constrained Willmore tori}~\ref{constrained Willmore tori}).
\index{Willmore!torus!constrained $\sim$}
The stereographic projections of minimal tori in $S^3$
are Willmore tori, and the area of these minimal tori is equal to the
Willmore functional. In spite of a quite satisfactory description of
the minimal tori in $S^3$ given in \cite{Hi} it was not clear which
minimal torus corresponds to the lowest value of the Willmore functional.
U.\ Pinkall was the first to construct Willmore tori which are not
stereographic projections of minimal tori in $S^3$
\cite{Pi1}. Afterwards many other Willmore tori were constructed
\cite{Ba,BB}, indicating that a general classification of the lowest
critical points seems to be even more difficult.

Thirdly, there exist several attempts to pass around these two
problems. In this direction first we mention the topological lower
bound $4\pi$ (compare with \cite[Section~7.6]{Wi2} and references therein),
which is sharp for immersions from $S^2$ into $\mathbb{R}^3$.
Furthermore, in \cite{LY} P.\ Li and S.\ T.\ Yau  proved several
lower bounds of the Willmore functional, one of which implies that all
immersions, whose conformal classes belong to a closed subset including the
conformal class of the \Em{Clifford torus}, have a Willmore functional
larger or equal to $2\pi^2$.

We want to present a new approach to both problems. The starting point
is an observation due to I.\ Taimanov and which initiated this project.
It states that the local Weierstra{\ss} representation \cite{Ei}
of immersions of surfaces in $\mathbb{R}^3$ can be extended
to a global Weierstra{\ss} representation.
More precisely, for all conformal immersions
\index{conformal!immersion}
of some flat torus $\mathbb{T}^2$ into
$\mathbb{R}^3$, there exists a unique
real potential $U$ on $\mathbb{T}^2$ and some spinor
$\left(\begin{smallmatrix}
\psi_1\\
\psi_2
\end{smallmatrix}\right)$ in the kernel of the Dirac operator
$\left(\begin{smallmatrix}
U & \partial\\
-\Bar{\partial} & U
\end{smallmatrix}\right)$
on $\mathbb{T}^2$ with potential $U$, such that the differential of
the immersion is the $\mathbb{R}^3$--valued one form
$(\Re(\psi_1^2dz-\psi_2^2d\Bar{z}),\Im(\psi_1^2dz-\psi_2^2d\Bar{z}),
\psi_1\Bar{\psi_2}dz+\Bar{\psi}_1\psi_2d\Bar{z})$ on $\mathbb{T}^2$
\cite{Ta1,Fr2}. With respect to the spinor this $\mathbb{R}^3$--valued
one--form is a quadratic form. Therefore, it is single--valued on
$\mathbb{T}^2$. We should remark that these spinors
$\left(\begin{smallmatrix}
\psi_1\\
\psi_2
\end{smallmatrix}\right)$ are in general not single--valued functions
on $\mathbb{T}^2$. The non--uniqueness, however, is soft in the sense
that the difference is only given by the factor $-1$.
More precisely, the universal covering $\mathbb{R}^2$ is a
principal bundle over $\mathbb{T}^2$ whose structure group is the
fundamental group of $\mathbb{T}^2$, which is isomorphic to
$\mathbb{Z}^2$. Any homomorphism of this group into the group
$\{1,-1\}\simeq \mathbb{Z}_2$, induces some
line bundle on $\mathbb{T}^2$. The spinor bundle is then given by
this line bundle tensored with the trivial vector bundle of rank two
over $\mathbb{T}^2$. Obviously there exist four different
homomorphisms from $\mathbb{Z}^2$ into $\mathbb{Z}_2$, and therefore
also four different spin structures on $\mathbb{T}^2$.
The spin structures of the global Weierstra{\ss} representations are
obviously homotopy invariants of the corresponding immersions, and in
\cite{Pi2} it was shown that they are the only homotopy invariants.

This global Weierstra{\ss} representation was generalized to
immersions into $\mathbb{R}^4$ by F.\ Pedit and U.\ Pinkall \cite{PP}.
For this purpose generalized Dirac operators with complex potentials
were considered. They act on the trivial vector bundle of rank two
over $\mathbb{T}^2$ tensored with the line bundles induced by a
real character of the fundamental group of $\mathbb{T}^2$.
These bundles may be considered as
quaternionic line bundles\index{quaternionic!line bundle}

Conversely, for all real potentials, whose Dirac operators have a
non--empty kernel, the corresponding quadratic form of any element of
this kernel is a closed $\mathbb{R}^3$--valued one--form on the torus
$\mathbb{T}^2$. The
integral of this one--form yields a mapping from the universal covering
$\mathbb{R}^2$ of $\mathbb{T}^2$ into $\mathbb{R}^3$.
If the spinor is everywhere non--vanishing,
this mapping is a conformal immersion.
Moreover, if this mapping is the composition of
the universal covering map from $\mathbb{R}^2$ to $\mathbb{T}^2$
with some mapping from $\mathbb{T}^2$ into $\mathbb{R}^3$,
or equivalently if the integrals of the
$\mathbb{R}^3$--valued one--form over all cycles of $\mathbb{T}^2$
vanish, then the corresponding Willmore functional is equal to
the square of the $\banach{2}$--norm of the potential times four.

\newtheorem{Periodicity condition}[Lemma]{Periodicity condition}
\index{periodicity condition}
\index{condition!periodicity $\sim$}
\begin{Periodicity condition}\label{periodicity condition}
The integrals of the two closed one forms
$\psi_1^2dz-\psi_2^2d\Bar{z}$ and
$\psi_1\Bar{\psi_2}dz+\Bar{\psi}_1\psi_2d\Bar{z}$ along all cycles of
$\mathbb{T}^2$ vanish.
\end{Periodicity condition}

An easy calculation
(compare with Corollary~\ref{involutions} and Lemma~\ref{branchpoints})
shows that this condition on any spinor
$\left(\begin{smallmatrix}
\psi_1\\
\psi_2
\end{smallmatrix}\right)$ on $\mathbb{R}^2$ in the kernel of the
Dirac operator
$\left(\begin{smallmatrix}
U & \partial\\
-\Bar{\partial} & U
\end{smallmatrix}\right)$
has another

\newtheorem{Equivalent form}[Lemma]{Equivalent form of the Periodicity
  condition}\index{periodicity condition!equivalent form of the $\sim$}
\begin{Equivalent form}\label{equivalent form}
$$\int\limits_{\mathbb{T}^2}
\psi_1^2(x)d^2x=0,
\int\limits_{\mathbb{T}^2}
\psi_2^2(x)d^2x=0\text{ and }
\int\limits_{\mathbb{T}^2}
\psi_1(x)\Bar{\psi}_2(x)d^2x=0.$$
\end{Equivalent form}

Therefore, we may transform the variational problem of the restriction
of the Willmore functional to some conformal class into the variational
problem of the square of the $\banach{2}$--norm times four
on the space of all real potentials on $\mathbb{T}^2$,
whose Dirac operator has some non--trivial spinor in the kernel
satisfying the \De{Periodicity condition}~\ref{periodicity condition}.
The conformal classes of the corresponding immersions are all equal
to the conformal class of the flat torus $\mathbb{T}^2$. For all
$\banach{2}$--potentials the resulting Dirac operator defines an
operator with compact resolvent on an appropriate Hilbert space of
spinors on $\mathbb{T}^2$. For all spinors in the kernel of this
Dirac operator both equivalent forms of the
\De{Periodicity condition}~\ref{periodicity condition}
are well defined. In fact, we shall see in Remark~\ref{sobolev embedding}
that the spinors in the kernel of belong to the Sobolev spaces
$\sobolev{1,p}(\mathbb{T}^2)\times\sobolev{1,p}(\mathbb{T}^2)$
with $1<p<2$. Due to the Sobolev embedding theorem \cite[5.4 Theorem]{Ad}
this implies that the restrictions of the spinors to a
one--dimensional plane in $\mathbb{T}^2$ are $\banach{2}$--spinors.
In this sense the
\De{Periodicity condition}~\ref{periodicity condition} is
well defined. The \Em{Equivalent form}~\ref{equivalent form}
is well defined since these spinors in the kernel of the Dirac operator
are $\banach{2}$--spinors. Therefore, it is natural
to extend the Willmore functional (defined by the square of the
$\banach{2}$--norm of the potential times four) to the space of all mappings
from $\mathbb{T}^2$ into $\mathbb{R}^3$, defined as the integrals
over the $\mathbb{R}^3$--valued one--forms corresponding to some spinor
in the kernel of the Dirac operator fulfilling the
\De{Periodicity condition}~\ref{periodicity condition}.
These mappings are always continuous, but in
general not smooth. Furthermore, they contain also smooth mappings,
which are not immersions, since the images have branch points.
We consider these mappings to be the maximal domain of definition of
the Willmore functional and therefore make the definitions

\newtheorem{Weierstrass potentials}[Lemma]{Weierstra{\ss} potentials}
\begin{Weierstrass potentials} \label{Weierstrass potentials}
\index{Weierstra{\ss}!potential}
Real potentials $U\in \banach{2}(\mathbb{T}^2)$, whose
Dirac operators contain in their kernel some non--trivial
spinor fulfilling the
\De{Periodicity condition}~\ref{periodicity condition},
are called \De{Weierstra{\ss} potentials}.
\end{Weierstrass potentials}

The main achievement of the observation of I.\ Taimanov is the
description of the space of conformal immersions from a flat torus
$\mathbb{T}^2$ into $\mathbb{R}^3$ as a subspace of the phase space of
the integrable system, whose Lax operator is the Dirac operator
with a pair of potentials on $\mathbb{T}^2$.
Consequently, we use the methods from the theory of integrable systems
related to Lax operators and their spectral curves.
The corresponding evolution is described by the Lax equation
$\frac{\partial L}{\partial t}=\left[L,A\right]$ of a Lax pair $(L,A)$.
This equation indicates that the spectrum of the Lax operator $L$ is
an integral of motion. Roughly speaking, the eigenvalues of the
Lax operator fits together to all integrals of motion of
a completely integrable system \cite{DKN,Kon}.
The application of these methods to differential geometry
accompanied their development almost from the very beginning.
But in contrast to other geometric applications
the present approach uses these methods in order to describe
the whole space of immersions, and not only some special immersions
(Willmore tori, CMC-tori, minimal tori etc.).
This makes the interplay between Hamiltonian systems and
variational analysis available for this problem.
Since Hamiltonian mechanics and variational analysis 
have as common origin Lagrangian mechanics,
it seems natural to apply the new methods of the former theory
to an unsolved problem of the latter.

These Dirac operators with potentials
are the Lax operators\index{Lax operator}
of the Davey--Stewartson equation \cite{Kon},
and the corresponding
spectral curves\index{spectral!curve}\index{curve!spectral $\sim$}
(i.\ e.\ a complete set of integrals of motion)
are the \Em{complex Fermi curves} of these
Dirac operators corresponding to vanishing energy.
These \Em{complex Fermi curves}
\index{Fermi curve!complex $\sim$ $\fermi$}
\index{curve!complex Fermi $\sim$ $\fermi$}
are the subsets of all complex characters of the
fundamental group of $\mathbb{T}^2$ such that the Dirac operators
$\left(\begin{smallmatrix}
V & \partial\\
-\Bar{\partial} & W
\end{smallmatrix}\right)$
acting on the trivial vector bundle of rank two tensored with the
line bundle induced by the character have a non--trivial kernel.
The Bloch theory investigates the common spectrum of a periodic
differential operator with the commuting space translations \cite{Ku}.
For three--dimensional periodic Schr\"odinger operators this common
spectrum is called band structure. The physical \Em{Fermi curves}
are the intersection of the band structure with a fixed energy level
(compare with \cite[\S18.-\S20.]{FKT}).\index{Fermi curve}
\index{curve!Fermi $\sim$} Our \Em{complex Fermi curves}
are the complex continuations of the analogs
for two--dimensional Dirac operators with periodic potentials
instead of three--dimensional periodic Schr\"odinger operators.

\newtheorem{Weierstrass curves}[Lemma]{Weierstra{\ss} curves}
\index{Weierstra{\ss}!curve}
\index{curve!Weierstra{\ss} $\sim$}
\begin{Weierstrass curves}\label{Weierstrass curves}
The \Em{complex Fermi curves} of the \De{Weierstra{\ss} potentials}
are called \De{Weierstra{\ss} curves}.
(This terminology was suggested by I.\ Taimanov.)
\end{Weierstrass curves}

The linear involution $(V,W)\mapsto (W,V)$ and the anti--linear
involution $(V,W)\mapsto (\Bar{V},\Bar{W})$ induce two involutions on
the space of \Em{complex Fermi curves}.
Therefore, the fixed points of these involutions
(i.\ e.\ the space of real potentials $U$)
is some reduction of this integrable system,
and called the modified Novikov--Veselov equation,
since it is a modification of the Novikov--Veselov equation
analogous to the modified Korteweg--de Vries equation \cite{Kon}.
The corresponding \Em{complex Fermi curves} have
a holomorphic and an anti--holomorphic involution.
Unfortunately, the reduction corresponding to the linear
involution destroys the Hamiltonian structure of this integrable
system. This is because the natural symplectic form vanishes on the
fixed point set of this involution
(compare with Remark~\ref{trivial symplectic form}).
For this reason we shall consider the unreduced system.
Moreover, due to the generalization of the
global Weierstra{\ss} representation to immersions into
$\mathbb{R}^4$, the corresponding \Em{complex Fermi curves} are
invariant only under the anti--holomorphic involution,
but in general not under the holomorphic involution.
In \cite{GS2} it was proven that the
\De{Weierstra{\ss} curves}~\ref{Weierstrass curves} are invariant under
conformal transformations of $\mathbb{R}^3$
acting on the corresponding immersions.
Finally, the Willmore functional is the
\Em{first integral}\index{integral!first $\sim$ $\willmore$}
of this integrable system.

In a second step we will investigate the
\De{Periodicity condition}~\ref{periodicity condition}.
I.\ Taimanov conjectured at the
beginning of this project that these
\De{Weierstra{\ss} curves}~\ref{Weierstrass curves} are
conformal invariants and he confirmed this numerically.
Usually, the real parts of the \Em{isospectral sets}
(i.\ e.\ all potentials fulfilling some reality condition,
whose spectral curves are equal to some given spectral curve)
are compact tori. On the other hand, the set of all
potentials corresponding to all conformal transformations of a given
immersion is in general not compact. For this reason the subset of
some \Em{isospectral set} containing all potentials fulfilling the
\De{Periodicity condition}~\ref{periodicity condition}
seems not to be closed.
But he proved that the isospectral flows preserve this subset \cite{Ta1,Ta2}.
This apparent inconsistency was the first hint that all
\De{Weierstra{\ss} curves}~\ref{Weierstrass curves} have some
non--trivial singularity.
In fact, in \cite{Sch} it was shown that for singular spectral curves
the real parts of the \Em{isospectral sets} are compact, but
decompose into a union of non--closed subsets of different
dimensions which are invariant under the isospectral flows.
We should remark that all \Em{complex Fermi curves} corresponding to
Dirac operators with non--trivial kernel have some kind of
singularity. In fact, each character corresponding to some spin
structure (which takes values in$\{1,-1\}$) is a fixed point of
the anti--holomorphic involution. Since this involution may not
have fixed points on the normalization of the \Em{complex Fermi curve},
the preimage of these elements of the \Em{complex Fermi curves}
under the normalization map has to contain at least two elements.
But in general these singularities do not give rise
to non compact components of the real part of the
\Em{isospectral sets}.

Hence the characterization of
\De{Weierstra{\ss} potentials}~\ref{Weierstrass potentials}
refers to the characterization of
\De{Weierstra{\ss} curves}~\ref{Weierstrass curves}.
But we emphasize that not all potentials,
whose \Em{complex Fermi curves} are
\De{Weierstra{\ss} curves}~\ref{Weierstrass curves},
are \De{Weierstra{\ss} potentials}~\ref{Weierstrass potentials}.
The \De{Weierstra{\ss} curves}~\ref{Weierstrass curves}
may be  characterized as those \Em{complex Fermi curves},
which obey the
\index{singularity!condition}
\index{condition!singularity $\sim$}
\De{Singularity condition}~\ref{singularity condition}.
This \De{Singularity condition}~\ref{singularity condition}
assumes, roughly speaking, that the
\Em{complex Fermi curve} has some fourth--order multiple point.
The exact form of this condition involves the first derivatives
of the regular functions. There exists a natural simplified
and weaker condition called

\newtheorem{Weak Singularity condition}[Lemma]{Weak Singularity condition}
\index{singularity!condition!weak $\sim$}
\index{condition!weak singularity $\sim$}
\begin{Weak Singularity condition}\label{weak singularity condition 1}
  (See also \ref{weak singularity condition}.)
  The \Em{complex Fermi curve} contains
  a singular point corresponding to characters with values in
  $\{1,-1\}$, whose preimage in the normalization contains either at
  least four points, or two points, which are zeroes of both
  quasi--momenta (i.\ e.\ the logarithmic derivatives of the
  values of the characters at two generators of the fundamental group
  of $\mathbb{T}^2$).
\end{Weak Singularity condition}

For immersions into $\mathbb{R}^4$ `quaternionic function theory'
\index{quaternionic!function theory|(}
developed by F.\ Pedit and U.\ Pinkall \cite{PP,BFLPP,FLPP} leads to the

\newtheorem{Quaternionic Singularity condition}[Lemma]{Quaternionic
  Singularity condition}
\index{singularity!condition!quaternionic $\sim$}
\index{condition!quaternionic singularity $\sim$}
\index{quaternionic!singularity condition}
\begin{Quaternionic Singularity condition}\label{quaternionic condition}
  The \Em{complex Fermi curve}
  contains a singular point corresponding to real characters
  (i.\ e.\ quaternionic line bundles), whose preimage in the
  normalization contains either at least four different points, or at
  least two different points, at which both quasi--momenta vanish.
\end{Quaternionic Singularity condition}

This \De{Quaternionic Singularity condition}~\ref{quaternionic condition}
is equivalent to the condition that the corresponding
\Em{isospectral set} contains some complex potential,
whose Dirac operator acting on some quaternionic line bundle
(i.\ e.\ the trivial vector bundle of rank two on $\mathbb{T}^2$
tensored with the line bundle induced by some real character of the
fundamental group of $\mathbb{T}^2$) has a kernel
of quaternionic dimension larger than one.
In general the quaternionic line bundles
may have non--vanishing degree. In this case the action of the
Dirac operator has to be modified. But we shall see that
the notion of \Em{complex Fermi curves} generalizes to these cases,
and the construction of the \Em{Baker--Akhiezer function}
shows that at least in the finite genus case the
\De{Quaternionic Singularity condition}~\ref{quaternionic condition}
is equivalent to the condition that the
\Em{isospectral set} contains some
Dirac operator acting on some quaternionic line bundle with at least
two--dimensional kernel.

The \De{Weak Singularity condition}~\ref{weak singularity condition 1}
is equivalent to the
condition that some \Em{complex Fermi curve} endowed with a
holomorphic and an anti--holomorphic involution obeys the
\De{Quaternionic Singularity condition}~\ref{quaternionic condition},
and that the singularity
is invariant under the holomorphic involution.
In particular, all \Em{complex Fermi curves}
obeying the
\De{Weak Singularity condition}~\ref{weak singularity condition 1}
correspond to some not necessarily injective
conformal mappings\index{conformal!mapping}
from $\mathbb{T}^2$ into $\mathbb{R}^4$.

The characterization of the
\De{Weierstra{\ss} potentials}~\ref{Weierstrass potentials}
in terms of their \Em{complex Fermi curves} suggests to consider
the slightly larger set of all potentials,
whose \Em{complex Fermi curves} satisfy the
\De{Singularity condition}~\ref{singularity condition}.
Instead of this set we now introduce the even larger set of all

\newtheorem{Generalized Weierstrass potentials}[Lemma]{Generalized
Weierstra{\ss} potentials}
\index{Weierstra{\ss}!generalized $\sim$ potential}
\index{generalized!Weierstra{\ss}!potential}
\begin{Generalized Weierstrass potentials}\label{generalized potentials}
Real potentials $U\in \banach{2}(\mathbb{T}^2)$, whose
\Em{complex Fermi curves} fulfill the
\De{Weak Singularity condition}~\ref{weak singularity condition 1},
are called \De{Generalized Weierstra{\ss} potentials}.
\end{Generalized Weierstrass potentials}

The set of all \De{Generalized Weierstra{\ss} potentials} is a
subvariety of the Hilbert space of all $\banach{2}$--potentials.
More precisely, this set is locally the zero set of at most
ten analytic functions. Unfortunately, these functions are not
weakly continuous. But the restriction of these functions to small
balls are weakly continuous. Indeed, we shall explicitly see
that these functions are not weakly continuous on balls,
whose radius correspond to Willmore functionals larger than $4\pi$.
More precisely, the moduli space\index{moduli space}
of all \Em{complex Fermi curves} is not complete.
There exist limits of \Em{complex Fermi curves},
whose \Em{first integrals} are equal to $4\pi$, but which
correspond to no real (or complex) potential.
Nevertheless we shall call these limits \Em{complex Fermi curves},
since they are the \Em{complex Fermi curves} of
\Em{finite rank perturbations} of the Dirac operators.
Moreover, they may be considered as \Em{complex Fermi curves}
of Dirac operators acting on line bundles of non--vanishing degree,
which occur in `quaternionic function theory' \cite{PP,BFLPP,FLPP}.

For this reason we transform in a last step the variational
problem of the $\banach{2}$--norm on the space of all
\De{Generalized Weierstra{\ss} potentials}\ref{generalized potentials}
to a variational problem of the natural extension of the
\Em{first integral} to all\vspace{.4cm}

\noindent
{\bf Generalized Weierstra{\ss} curves~\ref{generalized curves}.}
{\em Roughly speaking, these are those
\Em{complex Fermi curves} of \Em{finite rank perturbations}
of Dirac operators with real potentials $U\in\banach{2}(\mathbb{T}^2)$,
which obey the
\De{Weak Singularity condition}~\ref{weak singularity condition 1}.}
\vspace{.4cm}

To sum up we have the following chain of modifications of
variational problems starting with the Willmore functional:
\begin{description}
\item[1. Immersions.] On the space of smooth immersions from
  $\mathbb{T}^2$ into $\mathbb{R}^3$ we are given the
  Willmore functional. The space is an open subset of the vector space
  of smooth mappings from $\mathbb{T}^2$ into $\mathbb{R}^3$.
  But the functional is quite complicated involving second--order
  partial derivatives.
\item[2. \De{Weierstra{\ss} potentials}.]\index{Weierstra{\ss}!potential}
  Consider the subspace of all $\banach{2}$--potentials,
  whose Dirac operators contain in their kernel some spinor
  fulfilling the
  \De{Periodicity condition}~\ref{periodicity condition},
  and on this subspace the functional, which is the square of the
  $\banach{2}$--norm times four.
  This space is quite complicated, but the functional is very simple.
\item[3. \De{Generalized Weierstra{\ss} potentials}.]
  \index{Weierstra{\ss}!generalized $\sim$ potential}
  \index{generalized!Weierstra{\ss}!potential}
  Consider the real subspace of all $\banach{2}$--potentials,
  whose \Em{complex Fermi curves} fulfill the
  \De{Weak Singularity condition}~\ref{weak singularity condition 1},
  and on this subspace again the functional,
  which is the square of the $\banach{2}$--norm times four.
  The space is a real subvariety of $\banach{2}$.
  The functional is again very simple.
  The corresponding infimum yields only a lower bound
  for the Willmore functional.
\item[4. \De{Generalized Weierstra{\ss} curves}.]
  \index{Weierstra{\ss}!generalized $\sim$ curve}
  \index{generalized!Weierstra{\ss}!curve}
  Consider the subspace of the completion of the moduli space
  containing all \Em{complex Fermi curves} fulfilling the
  \De{Weak Singularity condition}~\ref{weak singularity condition 1},
  and the functional, which is the restriction of the
  \Em{first integral}. The corresponding infimum yields only a lower
  bound of the Willmore functional.
\end{description}
In a concluding remark we emphasize the motivation of
the modification of the second problem into the third problem.
From the point of view of `quaternionic function theory'
\index{quaternionic!function theory|)}
it seems to be more natural to consider the variational problem of
the Willmore functional on immersions into $\mathbb{R}^4$ instead of
$\mathbb{R}^3$. Comparing these two problems the latter has the
advantage of providing a smaller space, and the former has the
advantage of a  stronger condition on relative minimizers.
In some sense the
\De{Weak Singularity condition}~\ref{weak singularity condition 1}
uses this advantage of immersions into $\mathbb{R}^4$
without giving up the advantage of immersions into $\mathbb{R}^3$.
In fact, we shall provide some conditions on
relative minimizers of the fourth variational problem,
which singles out many Willmore tori
(in particular all Willmore tori described in \cite{BB}).

\subsection{The solution of the transformed problem}
\label{subsection solution}

Finally, we arrive at a variational problem of the \Em{first integral}
of an integrable system related to the two--dimensional Dirac operator
with constraints on the corresponding spectral curves.
Therefore, our main effort aims at investigating the moduli space,
i.\ e.\ the space of all generalized \Em{complex Fermi curves}.
Roughly speaking, we are interested in two properties:
\begin{description}
\item[Compactness] of the space of all \Em{complex Fermi curves},
  whose \Em{first integral} is bounded from above, and a
\item[Differentiable structure]\index{structure!differentiable $\sim$}
  of the moduli space.
\end{description}
The first property contains the essential part of the
\Em{Proof of the existence of a minimizer}.
The second property provides the main tools for the
\Em{Classification of the lowest critical points}.
Both properties are proven in full generality.

For the \Em{Compactness} we have to endow this moduli space
with some topology. We shall use the topology of the Hausdorff metric
on the space of all closed subsets of the one--point--compactification
of $\mathbb{C}^2$. With this topology the corresponding space of
closed subsets of the one--point--compactification of $\mathbb{C}^2$
becomes a compact Hausdorff space.
So in order to obtain compactness of some subset of the moduli space,
we only have to prove that this subset is closed.
First we shall prove this for the subsets of the moduli space
containing all generalized \Em{complex Fermi curves}
of bounded geometric genus and with bounded \Em{first integral}.
The main technical tools are a representation of the
\Em{complex Fermi curves} as two complex planes glued along
several cuts and the evaluation of the \Em{first integral} by some
Dirichlet integral of some linear combination of the quasi--momenta over
one of these planes. These tools should be considered as generalizations
of analogous tools for the KdV equations, which were
first used by V.\ A.\ Marchenko and I.\ V.\ Ostrowski \cite{MO} in their
investigation of the moduli space corresponding to real periodic
solutions of the KdV equation with the help of
conformal mappings\index{conformal!mapping}.
In the present situation we encounter two additional technical
complications.
\begin{description}
\item[1. Two--dimensionality of the system:] Integrable systems having
  Lax operators may be divided into three groups.
  \begin{description}
  \item[Two trivial flows.] These systems are finite--dimensional. The
    Lax operator is a spectral-parameter dependent matrix, and
    the spectral curves have two meromorphic functions with prescribed
    poles.
  \item[One trivial and one periodic flow.] The phase space of these
    systems is a space of functions depending on one variable. The
    Lax operator is some ordinary differential operator. The
    associated spectral curves have one meromorphic function and one
    transcendental function with prescribed singularities.
  \item[Two periodic flows.] The phase space of these systems is a
    space of functions depending on two variables. The Lax operator is
    some partial differential operator and the associated
    spectral curves have two transcendental functions with
    prescribed singularities.
  \end{description}
  The KdV and the NLS equation belong to the second group. In
  particular the corresponding spectral curves are hyperelliptic. In
  this case there exist natural parameters of the moduli space, which
  were introduced in \cite{MO} and studied further in \cite{GS1,Tk}.
  The integrable system corresponding to the two--dimensional
  Dirac operator belongs to the third group. In this case
  natural parameters of the moduli space do not exist. We overcome
  this problem by considering the moduli space as a closed subspace of
  some higher--dimensional space, which has natural parameters.
\item[2. A complicated reality condition:] Some reality conditions
  imply that the corresponding spectral curves are \Em{M--curves}
  \index{M--curve}
  (i.\ e.\ the real part has $g+1$ connected components, where $g$ is
  the geometric genus of the spectral curve).
  In this case the real cycles
  cut the spectral curves into pieces, which may be described by some
  conformal mappings. The KdV equation has this property, but the
  integrable system corresponding to the two--dimensional
  Dirac operator does not. It might also be that the
  \Em{complex Fermi curves} can be divided into pieces,
  which can be described by conformal mappings,
  but there does not exist a natural choice.
  This difficulty is overcome by dividing the
  \Em{complex Fermi curve} again into two pieces, where we make use of
  holomorphic functions instead of univalent (or schlicht) functions.
\end{description}

In a second and independent step we determine the limits of the
resolvents of a weakly convergent sequence of potentials.
It is contained in the technically most involved
Section~\ref{subsection limits} and is motivated by two results
on the spectral theory of Dirac operators in
Section~\ref{subsection resolvent} and on the compactifications
of isospectral sets of \Em{complex Fermi curves}
of finite \Em{geometric genus} in
Section~\ref{subsubsection compactified isospectral}.
The first result states that the resolvents are almost weakly continuous
with respect to the potentials. More precisely,
if for any $\varepsilon>0$ the $\banach{2}$--norms of the restrictions
to all $\varepsilon$--balls of a weakly convergent sequence
of potentials is bounded by the inverse Sobolev constant,
then the limit of the corresponding sequence of resolvents
is the resolvent of the limit. Therefore, in general the limits of the
resolvents of a weakly convergent sequence of potentials $U_n$
should differ from the resolvent of the limit only in those points
where the limit of the corresponding measures $U_n(x)\Bar{U}_n(x)d^2x$
has a point measure whose mass is not smaller than
the inverse of the square of the Sobolev constant.
The second result states that the compactified isospectral sets of
\Em{complex Fermi curves} of finite \Em{geometric genus} contain
potentials on the complement of a finite set of the torus
(i.\ e.\ the corresponding \Em{Baker Akhiezer functions}
are not eigenfunctions in small neighbourhoods of this finite set).
The corresponding \Em{first integrals} are equal to four times the
square of the $\banach{2}$--norms of these potentials
plus a multiple of $4\pi$. In particular, the isospectral sets
are compact, if the \Em{first integral} is smaller than $4\pi$.
This follows also from the first result with the optimal Sobolev constant
(compare with Remark~\ref{sobolev 1}).

Therefore, it seems natural to expect that in general
the limits of resolvents of weakly convergent sequences of potentials
are resolvents of perturbations of the Dirac operator corresponding
to the weak limit. More precisely, the support of the perturbation is
contained in a finite subset of $\mathbb{T}^2$ and therefore should be
described by a finite linear combination of such distributions.
This is proven by some kind of blowing up of the sequences
of potentials nearby the singular points to potentials
on $\mathbb{P}^1$. These finite rank perturbations of Dirac operators
can be considered as holomorphic structures of
quaternionic line bundles of non--vanishing degree,
which occur in `quaternionic function theory' \cite{PP,BFLPP,FLPP}.

As a consequence, we prove in
Corollary~\ref{existence of constrained minimizers}
the existence of minimizers of the restrictions of
the Willmore functional to all conformal classes.
We remark that the application of the Heine--Borel property
\cite[Chapter~I. 1.1~Theorem]{Str}
in this variational problem is possible,
since the transformation to a variational problem on the
\Em{complex Fermi curves} throws out the non--compact conformal symmetry.
Furthermore, the Willmore functional was extended to the
maximal domain of definition.

The main techniques for establishing the \Em{Differentiable structure}
is an application of deformation theory of compact complex spaces
combined with some general features of integrable systems having
Lax operators. The spectral curves of these integrable systems should
be considered as a combination of all integrals.
Moreover, the differentials of these integrals fit together
to some linear map from the tangent space of the phase space
into the space of holomorphic one--forms
of the corresponding spectral curve.
On the other hand, the real part of the connected component
of the identity of the Picard group
\index{Picard group}
acts freely and transitively on the \Em{isospectral sets}
\index{isospectral!set}
(i.\ e.\ the set of those elements of the phase space,
whose spectral curves are equal to some given spectral curve).
The differential of this action yields some mapping
from the Lie algebra of the Picard group into the tangent space.
\index{tangent!space of the phase space}
On compact Riemann surfaces S\'{e}rre duality \index{S\'{e}rre duality}
establishes a canonical
pairing between the Lie algebra of the Picard group and the space of
holomorphic one--forms. On the other hand, the symplectic form gives a
pairing on the tangent space of the phase space. Due to a
general rule for such integrable systems these two mappings are dual
to each other with respect to the pairings given by the symplectic
form and S\'{e}rre duality. As a starting point we establish this
structure for the integrable system, whose Lax operator is the
two--dimensional Dirac operator. Since both pairings (symplectic form
and S\'{e}rre duality) are non--degenerate, we conclude that the
tangent space of the moduli space may be identified with the space of
regular forms of the spectral curves. Afterwards we apply a theorem of
Grauert, which is the base of the generalization of
the deformation theory of compact complex manifolds initiated by
K.\ Kodaira and D.\ C.\ Spencer \cite{Ko} to compact complex spaces
(see \cite[Chapter~III. \S4.]{GPR}). This application shows that some
special subsets of the moduli space containing only
\Em{complex Fermi curves} of bounded geometric genus are
complex manifolds. All \Em{complex Fermi curves}
of finite geometric genus are contained in some of these subsets.

In a second step we extend this deformation theory to all
\Em{complex Fermi curves}. For this purpose we show that all
\Em{complex Fermi curves} are arbitrary small perturbations of
\Em{complex Fermi curves} of finite \Em{geometric genus},
and as a consequence we may extend the deformation theory
to these perturbations.
The essential tool of this step is the \Em{asymptotic analysis}
of \Em{complex Fermi curves}.
Usually this is done only for a restricted class of potentials.
More precisely, since the mapping from the potentials to the
spectral curves (i.\ e.\ the combined integrals of motion)
is a non--linear perturbation of the mapping
to the absolute values of the Fourier coefficients,
differentiability assumptions on the potentials
lead to decreasing Fourier coefficients and a simplification
of the asymptotic analysis (compare with \cite{FKT}).
We shall establish the asymptotic analysis in
Theorem~\ref{asymptotic analysis 1} and
Proposition~\ref{asymptotic analysis 2}
for all $\banach{2}$--potentials.
As a consequence in Theorem~\ref{l1-structure}
the moduli space is shown to be a Banach manifold modelled on
Banach space $\ell_1$. For short we use the notation $\ell_1$--manifold.
The corresponding topology differs from the topology used for the
\Em{Compactness}.

The \Em{Differentiable structure} allows us to prove
that all relative minimizers of the fourth variational problem
have dividing real parts
(i.\ e.\ the real part divides the normalization of the
\Em{complex Fermi curve} into two connected components).
Here we make use of local one--dimensional deformation families,
which for negative deformation parameters violate
the reality condition. More precisely, the corresponding
anti--holomorphic involutions have fixed points
on the normalization of the \Em{complex Fermi curves}.
Therefore, the corresponding first--order derivative
should correspond to a higher--order derivative of deformations of
immersions. For relative minimizers the derivative of the
\Em{first integral} with respect to the deformation parameter has to
be non--negative. This condition can be used with the help of a
calculation of intersection numbers of the real part to show that
either this real part divides the normalization or it has self
intersecting points. In the second case the corresponding double
points can again be used to construct local deformations with decreasing
\Em{first integral}. To sum up, the dividing property of the real part
is a consequence of the positivity of some derivative of
even--order, and is definitely not fulfilled at all critical points.
We remark that this dividing property of relative minimizers
is not true for the original Willmore functional.
Roughly speaking, the simplification of the
\De{Weak Singularity condition}~\ref{weak singularity condition 1}
results in a smaller set of relative minimizers
comparable to the smaller set of stereographic projections
of minimal tori in $S^3$.

Finally, we investigate the space of all
\De{Generalized Weierstra{\ss} curves}~\ref{generalized curves}
with dividing real parts and with the \Em{first integral} being
smaller than $8\pi$. Fortunately, this space is small enough to
allow a classification of all relative minimizers in this space.
Since the \Em{geometric genera} of all \Em{complex Fermi curves}
with dividing real parts are bounded in
terms of the \Em{first integral}, this space is finite--dimensional.
We shall first show that any \Em{complex Fermi curve}
of this space may be deformed within this space
into a \De{Generalized Weierstra{\ss} curve}~\ref{generalized curves},
whose normalization has two connected components
and whose \Em{first integral} is equal to $8\pi$.
One connected component of these
\De{Generalized Weierstra{\ss} curves}~\ref{generalized curves}
has to be the spectral curve of some elliptic solution of the
Kadomtsev--Petviashvili Equation. Moreover, due to an estimate
on the genus in terms of the \Em{first integral},
which follows from a formula of Krichever \cite{Kr2},
we will classify these
\De{Generalized Weierstra{\ss} curves}~\ref{generalized curves}.
Finally, this classification is extended to
the whole deformation family.
As a result for each conformal class we have
for all of the three non--trivial spin structures
exactly one relative minimizer of the fourth variational problem,
whose \Em{first integrals} are not larger than $8\pi$.
This agrees with the identification of the
homotopy invariants with the spin structures and a result of P.\ Li and
S.\ T.\ Yau \cite{LY}, which implies that all immersions corresponding
to the trivial spin structure have a self intersection, and therefore a
Willmore functional larger than $8\pi$
(compare with Remark~\ref{trivial spin structure}). The
\De{Generalized Weierstra{\ss} curves}~\ref{generalized curves} of
these relative minimizers fit together into one family of compact
Riemann surfaces whose geometric genera vary between zero and two.
Due to a special feature of this family
the absolute minimum is attained for a 
\De{Generalized Weierstra{\ss} curve}~\ref{generalized curves}
of geometric genus zero. The
\De{Generalized Weierstra{\ss} curves}~\ref{generalized curves}
of geometric genus zero and one are calculated explicitly, and the
absolute minimum is attained for the \Em{complex Fermi curve} of the
\Em{Clifford torus}.

Finally, we give a short summary of the content of the main sections.
Section~\ref{section bloch} contains the spectral theory
of the Dirac operator with periodic $\banach{2}$--potentials
on a flat torus. First we reproduce in the present situation
several classical result of Bloch theory.
We show in Section~\ref{subsection resolvent}
that these Dirac operators are direct integrals of families of
operators with compact resolvents.
Section~\ref{subsection fermi curve} introduces the
complex Bloch variety and the \Em{complex Fermi curve}
together with the first steps of the asymptotic analysis
of these varieties. After considering several symmetries
and reality conditions in Section~\ref{subsection reductions}
we present in Section~\ref{subsection spectral projections} a
remarkable relation of different representations
of these Dirac operators and the corresponding spectral projections.
The \De{Singularity conditions}
(Section~\ref{subsection weak singularity} as well as
Section~\ref{subsection singularity})
and the generalized Cauchy's integral formula
(Section~\ref{subsubsection finite rank perturbations})
are based on this structure. In the last
Subsection~\ref{subsection complex Fermi curves of finite genus}
we leave for a moment the main road to the proof and apply some
techniques of the `finite gap' theory of integrable systems.
This discussion motivates the complicated theorems in the following
Section~\ref{section moduli}.

Section~\ref{section moduli} essentially contains the proof of the 
\Em{Differentiable structure} of the moduli space and the
\Em{Compactness}. The former property is proven in three steps in
Section~\ref{subsection variation zero energy}--\ref{subsection deformation}
and Section~\ref{subsection compactified moduli}.
In the first step the tangent mapping of the phase space to the
spectral curves and its dual (with respect to the symplectic form and
S\'{e}rre duality) is constructed. In the second step an application
of Deformation theory establishes the \Em{Differentiable structure}
for the moduli spaces with bounded genus. Finally, the asymptotic
analysis extends this to the whole moduli space.
The latter property is proven first for the moduli spaces with
bounded genus (Section~\ref{subsection compactified bounded genus})
and independently in the most involved 
Section~\ref{subsection limits}.

After these preparations the proof of the Willmore conjecture is
contained in Section~\ref{section generalized}.
In Section~\ref{subsection weak singularity} we prove that all
\De{Weierstra{\ss} curves}~\ref{Weierstrass curves}
satisfy the
\De{Weak Singularity condition}~\ref{weak singularity condition 1}.
Then the dividing property of relative minimizers follows in
Section~\ref{subsection relative minimizers}. Finally the
classification of all local minimizers below $8\pi$ is presented in
Section~\ref{subsection absolute minimizer}.

The last Section contains an investigation of the
\De{Singularity condition}~\ref{singularity condition}
(Section~\ref{subsection singularity}) and a proof
that all conformal classes contain a minimum of finite type
(Section~\ref{subsection constrained}).

Whenever possible we have tried to quote new standard textbooks
for the tools we have used. These textbooks contain the references
to the original articles.

\noindent
{\bf Acknoledgements:} This work is dedicated to my teacher Robert
Schrader. His constant support gave me the opportunity to attack such
a difficult problem without the safety of having a permanent position.
In the beginning, Iskander Taimanov participated in this project and
later accompanied it with several discussions.
The collaboration inside the SFB~288
with the mathematicians Alexander Bobenko, Frank Duzaar,
Daniel Grieser, Frank Hausser, Franz Pedit and Ulrich Pinkall
was very stimulating and helpful. 
So results from `quaternionic function theory' provided useful
checks and led to new ideas. During several visits in Berlin Piotr
Grinevich was a good adviser. Elmar Vogt helped to recognize an impasse.
David Jerison pointed out to the author, that the arguments of
\cite[Proposition~2.6]{Wo}, where the analogous but weaker
Carleman inequality for the gradient term of the Laplace operator is
treated, carry over to the Dirac operator.
Sonya Miller, Robert Schrader and Elmar Vogt helped to improved
my English style.

\section{The Bloch theory}\label{section bloch}

Originally the Bloch theory investigates the spectral theory of a
Schr\"odinger operator with a periodic potential.
Bearing in mind that the Dirac operator is a natural relativistic
modification of the Schr\"odinger operator we continue to use this
notion of a Bloch theory in the context of Dirac operators with
periodic potentials. The resolvent of an unbounded operator is a
bounded operator, which contains all relevant information of the
spectral theory of the former operator. Hence we first introduce the
resolvents of the corresponding Dirac operators on the torus.

\subsection{The resolvent of the Dirac operator}\label{subsection resolvent}

Let $\lattice $ be a lattice in $\mathbb{R}^2$
\index{lattice!$\lattice$, $\lattice_\mathbb{C}$}
and $V$ and $W$ two
periodic functions from $\torus $ into $\mathbb{C}$. In
this section we investigate the resolvent of the
Dirac operator
\index{Dirac operator $\Op{D}$}
$$\Op{D}=
\begin{pmatrix}
V & \partial \\
-\Bar{\partial} & W
\end{pmatrix}
\text{, with }
\partial = \frac{1}{ 2}\left(\frac{\partial}{ \partial x_1}
-\sqrt{-1}\frac{\partial }{ \partial x_2}\right),\;
\Bar{\partial} = \frac{1}{ 2}\left(\frac{\partial}{ \partial x_1}
+\sqrt{-1}\frac{\partial }{ \partial x_2}\right).$$
\index{coordinate!$x$, $z$}
By abuse of notation the operator of multiplication with some
function is denoted by the same symbol as the function.
The operators $\partial$ and $\Bar{\partial}$ coincide with
holomorphic and anti--holomorphic derivation with respect to
$z=x_1+\sqrt{-1}x_2$. Therefore, we use sometimes the coordinates
$z$ and $\Bar{z}$ instead of $x$. Consequently $\lattice_\mathbb{C}$
denotes the corresponding sublattice of $\mathbb{C}$:
$$\lattice_\mathbb{C}=
\left\{\gamma_1+\sqrt{-1}\gamma_2\in\mathbb{C}\mid\gamma\in\lattice\right\}.$$ 
We use the natural scalar product on $\mathbb{C}^2$:
\begin{align*}
g(x,x')&=x_1 x'_1 + x_2 x'_2&\forall x,x'\in\mathbb{C}^2.&
\end{align*}
For all $k\in \mathbb{C}^2$ the map
$\gamma\mapsto \exp\left(2\pi\sqrt{-1}g(\gamma,k)\right)$ is a
\index{quasi--momenta!$k$}
one--dimensional representation of $\lattice$. Two of these
representations are equivalent if and only if the difference of the
corresponding elements in $\mathbb{C}^2$ is some element of the
dual lattice\index{lattice!dual $\sim$ $\lattice\dual$}
$$\lattice\dual :=\left\{ \kappa \in \mathbb{R}^2\mid
g(\gamma ,\kappa)\in \mathbb{Z}
\;\forall\gamma \in \lattice \right\} .$$
Hence the set of one--dimensional representations of $\lattice$ may be
identified with $\mathbb{C}^2/\lattice\dual$. All these
representations $[k]\in\mathbb{C}^2/\lattice\dual$ induce some
line bundle on the torus $\torus$, whose cross sections
are functions $\psi$ on $\mathbb{R}^2$ with the property
\begin{align*}
\psi(x+\gamma)&=\exp\left(2\pi\sqrt{-1}g(\gamma,k)\right)\psi (x)
&\text{for all }x \in \mathbb{R}^2\text{ and all }\gamma\in\lattice.&
\end{align*}
In the sequel we will often consider operators, which
act on the space of cross sections of these bundles. These
cross sections may be described in two different ways. We may either
use a Fundamental domain or a Trivialization.

\newtheorem{Fundamental domain}[Lemma]{Fundamental domain}
\index{fundamental domain}
\index{lattice!generator!$\xx{\gamma}$ and $\yy{\gamma}$}
\index{lattice!dual $\sim$ generator!$\xx{\kappa}$ and $\yy{\kappa}$}
\index{coordinate!$\xx{q}$ and $\yy{q}$}
\index{quasi--momenta!$\xx{p}$ and $\yy{p}$}
\begin{Fundamental domain}\label{fundamental domain}
  We chose some fundamental domain
  $$\Delta=\left\{x\in\mathbb{R}^2\mid
  x=x'+\xx{q}_\Delta\xx{\gamma}_\Delta+
  \yy{q}_\Delta\yy{\gamma}_\Delta
  \text{ with }
  (\xx{q}_\Delta,\yy{q}_\Delta)\in [0,1]^2\right\},$$
  of the action of $\lattice$ on $\mathbb{R}^2$. Here $\xx{\gamma}_\Delta$
  and $\yy{\gamma}_\Delta$ are generators of $\lattice$ and
  $x'$ is some element of $\mathbb{R}^2$. In the sequel we will
  sometimes pick out finitely many special elements of
  $\mathbb{R}^2/\lattice$, and investigate the behaviour of some
  functions nearby these elements. In this case we will always assume
  that $x'$ is chosen in such a way that the corresponding elements
  of $\Delta$ are contained in the interior of $\Delta$. The sections
  of the line bundles corresponding to some
  $[k]\in\mathbb{C}^2/\lattice\dual$ may be described by usual
  functions on $\Delta$, which on the boundary of $\Delta$ obey the
  conditions introduced above. Consequently the Hilbert space of
  cross sections of these line bundles are all equal to
  $\banach{2}(\Delta)$. In this case the dual
  generators of the dual lattice are denoted by $\xx{\kappa}_\Delta$ and
  $\yy{\kappa}_\Delta$, with
  $g(\xx{\gamma}_\Delta,\xx{\kappa}_\Delta)=1=
  g(\yy{\gamma}_\Delta,\yy{\kappa}_\Delta)$ and
  $g(\xx{\gamma}_\Delta,\yy{\kappa}_\Delta)=0=
  g(\yy{\gamma}_\Delta,\xx{\kappa}_\Delta)$. Therefore, we have natural
  coordinates of the space $\mathbb{R}^2$ denoted by
  $\xx{q}_\Delta=g(x,\xx{\kappa}_\Delta)$ and
  $\yy{q}_\Delta=g(x,\yy{\kappa}_\Delta)$, and the dual
  space denoted by $\xx{p}_\Delta=g(\xx{\gamma}_\Delta,k)$ and
  $\yy{p}_\Delta=g(\yy{\gamma}_\Delta,k)$.
  Moreover, we assume that the orientation of
  the generators of $\lattice$ is positive, which means that the real
  number $\xx{\gamma}_{\Delta,1}\yy{\gamma}_{\Delta,2}-
  \xx{\gamma}_{\Delta,2}\yy{\gamma}_{\Delta,1}$ is positive and
  therefore equal to the Euclidean volume
  $\vol(\torus)$ of
  the torus. Due to this assumption the modular group
  $SL(2,\mathbb{Z})$ acts freely and transitively on the set of all
  possible choices of such generators. If we want to emphasize the
  dependence of some objects on this choice, we will decorate this
  object with an index $\Delta$, but in general we shall omit this index.
\end{Fundamental domain}

\newtheorem{Trivialization}[Lemma]{Trivialization}
\index{trivialization}
\begin{Trivialization}\label{trivialization}
  The functions $\psi_{k}$ are global
  non--vanishing sections of the line bundles corresponding to
  $[k]$. Hence all cross sections may be described by usual functions
  times these trivializations. However, these Trivializations are different
  for different representatives $k$ and $k+\kappa$ of the class
  $[k]$. Therefore the representatives of the
  cross sections of the line bundles corresponding to
  $[k]\in\mathbb{C}^2/\lattice\dual$
  belong to a vector bundle over
  $\mathbb{C}^2/\lattice\dual$. If for example we consider the
  Hilbert spaces of $\banach{2}$--sections,
  then the cross sections belong to the
  induced Hilbert bundle on the quotient $\mathbb{C}^2/\lattice\dual$, 
  of the unitary representation $\kappa\mapsto \psi_{\kappa}$ of
  $\lattice\dual$ on $\banach{2}(\torus)$. The sections $\psi(k)$ of
  this induced bundle have to obey the relation
  $\psi(k+\kappa)=\psi_{-\kappa}\psi(k)$. If we use
  this approach the elements of this Hilbert bundle are called
  periodic parts of the eigenfunctions, and the corresponding
  operators, which act on this Hilbert bundle, will be decorated by
  some tilde.\index{Hilbert bundle}\index{accent tilde}
\end{Trivialization}

In the sequel the dependence on $k$, which is easier to describe with
the help of the \De{Fundamental domain}, will be
more important than the dependence on $\lattice$, which is easier to
describe with the help of the \De{Trivialization}.
Hence we will mostly use the \De{Fundamental domain}.
Obviously the operator of multiplication with
$\psi_k$ transforms the periodic parts of the eigenfunctions defined
with the help of the \De{Trivialization} into the elements of the
Hilbert space $\banach{2}(\Delta)$ over the \De{Fundamental domain}.
Therefore, the operators on the former Hilbert bundle are transformed
into the operators on the latter Hilbert space by taking the conjugate
with the operator $\psi_k$.

In addition to the hermitian scalar product the Hilbert spaces
$\banach{2}(\Delta)$ and $\banach{2}(\torus)$ have a natural complex
conjugation denoted by $\psi\mapsto\Bar{\psi}$ and a symmetric
bilinear form denoted by
$$\langle\phi,\psi\rangle=\int\limits_{\Delta}
\phi(x)\psi(x)d^2x.$$
Consequently, the hermitian scalar product is denoted by
$\langle\Bar{\phi},\psi\rangle$. 
Let $\partial_{[k]}$ and $\Bar{\partial}_{[k]}$ be
for all $k\in \mathbb{C}^2$ the unique closures of the
unbounded operators $\partial$ and $\Bar{\partial}$
on $\banach{2}(\Delta)$ with the domain of smooth functions $\psi$, so that
\begin{align*}
\psi(x+\gamma)&=\exp\left(2\pi\sqrt{-1}g(\gamma,k)\right)\psi(x)
&\text{for all }x\in\mathbb{R}^2\text{ and all }\gamma\in\lattice.&
\end{align*}
Obviously, these operators depend only on the
representative $[k]\in \mathbb{C}^2/\lattice\dual$.
The eigenfunctions of these operators are given by
$\psi_{k}$ with eigenvalues
\begin{align*}
\partial_{[k]}\psi_{k}&=\pi\left(k_2+\sqrt{-1}k_1\right)\psi_{k}
&\text{and }
\Bar{\partial}_{[k]}\psi_{k}&=
\pi\left(-k_2+\sqrt{-1}k_1\right)\psi_{k},
\end{align*}
where $k$ runs over all elements of the preimage of $[k]$ under the
natural map $\mathbb{C}^2\rightarrow
\mathbb{C}^2/\lattice\dual$. Let
$\Op{D}(V,W,k)$ denote the operator
$$\Op{D}(V,W,k)=\begin{pmatrix}
V & \partial_{[k]}\\
-\Bar{\partial}_{[k]} & W
\end{pmatrix}.$$ The exact definition as a closed operator on some
Banach space is part of the following

\begin{Theorem}\label{meromorph}
\index{Dirac operator $\Op{D}$!resolvent of the $\sim$ $\Op{R}$}
\index{resolvent!of the Dirac operator $\Op{R}$}
The resolvent
$\Op{R}(V,W,k,\lambda)=\left(\lambda\unity-\Op{D}(V,W,k)\right)^{-1}$
defines a meromorphic map from
$\banach{2}(\torus)\times\banach{2}(\torus)\times
\mathbb{C}^2\times\mathbb{C}$ into the set of compact operators on
$\banach{2}(\Delta)\times\banach{2}(\Delta)$. Moreover,
for all $(V_0,W_0,k_0,\lambda_0)\in
\banach{2}(\torus)\times\banach{2}(\torus)\times
\mathbb{C}^2\times\mathbb{C}$
there exists some open neighbourhood
$\Set{U}\subset \banach{2}(\torus)\times
\banach{2}(\torus)\times \mathbb{C}^2$
of $(V_0,W_0,k_0)$ and some holomorphic function $\Breve{\Op{P}}$
from $\Set{U}$ into the finite rank projections on
$\banach{2}(\Delta)\times\banach{2}(\Delta)$, such that
the following equation holds:
$$\left(\unity-\Breve{\Op{P}}(V,W,k)\right)\comp\Op{R}(V,W,k,\lambda)=
\Op{R}(V,W,k,\lambda)\comp\left(\unity-\Breve{\Op{P}}(V,W,k)\right)$$
for all $(V,W,k)\in \Set{U}$ and all $\lambda$ in some small
neighbourhood of $\lambda_0$. Both functions are holomorphic on
this domain. Finally, for all $2\leq q<\infty$ the images of the
values of the map $\Breve{\Op{P}}$ are contained in
$\banach{q}(\Delta)\times\banach{q}(\Delta)$.
\end{Theorem}

\begin{proof} We will prove the equivalent statements for the
  resolvent $\triv{\Op{R}}(V,W,k,\lambda)$ of the operator 
  $\triv{\Op{D}}(V,W,k)=
  \left(\begin{smallmatrix}
    \psi_{-k} & 0 \\
    0 & \psi_{-k}
    \end{smallmatrix}\right) \comp\Op{D}(V,W,k)\comp
    \left(\begin{smallmatrix}
    \psi_k & 0 \\
    0 & \psi_k
    \end{smallmatrix}\right),$
  which is an operator on $\banach{2}(\torus)\times 
  \banach{2}(\torus)$. First we prove two lemmata.

\begin{Lemma} \label{free resolvent}
Let $1<p<2$ and $q=\frac{2p}{2-p}$.
Then there exists some constant $S_p$ depending only on $p$
with the following property:
For any fixed $S'_p>S_p$ and for all $k\in\mathbb{C}^2$
there exists some $\delta>0$ such that for all real $\lambda>\delta$
the operator $\triv{\Op{R}}(0,0,k,\sqrt{-1}\lambda)$
from $\banach{p}(\torus)\times\banach{p}(\torus)$ into
$\banach{q}(\torus)\times\banach{q}(\torus)$ has norm
smaller than $S'_p$.
\end{Lemma}

\begin{Remark}\label{sobolev 1}\index{Sobolev!constant $S_p$}
  The constant $\sqrt{2}S_p$ is the Sobolev constant
  of the embbeding $\sobolev{1,p}(\torus)\hookrightarrow
  \banach{\frac{2p}{2-p}}(\torus)$.
  With the optimal Sobolev constants in
  \cite[Chapter~I. Section~4.6]{Str} we may set
  $S_p=\frac{p}{\sqrt{\pi}(2-p)}$.
\end{Remark}

\begin{proof} Let $\Func{G}_\lambda(x)=\frac{1}{2\pi}
  \Func{K}_0(\lambda\sqrt{g(x,x)})$
  \index{Green's function!$\sim$ $\Func{G}_{\lambda}$}
  be the Green's function of $(\lambda^2-\Delta)^{-1}$ on
  $\mathbb{R}^2$ (\cite[Chapter~V \S3.1 and \S6.5]{St} and
  \cite[Section~7.2]{GJ}), where $\Func{K}_0$
  denotes the modified Bessel function of zero--order. Then the
  resolvent $\left(\begin{smallmatrix}
  \sqrt{-1}\lambda & -\partial \\
  \Bar{\partial} & \sqrt{-1}\lambda
  \end{smallmatrix}\right)^{-1}$ has on $\mathbb{R}^2$ the integral kernel
  $$\begin{pmatrix}
  a_\lambda(x-x') & b_\lambda(x-x')\\
  c_\lambda(x-x') & d_\lambda(x-x')
  \end{pmatrix}d^2x' =-4
  \begin{pmatrix}
  \sqrt{-1}\lambda\Func{G}_{2\lambda}(x-x') & 
  \partial\Func{G}_{2\lambda}(x-x')\\
  -\Bar{\partial}\Func{G}_{2\lambda}(x-x') &
  \sqrt{-1}\lambda\Func{G}_{2\lambda}(x-x')
  \end{pmatrix}d^2x'.$$
The arguments given in the proof of \cite[Proposition~7.2.1]{GJ}
also show that
\begin{description}
\item[(i)] $|\Func{K}_0(|t|)|\leq
  \text{\bf{O}}(1) \exp(-|t|)/|t|$ and
  $|\Func{K}_0'(|t|)|\leq {\bf{O}}(1) \exp(-|t|)/|t|$
  for $t$ bounded away from zero.
\item[(ii)] $\frac{1}{2\pi}\Func{K}_0(|t|)\sim -\ln(|t|)$
  in some neighbourhood of $t=0$. 
\item[(iii)] $\frac{1}{2\pi}\Func{K}_0'(|t|)\sim -1/|t|$
  in some neighbourhood of $t=0$. 
\end{description}
Hence the infinite sum
$$\begin{pmatrix}
a_{\lambda,\lattice}(x) & b_{\lambda,\lattice}(x)\\
c_{\lambda,\lattice}(x) & d_{\lambda,\lattice}(x)
\end{pmatrix} =\sum\limits_{\gamma\in\lattice}
\begin{pmatrix}
a_\lambda(x+\gamma) & b_\lambda(x+\gamma)\\
c_\lambda(x+\gamma) & d_\lambda(x+\gamma)
\end{pmatrix} $$
\index{integral!kernel!of $\triv{\Op{R}}(0,0,0,\sqrt{-1}\lambda)$}
converges on $\mathbb{R}^2$ to the integral kernel of
$\triv{\Op{R}}(0,0,0,\sqrt{-1}\lambda)$. Moreover, due to the
Hardy--Littlewood--Sobolev theorem \cite[Chapter~V. \S1.2 Theorem~1]{St}
the resolvent $\left(\begin{smallmatrix}
\sqrt{-1}\lambda & \partial \\
-\Bar{\partial} & \sqrt{-1}\lambda
\end{smallmatrix}\right)^{-1}$ of the Dirac operator on
$\mathbb{R}^2$ is a bounded operator from
$\banach{p}(\mathbb{R}^2)\times\banach{p}(\mathbb{R}^2)$
into $\banach{q}(\mathbb{R}^2)\times\banach{q}(\mathbb{R}^2)$. Using the
equalities
$$\lambda^{\frac{2}{p}}\left(\int\limits_{\mathbb{R}^2}
|\psi(\lambda x)|^p d^2x\right)^{1/p} = 
\left(\int\limits_{\mathbb{R}^2}
|\psi(x)|^p d^2x\right)^{1/p} \text{ and }$$
$$\begin{pmatrix}
a_\lambda(x-x') & b_\lambda(x-x')\\
c_\lambda(x-x') & d_\lambda(x-x')
\end{pmatrix}d^2x'=
\begin{pmatrix}
a_1(\lambda x-\lambda x') & b_1(\lambda x-\lambda x')\\
c_1(\lambda x-\lambda x') & d_1(\lambda x-\lambda x')
\end{pmatrix}\frac{d^2\lambda x'}{\lambda}$$
we conclude that the norm of this operator does not depend on
$\lambda$. Moreover, this resolvent considered as an operator on
$\banach{p}(\mathbb{R}^2)\times\banach{p}(\mathbb{R}^2)$ or
$\banach{q}(\mathbb{R}^2)\times\banach{q}(\mathbb{R}^2)$
converges to zero, if $\lambda$ tends to infinity.
Due to the definition of the integral kernel
of $\triv{\Op{R}}(0,0,0,\sqrt{-1}\lambda)$ this implies the equality
$$\begin{pmatrix}
a_{\lambda,\lattice}(x) & b_{\lambda,\lattice}(x)\\
c_{\lambda,\lattice}(x) & d_{\lambda,\lattice}(x)
\end{pmatrix}=
\lambda\begin{pmatrix}
a_{1,\lambda\lattice}(\lambda x) & b_{1,\lambda\lattice}(\lambda x)\\
c_{1,\lambda\lattice}(\lambda x) & d_{1,\lambda\lattice}(\lambda x)
\end{pmatrix},$$
where $\lambda\lattice$ denotes the lattice
$\lambda\lattice=\{\lambda\gamma\mid\gamma\in\lattice\}$. Moreover,
for all $S'_p$ larger than the Sobolev constant $S_p$
there exists some $\delta>0$ such that for all real $\lambda>\delta$,
the resolvent $\triv{\Op{R}}(0,0,0,\sqrt{-1}\lambda)$
is an operator from
$\banach{p}(\torus)\times\banach{p}(\torus)$ into
$\banach{q}(\torus)\times\banach{q}(\torus)$,
whose norm is smaller than $S'_p$.
This shows the claim for $k=0$. Finally, this resolvent
$\triv{\Op{R}}(0,0,0,\sqrt{-1}\lambda)$ considered as an operator on
$\banach{p}(\mathbb{R}^2)\times\banach{p}(\mathbb{R}^2)$ or
$\banach{q}(\mathbb{R}^2)\times\banach{q}(\mathbb{R}^2)$
converges to zero, if $\lambda$ tends to infinity.
Obviously the difference
$\triv{\Op{D}}(0,0,k)-\triv{\Op{D}}(0,0,0)$ is the bounded
operator $\left(\begin{smallmatrix}
0 & \pi(k_2+\sqrt{-1}k_1\\
\pi(k_2-\sqrt{-1}k_1) & 0
\end{smallmatrix}\right)$. Hence the Neumann series
\begin{multline*}
\triv{\Op{R}}(0,0,k,\sqrt{-1}\lambda)-
\triv{\Op{R}}(0,0,0,\sqrt{-1}\lambda)=\\
=\sum\limits_{l=1}^\infty
\triv{\Op{R}}(0,0,0,\sqrt{-1}\lambda)
\left(\comp\left(\begin{smallmatrix}
0 & \pi(k_2+\sqrt{-1}k_1\\
\pi(k_2-\sqrt{-1}k_1) & 0
\end{smallmatrix}\right)\comp
\triv{\Op{R}}(0,0,0,\sqrt{-1}\lambda)\right)^l
\end{multline*}
of operators from
$\banach{p}(\torus)\times\banach{p}(\torus)$ into
$\banach{q}(\torus)\times\banach{q}(\torus)$
converges to zero for $\lambda\rightarrow\infty$.
\end{proof}

In the following lemma we make use of a family of norms on
$\banach{2}(\torus)$. 
The family is parameterized by a positive real number
$\varepsilon$. We consider
the balls $B(x,\varepsilon)$ in $\torus$, which for large
$\varepsilon$ becomes equal to the whole torus independent of $x$. But we
are mainly interested in small $\varepsilon$. More precisely,
if $2\varepsilon$ is smaller than the lengths of all non--zero periods,
the corresponding balls have volume $\pi\varepsilon^2$.
For all $U\in \banach{2}(\torus)$, the function
$$x\mapsto \left\|\left.U\right|_{B(x,\varepsilon)}\right\|_2$$
is a continuous function on $\torus$. Obviously, the
supremum of this function in dependence of $U$, denoted by
$\|U\|_{\varepsilon,2}$, defines a norm on $\banach{2}(\torus)$,
which is equivalent to the usual $\banach{2}$--norm.
In fact, Fubini's theorem \cite[Theorem~I.22]{RS1}
and the invariance of the Lebegues measure under translations
\cite[Chapter~14 Section~6 24.~Proposition]{Ro2} yields the identity
\begin{multline*}
\begin{aligned}
\int\limits_{\torus}\int\limits_{B(x,\varepsilon)}f(x')d^2x'd^2x&=
\int\limits_{\torus}\int\limits_{B(0,\varepsilon)}f(x+x')d^2x'd^2x&=&\\
\int\limits_{B(0,\varepsilon)}\int\limits_{\torus}f(x+x')d^2xd^2x'&=
\int\limits_{B(0,\varepsilon)}\int\limits_{\torus}f(x)d^2xd^2x'&=&
\end{aligned}\\
=\vol(B(0,\varepsilon))\int\limits_{\torus}f(x)d^2x
\;\forall f\in\banach{1}(\torus).
\end{multline*}
For $f(x)$ equal to the square of the former function we obtain
that the usual $\banach{2}$--norm of the former function
is equal to $\sqrt{\vol(B(x',\varepsilon))}$ times $\|U\|_2$. 
Hence we have
$$\sqrt{\frac{\vol(B(x',\varepsilon))}
{\vol(\torus)}}
\|U\|_2\leq\|U\|_{\varepsilon,2}\leq\|U\|_2.$$
Due to \cite[Appendix to Chapter~V. Theorem~1.]{Yo}
the convex (weakly closed) subsets
$$\left\{U\in \banach{2}(\torus)\mid
\|U\|_{\varepsilon,2}\leq C\right\}$$
of $\banach{2}(\torus)$ are compact with respect to the weak topology
for all $\varepsilon >0$ and $C>0$. Moreover, for any
$U\in \banach{2}(\torus)$ the limit
$\lim\limits_{\varepsilon\downarrow 0} \|U\|_{\varepsilon,2}$ is zero, but
not uniformly on bounded sets of $\banach{2}(\torus)$. Hence
for all $U$ and all $C >0$ there exists some small $\varepsilon >0$, such
that $\|U\|_{\varepsilon,2}\leq C$.

\begin{Lemma} \label{weakly continuous resolvent}
\index{resolvent!weak continuity of the $\sim$}
All $C_p<S_p^{-1}$
(compare with Lemma~\ref{free resolvent} and Remark~\ref{sobolev 1})
have the following property:
For all $\varepsilon>0$ and open bounded subsets
$\Set{O}$ of $\mathbb{C}^2$ there exists some $\delta>0$ such
that for all $\lambda>\delta$ the mapping
$(V,W,k)\mapsto \Op{R}(V,W,k,\sqrt{-1}\lambda)$
is a holomorphic weakly continuous map from
$\left\{(V,W)\in \banach{2}(\torus)\times
\banach{2}(\torus)\mid\|V\|_{\varepsilon,2}\leq C_p\text{ and }
\|W\|_{\varepsilon,2}\leq C_p\right\}\times\Set{O}$ into the compact
operators on $\banach{2}(\Delta)\times
\banach{2}(\Delta)$.
\end{Lemma}

\begin{proof} We decompose the free resolvent $\Op{R}(0,0,k,\lambda)=
\Op{R}_{\varepsilon'\text{\scriptsize\rm--near}}(0,0,k,\lambda)+
\Op{R}_{\varepsilon'\text{\scriptsize\rm--distant}}(0,0,k,\lambda)$
into a sum of two operators, where the
integral kernel of the first summand is the product of the
integral kernel of $\Op{R}(0,0,k,\lambda)$,
with the function depending on $x-x'$,
which is equal to zero, if $\|x-x'\|>\varepsilon'$, and otherwise
equal to $1$. Consequently the integral kernel of the second summand
is the product of the integral kernel of $\Op{R}(0,0,k,\lambda)$
with the function depending on $x-x'$, which is equal to zero,
if $\|x-x'\|\leq \varepsilon'$, and otherwise equal to $1$.
Due to H\"older's inequality \cite[Theorem~III.1~(c)]{RS1}
the operators of multiplication with $V$ and $W$
are bounded operators from
$\banach{q}(\torus)$ into $\banach{p}(\torus)$.
Now Lemma~\ref{free resolvent} implies
that all fixed $C_p<S_p^{-1}$ and all $(V,W,k)\in
\left\{(V,W)\in \banach{2}(\torus)\times
\banach{2}(\torus)\mid
\|V\|_{\varepsilon,2}< C_p\text{ and }
\|W\|_{\varepsilon,2}< C_p\right\}\times\Set{O}$ and all large real
$\lambda$ the norm of the operator
$\left(\begin{smallmatrix}
V & 0\\
0 & W
\end{smallmatrix}\right)\comp 
\Op{R}_{\varepsilon'\text{\scriptsize\rm--near}}(0,0,k,\sqrt{-1}\lambda)$
on $\banach{p}(\torus)\times\banach{p}(\torus)$
is for small $0<\varepsilon'<\varepsilon$
bounded by some constant smaller than $\left(
\frac{\vol(B(0,\varepsilon))}{\vol(B(0,\varepsilon-\varepsilon'))}
\right)^{1/p}=\left(1-\frac{\varepsilon'}{\varepsilon}\right)^{-2/p}$.
In fact, for all $x\in\torus$ the
$\banach{p}(\torus)\times\banach{p}(\torus)$--norm
of the restriction of $\left(\begin{smallmatrix}
V & 0\\
0 & W
\end{smallmatrix}\right)\comp 
\Op{R}_{\varepsilon'\text{\scriptsize\rm--near}}(0,0,k,\sqrt{-1}\lambda)
\psi$ to $B(x,\varepsilon-\varepsilon')$ is smaller than the norm of the
restriction of $\psi$ to $B(x,\varepsilon)$.
Therefore, the bound on
the operator $\left(\begin{smallmatrix}
V & 0\\
0 & W
\end{smallmatrix}\right)\comp 
\Op{R}_{\varepsilon'\text{\scriptsize\rm--near}}(0,0,k,\sqrt{-1}\lambda)$
follows from the identity
\begin{align*}
\int\limits_{\torus}\left\|\left.f\right|_{B(x,\varepsilon)}\right\|_p^pd^2x&=
\vol(B(0,\varepsilon))\|f\|_p^p&
\;\forall f\in\banach{p}(\torus),&
\end{align*}
which follows again from Fubini's theorem \cite[Theorem~I.22]{RS1}
and the invariance of the Lebegues measure under translations
\cite[Chapter~14 Section~6 24.~Proposition]{Ro2}.
With very small $\varepsilon'$ the norm of this operator
$\left(\begin{smallmatrix}
V & 0\\
0 & W
\end{smallmatrix}\right)\comp 
\Op{R}_{\varepsilon'\text{\scriptsize\rm--near}}(0,0,k,\sqrt{-1}\lambda)$
on $\banach{p}(\torus)\times\banach{p}(\torus)$
is smaller than one.
On the other hand we have seen in the proof of Lemma~\ref{free resolvent}
that under the natural isometries of $\banach{p}(\torus)$ and
$\banach{q}(\torus)$ onto $\banach{p}(\mathbb{R}^2/\lambda\lattice)$
and $\banach{q}(\mathbb{R}^2/\lambda\lattice)$, respectively, the resolvent
$\Op{R}(0,0,k,\sqrt{-1})$ considered as an operator from
$\banach{p}(\torus)\times\banach{p}(\torus)$ to
$\banach{q}(\torus)\times\banach{q}(\torus)$ is
transformed into the resolvent $\Op{R}(0,0,k,\sqrt{-1}\lambda)$
considered as an operator from
$\banach{p}(\mathbb{R}^2/\lambda\lattice)\times
\banach{p}(\mathbb{R}^2/\lambda\lattice)$ 
to
$\banach{q}(\mathbb{R}^2/\lambda\lattice)\times
\banach{q}(\mathbb{R}^2/\lambda\lattice)$.
Moreover, the asymptotic behaviour (i) of
the Green's function of the Laplace operator mentioned in the proof of
Lemma~\ref{free resolvent} implies that for large $\lambda$, the
integral kernel of this resolvent becomes arbitrary small, if the
distance of $x-x'$ in $\mathbb{R}^2/\lambda\lattice$ becomes
large. Hence for all $(V,W,k)\in
\left\{(V,W)\in \banach{2}(\torus)\times
\banach{2}(\torus)\mid
\|V\|_{\varepsilon,2}< C_p\text{ and }
\|W\|_{\varepsilon,2}< C_p\right\}\times\Set{O}$
the norm of the operator
$\left(\begin{smallmatrix}
V & 0\\
0 & W
\end{smallmatrix}\right)\comp
\Op{R}_{\varepsilon'\text{\scriptsize\rm--distant}}(0,0,k,\sqrt{-1}\lambda)$
on $\banach{p}(\torus)\times\banach{p}(\torus)$
converges in the limit $\lambda\rightarrow\infty$ uniformly to zero.
Consequently the Neumann series
$$\triv{\Op{R}}(V,W,k,\sqrt{-1}\lambda)
=\sum\limits_{l=0}^{\infty} \triv{\Op{R}}(0,0,k,\sqrt{-1}\lambda)
\left(\comp\left(\begin{smallmatrix}
V & 0\\
0 & W
\end{smallmatrix}\right)\comp
\triv{\Op{R}}(0,0,k,\sqrt{-1}\lambda)\right)^l$$
converges to a holomorphic function with values
in the bounded operators from
$\banach{p}(\torus)\times\banach{p}(\torus)$ to
$\banach{q}(\torus)\times\banach{q}(\torus)$.
The latter Banach space is naturally contained in the former
Banach space. Therefore, it remains to show, firstly that the values of
this function are compact operators on the former Banach space, and
secondly that these functions are weakly continuous.

In doing so we consider again the free resolvent
$\triv{\Op{R}}(0,0,k,\sqrt{-1}\lambda)$ as an operator from
$\banach{p}(\torus)\times\banach{p}(\torus)$ to
$\banach{q}(\torus)\times\banach{q}(\torus)$, but this
time we assume the inequality $q<\frac{2p}{2-p}$, and therefore
$1/p=1/q+1/r$ with some $r>2$.
  Due to the Hausdorff--Young inequality
  \cite[Theorem~IX.8]{RS2} the (not normalized) Fourier transform
  \index{Fourier transform $\Op{F}$}
  $$\Op{F}:
  \psi\mapsto
  \left(\fourier{\psi}(\kappa)\right)_{\kappa\in\lattice\dual},
  \text{ with } \psi(x)=\sum\limits_{\kappa\in\lattice\dual}
  \fourier{\psi}(\kappa)\exp(2\pi\sqrt{-1}g(x,\kappa)).$$
  is a bounded map from $\banach{p}(\torus)$ into
  $\ell_{\frac{p}{p-1}}(\lattice\dual)$
  and the inverse $\Op{F}^{-1}$ from
  $\ell_{\frac{q}{q-1}}(\lattice\dual)$ into
  $\banach{q}(\torus)$.
  If $\lambda$ does not belong to the spectrum of the free
  Dirac operator $\triv{\Op{D}}(0,0,k)$,
  then the Fourier transform
  $\Op{F}\comp\triv{\Op{R}}(0,0,k,\sqrt{-1}\lambda)
  \comp\Op{F}^{-1}$
  acts on $\ell_{\frac{p}{p-1}}(\lattice\dual)\times
  \ell_{\frac{p}{p-1}}(\lattice\dual)$ as the point--wise multiplication
  with the sequence of matrices
  $$\begin{pmatrix}
  \sqrt{-1}\lambda & -\pi(k_2+\kappa_2+\sqrt{-1}(k_1+\kappa_1))\\
  -\pi(k_2+\kappa_2-\sqrt{-1}(k_1+\kappa_1)) & \sqrt{-1}\lambda
  \end{pmatrix}^{-1}.$$
  Obviously, for all $k\in\mathbb{C}^2$, $\varepsilon>0$ and $r>2$ there
  exists some $\delta>0$ such that all real $\lambda>\delta$
  fulfill the following estimate: 
  $$\left(\sum\limits_{\kappa\in\lattice\dual}
  \left\|\begin{pmatrix}
  \sqrt{-1}\lambda & -\pi(k_2+\kappa_2+\sqrt{-1}(k_1+\kappa_1))\\
  -\pi(k_2+\kappa_2-\sqrt{-1}(k_1+\kappa_1)) & \sqrt{-1}\lambda
  \end{pmatrix}^{-1}\right\|^r\right)^{1/r}< \infty.$$
  With the help of H\"older's inequality
  \cite[Theorem~III.1~(c)]{RS1} we conclude that for all $k\in
  \mathbb{C}^2$ there exists some $\delta>0$ such that for all real
  $\lambda>\delta$ the resolvent
  $\triv{\Op{R}}(0,0,k,\sqrt{-1}\lambda)$
  is a bounded operator from
  $\banach{p}(\torus)\times\banach{p}(\torus)$ into
  $\banach{q}(\torus)\times\banach{q}(\torus)$.
  Let $\triv{\Op{P}}_n$ be the natural projection
  onto the span of all $\psi_{\kappa}$
  with $g(\kappa,\kappa)\leq n$. The same argument shows that for all
  $l\in\mathbb{N}$ and the same $V$ and $W$ as before the sequence
  $$\triv{\Op{R}}(0,0,k,\sqrt{-1}\lambda)
  \left(\comp\triv{\Op{P}}_n\comp
  \left(\begin{smallmatrix}
  V & 0\\
  0 & W
  \end{smallmatrix}\right)
  \comp\triv{\Op{P}}_n\comp
  \triv{\Op{R}}(0,0,k,\sqrt{-1}\lambda)\right)^l$$
  converges in the limit $n\rightarrow \infty$ uniformly
  with respect to $V$ and $W$ to the operator
  $$\triv{\Op{R}}(0,0,k,\sqrt{-1}\lambda)\left(\comp
  \left(\begin{smallmatrix}
  V & 0\\
  0 & W
  \end{smallmatrix}\right)
  \comp\triv{\Op{R}}(0,0,k,\sqrt{-1}\lambda)\right)^l$$
  on $\banach{p}(\torus)\times\banach{p}(\torus)$.
  Since all elements of the sequences are
  finite rank operators, the limits are compact operators
  \cite[Theorem~VI.12]{RS1}. Obviously these sequences are
  sequences of weakly continuous functions in dependence of $(V,W)$.
  Since the uniform limits of continuous functions on a compact
  Hausdorff space are again a continuous functions, the limits are weakly
  continuous on the sets
  $$\left\{(V,W)\in \banach{2}(\torus)\times
  \banach{2}(\torus)\mid\|V\|_{\varepsilon,2}\leq C_p\text{ and }
  \|W\|_{\varepsilon,2}\leq C_p\right\}\times\Set{O}.$$
  Due to the considerations above the Neumann series converges
  uniformly on these sets, and the lemma follows.
\end{proof}

\noindent
{\it Continuation of the proof of Theorem~\ref{meromorph}.}
In the last lemma we actually proved that for all $(V_0,W_0,k_0)\in
\banach{2}(\torus)\times\banach{2}(\torus)\times
\mathbb{C}^2$, there exists a complex number $\lambda_1$ and some
open neighbourhood $\Set{U}$ of $(V_0,W_0,k_0)$, so that the
map $(V,W,k)\mapsto \triv{\Op{R}}(V,W,k,\lambda_1)$
is a holomorphic map from $\Set{U}$ into the compact operators on
$\banach{2}(\torus)\times\banach{2}(\torus)$.
Thus, due to the Riesz--Schauder theorem \cite[Theorem~VI.15]{RS1}
the spectrum of the operator 
$\triv{\Op{R}}(V,W,k,\lambda_1)$ is discrete and 
the only limit point of the spectrum is zero.
Let $\Set{S}\subset \mathbb{C}$ denote the set 
$$\Set{S}=\left\{ \left. \lambda_1 -\frac{1}{ \lambda }\;\right|\;
\lambda \text{ is an eigenvalue of } 
\triv{\Op{R}}(V_0,W_0,k_0,\lambda_1) \right\} . $$
Then for all $\lambda \in \mathbb{C}\setminus \Set{S}$, the operator
$$\triv{\Op{R}}(V_0,W_0,k_0,\lambda)=
\frac{\triv{\Op{R}}(V_0,W_0,k_0,\lambda_1)}
     { \lambda_1 - \lambda}\comp\left(
\frac{\unity }{ \lambda_1 -\lambda }-
\triv{\Op{R}}(V_0,W_0,k_0,\lambda_1) \right)^{-1}$$
is compact and depends holomorphically on $\lambda $.

Let $\lambda_0\in \Set{S}$ be a pole of this map and choose a
positive number $\varepsilon >0$, such that there exists
no other pole in the disc
$\{ \lambda\mid |\lambda -\lambda_0| < 2 \varepsilon \}$.
Then there exists an open neighbourhood $\Set{U} \subset 
\banach{2}(\torus)\times\banach{2}(\torus)\times
\mathbb{C}^2$ of $(V_0,W_0,k_0)$ such that
$\triv{\Op{R}}(V,W,k,\lambda)$
is compact for all $(V,W,k)\in \Set{U}$ and all $\lambda $ on the
circle 
$\{ \lambda\mid |\lambda -\lambda_0| =\varepsilon \} $. 
Now we define a map $\Breve{\Op{P}}$ from $\Set{U}$
into the set of compact operators on
$\banach{2}(\torus)\times\banach{2}(\torus)$ by
$$\Breve{\Op{P}}(V,W,k)= -\frac{1}{ 2 \pi \sqrt{-1}} 
\oint\limits_{\{ \lambda\mid |\lambda -\lambda_0| =\varepsilon \}}
\triv{\Op{R}}(V,W,k,\lambda)d\lambda .$$
The combination of the results from \cite[VI.5]{RS1}
and \cite[Appendix to XII.1]{RS4} implies that for all
$(V,W,k)\in \Set{U}$, $\Breve{\Op{P}}(V,W,k)$ is a finite--dimensional
projector, which commutes with $\triv{\Op{D}}$.
Also the eigenvalues of the restriction to the
corresponding finite--dimensional subspace of
$\banach{2}(\torus)\times\banach{2}(\torus)$
are elements of the disc 
$\{ \lambda\mid |\lambda -\lambda_0| <\varepsilon \} $. Furthermore, all
eigenfunctions, which correspond to eigenvalues inside of this disc, are
elements of the finite--dimensional subspace mentioned above.
It remains to prove the last statement of the theorem. In
Lemma~\ref{weakly continuous resolvent} we actually proved that for
all $1<p<2$ and all $p\leq q<\frac{2p}{2-p}$ the resolvent
$\triv{\Op{R}}(V,W,k,\sqrt{-1}\lambda)$ is a compact operator from
$\banach{p}(\torus)\times\banach{p}(\torus)$ into
$\banach{q}(\torus)\times\banach{q}(\torus)$. This
implies that $\Breve{\Op{P}}(V,W,k)$
is a finite--dimensional projector on 
$\banach{q}(\torus)\times\banach{q}(\torus)$, for all
$2\leq q<\infty$.
\end{proof}

\begin{Remark}\label{sobolev embedding}
Due to \cite[Chapter~V \S.3.4 Lemma~3.]{St} the operators
$\partial\comp\left(\unity-\partial\Bar{\partial}\right)^{-1/2}$ and 
$\Bar{\partial}\comp\left(\unity-\partial\Bar{\partial}\right)^{-1/2}$
are bounded operators on $\banach{p}(\torus)$ with $1<p<2$.
Therefore any function $f\in\banach{p}(\torus)$ with
either $\partial f\in\banach{p}(\torus)$ or
$\Bar{\partial}f\in\banach{p}(\torus)$ belongs to the Sobolev space
$\sobolev{1,p}(\torus)$ \cite[Chapter~V \S.3.4 Theorem~3.]{St}.
Hence our arguments show that the eigenfunctions belong
to the intersection $\bigcap\limits_{1<p<2}
\sobolev{1,p}(\Delta)\times\sobolev{1,p}(\Delta)$ of Sobolev spaces.
\end{Remark}

\begin{Remark}\label{sobolev 2}
Unfortunately the sets described in
Lemma~\ref{weakly continuous resolvent}, by which the whole
complex Bloch varieties depend weakly continuous on the potentials,
do not contain balls, whose radius is larger or equal to
$S_p^{-1}$ (compare with Remark~\ref{sobolev 1}).
This failure will cause the most difficulties of this approach.
So let us elaborate a bit on this failure.
For all elements $(V_0,W_0)$ of the sets
$$\left\{(V,W)\in \banach{2}(\torus)\times
\banach{2}(\torus)\mid \|V\|_{\varepsilon,2}\leq C_p\text{ and }
\|W\|_{\varepsilon,2}\leq C_p\right\}\text{ with $C_p<S_p^{-1}$}$$
the balls $B(V_0,C_p')\times B(W_0,C_p')$,
with $C_p'<S_p^{-1}$ are contained in another set of this form
$$\left\{(V,W)\in \banach{2}(\torus)\times
\banach{2}(\torus)\mid \|V\|_{\varepsilon',2}\leq C_p''\text{ and }
\|W\|_{\varepsilon',2}\leq C_p''\right\}\text{ with $C_p''<S_p^{-1}$}$$
and with some small $\varepsilon'>0$. Moreover, all bounded
subsets of
$\banach{r}(\torus)\times\banach{r}(\torus)$ with $2<r$
are contained in a set of this form. Finally, we shall see in
Lemma~\ref{bounded point measures} that a weakly convergent
subsequence is contained in a set of this form, if the
weak limit of the measures $V(x)\Bar{V}(x)d^2x$ and
$W(x)\Bar{W}(x)d^2x$ exist and contain no point measures
with mass larger or equal to $S_p^{-2}$.
If the norms of these sequences are bounded, this can fail
only at finitely many points of $\torus$.
In this case the limits of the corresponding resolvents are the
resolvents of some perturbation of the Dirac operator acting on
eigenfunctions with poles
(or on line bundles of non--vanishing degree)
as described in Section~\ref{subsection limits}.
\end{Remark}

\index{unique continuation property|(}
We close this section with a proof of the
\Em{strong unique continuation property} of the eigenfunctions of the
Dirac operators. We use the classical Carleman method, which is based
on an improved Sobolev inequality (compare with \cite{Ca} and
\cite[Proposition~1.3]{Wo}).

\newtheorem{Carleman inequality}[Lemma]{Carleman inequality}
\index{Carleman inequality}
\begin{Carleman inequality}\label{carleman inequality}
There exists some constant $S_p$
(compare with Lemma~\ref{free resolvent} and Remark~\ref{sobolev 1})
depending only on $1<p<2$, such that for all $n\in\mathbb{Z}$
and all $\psi\in C^{\infty}_0(\mathbb{C}\setminus\{0\})\times
C^{\infty}_0(\mathbb{C}\setminus\{0\})$ the following inequality
holds:
$$\left\| |z|^{-n}\psi\right\|_{\frac{2p}{2-p}}\leq S_p
\left\| |z|^{-n}\left(\begin{smallmatrix}
0 & \partial\\
-\Bar{\partial} & 0
\end{smallmatrix}\right)\psi\right\|_{p}.$$
\end{Carleman inequality}

The literature \cite{Je,Ki1,Ma,Ki2} deals with the much more difficult
higher--dimensional case and does not treat our case.
David Jerison pointed out to the author, that the arguments of
\cite[Proposition~2.6]{Wo}, where the analogous but weaker statement
about the gradient term of the Laplace operator is treated,
carry over to the Dirac operator.

\begin{proof}
Dolbeault's Lemma \cite[Chapter~I Section~D 2.~Lemma]{GuRo}
implies for all smooth $\psi$ with compact support the equality
$$\psi(z)=\int\limits_{\mathbb{C}}
\left(\begin{smallmatrix}
0 & (z'-z)^{-1}\\
(\Bar{z}-\Bar{z}')^{-1} & 0
\end{smallmatrix}\right)\comp\left(\begin{smallmatrix}
0 & \partial\\
-\Bar{\partial} & 0
\end{smallmatrix}\right)\psi(z')
\frac{d\Bar{z}'\wedge dz'}{2\pi\sqrt{-1}}.$$
In fact, the components of the difference of the left hand side minus
the right hand side are holomorphic and anti--holomorphic functions on
$\mathbb{C}$, respectively, which vanish at $z=\infty$.
In particular, the integrals
$\int_\mathbb{C} z^n\Bar{\partial}\psi_1d\Bar{z}\wedge dz$ and
$\int_\mathbb{C}\Bar{z}^n\partial\psi_2d\Bar{z}\wedge dz$
with $n\in\mathbb{N}_0$ are proportional to the Taylor coefficients
of $\psi$ at $\infty$, which vanish.
Moreover, if the support of $\psi$
does not contain $0$ and therefore also a small neighbourhood of $0$,
then the integrals
$\int_\mathbb{C}z^{-n}\Bar{\partial}\psi_1d\Bar{z}\wedge dz$ and
$\int_\mathbb{C}\Bar{z}^{-n}\partial\psi_2d\Bar{z}\wedge dz$ with
$n\in\mathbb{N}$ are proportional to the Taylor coefficients of
$\psi$ at $0$, which in this case also
vanish. These cancellations follow also from partial integration.
We conclude that for all $n\in\mathbb{Z}$
$$\psi(z)=\int\limits_{\mathbb{C}}
\begin{pmatrix}
0 & \left(\frac{z}{z'}\right)^n\frac{1}{z'-z}\\
\left(\frac{\Bar{z}}{\Bar{z}'}\right)^n\frac{1}{\Bar{z}-\Bar{z}'} & 0
\end{pmatrix}\comp\left(\begin{smallmatrix}
0 & \partial\\
-\Bar{\partial} & 0
\end{smallmatrix}\right)\psi(z')
\frac{d\Bar{z}'\wedge dz'}{2\pi\sqrt{-1}}.$$
In fact, for negative $n$ the left hand side minus the left hand side
of the foregoing formula is equal to the Taylor polynomial of $\psi$
at $\infty$ up to order $|n|$, and for positive $n$ equal to the
Taylor polynomial of $\psi$ at $0$ up to order $n-1$.
Finally, the Hardy--Littlewood--Sobolev theorem
\cite[Chapter~V. \S1.2 Theorem~1]{St} implies that the operator with
integral kernel
$$\begin{pmatrix}
0 & \left(\frac{|z'|z}{|z|z'}\right)^n\frac{1}{z'-z}\\
\left(\frac{|z'|\Bar{z}}{|z|\Bar{z}'}\right)^n\frac{1}{\Bar{z}-\Bar{z}'}
& 0
\end{pmatrix}\frac{d\Bar{z}'\wedge dz'}{2\pi\sqrt{-1}}$$
from $\banach{p}(\mathbb{C})\times\banach{p}(\mathbb{C})$ into
$\banach{\frac{2p}{2-p}}(\mathbb{C})\times
\banach{\frac{2p}{2-p}}(\mathbb{C})$
is bounded by some constant $S_p$ not depending on $n$, and maps
$|z|^{-n}\left(\begin{smallmatrix}
0 & \partial\\
-\Bar{\partial} & 0
\end{smallmatrix}\right)\psi$ onto $|z|^{-n}\psi$.
\end{proof}

Due to a standard argument
(e.\ g.\ \cite[Proof of Theorem~5.1.4]{So} and
\cite[Section Carleman Method]{Wo})
this \De{Carleman inequality} implies the

\newtheorem{Unique continuation}[Lemma]{Strong unique continuation property}
\begin{Unique continuation}\label{strong unique continuation}
Let $V$ and $W$ be potentials in
$\banach{2}_{\text{\scriptsize\rm loc}}(\Set{O})$ on an
open connected set $0\ni\Set{O}\subset\mathbb{C}$ and
$\psi\in \sobolev{1,p}_{\text{\scriptsize\rm loc}}(\Set{O})\times
\sobolev{1,p}_{\text{\scriptsize\rm loc}}(\Set{O})$
an element of the kernel of $\left(\begin{smallmatrix}
V & \partial\\
-\Bar{\partial} & W
\end{smallmatrix}\right)$ on $\Set{O}$ with $1<p<2$.
If the $\banach{\frac{2p}{2-p}}$--norm of the restriction of
$\psi$ to the balls $B(0,\varepsilon)$ converges in the limit
$\varepsilon\downarrow 0$ faster to zero than any power of $\varepsilon$:
$$\left(\int\limits_{B(0,\varepsilon)}
(|\psi_1|^{\frac{2p}{2-p}}+|\psi_2|^{\frac{2p}{2-p}})d^2x
\right)^{\frac{2-p}{2p}}\leq
\text{\bf{O}}(\varepsilon^n)\;\forall n\in\mathbb{N},$$
then $\psi$ vanishes identically on $\Set{O}$.
\end{Unique continuation}

\begin{proof}
The question is local so we may assume that $V$ and $W$ are elements
of $\banach{2}$ rather than $\banach{2}_{\text{\scriptsize\rm loc}}$. We fix
$\varepsilon$ small enough that
$\max\left\{\|V\|_{\banach{2}(B(z,2\varepsilon))},\|W\|_{\banach{2}(B(z,2\varepsilon))}\right\}
\leq 1/(2S_p)$ for all $z$, where $S_p$ is the constant of the
\De{Carleman inequality}.
Let $\phi\in C^{\infty}$ be $1$ on $B(0,\varepsilon)$ and $0$ on
$\mathbb{C}\setminus B(0,2\varepsilon)$. A limiting argument using the
infinite order vanishing of $\psi$ and the equality
$\left(\begin{smallmatrix}
0 & \partial\\
-\Bar{\partial} & 0
\end{smallmatrix}\right)\psi=-\left(\begin{smallmatrix}
V & 0\\
0 & W
\end{smallmatrix}\right)\psi$ shows that the proof of the
\De{Carleman inequality} is also true for $\phi\psi$. So
$$\left\| |z|^{-n}\phi\psi\right\|_{\frac{2p}{2-p}}\leq
S_p\left\| |z|^{-n}\left(\begin{smallmatrix}
0 & \partial\\
-\Bar{\partial} & 0
\end{smallmatrix}\right)\phi\psi\right\|_p\leq
S_p\left\| |z|^{-n}\phi\left(\begin{smallmatrix}
0 & \partial\\
-\Bar{\partial} & 0
\end{smallmatrix}\right)\psi\right\|_p+
S_p\left\| |z|^{-n}E\right\|_p.$$
Here $E$ (for error) $=\left(\begin{smallmatrix}
\psi_2\partial\phi\\
-\psi_1\Bar{\partial}\phi
\end{smallmatrix}\right)$ is an $\banach{p}\times\banach{p}$
function supported in
$\{x\mid \varepsilon\leq|z|\leq 2\varepsilon\}$. Using the equality
$\left(\begin{smallmatrix}
0 & \partial\\
-\Bar{\partial} & 0
\end{smallmatrix}\right)\psi=-\left(\begin{smallmatrix}
V & 0\\
0 & W
\end{smallmatrix}\right)\psi$ and
H\"older's inequality \cite[Theorem~III.1~(c)]{RS1} yields
\begin{eqnarray*}
\left\| |z|^{-n}\phi\psi\right\|_{\frac{2p}{2-p}}&\leq&
S_p\left\| |z|^{-n}\phi\left(\begin{smallmatrix}
V & 0\\
0 & W
\end{smallmatrix}\right)\psi\right\|_p+S_p\left\| |z|^{-n}E\right\|_p\\
&\leq& S_p\max\left\{\|V\|_{\banach{2}(B(0,2\varepsilon))},
\|W\|_{\banach{2}(B(0,2\varepsilon))}\right\}
\left\| |z|^{-n}\phi\psi\right\|_{\frac{2p}{2-p}}+
S_p\left\| |z|^{-n}E\right\|_p.
\end{eqnarray*}
By the choice of $\varepsilon$ the first term can be absorbed
into a factor $2$
$$\left\| |z|^{-n}\phi\psi\right\|_{\frac{2p}{2-p}}\leq
2S_p\left\| |z|^{-n}E\right\|_p.$$
Now comes the crucial observation: $E$ is supported in
$\{z\mid \varepsilon\leq|z|\leq 2\varepsilon\}$, so
\begin{align*}
\left\| |z|^{-n}\phi\psi\right\|_{\frac{2p}{2-p}}&
\leq 2S_p\varepsilon^{-n}\left\|E\right\|_p
&\text{ and }\left\|\left(\frac{\varepsilon}{|z|}
\right)^n\phi\psi\right\|_{\frac{2p}{2-p}}
&\leq 2S_p\left\|E\right\|_p.
\end{align*}
Using the limit $n\rightarrow\infty$ we conclude that $\phi\psi$
vanishes on $B(0,\varepsilon)$, and therefore also $\psi$.
In other words, the set
$\{z\mid\psi\text{ vanishes to infinite order at }z\}$
is open, and in fact contains a ball of fixed radius $\varepsilon$
centered at any of its points.
So this set must be all of $\Set{O}$ and the proof is complete.
\index{unique continuation property|)}
\end{proof}

\subsection{The Bloch variety and the Fermi curve}
\label{subsection fermi curve}

For two given periodic potentials $V$ and $W$ the corresponding
Dirac operator commutes with all shifts by the lattice
vectors. Hence these shifts and the Dirac operator
may be diagonalized simultaneously. The common spectrum
of all these operators is called a Bloch variety.\index{Bloch!variety}
The proof in \cite[Theorem~XIII.97]{RS4},
which applies to Schr\"odinger operators with periodic
potential, carries over to the Dirac operators and shows that a
Dirac operator with periodic potentials $V$ and $W$ is the direct
integral of the operators $\triv{\Op{D}}(V,W,k)$
where $k$ runs through $\torus\dual$
(compare with \cite{Ku}).
The Bloch variety may be considered as the set of points
$([k],\lambda)\in\torus\dual\times \mathbb{C}$, so that 
$\lambda $ is an eigenvalue of $\Op{D}(V,W,k)$.
Let us call the collection of all points $(k,\lambda)
\in \mathbb{C}^2\times \mathbb{C}$, so that 
$\lambda $ is an eigenvalue of $\Op{D}(V,W,k)$,
complex Bloch variety and denote this set by
$\bloch(V,W)$\index{Bloch!complex $\sim$ variety $\bloch$}.
Due to Theorem~\ref{meromorph} for all $(k_0,\lambda_0)\in
\mathbb{C}^2\times\mathbb{C}$ there exists some open
neighbourhood
$\Set{U}\times\{ \lambda\mid |\lambda -\lambda_0| <\varepsilon \}
\subset \mathbb{C}^2\times \mathbb{C}$ of
$(k_0,\lambda_0)$ and a holomorphic $l\times l$ matrix valued
function $\Mat{A}$ on $\Set{U}$,
such that $\bloch(V,W)\cap
(\Set{U}\times\{ \lambda\mid |\lambda -\lambda_0| <\varepsilon \})$
is the zero set of the function
$$\Set{U}\times\{ \lambda\mid |\lambda -\lambda_0| <\varepsilon \}
\rightarrow \mathbb{C},\; (k,\lambda)\mapsto
\det\left(\lambda\unity - \Mat{A}(k)\right).$$
Hence $\bloch(V,W)$ is a subvariety of 
$\mathbb{C}^2\times \mathbb{C}$.
The action of the dual lattice
$\lattice\dual$ on $\mathbb{C}^2$ given by $\kappa . k = \kappa+k$
for all $\kappa \in \lattice\dual$ and all $k\in \mathbb{C}^2$,
lifts to an action of the dual lattice $\lattice\dual$
on this subspace $\bloch(V,W)$ 
of $\mathbb{C}^2\times \mathbb{C}$. The set
of orbits of this action is equal to the subspace 
$\bloch(V,W)/\lattice\dual$
of all elements $([k],\lambda)\in
\mathbb{C}^2/\lattice\dual\times \mathbb{C}$,
such that $\lambda $ is an eigenvalue of the operator
$\Op{D}(V,W,k)$.

We call the zero energy level of the
complex Bloch variety the \Em{complex Fermi curve} and denote it by
\index{Fermi curve!complex $\sim$ $\fermi$}
\index{curve!complex Fermi $\sim$ $\fermi$}
$$\fermi(V,W)=\{k\in\mathbb{C}^2\mid
(k,0)\in\bloch(V,W)\}.$$
Sometimes we shall not distinguish between the values of $k$
and the corresponding elements of $\fermi$.
However, mostly $k$ denotes the corresponding
$\mathbb{C}^2$--valued function on the \Em{complex Fermi curves}.
Due to Theorem~\ref{meromorph} the \Em{complex Fermi curve}
$\fermi(V,W)$ is locally
the zero set of one holomorphic function. Hence it is a pure
one--dimensional subvariety of $\mathbb{C}^2$, which is invariant
under the action of $\lattice\dual$.
Strictly speaking, the \Em{complex Fermi curve} is
the quotient space $\fermi(V,W)/\lattice\dual$,
but we shall also call $\fermi(V,W)$ a \Em{complex Fermi curve}.

These \Em{complex Fermi curves} are of special interest
for two reasons: Firstly the corresponding eigenfunctions belong to
the kernel of the Dirac operator and therefore give rise to
local immersions of surfaces into the three--dimensional
Euclidean space due to the local Weierstra{\ss} representation.
Secondly, these \Em{complex Fermi curves} are the spectral curves of
the corresponding integrable system described by the
Davey--Stewartson equation.
\index{Davey--Stewartson equation}

Due to a general feature of differential operators, in the
high energy limit the higher--order derivatives dominate the
lower oder derivatives. Moreover, typically the complex Bloch variety
converges in the high energy limit to the complex Bloch variety of the
free Dirac operator and has an asymptotic expansion for
smooth potentials, which can be calculated explicitly
(compare with \cite[Chapter~1]{Sch}).
In the remainder of this section we shall investigate the asymptotic
behaviour of the \Em{complex Fermi curves} for general potentials
$V,W\in \banach{2}(\torus)$.

For this purpose we use a covariance property of the
\Em{complex Fermi curves} under some unitary transformations
of the potentials. By definition of the operators $\partial_{[k]}$ and
$\Bar{\partial}_{[k]}$ they transform under translations by some $k'$
as follows:
\begin{align*}
\partial_{[k+k']} &=\psi_{-k'}\comp\partial_{[k]}\comp\psi_{k'}
&\text{if }k'_2+\sqrt{-1}k'_1=0,\\
\Bar{\partial}_{[k+k']} &=\psi_{-k'}\comp\Bar{\partial}_{[k]}\comp\psi_{k'}
&\text{if }k'_2-\sqrt{-1}k'_1=0.
\end{align*}
For all $\kappa$ there exist unique pairs
$(k^-_{\kappa},k^+_{\kappa})\in\mathbb{C}^2\times\mathbb{C}^2$,
whose first element obeys the former condition
and the second element the latter condition, and whose difference
$k^+_{\kappa}-k^-_{\kappa}$ is equal to $\kappa$.

We conclude that for all $\kappa\in\lattice\dual$ the Dirac operators
transform as follows:
$$\Op{D}(V,W,k+k^{\pm}_{\kappa})=
\begin{pmatrix}
\psi_{-k^-_{\kappa}} & 0\\
0 & \psi_{-k^+_{\kappa}}
\end{pmatrix}\comp
\Op{D}(\psi_{-\kappa}V,\psi_{\kappa}W,k)\comp
\begin{pmatrix}
\psi_{k^+_{\kappa}} & 0\\
0 & \psi_{k^-_{\kappa}}
\end{pmatrix}.$$
This implies the following

\begin{Lemma} \label{covariant transformation}
For all $\kappa\in\lattice\dual$ the translation $k\mapsto
k+k^+_{\kappa}=k+k^-_{\kappa}+\kappa$ induces a biholomorphic
isomorphism of the \Em{complex Fermi curves}
$\fermi(V,W)\simeq \fermi(\psi_{-\kappa}V,\psi_{\kappa}W)$.
\qed
\end{Lemma}

For all potentials $V,W\in \banach{2}(\torus)$
there exists an $\varepsilon>0$ such that
$\|V\|_{\varepsilon,2}< S_p^{-1}$ and $\|W\|_{\varepsilon,2}<S_p^{-1}$.
Consequently, for all $\kappa\in \lattice\dual$
also the potentials $\psi_{-\kappa}V$ and $\psi_{\kappa}W$
obey these estimates.
Moreover, the sequence of potentials
$(\psi_{-\kappa}V,\psi_{\kappa}W)$ converges in the limit
$g(\kappa,\kappa)\rightarrow\infty$ weakly to the zero potentials in
$\banach{2}(\torus)$.
More precisely, if we identify the Banach space
$\banach{1}(\torus)\simeq \banach{1}(\Delta)$ with the corresponding
natural subspace of $\banach{1}(\mathbb{R}^2)$,
then due to the Riemann--Lebesgue Lemma \cite[Theorem~IX.7]{RS2}
the Fourier coefficients indexed by $\kappa\in\lattice\dual$
of any element in $\banach{1}(\torus)$
converges in the limit $g(\kappa,\kappa)\rightarrow\infty$ to zero.
Therefore, for any $V,W\in \banach{2}(\torus)$
and any weakly open neighbourhood $\Set{U}$ of
$0\in \banach{2}(\torus)\times\banach{2}(\torus)$ 
there exists some $\delta>0$ such that for all
$\kappa\in\lattice\dual$ with $g(\kappa,\kappa)>\delta^{-2}$
the pair of potentials
$(\psi_{-\kappa}V,\psi_{\kappa}W)$ belongs to $\Set{U}$.
We conclude from Lemma~\ref{weakly continuous resolvent}
that for all $k$ in a bounded open subsets
$\Set{O}\subset\mathbb{C}^2\setminus\fermi(0,0)$
the resolvents $\Op{R}(\psi_{-\kappa}V,\psi_{\kappa}W,k,0)$
converge in the limit $g(\kappa,\kappa)\rightarrow\infty$
uniformly to $\Op{R}(0,0,k,0)$.
Obviously the eigenvalues of the operators
$\left(\begin{smallmatrix}
0 & \pm\unity\\
\unity & 0
\end{smallmatrix}\right)\comp\triv{\Op{D}}(V,W,k)$ are completely
determined by the \Em{complex Fermi curve} $\fermi(V,W)$.
Conversely, one resolvent of the family of resolvents
$\Op{R}(V,W,k,0)$ indexed by
$k\in\mathbb{C}^2\setminus\fermi(V,W)$ determines the spectrum
of the operators $\left(\begin{smallmatrix}
0 & \pm\unity\\
\unity & 0
\end{smallmatrix}\right)\comp\triv{\Op{D}}(V,W,k)$ and therefore
also all other resolvents of this family.
Therefore, the intersections of the \Em{complex Fermi curves}
$\fermi(\psi_{-\kappa}V,\psi_{\kappa}W)$
with all bounded open subsets $\Set{O}\subset\mathbb{C}^2$
converge in the limit $g(\kappa,\kappa)\rightarrow\infty$
to $\Set{O}\cap\fermi(0,0)$.
Moreover, the corresponding eigenfunctions converge to the
corresponding eigenfunctions of the free Dirac operator.
An application of Lemma~\ref{covariant transformation}
yields that the intersections of
$\fermi(V,W)\cap\left(\Set{O}+k^+_{\kappa}\right)=
\left(\fermi(V,W)\cap\left(\Set{O}+k^-_{\kappa}\right)\right)+\kappa$
converges to
$\fermi(0,0)\cap\left(\Set{O}+k^+_{\kappa}\right)=
\left(\fermi(0,0)\cap\left(\Set{O}+k^-_{\kappa}\right)\right)+\kappa$.
An easy calculation shows that the free \Em{complex Fermi curve}
is equal to
$$\fermi(0,0)=\bigcup\limits_{\kappa\in\lattice\dual}
\kappa+\left\{k\in\mathbb{C}^2\mid k_2+\sqrt{-1}k_1=0\right\}
\bigcup\limits_{\kappa\in\lattice\dual}
\kappa+\left\{k\in\mathbb{C}^2\mid k_2-\sqrt{-1}k_1=0\right\}.$$
Therefore, the free \Em{complex Fermi curve}
$\fermi(0,0)/\lattice\dual$ is isomorphic to the two
subvarieties $\{k\mid k_1\pm\sqrt{-1}k_2=0\}$ of $\mathbb{C}^2$,
glued at infinitely many ordinary double points indexed by
$\kappa\in\lattice\dual$ of the form $(k_{\kappa}^-,k_{\kappa}^+)$ with
\begin{eqnarray*}
k_{\kappa}^- &=&
\left(-\kappa_1/2-\sqrt{-1}\kappa_2/2,-\kappa_2/2+\sqrt{-1}\kappa_1/2\right)
\text{ and }\\
k_{\kappa}^+ &=&
\left(\kappa_1/2-\sqrt{-1}\kappa_2/2,\kappa_2/2+\sqrt{-1}\kappa_1/2\right)
=k_{\kappa}^-+\kappa.
\end{eqnarray*}
\index{ordinary double point!--s of the free \Em{complex Fermi curve}
       $(k^+_{\kappa},k^-_{\kappa})$}
Also on these components the eigenfunctions of the free Dirac operator
are equal to
\begin{align*}
\psi=\psi_{k}\left(\begin{smallmatrix}
0\\
1
\end{smallmatrix}\right)&
\text{ on $\left\{k\mid k_1-\sqrt{-1}k_2=0\right\}$ and}&
\psi=\psi_{k}\left(\begin{smallmatrix}
1\\
0
\end{smallmatrix}\right)&
\text{ on $\left\{k\mid k_1+\sqrt{-1}k_2=0\right\}$.}
\end{align*}
respectively. For positive $\varepsilon$ and $\delta$ let
$\Set{V}^{\pm}_{\varepsilon,\delta}$ denote the intersection of
$\fermi(V,W)$ with the following open subsets
$\Set{U}^{\pm}_{\varepsilon,\delta}$ of $\mathbb{C}^2$:
\index{neighbourhood!$\Set{U}^{\pm}_{\varepsilon,\delta}$}
\index{neighbourhood!$\Set{V}^{\pm}_{\varepsilon,\delta}$}
$$\Set{U}^{\pm}_{\varepsilon,\delta}=
\left\{k\in\mathbb{C}^2\mid
\left|k_1\pm\sqrt{-1}k_2\right|<\varepsilon,\|k\|>1/\delta
\text{ and }\|k-k^{\pm}_\kappa\|>\varepsilon
\;\forall \kappa\in\lattice\dual\setminus\{0\}\right\}.$$
Here $\|k\|$ denotes the natural hermitian norm
$\|k\|=\sqrt{g(k,\Bar{k})}$ on $\mathbb{C}^2$.
Since for small $\varepsilon$ and $\delta$ these open sets
$\Set{U}^{\pm}_{\varepsilon,\delta}$ intersect any
$\lattice\dual$--orbit at most once, we may consider these open
sets as subsets of $\mathbb{C}^2/\lattice\dual$.
Since the quotient of $\mathbb{C}^2$ modulo the sublattice
$\left\{k^+_{\kappa}+\kappa'\mid\kappa,\kappa'\in\lattice\dual\right\}
=\left\{k^-_{\kappa}+\kappa'\mid\kappa,\kappa'\in\lattice\dual\right\}$
is compact, the combination of Lemma~\ref{weakly continuous resolvent}
and Lemma~\ref{covariant transformation} implies the following

\begin{Theorem}\label{asymptotic analysis 1}
\index{handle}
\index{asymptotic analysis}
For all pairs of potentials
$V,W\in \banach{2}(\torus)$ and all $\varepsilon>0$
there exists a $\delta>0$ so that the two open sets
$\Set{V}^{\pm}_{\varepsilon,\delta}$ are connected complex
one--dimensional submanifolds of $\mathbb{C}^2$ with one--dimensional
kernels of the corresponding Dirac operators $\Op{D}(V,W,k)$.
Moreover, the relative complement of
$\Set{V}^+_{\varepsilon,\delta}\cup\Set{V}^-_{\varepsilon,\delta}$
in $\fermi(V,W)/\lattice\dual$
decomposes into a compact set contained in
$\{k\in\mathbb{C}^2\mid \|k\|\leq 1/\delta\}/\lattice\dual$
and infinitely many small handles indexed by
$\lattice\dual_{\delta}=
\{\kappa\in\lattice\dual\mid \|k^+_{\kappa}\|>1/\delta\}$.
\index{index set!of handles $\lattice\dual_{\delta}$}
The handle with index $\kappa$ is contained in
$\left\{k\in\mathbb{C}^2\mid \|k-k^+_{\kappa}\|\leq\varepsilon\right\}=
\left\{k\in\mathbb{C}^2\mid \|k-k^-_{\kappa}\|\leq\varepsilon\right\}+\kappa$
and connects the small disc around $k^+_{\kappa}$ excluded from
$\Set{V}^+_{\varepsilon,\delta}$ with with the small disc around
$k^-_{\kappa}$ excluded from $\Set{V}^-_{\varepsilon,\delta}$.
Finally, for all $q<\infty$ the Fourier components $\fourier{\psi}(\kappa)$
of the periodic parts
$$\triv{\psi}=\exp\left(-2\pi\sqrt{-1}g(x,k)\right)\psi=
\sum\limits_{\kappa\in\lattice\dual}=
\fourier{\psi}(\kappa)\exp\left(2\pi\sqrt{-1}g(x,\kappa)\right)$$
obey for some $\varepsilon$ (depending on $q$)
on these parts of $\fermi(V,W)$ the following estimates:
\begin{align*}
\left\|\triv{\psi}-\begin{pmatrix}
0\\
\fourier{\psi}_2(0)
\end{pmatrix}\right\|_q<\varepsilon\left|\fourier{\psi}_2(0)\right|
&\text{ on $\Set{V}^-_{\varepsilon,\delta}$, }&
\left\|\triv{\psi}-\begin{pmatrix}
\fourier{\psi}_1(0)\\
0
\end{pmatrix}\right\|_q<\varepsilon\left|\fourier{\psi}_1(0)\right|
&\text{ on $\Set{V}^+_{\varepsilon,\delta}$, and}
\end{align*}
$$\left\|\triv{\psi}-\begin{pmatrix}
\fourier{\psi}_1(0)\\
\fourier{\psi}_2(-\kappa)\exp\left(-2\pi\sqrt{-1}g(x,\kappa)\right)
\end{pmatrix}\right\|_q<\varepsilon
\sqrt{|\fourier{\psi}_2(0)|^2+|\fourier{\psi}_1(-\kappa)|^2}$$
on the small handle with index $\kappa$
with respect to the wave vector
$k$, which coincides at the border to
$\Set{V}^+_{\varepsilon,\delta}$ with
the wave vector of $\Set{V}^+_{\varepsilon,\delta}$, respectively.
\qed
\end{Theorem}

The \Em{complex Fermi curves} $\fermi(V,W)$ are
locally finite sheeted coverings over $\xx{p}\in\mathbb{C}$ or
$\yy{p}\in\mathbb{C}$. If we enlarge the imaginary part of
$\xx{p}$ (or $\yy{p}$) we may connect any element of
$\fermi(V,W)$ by a path either with some
$\Set{V}^-_{\varepsilon,\delta}$ or with some
$\Set{V}^+_{\varepsilon,\delta}$.
Therefore, Theorem~\ref{asymptotic analysis 1} implies

\begin{Corollary}\label{regular fermi}
\index{connected components!of the normalization}
\index{normalization!connected components of the $\sim$}
The regular part of $\fermi(V,W)/\lattice\dual$
has at most two connected
components. Each component contains one of the two sets
$\Set{V}^{\pm}_{\varepsilon,\delta}$ described in
Theorem~\ref{asymptotic analysis 1}.
\qed
\end{Corollary}

This implies that the normalization of a \Em{complex Fermi curve} also
has at most two connected components, and each component contains one
of the smooth open subsets $\Set{V}^{\pm}_{\varepsilon,\delta}$.
The fixed energy level of the \Em{complex Bloch variety}
$\bloch(V,W)$ corresponding to the eigenvalue $\lambda$ is equal
to the \Em{complex Fermi curve} $\fermi(V+\lambda,W+\lambda)$.
Since the subset  of the \Em{complex Bloch variety} corresponding to
degenerated eigenvalues is a subvariety,
Theorem~\ref{asymptotic analysis 1} implies

\begin{Corollary}\label{degenerate eigenvalues}
The complement of the subvariety of degenerated eigenvalues is open
and dense in $\bloch(V,W)$.
\qed
\end{Corollary}

We shall see in Section~\ref{subsection compactified bounded genus}
that if the normalization of a \Em{complex Fermi curve}
$\fermi(V,W)/\lattice\dual$ has disconnected normalization,
then the \Em{first integral}
$4\int\limits_{\torus}V(x)W(x)d^2x$ is a multiple of
$4\pi$. Therefore, Theorem~\ref{asymptotic analysis 1} implies finally
that the regular part of the \Em{complex Bloch variety}
$\bloch(V,W)/\lattice\dual$ is connected
(i.\ e.\ the \Em{complex Bloch variety}
$\bloch(V,W)/\lattice\dual$ is irreducible).

\subsection{Several reductions}\label{subsection reductions}

In the last two sections we investigated the spectrum of the
Dirac operator with two periodic complex--valued potentials
$V$ and $W$ and introduced the \Em{complex Fermi curve} being the
spectral curve of the integrable system described by the
Davey--Stewartson equation. The Weierstra{\ss} representation of
immersions into $\mathbb{R}^3$ yields a one--to--one correspondence
of Dirac operators $\left(\begin{smallmatrix}
U & \partial\\
-\Bar{\partial} & U
\end{smallmatrix}\right)$ with real potentials $U$ on
$\torus$, whose kernel contains non--trivial elements
obeying the \De{Periodicity condition}~\ref{periodicity condition}.
Hence we shall reduce the integrable system,
whose spectral curves are the \Em{complex Fermi curves}
$\fermi(V,W)/\lattice\dual$, to the sub--system, whose spectral
curves are the \Em{complex Fermi curves}
$\fermi(U,U)/\lattice\dual$ of a real potential $U$.
Typically the reduced system is the set of fixed points of one or
several involutions. Moreover, these involutions induce involutions of
the corresponding spectral curves, and the spectral curves of the
fixed points of the former involutions are invariant
under the latter involutions.

In the sequel we will also consider the transposed Dirac operator
$\Op{D}^{t}(V,W,k)$. We will mostly denote the eigenfunction
of the Dirac operator by
\index{eigenfunction ($\rightarrow$ Baker--Akhiezer function)}
\index{eigenfunction ($\rightarrow$ Baker--Akhiezer function)!$\psi$
  of the Dirac operator}
\index{Dirac operator $\Op{D}$!eigenfunction $\psi$ of the $\sim$}
$\psi=\left(\begin{smallmatrix}
\psi_1\\
\psi_2
\end{smallmatrix}\right)$
and the eigenfunction of the transposed Dirac operator by
\index{eigenfunction ($\rightarrow$ Baker--Akhiezer function)! $\phi$
  of the transposed Dirac operator}
$\phi=\left(\begin{smallmatrix}
\phi_1\\
\phi_2
\end{smallmatrix}\right)$.
The transposition is always the one induced by the symmetric bilinear form
$$\langle\langle\phi,\psi\rangle\rangle
=\int\limits_{\Delta} 
\left(\phi_1(x)\psi_1(x)+\phi_2(x)\psi_2(x)\right)d^2x$$
on $\banach{2}(\Delta)\times\banach{2}(\Delta)$.
Consequently the pointwise complex conjugation on this Hilbert space
is used. Obviously the differential operators $\partial_{[k]}$ and
$\Bar{\partial}_{[k]}$ satisfy the relations
$$\partial_{[k]}^{t}=-\partial_{[-k]},
\Bar{\partial}_{[k]}^{t}=-\Bar{\partial}_{[-k]}\text{ and }
\partial_{[k]}^{\ast}=-\Bar{\partial}_{[\Bar{k}]},
\Bar{\partial}_{[k]}^{\ast}=-\partial_{[\Bar{k}]}.$$
Hence we have
\begin{eqnarray*}
\Op{D}^{t}(V,W,k) &=&\Op{J}^{-1}\comp\Op{D}(W,V,-k)\comp\Op{J},\\
\Op{D}^{\ast}(V,W,k) &=&\Op{D}(\Bar{V},\Bar{W},\Bar{k})
\text{, and}\\
\Bar{\Op{D}}(V,W,k) &=&
\Op{J}^{-1}\comp\Op{D}(\Bar{W},\Bar{V},-\Bar{k})\comp\Op{J},\\
\end{eqnarray*}
where $\Op{J}$ denotes the operator
$\left(\begin{smallmatrix}
0 & \unity\\
-\unity & 0
\end{smallmatrix}\right)$. This implies the following

\begin{Lemma} \label{transformations}
\begin{description}
\item[(i)] \index{involution!$\sigma$}
The holomorphic involution 
$\sigma:\mathbb{C}^3\rightarrow \mathbb{C}^3,
(k,\lambda)\mapsto (-k,\lambda)$, induces a
biholomorphic isomorphism of the complex Bloch varieties
$\bloch(V,W)\simeq \bloch(W,V)$.
\item[(ii)]\index{involution!$\rho$}
The anti--holomorphic involution 
$\rho:\mathbb{C}^3\rightarrow \mathbb{C}^3,
(k,\lambda)\mapsto (\Bar{k},\Bar{\lambda})$,
induces a anti--biholomor\-phic isomorphism of the
complex Bloch varieties
$\bloch(V,W)\simeq \bloch(\Bar{V},\Bar{W})$.
\item[(iii)]\index{involution!$\eta$}
The anti--holomorphic involution 
$\eta:\mathbb{C}^3\rightarrow \mathbb{C}^3,
(k,\lambda)\mapsto (-\Bar{k},\Bar{\lambda})$,
induces a anti--biholomorphic isomorphism of the
complex Bloch varieties
$\bloch(V,W)\simeq \bloch(\Bar{W},\Bar{V})$.
\qed
\end{description}
\end{Lemma}

\begin{Corollary} \label{involutions}
\index{involution!of the complex Bloch variety}
\begin{description}
\item[(i)]
The complex Bloch variety $\bloch(U,U)$ has a
holomorphic involution $\sigma: (k,\lambda)\mapsto (-k,\lambda)$.
Moreover, if $\psi$ is an eigenfunction corresponding to the
element $(k,\lambda)$ of $\bloch(U,U)$ then $\phi=\Op{J}\psi$ is an
eigenfunction of the transposed Dirac operator corresponding to the
element $\sigma(k,\lambda)=(-k,\lambda)$, and vice versa.
\item[(ii)]
If $V$ and $W$ are real valued periodic functions, the
complex Bloch variety $\bloch(V,W)$ has an
anti--holomorphic involution
$\rho:(k,\lambda)\mapsto(\Bar{k},\Bar{\lambda})$. Moreover, if
$\psi$ is an eigenfunction corresponding to the element $(k,\lambda)$
of $\bloch(V,W)$ then $\phi=\Bar{\psi}$ is an eigenfunction of
the transposed Dirac operator corresponding to the element
$\rho(k,\lambda)=(\Bar{k},\Bar{\lambda})$, and vice versa.
\item[(iii)]
The complex Bloch variety $\bloch(U,\Bar{U})$ has an
anti--holomorphic involution
$\eta:(k,\lambda)\mapsto(-\Bar{k},\Bar{\lambda})$. Moreover, if
$\psi$ is an eigenfunction corresponding to the element $(k,\lambda)$
of $\bloch(U,\Bar{U})$ then $\Op{J}\Bar{\psi}$
is an eigenfunction of the complex conjugate of the Dirac operator
corresponding to the element
$\eta(k,\lambda)=(\Bar{k},\Bar{\lambda})$, and vice versa.
\qed
\end{description}
\end{Corollary}

The holomorphic involution $\sigma$ induces a holomorphic involution
of $\fermi(U,U)$, $\fermi(U,U)/\lattice\dual$, and of the
normalization of $\fermi(U,U)/\lattice\dual$.
All of them are also denoted by $\sigma$.
If both potentials  $V$ and $W$ are real, then the
anti--holomorphic involution $\rho$
induces an involution of $\fermi(V,W)$,
$\fermi(V,W)/\lattice\dual$
and the normalization of
$\fermi(V,W)/\lattice\dual$.
Again these involutions are denoted by $\rho$.
We call the fixed points of the involution $\rho$
the real part\index{real!part} of the \Em{complex Fermi curve}.
Consequently the real part of the
normalization of the \Em{complex Fermi curve} consists of the
fixed points of the corresponding involution on the normalization.

\begin{Corollary} \label{fixed points}
\begin{description}
\item[(i)]
If $V$ and $W$ are real potentials, then all elements of the
preimage of $\fermi(V,W)/\lattice\dual\cap
\torus\dual$ under the normalization map, which are
not fixed points of the involution $\rho$, are isolated points of the
preimage of $\fermi(V,W)/\lattice\dual\cap
\torus\dual$ under the normalization map.
\item[(ii)]
The involution $\eta$
of the normalization of $\fermi(U,\Bar{U})/\lattice\dual$ does
not have any fixed point. More precisely, for all elements $y$ of the
normalization of $\fermi(U,U)/\lattice\dual$ the eigenfunctions
corresponding to $y$ and $\eta(y)$ are linearly independent.
\item[(iii)]
If $\fermi(U,\Bar{U})/\lattice\dual$
contains an element of the form $[k]$ with
$2\Re(k)\in\lattice\dual$, then the
dimension of the eigenspace corresponding to this element is
even.
\end{description}
\end{Corollary}

\begin{proof} All elements of the real part of
  $\fermi(V,W)/\lattice\dual$ are obviously fixed points of the
  involution $\rho$. In particular, all regular points of the real part of
  $\fermi(V,W)/\lattice\dual$ are fixed points of $\rho$. Since the
  singular points of $\fermi(V,W)/\lattice\dual$ are isolated, this
  implies (i). Let $y$ be some element of the
  normalization of $\fermi(U,\Bar{U})/\lattice\dual$,
  which is a fixed point of the involution $\eta$. Due to
  Lemma~\ref{involutions} the eigenfunction $\psi$ corresponding to
  this point has to be a multiple of $\Op{J}\Bar{\psi}$. But the equation
  $\psi=\alpha \Op{J}\Bar{\psi}$ implies $\psi=\alpha \Bar{\alpha}
  \Op{J}\overline{\left(\Op{J}\Bar{\psi}\right)}$. This would imply
  $\alpha\Bar{\alpha}=-1$, which is impossible. This proves (ii).
  All elements of $\fermi(U,\Bar{U})/\lattice\dual$
  of the form $[k]$ with $2\Re(k)\in \lattice\dual$ are
  fixed points of the involution $\eta$.
  Due to Lemma~\ref{involutions} the
  eigenspace corresponding to this element $[k]$ is invariant
  under the anti--unitary map $\psi\mapsto\Op{J}\Bar{\psi}$.
  The square of this map is equal to $-\unity$.
  Hence for all elements $\psi$ of the
  eigenspace corresponding to $[k]$ the subspace spanned by
  $\psi$ and $\Op{J}\Bar{\psi}$ is invariant under this map and
  two--dimensional, because $\Op{J}\Bar{\psi}=\alpha\psi$ would imply
  $\alpha\Bar{\alpha}=-1$, which is impossible. Since the eigenspace
  corresponding to $[k]$ has to be a direct sum of subspaces of
  this form, its dimension is even.
\end{proof}

Besides these involutions there exists some transformations
of the potentials, which do not change the complex Bloch variety.
In particular, all automorphisms of the group
$\torus$ induce transformations of the potentials, which
do not change the complex Bloch variety. Moreover, the translations
also induce transformations of the potentials, which do not change
the complex Bloch variety (the corresponding transformations of
the eigen bundle plays an important role in the
inverse spectral theory (see e.\ g.\ \cite{Kr1,DKN}).
Furthermore, for all $z\in\mathbb{C}^{\ast}$ the
Dirac operator transforms like
$$\Op{D}(z^2V,z^{-2}W,k)=
\left(\begin{smallmatrix}
z & 0\\
0 & z^{-1}
\end{smallmatrix}\right)\comp
\Op{D}(V,W,k)\comp
\left(\begin{smallmatrix}
z & 0\\
0 & z^{-1}
\end{smallmatrix}\right).$$
\index{isospectral!transformation}
Therefore, the \Em{complex Fermi curves} $\fermi(zV,z^{-1}W)$ are
independent of $z\in\mathbb{C}^{\ast}$.
This action of the group $\mathbb{C}^{\ast}$ on the
\Em{isospectral sets} (i.\ e.\ the set of potentials, which have the same
\Em{complex Fermi curves}) is only the simplest example of a hierarchy of
isospectral transformations \cite{GS2}, which may be constructed
with the tools of soliton theory \cite{Kr1,DKN}.

Finally, we remark that all affine symmetries of the torus
$\torus$ induce symmetries of the corresponding complex Bloch varieties.
Besides the translations there exist in general only one,
namely the parity transformation:
\index{parity transformation $\parity$|(}
$$\torus\rightarrow \torus,\;
x\mapsto-x.$$
This symmetry induces an isometry on $\banach{p}(\torus)$ and
$\banach{p}(\Delta)$, which is denoted by $\parity$.
By definition the operators $\partial_{[k]}$ and
$\Bar{\partial}_{[k]}$ transform as follows under this symmetry:
$$\parity\comp\partial_{[k]}\comp\parity=-\partial_{[-k]}=
\partial_{[k]}^{t} \text{ and }
\parity\comp\Bar{\partial}_{[k]}\comp\parity=
-\Bar{\partial}_{[-k]}=\Bar{\partial}_{[k]}^{t}.$$
This implies the following relations:
\begin{eqnarray*}
\left(\begin{smallmatrix}
0 & \parity\\
\parity & 0
\end{smallmatrix}\right)
\comp\Op{D}(V,W,k)\comp
\left(\begin{smallmatrix}
0 & \parity\\
\parity & 0
\end{smallmatrix}\right) &=&
\Op{D}^{t}(\parity(W),\parity(V),k),\\
\left(\begin{smallmatrix}
0 & \parity\\
\parity & 0
\end{smallmatrix}\right)
\comp\Bar{\Op{D}}(V,W,k)\comp
\left(\begin{smallmatrix}
0 & \parity\\
\parity & 0
\end{smallmatrix}\right) &=&
\Op{D}(\parity(\Bar{W}),\parity(\Bar{V}),\Bar{k}),\\
\left(\begin{smallmatrix}
\parity & 0\\
0 & -\parity
\end{smallmatrix}\right)
\comp\Op{D}(V,W,k)\comp
\left(\begin{smallmatrix}
\parity & 0\\
0 & -\parity
\end{smallmatrix}\right) &=&
\Op{D}(\parity(V),\parity(W),-k),\\
\left(\begin{smallmatrix}
\parity & 0\\
0 & -\parity
\end{smallmatrix}\right)
\comp\Bar{\Op{D}}(V,W,k)\comp
\left(\begin{smallmatrix}
\parity & 0\\
0 & -\parity
\end{smallmatrix}\right) &=&
\Op{D}^{t}(\parity(\Bar{V}),\parity(\Bar{W}),-\Bar{k}).
\end{eqnarray*}
Hence we obtain the following

\begin{Lemma}\label{parity}
The involutions $\unity$, $\rho$, $\sigma$ and $\eta$ induce
isomorphisms of the complex Bloch varieties
$\bloch(V,W)\simeq
\bloch(\parity(W),\parity(V))$,
$\bloch(V,W)\simeq
\bloch(\parity(\Bar{W}),\parity(\Bar{V}))$,
$\bloch(V,W)\simeq
\bloch(\parity(V),\parity(W))$ and
$\bloch(V,W)\simeq
\bloch(\parity(\Bar{V}),\parity(\Bar{W}))$,
respectively. More precisely, if $\psi$ is an eigenfunction of
$\Op{D}(V,W,k)$ with eigenvalue $\lambda$, then
$\left(\begin{smallmatrix}
0 & \parity\\
\parity & 0
\end{smallmatrix}\right)
\Bar{\psi}$ and
$\left(\begin{smallmatrix}
\parity & 0\\
0 & -\parity
\end{smallmatrix}\right)
\psi$ are eigenfunctions of the operators
$\Op{D}(\parity(\Bar{W}),\parity(\Bar{V}),\Bar{k})$ and 
$\Op{D}(\parity(V),\parity(W),-k)$ with eigenvalues
$\Bar{\lambda}$ and $\lambda$, respectively. Moreover,
$\left(\begin{smallmatrix}
0 & \parity\\
\parity & 0
\end{smallmatrix}\right)
\psi$ and
$\left(\begin{smallmatrix}
\parity & 0\\
0 & -\parity
\end{smallmatrix}\right)
\Bar{\psi}$ are transposed eigenfunctions of the operators
$\Op{D}(\parity(W),\parity(V),k)$ and 
$\Op{D}(\parity(\Bar{V}),\parity(\Bar{W}),-\Bar{k})$
with eigenvalues $\lambda$ and $\Bar{\lambda}$, respectively.
In particular, if the pair of potentials is invariant under the
transformation
$(V,W)\mapsto(\parity(W),\parity(V))$,
$(V,W)\mapsto(\parity(\Bar{W}),\parity(\Bar{V}))$,
$(V,W)\mapsto(\parity(V),\parity(W))$ and
$(V,W)\mapsto(\parity(\Bar{V}),\parity(\Bar{W}))$,
then the complex Bloch variety is invariant under the corresponding
involution $\unity$, $\rho$, $\sigma$ and $\eta$,
respectively.\qed
\index{parity transformation $\parity$|)}
\end{Lemma}

\subsection{Spectral projections}\label{subsection spectral projections}

In Section~\ref{subsection fermi curve} the \Em{complex Fermi curve} was
introduced as a subvariety of $\mathbb{C}^2$. In this section we will
work out a slightly different representation of this
\Em{complex Fermi curve}, which in
Section~\ref{subsection weak singularity} and
Section~\ref{subsection singularity} 
will turn out to be the key for the understanding
of the \De{Periodicity condition}~\ref{periodicity condition}.
As a preparation we investigate the spectral projections of the
Dirac operator in some more detail.

Due to Theorem~\ref{meromorph} and
Corollary~\ref{degenerate eigenvalues}
locally there exist two meromorphic
functions $\psi$ and $\phi$ on the complex Bloch variety
$\bloch(V,W)/\lattice\dual$, which map the elements
$([k],\lambda)$ onto
an eigenfunction of the Dirac operator $\Op{D}(V,W,k)$ and an
eigenfunction of the transposed Dirac operator $\Op{D}^{t}(V,W,k)$,
respectively. Moreover, these functions are unique up to
multiplication by some invertible meromorphic function.
Now we define a meromorphic function on
$\bloch(V,W)/\lattice\dual$ with values in the
finite rank operators on
$\banach{2}(\Delta)\times\banach{2}(\Delta)$ by
\index{spectral!projection!$\Op{P}$}
  $$\Op{P}([k],\lambda):\begin{pmatrix}
  \chi_1\\
  \chi_2
  \end{pmatrix}\longmapsto
  \frac{\langle\langle\phi([k],\lambda),
  \chi\rangle\rangle}
  {\langle\langle\phi([k],\lambda),
  \psi([k],\lambda)\rangle\rangle}
  \begin{pmatrix}
  \psi_1([k],\lambda)\\
  \psi_2([k],\lambda)
  \end{pmatrix}.$$
Obviously this definition does not depend on the normalization of the
functions $\psi$ and $\phi$. In fact, if we multiply both functions
with some non--vanishing meromorphic functions, then the operator
$\Op{P}$ does not change.
Moreover, due to Corollary~\ref{degenerate eigenvalues} the
denominator does not vanish identically on the
complex Bloch variety. Thus $\Op{P}$ is a well defined
global meromorphic function on $\bloch(V,W)/\lattice\dual$.
We will see that in some sense this projection can be considered
as a projection--valued regular form.

On varieties with singularities the suitable generalization of
holomorphic forms are the regular forms
(\cite[Chapter~IV \S3.]{Se} and \cite{Kun}).
In general the corresponding sheaf is not
the generalized cotangent sheaf
\cite[Chapter~II. \S1.2.]{GPR}.
Of particular interest are the regular forms of degree equal to
the dimension of the variety. The sheaf of these regular forms is called
dualizing sheaf \cite[Chapter~II. \S5.3.]{GPR}
\index{sheaf!dualizing $\sim$}. In our case there
exists a simple definition of the regular forms of degree equal
to the dimension of the variety. In order to keep this paper
self--contained we include a proof of the following lemma,
which is a special case of a more general statement
\cite[Chapter~II. Lemma~5.23]{GPR}.

\begin{Remark}\label{local sum}
In the sequel we shall meet quite often complex spaces, which are
locally biholomorphic to finite sheeted coverings over open subsets of
$\mathbb{C}^n$. If we restrict these coverings to the preimage of
small open balls, then different sheets without branch points are not
connected with each other. However, in arbitrary small neighbourhoods
of branch points several sheets are connected. In the sequel we shall
call those sheets, whose restrictions to arbitrary small
neighbourhoods of a given element contain this element,
the local sheets, which contain this element.
\end{Remark}

\begin{Lemma} \label{dualizing sheaf}
Let $R(z_0,\ldots,z_n)$ be a holomorphic function on some open subset
$U\subset \mathbb{C}^{n+1}$, whose partial derivative
$\partial R/\partial z_0$ is not
identically zero on the connected components of the subvariety
of $U$ defined by the equation $R(z_0,\ldots,z_n)=0$. Therefore
this subvariety locally may be considered as a covering space over
$(z_1,\ldots,z_n)\in\mathbb{C}^n$. The following conditions on a
meromorphic function $f$ on this subvariety are equivalent:
\begin{description}
\item[(i)] The function $f$ is holomorphic.
\item[(ii)] For all holomorphic functions $g$ the local sum of
  $gf/(\partial R/\partial z_0)$ over all
  sheets of the subvariety considered as a covering space over
  $(z_1,\ldots,z_n)\in\mathbb{C}^n$, which contain an arbitrary
  element (compare with Remark~\ref{local sum}), is a holomorphic function.
\end{description}
\end{Lemma}

\begin{proof} Let us choose some element $u$ of the subvariety. Due to
  the Weierstra{\ss} Preparation theorem
  \cite[Chapter~I. Theorem~1.4]{GPR} locally there exists a
  polynomial $Q$ with respect to $z_0$, whose coefficients are
  holomorphic functions depending on $z_1,\ldots,z_n$, such that $R/Q$
  is locally holomorphic. Moreover, the degree $d$ of $Q$ may be
  chosen to be equal to the number of all sheets of the subvariety
  considered as a covering space over
  $(z_1,\ldots,z_n)\in\mathbb{C}^n$, which contain the element
  $u$. Near this element $u$ the condition (ii) is
  obviously equivalent to an analogous condition, where $R$ is
  replaced by $Q$. Let $Q(z_0)$ be an arbitrary polynomial of
  degree $d$, whose coefficients are complex numbers, and let
  $u_{0,1},\ldots,u_{0,d}$ be the zeroes of this polynomial. For an
  arbitrary polynomial $g(z_0)$ of degree less than $d$ and with
  complex--valued coefficients we have the identity
  $$g(z_0) = \sum\limits_{i=1}^{d}\frac{g(u_{0,i})}
  {\partial Q(u_{0,i})/\partial z_0}
  \prod\limits_{j\neq i} (z_0-u_{0,j}).$$
  In fact, if all zeroes of $Q$ are pairwise different, the polynomials
  on both sides take the same values at all these zeroes and therefore
  have to be equal. In particular, the sum $\sum\limits_{i}^{d}
  \frac{g(u_{0,i})}{\partial Q(u_{0,i})/\partial z_0}$
  is equal to the coefficient of the monomial $z_0^{d-1}$ in $g(z_0)$.
  For the monomials $g(z_0)=z_0^l$ this implies that
  $$\sum\limits_{i=1}^{d}\frac{u_{0,i}^l}
  {\partial Q(u_{0,i})/\partial z_0} =
  \begin{cases} 0 & \text{if } l<d-1 \\
  1 & \text{if } l=d-1.
  \end{cases} $$
  We conclude that the same is true if the coefficients of $Q$ are
  holomorphic functions depending on $z_1,\ldots,z_n$, such that the
  discriminant of $Q$ does not vanish identically.
  Due to the Weierstra{\ss} Preparation theorem
  \cite[Chapter~I. Theorem~1.4]{GPR} each meromorphic function
  on the subvariety defined by the equation $Q=0$
  may be written locally
  as a polynomial with respect to $z_0$ of degree $d-1$,
  whose coefficients are meromorphic functions
  depending on $z_1,\ldots,z_n$.
  Hence the local sum of this function
  over all sheets of the subvariety considered as a covering space
  over $(z_1,\ldots,z_n)\in\mathbb{C}^n$, which contain  $u$,
  is locally a holomorphic function,
  if and only if the coefficient corresponding to
  $z_0^{d-1}$ of the polynomial is locally holomorphic.
  This shows that (i) and (ii) are equivalent.
\end{proof}

If besides the partial derivative $\partial R/\partial z_0$ some
other partial derivative $\partial R/\partial z_i$ is also not
identically zero on the connected components of the subvariety
of $U$ defined by the equation $R(z_0,\ldots,z_n)=0$,
then we may consider this subvariety locally either as a covering
space over 
$(z_1,\ldots,z_n)\in\mathbb{C}^n$ or as a covering space over
$(z_0,\ldots,z_{i-1},z_{i+1},\ldots,z_n)\in\mathbb{C}^n$.
The equation $R(z_0,\ldots ,z_n)=0$ implies the relation
$(\partial R/\partial z_0) dz_0 + \ldots + 
(\partial R/\partial z_n) dz_n = 0$. We conclude that the $n$--form
$\frac{1}{\partial R/\partial z_0} dz_1 \wedge \ldots \wedge dz_n$ is
equal to the $n$--form
$\frac{(-1)^i}{\partial R/\partial z_i}
dz_0 \wedge \ldots \wedge dz_{i-1} \wedge dz_{i+1} \wedge \ldots
\wedge dz_{n-1}$.
Consequently we call a $n$--form
\begin{multline*}
\omega =
\frac{f(z_0,\ldots,z_n)}{\partial R(z_0,\ldots,z_n)/\partial z_0}
dz_1\wedge \ldots \wedge dz_n = \\
(-1)^i
\frac{f(z_0,\ldots,z_n)}{\partial R(z_0,\ldots,z_n)/\partial z_i}
dz_0 \wedge \ldots \wedge dz_{i-1} \wedge dz_{i+1} \wedge \ldots
\wedge dz_{n-1}
\end{multline*}
regular\index{form!regular $\sim$},
if and only if $f$ satisfies one of the equivalent
conditions (i), (ii) or the analogous condition to (ii), when $z_0$ is
replaced by $z_i$ (\cite[Chapter~IV \S3.]{Se},
\cite{Kun} and \cite[Chapter~II. Definition~5.22]{GPR}).

\begin{Lemma} \label{projection 1}
\index{spectral!projection!$\Op{P}$}
The function $\Op{P}$ has the following properties:
\begin{description}
\item[(i)] The values of $\Op{P}$ are projections of rank one.
\item[(ii)] If $([k],\lambda)$ and $([k],\lambda')$ are two elements
  of the complex Bloch variety with $\lambda\neq\lambda'$, then
  $\Op{P}([k],\lambda)\comp\Op{P}([k],\lambda') = 0 =
  \Op{P}([k],\lambda')\comp\Op{P}([k],\lambda)$.
\item[(iii)] Due to Theorem~\ref{meromorph} the
  complex Bloch variety is locally a Weierstra{\ss} covering
  (see e.g. \cite[Chapter~I. \S12.3]{GPR}) over $k\in\mathbb{C}^2$.
  The local sum of $\Op{P}$ over all sheets of this covering,
  which contain some element
  $([k],\lambda)\in \bloch(V,W)/\lattice\dual$
  (compare with Remark~\ref{local sum}),
  is a holomorphic function on some neighbourhood of
  $k\in\mathbb{C}^2$. Moreover, the values of this function are
  projections,  whose rank is equal to the number of sheets, over
  which that sum is taken. Finally, the value of this function at
  $([k],\lambda)$ is equal to the spectral projection
  $\Breve{\Op{P}}([k],\lambda)$ of the Dirac operator
  $\Op{D}(V,W,k)$ onto the
  generalized eigenspace associated with eigenvalue $\lambda$
  (see e.\ g.\ \cite[Appendix to XII.1]{RS4}).
\item[(iv)] The projection--valued form
$\Op{P} dk_1\wedge dk_2$ is regular.
\end{description}
\end{Lemma}

\begin{proof} The first statement is obvious. Due to
  Theorem~\ref{meromorph} locally we may replace the
  Dirac operator $\Op{D}(V,W,k)$ by some holomorphic
  matrix--valued function
  $\Mat{A}(k)$ on some open subset $k\in\mathbb{C}^2$. Due to
  Corollary~\ref{degenerate eigenvalues}
  this function generically has pairwise different eigenvalues.
  If $\Mat{A}$ is a $l\times l$ matrix with
  pairwise different eigenvalues $\lambda_1,\ldots,\lambda_l$,
  then the eigenvectors $\psi_1,\ldots,\psi_l$ to the eigenvalues
  $\lambda_1,\ldots,\lambda_l$ form a basis of $\mathbb{C}^2$.
  Moreover, the dual basis $\phi_1,\ldots,\phi_l$ gives the
  eigenvectors of $\Mat{A}^{t}$ to the eigenvalues
  $\lambda_1,\ldots,\lambda_l$. In fact, if $\Tilde{\phi}_i$ is some
  eigenvector of $\Mat{A}^{t}$ with eigenvalue $\lambda_i$, then we have
  $\Tilde{\phi}_i^{t}\psi_j(\lambda_i-\lambda_j)=
  \left(\Mat{A}^{t}\Tilde{\phi}_i\right)^{t}\psi_j -
  \Tilde{\phi}_i^{t}\Mat{A}\psi_j = 0$ for all $1\leq i,j\leq l$ and
  $\Tilde{\phi}_i$ has to be proportional to $\phi_i$. We conclude that
  for all $1\leq i\neq j\leq l$ the eigen projections
  $$\Op{P}_i = \frac{\psi_i\phi_i^{t}}{\phi_i^{t}\psi_i}
  \; i=1,\ldots,l$$
  satisfy the relations $\Op{P}_i\comp\Op{P}_j = 0 =
  \Op{P}_j\comp\Op{P}_i$. This proves
  (ii). On the other hand $\Op{P}_1+\ldots+\Op{P}_l$
  is equal to the identity matrix and (iii) is proven.
  Due to Theorem~\ref{meromorph} the
  complex Bloch variety is locally a Weierstra{\ss} covering (see
  e.g. \cite[Chapter~I. \S12.3]{GPR}) over $k\in\mathbb{C}^2$. Hence
  any holomorphic function $f$ on the complex Bloch variety is locally
  equal to some polynomial with respect to $\lambda$, whose
  coefficients are holomorphic functions on some open
  subset of $k\in\mathbb{C}^2$.
  Again due to Theorem~\ref{meromorph} locally the sum over the sheets
  of this Weierstra{\ss} covering of $\Op{P}$
  times any holomorphic function is a holomorphic function on some
  open subset of $k\in\mathbb{C}^2$ with values in the
  finite rank operators on
  $\banach{2}(\Delta)\times\banach{2}(\Delta)$. Hence
  the function $\Op{P}$ satisfies condition (ii) of
  Lemma~\ref{dualizing sheaf}.
\end{proof}

Since the Dirac operator is a holomorphic unbounded--operator--valued
function depending on the wave vectors $k\in\mathbb{C}^2$, the
complex Bloch variety has naturally the structure of a
covering space over $k\in \mathbb{C}^2$. But in order to understand the
\Em{complex Fermi curve}, which is our main interest, it would be more
appropriate to consider the complex Bloch variety as a
covering space over $\lambda$ and
some component of $k$, because such a covering would just have to be
restricted to the plane $\lambda=0$ to obtain the
\Em{complex Fermi curve}. In the remainder of this section
this point of view will be
worked out. The starting point is the following observation, which
generalizes the equivalence of the
\De{Periodicity condition}~\ref{periodicity condition}. Here we use the
notations introduced in the context of the
\De{Fundamental domain}~\ref{fundamental domain}, but we shall omit
the index $\Delta$.

\begin{Lemma} \label{branchpoints}
Let $\psi(x)=\left(\begin{smallmatrix}
\psi_1(x)\\
\psi_2(x)
\end{smallmatrix}\right)$ be an eigenfunction corresponding to an
element $(k,\lambda)$ of $\bloch(V,W)$ and
$\phi(x)=\left(\begin{smallmatrix}
\phi_1(x)\\
\phi_2(x)
\end{smallmatrix}\right)$ an eigenfunction of the transposed
Dirac operator corresponding to an element $(k',\lambda)\in
\bloch(V,W)$. Then $\phi_2\psi_1 dz + \phi_1\psi_2
d\Bar{z}$, with $dz = dx_1+\sqrt{-1}dx_2$ and
$d\Bar{z}=dx_1-\sqrt{-1}dx_2$, is a closed form. If the
two components $\xx{p}=g(\xx{\gamma},k)$ and
$\xx{p}'=g(\xx{\gamma},k')$ of $k$ and $k'$ coincide,
then the integral of this form from some point
$x\in\mathbb{R}^2$ to $x+\xx{\gamma}$ does not depend on $x$,
and is equal to
$$\frac{\left(\xx{\gamma}_1+\sqrt{-1}\xx{\gamma}_2\right)
\langle\phi_2,\psi_1\rangle +
\left(\xx{\gamma}_1-\sqrt{-1}\xx{\gamma}_2\right)
\langle\phi_1,\psi_2\rangle}
{\xx{\gamma}_1\yy{\gamma}_2-
\xx{\gamma}_2\yy{\gamma}_1}.$$
If in addition the elements $[k]$ and $[k']$ of
$\mathbb{C}^2/\lattice\dual$ are different, then both integrals vanish.
\end{Lemma}

\begin{proof} Since $\psi$ is an eigenfunction of the Dirac operator
  and $\phi$ is a eigenfunction of the transposed Dirac operator with
  eigenvalues zero, we have
\begin{align*}
V\psi_1 + \partial \psi_2 &=\lambda \psi_1 & 
V\phi_1 + \Bar{\partial}\phi_2 &=\lambda \phi_1\\
-\Bar{\partial}\psi_1 + W\psi_2 &=\lambda \psi_2 &
-\partial\phi_1 + W\phi_2 &=\lambda \phi_2.
\end{align*}
Thus
$$d\left(\phi_2\psi_1 dz + \phi_1\psi_2 d\Bar{z} \right)=
\left(\left(\Bar{\partial}\phi_2\right)\psi_1 +
\phi_2\left(\Bar{\partial}\psi_1\right) -
\left(\partial\phi_1\right)\psi_2 -
\phi_1\left(\partial\psi_2\right) \right) d\Bar{z}\wedge dz= 0.$$
If the two components $\xx{p}$ and $\xx{p}'$ of $k$ and $k'$
coincide, then this closed
form is invariant under the shift by the period
$\xx{\gamma}$. This implies that the integral of this form
along some path from some arbitrary point $x\in\mathbb{R}^2$ to
$x+\xx{\gamma}$, does not depend on $x$. In particular, this
integral is equal to the integral
\begin{multline*}
\int\limits_{[0,1]^2}
\left(\left(\xx{\gamma}_1+\sqrt{-1}\xx{\gamma}_2\right)
\phi_2(x)\psi_1(x) +
\left(\xx{\gamma}_1-\sqrt{-1}\xx{\gamma}_2\right)
\phi_1(x)\psi_2(x)\right)
d\xx{q}\wedge d\yy{q} =\\
\frac{\left(\xx{\gamma}_1+\sqrt{-1}\xx{\gamma}_2\right)
\langle\phi_2,\psi_1\rangle +
\left(\xx{\gamma}_1-\sqrt{-1}\xx{\gamma}_2\right)
\langle\phi_1,\psi_2\rangle}
{\xx{\gamma}_1\yy{\gamma}_2-
\xx{\gamma}_2\yy{\gamma}_1}\; .
\end{multline*}
If the components with respect to the other period
$\yy{\gamma}$ are different
$\yy{p}\neq \yy{p} \mod \mathbb{Z}$, then the shift by
$\yy{\gamma}$ changes the form by the factor
$\exp\left(2\pi\sqrt{-1}(\yy{p}-\yy{p}')\right)\neq 1$.
Hence the integral of this form along $\xx{\gamma}$ has to vanish.
\end{proof}

This lemma motivates the following definition: For all
$\mu\in\mathbb{C}^2$, let
$\langle\langle\cdot,\cdot\rangle\rangle_{\mu}$ be
the following bilinear form on the Hilbert spaces
$\banach{2}(\torus)\times\banach{2}(\torus)$
and $\banach{2}(\Delta)\times\banach{2}(\Delta)$, respectively:
$$\langle\langle\phi,\psi\rangle\rangle_{\mu} := 
\left(\mu_1+\sqrt{-1}\mu_2\right)\langle\phi_2,\psi_1\rangle +
\left(\mu_1-\sqrt{-1}\mu_2\right)\langle\phi_1,\psi_2\rangle .$$
The complex Bloch variety is locally a finite--sheeted covering over
$(\xx{p},\lambda)\in\mathbb{C}^2$, and the zeroes of
the differential $d\xx{p}\wedge d\lambda$ are the branch points of this
covering. Now we define another projection, which is related to this
covering. Again we use the notations introduced in the context of the
\De{Fundamental domain}~\ref{fundamental domain}.
Furthermore, we again use the meromorphic functions
$\psi=\left(\begin{smallmatrix}
\psi_1\\
\psi_2
\end{smallmatrix}\right)$ and
$\phi=\left(\begin{smallmatrix}
\phi_1\\
\phi_2
\end{smallmatrix}\right)$ on $\bloch(V,W)/\lattice\dual$.
Let $\Op{P}_{\xx{\gamma}}$ be the meromorphic map from
$\bloch(V,W)/\lattice\dual$ into the finite rank operators on the
Hilbert space $\banach{2}(\Delta)\times\banach{2}(\Delta)$, which for
all elements $([k],\lambda)$ of the complex Bloch variety is defined by
\index{spectral!projection!$\Op{P}_{\xx{\gamma}}$}
$$\Op{P}_{\xx{\gamma}}([k],\lambda):
\begin{pmatrix}
\chi_1\\
\chi_2
\end{pmatrix}\longmapsto
\frac{\langle\langle\phi([k],\lambda),
\chi\rangle\rangle_{\xx{\gamma}}}
{\langle\langle\phi([k],\lambda),
\psi([k],\lambda)\rangle\rangle_{\xx{\gamma}}}
\begin{pmatrix}
\psi_1([k],\lambda)\\
\psi_2([k],\lambda)
\end{pmatrix}.$$
Obviously this definition does not depend on the normalization of the
functions $\psi$ and $\phi$. In fact, if we multiply both functions
with some non--vanishing meromorphic functions, then the operator
$\Op{P}_{\xx{\gamma}}$ does not change. Moreover, due to
proof of Theorem~\ref{asymptotic analysis 1} and
Corollary~\ref{regular fermi} the denominator
does not vanish identically on the complex Bloch variety.
Thus $\Op{P}_{\xx{\gamma}}$ is a well defined
global meromorphic function on
$\bloch(V,W)/\lattice\dual$.

\begin{Lemma} \label{projection 2}
\index{spectral!projection!$\Op{P}_{\xx{\gamma}}$}
The function $\Op{P}_{\xx{\gamma}}$ has the following properties:
\begin{description}
\item[(i)] The values of $\Op{P}_{\xx{\gamma}}$ are
  projections of rank one.
\item[(ii)] If $([k],\lambda)$ and $([k'],\lambda)$ are two
  different elements of the complex Bloch variety with
  $g(\xx{\gamma},k)= g(\xx{\gamma},k')\mod
  \mathbb{Z}$, then $\Op{P}_{\xx{\gamma}}([k],\lambda)\comp
  \Op{P}_{\xx{\gamma}}([k'],\lambda) = 0 =
  \Op{P}_{\xx{\gamma}}([k'],\lambda)\comp
  \Op{P}_{\xx{\gamma}}([k],\lambda)$.
\item[(iii)] Locally the complex Bloch variety is a
  finite--sheeted covering over
  $(\xx{p},\lambda)\in\mathbb{C}^2$. The local sum of
  $\Op{P}_{\xx{\gamma}}$ over all sheets, which contain one element
  $([k],\lambda)$ (compare with Remark~\ref{local sum}),
  is a holomorphic function on some open subset of
  $(\xx{p},\lambda)\in\mathbb{C}^2$, with values in the
  finite rank projections on $\banach{2}(\Delta)\times\banach{2}(\Delta)$.
  Moreover, the rank of any of these projections is equal to the number of
  sheets, over which the sum is taken.
\item[(iv)] The $2$--form $\Op{P}_{\xx{\gamma}}
  d\xx{p}\wedge d\lambda$ is
  a regular form on the complex Bloch variety and the
  one--form $\Op{P}_{\xx{\gamma}} d\xx{p}$
  is a regular form\index{form!regular $\sim$}
  on the \Em{complex Fermi curve}.
\item[(v)] Due to (iii) for all elements $([k],\lambda)$ of
  $\bloch(V,W)/\lattice\dual$ there exists a unique projection
  $\Breve{\Op{P}}_{\xx{\gamma}}([k],\lambda)$,
  which is the value at
  $([k],\lambda)$ of the local sum of $\Op{P}_{\xx{\gamma}}$
  over all sheets of
  of $\bloch(V,W)/\lattice\dual$, which contain $([k],\lambda)$.
  For two different elements $([k],\lambda)$ and $([k'],\lambda)$ of
  $\bloch(V,W)/\lattice\dual$ with equal components
  $g(\xx{\gamma},k)=g(\xx{\gamma},k') \mod \mathbb{Z}$ the corresponding
  projections are disjoint:
  $\Breve{\Op{P}}_{\xx{\gamma}}([k],\lambda)\comp
  \Breve{\Op{P}}_{\xx{\gamma}}([k'],\lambda)=0=
  \Breve{\Op{P}}_{\xx{\gamma}}([k'],\lambda)\comp
  \Breve{\Op{P}}_{\xx{\gamma}}([k]\lambda)$. 
\item[(vi)] If $\chi$ is a proper eigenfunction of the Dirac operator
  corresponding to some element $([k],\lambda)$
  of the complex Bloch variety,
  then the range of $\Breve{\Op{P}}_{\xx{\gamma}}([k],\lambda)$
  contains $\chi$.
\end{description}
\end{Lemma}

We remark that the range of the spectral projections
contain all elements in the kernels of
$\left(\lambda\unity-\Op{D}(V,W,k)\right)^{l}$ with $l\in\mathbb{N}$.

\begin{proof} The first statement is obvious.
  Statement (ii) is a direct consequence of Lemma~\ref{branchpoints}.
  If we consider the Dirac operator as an operator
  acting on the Hilbert bundle introduced in the context of the
  \De{Trivialization}~\ref{trivialization},
  then the partial derivative $\partial
  \triv{\Op{D}}(V,W,k)/\partial \yy{p}$ is equal to
  $$\frac{\partial \triv{\Op{D}}(V,W,k)}{\partial \yy{p}} =
  \pi\begin{pmatrix}
  0 & \left(\yy{\kappa}_2+\sqrt{-1}\yy{\kappa}_1\right)\unity\\
  \left(\yy{\kappa}_2-\sqrt{-1}\yy{\kappa}_1\right)\unity
  \end{pmatrix}.$$ We conclude that
  $$\frac{\partial \lambda}{\partial \yy{p}}
  \langle\langle \phi([k],\lambda),
  \psi([k],\lambda)\rangle\rangle =
  \frac{\pi}{\xx{\gamma}_1\yy{\gamma}_2-\xx{\gamma}_2\yy{\gamma}_1}
  \langle\langle \phi([k],\lambda),
  \psi([k],\lambda)\rangle\rangle_{\xx{\gamma}},$$
  because $\xx{\kappa}$ and $\yy{\kappa}$ form the dual basis of
  $\xx{\gamma}$ and $\yy{\gamma}$, and therefore $\yy{\kappa}$
  is equal to
  $\frac{1}{\xx{\gamma}_1\yy{\gamma}_2-
  \xx{\gamma}_2\yy{\gamma}_1}(-\xx{\gamma}_2,\xx{\gamma}_1)$.
  Since the zeroes and poles of the numerator of
  $\Op{P}_{\xx{\gamma}}$ and $\Op{P}$
  coincide, condition (iv) of Lemma~\ref{projection 1} implies that
  $\Op{P}_{\xx{\gamma}} d\xx{p}\wedge d\lambda =
  (\partial \lambda/\partial 
  \yy{p})\Op{P}_{\xx{\gamma}}d\xx{p}\wedge d\yy{p}$
  is a regular form on the complex Bloch variety with values
  in the projections of rank one on
  $\banach{2}(\Delta)\times\banach{2}(\Delta)$. This shows 
  that the local sum of $\Op{P}_{\xx{\gamma}}$ over all those sheets
  of $\bloch(V,W)/\lattice\dual$
  considered as a covering space over
  $(\xx{p},\lambda)\in\mathbb{C}^2$,
  which contain $([k],\lambda)$, is holomorphic. Finally,
  Lemma~\ref{branchpoints} implies that the values of this sum are
  projections, whose rank is equal to the number of sheets. Due to
  Theorem~\ref{asymptotic analysis 1} and Corollary~\ref{regular fermi}
  the restrictions of the sheets of the complex Bloch variety
  considered as a covering space over
  $(\xx{p},\lambda)\in\mathbb{C}^2$ to the
  \Em{complex Fermi curve} are locally different. Hence the fourth
  statement (iv) is a consequence of (iii) and
  Lemma~\ref{dualizing sheaf}. Condition (v) is a direct consequence
  of (iii). As for the proof of (vi), let $\chi$ be an eigenfunction
  of the Dirac operator $\Op{D}(V,W,k)$ with eigenvalue $\lambda$.
  The exterior derivative of the form
  $\phi_2([k'],\lambda')\chi_1 dz +
  \phi_1([k'],\lambda')\chi_2 d\Bar{z}$
  is equal to
  $$d\left(\phi_2([k'],\lambda')\chi_1 dz +
  \phi_1([k'],\lambda')\chi_2 d\Bar{z}\right) = 
  (\lambda'-\lambda)
  \left(\phi_1([k'],\lambda')\chi_1+\phi_2([k'],\lambda')\chi_2\right)
  d\Bar{z} \wedge dz.$$ 
  If furthermore $g(\xx{\gamma},k)=g(\xx{\gamma},k')\mod \mathbb{Z}$,
  we may apply Lemma~\ref{branchpoints} and conclude
  \begin{multline*}
  \frac{\exp\left(2\pi\sqrt{-1}g(\yy{\gamma},k'-k)\right)-1}
  {\xx{\gamma}_1\yy{\gamma}_2-\xx{\gamma}_2\yy{\gamma}_1}
  \langle\langle\phi([k'],\lambda'),
  \chi\rangle\rangle_{\xx{\gamma}} =\\
  \left(\exp\left(2\pi\sqrt{-1}g(\yy{\gamma},k'-k)\right)-1\right)
  \int\limits_{x}^{x+\xx{\gamma}}
  \phi_2([k'],\lambda')\chi_1 dz +
  \phi_1([k'],\lambda')\chi_2 d\Bar{z} =\\
  2\sqrt{-1}(\lambda'-\lambda)
  \langle\langle\phi([k'],\lambda'),
  \chi\rangle\rangle.
  \end{multline*}
  Hence, due to the relation between the denominators of the
  projections $\Op{P}$ and $\Op{P}_{\xx{\gamma}}$
  established in the proof of (iii),
  the restrictions of the following two forms to the subvariety
  $\xx{p}'=\xx{p}=g(\xx{\gamma},k)\mod \mathbb{Z}$ of
  $\bloch(V,W)/\lattice\dual$ coincide:
  $$\frac{2\pi\sqrt{-1}d\yy{p}'}
     {\exp\left(2\pi\sqrt{-1}(\yy{p}'-\yy{p})\right)-1}
     \Op{P}([k'],\lambda')\chi =
  \frac{d\lambda'}{\lambda'-\lambda}
  \Op{P}_{\xx{\gamma}}([k'],\lambda')\chi.$$
The subvariety $\left\{(\xx{p}',\lambda')\in\mathbb{C}^2 \mid
(\xx{p}'\xx{\kappa}+\yy{p}\yy{\kappa},\lambda')
\in\bloch(V,W) \right\}$  of the complex Bloch variety may be
considered locally either as
a covering space over $\xx{p}'\in\mathbb{C}$ or as a covering
space over $\lambda'\in\mathbb{C}$.
The residue of the left hand side is the value of the local sum of
$\Op{P}\chi$ over all sheets, which contains $(\xx{p},\lambda)$, of the
subvariety considered as a covering over $\xx{p}'\in\mathbb{C}$,
and the residue of the right hand side is the value of the local sum of
$\Op{P}_{\xx{\gamma}}\chi$ over all sheets, which contains
$(\xx{p},\lambda)$, of the subvariety considered as a covering over
$\lambda'\in\mathbb{C}$. Obviously the restrictions of the
sheets of $\bloch(V,W)$ considered as a covering space over
$k\in\mathbb{C}^2$ to the subvariety $g(\xx{\gamma},k)=$constant are
locally different and also the restrictions of the sheets of
$\bloch(V,W)$ considered as a covering space over
$(\xx{p},\lambda)\in\mathbb{C}^2$ to the subvariety
$\xx{p}=$constant are locally different. Hence the residue of the
left hand side is equal to $\Breve{\Op{P}}([k],\lambda)\chi$ and the
residue of the right hand side is equal to
$\Breve{\Op{P}}_{\xx{\gamma}}([k],\lambda)\chi$. This proves (vi).
\end{proof}

These Lemmata have some remarkable consequences. If we fix the
component $\xx{p}=\xx{p}'$ of the wave vectors $k$ and the energy
$\lambda=\lambda'$, then all eigenfunctions $\chi$ of the
Dirac operator with eigenvalue $\lambda=\lambda'$ and boundary
condition $\chi(\xx{q},\yy{q})=
\exp\left(2\pi\sqrt{-1}\xx{q}\xx{p}'\right)\chi(\xx{q}+1,\yy{q})$
are uniquely determined by their restriction to some line
$\yy{q}=\yy{q}'$. In fact, if the component $\xx{p}$ of $k$ is
equal to $\xx{p}'$, then due to Lemma~\ref{branchpoints}
the action of the projection $\Breve{\Op{P}}_{\xx{\gamma}}([k])$
defined in (v) of Lemma~\ref{projection 2} on
$\chi$ depends only on the restriction to $\yy{q}=\yy{q}'$.

For reasons of simplicity we restrict ourself in the following
discussion to the \Em{complex Fermi curve},
which is our main interest. The restrictions of
the functions $\Op{P}$ and $\Op{P}_{\xx{\gamma}}$ to the
\Em{complex Fermi curve} are also denoted by $\Op{P}$ and
$\Op{P}_{\xx{\gamma}}$, respectively. If we identify the Hilbert space
$\banach{2}(\mathbb{R}/\mathbb{Z})\times\banach{2}(\mathbb{R}/\mathbb{Z})$
with all $\banach{2}$ functions on some line
$\xx{q}\in\mathbb{R}/\mathbb{Z},\yy{q}=\yy{q}'$,
then the composition of the operator
$$\begin{pmatrix}
\chi_1\\
\chi_2
\end{pmatrix}\longmapsto
\frac{\int\limits_{\xx{q}\in\mathbb{R}/\mathbb{Z},\yy{q}=\yy{q}'}
\phi_2([k])\chi_1 dz + \phi_1([k])\chi_2 d\Bar{z}}
{\int\limits_{\xx{q}\in\mathbb{R}/\mathbb{Z},\yy{q}=\yy{q}'}
\phi_2([k])\psi_1([k]) dz + \phi_1([k])\psi_2([k]) d\Bar{z}}
\begin{pmatrix}
\psi_1([k])\\
\psi_2([k])
\end{pmatrix}$$
from $\banach{2}(\mathbb{R}/\mathbb{Z})\times\banach{2}(\mathbb{R}/\mathbb{Z})$
into $\banach{2}(\Delta)\times\banach{2}(\Delta)$
and the operator, which restricts the functions to the line
$\xx{q}\in\mathbb{R}/\mathbb{Z},\yy{q}=\yy{q}'+1$
defines an operator on the Hilbert space
$\banach{2}(\mathbb{R}/\mathbb{Z})\times\banach{2}(\mathbb{R}/\mathbb{Z})$.
This restriction is not
defined for all $\banach{2}$--functions. But for all eigenfunctions of
the Dirac operators it is. In fact, due to Remark~\ref{sobolev embedding}
the eigenfunctions belong to the Sobolev spaces
$\sobolev{1,p}(\Delta)\times\sobolev{1,p}(\Delta)$ with
$1<p<2$. Due to the Sobolev embedding theorem \cite[5.4 Theorem]{Ad}
this implies that the restrictions of the eigenfunctions to a
one--dimensional plane in $\torus$ are $\banach{q}$--functions
for all $q<\infty$. The result is some finite rank operator--valued
function on the \Em{complex Fermi curve}, whose eigenvalues are equal
to
$\exp\left(2\pi\sqrt{-1}g(\yy{\gamma},k)\right)=
\exp\left(2\pi\sqrt{-1}\yy{p}\right)$.
Moreover, due to Lemma~\ref{projection 2} the local sum of this
function over all sheets of $\fermi(V,W)$ considered as a 
covering space over $\xx{p}\in\mathbb{C}$, which contain an element
$k'\in \fermi(V,W)$ (compare with Remark~\ref{local sum}),
is locally a holomorphic function depending on  $\xx{p}$.
It can be proven that the sum of these operators over all sheets of
$\fermi(V,W)$ considered as a covering space over
$\xx{p}\in\mathbb{C}$ defines a holomorphic function from $\mathbb{C}$
into the unbounded closed operators on
$\banach{2}(\mathbb{R}/\mathbb{Z})\times\banach{2}(\mathbb{R}/\mathbb{Z})$,
which is periodic with period $1$. The \Em{complex Fermi curve}
is the complex Bloch variety of this entire function.
Due to Lemma~\ref{branchpoints}, the restriction of the projection
$\Op{P}_{\xx{\gamma}}$ to the Hilbert space
$\banach{2}(\mathbb{R}/\mathbb{Z})\times\banach{2}(\mathbb{R}/\mathbb{Z})$
is equal to the composition of the operator
$$\begin{pmatrix}
\chi_1\\
\chi_2
\end{pmatrix}\longmapsto
\frac{\int\limits_{\xx{q}\in\mathbb{R}/\mathbb{Z},\yy{q}=\yy{q}'}
\phi_2([k])\chi_1 dz + \phi_1([k])\chi_2 d\Bar{z}}
{\int\limits_{\xx{q}\in\mathbb{R}/\mathbb{Z},\yy{q}=\yy{q}'}
\phi_2([k])\psi_1([k]) dz + \phi_1([k])\psi_2([k]) d\Bar{z}}
\begin{pmatrix}
\psi_1([k])\\
\psi_2([k])
\end{pmatrix}$$
from $\banach{2}(\mathbb{R}/\mathbb{Z})\times\banach{2}(\mathbb{R}/\mathbb{Z})$
into $\banach{2}(\Delta)\times\banach{2}(\Delta)$
and the operator, which restricts the functions to the line
$\xx{q}\in\mathbb{R}/\mathbb{Z},\yy{q}=\yy{q}'$.
These projections may be considered as the spectral projection of the
operator--valued entire function depending on $\xx{p}$ like
$\Op{P}$ is the spectral projection of the Dirac operator.
In particular, the sum of $\Op{P}_{\xx{\gamma}}$
considered as an projection on
$\banach{2}(\mathbb{R}/\mathbb{Z})\times\banach{2}(\mathbb{R}/\mathbb{Z})$
over all sheets of $\fermi(V,W)$ considered as a covering space
over $\xx{p}\in\mathbb{C}$ converges in the strong operator topology
to the identity.

All operator--valued holomorphic functions
depending on $\xx{p}\in\mathbb{C}$,
which commute pointwise with this unbounded--operator--valued entire
function, may be diagonalized simultaneously with this entire
function. The eigenvalues of such
commuting operator--valued functions defines the structure sheaf of a
one--sheeted covering over the \Em{complex Fermi curve}:

\newtheorem{Structure sheaf}[Lemma]{Structure sheaf}
\index{sheaf!structure $\sim$}
\index{structure!sheaf}
\begin{Structure sheaf}\label{structure sheaf}
Due to \cite[Chapter~8. \S3.]{GrRe} there exists an up to isomorphism
unique one--sheeted covering of the \Em{complex Fermi curve},
whose direct image of the \De{Structure sheaf} yields the following
coherent subsheaf of the normalization sheaf
(which is the direct image of the structure sheaf
of the normalization under the normalization map):
A germ $f$ of the normalization sheaf belongs to the
direct image of the \De{Structure sheaf},
if and only if the local sum of $f\Op{P}_{\xx{\gamma}}$
over all local sheets of the covering map $\xx{p}$
of the \Em{complex Fermi curve}, which contain the base point of the
corresponding stalk (compare with Remark~\ref{local sum}),
is a germ of the sheaf of operator--valued holomorphic functions
depending on $\xx{p}\in\mathbb{C}$.
\end{Structure sheaf}

Obviously the normalization of this one sheeted covering is isomorphic
to the normalization of the \Em{complex Fermi curve}.
To sum up the complex space with
\De{Structure sheaf}~\ref{structure sheaf} is a one--sheeted covering
over  the \Em{complex Fermi curve}
considered as a subvariety of $\mathbb{C}^2/\lattice\dual$,
and the normalization is a one--sheeted covering
over the former space. Therefore, we distinguish between the
\begin{description}
\item[geometric genus] of the normalization, the
\index{genus!$\rightarrow$ geometric $\sim$}
\index{geometric genus}
\item[arithmetic genus of the \De{Structure sheaf}~\ref{structure
    sheaf}], and, finally, the
\index{arithmetic genus!of the $\rightarrow$ structure sheaf}
\index{structure!sheaf!arithmetic genus of the $\sim$}
\item[arithmetic genus] of the \Em{complex Fermi curve}.
\index{genus!$\rightarrow$ arithmetic $\sim$}
\index{arithmetic genus}
\end{description}
We shall see in Proposition~\ref{asymptotic analysis 2} that the last
genus is always infinite and in some sense locally preserved.
Moreover, in Lemma~\ref{singularity} we shall see
that for pairs of potentials of the form $(U,\Bar{U})$ the difference of the
\Em{arithmetic genus} of the \De{Structure sheaf}~\ref{structure sheaf}
minus the \Em{geometric genus} is finite.
Potentials, whose \De{Structure sheaf}~\ref{structure sheaf}
have finite arithmetic genus, are called
\Em{finite type potentials}
(compare with Section~\ref{subsection complex Fermi curves of finite genus}).

It can be shown with similar arguments to those  used in the proof of
Lemma~\ref{projection 2}~(vi), that this \De{Structure sheaf}
does not depend on the choice of $\xx{\gamma}$ and $\yy{\gamma}$.
This one--sheeted covering would be a more appropriate definition
of the \Em{complex Fermi curve}. For example, the degree of the dual
eigen bundle is equal to the arithmetic genus of this curve plus one
\cite{Sch}.

So we have investigated the structure of the \Em{complex Fermi curve}
as a covering space over some coordinate  $\xx{p}=g(\xx{\gamma},k)$.
Let us generalize this construction to the more general case,
when the period $\xx{\gamma}$ is replaced by some arbitrary element
$\mu\in\mathbb{C}^2$.
The starting point is again some generalization of the equivalence of
the \De{Periodicity condition}~\ref{periodicity condition}.
For this purpose, however, we use
the \De{Trivialization}~\ref{trivialization}
and the eigenfunctions are described by
their periodic parts in the corresponding Hilbert bundle.

\begin{Lemma} \label{general branchpoints}
Let $\psi=\left(\begin{smallmatrix}
\psi_1(x)\\
\psi_2(x)
\end{smallmatrix}\right)$ be an eigenfunction corresponding to an
element $(k,\lambda)$ of $\bloch(V,W)$ and
$\phi=\left(\begin{smallmatrix}
\phi_1(x)\\
\phi_2(x)
\end{smallmatrix}\right)$ an eigenfunction of the transposed
Dirac operator corresponding to an element $(k',\lambda)\in
\bloch(V,W)$, and let
$\triv{\psi}(x)=\exp\left(-2\pi\sqrt{-1}g(x,k)\right)\psi(x)$ and
$\triv{\phi}(x)=\exp\left(2\pi\sqrt{-1}g(x,k')\right)\phi(x)$ the
corresponding periodic parts. Then the exterior derivative of the form
$\triv{\phi}_2\triv{\psi}_1 dz + \triv{\phi}_1\triv{\psi}_2
d\Bar{z}$ is equal to
$$\pi\left(\left(k_2-k_2'-\sqrt{-1}(k_1-k_1')\right)
\triv{\phi}_2\triv{\psi}_1 +
\left(k_2-k_2'+\sqrt{-1}(k_1-k_1')\right)
\triv{\phi}_1\triv{\psi}_2\right) d\Bar{z}\wedge dz.$$
Therefore, the bilinear form
$\langle\langle\triv{\phi},
\triv{\psi}\rangle\rangle_{\mu}$ vanishes, whenever
$g(\mu,k-k')=0$ but $k-k'\neq 0$.
\end{Lemma}

\begin{proof} Since $\psi$ is an eigenfunction of the Dirac operator
  and $\phi$ is a eigenfunction of the transposed Dirac operator with
  eigenvalues zero, we have
\begin{align*}
V\triv{\psi}_1 + \left(\pi(k_2+\sqrt{-1}k_1)+\partial\right)
\triv{\psi}_2 &=\lambda \triv{\psi}_1 & 
V\triv{\phi}_1 + \left(\pi(k_2'-\sqrt{-1}k_1')+\Bar{\partial}\right)
\triv{\phi}_2 &=\lambda \triv{\phi}_1\\
\left(\pi(k_2-\sqrt{-1}k_1)-\Bar{\partial}\right)\triv{\psi}_1 +
W\triv{\psi}_2 &=\lambda \triv{\psi}_2 &
\left(\pi(k_2'+\sqrt{-1}k_1')-\partial\right)\triv{\phi}_1 +
W\triv{\phi}_2 &=\lambda \triv{\phi}_2.
\end{align*}
Thus the exterior derivative of 
$d\left(\triv{\phi}_2\triv{\psi}_1 dz + 
\triv{\phi}_1\triv{\psi}_2 d\Bar{z} \right)$ is equal to
$$\pi\left(\left(k_2-k_2'-\sqrt{-1}(k_1-k_1')\right)
\triv{\phi}_2\triv{\psi}_1 +
\left(k_2-k_2'+\sqrt{-1}(k_1-k_1')\right)
\triv{\phi}_1\triv{\psi}_2\right) d\Bar{z}\wedge dz.$$
Since both forms are periodic, the integral of the right hand side
over the torus $\torus$ vanishes.
This integral is equal to
$2\pi\sqrt{-1}\langle\langle\triv{\phi},
\triv{\psi}\rangle\rangle_{\mu}$ if $\mu$ is equal to
$(k'_2-k_2,k_1-k'_1)$, which implies that $g(\mu,k-k')=0$.
\end{proof}

This lemma implies that whenever for some element
$\mu\in\mathbb{C}^2$ the  components $g(\mu,k)$ and $g(\mu,k')$
coincide, the bilinear form $\langle\langle\triv{\phi},
\triv{\psi}\rangle\rangle_{\mu}$ vanishes. But for this
generalization we have to pay a price. Since instead of the
eigenfunctions we used the periodic part of the eigenfunctions,
we lose invariance under the action of the dual lattice. Now we can also
define for general $\mu\in\mathbb{C}^2$ a projection valued function
$\triv{\Op{P}}_{\mu}$,
but this function should be considered to act not on
the Hilbert space, but on the Hilbert bundle introduced in the
\De{Trivialization}~\ref{trivialization}.
Let $\triv{\Op{P}}_{\mu}$ be the meromorphic map from
$\bloch(V,W)/\lattice\dual$ into the finite rank operators on the
Hilbert space
$\banach{2}(\torus)\times\banach{2}(\torus)$, which for
all elements $([k],\lambda)$ of the
complex Bloch variety is defined by
\index{spectral!projection!$\triv{\Op{P}}_{\mu}$}
$$\triv{\Op{P}}_{\mu}([k],\lambda):
\begin{pmatrix}
\triv{\chi}_1\\
\triv{\chi}_2
\end{pmatrix}\longmapsto
\frac{\langle\langle\phi([k],\lambda),
\psi_{k}\triv{\chi}\rangle\rangle_{\mu}}
{\langle\langle\phi([k],\lambda),
\psi([k],\lambda)\rangle\rangle_{\mu}}
\begin{pmatrix}
\psi_{-k}\psi_1([k],\lambda)\\
\psi_{-k}\psi_2([k],\lambda)
\end{pmatrix}.$$
Here $\psi_{k}$ is the operator--valued holomorphic function on the
complex Bloch variety, which is given by the multiplication with
the function $\psi_{k}$ defined in
Section~\ref{subsection resolvent}.
Lemma~\ref{projection 2} carries over to this function,
in particular $\triv{\Op{P}}_{\mu} dg(\mu,k)\wedge d\lambda$ is a
regular $2$--form on the complex Bloch variety and all
eigenfunctions $\chi$, which correspond to some element
$(k,\lambda)$ of the complex Bloch variety belong to the
range of the projection
$\psi_{k}\comp\Breve{\triv{\Op{P}}}_{\mu}(k,\lambda)$, which is
analogous to the projection
$\Breve{\Op{P}}_{\xx{\gamma}}(k,\lambda)$ defined in (v) of
Lemma~\ref{projection 2}. It might happen, however, that the
restriction of $\triv{\Op{P}}_{\mu}$ to the \Em{complex Fermi curve}
is not well defined, because the denominator might vanish identically
on a connected component of the \Em{complex Fermi curve}.
In these cases, the second part of the statement~(iv) of
Lemma~\ref{projection 2} cannot be generalized to this projection.
But Theorem~\ref{asymptotic analysis 1} and
Corollary~\ref{regular fermi} imply that this can only happen,
if $\mu$ is proportional to $(1,\pm\sqrt{-1})$
and if the \Em{complex Fermi curve} is the same as the
\Em{complex Fermi curve} of the zero potential.

The restriction of this projection to the
\Em{complex Fermi curve} is again the projection of
some operator--valued holomorphic function, whose
complex Bloch variety is equal to the
\Em{complex Fermi curve} $\fermi(V,W)$.
Let $\xx{\mu}$ and $\yy{\mu}$ be an arbitrary bases of
$\mathbb{C}^2$ and $\xx{\nu}$ and $\yy{\nu}$ the dual basis:
$$g(\xx{\mu},\xx{\nu})=1=g(\yy{\mu},\yy{\nu}) \text{ and }
g(\xx{\mu},\yy{\nu})=0=g(\yy{\mu},\xx{\nu}).$$
In analogy to the notation introduced in the context of the
\De{Fundamental domain}~\ref{fundamental domain}
we parameterize the momentum space by
$\xx{p}=g(\xx{\mu},k)$ and $\yy{p}=g(\yy{\mu},k)$,
which is equivalent to $k=\xx{p}\xx{\nu}+\yy{p}\yy{\nu}$.
If we assume that $g(\yy{\nu},\yy{\nu})$ is not equal to zero,
then the matrix
$\left(\begin{smallmatrix}
0 & \yy{\nu}_2+\sqrt{-1}\yy{\nu}_1\\
\yy{\nu}_2-\sqrt{-1}\yy{\nu}_1 & 0
\end{smallmatrix}\right)$ is invertible, and the inverse equals
$\left(\begin{smallmatrix}
0 & 1/(\yy{\nu}_2-\sqrt{-1}\yy{\nu}_1)\\
1/(\yy{\nu}_2+\sqrt{-1}\yy{\nu}_1) & 0
\end{smallmatrix}\right)$. The function
$\psi=\psi_{\xx{p}\xx{\nu}+\yy{p}\yy{\nu}}\triv{\psi}$
is an eigenfunction of the Dirac operator
$\Op{D}(V,W,\xx{p}\xx{\nu}+
\yy{p}\yy{\nu})$ with eigenvalue zero if and only if
the function $\triv{\psi}$ is an
eigenfunction of the operator
\begin{eqnarray*}
\triv{\Op{D}}_{\xx{\mu},\yy{\mu}}(V,W,\xx{p})&=&\left(\begin{smallmatrix}
0 & \yy{\nu}_2+\sqrt{-1}\yy{\nu}_1\\
\yy{\nu}_2-\sqrt{-1}\yy{\nu}_1 & 0
\end{smallmatrix}\right)^{-1}
\left(\Op{D}(V,W,0)+\xx{p}\pi\left(\begin{smallmatrix}
0 & \xx{\nu}_2+\sqrt{-1}\xx{\nu}_1\\
\xx{\nu}_2-\sqrt{-1}\xx{\nu}_1 & 0
\end{smallmatrix}\right)\right)\\
&=&\begin{pmatrix}
\frac{\xx{p}\pi(\xx{\nu}_2-\sqrt{-1}\xx{\nu}_1)-\Bar{\partial}_{0}}
{\yy{\nu}_2-\sqrt{-1}\yy{\nu}_1} & 
\frac{W}{\yy{\nu}_2-\sqrt{-1}\yy{\nu}_1}\\
\frac{V}{\yy{\nu}_2+\sqrt{-1}\yy{\nu}_1} &
\frac{\xx{p}\pi(\xx{\nu}_2+\sqrt{-1}\xx{\nu}_1)+\partial_{0}}
{\yy{\nu}_2+\sqrt{-1}\yy{\nu}_1} 
\end{pmatrix}
\end{eqnarray*}
with eigenvalue $\yy{p}\pi$. Thus the \Em{complex Fermi curve}
$\fermi(V,W)$ is equal to the set of all
$$\left\{k=\xx{p}\xx{\nu}+\yy{p}\yy{\nu}\in\mathbb{C}^2\mid
\yy{p}\pi \text{ is an eigenvalue of }
\triv{\Op{D}}_{\xx{\mu},\yy{\mu}}(V,W,\xx{p})\right\},$$
and due to Lemma~\ref{general branchpoints} for all
$k=\xx{p}\xx{\nu}+\yy{p}\yy{\nu}\in \fermi(V,W)$ 
the projection $\Breve{\triv{\Op{P}}}_{\xx{\mu}}(k)$ is the
spectral projection of
$\triv{\Op{D}}_{\xx{\mu},\yy{\mu}}(V,W,\xx{p})$
onto the generalized eigenspace associated with eigenvalue
$\yy{p}\pi$, which is equal to the residue of
the resolvent of $\triv{\Op{D}}_{\xx{\mu},\yy{\mu}}(V,W,\xx{p})$
at the eigenvalue $\yy{p}\pi$
(see e.\ g.\ \cite[Appendix to XII.1]{RS4}).

\subsection{\Em{Complex Fermi curves} of finite genus}
\label{subsection complex Fermi curves of finite genus}

If there exists some finite subset of $\lattice\dual$, such that for
all $\kappa$ in the complement of this finite subset the
normalizations of the handle  with index $\kappa$ described in
Theorem~\ref{asymptotic analysis 1} decomposes
into two connected components,
one of which is a disc excluded from $\Set{V}^+_{\varepsilon,\delta}$
around the point $k_{\kappa}^+$ and the other is a disc excluded
from $\Set{V}^-_{\varepsilon,\delta}$ around the point $k_{\kappa}^-$,
then the normalization of $\fermi(V,W)/\lattice\dual$ is
biholomorphic to a compact Riemann surface with two points at infinity
removed. In fact, in this case the normalization of
$\fermi(V,W)/\lattice\dual$
contains two open subsets, which are biholomorphic to the open subsets
\begin{equation*}
\begin{split}
\left\{ k\in \mathbb{C}^2\;\left|\; k_1-\sqrt{-1}k_2=0
\text{ and } \|k\|>1/\delta\right.\right\} &\text{ and }\\
\left\{ k\in \mathbb{C}^2\;\left|\; k_1+\sqrt{-1}k_2=0
\text{ and } \|k\|>1/\delta\right. \right\}
\end{split}
\end{equation*}
of the normalization of $\fermi(0,0)$, with some
$\varepsilon>0$. Due to Theorem~\ref{asymptotic analysis 1}
the complement of these two open subsets of the normalization of
$\fermi(V,W)/\lattice\dual$ is compact. Hence we may compactify
the normalization of $\fermi(U,U)/\lattice\dual$ by adding to
the first subset the point
$\infty^-=\lim\limits_{t\rightarrow\infty}(t,-\sqrt{-1}t)$ and to
the second subset the point
$\infty^+=\lim\limits_{t\rightarrow\infty}(t,\sqrt{-1}t)$.

If we normalize the eigenfunctions $\psi$ of the Dirac operator by
some linear condition, then the eigenfunction is
a meromorphic function on $\fermi(V,W)/\lattice\dual$. The
following normalization of the eigenfunction
\index{normalization!of the eigenfunction}
\index{eigenfunction ($\rightarrow$ Baker--Akhiezer function)!normalization
  of the $\sim$}
$\psi(x)=\left(\begin{smallmatrix}
\psi_1(x)\\
\psi_2(x)
\end{smallmatrix}\right)$ and the eigenfunction of the transposed
Dirac operator 
$\phi(x)=\left(\begin{smallmatrix}
\phi_1(x)\\
\phi_2(x)
\end{smallmatrix}\right)$ is convenient:
$$\psi_1(x=0)+\psi_2(x=0)=\sqrt{2} \text{ and }
\phi_1(x=0)-\phi_2(x=0)=\sqrt{2},$$
respectively. If $h_0$ denotes the matrix
$h_0=1/\sqrt{2}\left(\begin{smallmatrix}
1 & 1\\
1 & -1
\end{smallmatrix}\right)$, then it may be written in the form
$\left(h_0\psi(x=0)\right)_1=1$ and $\left(\Op{J}h_0\phi(x=0)\right)_1=1$.
Although these eigenfunctions are always meromorphic
functions on the normalization of $\fermi(V,W)/\lattice\dual$,
it might happen that this normalization is biholomorphic to a
compact Riemann surface with two points removed,
but the eigenfunctions have an infinite number of singularities
on the normalization of $\fermi(V,W)/\lattice\dual$ \cite{Sch}.
In this case they are not meromorphic functions on the compact Riemann
surface with two essential singularities at the two points at
infinity. Consequently we call a pair of potentials $(V,W)$
to be of finite type, if they obey the following conditions:
\begin{description}\index{finite!type potentials|(}
\item[Finite type potentials (i)]\index{normalization!of finite genus}
\index{finite!geometric genus}
\index{geometric genus!finite $\sim$}
The normalization of
$\fermi(V,W)/\lattice\dual$ can be compactified to a compact
Riemann surface  by adding the two points
$\infty^-$ and $\infty^+$ at infinity.
\item[Finite type potentials (ii)] The pullback of the
normalized eigenfunctions under the normalization map have only a
finite number of poles on the normalization of
$\fermi(V,W)/\lattice\dual$.
\index{finite!type potentials|)}
\end{description}
The entries of the normalized eigenfunctions generate a coherent
subsheaf of the sheaf of meromorphic functions.
Due to an improved version of statement (iv) in Lemma~\ref{projection 2}
(compare with \cite[Chapter~9]{Sch})
the number of poles of the normalized eigenfunction is equal to the
arithmetic genus of the \Em{Structure sheaf}~\ref{structure sheaf}
plus one (i.\ e.\ the arithmetic genus of the modified
\Em{Structure sheaf}~\ref{structure sheaf},
where the two infinities are identified to an ordinary double point).
Therefore, condition \Em{Finite type potentials}~(ii) is equivalent to
the condition that the \Em{arithmetic genus} of the
\index{finite!arithmetic genus}
\Em{Structure sheaf}~\ref{structure sheaf} is finite.
Moreover, condition \Em{Finite type potentials}~(i) is equivalent
to the condition that the \Em{geometric genus} is finite.
Therefore, condition \Em{Finite type potentials}~(ii) implies
condition \Em{Finite type potentials}~(i).
We should remark that these eigenfunctions may be characterized by
some abstract properties and are called \Em{Baker--Akhiezer functions}
\index{Baker--Akhiezer function ($\rightarrow$ eigenfunction)}
(compare with \cite[Chapter~2 \S2.]{DKN} and conditions
\Em{Baker--Akhiezer function}~(i)--(ii)
in the proof of Lemma~\ref{compactification of isospectral sets}).
Obviously these functions uniquely
determine the corresponding potentials, which therefore have to be
analytic. Finally, all sheaves of this form may be characterized
(compare with \cite{Sch}).

\begin{Lemma} \label{singularity}
If a pair of potentials of the form $(U,\Bar{U})$ obeys
the condition \Em{Finite type potentials}~(i), then it satisfies also
the condition \Em{Finite type potentials}~(ii).
\end{Lemma}

\begin{proof} For small $\varepsilon>0$ the open subsets
  $\Set{V}^{\pm}_{\varepsilon,\delta}$ of
  $\fermi(0,0)$ are one--sheeted coverings over some open subset of
  $k_1\in\mathbb{C}$. Thus, due to the proof of
  Theorem~\ref{asymptotic analysis 1}, the same is true
  for arbitrary potentials $V$ and $W$.
  Furthermore, with the exception of finitely many
  $\kappa\in\lattice\dual$ all handles with index $\kappa$ are
  two--sheeted coverings over $\xx{p}\in\mathbb{C}$. Hence the only
  possible singularities of these handles are double points. If the
  normalization of $\fermi(U,\Bar{U})/\lattice\dual$ is a
  compact Riemann surface with two points removed, then with the
  exception of finitely many $\kappa\in\lattice\dual$ the
  normalization of the handles with index $\kappa$ have two connected
  components. Since the involution $\eta$ permutes the points
  $k_{\kappa}^{\pm}$, this involution also permutes the two connected
  components of the normalization of these handles. This implies that
  with the exception of finitely many $\kappa\in\lattice\dual$
  the only possible singularities of the handle with index $\kappa$
  are double points, which are invariant under $\eta$.
  Now Corollary~\ref{fixed points} and condition~(iii) of
  Lemma~\ref{projection 2} imply that the
  projection valued function $\Op{P}_{\xx{\gamma}}$
  has no poles on these handles. Thus due to the Maximum principle
  (see e.g. \cite[Chapter~4. Theorem~12.]{Ah}), the estimates of the
  eigenfunctions from Theorem~\ref{asymptotic analysis 1} on the open sets
  $\Set{V}^{\pm}_{\varepsilon,\delta}$ extend to these handles. This
  proves that the normalized eigenfunctions have only a finite number
  of poles on the normalization of
  $\fermi(U,\Bar{U})/\lattice\dual$.
\end{proof}

\subsubsection{The inverse problem}
\label{subsubscetion inverse problem}

As we have seen, in some cases we may associate to potentials
a compact Riemann surface, and moreover we found several properties of
theses Riemann surface. Let us now try to invert this procedure and to
associate to some Riemann surface with some specified properties
a pair of potentials, such that the compactification of the
normalization of the corresponding \Em{complex Fermi curve} is equal
to the given Riemann surface.
First we introduce the data of tuples
$(\Spa{Y},\infty^-,\infty^+,k)$.
Here $\Spa{Y}$ is a compact pure one--dimensional
complex space with two marked non--singular points $\infty^{\pm}$
with the property that each connected component
of the subspace of nonsingular points
contains at least one of these points,
and $k$ is a multi--valued meromorphic
function from $\Spa{Y}$ into $\mathbb{C}^2$,
which is holomorphic on $\Spa{Y}\setminus\{\infty^-,\infty^+\}$.
The complex space $\Spa{Y}\setminus\{\infty^-,\infty^+\}$ will be a
candidate for the one--sheeted covering of a \Em{complex Fermi curve},
whose  structure sheaf is defined in
\De{Structure sheaf}~\ref{structure sheaf}. In the sequel
it mostly can be chosen to be a compact smooth Riemann surface.
Furthermore, we assume the following conditions:
\begin{description}
\index{condition!quasi--momenta  (i)--(iv)|(}
\index{quasi--momenta!condition $\sim$  (i)--(iv)|(}
\item[Quasi--momenta (i)]
  The function $k$ has poles of first--order at the two
  marked points and for some branch of this function the
  $\mathbb{C}$--valued functions $k_1-\sqrt{-1}k_2$ and 
  $k_1+\sqrt{-1}k_2$ vanish at $\infty^-$ and $\infty^+$,
  respectively.
\item[Quasi--momenta (ii)]
  The difference of two arbitrary branches of the function
  $k$ is some element of $\lattice\dual$. Conversely, for all
  $\kappa\in\lattice\dual$ and all branches of $k$, $k+\kappa$ is
  some other branch of $k$.
\end{description}
Here we do not assume that we may pass
from one branch to all other branches by continuation along closed
cycles of $\Spa{Y}$ (compare with Remark~\ref{effective lattice}).

For all such data $(\Spa{Y},\infty^-,\infty^+,k)$ the image of
$\Spa{Y}\setminus\{\infty^-,\infty^+\}$ under the mapping $k$ is a pure
one--dimensional complex subvariety of $\mathbb{C}^2$, which is
invariant under translations by elements of $\lattice\dual$.
As shown in the following lemma these images are
\Em{complex Fermi curves} and consequently denoted by
$\fermi(\Spa{Y},\infty^-,\infty^+,k)$ or just $\fermi$.
The \Em{complex Fermi curves} we are interested in are endowed with
the anti--holomorphic involutions $\rho$ and $\eta$ and the
holomorphic involution $\sigma$. In the sequel we therefore impose
sometimes the additional conditions:
\begin{description}
\item[Quasi--momenta (iii)]
  The complex space $\Spa{Y}$ is endowed with an
  anti--holomorphic involution $\eta$ without fixed points. Moreover,
  any branch of the transformed function $\eta^{\ast} k$ is equal to
  some branch of the function $-\Bar{k}$.
\item[Quasi--momenta (iv)]
  The complex space $\Spa{Y}$ is endowed with an
  holomorphic involution $\sigma$. Moreover,
  any branch of the transformed function $\sigma^{\ast}k$ is equal to
  some branch of the function $-k$.
\index{condition!quasi--momenta  (i)--(iv)|)}
\index{quasi--momenta!condition $\sim$  (i)--(iv)|)}
\end{description}

\begin{Lemma}\label{existence of potentials}
\index{existence!of a potential}
Let $\Spa{Y}$ be a smooth connected compact smooth Riemann surface with two
marked points $\infty^-$ and $\infty^+$ and let $k$ be a multi--valued
meromorphic function, which is holomorphic on
$\Spa{Y}\setminus\{\infty^-,\infty^+\}$.
If $(\Spa{Y},\infty^-,\infty^+,k)$ satisfies the conditions
\Em{Quasi--momenta}~(i)--(iii),
then there exists an analytic complex potential $U$,
whose  \Em{complex Fermi curve} is equal to
$\fermi(\Spa{Y},\infty^-,\infty^+,k)$.
Moreover, the pullback of the eigenfunction under the normalization map
takes at all pairwise different elements of
$\Spa{Y}\setminus\{\infty^-,\infty^+\}$ linearly independent values.
Finally, for any choice of finitely many points $y_1,\ldots,y_l$ of
$\Spa{Y}\setminus\{\infty^-,\infty^+\}$ there exists one of these potentials,
such that the values of the pullback of the  corresponding
eigenfunction at these points $y_1,\ldots,y_l$ has no zeroes
considered as a function from $\Delta$ to $\mathbb{C}^2$.
\end{Lemma}

Due to \cite[Lemma~2.3]{GS1} the assumption of connectedness of $Y$
is necessary. In fact, if $\Spa{Y}$ has two connected
components, the forms $dk_1\pm\sqrt{-1}dk_2$ are holomorphic on the
component containing $\infty^{\pm}$, respectively. This implies that
the intersection form of those two cycles does not vanish, whose
intersection forms with any other cycle is equal to the integral of
$\xx{p}$ and $\yy{p}$ along this cycle.

\begin{proof}
We use the methods of inverse spectral theory of Lax operators
(see e.\ g.\ \cite{Kr1,DKN}). For this purpose it is convenient to
introduce the compact complex space $\Tilde{\Spa{Y}}$, which is obtained
form $\Spa{Y}$ by identifying the two marked points $\infty^{\pm}$ to an
\index{ordinary double point!at infinity $(\infty^-,\infty^+)$}
ordinary double point. Furthermore, we shall introduce the real
two--dimensional subgroup of the Picard group corresponding to the
translations by $\torus$. For all $x\in\mathbb{R}^2$,
let $\Tilde{\Sh{L}}(x)$ denote the
locally free sheaf of rank one on $\Tilde{\Spa{Y}}$,
which is defined by a co--cycle with respect to the covering
$\Tilde{\Spa{Y}}=\Tilde{\Spa{Y}}\setminus\{\infty^-,\infty^+\}\cup$
some small neighbourhood of the double point $(\infty^-,\infty^+)$.
On the intersection of the two open sets
the transition function of this co--cycle is equal to
$$\exp\left(2\pi\sqrt{-1}g(x,k)\right)=
\exp\left(2\pi\sqrt{-1}
\left(g(x,\xx{\kappa})\xx{p}+g(x,\yy{\kappa})\yy{p}\right)\right)=
\exp\left(2\pi\sqrt{-1}
\left(\xx{q}\xx{p}+\yy{q}\yy{p}\right)\right)$$
for any $k\in\mathbb{C}^2$.
If $x$ belongs to the lattice $\lattice$, then this transition
function extends to a global non--vanishing holomorphic function on
$\Tilde{\Spa{Y}}\setminus\{\infty^-,\infty^+\}$, and the corresponding
sheaf is equal to the structure sheaf of $\Tilde{\Spa{Y}}$.
Due to the general construction (see e.\ g.\ \cite{DKN,Sch}) for any
divisor $\Tilde{D}$ on $\Tilde{\Spa{Y}}$, which has the following properties,
there is associated a unique \Em{Baker--Akhiezer function}.
\begin{description}
\index{divisor!of the Baker--Akhiezer function|(}
\index{Baker--Akhiezer function ($\rightarrow$ eigenfunction)!divisor of the
  $\sim$|(}
\index{condition!divisor (i)--(iii)|(}
\index{divisor!condition $\sim$ (i)--(iii)|(}
\item[Divisor (i)]
  $\Tilde{D}$ is an integral divisor, whose support
  is contained in $\Tilde{\Spa{Y}}\setminus\{\infty^-,\infty^+\}$. Moreover,
  the degree of $\Tilde{D}$ is equal to the genus of $\Spa{Y}$ plus one
  (i.\ e.\ the arithmetic genus of $\Tilde{\Spa{Y}}$).
\item[Divisor (ii)]
  The difference $\eta(\Tilde{D})-\Tilde{D}$ is the
  principal divisor of a function $\func{f}$,
  which takes the values $\pm 1$ at the marked points $\infty^{\pm}$.
  In particular, this function obeys
  $\func{f}\eta^{\ast}\Bar{\func{f}}=-1$.
\item[Divisor (iii)]
  For all $x\in\mathbb{R}^2$ the sheaf
  $\Tilde{\Sh{L}}(x)\otimes\Sh{O}_{\Tilde{D}}$ on
  $\Tilde{\Spa{Y}}$ has no non--trivial global section, which vanishes at
  the double point $(\infty^-,\infty^+)$.
\index{divisor!of the Baker--Akhiezer function|)}
\index{Baker--Akhiezer function ($\rightarrow$ eigenfunction)!divisor of the
  $\sim$|)}
\index{condition!divisor (i)--(iii)|)}
\index{divisor!condition $\sim$ (i)--(iii)|)}
\end{description}
For given $\Tilde{D}$ let $D$ denote the corresponding divisors on
$\Spa{Y}$. With the help of the
Riemann--Roch theorem \cite[Theorem~16.9]{Fo} it is quite easy to see
that for any divisor $D$ on $\Spa{Y}$,
which has the following properties \Em{Divisor}~(i')--(iii'),
there exists a family $\Tilde{D}_{t}$ parameterized by $t\in S^1$
of divisors on $\Tilde{\Spa{Y}}$, fulfilling conditions
\Em{Divisor}~(i)--(ii), whose corresponding divisors $D_{t}$
are all isomorphic to $D$. The degree of freedom of this family
corresponds to the isospectral transformations described at the end of
Section~\ref{subsection reductions}.

\begin{description}
\item[Divisor (i')] The degree of $D$ is equal to the genus of $\Spa{Y}$ plus one.
\item[Divisor (ii')] The difference $\eta(D)-D$ is a principal divisor of a
  function $\func{f}$,
  which obeys $\func{f}\eta^{\ast}\Bar{\func{f}}<0$.
\item[Divisor (iii')] The sheaf $\Sh{O}_{D}$ has no global non--trivial
  section, which vanishes at the two marked points $\infty^-$ and
  $\infty^+$. 
\end{description}

\begin{Remark}\label{quaternionic line bundle}
\index{quaternionic!complex line bundle}
Complex line bundles, whose divisors
obey condition~(ii) are called quaternionic, since the space
of sections is a quaternionic space \cite[p.\ 667]{Hi}.
In order to distinguish this notion of quaternionic line bundles
from the notion of quaternionic line bundles occurring in
`quaternionic function theory', we shall emphasize that the first
notion describes a subclass of complex line bundles. Consequently we
shall speak of  quaternionic complex line bundles.
\end{Remark}

For all $x\in\mathbb{R}^2$ let $\Sh{L}(x)$ denotes the analogous
locally free sheaf of rank one on $\Spa{Y}$. Obviously these sheaves correspond
to real line bundles of degree zero. In particular, the tensor product
of these line bundles with any quaternionic complex line bundle gives
again quaternionic complex line bundles. Hence we have to show that
there exists a divisor $D$, whose degree is equal to the genus of $\Spa{Y}$
plus one, such that the sheaves $\Sh{L}(x)\otimes\Sh{O}_D$
have no non--trivial global section, which vanishes at the two marked
points $\infty^-$ and $\infty^+$. Obviously all these sheaves
correspond to quaternionic complex line bundles. Hence the dimension of
the space of global sections is divisible by two.
Due to \cite[Proposition~3.3]{GH} the real part of the component of the
Picard group of $\Spa{Y}$ corresponding to degree $g\pm 1$ has a component,
which contains only quaternionic complex line bundles.
Due to Marten's theorem \cite[Chapter~IV (5.1)~Theorem]{ACGH} the
subvariety of the Picard group of divisors of degree $g-1$, whose
space of sections is at least two--dimensional, has co--dimension of at
least three. In Marten's theorem it is assumed that the genus
is not smaller than three. If the genus is smaller than three, this
statement is trivial, because in these cases no line bundle, whose
degree is equal to the genus minus one, can have a two--dimensional
space of global sections. Therefore, the subvariety of the real
part of the Picard group of $\Spa{Y}$ corresponding to
quaternionic complex line bundles of degree $g-1$, which have
non--trivial global sections, has co--dimension larger than two. We
conclude that there exists a divisor $\Tilde{D}$ fulfilling conditions
\Em{Divisor}~(i)--(iii), and consequently a corresponding
\Em{Baker--Akhiezer function} and a corresponding
pair of potentials of the form $(U,\Bar{U})$ whose
\Em{complex Fermi curve} is equal to
$\fermi(\Spa{Y},\infty^-,\infty^+,k)$.

If the eigenfunctions corresponding to two different elements
$y\neq y'$  of $\Spa{Y}\setminus\{\infty^-,\infty^+\}$ are linearly dependent,
then the function $k$ has to take the same values $\mod\lattice\dual$
at these elements. Due to Lemma~\ref{general branchpoints} this
implies that the projection $\Op{P}_{\xx{\gamma}}d\xx{p}$ is not
an operator--valued holomorphic form on
$\Spa{Y}\setminus\{\infty^-,\infty^+\}$.
More precisely, the local sum over all sheets of $\Spa{Y}$
considered as a covering over $\xx{p}\in\mathbb{C}$,
which contain $y$ or $y'$ (compare with Remark~\ref{local sum}), respectively,
will then not be a holomorphic function with respect to $\xx{p}$.
On the other hand the corresponding
\De{Structure sheaf}~\ref{structure sheaf} is equal to the
normalization sheaves. Therefore,
due to an improved version of statement (iv) in Lemma~\ref{projection 2}
(compare with \cite[Chapter~9]{Sch}),
this projection valued form $\Op{P}_{\xx{\gamma}}d\xx{p}$
is holomorphic, and the values of the \De{Baker--Akhiezer function}
at different points are linearly independent.

Finally, we remark that due to the construction of the
\Em{Baker--Akhiezer function}, the eigenfunction corresponding to some
element $y\in \Spa{Y}\setminus\{\infty^-,\infty^+\}$ considered as a
function from $\Delta$ to $\mathbb{C}^2$ has a  zero at $x\in\Delta$,
if and only if
\begin{description}
\item[(i)] the element $y$ belongs to the support of $\Tilde{D}(x)$,
  and 
\item[(ii)] the unique global section $\func{f}$ of
  $\Sh{O}_{\Tilde{D}(x)}$ on $\Spa{Y}$, which takes the values $\pm 1$
  at the two marked points $\infty^{\pm}$, is also a section of
  $\Sh{O}_{\Tilde{D}(x)-y}$.
\end{description}
Here $\Tilde{D}(x)$ denotes the unique integral divisor of
$\Tilde{\Spa{Y}}$, whose sheaf $\Sh{O}_{\Tilde{D}(x)}$
is isomorphic to
$\Tilde{\Sh{L}}(x)\otimes\Sh{O}_{\Tilde{D}}$.
Due to condition \Em{Divisor}~(ii) this implies that $\eta(y)$
also belongs to the support of $\Tilde{D}(x)$,
and that $\func{f}$ is even a
global section of $\Sh{O}_{\Tilde{D}(x)-y-\eta(y)}$.
In particular, the corresponding divisor $D(x)-y-\eta(y)$
of $\Spa{Y}$ is an integral divisor of degree $g-1$,
whose space of global sections is at least two--dimensional.
Again  Marten's theorem \cite[Chapter~IV (5.1)~Theorem]{ACGH}
shows that for generic $y$ this is not the case
for all $x\in\mathbb{R}^2$. 
\end{proof}

\subsubsection{The B\"acklund transformations}
\label{subsubsection baecklund}

The one--sheeted coverings of the \Em{complex Fermi curve} of these
potentials, whose structure sheaf is defined in
\De{Structure sheaf}~\ref{structure sheaf}, are equal to
the normalization of $\Spa{Y}\setminus\{\infty^-,\infty^+\}$.
This construction can be generalized to all one--sheeted coverings
of $\fermi(\Spa{Y},\infty^-,\infty^+,k)/\lattice\dual$
of finite arithmetic genus,
on which the involution $\eta$ does not have any fixed points
(i.\ e.\ all double points of the form $(y,\eta(y))$ of the
\Em{complex Fermi curve} have to be removed). Moreover, the
\Em{B\"acklund transformation}
\index{B\"acklund transformation|(}
provides a tool for transforming a
potential corresponding to a one--sheeted covering of the
\Em{complex Fermi curve} into a potential corresponding to
another one--sheeted covering, such that the former is a
one--sheeted covering of the latter
(compare with \cite{EK},\cite[Chapter~6.]{MS} and \cite[Section~4]{LMcL}).
We shall now derive the formulas of these
\Em{B\"acklund transformations} in the present context.
First we remark that for all pairs $(y,y')$ of the
\Em{complex Fermi curve} the transformation
$D\mapsto D+\infty^-+\infty^+-y-y'$ preserves the degree of the
divisor and therefore should induce an isospectral transformation of
the corresponding potentials, in case they exist.
Moreover, due to the properties of the \Em{Baker--Akhiezer function},
the $n$--th partial derivatives of the eigenfunctions
yield cross sections of the sheaves
$\Sh{O}_{D(x)+n\infty^-n\infty^+}$.
Therefore, we may calculate the \Em{Baker--Akhiezer functions}
corresponding to the divisors $D(x)+\infty^-+\infty^+-y-y'$ in terms of
the \Em{Baker--Akhiezer function} corresponding to the divisors $D(x)$.
Since we are only interested in pairs of potentials of the form
$(U,\Bar{U})$ the pair $(y,y')$ has to be invariant under $\eta$. In fact,
the transformation $D\mapsto D+\infty^-+\infty^+-y-\eta(y)$ preserves
the divisors corresponding to quaternionic complex line bundles.
Let $\psi$ be the \Em{Baker--Akhiezer function} of the pair of
complex potentials $(U,\Bar{U})$ and $\chi$ and $\Op{J}\Bar{\chi}$
the values of $\psi$ at the pair $(y,\eta(y))$ of points
of the corresponding \Em{complex Fermi curve}.
In the following considerations $\chi$ may be multiplied with some
complex number without changing the transformation. The function
$$\psi'=\begin{pmatrix}
\partial+a & b\\
c & \Bar{\partial}+d
\end{pmatrix}\psi$$
with periodic functions $a,b,c,d$ on $\torus$ has the
correct asymptotic behaviour of the \Em{Baker--Akhiezer function}
corresponding to $D(x)+\infty^-+\infty^+-y-\eta(y)$.
If we impose the condition, that $\psi'$ vanishes at $y$ and
$\eta(y)$, then these functions are given by
\begin{align*}
a&=-\frac{(\partial \chi_1)\Bar{\chi}_1+\chi_2(\partial \Bar{\chi}_2)}
        {\chi_1\Bar{\chi}_1+\chi_2\Bar{\chi}_2} &
b&=\frac{\chi_1(\partial \Bar{\chi}_2)-(\partial \chi_1)\Bar{\chi}_2}
        {\chi_1\Bar{\chi}_1+\chi_2\Bar{\chi}_2}\\
c&=-\Bar{b}=
   \frac{\chi_2(\Bar{\partial}\Bar{\chi}_1)-
                   (\Bar{\partial}\chi_2)\Bar{\chi}_1}
        {\chi_1\Bar{\chi}_1+\chi_2\Bar{\chi}_2} &
d&=\Bar{a}=
   -\frac{\chi_1(\Bar{\partial}\Bar{\chi}_1)+
                   (\Bar{\partial}\chi_2)\Bar{\chi}_2}
        {\chi_1\Bar{\chi}_1+\chi_2\Bar{\chi}_2}.
\end{align*}
In order to verify that the corresponding
\Em{Baker--Akhiezer function} $\psi'$ is indeed an eigenfunction of
some Dirac operator we calculate the derivatives
\begin{align*}
\Bar{\partial} a&=U\Bar{U}-b\Bar{b} &
\Bar{\partial} b&=-\partial \Bar{U}-\Bar{U}a+b\Bar{a}\\
\partial c&=-\partial\Bar{b}=\Bar{\partial}U+U\Bar{a}-a\Bar{b} &
\partial d&=\partial\Bar{a}=U\Bar{U}-b\Bar{b}.
\end{align*}
Now a direct calculation shows that $\psi'$ is indeed an eigenfunction
corresponding to the pair of potentials $(U',\Bar{U}')=(\Bar{b},b)$,
if $\chi$ has no zero considered as a function from $\Delta$ to
$\mathbb{C}^2$. In this case the \Em{complex Fermi curve}
$\fermi(U',\Bar{U}')$ coincides with $\fermi(U,\Bar{U})$.
We should remark that the condition, that the value $\chi$ of the original
\Em{Baker--Akhiezer function} at $y$ (or $\eta(y)$) has no zero
considered as a function from $\Delta$ to $\mathbb{C}^2$, fits
perfectly with Lemma~\ref{existence of potentials}.
Obviously these formulas extend to the situation,
where $\chi$ is any element of the kernel of the Dirac operator
$\Op{D}(U,\Bar{U},k)$ with $k\in\fermi(U,\Bar{U})$,
which has no zero considered as a function from $\Delta$ to $\mathbb{C}^2$.
In case that $k$ is a singularity of the
\Em{complex Fermi curve} there might exist several values of the
pullback of the \Em{Baker--Akhiezer function} to the normalization and
more generally the eigenspace might have dimension larger than one.
However, only those elements of the range corresponding
spectral projection $\Breve{\Op{P}}([k])$
may be used in order to obtain a periodic potential,
which are proper eigenfunctions and have no zero considered
as a function from $\Delta$ to $\mathbb{C}^2$.

In general these formulas yield for two given
(not necessarily quasi--periodic) eigenfunctions $\chi$ and $\psi$ of
the Dirac operator
$\left(\begin{smallmatrix}
U & \partial\\
-\Bar{\partial} &\Bar{U}
\end{smallmatrix}\right)$
on some domain of $\mathbb{C}$ another eigenfunction $\psi'$
corresponding to another Dirac operator
$\left(\begin{smallmatrix}
U' & \partial\\
-\Bar{\partial} &\Bar{U}'
\end{smallmatrix}\right)$.
We would have obtained more general formulas, if we had started
with two elements of the complex Bloch variety $(y,y')$, which are
not invariant under $\eta$. In this context they are called
B\"acklund transformations \cite{EK}.

We shall derive also the formulas for the inverse transformation
corresponding to the transformation
$D\mapsto D+y+\eta(y)-\infty^--\infty^+$. For this purpose we remark
that the arguments of the proof of Lemma~\ref{projection 2}~(iv) shows
that the divisor $D$ of the \Em{Baker--Akhiezer function} and the
divisor $D^{t}$ of the transposed \Em{Baker--Akhiezer function} obey the
relation $\Sh{O}_{D}\otimes\Sh{O}_{D^{t}}\simeq
\Sh{O}_{K+2\infty^-+2\infty^+}$, where $\Sh{O}_{K}$ is the
dualizing sheaf of the one--sheeted covering described in
\De{Structure sheaf}~\ref{structure sheaf}
(for smooth one sheeted coverings $K$ is the canonical divisor).
More precisely, there exists a more restrictive version of
this relation, which makes use of the corresponding
one--sheeted coverings of the \Em{complex Fermi curve},
whose infinities are identified to an ordinary double point
(compare with \cite{Sch}). Therefore, the
transformation $D\mapsto D+y+\eta(y)-\infty^--\infty^+$ corresponds to
the transformation $D^{t}\mapsto D^{t}+\infty^-+\infty^+-y-\eta(y)$,
i.\ e.\ to the transformation described above applied to the
transposed operators. However, we may use Lemma~\ref{branchpoints} in
order to obtain a formula for the \Em{Baker--Akhiezer functions}
instead of the transposed \Em{Baker--Akhiezer functions}.
Let $\psi$ denote the \Em{Baker--Akhiezer function} of the
pair of potentials $(U,\Bar{U})$, and let $\phi$ and $\Op{J}\phi$ denote
the values of the transposed \Em{Baker--Akhiezer function} at the
two elements $k'$ and $-\Bar{k}'$ of the \Em{complex Fermi curve}.
Again in the following formulas $\phi$ may be multiplied
with some complex number without changing the transformation.
Due to Lemma~\ref{branchpoints} the two functions
$$\frac{\int\limits_{x}^{x+\xx{\gamma}}
                   \phi_2\psi_1dz+\phi_1\psi_2d\Bar{z}}
       {\exp(2\pi\sqrt{-1}g(\xx{\gamma},k-k'))-1}
\text{ and }
\frac{\int\limits_{x}^{x+\xx{\gamma}}
                   \Bar{\phi}_1\psi_1dz-\Bar{\phi}_2\psi_2d\Bar{z}}
     {\exp(2\pi\sqrt{-1}g(\xx{\gamma},k+\Bar{k}'))-1}$$
define for all $x\in\mathbb{R}^2$ meromorphic functions on the
\Em{complex Fermi curve} with  first--order poles at the elements
$k'$ and $-\Bar{k}'$, respectively. Moreover, due to the
normalization factor, in which $k$ denotes the corresponding function
on the \Em{complex Fermi curve}, these functions are cross sections
of the bundle on $\torus$ corresponding to $[k-k']$ and
$[k+\Bar{k}']$, respectively. Furthermore, due to the
normalization factor these functions do not depend
on the period $\xx{\gamma}$.
Their asymptotic behaviour near $\infty^{\pm}$ may be calculated as
\begin{equation*}\begin{split}
\frac{\int\limits_{x}^{x+\xx{\gamma}} \phi_2\psi_1dz+\phi_1\psi_2d\Bar{z}}
{\exp(2\pi\sqrt{-1}g(\xx{\gamma},k-k'))-1}&=
\begin{cases}
\frac{(\xx{\gamma}_1-\sqrt{-1}\xx{\gamma}_2)
                 \phi_1(x)\exp(2\pi\sqrt{-1}g(x,k))}
     {2\pi\sqrt{-1}g(\xx{\gamma},k)}
\left(1+\text{\bf{O}}(1/g(\xx{\gamma},k))\right)
&\text{at }\infty^-\\
\frac{(\xx{\gamma}_1+\sqrt{-1}\xx{\gamma}_2)
                 \phi_2(x)\exp(2\pi\sqrt{-1}g(x,k))}
     {2\pi\sqrt{-1}g(\xx{\gamma},k)}
\left(1+\text{\bf{O}}(1/g(\xx{\gamma},k))\right)
&\text{at }\infty^+.
\end{cases}\\
\frac{\int\limits_{x}^{x+\xx{\gamma}} 
\Bar{\phi}_1\psi_1dz-\Bar{\phi}_2\psi_2d\Bar{z}}
{\exp(2\pi\sqrt{-1}g(\xx{\gamma},k+\Bar{k}'))-1}&=
\begin{cases}
\frac{-(\xx{\gamma}_1-\sqrt{-1}\xx{\gamma}_2)
                  \Bar{\phi}_2(x)\exp(2\pi\sqrt{-1}g(x,k))}
     {2\pi\sqrt{-1}g(\xx{\gamma},k)}
\left(1+\text{\bf{O}}(1/g(\xx{\gamma},k))\right)
&\text{at }\infty^-\\
\frac{(\xx{\gamma}_1+\sqrt{-1}\xx{\gamma}_2)
                  \Bar{\phi}_1(x)\exp(2\pi\sqrt{-1}g(x,k))}
     {2\pi\sqrt{-1}g(\xx{\gamma},k)}
\left(1+\text{\bf{O}}(1/g(\xx{\gamma},k))\right)
&\text{at }\infty^+.
\end{cases}
\end{split}\end{equation*}
We conclude that
$$\psi'=\frac{\begin{pmatrix}
      \Bar{\phi}_2\\
      \Bar{\phi}_1
      \end{pmatrix}}
     {\phi_1\Bar{\phi}_1+\phi_2\Bar{\phi}_2}
\frac{\int\limits_{x}^{x+\xx{\gamma}} \phi_2\psi_1dz+\phi_1\psi_2d\Bar{z}}
{\exp(2\pi\sqrt{-1}g(\xx{\gamma},k-k'))-1}
+\frac{\begin{pmatrix}
       \phi_1\\
       -\phi_2
       \end{pmatrix}}
      {\phi_1\Bar{\phi}_1+\phi_2\Bar{\phi}_2}
\frac{\int\limits_{x}^{x+\xx{\gamma}} 
\Bar{\phi}_1\psi_1dz-\Bar{\phi}_2\psi_2d\Bar{z}}
{\exp(2\pi\sqrt{-1}g(\xx{\gamma},k+\Bar{k}'))-1}$$
has the properties of the \Em{Baker--Akhiezer function} corresponding
to the divisor $D+y+\eta(y)-\infty^--\infty^+$.
The derivatives of these functions are given by
\begin{align*}
\frac{\partial\int\limits_{x}^{x+\xx{\gamma}}
                   \phi_2\psi_1dz+\phi_1\psi_2d\Bar{z}}
{\exp(2\pi\sqrt{-1}g(\xx{\gamma},k-k'))-1}&=
\phi_2(x)\psi_1(x), &
\frac{\partial\int\limits_{x}^{x+\xx{\gamma}} 
\Bar{\phi}_1\psi_1dz-\Bar{\phi}_2\psi_2d\Bar{z}}
{\exp(2\pi\sqrt{-1}g(\xx{\gamma},k+\Bar{k}'))-1}&=
\Bar{\phi}_1(x)\psi_1(x),\\
\frac{\Bar{\partial}\int\limits_{x}^{x+\xx{\gamma}}
                    \phi_2\psi_1dz+\phi_1\psi_2d\Bar{z}}
{\exp(2\pi\sqrt{-1}g(\xx{\gamma},k-k'))-1}&=
\phi_1(x)\psi_2(x), &
\frac{\Bar{\partial}\int\limits_{x}^{x+\xx{\gamma}} 
\Bar{\phi}_1\psi_1dz-\Bar{\phi}_2\psi_2d\Bar{z}}
{\exp(2\pi\sqrt{-1}g(\xx{\gamma},k+\Bar{k}'))-1}&=
-\Bar{\phi}_2(x)\psi_2(x).
\end{align*}
Furthermore, the values of $\psi'$ at the points $k'$ and $-\Bar{k}'$
satisfy the equations
\begin{align*}
\begin{pmatrix}
U' & \partial\\
-\Bar{\partial} & \Bar{U}'
\end{pmatrix}
\frac{\begin{pmatrix}
      \Bar{\phi}_2\\
      \Bar{\phi}_1
      \end{pmatrix}}
     {\phi_1\Bar{\phi}_1+\phi_2\Bar{\phi}_2}&=0, &
\begin{pmatrix}
U' & \partial\\
-\Bar{\partial} & \Bar{U}'
\end{pmatrix}
\frac{\begin{pmatrix}
      \phi_1\\
      -\phi_2
      \end{pmatrix}}
     {\phi_1\Bar{\phi}_1+\phi_2\Bar{\phi}_2}&=0,
\end{align*} 
with the transformed potential
$$U'=\frac{(\partial\phi_2)\Bar{\phi}_1-\phi_2(\partial\Bar{\phi}_1)}
          {\phi_1\Bar{\phi}_1+\phi_2\Bar{\phi}_2}.$$
All these equations imply that $\psi'$ is indeed the
\Em{Baker--Akhiezer function} of the pair of potentials $(U',\Bar{U}')$,
if $\phi$ considered as a function from $\Delta$ to $\mathbb{C}^2$ has
no zero. Moreover, if in addition $\chi$ is an
element of the kernel of $\Op{D}(U,\Bar{U},k')$, such that the integral
$\int\limits_{x'}^{x} \phi_2\chi_1dz+\phi_1\chi_2d\Bar{z}$
yields a periodic function with respect to $x$, then for all
$x'\in\mathbb{R}^2$ the function
$$\chi'=\frac{\begin{pmatrix}
      \Bar{\phi}_2\\
      \Bar{\phi}_1
      \end{pmatrix}}
     {\phi_1\Bar{\phi}_1+\phi_2\Bar{\phi}_2}
\int\limits_{x'}^{x} \phi_2\chi_1dz+\phi_1\chi_2d\Bar{z}
+\frac{\begin{pmatrix}
       \phi_1\\
       -\phi_2
       \end{pmatrix}}
      {\phi_1\Bar{\phi}_1+\phi_2\Bar{\phi}_2}
\frac{\int\limits_{x}^{x+\xx{\gamma}} 
\Bar{\phi}_1\chi_1dz-\Bar{\phi}_2\chi_2d\Bar{z}}
{\exp(2\pi\sqrt{-1}g(\xx{\gamma},k'+\Bar{k}'))-1}$$
belongs to the kernel of the transformed Dirac operator
$\Op{D}(U',\Bar{U}',k')$. Due to Lemma~\ref{branchpoints} this
condition on $\chi$ is equivalent to the condition
$\left\langle\phi_2,\chi_1\right\rangle =0\text{ and }
\left\langle\phi_1,\chi_2\right\rangle =0.$

A direct calculation shows that for any element $y$ of the
normalization the two \Em{B\"acklund transformations} corresponding
to the values of the pullback of the (transposed)
\Em{Baker--Akhiezer function} to the normalization are inverse to each
other, in agreement with the corresponding transformations of the
divisors: $D\mapsto D+\infty^-+\infty^+-y-\eta(y)$ and
$D\mapsto D+y+\eta(y)-\infty^--\infty^+$.
We collect our calculations and observations in the following

\begin{Lemma}\label{Baecklund transformation}
Let $\psi$ be the \Em{Baker--Akhiezer function}
of the pair of potentials $(U,\Bar{U})$.
If $\chi$ belongs to the kernel of $\Op{D}(U,\Bar{U},k)$
and has no zero considered as a function
from $\Delta$ to $\mathbb{C}^2$, then the function
$$\psi'=\left(\frac{1}{\chi_1\Bar{\chi}_1+\chi_2\Bar{\chi}_2}
\begin{pmatrix}
-(\partial \chi_1)\Bar{\chi}_1+\chi_2(\partial \Bar{\chi}_2) &
\chi_1(\partial \Bar{\chi}_2)-(\partial \chi_1)\Bar{\chi}_2\\
\chi_2(\Bar{\partial}\Bar{\chi}_1)-(\Bar{\partial}\chi_2)\Bar{\chi}_1 &
-\chi_1(\Bar{\partial}\Bar{\chi}_1)+(\Bar{\partial}\chi_2)\Bar{\chi}_2
\end{pmatrix}+
\begin{pmatrix}
\partial & 0\\
0 & \Bar{\partial}
\end{pmatrix}\right)\psi$$
is the \Em{Baker--Akhiezer function}
of the pair of potentials $(U',\Bar{U}')$ with
$$U'=\frac{(\Bar{\partial}\chi_2)\Bar{\chi}_1-
                     \chi_2(\Bar{\partial}\Bar{\chi}_1)}
          {\chi_1\Bar{\chi}_1+\chi_2\Bar{\chi}_2},$$
and the \Em{complex Fermi curves} $\fermi(U,\Bar{U})$
and $\fermi(U',\Bar{U}')$ coincide. In particular, if $\chi$ is
the value of the pullback of $\psi$ under the normalization map at some
element $y$ of the normalization, then the corresponding divisors obey
the relation $D'\sim D+\infty^-+\infty^+-y-\eta(y)$. Conversely, if
$\psi$ belongs to the kernel of
$\Op{D}^{t}(U,\Bar{U},k)$ and has no zero considered as a function
from $\Delta$ to $\mathbb{C}^2$, then the function
$$\psi'=\frac{\begin{pmatrix}
      \Bar{\phi}_2\\
      \Bar{\phi}_1
      \end{pmatrix}}
     {\phi_1\Bar{\phi}_1+\phi_2\Bar{\phi}_2}
\frac{\int\limits_{x}^{x+\xx{\gamma}} \phi_2\psi_1dz+\phi_1\psi_2d\Bar{z}}
{\exp(2\pi\sqrt{-1}g(\xx{\gamma},k-k'))-1}
+\frac{\begin{pmatrix}
       \phi_1\\
       -\phi_2
       \end{pmatrix}}
      {\phi_1\Bar{\phi}_1+\phi_2\Bar{\phi}_2}
\frac{\int\limits_{x}^{x+\xx{\gamma}} 
\Bar{\phi}_1\psi_1dz-\Bar{\phi}_2\psi_2d\Bar{z}}
{\exp(2\pi\sqrt{-1}g(\xx{\gamma},k+\Bar{k}'))-1}$$
is again the \Em{Baker--Akhiezer function}
of the pair of potentials $(U',\Bar{U}')$ with
$$U'=\frac{(\partial\phi_2)\Bar{\phi}_1-\phi_2(\partial\Bar{\phi}_1)}
          {\phi_1\Bar{\phi}_1+\phi_2\Bar{\phi}_2}.$$
Again the \Em{complex Fermi curves} $\fermi(U,\Bar{U})$
and $\fermi(U',\Bar{U}')$ coincide.
If in particular $\phi$ is the value of the pullback of the transposed
\Em{Baker--Akhiezer function} under the normalization map
at some element $y$ of the normalization, then the corresponding
divisors obey the relation $D'\sim D+y+\eta(y)-\infty^--\infty^+$.
Finally, the special cases of these two \Em{B\"acklund transformations},
which correspond to the values of the pullback of the (transposed)
\Em{Baker Akhiezer function} at some element $y$ of the
normalization, are inverse to each other.\qed
\end{Lemma}
\index{B\"acklund transformation|)}

\subsubsection{A compactification of the \Em{isospectral sets}}
\label{subsubsection compactified isospectral}

For complex potentials $U$, the values of the \Em{first integral}
is equal to the square of the $\banach{2}$--norm times four.
Hence, due to Lemma~\ref{weakly continuous resolvent}
the \Em{isospectral set} of a \Em{complex Fermi curve}
corresponding to small values of the
\Em{first integral} is weakly compact. In the proof of
Lemma~\ref{existence of potentials} we
actually showed that the isospectral set of complex potentials $U$ is
isomorphic to an open dense subset of one connected real component of
the Picard group. Therefore, for small values of the
\Em{first integral} all elements of this connected real component of
the Picard group obey condition \Em{Divisor}~(iii)
(i.\ e.\ they are non special divisors of $\Tilde{\Spa{Y}}$).
For two other classes of \Em{complex Fermi curves} this is also
true: in \cite[Proposition~7.15]{Hi} it is proven that for all
\Em{real--$\sigma$--hyperelliptic} \Em{complex Fermi curves}
(i.\ e.\ the normalization of the quotient $\Spa{Y}_{\sigma}$ of $\Spa{Y}$
modulo the holomorphic involution $\sigma$
is isomorphic to $\mathbb{P}^1$ endowed with the
anti--holomorphic involution with one real cycle)
again all quaternionic complex line bundles of degree equal
to the genus plus one are non--special.
If a linear combination of the
functions $\xx{p}$ and $\yy{p}$ extends to a single--valued
meromorphic functions on the Riemann surface $\Spa{Y}$,
the same can be proven with the help of the arguments in
\cite[Theorem~8.3]{Sch}. In all these cases the proof of
Lemma~\ref{existence of potentials} actually shows that the
eigenfunctions corresponding to any element of the normalization of
the \Em{complex Fermi curve} has no zeroes considered as a
$\mathbb{C}^2$--valued function on $\Delta$.

We shall see that not
for all compact smooth Riemann surfaces $\Spa{Y}$ obeying conditions
\Em{Quasi--momenta}~(i)--(iii), the \Em{isospectral sets} of all complex
potentials $U$ are compact. In contrast for example to the case where
a linear combination of the functions $\xx{p}$ and $\yy{p}$
extends to a single--valued meromorphic function, the $\banach{2}$--norms of the
partial derivatives of the potential $U$, cannot be estimated in terms
of the higher integrals. Moreover, later we shall consider
non--connected compact Riemann surfaces $\Spa{Y}$ fulfilling
conditions \Em{Quasi--momenta}~(i)--(iii),
which are not \Em{complex Fermi curves},
i.\ e.\ the corresponding \Em{isospectral sets} are empty.
But the \Em{isospectral sets} have a natural compactification,
which we shall construct now.

\begin{Lemma}\label{compactification of isospectral sets}
\index{isospectral!set!compactification of an $\sim$}
\index{compactification!of an isospectral set}
Let $\Spa{Y}$ be a compact smooth Riemann surface with two
marked points $\infty^-$ and $\infty^+$
and $k$ a multi--valued meromorphic function,
which is holomorphic on $\Spa{Y}\setminus\{\infty^-,\infty^+\}$.
If $\Spa{Y},k,\infty^-,\infty^+)$ satisfies the conditions
\Em{Quasi--momenta}~(i)--(iii),
then the \Em{isospectral set} has a natural compactification
isomorphic to a real torus,
whose dimension is equal to the genus of $\Spa{Y}$ plus one.
The elements of this compactified \Em{isospectral sets}
are integral divisors on $\torus$
together with analytic complex potentials $U$ on the complement of the
support of the divisor.
The corresponding \Em{Baker--Akhiezer functions}
yield a family of elements in the kernel of the restriction of
the Dirac operator $\Op{D}(U,\Bar{U},k)$
to the complement of the support of the divisor,
which is parameterized by the corresponding \Em{complex Fermi curve}.
Finally, the \Em{first integral} is equal to the square of the
$\banach{2}(\torus)$--norm of the potential times four plus
$4\pi$ times the degree of the divisor.
\end{Lemma}

The corresponding \Em{Baker--Akhiezer} functions belong to the kernel
of \De{Finite rank Perturbations}~\ref{finite rank perturbations}
of the Dirac operators $\Op{D}(U,\Bar{U},k)$, whose support is the
support of the divisor
(compare with Section~\ref{subsubsection finite rank perturbations}). 

\begin{proof}
We extend the construction of the \Em{Baker--Akhiezer function}
in Lemma~\ref{existence of potentials} to all divisors of $\Tilde{\Spa{Y}}$,
which obey conditions \Em{Divisor}~(i)--(ii).
Here $\Tilde{\Spa{Y}}$ is again obtained from $\Spa{Y}$
by identifying the two marked points $\infty^{\pm}$ to an
ordinary double point.

First we claim that these divisors obey condition
\Em{Divisor}~(iii) for all $x\in\torus$ with the
exception of finitely many $x_1,\ldots,x_L\in\torus$.
Due to the consideration in the proof of
Lemma~\ref{existence of potentials}
it suffices to show, that for all divisors $D$ on $\Spa{Y}$ fulfilling
conditions \Em{Divisor}~(i')--(ii'),
there exists only finitely many
$x_1,\ldots,x_L\in\torus$, for which the sheaf
$\Sh{L}(x)\otimes\Sh{O}_{D-\infty^--\infty^+}$
has non--trivial global sections.
Due to the Brill--Noether theory \cite[Chapter~IV]{ACGH}
the set of all divisors on $\Spa{Y}$ fulfilling condition
\Em{Divisor}~(i'), which violate condition \Em{Divisor}~(iii'),
is a subvariety of the $(g+1)$--th symmetric product of $\Spa{Y}$.
We conclude that for all divisors $D$ on $\Spa{Y}$ fulfilling conditions
\Em{Divisor}~(i')--(ii') the set of all $x\in\torus$,
for which the sheaf
$\Sh{L}(x)\otimes\Sh{O}_{D-\infty^--\infty^+}$
has non--trivial global sections is a real subvariety.
Therefore, either this set is finite, or it contains
locally an analytic one--dimensional submanifold.
We exclude the latter case by contradiction.
Without loss of generality we may assume
that $x=0$ belongs to the analytic submanifold, and that
$\dim H^0(\Spa{Y},\Sh{L}(x)\otimes\Sh{O}_{D-\infty^--\infty^+})$
is locally constant.
In this case we have a flat family of locally free sheaves on $\Spa{Y}$
parameterized by some open disc in $\mathbb{C}$.
Due to \cite[Chapter~III. \S4. Theorem~4.7]{GPR}
the direct image of this family has a section on the disc.
The $n$--th derivatives of this section
with respect to the deformation parameter yield global sections of
$\Sh{O}_{D+(n-1)\infty^-+(n-1)\infty^+}$,
which are not global sections of
$\Sh{O}_{D+(n-2)\infty^-+(n-2)\infty^+}$.
Since for all $n\in\mathbb{N}_0$ the divisor
$D+(n-1)\infty^-+(n-1)\infty^+$ corresponds to a
quaternionic complex line bundle,
this implies that for all $n\in\mathbb{N}_0$
$$\dim H^0(\Spa{Y},\Sh{O}_{D+(n-1)\infty^-+(n-1)\infty^+})=
2n+\dim H^0(\Spa{Y},\Sh{O}_{D-\infty^--\infty^+})\geq 2n+2.$$
This contradicts to the Riemann--Roch theorem
\cite[Theorem~16.9 and Theorem~17.16]{Fo}.

We conclude that the \Em{Baker--Akhiezer function}
exists for all divisors $\Tilde{D}$ on $\Tilde{\Spa{Y}}$,
which satisfy conditions \Em{Divisor}~(i)--(ii).
But in general these \Em{Baker--Akhiezer functions}
have besides the poles on $\Tilde{\Spa{Y}}$,
which are independent of $x\in\torus$,
in addition finitely many poles $x_1,\ldots,x_L$ on $\torus$,
which are independent of $y\in\Tilde{\Spa{Y}}$.
As in the usual case, on the complements of the second type of poles
the asymptotic behaviour of the \Em{Baker--Akhiezer function}
determines a complex potential $U$,
such that the value of the \Em{Baker--Akhiezer function}
at any element of $\Tilde{\Spa{Y}}\setminus\{\infty^-,\infty^+\}$
restricted to this complement
$\left(\mathbb{R}^2\setminus\{x_1,\ldots,x_L\}\right)/\lattice$
is an element of the kernel of the Dirac operator
$\Op{D}(U,\Bar{U},k)$,
where $k$ is the value of the function $k$
at the same element of $\Tilde{\Spa{Y}}\setminus\{\infty^-,\infty^+\}$.
Furthermore, the $\banach{2}$--norms of these potentials are bounded,
since they are limits of potentials with bounded $\banach{2}$--norm.

If the divisor $\Tilde{D}$ satisfies conditions
\Em{Divisor}~(i)--(iii),
then the \Em{Baker--Akhiezer function} $\psi(x,y)=
\left(\begin{smallmatrix}
\psi_1(x,y)\\
\psi_2(x,y)
\end{smallmatrix}\right)$
is the unique function with the following properties:
\begin{description}
\index{condition!Baker--Akhiezer function (i)--(ii)|(}
\index{Baker--Akhiezer function ($\rightarrow$ eigenfunction)!condition
  $\sim$ (i)--(ii)|(}
\item[Baker--Akhiezer function (i)] For all $x\in\mathbb{R}^2$ the restriction
  of $\psi(x,\cdot)$ to $\Spa{Y}\setminus\{\infty^-,\infty^+\}$
  is a section of $\Sh{O}_{D}$.
\item[Baker--Akhiezer function (ii)] For all $x\in\mathbb{R}^2$ the function
  $\exp\left(-2\pi\sqrt{-1}g(x,k)\right)\psi(x,\cdot)$ is holomorphic
  in some neighbourhood of $\infty^-$ and $\infty^+$. At these points
  it takes the values $\left(\begin{smallmatrix}
  1\\
  0
  \end{smallmatrix}\right)$ and $\left(\begin{smallmatrix}
  0\\
  1
  \end{smallmatrix}\right)$, respectively.
\end{description}
\index{condition!Baker--Akhiezer function (i)--(ii)|)}
\index{Baker--Akhiezer function ($\rightarrow$ eigenfunction)!condition
  $\sim$ (i)--(ii)|)}
In some neighbourhood of $\infty^{\pm}$ we may expand the product
$\exp\left(-2\pi\sqrt{-1}g(x,k)\right)\psi(x,y)$
with respect to the local parameter $1/k_1$:
$$\exp\left(-2\pi\sqrt{-1}g(x,k)\right)\psi(x,y)=
\begin{cases}
\begin{pmatrix}
0\\
1
\end{pmatrix}+
\begin{pmatrix}
\triv{\psi}_{1,1}^-(x)/k_1\\
\triv{\psi}_{2,1}^-(x)/k_1
\end{pmatrix}
+\text{\bf{O}}(1/k_1^2)
&\text{at }\infty^-\\
\begin{pmatrix}
1\\
0
\end{pmatrix}+
\begin{pmatrix}
\triv{\psi}_{1,1}^+(x)/k_1\\
\triv{\psi}_{2,1}^+(x)/k_1
\end{pmatrix}
+\text{\bf{O}}(1/k_1^2)
&\text{at }\infty^+.
\end{cases}$$
Due to Lemma~\ref{Willmore functional} the function $k_2$ has the
asymptotic expansion
$$k_2=\pm\sqrt{-1}\left(k_1-
\frac{\willmore}{8\pi^2\vol(\torus)}
1/k_1\right)+\text{\bf{O}}(1/k_1^2)$$
 at $\infty^{\pm}$.
Since the \Em{Baker--Akhiezer function} lies in the kernel of the
Dirac operator $\Op{D}(U,\Bar{U},k)$,
we obtain the equations
\begin{align*}
U(x)+2\pi\sqrt{-1}\triv{\psi}_{2,1}^+(x)&=0
&U(x)\triv{\psi}_{1,1}^-(x)+
\frac{\sqrt{-1}\willmore}{8\pi\vol(\torus)}+
\partial\triv{\psi}_{2,1}^-(x)&=0\\
\Bar{U}(x)\triv{\psi}_{2,1}^+(x)-
\frac{\sqrt{-1}\willmore}{8\pi\vol(\torus)}-
\Bar{\partial}\triv{\psi}_{1,1}^+(x)&=0
&\Bar{U}-2\pi\sqrt{-1}\triv{\psi}_{1,1}^-(x)&=0
\end{align*}
These equation show that the \Em{first integral} is equal to
$\willmore=4\int\limits_{\torus}U(x)\Bar{U}(x)d^2x$.

We shall modify this calculation of the \Em{Baker--Akhiezer function}
in the case of poles with respect to $x$.
In this case the \Em{Baker--Akhiezer function} satisfies conditions
\Em{Baker--Akhiezer function}~(i)--(ii) only for all
$x$ in the complement of $(x_1+\lattice)\cup\ldots\cup (x_L+\lattice)$.
This function has a unique meromorphic extension to the product of
a small neighbourhood of $\mathbb{R}^2$ in $\mathbb{C}^2$ with
$\Spa{Y}\setminus\{\infty^-,\infty^+\}$. Moreover, the product
$\exp\left(-2\pi\sqrt{-1}g(x,k)\right)\psi(x,y)$
has also a meromorphic extension to the product of
a small neighbourhood of $\mathbb{R}^2$ in $\mathbb{C}^2$
with some neighbourhood of $\infty^{\pm}$.
We shall see that the functions $\triv{\psi}_{1,1}^+$ and
$\triv{\psi}_{2,1}^-$ have only first order poles
at the points $x_1,\ldots,x_L$. Dolbeault's Lemma
\cite[Chapter~I Section~D 2.~Lemma]{GuRo}
and the equations above imply
that in this case the \Em{first integral} is equal to
$\willmore=4\int\limits_{\Delta\setminus\{x_1,\ldots,x_L\}}U(x)\Bar{U}(x)d^2x$
minus $8\pi^2\sqrt{-1}$ times the sum over all Laurent coefficients of
$\triv{\psi}_{1,1}^+(x)$ and $\triv{\psi}_{2,1}^-(x)$
at $x_1,\ldots,x_L$ corresponding to the singular term
$\left((x-x_l)_1\pm\sqrt{-1}(x-x_l)_2\right)^{-1}$.
In order to calculate these Laurent coefficients it is convenient to
consider the limit
$$\psi_l^{\pm}(x,k_1)=\lim\limits_{\lambda\rightarrow\infty}
\exp\left(-2\pi\sqrt{-1}g(x_l,k)\right)\psi(x_l+x/\lambda,k_1\lambda)$$
in the neighbourhood of $\infty^{\pm}$, respectively.
Here we parameterize small neighbourhoods of $\infty^{\pm}$
by the local parameters $1/k_1$.
More precisely, the \Em{Baker--Akhiezer function} extends to a
function on a neighbourhood of $\mathbb{R}^2\times\infty^{\pm}$
in $\mathbb{C}^2\times\Spa{Y}$
with prescribed essential singularity at $\infty^{\pm}$,
which is meromorphic with respect to $x\in\mathbb{C}^2$.
The functions $\psi_{l}^{\pm}(x,k_1)$ are the pullbacks of the
products of $\exp\left(-2\pi\sqrt{-1}g(x_l,k)\right)$ with the
\Em{Baker--Akhiezer function} to the blowing up of
$\mathbb{C}^2\times\Spa{Y}$ at the points $x_l\times\infty^{\pm}$
(compare with \cite[Chapter~I Section~4]{Har}).
The condition \Em{Baker--Akhiezer function}~(ii)
determines the values of the products of these functions with
$\exp\left(-2\pi\sqrt{-1}g(x,k)\right)$ at $k_1=\infty$.
This implies that these products are meromorphic. 
Due to Lemma~\ref{existence of potentials} the potential $U$ is a
limit of potentials with finite $\banach{2}$--norms.
Therefore $U$ has finite $\banach{2}$--norm. Moreover, the limits
$\lim\limits_{\lambda\rightarrow\infty}U(x_l+x/\lambda)/\lambda$
are equal to zero for $x\neq 0$. Furthermore, since on the complement of
$(x_1+\lattice)\cup\ldots\cup (x_L+\lattice)$ the
\Em{Baker--Akhiezer functions} belong to the kernel of the
Dirac operator, for $x\neq 0$ the functions $\psi_l^{\pm}$
belong to the kernel of $\left(\begin{smallmatrix}
0 & \partial\\
-\Bar{\partial} & 0
\end{smallmatrix}\right)$. This implies that
the functions
$\psi_l^{\pm}(x,k_1)$ are of the form
$$\psi_l^{\pm}(x,k_1)=q_l(1/z^{\pm})\exp(z^{\pm})
\left(\begin{smallmatrix}
(1\pm 1)/2\\
(1\mp 1)/2
\end{smallmatrix}\right)
\text{ with }z^{\pm}=2\pi\sqrt{-1}k_1(x_1\pm\sqrt{-1}x_2).$$
Here $q_l$ is a rational function
with respect to $1/z^{\pm}\in\mathbb{P}^1$,
which is at $z^{\pm}=\infty$ equal to $1$,
and has poles only at $z^{\pm}=0$.

\begin{Remark}\label{rational KP solutions}
  The non--trivial components of these functions
  can be considered as \Em{Baker--Akhiezer functions} of rational
  solutions of the Kadomtsev--Petviashvili Equation.
  For example the non--trivial components of the functions
  $\psi_l^+(x,k)$ corresponding
  to the cases $\Set{N}_l=\{n\in\mathbb{Z}\mid-d\leq n\}\setminus
  \{-d+1,-d+3,\ldots,d-1\}$ are the \Em{Baker--Akhiezer functions}
  of the rational solutions of the KdV equation corresponding
  to the Lax operators $(\partial)^2-d(d+1)/(x_1+\sqrt{-1}x_2)^2$
  \cite{AMcKM}.
  The blowing up extracts from the \Em{Baker--Akhiezer function}
  exactly that part, which describes the singularity at $x_l$,
  and shifts the potential to infinity. Therefore, these
  \Em{Baker--Akhiezer functions} are elements of the kernel
  of finite rank perturbations of the Dirac operator
  without potential on $\mathbb{P}^1$
  (compare with Section~\ref{subsubsection finite rank perturbations}).
  The corresponding spectral curves have disconnected normalization
  similar to the \Em{complex Fermi curves}
  with disconnected normalization. The former
  \Em{Baker--Akhiezer functions}, however, correspond to
  rational solutions of the Kadomtsev--Petviashvili Equation,
  and the latter to elliptic solutions.
\end{Remark}

For all $x_l$ let $\Set{N}_l\subset\mathbb{Z}$ be the unique
subset of the integer numbers with the property that
$$\dim H^0\left(\Spa{Y},\Sh{O}_{D(x_l)+n(\infty^-+\infty^+)}\right)- 
\dim H^0\left(\Spa{Y},\Sh{O}_{D(x_l)+(n-1)(\infty^-+\infty^+)}\right)
=\begin{cases}
2&\text{if }n\in \Set{N}_l\\
0&\text{if }n\notin \Set{N}_l.
\end{cases}$$
We should remark that the space of global sections of the sheaves
$\Sh{O}_{D(x_l)+n\infty^-+n\infty^+}$ have even dimension,
since they correspond to quaternionic complex line bundles.
Due to the Riemann--Roch theorem
\cite[Theorem~16.9 and Theorem~17.16]{Fo}
$\Set{N}_l$ has as many negative elements
as there exists non--negative natural numbers,
which does not belong to $\Set{N}_l$.
Therefore, the Laurent expansion of $\psi_l^{\pm}$
at $z^{\pm}=\infty$ has no term of the form
$(z^{\pm})^n$, if $n\notin \Set{N}_l$.
These properties uniquely determine the functions $\psi_l^{\pm}$.
The corresponding rational functions $q_l$
are polynomials with integer coefficients.
In fact, if $-d_l$ denotes the smallest element of $\Set{N}_l$,
then $q_l$ has to be a polynomial of degree $d_l$,
whose lowest coefficient is equal to one.
Moreover, since $\Set{N}_l$ does not contain exactly $d_l$ elements of the
subset of all natural numbers, which are larger than $-d_l$,
this polynomial has to obey exactly $d_l$ equations.
If $1,a_1,\ldots,a_{d_l}$ are the coefficients
of the polynomial $q_l$, then the polynomial
$$p_l(z)=(z+1)(z+2)\ldots(z+d_l)+
a_1\left((z+2)\ldots(z+d_l)\right)+\ldots+a_{d_l}$$
has to vanish at all
$n\in\{n\in\mathbb{Z}\mid-d_l\leq n\}\setminus \Set{N}_l$.
Since $\{n\in\mathbb{Z}\mid-d_l\leq n\}\setminus \Set{N}_l$
has exactly $d_l$ elements,
the polynomial $p_l$ is the unique polynomial with highest
coefficient $1$, which vanishes at all elements of
$\{n\in\mathbb{Z}\mid-d_l\leq n\}\setminus \Set{N}_l$.
Obviously the coefficients of $q_l$
are integer linear combinations of the coefficients of $p_l$.
Therefore, the set $\Set{N}_l$ uniquely determines the polynomial $q_l$,
and these polynomials have integer coefficients.
Moreover, $a_1$ is equal to $-d_l(d_l+1)/2$
minus the sum over all elements of
$\{n\in\mathbb{Z}\mid-d_l\leq n\}\setminus \Set{N}_l$.
This natural number is equal to the sum over all elements in
$\Set{N}_l\cap-\mathbb{N}$ minus the sum over all elements in
$\mathbb{N}_0$, which does not belong to $\Set{N}_l$.
We define the local contribution to the Willmore energy
(compare with \cite{PP,BFLPP,FLPP})
at $x_l$ of the \Em{Baker--Akhiezer function}
$\willmore_{\text{\scriptsize\rm sing},x_l}(\psi)$
as this natural number times $-4\pi$,
which is always non--negative,
and zero if and only if $\Set{N}_l=\mathbb{N}_0$.
We conclude that the \Em{first integral} is equal to
\begin{align*}
\willmore&=4\int\limits_{\torus}U(x)\Bar{U}(x)d^2x
+\willmore_{\text{\scriptsize\rm sing}}&\text{with }
\willmore_{\text{\scriptsize\rm sing}}&=\sum\limits_{l=1}^{L}
\willmore_{\text{\scriptsize\rm sing},x_l}(\cdot).
\end{align*}
\index{local contribution!to the Willmore energy
$\willmore_{\text{\scriptsize\rm sing}}$}
\index{singularity!of the holomorphic structure}
We consider the points $x_1,\ldots,x_L$ as singularities of the
holomorphic structure of the corresponding quaternionic line bundle
in the sense of `quaternionic function theory'
(compare with \cite{PP,BFLPP,FLPP}).
In analogy to the
\De{Local contribution to the arithmetic genus}~\ref{local contribution}
of singularities in the sense of complex function theory
these singularities of `quaternionic function theory' have in
addition a local contribution to the Willmore energy.
\end{proof}

Due to this formula the \Em{isospectral sets} of all
\Em{complex Fermi curves} of finite \Em{geometric genus},
whose \Em{first integral} is smaller than $4\pi$ are compact
(compare with Remark~\ref{sobolev 1} and Remark~\ref{sobolev 2}).
Moreover, due to the explicit calculation of
$\willmore_{\text{\scriptsize\rm sing},x_l}(\cdot)$
in terms of $\Set{N}_l$ it is possible to classify all $\Set{N}_l$, with the
property that the corresponding
$\willmore_{\text{\scriptsize\rm sing},x_l}(\cdot)\leq 4n\pi$.
For some given $\Set{N}_l$
the local contribution to the Willmore energy
$\willmore_{\text{\scriptsize\rm sing},x_l}(\cdot)$ is equal to
the sum over all elements of $\mathbb{N}_0$,
which do not belong to $\Set{N}_l$,
minus the sum over all elements of $\Set{N}_l\cap -\mathbb{N}$.
For some given $d_l$ the corresponding
$\willmore_{\text{\scriptsize\rm sing},x_l}(\cdot)$ satisfies
the inequality $\willmore_{\text{\scriptsize\rm sing},x_l}(\cdot)\geq 4d_l\pi$,
and equality holds if and only if $\Set{N}_l=\mathbb{N}\cup\{-d_l\}$.
The possible $\Set{N}_l$'s with some given
$\willmore_{\text{\scriptsize\rm sing},x_l}(\cdot)$
decompose into groups with $m=1,\ldots,M$ negative elements, where
$$2M(M+1)\pi+2M(M-1)\pi=4M^2\pi\leq
\willmore_{\text{\scriptsize\rm sing},x_l}(\cdot).$$
Hence for each $m$ with
$4m^2\pi\leq\willmore_{\text{\scriptsize\rm sing},x_l}(\cdot)$
the sum $\Sigma_l^-$ of all negative elements of $\Set{N}_l$
and the sum $\Sigma_l^+$ of all elements of $\mathbb{N}_0$,
which do not belong to $\Set{N}_l$ satisfy
\begin{align*}
\willmore_{\text{\scriptsize\rm sing},x_l}(\cdot)
&=\Sigma_l^+-\Sigma_l^-&\text{and }
&-\frac{\willmore_{\text{\scriptsize\rm sing},x_l}(\cdot)}{4\pi}
+m(m-1)/2\leq\Sigma_l^-\leq-m(m+1)/2.
\end{align*}
Hence for each $m$ with $4m^2\pi\leq
\willmore_{\text{\scriptsize\rm sing},x_l}(\cdot)$
there are exactly 
$\willmore_{\text{\scriptsize\rm sing},x_l}(\cdot)/(4\pi)-m^2+1$
possible values of $\Sigma_l^-$ and $\Sigma_l^+$.
Moreover, for some given $m$, $\Sigma_l^-$ and $\Sigma_l^+$,
there are exactly so many $\Set{N}_l$'s, as the number of choices
of $m$ different elements of $-\mathbb{N}$,
whose sum is equal to $\Sigma_l^-$,
times the number of choices
of $m$ different elements of $\mathbb{N}_0$,
whose sum is equal to $\Sigma_l^+$.
In the following table we state all possible $\Set{N}_l$'s,
whose local contribution to the Willmore energy
$\willmore_{\text{\scriptsize\rm sing},x_l}(\cdot)$
is not larger than $20\pi$.
Instead of $\Set{N}_l$ we state all negative elements
of $\Set{N}_l$ together with all elements of $\mathbb{N}_0$,
which does not belong to $\Set{N}_l$:
$$\begin{array}{|c|ccccc|ccc|}
\hline
\willmore_{\text{\scriptsize\rm sing},x_l}(\cdot) & & & m=1 & & & & m=2 &\\
\hline
4\pi & -1,0 & & & & & & &\\
8\pi & -1,1 & -2,0 & & & & & &\\
12\pi & -1,2 & -2,1 & -3,0 & & & & &\\
16\pi & -1,3 & -2,2 & -3,1 & -4,0 & & -2,-1,0,1 & &\\
20\pi & -1,4 & -2,3 & -3,2 & -4,1 & -5,0 & -2,-1,0,2 & & -3,-1,0,1 \\
\hline
\end{array}$$

On Riemann surfaces there exists a one--to--one correspondence between
line bundles and divisors \cite[Chapter~III. \S29.11]{Fo},
which is compatible with the degrees of line bundles and divisors.
With the help of this correspondence those elements
in the compactification of the \Em{isospectral sets},
which do not belong to the \Em{isospectral sets},
can be shown to correspond to Dirac operators with potentials
acting on topological non--trivial bundles.
These generalized Dirac operators
\index{generalized!Dirac operator}
\index{Dirac operator $\Op{D}$!generalized $\sim$}
occur also in `quaternionic function theory'
developed by F.\ Pedit and U.\ Pinkall,
which describes immersions into $\mathbb{R}^4$.

If a \Em{complex Fermi curve} has singularities, then the
\Em{isospectral set} contains besides the potentials considered in
Lemma~\ref{existence of potentials},
whose corresponding \De{Structure sheaf}~\ref{structure sheaf}
is the normalization sheaf of the \Em{complex Fermi curve}
in addition potentials, whose corresponding
\De{Structure sheaf}~\ref{structure sheaf} differs from the
normalization sheaf. Due to Corollary~\ref{fixed points}~(ii),
on the one--sheeted coverings introduced in
\De{Structure sheaf}~\ref{structure sheaf} corresponding to
potentials of the form $(U,\Bar{U})$ the involution $\eta$ does not
have fixed points. Due to Lemma~\ref{singularity} all
one--sheeted coverings of the \Em{complex Fermi curve}
of finite \Em{geometric genus} have finite arithmetic genus.
The proofs of Lemma~\ref{existence of potentials} and
Lemma~\ref{compactification of isospectral sets} carry over to these
one--sheeted coverings and show that for any such one--sheeted covering
the real component of the Picard group, whose degree is equal to the
arithmetic genus of this one--sheeted covering plus one, parameterizes a
family of generalized potentials in the sense of
Lemma~\ref{compactification of isospectral sets}.
More precisely, sue to \cite[Chapter~9]{Sch} all sheaves, which
correspond to some potential with a given
\De{Structure sheaf}~\ref{structure sheaf}, are characterized as
coherent subsheaves of the sheaf of meromorphic functions,
whose degree is equal to the corresponding arithmetic genus,
if the two infinities are identified to an ordinary double point. But
for our purpose it is enough to consider only those potentials,
which correspond to locally free rank--one sheaves.
For any data $(\Spa{Y},\infty^-,\infty^+,k)$
fulfilling \Em{Quasi--momenta}~(i)--(iii) let
$\willmore_{\text{\scriptsize\rm sing},\max}(\Spa{Y},\infty^-,\infty^+,k)$
denote the maximal contribution of the singularities
to the Willmore energy of the corresponding compactified
\Em{isospectral set} described analogous to those described in
Lemma~\ref{compactification of isospectral sets},
whose \De{Structure sheaf}~\ref{structure sheaf} coincides
with the structure sheaf of $Y$.
Due to Lemma~\ref{compactification of isospectral sets}
this maximal contribution of the singularities to the Willmore energy
obeys the inequality
$$\willmore_{\text{\scriptsize\rm sing},\max}(\Spa{Y},\infty^-,\infty^+,k)
\leq\willmore(\fermi)=\willmore(\Spa{Y},\infty^-,\infty^+,k).$$
A \Em{Baker--Akhiezer function}
has a pole at $x\in\torus$,
if and only if the divisor $D(x)-\infty^--\infty^+$
is the pole divisor of a function $\func{f}$
fulfilling $\func{f}\eta^{\ast}\Bar{\func{f}}=-1$,
whose degree is equal to the arithmetic genus $g$ of $\Spa{Y}$
minus a positive odd number.
More precisely, the minimum of the degrees of all functions
$\func{f}$ fulfilling $\func{f}\eta^{\ast}\Bar{\func{f}}=-1$,
whose degree is equal to $g+1$ modulo two,
is larger or equal to
$g+1-2\willmore_{\text{\scriptsize\rm sing},\max}(\Spa{Y},\infty^-,\infty^+,k)
/(4\pi)$.
Moreover, the maximal contribution of the singularities to the
Willmore energy does not vanish,
if and only if there exists a function $\func{f}$
fulfilling $\func{f}\eta^{\ast}\Bar{\func{f}}=-1$,
whose degree is equal to $g$ minus a positive odd number.

If in addition ($\Spa{Y},\infty^-,\infty^+,k)$ fulfills condition
\Em{Quasi--momenta}~(iv), then the following condition implies
that the corresponding (generalized) potentials are real:
\begin{description}
\index{condition!divisor (iv)}
\index{divisor!condition $\sim$ (iv)}
\item[Divisor (iv)] There exists a function $\func{g}$ on
  $\Tilde{\Spa{Y}}$ which takes the same values at $\infty^{\pm}$, whose
  divisor is equal to the zero divisor of
  $dk_2$ minus $\Tilde{D}+\sigma(\Tilde{D})$.
\end{description}
Therefore, the real part of the compactification of the
\Em{isospectral sets} are real tori,
which are isomorphic to the Prym varieties of
$\Tilde{\Spa{Y}}$ together with the holomorphic involution $\sigma$
(see e.\ g.\ \cite[Chapter~12]{LB}).
In particular, if the \Em{isospectral sets} are compact
(e.\ g.\ if the \Em{complex Fermi curve} is
\Em{real--$\sigma$--hyperelliptic},
or if a linear combination of $\xx{p}$ and $\yy{p}$
extends to a single--valued meromorphic function,
or if the \Em{first integral} is smaller than $4\pi$)
then the real part of the \Em{isospectral} is not empty.
We do not know whether there exist \Em{complex Fermi curves},
whose \Em{isospectral sets} are not empty,
but whose real parts are empty.

\section{The moduli space}\label{section moduli}
\subsection{Variations of \Em{complex Fermi curves}}
\label{subsection variation zero energy}

In this section we want to investigate the behaviour of the map
$(V,W)\mapsto \fermi(V,W)$ under small variations of the
potentials.
In particular, we will characterize those variations
$(\var V,\var W)$ of the potentials, which induce no
variation of the \Em{complex Fermi curve}. For this purpose we use
the \De{Fundamental domain}~\ref{fundamental domain} and the
corresponding notations. If $\var V$ and $\var W$ are variations
of the potentials $V$ and $W$, then we may define a one--form on the
\Em{complex Fermi curve} $\fermi(V,W)$ by the formula
$$\Omega_{V,W}(\var V,\var W) = \var\yy{p}\cdot d\xx{p} -
\var\xx{p}\cdot d\yy{p}.
\label{omega}$$
Since in general the \Em{complex Fermi curves} of deformed potentials
are also deformed, it is not clear a priori whether this definition
makes sense. For example, either the condition $\var\xx{p}=0$
or the condition $\var\yy{p}=0$ may be used to identify locally
the \Em{complex Fermi curve} of the deformed potential with
the undeformed \Em{complex Fermi curve}. In the first
case we may consider $\yy{p}$ locally as a function of $\xx{p}$
and a deformation parameter $\parameter{t}$ and in the second case
$\xx{p}$ locally as a function of $\yy{p}$ and 
$\parameter{t}$. In the first case $\Omega_{V,W}(\var V,\var W)$ equals 
$\frac{\partial\yy{p}}{\partial \parameter{t}} d\xx{p}$ and in the
second case it equals
$-\frac{\partial\xx{p}}{\partial \parameter{t}} d\yy{p}$.
But as shown in \cite[Section~1.3]{GS1} the form does not depend
in which way we identify the different \Em{complex Fermi curves}.
In fact, due to Theorem~\ref{meromorph} locally
the \Em{complex Fermi curves} may be described by
an equation of the form $R(\xx{p},\yy{p})=0$.
Moreover, if $d$ is the number of sheets
of the covering map $\xx{p}$, then
due to the Weierstra{\ss} Preparation theorem
\cite[Chapter~I. Theorem~1.4]{GPR} $R$ can be chosen to be a
polynomial with respect to $\yy{p}$ of degree $d$,
whose coefficients are holomorphic functions depending on $\xx{p}$,
and whose highest coefficient is equal to $1$.
The equations
\begin{align*}
dR(\xx{p},\yy{p})=
\var R(\xx{p},\yy{p})+
\frac{\partial R(\xx{p},\yy{p})}{\partial \xx{p}}\var\xx{p}+
\frac{\partial R(\xx{p},\yy{p})}{\partial \yy{p}}\var\yy{p}
&=0\text{ and}\\
\frac{\partial R(\xx{p},\yy{p})}{\partial \xx{p}}d{\xx{p}}+
\frac{\partial R(\xx{p},\yy{p})}{\partial \yy{p}}d{\yy{p}}&=0
\end{align*}
imply the equation
$$\var\yy{p}d\xx{p}-\var\xx{p}d\yy{p}=
-\var R(\xx{p},\yy{p})\frac{d\xx{p}}{
\frac{\partial R(\xx{p},\yy{p})}{\partial \yy{p}}}=
\var R(\xx{p},\yy{p})\frac{d\yy{p}}{
\frac{\partial R(\xx{p},\yy{p})}{\partial \xx{p}}}.$$
Hence due to Lemma~\ref{dualizing sheaf}
$\Omega_{V,W}(\var V,\var W)$ is a well defined regular form.
If we consider $\lambda$ locally as a function of
$\xx{p},\yy{p}$, and $t$, then the choice $\var\xx{p}=0$ yields
\begin{multline*}
\Omega_{V,W}(\var V,\var W) = \frac{\partial\lambda}{\partial \parameter{t}}
\frac{d\xx{p}}{\frac{\partial\lambda}{\partial\yy{p}}} =
\tr\left(\Op{P}\comp
\begin{pmatrix}
\var V & 0\\
0 & \var W
\end{pmatrix}\right)\frac{d\xx{p}}
       {\frac{\partial\lambda}{\partial\yy{p}}}\\
=\tr\left(\Op{P}_{\xx{\gamma}}\comp
  \begin{pmatrix}
  0 & \frac{\var W}{\yy{\kappa}_2-\sqrt{-1}\yy{\kappa}_1}\\
  \frac{\var V}{\yy{\kappa}_2+\sqrt{-1}\yy{\kappa}_1} & 0
  \end{pmatrix}\right) \frac{d\xx{p}}{\pi}.
\end{multline*}
Here we have used the relation between
$\Op{P}$ and $\Op{P}_{\xx{\gamma}}$
established in the proof of Lemma~\ref{projection 2}.

\begin{Lemma} \label{regular form}
\index{form!regular $\sim$}
Let $\var V\in \banach{2}(\torus)$ and
$\var W\in \banach{2}(\torus)$ be variations of
$V$ and $W$, respectively. Then
$\Omega_{V,W}(\var V,\var W)=
\tr\left(\Op{P}_{\xx{\gamma}}\comp
  \begin{pmatrix}
  0 & \frac{\var W}{\yy{\kappa}_2-\sqrt{-1}\yy{\kappa}_1}\\
  \frac{\var V}{\yy{\kappa}_2+\sqrt{-1}\yy{\kappa}_1} & 0
  \end{pmatrix}\right)\displaystyle{\frac{d\xx{p}}{\pi}}$
is a regular one--form, on
$\fermi(V,W)$, independent of the choice of $\xx{\gamma}$ and
$\yy{\gamma}$.
\end{Lemma}

\begin{proof}
By Lemma~\ref{projection 2} (iv),
$\Op{P}_{\gamma} d\xx{p}$ is a regular one--form on
$\fermi(V,W)$. Thus also $\Omega_{V,W}(\var V,\var W)$ is a
regular one--form. All matrices
$\left(\begin{smallmatrix}
a & b\\
c & d
\end{smallmatrix}\right) \in SL(2,\mathbb{Z})$ induce
the modular transformations:
\begin{align*}
d\xx{p}_{\Delta'} & = a\;d\xx{p}_{\Delta}
+b\;d\yy{p}_{\Delta} &
\frac{\partial}{\partial\xx{p}_{\Delta'}}=
d\frac{\partial}{\partial\xx{p}_{\Delta}}
-c\frac{\partial}{\partial\yy{p}_{\Delta}} \\
d\yy{p}_{\Delta'} & = c\;d\xx{p}_{\Delta}
+d\;d\yy{p}_{\Delta} &
\frac{\partial}{\partial\yy{p}_{\Delta'}}
=-b \frac{\partial}{\partial\xx{p}_{\Delta}}
+ a\frac{\partial}{\partial\yy{p}_{\Delta}}
\end{align*}
On the \Em{complex Fermi curve} we have the following identity:
$$
\frac{\partial\lambda}{\partial\xx{p}} d\xx{p} +
  \frac{\partial\lambda}{\partial\yy{p}} d\yy{p} = 0.
$$
Therefore we obtain
$$
\frac{d\xx{p}_{\Delta'}}
{\frac{\partial\lambda}{\partial\yy{p}_{\Delta'}}} =
  \frac{a\cdot d\xx{p}_{\Delta} + b\cdot d\yy{p}_{\Delta}}
       {a\frac{\partial\lambda}{\partial\yy{p}_{\Delta}}-
        b\frac{\partial\lambda}{\partial\xx{p}_{\Delta}}} =
  \frac{a\frac{\partial\lambda}
              {\partial\yy{p}_{\Delta}}d\xx{p}_{\Delta} +
        b\frac{\partial\lambda}
       {\partial\yy{p}_{\Delta}}d\yy{p}_{\Delta}}
       {\frac{\partial\lambda}{\partial\yy{p}_{\Delta}}\left(
        a\frac{\partial\lambda}{\partial\yy{p}_{\Delta}}-
        b\frac{\partial\lambda}{\partial\xx{p}_{\Delta}}\right)} =
  \frac{\left(
        a\frac{\partial\lambda}{\partial\yy{p}_{\Delta}}-
        b\frac{\partial\lambda}{\partial\xx{p}_{\Delta}}\right)
                d\xx{p}_{\Delta}}
       {\frac{\partial\lambda}{\partial\yy{p}_{\Delta}}\left(
        a\frac{\partial\lambda}{\partial\yy{p}_{\Delta}}-
        b\frac{\partial\lambda}{\partial\xx{p}_{\Delta}}\right)} =
   \frac{d\xx{p}_{\Delta}}
              {\frac{\partial\lambda}{\partial\yy{p}_{\Delta}}}.
$$
On the other hand
it follows from the proof of (iii) of
Lemma~\ref{projection 2} that
$\Omega_{V,W}(\var V,\var W)$ is equal to
$$
\Omega_{V,W}(\var V,\var W)=
  \tr\left(\Op{P}\comp\begin{pmatrix}
  \var V & 0\\
  0 & \var W
  \end{pmatrix} \right)
  \frac{d\xx{p}}{\frac{\partial\lambda}{\partial\yy{p}}},
$$
and does not change under the modular transformations.
\end{proof}

This map $\Omega_{V,W}$ should be considered as the derivative of the
map $(V,W)\mapsto \fermi(V,W)$. Due to a general feature of the
integrable systems having Lax operators, the dual of this map is the
\index{isospectral!transformation!infinitesimal $\sim$}
infinitesimal action of the Lie algebra of the Picard group on the
\Em{isospectral sets}. More precisely, S\'{e}rre duality is a pairing
between the space of holomorphic one--forms and the Lie algebra of the
Picard group \cite[\S17. and \S21.6]{Fo},
and on the other hand the symplectic form is a pairing
on the tangent space of the phase space of the integrable system.
\index{tangent!space of the phase space}
Therefore, the dual of the map $\Omega_{V,W}$ with respect to
these pairings is a map from the Lie algebra of the Picard group into
the tangent space of the phase space. Furthermore, the
\Em{isospectral sets} correspond to open dense subsets of one
component of the Picard group (with specified degree). Therefore, the
Picard group acts on the \Em{isospectral sets}, and the derivative of
this action yields a map from the Lie algebra of the Picard group into
the tangent space of the phase space. The Hamiltonian structure fits
together with the function theoretic properties of the
spectral curves, and the second map is the dual of the first map.
To work out this structure completely implies a lot of analytical work
(compare with \cite{Sch}). In the following lemma we define the
infinitesimal action of the Lie algebra of the Picard group on the
\index{Picard group}\index{isospectral!transformation}
\Em{isospectral sets}. We show that this map is the dual map of
$\Omega_{V,W}$ with respect to
S\'{e}rre duality\index{S\'{e}rre duality} and the
symplectic form.
For this purpose we represent the elements of the
Lie algebra of the Picard group by those data, which enter in the
first Cousin Problem (see e.\ g.\ \cite[Chapter~5 2.1.]{Ah} and
\cite[Chapter~III. \S3.3.]{GPR}). These data are the
singular parts of local meromorphic functions, sometimes called
Mittag-Leffler distributions. If we identify those data with zero,
which are the singular parts of global meromorphic functions, the
resulting space is canonically isomorphic with the Lie algebra of the
Picard group \cite[\S18.2]{Fo}.

We shall use the notations introduced in the context of the
\De{Fundamental domain}~\ref{fundamental domain}, but the operators
under consideration act on the Hilbert bundle corresponding to the
\De{Trivialization}~\ref{trivialization}.

\begin{Lemma} \label{residue}
\index{residue}
Let $f$ be a meromorphic function with finitely many
poles on some open neighbourhood $\Set{U}$ of
$\fermi(V,W)/\lattice\dual$.
Then there exists a meromorphic function
$\triv{\Op{A}}_{f}(\xx{p})$ from the complex plane
$\xx{p}\in\mathbb{C}$ into the bounded operators from the
Banach spaces $\banach{p}(\torus)\times\banach{p}(\torus)$
into $\banach{p'}(\torus)\times\banach{p'}(\torus)$
for all $1<p'<p<\infty$ with the following properties:
\begin{description}
\item[(i)] For all $(n_1,n_2)\in\mathbb{Z}^2$ and all
  $\xx{p}\in\mathbb{C}$ the following identity is valid:
  $$\triv{\Op{A}}_{f}(\xx{p}) \comp
  \psi_{n_1\xx{\kappa}+n_2\yy{\kappa}} = 
  \psi_{n_1\xx{\kappa}+n_2\yy{\kappa}}\comp
  \triv{\Op{A}}_{f}(\xx{p}+n_1).$$
\item[(ii)] The commutator
  $\left[\triv{\Op{A}}_{f}(\xx{p}),
  \triv{\Op{D}}_{\xx{\gamma},\yy{\gamma}}(V,W,\xx{p})\right]$
  does not depend on $\xx{p}$ and is equal to an operator of the form
  $\var\triv{\Op{D}}_{\xx{\gamma},\yy{\gamma}}(V,W,\cdot)=
  \begin{pmatrix}
  0 & \frac{\var W_{f}}{\yy{\kappa}_2-\sqrt{-1}\yy{\kappa}_1}\\
  \frac{\var V_{f}}{\yy{\kappa}_2+\sqrt{-1}\yy{\kappa}_1} & 0
  \end{pmatrix}$, with some variations
  $\var V_{f},\var W_{f}\in
  \bigcap\limits_{q<\infty}\banach{q}(\torus)$. Moreover,
  this variation does not change the \Em{complex Fermi curve}.
  Equivalently, the form $\Omega_{V,W}(\var V_{f},\var W_{f})$ is
  identically zero.
\item[(iii)] For all elements $\triv{\chi}$ and $\triv{\xi}$ of
  $\bigcup\limits_{r>2}
  \banach{r}(\torus)\times\banach{r}(\torus)$ and all
  $\varepsilon,\delta>0$ there exists a $\delta'>0$, such that for all
  $\xx{p}$ in
  $$\left\{\xx{p}\in\mathbb{C}\mid
  |\xx{p}-\xx{p}'|>\delta\mod\mathbb{Z}
  \text{ for all values }\xx{p}'\text{ of the function }\xx{p}
  \text{ at the poles of }f\right\}$$
  with $|\xx{p}|>\delta'$ the function
  $\left|\langle\langle\triv{\chi},
  \triv{\Op{A}}(\xx{p})\triv{\xi}\rangle\rangle\right|$
  is bounded by $\varepsilon$.
\item[(iv)] For all variations $\var V,\var W
  \in \banach{2}(\torus)$ the total residue of the
  form $f\cdot\Omega_{V,W}(\var V,\var W)$
  on $\Set{U}$ is equal to the value of a symplectic form at
  $(\var V,\var W)$ and $(\var V_{f},\var W_{f})$:
  \index{symplectic form}
  $$\res\limits_{\Set{U}}
  \left(f\cdot\Omega_{V,W}(\var V,\var W)\right) =
  \frac{1}{2\pi^2\sqrt{-1}}
  \int\limits_{\torus}
  \left(\var V\var W_{f}-\var W\var V_{f}\right)d^2x.$$
\end{description}
\end{Lemma}

\begin{proof} If $f$ is a linear combination of two functions with
  two associated functions fulfilling conditions (i)--(iv), then the
  linear combination $\triv{\Op{A}}_{f}$
  of the two associated functions
  has the required properties. Hence we may assume that $f$ has a pole
  at one point $[k']\in \fermi(V,W)/\lattice\dual$. Let
  $\triv{\Op{F}}(\cdot)$ denote the local sum of
  $f\cdot\triv{\Op{P}}_{\xx{\gamma}}$ over all sheets of
  $\fermi(V,W)/\lattice\dual$ considered as a covering space
  over $\xx{p}\in\mathbb{C}$, which contain the element $[k']$
  (compare with Remark~\ref{local sum}).
  By definition this is a meromorphic function from some neighbourhood
  of $\xx{p}'$ into the finite rank operators on
  $\banach{2}(\torus)\times\banach{2}(\torus)$.
  Therefore, the singular part
  $\triv{\Op{F}}_{\text{\scriptsize\rm singular}}(\cdot)$
  is a meromorphic function
  on the whole plane $\xx{p}\in\mathbb{C}$.
  Now we claim that the infinite sum
  $$\triv{\Op{A}}_{f}(\xx{p}) = \sum\limits_{(n_1,n_2)\in\mathbb{Z}^2}
  \psi_{n_1\xx{\kappa}+n_2\yy{\kappa}}\comp
  \triv{\Op{F}}_{\text{\scriptsize\rm singular}}(\xx{p}+n_1) \comp
  \psi_{-n_1\xx{\kappa}-n_2\yy{\kappa}}$$
  converges in the strong operator topology. Due to
  Theorem~\ref{meromorph} the operator--valued function
  $\triv{\Op{F}}_{\text{\scriptsize\rm singular}}(\xx{p})$
  is a finite sum of operators of the form
  $\frac{|\triv{\psi}\rangle\rangle\langle\langle\triv{\phi}|}
  {(\xx{p}-\xx{p}')^{l}}:
  \triv{\chi}\mapsto \langle\langle\triv{\phi},
  \triv{\chi}\rangle\rangle (\xx{p}-\xx{p}')^{-l}\triv{\psi}$
  with elements $\triv{\phi}$ and $\triv{\psi}$ of the Banach spaces
  $\bigcap\limits_{q<\infty}
  \banach{q}(\torus)\times\banach{q}(\torus)$. For all
  $l\in\mathbb{N}$ and $\xx{q}\in [0,1]$ the sum
  $d_l(\xx{p},\xx{q})=\sum\limits_{n\in\mathbb{Z}}(\xx{p}-\xx{p}'+n)^{-l}
  \exp\left(2n\pi\sqrt{-1}\xx{q}\right)$
  converges. Let us calculate these limits.
  Fourier decomposition gives the relation 
  \begin{align*}
  d_1(\xx{p},\xx{q})&=\frac{2\pi\sqrt{-1}
  \exp\left(2\pi\sqrt{-1}(\xx{p}'-\xx{p})\xx{q}\right)}
  {1-\exp\left(2\pi\sqrt{-1}(\xx{p}'-\xx{p})\right)}
  &&\;\forall \xx{q}\in [0,1].
  \end{align*}
  Moreover, they satisfy the recursion equation
  \begin{align*}
  2\pi\sqrt{-1}(\xx{p}-\xx{p}')d_{l+1}(\xx{p},\xx{q})+
  \frac{\partial d_{l+1}(\xx{p},\xx{q})}{\partial \xx{q}}&=
  2\pi\sqrt{-1}d_{l}(\xx{p},\xx{q}) 
  &&\;\forall\xx{q}\in\mathbb{R}/\mathbb{Z}.
  \end{align*}
  Due to the Riemann--Lebesgue Lemma \cite[Theorem~IX.7]{RS2}
  for $l>1$ the functions $d_l(\xx{p},\xx{q})$ are continuous and
  periodic with period $1$ with respect to $\xx{q}$.
  This implies the following recursion relation for the functions
  $\Tilde{d}_{l+1}(\xx{p},\xx{q})=
  \exp\left(2\pi\sqrt{-1}(\xx{p}-\xx{p}')\xx{q}\right)
  d_l(\xx{p},\xx{q})$:
  $$\Tilde{d}_l(\xx{p},\xx{q}) =
  2\pi\sqrt{-1}\left(\int\limits_{0}^{\xx{q}}
  \Tilde{d}_{l}(\xx{p},\xx{q}')d\xx{q}'
  +\frac{1}{\left(\exp\left(2\pi\sqrt{-1}(\xx{p}-\xx{p}')\right)-1\right)}
  \int\limits_{0}^{1} \Tilde{d}_{l}(\xx{p},\xx{q}')
  d\xx{q}'\right).$$
  Since $\Tilde{d}_1(\xx{p},\xx{q})$ is equal to
  $\left(2\pi\sqrt{-1}\right)
  \left(1-\exp\left(2\pi\sqrt{-1}(\xx{p}'-\xx{p})\right)\right)^{-1}$,
  these functions $\Tilde{d}_l(\xx{p},\xx{q})$ are 
  polynomials with respect to $\xx{q}$ and
  $\left(\exp\left(2\pi\sqrt{-1}(\xx{p}'-\xx{p})\right)-1\right)^{-1}$
  times $\Tilde{d}_1(\xx{p},\xx{q})$.
  The unique solution of this recursion is given by the generating
  function
  $$\sum\limits_{l\in\mathbb{N}} t^{l-1}\Tilde{d}_{l}(\xx{p},\xx{q}) =
  \frac{\left(2\pi\sqrt{-1}\right)
  \left(1-\exp\left(2\pi\sqrt{-1}(\xx{p}'-\xx{p})\right)\right)^{-1}
  \exp\left(2\pi\sqrt{-1}\xx{q}t\right)}
  {1-\left(\exp\left(2\pi\sqrt{-1}(\xx{p}'-\xx{p})\right)-1\right)^{-1}
  \left(\exp\left(2\pi\sqrt{-1}t\right)-1\right)}.$$
  In fact, a transformation of the recursion relation
  \begin{align*}
  f_{l+1}(x)&=a\left(\int\limits_0^x
  f_l(x')dx'+b\int\limits_0^1f_l(x')dx'\right)
  &\text{with }f_1(x)&=c
  \end{align*}
  into an integral equation of the corresponding generating function
  yields the unique solution
  $$\sum\limits_{l\in\mathbb{N}}t^{l-1}f_l(x)=
  \displaystyle{\frac{c\exp(atx)}{1-b(\exp(at)-1)}}.$$
  Hence the functions $d_{l}(\xx{p},\xx{q})$ may be considered as
  meromorphic functions with respect to $\xx{p}\in\mathbb{C}$, whose
  values are bounded functions of $\xx{q}\in
  [0,1]$. Finally, using the identity $\sum\limits_{n\in\mathbb{Z}}
  \exp\left(2n\pi\sqrt{-1}\left(\yy{q}-\yy{q}'\right)\right)
  =\delta\left(\yy{q}-\yy{q}'\right)$ we conclude that
  the following infinite sum of finite rank operators
  converges in the strong operator topology to an operator, which is
  given by an integral kernel with respect to $\xx{q}\in [0,1]$
  depending on $\yy{q}\in \mathbb{R}/\mathbb{Z}$:
  \begin{multline*}
  \left(\sum\limits_{(n_1,n_2)\in\mathbb{Z}^2}
  \psi_{n_1\xx{\kappa}+n_2\yy{\kappa}}\comp
  \frac{|\triv{\psi}\rangle\rangle\langle\langle\triv{\phi}|}
       {(\xx{p}-\xx{p}'+n_1)^{l}}
  \comp\psi_{-n_1\xx{\kappa}-n_2\yy{\kappa}}
  \triv{\chi}\right)(\xx{q},\yy{q})=\\
  \vol(\torus)
  \int\limits_{\xx{q}-1}^{\xx{q}}
  d_l(\xx{p},\xx{q}-\xx{q}')
  \begin{pmatrix}
  \triv{\psi}_1(\xx{q},\yy{q})\triv{\phi}_1(\xx{q}',\yy{q})
  & \triv{\psi}_1(\xx{q},\yy{q})\triv{\phi}_2(\xx{q}',\yy{q})\\
  \triv{\psi}_2(\xx{q},\yy{q})\triv{\phi}_1(\xx{q}',\yy{q}) 
  & \triv{\psi}_2(\xx{q},\yy{q})\triv{\phi}_2(\xx{q}',\yy{q})
  \end{pmatrix} \begin{pmatrix}
  \triv{\chi}_1(\xx{q}',\yy{q})\\
  \triv{\chi}_2(\xx{q}',\yy{q})
  \end{pmatrix} d\xx{q}'.
  \end{multline*}
  This proves the claim. Thus, due to
  H\"older's inequality \cite[Theorem~III.1~(c)]{RS1} the function
  $\triv{\Op{A}}_f(\cdot)$ is a meromorphic function from the complex
  plane $\xx{p}\in\mathbb{C}$ into the bounded operators from
  $\banach{p}(\torus)\times\banach{p}(\torus)$ to
  $\banach{p'}(\torus)\times\banach{p'}(\torus)$
  for all $1<p'<p<\infty$. By definition this function
  $\triv{\Op{A}}_{f}(\cdot)$ satisfies condition (i).

  In order to show the existence of the commutator of this
  function with the unbounded-operator valued function
  $\triv{\Op{D}}_{\xx{\gamma},\yy{\gamma}}(V,W,\cdot)$
  in (ii) we have to ensure that the corresponding operators
  maps the domain of
  $\triv{\Op{D}}_{\xx{\gamma},\yy{\gamma}}(V,W,\cdot)$
  (which does not depend on $\xx{p}$) into this domain.
  With the help of Theorem~\ref{meromorph} it is quite straightforward
  to chose appropriate Banach spaces and to show this property. For all
  $(n_1,n_2)\in\mathbb{Z}^2$ and all $\xx{p}\in\mathbb{C}$ we have 
  $$\triv{\Op{D}}_{\xx{\gamma},\yy{\gamma}}(V,W,\xx{p})\comp
  \psi_{n_1\xx{\kappa}+n_2\yy{\kappa}}
  =\psi_{n_1\xx{\kappa}+n_2\yy{\kappa}}\comp
  \left(\triv{\Op{D}}_{\xx{\gamma},\yy{\gamma}}(V,W,\xx{p}+n_1)
  +n_2\pi\unity\right).$$
  Hence for all $\xx{p}\in\mathbb{C}$ the commutator is equal to
  \begin{eqnarray*}
  \left[\triv{\Op{A}}_{f}(\xx{p}),
  \triv{\Op{D}}_{\xx{\gamma},\yy{\gamma}}(\xx{p})\right]&=&
  \sum\limits_{(n_1,n_2)\in\mathbb{Z}^2}
  \left[\psi_{n_1\xx{\kappa}+n_2\yy{\kappa}}\comp
  \triv{\Op{F}}_{\text{\scriptsize\rm singular}}(\xx{p}+n_1) \comp
  \psi_{-n_1\xx{\kappa}-n_2\yy{\kappa}},
  \triv{\Op{D}}_{\xx{\gamma},\yy{\gamma}}(\xx{p})\right]
  \\&=&\sum\limits_{(n_1,n_2)\in\mathbb{Z}^2}
  \psi_{n_1\xx{\kappa}+n_2\yy{\kappa}}\comp
  \left[\triv{\Op{F}}_{\text{\scriptsize\rm singular}}(\xx{p}+n_1),
  \triv{\Op{D}}_{\xx{\gamma},\yy{\gamma}}(\xx{p}+n_1)\right] \comp
  \psi_{-n_1\xx{\kappa}-n_2\yy{\kappa}}.\end{eqnarray*}
  The operator--valued functions $\triv{\Op{F}}(\cdot)$ and
  $\triv{\Op{D}}_{\xx{\gamma},\yy{\gamma}}(\cdot)$
  commute pointwise, because the projection
  $\triv{\Op{P}}_{\xx{\gamma}}$ is the spectral projection of
  $\triv{\Op{D}}_{\xx{\gamma},\yy{\gamma}}(\cdot)$.
  If we insert the explicit form
  of the operator--valued function
  $\triv{\Op{D}}_{\xx{\gamma},\yy{\gamma}}(\xx{p})$,
  then we obtain
  $$\left[\triv{\Op{A}}_{f}(\xx{p}),
  \triv{\Op{D}}_{\xx{\gamma},\yy{\gamma}}(\xx{p})\right] =
  \pi\sum\limits_{\kappa\in\lattice\dual}
  \psi_{\kappa}\comp
  \left[\triv{\Op{F}}_{-1},\begin{pmatrix}
  \frac{\xx{\kappa}_2-\sqrt{-1}\xx{\kappa}_1}
  {\yy{\kappa}_2-\sqrt{-1}\yy{\kappa}_1} & 0\\
  0 & \frac{\xx{\kappa}_2+\sqrt{-1}\xx{\kappa}_1}
  {\yy{\kappa}_2+\sqrt{-1}\yy{\kappa}_1} 
  \end{pmatrix} \right] \comp
  \psi_{-\kappa},$$
  where $\triv{\Op{F}}_{-1}$ is the residue
  of the operator--valued form
  $\triv{\Op{F}}(\xx{p})d\xx{p}$ at the pole $\xx{p}=\xx{p}'$.
  This shows that the commutator does not depend on $\xx{p}$ and that
  the diagonal elements with respect to the decomposition
  $\banach{2}(\torus)\times\banach{2}(\torus)$
  vanishes. Moreover, since
  $\sum\limits_{\kappa\in\lattice\dual}\psi_{\kappa}(x)$ is equal to
  $\vol(\torus)$ times
  Dirac's $\delta$--function $\delta(x)$,
  for all finite rank operators of the form
  $\triv{\chi}\mapsto \langle\langle\triv{\phi},
  \triv{\chi}\rangle\rangle \triv{\psi}$
  with elements $\triv{\phi}$ and $\triv{\psi}$ of
  $\bigcap\limits_{q<\infty}\banach{q}(\torus)$, the infinite sum
  $$\triv{\chi}\longmapsto \sum\limits_{\kappa\in\lattice\dual}
  \langle\langle\triv{\phi},
  \psi_{-\kappa}\triv{\chi}\rangle\rangle
  \psi_{\kappa} \triv{\psi}$$
  converges in the strong operator topology to the operator
  $\vol(\torus)$ times
  $\left(\begin{smallmatrix}
  \triv{\psi}_1\triv{\phi}_1 & \triv{\psi}_1\triv{\phi}_2\\
  \triv{\psi}_2\triv{\phi}_1 & \triv{\psi}_2\triv{\phi}_2
  \end{smallmatrix}\right)$, where the functions
  $\triv{\psi}_1,\triv{\psi}_2,\triv{\phi}_1$ and $\triv{\phi}_2$
  are considered as the operators
  of multiplication with these functions.
  Hence the commutator is equal to a variation
  $\var\triv{\Op{D}}_{\xx{\gamma},\yy{\gamma}}$ associated
  to unique variations $\var V_{f}$ and $\var W_{f}$ of the
  potentials. Due to Lemma~\ref{regular form} the form
  $\Omega_{V,W}(\var V_{f},\var W_{f})$ is equal to
  $$\Omega_{V,W}(\var V_{f},\var W_{f})=\tr \left(
  \triv{\Op{P}}_{\xx{\gamma}}\comp 
  \var\triv{\Op{D}}_{\xx{\gamma},\yy{\gamma}}\right)\frac{d\xx{p}}{\pi}=
  \tr \left(\triv{\Op{P}}_{\xx{\gamma}}\comp 
  \left[\triv{\Op{A}}_{f}(\xx{p}),\triv{\Op{D}}_{\xx{\gamma},
  \yy{\gamma}}(\xx{p})\right]\right)\frac{d\xx{p}}{\pi}$$
  for all $\xx{p}\neq \xx{p}'$.
  If we choose $\xx{p}$ to be equal to the
  corresponding value of the point,
  where $\triv{\Op{P}}_{\xx{\gamma}}$ is
  evaluated, then the projection
  $\triv{\Op{P}}_{\xx{\gamma}}$ commutes with
  the operator
  $\triv{\Op{D}}_{\xx{\gamma},\yy{\gamma}}(\xx{p})$
  and the trace vanishes. Thus the form vanishes on an
  open dense subset of the \Em{complex Fermi curve} and therefore on
  the whole \Em{complex Fermi curve}. This proves (ii).

  For all elements $\triv{\chi}\in \banach{1}(\mathbb{R}/\mathbb{Z})$
  the Riemann--Lebesgue Lemma \cite[Theorem~IX.7]{RS2}
  shows that the estimate
  $$\left|\int\limits_{0}^{1}\chi(\xx{q}')
  \exp\left(2\pi\sqrt{-1}(\xx{p}'-\xx{p})\xx{q}'\right)
  d\xx{q}'\right|\leq
  \exp\left(2\pi\Im(\xx{p}-\xx{p}')\right)\cdot
  |\Check{f}(\Re(\xx{p}'-\xx{p})|$$ holds with some continuous function
  $\Check{f}$ on $\mathbb{R}$, which vanishes at infinity. Due to the
  form of the functions $d_l(\xx{p},\xx{q})$, this implies the
  estimate~(iii).

  The total residue of the form
  $f\cdot \Omega_{V,W}$ on $\Set{U}$ is equal to
  \begin{eqnarray*}
  \res\limits_{\Set{U}}\left(f\cdot\Omega_{V,W}\right) &=&
  \res\limits_{\Set{U}}\left(\tr\left(
  f\cdot \triv{\Op{P}}_{\xx{\gamma}}\comp 
  \begin{pmatrix}
  0 & \frac{\var W}{\yy{\kappa}_2-\sqrt{-1}\yy{\kappa}_1}\\
  \frac{\var V}{\yy{\kappa}_2+\sqrt{-1}\yy{\kappa}_1} & 0
  \end{pmatrix}\right)\frac{d\xx{p}}{\pi}\right)\\ &=&
  \tr\left(\res\limits_{\xx{p}=\xx{p}'}\left(
  \triv{\Op{F}}(\xx{p})\comp
  \begin{pmatrix}
  0 & \frac{\var W}{\yy{\kappa}_2-\sqrt{-1}\yy{\kappa}_1}\\
  \frac{\var V}{\yy{\kappa}_2+\sqrt{-1}\yy{\kappa}_1} & 0
  \end{pmatrix}\frac{d\xx{p}}{\pi}\right)\right) \\ &=&
  \frac{1}{\pi}\tr\left(\triv{\Op{F}}_{-1}\comp
  \begin{pmatrix}
  0 & \frac{\var W}{\yy{\kappa}_2-\sqrt{-1}\yy{\kappa}_1}\\
  \frac{\var V}{\yy{\kappa}_2+\sqrt{-1}\yy{\kappa}_1} & 0
  \end{pmatrix}\right).
  \end{eqnarray*}
  Since $\triv{\Op{F}}_{-1}$ is a finite rank operator,
  it is an integral operator.
  The trace of the composition of an integral operator with
  the operator of multiplication with some function is equal to
  the integral over the product of the diagonal part of the
  integral kernel and the function. On the other hand the
  infinite sum 
  $\sum\limits_{\kappa\in\lattice\dual}
  \psi_{\kappa}\comp\triv{\Op{F}}_{-1}\comp
  \psi_{-\kappa}$ converges in the
  strong operator topology to the operator of
  $\vol(\torus)$ times the
  multiplication with the diagonal part of
  the integral kernel of $\triv{\Op{F}}_{-1}$.
  By definition of $\var V_f$ and $\var W_f$ we have
  $$\left[\sum\limits_{\kappa\in\lattice\dual}
  \psi_{\kappa}\comp\triv{\Op{F}}_{-1} \comp
  \psi_{-\kappa},\begin{pmatrix}
  \frac{\xx{\kappa}_2-\sqrt{-1}\xx{\kappa}_1}
  {\yy{\kappa}_2-\sqrt{-1}\yy{\kappa}_1} & 0\\
  0 & \frac{\xx{\kappa}_2+\sqrt{-1}\xx{\kappa}_1}
  {\yy{\kappa}_2+\sqrt{-1}\yy{\kappa}_1} 
  \end{pmatrix} \right]=\frac{1}{\pi}
  \begin{pmatrix}
  0 & \frac{\var W_{f}}{\yy{\kappa}_2-\sqrt{-1}\yy{\kappa}_1}\\
  \frac{\var V_{f}}{\yy{\kappa}_2+\sqrt{-1}\yy{\kappa}_1} & 0
  \end{pmatrix}.$$
  This implies
  $\displaystyle{\frac{1}{\pi\vol(\torus)}\sum\limits_{\kappa\in\lattice\dual}
  \psi_{\kappa}\comp\triv{\Op{F}}_{-1} \comp
  \psi_{-\kappa}=}$
  \begin{equation*}\begin{split}
  &=\frac{1}{\pi^2\vol(\torus)}\begin{pmatrix}
  \cdot & \left(\frac{\xx{\kappa}_2+\sqrt{-1}\xx{\kappa}_1}
  {\yy{\kappa}_2+\sqrt{-1}\yy{\kappa}_1}-
  \frac{\xx{\kappa}_2-\sqrt{-1}\xx{\kappa}_1}
  {\yy{\kappa}_2-\sqrt{-1}\yy{\kappa}_1}\right)^{-1}
  \frac{\var W_{f}}{\yy{\kappa}_2-\sqrt{-1}\yy{\kappa}_1}\\
  \left(\frac{\xx{\kappa}_2-\sqrt{-1}\xx{\kappa}_1}
  {\yy{\kappa}_2-\sqrt{-1}\yy{\kappa}_1}-
  \frac{\xx{\kappa}_2+\sqrt{-1}\xx{\kappa}_1}
  {\yy{\kappa}_2+\sqrt{-1}\yy{\kappa}_1}\right)^{-1}
  \frac{\var V_{f}}{\yy{\kappa}_2+\sqrt{-1}\yy{\kappa}_1} & \cdot
  \end{pmatrix}\\
  &=\frac{g(\yy{\kappa},\yy{\kappa})}{2\pi^2\sqrt{-1}}
  \begin{pmatrix}
  \cdot & \frac{\var W_{f}}{\yy{\kappa}_2-\sqrt{-1}\yy{\kappa}_1}\\
  \frac{-\var V_{f}}{\yy{\kappa}_2+\sqrt{-1}\yy{\kappa}_1} & \cdot
  \end{pmatrix},
  \end{split}\end{equation*}
  and the formula (iv) follows.
\end{proof}

The integral on the left hand side of the formula~(iv)
is equal to the value of some unique symplectic form
on the Hilbert space of all
potentials evaluated at the variations $(\var V,\var W)$ and
$(\var V_{f}, \var W_{f})$. Due to (ii) the kernel of the map
$f\mapsto (\var V_{f}, \var W_{f})$ contains data $f$ entering in
first Cousin Problem, whose corresponding operator--valued function
$\triv{\Op{A}}_f(\cdot)$ commute with the unbounded-operator function
$\triv{\Op{D}}_{\xx{\gamma},\yy{\gamma}}(V,W,\cdot)$. The
eigenvalue of this function $\triv{\Op{A}}_f(\cdot)$ yields a global
meromorphic function on the \Em{complex Fermi curve}. Due to the
definition of the function $\triv{\Op{A}}_f(\cdot)$, the singular part
of this meromorphic function coincides with the singular part of $f$
on $\Set{U}$, and has no other poles. Therefore, this function is a
solution of the corresponding first Cousin Problem. More precisely, it
can be shown that on the complement of the singular part described in
the initial data $f$ this function is holomorphic with respect to the
structure sheaf of all eigenvalues of holomorphic functions commuting
with the function
$\triv{\Op{D}}_{\xx{\gamma},\yy{\gamma}}(V,W,\cdot)$.
Furthermore, the points $\infty^{\pm}$ are
identified to an ordinary double point
(compare with \cite{Sch}). In fact property (iii) implies that
these functions vanish at $\infty^{\pm}$. If we identify all data $f$
entering in the first Cousin problem with zero, which are the
singular parts of global meromorphic functions,
then the quotient space is naturally isomorphic to the Lie algebra
of the Picard group \cite[\S18.2]{Fo}.
We conclude that the map
$f\mapsto (\var V_{f}, \var W_{f})$ induces a map from the
Lie algebra of the Picard group into the tangent space of the
phase space. The pairing on the left hand side of the formula~(iv)
between the space of regular one--forms and the Lie algebra of the
Picard group is given by S\'{e}rre Duality \cite[\S17]{Fo}.
As a remarkable consequence of the formula~(iv) the image of the map
$\Omega_{V,W}$ described in Lemma~\ref{regular form} is surjective, if
the \Em{complex Fermi curve} has finite genus
(compare with Proposition~\ref{smooth moduli}).
We expect that a similar statement holds also in
the general case (compare with \cite[Chapter~7]{Sch}).

\begin{Remark}\label{trivial symplectic form}
The reduction to the fixed points of the involution $(V,W)\mapsto(W,V)$
(compare with Section~\ref{subsection reductions})
destroys this structure. More precisely, for geometric algebraic solutions
(corresponding to \Em{complex Fermi curves} of finite geometric genus)
the number of independent integrals is equal to the geometric genus
of the quotient of the \Em{complex Fermi curve}
modulo the holomorphic involution $\sigma$.
On the other hand the number of degrees of freedom of the
\Em{isospectral sets} (i.\ e.\ the set of all real potentials $U$,
whose \Em{complex Fermi curves} are equal to some given
\Em{complex Fermi curve}) is equal to the dimension of the
Prym varieties associated to the \Em{complex Fermi curves} together
with the holomorphic involution  $\sigma$ 
(see e.\ g.\ \cite[Chapter~12]{LB}).
Due to the Riemann--Hurwitz formula the
difference of these numbers is proportional to the number of
fixed points of the holomorphic involution.
Any bound on this difference would yield a bound
on the geometric genus of spectral curves
corresponding to minimal tori in $S^3$. 
\end{Remark}

Let us now investigate the compatibility of the involutions and the
corresponding relations with the map $f\mapsto
(\var V_{f},\var W_{f})$. The operator $\left(\begin{smallmatrix}
0 & \yy{\kappa}_2+\sqrt{-1}\yy{\kappa}_1\\
\yy{\kappa}_2-\sqrt{-1}\yy{\kappa}_1 & 0
\end{smallmatrix}\right)$ is self--adjoint and satisfies the relation
$$\Op{J}\comp\left(\begin{smallmatrix}
0 & \yy{\kappa}_2+\sqrt{-1}\yy{\kappa}_1\\
\yy{\kappa}_2-\sqrt{-1}\yy{\kappa}_1 & 0
\end{smallmatrix}\right) +\left(\begin{smallmatrix}
0 & \yy{\kappa}_2+\sqrt{-1}\yy{\kappa}_1\\
\yy{\kappa}_2-\sqrt{-1}\yy{\kappa}_1 & 0
\end{smallmatrix}\right)^{t}\comp\Op{J} = 0.$$
Hence the operator
$\triv{\Op{D}}_{\xx{\gamma},\yy{\gamma}}(\xx{p})$
satisfies the relations
\begin{eqnarray*}
-\triv{\Op{D}}_{\xx{\gamma},\yy{\gamma}}^{t}(-\xx{p})
&=&-\left(\begin{smallmatrix}
0 & \yy{\kappa}_2+\sqrt{-1}\yy{\kappa}_1\\
\yy{\kappa}_2-\sqrt{-1}\yy{\kappa}_1 & 0
\end{smallmatrix}\right)^{-1}\comp
\Op{J}^{-1}\comp\triv{\Op{D}}_{\xx{\gamma},\yy{\gamma}}
(\xx{p})\comp\left(\begin{smallmatrix}
0 & \yy{\kappa}_2+\sqrt{-1}\yy{\kappa}_1\\
\yy{\kappa}_2-\sqrt{-1}\yy{\kappa}_1 & 0
\end{smallmatrix}\right)^{t}\comp\Op{J} \\
&=&\left(\begin{smallmatrix}
0 & \yy{\kappa}_2+\sqrt{-1}\yy{\kappa}_1\\
\yy{\kappa}_2-\sqrt{-1}\yy{\kappa}_1 & 0
\end{smallmatrix}\right)^{-1}\comp
\Op{J}^{-1}\comp\triv{\Op{D}}_{\xx{\gamma},\yy{\gamma}}
(\xx{p})\comp\Op{J}\comp
\left(\begin{smallmatrix}
0 & \yy{\kappa}_2+\sqrt{-1}\yy{\kappa}_1\\
\yy{\kappa}_2-\sqrt{-1}\yy{\kappa}_1 & 0
\end{smallmatrix}\right) \\ &=&
\left(\begin{smallmatrix}
-\yy{\kappa}_2+\sqrt{-1}\yy{\kappa}_1 & 0\\
0 & \yy{\kappa}_2+\sqrt{-1}\yy{\kappa}_1
\end{smallmatrix}\right)^{-1}\comp
\triv{\Op{D}}_{\xx{\gamma},\yy{\gamma}}(\xx{p})\comp
\left(\begin{smallmatrix}
-\yy{\kappa}_2+\sqrt{-1}\yy{\kappa}_1 & 0\\
0 & \yy{\kappa}_2+\sqrt{-1}\yy{\kappa}_1
\end{smallmatrix}\right),\\
\triv{\Op{D}}_{\xx{\gamma},\yy{\gamma}}^{\ast}(\Bar{\xx{p}})
&=&\left(\begin{smallmatrix}
0 & \yy{\kappa}_2+\sqrt{-1}\yy{\kappa}_1\\
\yy{\kappa}_2-\sqrt{-1}\yy{\kappa}_1 & 0
\end{smallmatrix}\right)^{-1}\comp
\triv{\Op{D}}_{\xx{\gamma},\yy{\gamma}}(\xx{p})\comp
\left(\begin{smallmatrix}
0 & \yy{\kappa}_2+\sqrt{-1}\yy{\kappa}_1\\
\yy{\kappa}_2-\sqrt{-1}\yy{\kappa}_1 & 0
\end{smallmatrix}\right),\\
-\Bar{\triv{\Op{D}}}_{\xx{\gamma},\yy{\gamma}}(-\Bar{\xx{p}})
&=&\Op{J}^{-1}\comp
\triv{\Op{D}}_{\xx{\gamma},\yy{\gamma}}(\xx{p})\comp\Op{J}.
\end{eqnarray*}
Therefore, the projection
$\triv{\Op{P}}_{\xx{\gamma}}$ satisfies the relations
\begin{eqnarray*}
\triv{\Op{P}}_{\xx{\gamma}} &=&
\left(\begin{smallmatrix}
-\yy{\kappa}_2+\sqrt{-1}\yy{\kappa}_1 & 0\\
0 & \yy{\kappa}_2+\sqrt{-1}\yy{\kappa}_1
\end{smallmatrix}\right)
\comp\sigma^{\ast}\left(\triv{\Op{P}}_{\xx{\gamma}}^{t}\right)\comp
\left(\begin{smallmatrix}
-\yy{\kappa}_2+\sqrt{-1}\yy{\kappa}_1 & 0\\
0 & \yy{\kappa}_2+\sqrt{-1}\yy{\kappa}_1
\end{smallmatrix}\right)^{-1}\\
&=&\left(\begin{smallmatrix}
0 & \yy{\kappa}_2+\sqrt{-1}\yy{\kappa}_1\\
\yy{\kappa}_2-\sqrt{-1}\yy{\kappa}_1 & 0
\end{smallmatrix}\right)
\comp\rho^{\ast}\left(\triv{\Op{P}}_{\xx{\gamma}}^{\ast}\right)\comp
\left(\begin{smallmatrix}
0 & \yy{\kappa}_2+\sqrt{-1}\yy{\kappa}_1\\
\yy{\kappa}_2-\sqrt{-1}\yy{\kappa}_1 & 0
\end{smallmatrix}\right)^{-1}\\
&=&\Op{J}\comp\eta^{\ast}\left(
\Bar{\triv{\Op{P}}}_{\xx{\gamma}}\right)\comp
\Op{J}^{-1}.
\end{eqnarray*}
The following lemma is verified by some easy calculations.

\begin{Lemma} \label{compatibility}
The map $f\mapsto \triv{\Op{A}}_{f}$ has the following properties:
\begin{description}
\item[(i)] If the pair of potentials is of the form $(U,\Bar{U})$
  and $f$ is a meromorphic function with finitely many poles
  on some open subset $\Set{U}$ of the \Em{complex Fermi curve},
  then $\sigma^{\ast}f$ is a meromorphic function
  with finitely many poles on $\sigma(\Set{U})$.
  The associated operator--valued functions satisfy
  for all $\xx{p}\in\mathbb{C}$:
  $$\left(\begin{smallmatrix}
  -\yy{\kappa}_2+\sqrt{-1}\yy{\kappa}_1 & 0\\
  0 & \yy{\kappa}_2+\sqrt{-1}\yy{\kappa}_1
  \end{smallmatrix}\right)
  \comp\triv{\Op{A}}_{\sigma^{\ast}f}^{t}(-\xx{p})\comp
  \left(\begin{smallmatrix}
  -\yy{\kappa}_2+\sqrt{-1}\yy{\kappa}_1 & 0\\
  0 & \yy{\kappa}_2+\sqrt{-1}\yy{\kappa}_1
  \end{smallmatrix}\right)^{-1} =\triv{\Op{A}}_{f}(\xx{p}),$$
  and the associated variation
  $(\var V_{\sigma^{\ast}f},
  \var W_{\sigma^{\ast}f})$ is equal to
  $(-\var W_{f},-\var V_{f})$.
\item[(ii)] If the potentials are real and $f$ is a
  meromorphic function with finitely many poles on some open subset
  $\Set{U}$ of the \Em{complex Fermi curve}, then
  $\rho^{\ast}\Bar{f}$ is a
  meromorphic function with finitely many poles on
  $\rho(\Set{U})$. For all $\xx{p}\in\mathbb{C}$
  the associated operator--valued functions satisfy
  $$\left(\begin{smallmatrix}
  0 & \yy{\kappa}_2+\sqrt{-1}\yy{\kappa}_1\\
  \yy{\kappa}_2-\sqrt{-1}\yy{\kappa}_1 & 0
  \end{smallmatrix}\right)
  \comp\triv{\Op{A}}_{\rho^{\ast}\Bar{f}}^{\ast}(\Bar{\xx{p}})\comp
  \left(\begin{smallmatrix}
  0 & \yy{\kappa}_2+\sqrt{-1}\yy{\kappa}_1\\
  \yy{\kappa}_2-\sqrt{-1}\yy{\kappa}_1 & 0
  \end{smallmatrix}\right)^{-1} =\triv{\Op{A}}_{f}(\xx{p}).$$
  The associated variation
  $(\var V_{\rho^{\ast}\Bar{f}},
  \var W_{\rho^{\ast}\Bar{f}})$ is
  equal to $(-\var \Bar{V}_{f},-\var \Bar{W}_{f})$.
\item[(iii)] If the pair of potentials is of the form $(U,\Bar{U})$
  and $f$ is a
  meromorphic function with finitely many poles on some open subset
  $\Set{U}$ of the \Em{complex Fermi curve}, then
  $\eta^{\ast}\Bar{f}$ is a
  meromorphic function with finitely many poles on
  $\eta(\Set{U})$. The associated operator--valued functions satisfy
  for all $\xx{p}\in\mathbb{C}$:
  $$\Op{J}\comp\triv{\Op{A}}_{\eta^{\ast}\Bar{f}}^{\ast}
  (-\Bar{\xx{p}})\comp\Op{J}^{-1}
  =\triv{\Op{A}}_{f}(\xx{p}),$$
  and the associated variation
  $(\var V_{\eta^{\ast}\Bar{f}},
  \var W_{\eta^{\ast}\Bar{f}})$ is
  equal to $(\var \Bar{W}_{f},\var \Bar{V}_{f})$
\item[(iv)] The constant operator
  $\triv{\Op{A}}_{f}(\xx{p})=\left(\begin{smallmatrix}
  \unity & 0\\
  0 & -\unity
  \end{smallmatrix}\right)$ fulfills condition~(i)--(ii) of
  Lemma~\ref{residue}. The corresponding variation is equal to
  $\var V_{f}=-2V$ and $\var W_{f}=2W$.
\qed
\end{description}
\end{Lemma}

The corresponding function $f$ in (iv) would be a function,
which is in some neighbourhood of $\infty^-$ equal to $1$
and on some neighbourhood of $\infty^+$ equal to $-1$.
Since the two infinities should be considered as one double point,
this function is not holomorphic (compare with \cite{Sch}).
Moreover, the forms $\Omega_{V,W}(\var V,\var W)$ have in general
two simple poles at the two infinities, and the residue is
proportional to the  variation of the functional $\int VWd^2x$.

The foregoing Lemmata have an application,
which will be an essential tool
in the classification of relative minimizers.
By the residue formula
$\sum\limits_{\text{poles of }f\omega}\res(f\omega)=0$
on a compact Riemann surface every meromorphic function $f$
induces a relation of the values of a holomorphic (meromorphic)
form $\omega$ at the poles of $f$ (and the poles of $\omega$).
Due to S\'{e}rre duality \cite[\S17.]{Fo} all relations on the
values of all holomorphic forms are of this form.
We shall see now that the same is true for general
\Em{complex Fermi curves}.

\begin{Proposition}\label{global meromorphic function}
Let $\fermi(V,W)$ be a \Em{complex Fermi curve} with the
property that for all variations $\var V$ and $\var W$, whose
regular forms $\Omega_{V,W}(\var V,\var W)$
vanish at finitely many points
$y_1,\ldots,y_L\in\fermi(V,W)/\lattice\dual$,
the corresponding variation
$4\int\limits_{\torus}
\left(V(x)\var W(x)+\var V(x)W(x)\right)d^2x$ of the
\Em{first integral} vanishes. Then this \Em{complex Fermi curve}
$\fermi(V,W)/\lattice\dual$ can be compactified to a compact
Riemann surface by adding the points $\infty^{\pm}$
(compare with condition \Em{Finite type potentials}~(i)).
Moreover, this compactified \Em{complex Fermi curve}
has a global meromorphic function whose poles are contained in
\{$y_1,\ldots,y_l\}$ and which takes the values $\pm 1$ at $\infty^{\pm}$.
\end{Proposition}

\begin{proof}
The form $\Omega_{V,W}(\var V,\var W)$ vanishes at $y_l$,
if and only if the residue of this form multiplied with a
meromorphic function $f_l$ on a neighbourhood $\Set{U}$ of $y_l$
vanishes, where $f$ has a first order pole at $y_l$ and no other poles.
Therefore, for all variations $(\var V,\var W)$, whose residues
$$\res\limits_{y_1}\left(f_1\Omega_{V,W}(\var V,\var W)\right)
,\ldots,
\res\limits_{y_L}\left(f_L\Omega_{V,W}(\var V,\var W)\right)$$
vanish, the variation of the \Em{first integral} has to vanish.
Due to Lemma~\ref{residue}~(iv) the former residues vanish,
if and only if the symplectic form of $(\var V,\var W)$
and the variations
$$(\var V_{f_1},\var W_{f_1}),\ldots,(\var V_{f_L},\var W_{f_L})$$
vanish. Since the variation of the \Em{first integral} is equal to
four times the symplectic form of $(V,-W)$ and $(\var V,\var W)$,
we conclude that $(-V,W)$ is a complex linear combination of
$(V_{f_1},W_{f_1}),\ldots,(V_{f_L},W_{f_L})$. Due to Lemma~\ref{residue}~(ii)
this implies that the sum of the operator mentioned in
Lemma~\ref{compatibility}~(iv) and a complex linear combination
of the operators
$$\triv{\Op{A}}_{f_1}(\cdot),\ldots,\triv{\Op{A}}_{f_L}(\cdot)$$
commutes with $\triv{\Op{D}}_{\xx{\gamma},\yy{\gamma}}(\cdot)$.
The eigenvalues of this family of operators commuting with
$\triv{\Op{D}}_{\xx{\gamma},\yy{\gamma}}(\cdot)$
defines a meromorphic function on $\fermi(V,W)$.
Due to Lemma~\ref{residue}~(i) this function
induces a meromorphic function on $\fermi(V,W)/\lattice\dual$.
By construction this function solves the first Cousin Problem
(see e.\ g.\ \cite[Chapter~5 2.1.]{Ah} and \cite[Chapter~III. \S3.3.]{GPR})
of the data given by the corresponding complex linear combination of
$f_1,\ldots,f_L$. Moreover, due to Theorem~\ref{asymptotic analysis 1}
Lemma~\ref{residue}~(iii) and Lemma~\ref{compatibility}~(iv)
this function takes values in $\{z\in\mathbb{C}\mid|z+1|<\varepsilon\}$
on any open set of the form $\Set{V}^-_{\varepsilon,\delta}$
($\varepsilon$ arbitrary) and values in
$\{z\in\mathbb{C}\mid|z-1|<\varepsilon\}$ on any set of the form
$\Set{V}^+_{\varepsilon,\delta}$ ($\varepsilon$ arbitrary).
If $\varepsilon<1$, then the normalizations of the corresponding handles
described in Theorem~\ref{asymptotic analysis 1}
have two connected components. Therefore, the addition of the points
$\infty^{\pm}$ compactifies the \Em{complex Fermi curve}
$\fermi(V,W)/\lattice\dual$ to a compact Riemann surface
of finite genus. The solution of the first Cousin problem yields a
meromorphic function on this compact Riemann surface
whose poles are contained in \{$y_1,\ldots,y_l\}$
and which takes at $\infty^{\pm}$ the values $\pm 1$.
\end{proof}

\begin{Remark}\label{more general conditions}
The proof carries over to more general situations,
where the vanishing of $\Omega_{V,W}(\var V,\var W)$ at
$y_1,\ldots,y_L$ with orders $m_1-1,\ldots,m_L-1$ implies that
the corresponding variation of the \Em{first integral} vanishes.
Furthermore, due to the Separating hyperplane theorem
\cite[Theorem~V.4]{RS1},
the conclusion of the proposition holds as well,
if for variations $(\var V,\var W)$, whose forms
$\Omega_{V,W}(\var V,\var W)$ takes values in a convex sub--cone of
$\mathbb{C}$ at the points $y_1,\ldots,y_L$,
the corresponding variation of the \Em{first integral} takes values
in the complement of that convex sub--cone.
Finally, by construction of the variations $(\var V_f,\var W_f)$
in the proof of Lemma~\ref{residue}, the conclusion of the proof of
this proposition $(-V,W)$ being a linear combination of the variations
$(\var V_{f_1},\var W_{f_1}),\ldots,(\var V_{f_L},\var W_{f_L})$
implies that these potentials are smooth.
In fact, due to Remark~\ref{sobolev embedding} the resolvent of
the free Dirac operator is a bounded operator from the Sobolev spaces
$\sobolev{m-1,p}(\torus)\times\sobolev{m-1,p}(\torus)$ onto
$\sobolev{m,p}(\torus)\times \sobolev{m,p}(\torus)$
for all $m\in\mathbb{N}$ and $1<p<\infty$.
Hence, due to the proof of Theorem~\ref{meromorph},
for all $m\in\mathbb{N}$ the eigenfunctions
of potentials in the Sobolev spaces
$\bigcap\limits_{q<\infty}\sobolev{m-1,q}(\torus)$
belong to the Sobolev spaces
$\bigcap\limits_{q<\infty}\sobolev{m,q}(\torus)$.
We conclude inductively that the potential of a critical point
is indeed smooth.
If the pair of potentials is of the form $(U,\Bar{U})$,
then due to Lemma~\ref{singularity} $U$ has to be of \Em{finite type}.
Consequently, due to the obvious generalization of
Lemma~\ref{existence of potentials} to compact
one--dimensional complex spaces $\Spa{Y}$
(instead of compact smooth Riemann surfaces) $U$ is actually analytic.
\end{Remark}

We shall see that this proposition applies to critical points of
the Willmore functional, and yields a function theoretic
characterization of these critical points.
The corresponding global meromorphic function induces a holomorphic
mapping from the spectral curve onto $\mathbb{P}^1$.
It is well known \cite{DKN}
that if the preimage of one element of
$\mathbb{P}^1$ is contained in $\{\infty^-,\infty^+\}$,
then the corresponding potentials solve a differential equation,
and belong to a reduction of the integrable system.
The known descriptions of Willmore surfaces
\cite{Ba,BB,Bo} and of conformal minimal immersion into $S^3$
\cite{Hi,PS} are exactly of this form.

\subsection{Deformations of \Em{complex Fermi curves}}
\label{subsection deformation}

In this section we shall prove the smoothness of some
subsets of the moduli space\index{moduli space}.
This result will become the essential tool
in the investigation of relative minimizers
of the transformed variational problem of the Willmore functional.
Due to condition \Em{Quasi--momenta}~(i)
the complex space $\Spa{Y}\setminus\{\infty^-,\infty^+\}$
is a one-sheeted covering of $\fermi(\Spa{Y},\infty^-,\infty^+,k)$.
The elements of these subsets are locally isomorphic to
$\fermi(\Spa{Y},\infty^-,\infty^+,k)$ and therefore
local complete intersections of $\mathbb{C}^2$
\cite[Chapter~VI. \S2.1.]{GPR}
with fixed arithmetic genus.
To make this precise we introduce another condition
on the data $(\Spa{Y},\infty^-,\infty^+,k)$.
\begin{description}
\item[Quasi--momenta (v)]
  \index{condition!quasi--momenta  (v)}
  \index{quasi--momenta!condition $\sim$  (v)}
  There exists an open covering of
  $\Spa{Y}\setminus\{\infty^-,\infty^+\}$ with the property that the
  restriction of any branch of $k$ to any neighbourhood of this
  covering is an isomorphism of this neighbourhood onto some
  one--dimensional complex subspace of $\mathbb{C}^2$.
\end{description}
Let $\moduli_{g,\lattice}$
\index{moduli space!$\moduli_{g,\lattice}$}
be the set of all data fulfilling
conditions \Em{Quasi--momenta}~(i)--(ii) and \Em{Quasi--momenta}~(v) with the
property, that the
arithmetic genus of $\Spa{Y}$ is equal to $g$.
For all data $(\Spa{Y}',\infty^-,\infty^+,k)$,
where $\Spa{Y}'$ are smooth compact Riemann surfaces, there exist some data
$(\Spa{Y},\infty^-,\infty^+,k)$ in some $\moduli_{g,\lattice}$, where the
normalization of $\Spa{Y}$ is equal to $\Spa{Y}'$ and
the pullback of the function $k$ on $\Spa{Y}$
under the normalization map is equal to the corresponding
function of $k$ on $\Spa{Y}'$. Hence the union
$\bigcup\limits_{g\in\mathbb{N}_0}\moduli_{g,\lattice}$ contains
all \Em{complex Fermi curves} of finite \Em{geometric genus}.

\begin{Proposition}\label{smooth moduli}
\index{deformation!of a complex Fermi curve}
The set $\moduli_{g,\lattice}$ is a
  $g+1$--dimensional complex manifold.
For all elements $(\Spa{Y},\infty^-,\infty^+,k)$ the tangent space
\index{moduli space!tangent space of the $\sim$}
\index{tangent!space of the moduli space}
of this manifold at this point is isomorphic to the space of all
regular forms of the complex space $\Tilde{\Spa{Y}}$, which is obtained from
$\Spa{Y}$ by identifying $\infty^-$ and $\infty^+$
to an ordinary double point. More precisely,
for any such regular form $\omega$ there exists
a one--dimensional complex family of data, such that
$\omega=\Dot{\yy{p}}d\xx{p}$.
\end{Proposition}

\begin{Remark}\label{disconnected moduli}
Here we assume that $\Spa{Y}$ is connected.
More precisely, we assume that if $\Spa{Y}$ has a
disconnected normalization, then $\Spa{Y}$ has
a singularity containing elements of both components.
Otherwise the regular forms of $\Tilde{\Spa{Y}}$ coincide with the
regular forms of $\Spa{Y}$, and $\moduli_{g,\lattice}$ is a
$g$--dimensional manifold. Moreover, in this case the
\Em{first integral} is locally constant on this manifold.
\end{Remark}

Here we have used the notation introduced
in \De{Fundamental domain}~\ref{fundamental domain}. Furthermore, we
make the convention that the function $\xx{p}$ does not depend on the
deformation.

\begin{proof} For given data
$(\Spa{Y},\infty^-,\infty^+,k)\in\moduli_{g,\lattice}$ we shall first
construct an open neighbourhood in some complex manifold, whose
elements correspond to deformations of $\Spa{Y}$. More precisely, we
consider all deformations of compact one--dimensional complex spaces of
arithmetic genus $g$, on which the multi--valued function
$\xx{p}=g(\xx{\gamma},k)$ exists with the analogous properties
\Em{Quasi--momenta}~(i)--(ii), and which is a
local complete intersection in $\mathbb{C}^2$
\cite[Chapter~VI. \S2.1.]{GPR}.
Here we carry over the methods of \cite{Ko} to the
case of complex spaces instead of complex manifolds. The tools are
essentially developed by Grauert (see e.\ g.\
\cite{Gr,BaSt,GrRe,GK,GPR})\index{deformation!of a complex space}.
Afterwards we shall extend any regular form
of $\Tilde{\Spa{Y}}$ to a holomorphic family of regular forms on the
corresponding family of complex spaces. Then we shall show that each
of these holomorphic families of regular forms defines some flow
on the manifold, which preserves the subset corresponding to
elements of $\moduli_{g,\lattice}$. Finally, we show that the
regular forms on the family of deformed complex spaces
may be considered as a holomorphic sub--bundle
of the tangent bundle of the manifold of
deformations constructed before, and that this sub--bundles is the
tangent bundle along the submanifold corresponding to elements of
$\moduli_{g,\lattice}$.
Along these lines the proof is achieved in ten steps.

\noindent
{\bf 1.} Since the arithmetic genus is finite,
the union of the singular points of $\Spa{Y}$ and the zeroes of
$d\xx{p}$ is a finite subset of $\Spa{Y}$. Now we choose a
finite covering $\Set{U}_0,\ldots,\Set{U}_L$ of $\Spa{Y}$, such that
$\Set{U}_1,\ldots,\Set{U}_L$
are disjoint open neighbourhoods, each of which contains either
one singular point or one zero of $d\xx{p}$.
Without loss of generality we may assume that the function
$\xx{p}$ maps all these open subsets onto Weierstra{\ss} coverings
over open subsets $\Set{V}_1,\ldots,\Set{V}_L$ of $\xx{p}\in\mathbb{C}$.
We in addition assume that the subset $\Set{U}_0$ is contained
in the complement of some neighbourhood of the union of
all singular points and all zeroes of $d\xx{p}$.
Finally, the two marked points $\infty^-$ and $\infty^+$ are 
assumed to be elements of the complement of the closure of
$\Set{U}_1\cup\ldots\cup \Set{U}_L$ in $\Set{U}_0$.
We conclude that the restriction of
the function $\xx{p}$ to $\Set{U}_0$ is locally biholomorphic,
and therefore
$\Set{U}_0$ may be identified with some union of open subsets of
$\xx{p}\in\mathbb{P}^1$. Due to condition
\Em{Quasi--momenta}~(v) any $\Set{U}_l$ of the subsets
$\Set{U}_1,\ldots,\Set{U}_L$ is
isomorphic to the complex spaces defined by some equation
$R(\xx{p},\yy{p})=0$
over some open subset $\Set{V}_l$ of $\mathbb{C}^2$.
Moreover, we may choose this function $R(\xx{p},\yy{p})$
to be a polynomial with respect to $\yy{p}$,
whose degree $d$ is equal to the number of sheets of
this Weierstra{\ss} covering with respect to the function $\xx{p}$,
and whose highest coefficient is $1$.
It would be more accurate to decorate the functions $R$ and their
degrees $d$ with an index $l=1,\ldots,L$. We omit these indices,
and hope that from the context it will be clear, which objects are
associated to which member of the covering.

\noindent
{\bf 2.} Let us simplify these local descriptions to polynomial equations
$Q(\xx{p},\Breve{p})=0$ together with an equation
$\yy{p}=P(\xx{p},\Breve{p})$,
which is polynomial with respect to $\Breve{p}$ of degree $d-1$,
and whose coefficients are holomorphic functions depending on
$\xx{p}\in \Set{V}_l$.
If we choose the neighbourhoods $\Set{U}_1,\ldots,\Set{U}_L$
small enough, then there always exists such a simplification,
where the coefficients of $Q(\xx{p},\Breve{p})$ are equal to the
corresponding lowest non--vanishing Taylor coefficients of
$R(\xx{p},\Breve{p})$. More precisely, if the preimage
of the corresponding zero of $d\xx{p}$ in the normalization
is a single point, then $Q(\xx{p},\Breve{p})$ can be chosen
to be of the form $\Breve{p}^d-\xx{p}^m$, where $c\xx{p}^m$
is the lowest non--vanishing coefficient of $R$ of this form.
In this case $d$ and $m$ have to be co--prime,
and we define the degree of $\xx{p}$ as $d$
and the degree of $\Breve{p}$ (and $\yy{p}$) as $m$.
In general the
\De{Local contribution to the arithmetic genus}~\ref{local contribution}
of the equation $R(\xx{p},\yy{p})=0$
at the corresponding zero of $d\xx{p}$ has to be the same as the
corresponding
\De{Local contribution to the arithmetic genus}~\ref{local contribution}
of the equation $Q(\xx{p},\Breve{p})=0$.

\noindent
{\bf 3.} If $\Dot{Q}(\xx{p},\Breve{p})$
is any polynomial with respect to $\xx{p}$ and $\Breve{p}$, whose
degree with respect to $\Breve{p}$ is smaller than $d$, the equation
$Q(\xx{p},\Breve{p},t)=Q(\xx{p},\Breve{p})+t\Dot{Q}(\xx{p},\Breve{p})$
defines a deformation of the complex
space $\Set{U}_l$. Also for small $|t|$ the function $\xx{p}$ of this
deformation defines a Weierstra{\ss} covering over $\Set{V}_l$,
whose restriction to the corresponding deformed sets
$\Set{U}_l\cap \Set{U}_0$ is biholomorphic to the undeformed set.
Hence this deformation defines a deformation of the entire
compact complex space. If we consider $\Breve{p}$ as a function
depending on $\xx{p}$ and $t$, the derivative
$\partial \Breve{p}/\partial t$ is
equal to $-\Dot{Q}(\xx{p},\Breve{p})/
(\partial Q(\xx{p},\Breve{p},t)/\partial \Breve{p})$.
If for small $|t|$ this quotient is a holomorphic function
on the deformed space,
then we conclude that at least for small values of $|t|$ the
deformed function
$\Breve{p}$ may be expressed by the undeformed functions $\xx{p}$ and
$\Breve{p}$ and
vice versa. Therefore, the family of complex spaces are all isomorphic.
Consequently we choose finitely many monomials $Q_1,\ldots,Q_r$
with respect to $\xx{p}$ and $\Breve{p}$,
whose degree with respect to $\Breve{p}$ is smaller than $d$,
and which induce a basis of the space of global sections of
the sheaf $\Sh{O}/\Sh{O}
(\partial Q(\xx{p},\Breve{p})/\partial \Breve{p})$ on
the subset $\Set{U}_l$. Let $\Tilde{\moduli}_l\subset\mathbb{C}^r$ be
some small open neighbourhood of $0$, and let $\Spa{Z}_l$ be the complex
subspace of $(\xx{p},\Breve{p},t_1,\ldots,t_r)\in
\Set{V}_l\times\mathbb{C}\times\Tilde{\moduli}_l$ defined by the
equation
$$Q(\xx{p},\Breve{p},t_1,\ldots,t_r)=
Q(\xx{p},\Breve{p})
+t_1Q_1(\xx{p},\Breve{p})+\ldots+t_rQ_r(\xx{p},\Breve{p})=0.$$
There exists a natural map $f_l:\Spa{Z}_l\rightarrow
\Tilde{\moduli}_l$, which is induced by the projection
$\Set{V}_l\times\mathbb{C}\times\Tilde{\moduli}_l\rightarrow
\Tilde{\moduli}_l$. For all elements
$t_1\ldots,t_r\in\Tilde{\moduli}_l$
the complex subspaces defined by the pullback under $f_l$ under the
unique maximal ideal corresponding to $t_1,\ldots,t_r$ is obviously the
deformation of $\Set{U}_l$ defined by the equation
$Q(\xx{p},\Breve{p},t_1,\ldots,t_r)=0.$
Hence $\Spa{Z}_l$ is some ``universal family'' of deformations of $\Set{U}_l$.
Now we claim that the map $f_l$ is flat. In fact, $f_l$ is the
composition of some Weierstra{\ss} map
$\Spa{Z}_l\rightarrow \Set{V}_l\times\Tilde{\moduli}_l$
and the projection 
$\Set{V}_l\times\Tilde{\moduli}_l\rightarrow\Tilde{\moduli}_l$. The
first map is flat due to \cite[Chapter~II Proposition~2.10]{GPR} and
the Weierstra{\ss} Isomorphism \cite[Chapter~2. \S4.2.]{GrRe}. The
second map is flat due to \cite[Chapter~II Corollary~2.7]{GPR}. Hence
\cite[Chapter~II Proposition~2.6]{GPR} shows the claim. If
$\Tilde{\moduli}_l$ is chosen small enough, the deformations of
the subsets $\Set{U}_l\cap \Set{U}_0$ do not contain any zero
of the function
$\partial Q(\xx{p},\Breve{p},t_1,\ldots,t_r)/\partial \Breve{p}$.
In this case we may
glue $\Spa{Z}_l$ with the complex space $\Set{U}_0\times\Tilde{\moduli}_l$
along the deformations of the subset $\Set{U}_l\cap \Set{U}_0$.
If we repeat this for all $l=1,\ldots,L$ we obtain a universal family
$f:\Spa{Z}\rightarrow \Tilde{\moduli}=
\Tilde{\moduli}_1\times\ldots\times
\Tilde{\moduli}_L$,
of deformations of the compact complex space $\Spa{Y}$.
We conclude that the map $f$ is proper and flat.

\noindent
{\bf 4.} Now we shall prove that for all
$l=1,\ldots,L$ the corresponding polynomials $Q_1,\ldots,Q_r$ form a
basis of all sections of the sheaf
$\Sh{O}/\Sh{O}
(\partial Q(\xx{p},\Breve{p},t_1,\ldots,t_r)/\partial \Breve{p})$
on $\Spa{Z}_l$. The support of this sheaf is the complex subspace
$\Check{\Spa{Z}}_l$ of $\Spa{Z}_l$ defined by the equation
$\partial Q(\xx{p},\Breve{p},t_1,\ldots,t_r)/\partial \Breve{p}=0$.
The union of these
subspaces is a subspace $\Check{\Spa{Z}}$ of $\Spa{Z}$, and the map $f$ induces a
map $\Check{f}:\Check{\Spa{Z}}\rightarrow \Tilde{\moduli}$. This map is
again proper and flat. In fact, it is a finite map, and the dimension
of the stalk of the preimage of any point of $\Tilde{\moduli}$
considered as a $\mathbb{C}$--vector space is equal to the number of
zeroes of the form $d\xx{p}$ on the corresponding deformation of $\Spa{Y}$.
Due to the
Semi-Continuity theorem \cite[Chapter~III Theorem~4.7]{GPR}
the arithmetic genus of all these deformations is equal to
$g$. Moreover, since all deformations are
local complete intersections in $\mathbb{C}^2$
\cite[Chapter~VI. \S2.1.]{GPR},
Lemma~\ref{dualizing sheaf} shows that the sheaf of
regular forms are locally free. Therefore,
the degree of the canonical divisor is equal to $2g-2$
\cite[Chapter~IV. \S3.11]{Se}.
Since $d\xx{p}$ has two poles at the
non--singular points $\infty^-$ and $\infty^+$, the number of zeroes of
$d\xx{p}$ is equal to $2g+2$ and therefore constant. Again
\cite[Chapter~II Proposition~2.10]{GPR}
shows that $\Check{f}$ is flat. Now the
Semi-Continuity theorem \cite[Chapter~III Theorem~4.7]{GPR}
implies that the direct image sheaf of the
structure sheaf of $\Check{\Spa{Z}}$ under $\Check{f}$ defines
some holomorphic vector bundles on $\Tilde{\moduli}$ of
rank $2g+2$. Hence for any choice of sections,
which form a basis in the fibre over
$0\in\Tilde{\moduli}$, these sections form a basis of all fibres
in some neighbourhood. If $\Tilde{\moduli}_l$ is chosen small
enough, then the corresponding elements $Q_1,\ldots,Q_r$ form a basis
of all global sections of the sheaf
$\Sh{O}/\Sh{O}
(\partial Q(\xx{p},\Breve{p},t_1,\ldots,t_r)/\partial \Breve{p})$
on $\Spa{Z}_l$.

\noindent
{\bf 5.} Due to Lemma~\ref{dualizing sheaf} the sheaf of regular forms on the
family of deformations of $\Spa{Y}$ is locally free. Therefore,
\cite[Chapter~II. Proposition~2.6]{GPR} implies that this sheaf is
also $f$--flat, and the
Semi-Continuity theorem \cite[Chapter~III Theorem~4.7]{GPR}
shows that the direct image of the sheaf of sections of the
regular forms of the deformations of $\Spa{Y}$ under $f$ defines
some vector bundle on $\Tilde{\moduli}$ of rank $g$.
The same argumentation carries over to the direct image of the sheaf
of regular forms of the corresponding deformations of
$\Tilde{\Spa{Y}}$, which is defined by identifying the
non--singular points $\infty^-$ and $\infty^+$ to an
ordinary double point. If $\Spa{Y}$ is connected
(i.\ e.\ either the normalization of $\Spa{Y}$ is connected
or the two connected components of the normalization of $\Spa{Y}$
are connected by some singularity), then this yields a
holomorphic vector bundle of rank $g+1$ on $\Tilde{\moduli}$,
otherwise a holomorphic vector bundle of rank $g$
(compare with Remark~\ref{disconnected moduli}).
We conclude that any basis of regular forms of $\Tilde{\Spa{Y}}$,
extends to a holomorphic family of regular forms on the
preimage of some open neighbourhood of $0\in\Tilde{\moduli}$.

\noindent
{\bf 6.} If for the deformations of the sets $\Set{U}_l$ we make
the ansatz $\yy{p}=P(\xx{p},\Breve{p})$
and conceive the functions $\yy{p}$ and $\Breve{p}$
as functions depending on $\xx{p}$ and the moduli,
then the differential equation $\Dot{\yy{p}}d\xx{p}=\omega$,
where $\omega$ is a section
of the vector bundle on $\Tilde{M}$ constructed above,
transforms into the differential equation
\begin{eqnarray*}
\frac{\omega}{\partial P(\xx{p},\Breve{p})/\partial \Breve{p}}&=&
\frac{\Dot{\yy{p}}d\xx{p}}
       {\partial P(\xx{p},\Breve{p})/\partial \Breve{p}}=
\left(\frac{\Dot{P}(\xx{p},\Breve{p})}
       {\partial P(\xx{p},\Breve{p})/\partial \Breve{p}}
+\Dot{\Breve{p}}\right)d\xx{p}\\
&=&\left(\frac{\Dot{P}(\xx{p},\Breve{p})}
           {\partial P(\xx{p},\Breve{p})/\partial \Breve{p}}-
\frac{\Dot{Q}(\xx{p},\Breve{p})}
     {\partial Q(\xx{p},\Breve{p})/\partial \Breve{p}}\right)d\xx{p}.
\end{eqnarray*}
The partial derivative
$\partial P(\xx{p},\Breve{p})/\partial\Breve{p}$
of the initial $P(\xx{p},\Breve{p})$ has a holomorphic inverse
on $\Set{U}_l$. By reduction of $\Tilde{\moduli}$ we may achieve
that the same is true for all deformed $P(\xx{p},\Breve{p})$.
Due to Lemma~\ref{dualizing sheaf} any regular form on some deformed
$\Spa{Y}$ of the family constructed above may be written on each
$\Set{U}_l$ of the deformed subsets
$\Set{U}_1,\ldots,\Set{U}_L$ as the quotient of some
holomorphic function $f(\xx{p},\Breve{p})$ divided by
$\partial Q(\xx{p},\Breve{p},t_1,\ldots,t_r)/\partial \Breve{p}$.
Therefore, the same is true for
$\frac{\omega}{\partial P(\xx{p},\Breve{p})/\partial\Breve{p}}$.
Since the polynomials
$Q_1,\ldots,Q_r$ forms a basis of the sheaf
$\Sh{O}/\Sh{O}
(\partial Q(\xx{p},\Breve{p},t_1,\ldots,t_r)/\partial \Breve{p})$, 
any such regular form may be written locally as
$$\frac{\omega}{\partial P(\xx{p},\Breve{p})/\partial\Breve{p}}=
\frac{c_1Q_1(\xx{p},\Breve{p})+\ldots+c_rQ_r(\xx{p},\Breve{p})}
{\partial Q(\xx{p},\Breve{p},t_1,\ldots,t_r)/\partial \Breve{p}}
d\xx{p}+g(\xx{p},\Breve{p})d\xx{p},$$
where $c_1,\ldots,c_r$ are uniquely defined complex numbers (depending
holomorphically on the moduli in $\Tilde{\moduli}$)
and $g(\xx{p},\Breve{p})$ is a unique holomorphic function
on $\Set{U}_l$ (also depending holomorphically on the moduli).
Hence we obtain the differential equation
$$\Dot{t}_1=-c_1,\ldots,\Dot{t}_r=-c_r \text{ and }
\Dot{P}(\xx{p},\Breve{p})=g(\xx{p},\Breve{p})
\frac{\partial P(\xx{p},\Breve{p})}{\partial \Breve{p}}.$$

\noindent
{\bf 7.} With the help of the Weierstra{\ss} Isomorphism
\cite[Chapter~2.\ \S4.2.]{GrRe} the holomorphic functions
$P(\xx{p},\Breve{p})$, $f(\xx{p},\Breve{p})$
and $g(\xx{p},\Breve{p})$ may be uniquely
represented as polynomials with respect to $\Breve{p}$ of degree $d-1$
(where $d$ denotes the number of sheets of the Weierstra{\ss} covering
$\Set{U}_l\rightarrow \Set{V}_l$), whose coefficients are holomorphic functions
depending on $\xx{p}$. It is quite easy to see
that with this choice the holomorphic function
$\frac{f(\xx{p},\Breve{p})}
      {\partial P(\xx{p},\Breve{p})/\partial \Breve{p}}$
is also of this form, and the coefficients of this polynomial
of degree $d-1$ with respect to $\Breve{p}$ are rational functions
of the corresponding coefficient of $P(\xx{p},\Breve{p})$
and $f(\xx{p},\Breve{p})$.
We conclude that the numbers $c_1,\ldots,c_r$
and the coefficient of the polynomial $g(\xx{p},\Breve{p})$
with respect to $\Breve{p}$ depend holomorphically
on the coefficients of the polynomials $P(\xx{p},\Breve{p})$
and $f(\xx{p},\Breve{p})$ with respect to $\Breve{p}$.
On $\Set{U}_0$ we have the ordinary differential equation
$\Dot{P}(\xx{p})=\Dot{\yy{p}}=\omega/d\xx{p}$.
We conclude that the section $\omega$ defines
a holomorphic vector field on the product of
$\Tilde{M}$ with the space of polynomials
$P(\xx{p},\Breve{p})$ on $\Set{U}_0\times\ldots\times\Set{U}_L$
of the corresponding degrees $d-1$ with respect to $\Breve{p}$,
whose coefficients depend holomorphically on
$\xx{p}\in \Set{V}_0\times\ldots\times\Set{V}_L$.
Since these vector fields are holomorphic they may be integrated to
some holomorphic flow at least in some neighbourhood of
$0\in\Tilde{\moduli}$ \cite[Chapter~IV \S1.]{La},
which have to be chosen small enough in order to guarantee the
assumption that $\partial P(\xx{p},\Breve{p})/\partial \Breve{p}$
has a holomorphic inverse on $\Set{U}_l$ for $l=1,\ldots,L$.

\noindent
{\bf 8.} This implies that in some neighbourhood of
$0\in\Tilde{\moduli}$ on the whole family of the integrated flow
two functions $\xx{p}$ and $\yy{p}+P(\xx{p},\Breve{p})$
are defined,  which may be combined to a function $k$ with
$\xx{p}=g(\xx{\gamma},k)$ and $\yy{p}=g(\yy{\gamma},k)$.
This function $k$ fulfills conditions \Em{Quasi--momenta}~(i)--(ii).
Moreover, due to the assumption that
$\partial P(\xx{p},\Breve{p})/\partial \Breve{p}$
has a holomorphic inverse on $\Set{U}_l$, the map
$(\xx{p},\Breve{p})\mapsto(\xx{p},\yy{p})$ is locally
a biholomorphic map. Therefore, also condition \Em{Quasi--momenta}~(v) is
fulfilled for small moduli. This shows that any holomorphic family of
regular forms on the family of deformations of $\Spa{Y}$ constructed above
define a vector field on the manifold of moduli, whose flow
preserves the subspace of deformations contained in
$\moduli_{g,\lattice}$.

\noindent
{\bf 9.} Now we shall prove that the subset of $\Tilde{\moduli}$,
on whose elements besides the function $\xx{p}$
there exists another function $\yy{p}$ such that the
resulting function $k$ fulfills conditions
\Em{Quasi--momenta}~(i)--(ii), is a complex subspace of this complex
manifold. It is part of the Krichever construction
\cite{Kr1} that the space translations of these soliton equations
are induced by the action of some linear subgroup of the
Picard group. From this point of view the existence of the function
$k$ with the properties \Em{Quasi--momenta}~(i)--(ii) is equivalent to the
existence of two periods of a natural two--dimensional subgroup of
the Picard group associated to the choice of the two marked points
$\infty^-$ and $\infty^+$. This periodicity is equivalent to
the assumption that the two line bundles corresponding to the
translation by $\xx{\gamma}$ and $\yy{\gamma}$ are trivial.
Since the function $\xx{p}$ exists on the whole universal family
$\Tilde{\Spa{Z}}$, the line bundle corresponding to the translation by
$\xx{\gamma}$ is trivial for the whole family.
Now we construct the line bundle corresponding to the shift by
$\yy{\gamma}$. The corresponding global sections are of the form
$c\exp(2\pi\sqrt{-1}\yy{p})$, where $c$ is some complex number and
$\yy{p}=g(\yy{\gamma},k)$ together with $\xx{p}=g(\xx{\gamma},k)$
defines a holomorphic function $k$ on $\Spa{Y}$ fulfilling
conditions \Em{Quasi--momenta}~(i)--(ii).
Let $\Set{U}^+$ and $\Set{U}^-$ denote small open disjoint
neighbourhoods of $\infty^-$ and $\infty^+$,
which are contained in the complement of
$\Set{U}_1,\ldots,\Set{U}_L$.
Due to condition \Em{Quasi--momenta}~(i) some branch of
the function $\yy{p}$ is proportional to
$\frac{\yy{\gamma}_1-\sqrt{-1}\yy{\gamma}_2}
{\xx{\gamma}_1-\sqrt{-1}\xx{\gamma}_2}$
times $\xx{p}$ at $\infty^-$ and proportional to 
$\frac{\yy{\gamma}_1+\sqrt{-1}\yy{\gamma}_2}
{\xx{\gamma}_1+\sqrt{-1}\xx{\gamma}_2}$
times $\xx{p}$ at $\infty^+$ up to order
$\text{\bf{O}}(1/\xx{p})$. Hence the
line bundle, which describes the lifting of the translation by
$\yy{\gamma}$, may be defined by the trivial bundle on
$\Spa{Y}\setminus\{\infty^-,\infty^+\}$, $\Set{U}^+$ and $\Set{U}^-$,
which are glued by the function
$\exp\left(2\pi\sqrt{-1}\xx{p}
\frac{\yy{\gamma}_1+\sqrt{-1}\yy{\gamma}_2}
{\xx{\gamma}_1+\sqrt{-1}\xx{\gamma}_2}\right)$
on $(\Spa{Y}\setminus\{\infty^+\})\cap \Set{U}^+$ and by the function
$\exp\left(2\pi\sqrt{-1}\xx{p}
\frac{\yy{\gamma}_1-\sqrt{-1}\yy{\gamma}_2}
{\xx{\gamma}_1-\sqrt{-1}\xx{\gamma}_2}\right)$
on $(\Spa{Y}\setminus\{\infty^+\})\cap \Set{U}^-$. Due to the
Direct Image theorem \cite[Chapter~10. \S4.6.]{GrRe}
the direct image of the sheaf of sections of this line bundle under
$\Tilde{f}$ is a coherent sheaf on $\Tilde{\moduli}$.
The subset of this manifold, on which another
function $\yy{p}$ is defined and whose combination $k$ fulfills
conditions \Em{Quasi--momenta}~(i)--(ii) is exactly the support of this
direct image sheaf. In fact, any non--zero germ of this
direct image sheaf defines a global holomorphic function on
$\Spa{Y}\setminus\{\infty^-,\infty^+\}$ with
essential singularity at the two marked points. Since the degree of
the family of line bundles on the family of complex spaces is zero,
this global holomorphic function does not vanish, and the logarithm of
this function divided by some uniquely defined complex number yields a
function $\yy{p}$ with the desired properties. Since the support of
any coherent sheaf is an analytic set \cite[Annex. \S4.5.]{GrRe} this
establishes the claim.

\noindent
{\bf 10.} This last result implies that the subset of
$\Tilde{\moduli}$ corresponding to elements of $\moduli_{g,\lattice}$
contains an open dense subset of nonsingular points,
and this open dense subset forms a submanifold of
$\Tilde{\moduli}$. Moreover, the derivative of $\yy{p}$ with
respect to any vector field along this submanifold defines some
regular form on the corresponding family of complex subspaces.
We conclude that the dimension of this complex space is equal to
$g+1$, and that the tangent sheaf extends to
a locally free sheaf of rank $g+1$. Now
the same arguments as in the proof of
\cite[Chapter~II. \S1.4. Theorem~1.12.]{GPR}
(i.\ e. \cite[Chapater~II. \S1.4. Fact~1.13.]{GPR} combined with
\cite[Chapter~I. \S10.5. (10.20)]{GPR}) imply that this subset is smooth.
\end{proof}

For the proof of the existence of the minimizer
we shall only consider \Em{complex Fermi curves} of bounded genus.
Therefore, we should investigate the subspace of this moduli space
containing only \Em{complex Fermi curves}
of finite \Em{geometric genus}. Obviously, a given
\Em{complex Fermi curve} of \Em{geometric genus} $g$
corresponds to some element in $\moduli_{g,\lattice}$
if and only if the normalization does not contain zeroes
of both differentials $d\xx{p}$ and $d\yy{p}$.
In this case $\Spa{Y}$ may be chosen to be the normalization.
Otherwise there exists a unique representative element
$(\Spa{Y},\infty^-,\infty^+,k)\in\moduli_{g',\lattice}$ of minimal
arithmetic genus $g'$. In this case $\Spa{Y}$ is locally isomorphic to
the images of the normalization under $k$,
and the arithmetic genus is equal to the genus
of the normalization plus the
\De{Local contribution to the arithmetic genus}~\ref{local contribution}
of all these common zeros of $d\xx{p}$ and $d\yy{p}$,
which are called \Em{cusps}. Unfortunately, the subspace of
$\moduli_{g',\lattice}$ containing all \Em{complex Fermi curves}
of \Em{geometric genus} $g$,
has itself \Em{cusps} at \Em{complex Fermi curves},
whose normalization contains \Em{cusps}.

\begin{Proposition}\label{modified smooth moduli}
\index{moduli space!cusp of the $\sim$}
\index{cusp!of the moduli space}
\index{cusp}
For any \Em{complex Fermi curve}
$\fermi(\Spa{Y},\infty^-,\infty^+,k)\in\moduli_{g',\lattice}$,
whose normalization contains a \Em{cusp} and has genus $g<g'$,
and all meromorphic forms $\omega$ on the normalization of $\Spa{Y}$
with at most first--order poles at $\infty^{\pm}$ and no other poles,
there exists a one--dimensional complex subspace
of some neighbourhood of $\fermi$ in $\moduli_{g',\lattice}$,
which contains $\fermi$,
and which is contained in the set of the \Em{complex Fermi curves}
of \Em{geometric genus} $g$.
Moreover, this one--dimensional subspace
contains a \Em{cusp} at $\fermi$ in tangent direction
of the line corresponding to $\omega$.
More precisely, the preimage of $\fermi$
in the normalization of this one--dimensional subspace
contains a single point. Furthermore, with respect to
any local parameter of the normalization at this point
the lowest non--vanishing derivative is proportional to $\omega$,
considered as an element of the tangent space of
$\moduli_{g',\lattice}$ over $\fermi$.
\end{Proposition}

\begin{proof} This proposition is proven in nine steps.

\noindent
{\bf 1.} First we remark that due to
Proposition~\ref{smooth moduli}
the flow of a vector field on some neighbourhood
of some \Em{complex Fermi curve} without \Em{cusps},
which takes values in the subspace of the tangent space
corresponding to meromorphic forms on the normalization of the
corresponding \Em{complex Fermi curve} with at most
first--order poles at $\infty^{\pm}$ and no other poles,
preserves the genus of the normalization
of the corresponding \Em{complex Fermi curves}.

\noindent
{\bf 2.} We shall see that again we may transform a \Em{cusp}
of the \Em{complex Fermi curve}
into \Em{multiple points} (i.\ e.\ singularities, whose
preimage in the normalization does not contain
zeroes of both differentials $d\xx{p}$ and $d\yy{p}$)
by blowing up the moduli space.
We choose the unique representative $(\Spa{Y},\infty^-,\infty^+,k)$
of minimal arithmetic genus described above.
Without loss of generality we may choose
the local descriptions of the elements of $\Tilde{M}$
to be locally zero sets of equations
$\Breve{p}^d-\xx{p}^m=0$ with two co--prime natural numbers
$d$ and $m$. Since we are interested in deformations,
which preserve the \Em{geometric genus},
we change the description of $\Tilde{M}_l$ near the \Em{cusps}.
All these deformations correspond locally to those deformations
in $\Tilde{M}_l$, whose equations
$Q(\xx{p},\Breve{p},t_1,\ldots,t_r)=0$
describe algebraic curves birational isomorphic to $\mathbb{P}^1$.
In particular, the functions $\xx{p}$ and $\Breve{p}$
may be written as polynomials
$$\xx{p}=z^d+\Hat{t}_1z^{d-1}+\Hat{t}^2_2z^{d-2}
+\ldots+\Hat{t}^d_{d} \text{ and }
\Breve{p}= z^m+\Breve{t}_1z^{m-1}+\Breve{t}^2_2z^{m-2}
+\ldots+\Breve{t}^m_{n}.$$
The coefficients are written in terms of powers of the moduli
$(\Hat{t}_1,\ldots,\Hat{t}_d,\Breve{t}_1,\ldots,\Breve{t}_m)
\in\Breve{\moduli}_l$
in such a way, that the re-scaling $z\mapsto\lambda z$
acts projectively on the moduli.
For all these moduli there exists an unique polynomial with
$$Q(\xx{p},\Breve{p},
\Hat{t}_1,\ldots,\Hat{t}_d,\Breve{t}_1,\ldots,\Breve{t}_m)=0.$$
If we define the degree of $\xx{p}$ as $d$
and the degree of $\Breve{p}$ as $m$ this polynomial
has common degree $dm$, and the highest coefficient is equal to
$\Breve{p}^d-\xx{p}^m$.

\noindent
{\bf 3.} We should remark that the translations $z\mapsto z+z_0$
do not change $Q$, but change the moduli.
Furthermore, a finite groups acts on the moduli without changing
the coefficients of the polynomials $\xx{p}$ and $\Breve{p}$
with respect to $z$. For this reason the moduli space
should be defined as the quotient by the action of the translations
and the finite group on the moduli.
The orbits of this group action define one--dimensional curves
in the parameter space parameterized by $z_0$.
At some points these curves might have singularities
(i.\ e.\ the derivatives of all parameters with respect to $z_0$ vanish),
but obviously for all orbits this happens only
at finitely many values of $z_0$.
Therefore, we may choose some non--singular representative
and use the implicit function theorem
\cite[Supplementary material~V.5 Theorem~S.11]{RS1}
in order to obtain locally smooth parameters
of the quotient by this group action.
In the sequel we will ignore this ambiguity of the moduli
$(\Hat{t}_1,\ldots,\Hat{t}_d,\Breve{t}_1,\ldots,\Breve{t}_m)$
and assume that the real moduli space is locally
the quotient by this group action.

\noindent
{\bf 4.} If we glue the local normalizations of the deformed spaces,
parameterized again by $z$ on $\Set{U}_0$ in the natural way,
then the arguments in the proof of Proposition~\ref{smooth moduli}
carry over to the modified moduli space $\Breve{\moduli}$.
They show that the space of meromorphic forms on the normalizations,
which have at most first--order poles at $\infty^{\pm}$
and no other poles, define a holomorphic vector bundle
on $\Breve{\moduli}$ of rank $g+1$ or $g$
(compare with Remark~\ref{disconnected moduli}).
Therefore, all these forms on the normalization of the initial
\Em{complex Fermi curve} extend to holomorphic sections of this
vector bundle. Now we shall construct the deformations corresponding
to the sections of this vector bundle by blowing up the moduli space
$\Breve{\moduli}$ at the origin corresponding to the initial
\Em{complex Fermi curve}. The blowing up re-parameterizes
\index{blowing up!of the moduli space}
\index{moduli space!blowing up of the $\sim$}
locally the moduli as elements of the tautological bundle
over the corresponding projective space \cite[pp. 28-31]{Har}.
Locally  this blowing up is the product of the projective space with 
a one--dimensional space parameterizing
the fibre of the tautological bundle.
We should choose an element
$(\Hat{t}_1,\ldots,\Hat{t}_d,\Breve{t}_1,\ldots,\Breve{t}_m)
\in\mathbb{P}^{d+m-1}$ of the projective space,
such that the deformation changes in lowest--order
the parameter of the fibre.
All coefficients of the polynomial $Q(\xx{p},\Breve{p})$
corresponding to some moduli are homogeneous polynomials of the moduli,
whose degree is equal to $dm$ minus the common degree
of the corresponding monomial of $\xx{p}$ and $\Breve{p}$.
Due to Lemma~\ref{dualizing sheaf} the regular form $\omega$
may be locally uniquely written as
$$\omega=\frac{\Dot{Q}(\xx{p},\Breve{p})}
{\partial Q(\xx{p},\Breve{p})/\partial \Breve{p}}d\xx{p},$$
where $\Dot{Q}$ is a polynomial with respect to $\Breve{p}$
of degree $d-1$, whose coefficients are holomorphic functions
depending on $\xx{p}$.
We conclude that besides the coefficient of degree zero
(corresponding to $\Breve{p}^d-\xx{p}^m$)
those coefficients of $Q(\xx{p},\Breve{p})$
corresponding to the projective element 
$(\Hat{t}_1,\ldots,\Hat{t}_d,\Breve{t}_1,\ldots,\Breve{t}_m)
\in\mathbb{P}^{d+m-1}$ we start with should be zero,
which are of lower degree than the degree of the unique coefficient,
whose monomial corresponds to the lowest non--vanishing Taylor
coefficient of $\omega$. Furthermore, the coefficient of
$Q(\xx{p},\Breve{p})$ corresponding to the projective element
we start with, whose monomial
corresponds to the lowest non--vanishing Taylor coefficient,
should not vanish.

\noindent
{\bf 5.} Now we shall determine those projective moduli
$(\Hat{t}_1,\ldots,\Hat{t}_d,\Breve{t}_1,\ldots,\Breve{t}_m)$
we have to start with, if $\omega$ does not vanish
at the corresponding \Em{cusp}.
The common degree of the unique monomial $\Dot{Q}(\xx{p},\Breve{p})$,
corresponding locally to a holomorphic non--vanishing form,
is equal to $(d-1)(m-1)$. Therefore, besides the coefficient
of degree zero all coefficients of degree lower than $dm-(d-1)(m-1)=d+m-1$
of the polynomial $Q(\xx{p},\Breve{p})$,
which corresponds to those projective moduli we start with,
should vanish. There are exactly $d+m-2$ such coefficients.
Hence we have to determine the zero set of $d+m-2$
homogeneous polynomials on $\mathbb{P}^{d+m-1}$.
Due to \cite[Chapter~I Theorem~7.2]{Har} this set
is at least one--dimensional. In order to investigate these solutions
let us first determine the derivative of the coefficients
with respect to some infinitesimal change
$(\Dot{\xx{p}},\Dot{\Breve{p}})$
of the polynomials $\xx{p}$ and $\Breve{p}$. The equations
\begin{align*}
\frac{dQ(\xx{p},\Breve{p})}{dt}=
\Dot{Q}(\xx{p},\Breve{p})+
\frac{\partial Q(\xx{p},\Breve{p})}{\partial \xx{p}}\Dot{\xx{p}}+
\frac{\partial Q(\xx{p},\Breve{p})}
       {\partial \Breve{p}}\Dot{\Breve{p}}&=0&
\frac{\partial Q(\xx{p},\Breve{p})}{\partial \xx{p}}d\xx{p}+
\frac{\partial Q(\xx{p},\Breve{p})}
       {\partial \Breve{p}}d\Breve{p}&=0
\end{align*}
imply the identity
$$\Dot{Q}(\xx{p},\Breve{p})
\frac{d\xx{p}}{\partial Q(\xx{p},\Breve{p})/\partial \Breve{p}}=
-\Dot{Q}(\xx{p},\Breve{p})
\frac{d\Breve{p}}{\partial Q(\xx{p},\Breve{p})/\partial \xx{p}}=
\Dot{\xx{p}}d\Breve{p}-\Dot{\Breve{p}}d\xx{p}.$$
For the initial curve defined by $\Breve{p}^d-\xx{p}^m=0$
the right hand side divided by $\Dot{Q}(\xx{p},\Breve{p})$
is equal to $z^{-(d-1)(m-1)}dz$. If we consider the projective deformation
from the initial curve (defined by $\Breve{p}^d-\xx{p}^m=0$) into one
of these solutions, then the derivative with respect to
the lowest power of the deformation parameter
(the parameter of the fibre of the tautological bundle),
which does not vanish, yields a meromorphic one--form
with poles only at $z=\infty$.
We conclude that the corresponding $\Dot{Q}(\xx{p},\Breve{p})$
has to be a polynomial with respect to $z$, whose degree is larger
than $(d-1)(m-1)$. In particular, the coefficients of degree $d+m-1$
of the polynomials $Q(\xx{p},\Breve{p})$,
which correspond to these solutions, cannot vanish.
Therefore, for these deformations the derivative with respect to the
$(d+m-1)$--th power of the deformation parameter at the initial curve
corresponds to a one--form proportional to $dz$.

\noindent
{\bf 6.} Now we claim that the normalizations
of these curves do not have a \Em{cusp}. In fact,
otherwise we shift the \Em{cusp} into the origin ($z=0$).
In this case, due to the foregoing formula,
the derivative with respect to the lowest non--vanishing
(i.\ e.\ the $(d+m-1)$--th) power of the deformation parameter,
corresponds to a one--form, which vanishes at $z=0$.
This contradicts the considerations in step~5.

\noindent
{\bf 7.} Now we show that,
if we choose some non--singular representative
$(\Hat{t}_1,\ldots,\Hat{t}_d,\Breve{t}_1,\ldots,\Breve{t}_m)$
in the orbit under the group described above,
then the derivative of the map,
which maps the moduli onto the coefficients
of the corresponding polynomial $Q(\xx{p},\Breve{p})$,
has rank $d+m-1$. This statement is equivalent to the statement,
that the derivative of the map,
which maps the coefficients of the polynomials
$\xx{p}$ and $\Breve{p}$ with respect to $z$
onto the coefficients of the corresponding polynomial
$Q(\xx{p},\Breve{p})$, has rank $d+m-1$.
The foregoing formula shows that a tangent vector
$(\Dot{\xx{p}},\Dot{\Breve{p}})$
in the space of polynomials $\xx{p}$ and $\Breve{p}$
with respect to $z$ does not change $Q$,
if and only if the right hand side vanishes.
Since the curve defined by the equation $Q(\xx{p},\Breve{p})=0$
does not have \Em{cusps}, $\Dot{\xx{p}}$ has to be zero
at the zeroes of $d\xx{p}$, and $\Dot{\Breve{p}}$ has to be zero at
the zeroes of $d\Breve{p}$. We conclude that
$(\Dot{\xx{p}},\Dot{\Breve{p}})$ has to be proportional to
$(d\xx{p}/dz,d\Breve{p}/dz)$,
which corresponds exactly to translations.

\noindent
{\bf 8.} With the help of the implicit function theorem
\cite[Supplementary material~V.5 Theorem~S.11]{RS1}
we conclude from the last step that
in some arbitrary small neighbourhood of those non--singular moduli,
we start with if $\omega$ does not vanish at the \Em{cusp},
there also exist moduli, whose $Q(\xx{p},\Breve{p})$ has besides the
highest term $\Breve{p}^d-\xx{p}^m$
vanishing coefficients of degree larger than those corresponding to
the lowest non--vanishing coefficient of $\omega$.
Therefore, we may choose these moduli in such a way,
that the corresponding curves $Q(\xx{p},\Breve{p})=0$
also have no \Em{cusps}, and that the derivative of the map,
which maps the moduli onto the coefficients of $Q(\xx{p},\Breve{p})$,
has rank $d+m-1$. This in general defines the projective moduli,
we start with.

\noindent
{\bf 9.} Finally, we multiply the pullback of the section
of the vector bundle of meromorphic forms on the normalizations
with poles of at most first--order at $\infty^{\pm}$
and no other poles, to the blowing up of the moduli space
$\Breve{\moduli}$,
with the parameter of the fibre of the tautological bundle to the
power, which is equal to the order of the zero of $\omega$ plus $d+m-2$,
and which is also equal to the degree of the corresponding coefficient
of $Q$ minus one. With the help of the construction of
Proposition~\ref{smooth moduli} the flow defined by
$\Dot{\yy{p}}d\xx{p}=\omega$ can be locally integrated along the
vector field on the blowing up of the moduli space $\Breve{\moduli}$.
The image of the corresponding curve in $\Breve{\moduli}$
is contained in $\moduli$ and yields the desired curve.
\end{proof}

This proposition can be generalized to deformations
of one--sheeted coverings of the \Em{complex Fermi curves},
which are locally not complete intersection in $\mathbb{C}^2$.

If we consider deformations,
which preserve an anti--holomorphic involution with fixed points,
then the real regular forms describe the corresponding tangent space.
In case of a real cusp of the \Em{complex Fermi curve},
the construction of Proposition~\ref{modified smooth moduli}
has to be done over the real numbers.
In particular, the moduli we start with, should be real.
More precisely, it suffices to find two polynomials
$$\xx{p}(z)=z^d+a_1z^{d-1}+\ldots+a_d
\text{ and }
\Breve{p}(z)=z^m+b_1z^{m-1}+\ldots+b_m$$
with real coefficients $a_1,\ldots,a_d$ and $b_1\ldots,b_m$,
whose coefficients of the corresponding polynomial
$Q(\xx{p},\Breve{p})$ of degree $1,\ldots,d+m-2$ vanish.
Due to the considerations in part~5. of the forgoing proof
this is equivalent to the relation
$$\xx{p}'(z/t)t^d\frac{d}{dt}\Breve{p}(z/t)t^m-
\Breve{p}'(z/t)t^m\frac{d}{dt}\xx{p}(z/t)t^d=
t^{d+m-2}\text{\bf{O}}(1)\text{ with respect to }z.$$
This relation in turn is equivalent to the relation
$$m\xx{p}'(z)\Breve{p}(z)-d\Breve{p}'(z)\xx{p}(z)=
\text{\bf{O}}(1)\text{ with respect to }z.$$
First we remark that due to the invariance under the
translations $z\mapsto z+z_0$, we may assume that $a_1=0$.
Moreover, due to the structure of this relation,
the corresponding equations of order $z^{m+d-1},\ldots,z^{d-1}$
are linear in $b_1,\ldots,b_m$.
Therefore, we may use these equations to obtain $b_1,\ldots,b_m$
as polynomials with respect to $a_2,\ldots,a_d$.
The corresponding polynomial $\Breve{p}(z)$ is determined by
the $m$ lowest Taylor coefficients of
$(z^{-d}\xx{p}(z))^{m/d}=(1+a_1z^{-1}+\ldots+a_dz^{-d})^{m/d}$
at infinity:
$$\Breve{p}(z)=z^m\left(z^{-d}\xx{p}(z)\right)^{m/d}+
\text{\bf{O}}(z^{-1}).$$
It remains to find zeroes of $d-2$ homogeneous polynomials
with respect to $a_2,\ldots,a_d$ of common degree
$m+1,\ldots,m+d-2$,
where the degree of $a_l$ is equal to $l$ for all $l=2,\ldots,d$.
We conclude that for $d=2$
there always exists one unique real solution. 

If $d=3$ and $m$ is even, then there also exists a real solution.
In fact in this case $a_3=0$ is always a solution.
If $m$ is odd and $m+1$ is not divisible by $12$
then there again exists a real solution, since a polynomial,
whose degree is even but not divisible by $12$
with respect to $a_2$ and $a_3$
is either zero for $a_2=0$ or $a_3=0$
or is an odd polynomial either with respect to $a_2$ or $a_3$.
A direct calculation shows that for $m=11$
there also exists a real solution.

If $d=4$, the choice $a_3=0$ reduces the number of equations by one.
Therefore, we have to find a real zero of a polynomial of degree
$m+1$ with respect to $a_2$ and $a_4$.
If $m+1$ is not divisible by $8$ such a polynomial
has always a real solution.

We conjecture that there always exists a real solution
of these equations.
But if $d+m-1$ is even (and therefore not both degrees are odd),
then the real part
of the one--dimensional subspace of the moduli space,
which is constructed in Proposition~\ref{modified smooth moduli}
with the help of a real polynomials constructed above,
covers only one half of the real line in the tangent space
of the moduli space corresponding to any real regular form,
which does not vanish at the cusp.
Since in the simplest case $d=2$ and $m=3$ there exists only one
real solution,
the subspace of the moduli space containing
\Em{complex Fermi curves}, which are invariant under the
anti--holomorphic involution, and which have the same
\Em{geometric genus} as the original \Em{complex Fermi curve}
with the real cusp,
covers only a half-plane of the real tangent space.

\subsection{The compactified moduli spaces with bounded genus}
\label{subsection compactified bounded genus}

In Section~\ref{subsection weak singularity} we shall find some
necessary condition for a \Em{complex Fermi curve} to correspond
to some immersion. Since the Willmore functional is equal to the
\Em{first integral} of the corresponding
\De{Weierstra{\ss} curve}~\ref{Weierstrass curves},
we may transform the variational problem of the
Willmore functional into a variational problem of the
\Em{first integral} on the space of all \Em{complex Fermi curves},
which obey the
\De{Weak Singularity condition}~\ref{weak singularity condition}.
For this purpose we should endow the
set of all \Em{complex Fermi curves} with some topology. By definition all
\Em{complex Fermi curves} are closed subsets of $\mathbb{C}^2$.
\index{topology!of the moduli space}
Consequently we shall use some topology on the set of closed subsets
of $\mathbb{C}^2$. The set of closed subsets of some compact
Hausdorff space has a natural topology in contrast to the space of
closed subsets of a non--compact Hausdorff space \cite{Mi}.
Obviously the closures of all \Em{complex Fermi curves}
considered as subsets of the one--point--compactification
$\overline{\mathbb{C}^2}$ of $\mathbb{C}^2$ are give by the union of
the \Em{complex Fermi curves} (contained in $\mathbb{C}^2$) with the
one point set $\{\infty\}$ containing the single point at infinity.
In the sequel we shall identify the closed subsets of $\mathbb{C}^2$
with the corresponding closures in $\overline{\mathbb{C}^2}$.

Let us therefore introduce a topology on the space of
all closed subsets of some topological space
(\cite[Chapter~5.4]{Ho} and \cite{Mi}).
In \cite{Mi} several topologies are given, and we will use
\index{topology!finite $\sim$}
the so called \Em{finite topology}.
For all open sets $\Set{O}$ the subsets
of the set of closed subsets of the form
$\{\Set{A}\mid\Set{A}\cap\Set{O}\neq \emptyset\}$
and $\{\Set{A}\mid\Set{A}\subset\Set{O}\}$
are open and generate the set of open subsets
of the \Em{finite topology}. We recall some of the properties of
this topology proven in \cite{Mi}. If the underlying space
is a compact Hausdorff space, the space of closed subsets endowed with
this \Em{finite topology} is again a compact Hausdorff space.
Moreover, if the underlying space is in addition metrizable the
space of all closed subsets is also metrizable and coincides with the
topology of the Hausdorff metric\index{Hausdorff metric}.
Finally, the subspace of all finite subsets is dense,
therefore the space of all closed subsets is
separable if the underlying space is separable. The following lemma
is a special case of results proven in \cite{Mi}.

\begin{Lemma} \label{compact metric}
The space of all closed subsets of $\overline{\mathbb{C}^2}$ is a
separable compact metrizable space.\qed
\end{Lemma}

Due to Theorem~\ref{meromorph}, for any pair of potentials $(V,W)$
the mapping $k\mapsto\Op{R}(V,W,k,0)$ is a meromorphic mapping from
$\mathbb{C}^2$ into the compact operators on
$\banach{2}(\Delta)\times\banach{2}(\Delta)$.
Moreover, this mapping uniquely determines the pair of potentials.
The topology of uniform convergence on compact subsets of these
mappings  \cite[Chapter~VII. \S1.]{Co1}
yields a topology on $(V,W)\in\banach{2}(\torus)\times\banach{2}(\torus)$,
which  we call the \Em{compact open topology}
of the resolvents $\Op{R}(V,W,k,0)$.
\index{topology!compact open $\sim$}
Obviously this topology is equivalent to the
\Em{compact open topology} of the resolvents of
$\triv{\Op{D}}_{\xx{\gamma},\yy{\gamma}}(V,W,\xx{p})$
(compare with Section~\ref{subsection spectral projections}).

We claim that the mapping $(V,W)\mapsto\Bar{\fermi}(V,W)$ is continuous
with respect to this topology and the \Em{finite topology}.
In fact, since the two families of subsets
$\{\Bar{\fermi}\subset\overline{\mathbb{C}^2}\mid
\Bar{\fermi}\cap\Set{O}\neq\emptyset\}$
and $\{\Bar{\fermi}\subset\overline{\mathbb{C}^2}\mid
\Bar{\fermi}\subset\Set{O}\}$
indexed by the open subsets $\Set{O}$ of $\overline{\mathbb{C}^2}$
form a basis of the open subsets of the \Em{finite topology},
it suffices to show that the corresponding preimages are open.
If the closure $\Bar{\fermi}$ of a \Em{complex Fermi curve}
is contained in an open subset of the former family,
then there exists a finite $k\in\fermi\cap\Set{O}$.
Furthermore, due to the application of the arguments in the proof of
Theorem~\ref{meromorph} concerning the spectral projections
to the corresponding spectral projections of
$\triv{\Op{D}}_{\xx{\gamma},\yy{\gamma}}(V,W,\xx{p})$,
the same is true for all pairs of potentials in a neighbourhood with
respect to the \Em{compact open topology}.
If the closure $\Bar{\fermi}$ is contained in an open subset of the latter family,
then $\Set{O}$ has to contain $\infty$. Furthermore,
all \Em{complex Fermi curves} of pairs of potentials in a neighbourhood
(with respect to the \Em{compact open topology})
are also disjoint with the complement of $\Set{O}$,
since this complement is a compact subset of $\mathbb{C}^2$.
This proves

\begin{Lemma}\label{continuity of fermi curves}
The mapping $(V,W)\mapsto\Bar{\fermi}(V,W)$ is a continuous mapping from
$(V,W)\in\banach{2}(\torus)\times\banach{2}(\torus)$
with the \Em{compact open topology} of $\Op{R}(V,W,k,0)$
into the closed subsets of $\overline{\mathbb{C}^2}$
with the \Em{finite topology}.\qed
\end{Lemma}

Let $\moduli_{\lattice,\eta}$ and $\moduli_{\lattice,\eta,\willmore}$
\index{moduli space!$\moduli_{\lattice,\eta},
                     \moduli_{\lattice,\eta,\willmore},
                     \moduli_{\lattice,\eta,\sigma},
                     \moduli_{\lattice,\eta,\sigma,\willmore}$|(}
denote the moduli spaces
\begin{eqnarray*}
\moduli_{\lattice,\eta}&=&
\left\{ \left. \fermi(U,\Bar{U}) \right|
U\in \banach{2}(\torus) \right\}\\
\moduli_{\lattice,\eta,\willmore}&=&\left\{ \left.
\fermi(U,\Bar{U}) \right|
U\in \banach{2}(\torus)\text{ and }
4\|U\|^2\leq \willmore \right\}.
\end{eqnarray*}
Both sets are are subsets of the space of closed subsets of
$\overline{\mathbb{C}^2}$.
The subset of these moduli spaces, which contain the corresponding
\Em{complex Fermi curves} of real potentials are denoted by
$\moduli_{\lattice,\eta,\sigma}$
and
$\moduli_{\lattice,\eta,\sigma,\willmore}$,
respectively.
\index{moduli space!$\moduli_{\lattice,\eta},
                     \moduli_{\lattice,\eta,\willmore},
                     \moduli_{\lattice,\eta,\sigma},
                     \moduli_{\lattice,\eta,\sigma,\willmore}$|)}

We will see that neither $\moduli_{\lattice,\eta}$ nor
$\moduli_{\lattice,\eta,\willmore}$ is compact,
if $\willmore$ is not smaller than $4\pi$
(compare with Remark~\ref{sobolev 1} and Remark~\ref{sobolev 2}). In fact, in
Section~\ref{subsection complex Fermi curves of finite genus} we have
introduced conditions on compact one--dimensional complex spaces,
which for the dense subclass of \Em{connected} complex spaces implies,
that they are \Em{complex Fermi curves}.
In Lemma~\ref{disconnected} we will investigate the subclass of
complex spaces with two connected components.
It will turn out that the simplest example is not a
\Em{complex Fermi curve} of any potential, and that the corresponding
\Em{first integral} is equal to $4\pi$. But they are the
\Em{complex Fermi curves} of perturbations of Dirac operators
acting on direct sums of line bundles of degree $\pm1$.
In the next section we shall prove that the set of all
\Em{complex Fermi curves} of higher-degree perturbations of
Dirac operators, whose \Em{first integral} is bounded by any
$\willmore>0$ is compact.
Higher-degree perturbations of Dirac operators
also occur in `quaternionic function theory' \cite{PP}.

In this section we restrict ourself to the moduli spaces of
\Em{complex Fermi curves} of bounded \Em{geometric genus}:
Let $\moduli_{g,\lattice,\eta}$ and $\moduli_{g,\lattice,\eta,\willmore}$
\index{moduli space!$\moduli_{g,\lattice,\eta},
                     \moduli_{g,\lattice,\eta,\willmore},
                     \moduli_{g,\lattice,\eta,\sigma},
                     \moduli_{g,\lattice,\eta,\sigma,\willmore}$|(}
denote the moduli spaces
\begin{eqnarray*}
\moduli_{g,\lattice,\eta}&=&\left\{\fermi\in\moduli_{\lattice,\eta}\mid
\text{ geometric genus of $\fermi$ is not larger than $g$}\right\}\\ 
\moduli_{g,\lattice,\eta,\willmore}&=&
\moduli_{\lattice,\eta,\willmore}\cap
\moduli_{g,\lattice,\eta}.
\end{eqnarray*}
The corresponding subsets of these moduli spaces, which contain only
\Em{complex Fermi curves} of real potentials are denoted by
$\moduli_{g,\lattice,\eta,\sigma}$
and
$\moduli_{g,\lattice,\eta,\sigma,\willmore}$,
respectively.
\index{moduli space!$\moduli_{g,\lattice,\eta},
                     \moduli_{g,\lattice,\eta,\willmore},
                     \moduli_{g,\lattice,\eta,\sigma},
                     \moduli_{g,\lattice,\eta,\sigma,\willmore}$|)}

In Section~\ref{subsection complex Fermi curves of finite genus} we
characterized these \Em{complex Fermi curves} of finite genus by the
conditions \Em{Quasi--momenta}~(i)--(iii) of
Section~\ref{subsection complex Fermi curves of finite genus}.
There exists an explicit
formula of the \Em{first integral}, whose restriction to
\De{Weierstra{\ss} curves}~\ref{Weierstrass curves} is equal to the
Willmore functional, in terms of the data $(\Spa{Y},\infty^-,\infty^+,k)$.

\begin{Lemma} \label{Willmore functional}
\index{integral!first $\sim$ $\willmore$}
If the \Em{complex Fermi curve}
$\fermi(\Spa{Y},\infty^-,\infty^+,k)$ of some
data $(\Spa{Y},\infty^-,\infty^+,k)$, which fulfills conditions
\Em{Quasi--momenta}~(i)--(ii)
of Section~\ref{subsection complex Fermi curves of finite genus},
is equal to the \Em{complex Fermi curve} of
some pair of potentials $(V,W)$, then the corresponding
\Em{first integral} is equal to the residue\index{residue}
$$4\int\limits_{\torus} V(x)W(x)d^2x=
\willmore(\Spa{Y},\infty^-,\infty^+,k)=8\pi^2\sqrt{-1}
\vol(\torus)
\res\limits_{\infty^+}
\left(k_1dk_2\right).$$
\end{Lemma}

Obviously the residue on the right hand side does not depend on the
choice of the branch of $k$. Below we will show, that whenever the
data fulfill condition \Em{Quasi--momenta}~(iii)
of Section~\ref{subsection complex Fermi curves of finite genus}
these numbers are non--negative.

\begin{proof}
The \Em{complex Fermi curve} corresponding to the zero potentials is
given by the union of the solutions of the equations
$g(k+\kappa,k+\kappa)=0$ for $\kappa\in\lattice\dual$.
An easy calculation shows that in this
case the formula is correct. If we apply formula
Lemma~\ref{residue}~(iv) to the variation described in
Lemma~\ref{compatibility}~(iv) the right hand side is equal to
$1/(\pi^2\sqrt{-1})$ times the variation of the \Em{first integral}.
Since this variation corresponds to the function, which in some
neighbourhood of $\infty^+$ is equal to $1$ and in some neighbourhood of
$\infty^-$ equal to $-1$, the left hand side is equal to the residue
of $\Omega_{V,W}(\cdot,\cdot)$ at $\infty^+$ minus the residue of this
form at $\infty^-$. If we consider the function $\yy{p}$ as a
function depending on $\xx{p}$ and the deformation parameters, we
conclude from the definition of $\Omega$ that the left hand side is
equal to the residue of the form $\var\yy{p}d\xx{p}$ at $\infty^+$.
This implies that the \Em{first integral}
is equal to the residue at $\infty^+$ of the form
$2\pi^2\sqrt{-1}\yy{p}d\xx{p}$. The right hand side of the formula
of the lemma is invariant under linear conformal transformations of
the lattice $\lattice$. If we choose $\gamma=(1,0)$ the formula
follows.
\end{proof}

Consequently, define $\Bar{\moduli}_{g,\lattice,\eta}$
\index{moduli space! $\Bar{\moduli}_{g,\lattice,\eta},
                      \Bar{\moduli}_{g,\lattice,\eta,\willmore},
                      \Bar{\moduli}_{g,\lattice,\eta,\sigma},
                      \Bar{\moduli}_{g,\lattice,\eta,\sigma,\willmore}$|(}
as the set
of all \Em{complex Fermi curves} $\fermi(\Spa{Y},\infty^-,\infty^+,k)$,
whose data fulfill conditions \Em{Quasi--momenta}~(i)--(iii)
of Section~\ref{subsection complex Fermi curves of finite genus}
and whose \Em{geometric genus} is not larger than $g$.
Analogously we denote the subset of all varieties
$\fermi(\Spa{Y},\infty^-,\infty^+,k)$ of data, whose
\Em{first integral}
$\willmore(\Spa{Y},\infty^-,\infty^+,k)$ is not larger than
$\willmore$ by $\Bar{\moduli}_{g,\lattice,\eta,\willmore}$.
The subsets of these moduli spaces containing the corresponding
\Em{complex Fermi curves}, which are invariant under the involution
$\sigma$, are denoted by $\Bar{\moduli}_{g,\lattice,\eta,\sigma}$
and $\Bar{\moduli}_{g,\lattice,\eta,\sigma,\willmore}$,
respectively. Now we can state the main theorem of this section.
\index{moduli space! $\Bar{\moduli}_{g,\lattice,\eta},
                      \Bar{\moduli}_{g,\lattice,\eta,\willmore},
                      \Bar{\moduli}_{g,\lattice,\eta,\sigma},
                      \Bar{\moduli}_{g,\lattice,\eta,\sigma,\willmore}$|)}

\begin{Theorem} \label{compactification}
\index{moduli space!compactification of the $\sim$}
\index{compactification!of the moduli space}
The spaces $\Bar{\moduli}_{g,\lattice,\eta,\willmore}$ are compact
and the completions of $\moduli_{g,\lattice,\eta,\willmore}$.
\end{Theorem}

For the proof of this theorem we need some preparation. Let
$\Bar{\moduli}_{g,\eta}$
\index{moduli space! $\Bar{\moduli}_{g,\eta}$}
denote the set of all varieties
$\fermi(\Spa{Y},\infty^-,\infty^+,k)$, whose data obey conditions
\Em{Quasi--momenta}~(i) and (iii)
of Section~\ref{subsection complex Fermi curves of finite genus},
but instead of condition \Em{Quasi--momenta}~(ii) we assume condition
\begin{description}
\item[Quasi--momenta (ii')]
  \index{condition!quasi--momenta (ii')}
  \index{quasi--momenta!condition $\sim$ (ii')}
  The difference of two arbitrary branches of the
  function $k$ is some element of $\mathbb{R}^2$.
\end{description}
In this case $d\xx{p}$ and $d\yy{p}$ are meromorphic forms, and
the multi--valued function $k$ is determined by this form up to
the addition of some complex numbers. So let us investigate for a
given data $(\Spa{Y},\infty^-,\infty^+,k)$ the space of such forms
$d\xx{p}$ and $d\yy{p}$.

\begin{Lemma} \label{meromorphic forms}
Let $\Spa{Y}$ be a compact smooth Riemann surface of genus $g$ with two
marked points $\infty^-$ and $\infty^+$ with the property that each
connected component of $\Spa{Y}$ contains at least one point. Moreover, let
$\eta$ be an anti--holomorphic involution, which interchanges the two
marked points. Then for any local parameter $z$ at $\infty^+$, which
takes the value zero at $\infty^+$, there exists a unique meromorphic
form $dp$ with the following properties:
\begin{description}
\index{condition!form (i)--(iv)|(}
\index{form!condition $\sim$ (i)--(iv)|(}
\item[Form (i)]
  $dp-d(1/z)$ extends to a holomorphic form on some
  neighbourhood of $\infty^+$.
\item[Form (ii)]
  The integral of $dp$ along any cycle in $H_1(\Spa{Y},\mathbb{Z})$ is real.
\item[Form (iii)]
  $dp$ transforms under $\eta$ as
  $\eta^{\ast}dp=-\Bar{dp}$.
\item[Form (iv)]
The integral of $dp$ along any real cycle
(i.\ e.\ the cycle is homolog to the image under $\eta$)
vanishes.
\index{condition!form (i)--(iv)|)}
\index{form!condition $\sim$ (i)--(iv)|)}
\end{description}
\end{Lemma}

\begin{proof} Obviously conditions \Em{Form}~(i)--(iii) are
  equivalent to conditions \Em{Form}~(i) together with
  the conditions \Em{Form}~(iii)--(iv). Hence we may assume this
  last triple of conditions.
  The Riemann--Roch theorem
  \cite[Theorem~16.9 and Theorem~17.16]{Fo} implies that
  the space of meromorphic forms, which have poles of at most second
  order without residues at the two marked points and no other poles
  is a $(g+2)$--dimensional complex space. Since the involution $\eta$
  acts on this space as an anti--linear involution, the eigenspace of
  this action corresponding to the eigenvalue $-1$ is a
  $(g+2)$--dimensional real space. Hence  there exists some $dp$
  obeying conditions \Em{Form}~(i) and (ii). This form is determined
  up to the addition of some holomorphic form $\omega$ with the property
  that $\eta^{\ast}\omega=-\Bar{\omega}$. Since the real part of
  $H_1(\Spa{Y},\mathbb{Z})$ forms a Lagrangian subgroup with respect to the
  intersection form (i.\ e.\ a subgroup, which coincides with the group
  of all elements having zero intersection form with all elements of
  this subgroup), any set of generators of this real part gives rise
  to a dual base of holomorphic differential forms $\omega$ with the
  property that $\eta^{\ast}\omega=-\Bar{\omega}$
  \cite[III.3.3.~Proposition]{FK}. This proves the Lemma.
\end{proof}

We conclude that for all data
$(\Spa{Y},\infty^-,\infty^+,\eta,k)\in\Bar{\moduli}_{g,\eta}$
the components $\xx{p}$ completely determine the other components, and
more generally, any real component of $k$ determines the function $k$.
So in the sequel we will associate to the elements of
$\Bar{\moduli}_{g,\lattice,\eta}$ data of the form
$(\Spa{Y},\infty^-,\infty^+,p)$ where $p$ is a multi--valued in
$\mathbb{C}$ instead of $\mathbb{C}^2$ (either equal to
$\xx{p}$ or any other real component of $k$).
In this case we assume that these data obey the conditions analogous
to \Em{Quasi--momenta}~(i), (ii') and (iii).
Obviously, $(\Spa{Y},\infty^-,\infty^+,\xx{p})$ and
$(\Spa{Y},\infty^-,\infty^+,\yy{p})$
are data of this form. Also the formula of
Lemma~\ref{Willmore functional} for the \Em{first integral} extends
to elements of $\Bar{\moduli}_{g,\eta}$. On these larger
moduli spaces we may \Em{renormalize} the functions $k$ by
multiplication with some positive number.
The \Em{first integral} changes under these
transformations by multiplication with the square of this positive number.
It is quite easy to see that the quotient of the moduli spaces
$\Bar{\moduli}_{g,\eta}$ modulo these transformations and modulo
shifts of $k$ by some additive constant, is the moduli space of
all compact Riemann surfaces of genus $g$ with two marked points
$\infty^{\pm}$, at most two connected components, each containing at
least one of the two marked points, and the
anti--holomorphic involution $\eta$ without fixed points
interchanging $\infty^{\pm}$. Therefore, the description
of these Riemann surfaces in Lemma~\ref{gluing rule} yields a
parameterization of this moduli space. A simple count of degrees of
freedom shows that $\Bar{\moduli}_{g,\eta}$ has real dimension
$3g+2$. This coincides with the dimension of the moduli space of the
corresponding Riemann surfaces with two marked points ($3g-1$) plus
the degree of freedom of renormalization of $k$
and the two additive degrees of freedom given by the complex
additive constant just mentioned.

Let us now show that for all these data
$(\Spa{Y},\infty^-,\infty^+,p)$
the Riemann surface $\Spa{Y}$ may be realized as two copies of
$\mathbb{P}^1$, which are glued along several cuts.
In the sequel we will use this
construction for different choices of the function $p$, therefore we
shall decorate the corresponding symbols with an index $p$.
First we choose a branch of the
function $p$ in some small neighbourhood of $\infty^+$.
Obviously there exists an open neighbourhood $\Set{U}_{p}$ of $\infty^+$
obeying the following conditions:
\begin{description}
\index{condition!neighbourhood (i)--(ii)|(}
\index{neighbourhood!condition $\sim$ (i)--(ii)|(}
\item[Neighbourhood (i)]
  The fixed branch of $p$ near $\infty^+$ maps
  $\Set{U}_{p}$ biholomorphically onto the complement
  in $\mathbb{P}^1$ of
  finitely many disjoint compact convex sets contained in
  $\mathbb{C}\subset\mathbb{P}^1$.
  In particular, $\Set{U}_{p}$ contains no
  zero of $dp$.
\item[Neighbourhood (ii)]
  Each connected component of the complement
  of $\Set{U}_{p}\cup\eta(\Set{U}_{p})$
  contains at most one connected component
  of the boundary of the closures of $\Set{U}_{p}$.
\index{condition!neighbourhood (i)--(ii)|)}
\index{neighbourhood!condition $\sim$ (i)--(ii)|)}
\end{description}
Due to condition \Em{Neighbourhood}~(i) the connected components of 
the boundary of the closure of $\Set{U}_{p}$ are in one to one
correspondence with the compact convex subsets mentioned there.

\begin{Lemma} \label{neighbourhood}
\index{neighbourhood!$\Set{U}_{p}^{\pm}$}
There exists a unique maximal open neighbourhood $\Set{U}_{p}^+$ of
$\infty^+$ obeying conditions \Em{Neighbourhood}~(i)--(ii).
Moreover, the set $\Set{U}_{p}^-=\eta(\Set{U}_{p}^+)$ is the
unique maximal open neighbourhood of $\infty^-$
obeying the corresponding conditions
for the branch $-\eta^{\ast}\Bar{p}$ near $\infty^-$.
Both sets $\Set{U}_{p}^+$ and $\Set{U}_{p}^-$ are disjoint.
\end{Lemma}

\begin{proof} Let us first prove that the union of two open
neighbourhoods $\Set{U}_{p}$ and $\Set{U}_{p}'$ of $\infty^+$ obeying
these conditions obeys also these conditions. Obviously the union of
two complements of finitely many disjoint compact convex sets is again
a complement of finitely many disjoint compact convex sets. 
Due to condition \Em{Quasi--momenta}~(ii') the subset of all elements in
the intersection of the images of $\Set{U}_{p}$ and
$\Set{U}_{p}'$ under the corresponding branches of $p$,
whose preimages are the same element
of $\Set{U}_{p}\cap \Set{U}_{p}'$, is open and closed.
Hence it is a connected component of
$\Set{U}_{p}\cap \Set{U}_{p}'$.
But the image of $\Set{U}_{p}\cap \Set{U}_{p}'$
under the corresponding branches of $p$ is the complement of finitely
many not necessary disjoint compact convex subsets of $\mathbb{P}^1$
and therefore connected. Hence the branches of $p$ corresponding to
$\Set{U}_{p}$ and $\Set{U}_{p}'$, coincide,
and $\Set{U}_{p}\cup \Set{U}_{p}'$ obeys condition
\Em{Neighbourhood}~(i). Since the complement of
$\Set{U}_{p}\cup \Set{U}_{p}'\cup\eta(\Set{U}_{p}\cup
\Set{U}_{p}')$ is equal to the
intersection of the complements of
$\Set{U}_{p}\cup\eta(\Set{U}_{p})$ and
$\Set{U}_{p}'\cup\eta(\Set{U}_{p}')$,
the union $\Set{U}_{p}\cup \Set{U}_{p}'$ obeys also condition
\Em{Neighbourhood}~(ii). This argumentation extends to the
union of an arbitrary set of neighbourhoods obeying these
conditions. Therefore, the union of all neighbourhoods obeying these
conditions is the unique maximal $\Set{U}_{p}^+$. Obviously, the set
$\Set{U}_{p}^-=\eta(\Set{U}_{p}^+)$ has the corresponding properties
(with $\infty^-$ interchanged by $\infty^+$).
Arguments similar to those used above show that
$\Set{U}_{p}^+\cap \Set{U}_{p}^-$ has to be open
and closed. Therefore, it has to be empty,
since it does not contain $\infty^{\pm}$.
\end{proof}

In the sequel we want to extend the open sets
$\Set{U}_{p}^+$ and $\Set{U}_{p}^-$
to open disjoint sets, which are mapped by the corresponding
branches of $p$ biholomorphically onto the complement of
finitely many finite parts of straight lines in
$\mathbb{C}\subset\mathbb{P}^1$. The copy of
$\mathbb{P}^1$ containing the image of $\Set{U}_{p}^+$ is denoted by
$\mathbb{P}_{p}^+$ and the copy containing the image of $\Set{U}_{p}^-$
is denoted by $\mathbb{P}_{p}^-$. It will turn out that the
Riemann surface $\Spa{Y}$ may be realized by gluing $\mathbb{P}_{p}^+$
and $\mathbb{P}_{p}^-$ along several cuts, and the
gluing map is always of the form $p\mapsto p$ plus some element of
$\mathbb{R}$. We shall describe such surgeries by figures.
On these figures small letters will denote the values
of the functions $p$ at the zeroes of $dp$.
If a sheet of $\mathbb{P}_{p}^+$ is glued with some
sheets of $\mathbb{P}_{p}^-$ these two sheets will be denoted by a
capital letter with exponent $\pm$. If two sheets of
$\mathbb{P}_{p}^+$ or $\mathbb{P}_{p}^-$ are these
two sheets are denoted by a capital letter with indices $1$ and $2$
and exponent $\pm$. Let us first give three typical examples:
\begin{description}
\item[A horizontal cut.]
\index{cut!horizontal $\sim$}
  In this case $dp$ has two simple zeroes at the
  values $p=a$ and $p=b$ with respect to the branch of $p$ near
  $\infty^+$ and the values $p=-\Bar{b}$ and $p=-\Bar{a}$
  with respect to the branch of $p$ near $\infty^-$, respectively,
  where $a$ and $b$ are complex numbers with $\Im(a)=\Im(b)$ and
  $\Re(a)<\Re(b)$.
  The difference between the branch near $\infty^+$ minus the branch near
  $\infty^-$ is therefore equal to $a+\Bar{b}=b+\Bar{a}$. Finally, the
  Riemann surface $\Spa{Y}$ may be realized by gluing the slit along the
  straight line connecting $a$ and $b$ of $\mathbb{P}_{p}^+$ with the
  slit along the straight line connecting $-\Bar{b}$ and $-\Bar{a}$ of
  $\mathbb{P}_{p}^-$. More precisely, for the gluing we use the map
  $p\mapsto p+a+\Bar{b}$, which transforms the branch near of
  $\mathbb{P}_{p}^+$ into the branch of $\mathbb{P}_{p}^-$.

\setlength{\unitlength}{1cm}
\begin{picture}(15,3)
\put(0,0){\begin{picture}(5,3)
\put(1,2){\line(1,0){3}}
\put(1,2){\circle*{.07}}
\put(4,2){\circle*{.07}}
\put(2.3,2.2){\makebox(.5,.5)[b]{$A^+$}}
\put(.8,1.3){\makebox(.5,.5)[t]{$a$}}
\put(2.3,1.3){\makebox(.5,.5)[t]{$B^+$}}
\put(3.8,1.3){\makebox(.5,.5)[t]{$b$}}
\put(,0){\makebox(3,.5)[t]{figure of $\mathbb{P}_{p}^+$}} 
\end{picture} }

\put(10,0){\begin{picture}(5,3)
\put(1,2){\line(1,0){3}}
\put(1,2){\circle*{.07}}
\put(4,2){\circle*{.07}}
\put(2.3,2.2){\makebox(.5,.5)[b]{$B^-$}}
\put(.8,1.3){\makebox(.5,.5)[t]{$-\Bar{b}$}}
\put(2.3,1.3){\makebox(.5,.5)[t]{$A^-$}}
\put(3.8,1.3){\makebox(.5,.5)[t]{$-\Bar{a}$}}
\put(1,0){\makebox(3,.5)[t]{figure of $\mathbb{P}_{p}^-$}} 
\end{picture} }

\end{picture}
\item[A pair of parallel cuts.]
\index{cut!pair of parallel $\sim$s}
  In this case $dp$ has four simple
  zeroes at the values $p=a_1,a_2,b_1,b_2$ on $\mathbb{P}_{p}^+$, where
  $\Im(a_1)=\Im(a_2)$, $\Im(b_1)=\Im(b_2)$, $\Im(a_1)>\Im(b_1)$,
  $\Re(a_1)<\Re(a_2)$ and $\Re(b_1)-\Re(b_2)=\Re(a_1)-\Re(a_2)$.
  Also there are four other zeroes at the
  values $p=-\Bar{a}_2,-\Bar{a}_1,-\Bar{b}_2,-\Bar{b}_1$. In this case
  the slit along the straight lines of $\mathbb{P}_{p}^+$ connecting
  $a_1$ and $b_1$ is glued with the slit along straight lines
  connecting $a_2$ and $b_2$, and the analog slits of $\mathbb{P}_{p}^-$
  are also glued. Hence on $\mathbb{P}_{p}^+$ and on
  $\mathbb{P}_{p}^-$  two pairs of zeroes of $dp$ are identified,
  and all together we have four zeroes of $dp$. Besides this two pairs of
  parallel slits we have two horizontal slits. One connects the two
  slits of $\mathbb{P}_{p}^+$ and the other connects the two slits of
  $\mathbb{P}_{p}^-$. Both are part of a line $\Im(p)=c$ with some
  arbitrary constant in the interval $\Im(a_1)\geq c\geq
  \Im(b_1)$. Since the two slits at the end of each of these
  horizontal slits are glued, these horizontal slits are
  glued to form circles. Finally, these two horizontal slits of
  $\mathbb{P}_{p}^+$ and $\mathbb{P}_{p}^-$ are glued by
  some map $p\mapsto p$ plus some real constant. But if we add to this
  constant the period $\Re(a_2)-\Re(a_1)$ the identification of the
  two circles does not change. Obviously we may shift these two slits
  (this means that we change $c$  in the interval
  $\Im(a_1)\geq c\geq \Im(b_1)$)
  and simultaneously shift the constant of the gluing map in some
  specific way, without changing the data
  $(\Spa{Y},\infty^-,\infty^+,\eta,p)$ and  the branch of $p$ at $\infty^+$.

\setlength{\unitlength}{1cm}
\begin{picture}(15,6)
\put(0,0){\begin{picture}(6,6)

\put(1,2){\line(1,3){1}}
\put(1,2){\circle*{.07}}
\put(2,5){\circle*{.07}}
\put(.8,3.2){\makebox(.5,.5)[b]{$A_1^+$}}
\put(.8,1.3){\makebox(.5,.5)[t]{$b_1$}}
\put(1.8,5.2){\makebox(.5,.5)[b]{$a_1$}}
\put(1.4,2.4){\makebox(.5,.5)[t]{$C_1^+$}}
\put(1.9,4.1){\makebox(.5,.5)[t]{$B_1^+$}}

\put(3,2){\line(1,3){1}}
\put(3,2){\circle*{.07}}
\put(4,5){\circle*{.07}}
\put(3.7,3.2){\makebox(.5,.5)[b]{$A_2^+$}}
\put(2.8,1.3){\makebox(.5,.5)[t]{$b_2$}}
\put(3.8,5.2){\makebox(.5,.5)[b]{$a_2$}}
\put(2.6,2.4){\makebox(.5,.5)[t]{$C_2^+$}}
\put(3.2,4.1){\makebox(.5,.5)[t]{$B_2^+$}}

\put(1.5,3.5){\line(1,0){2}}
\put(2.2,3.6){\makebox(.5,.5)[b]{$D^+$}}
\put(2.2,2.9){\makebox(.5,.5)[t]{$E^+$}}

\put(1,0){\makebox(3,.5)[t]{figure of $\mathbb{P}_{p}^+$}} 
\end{picture} }

\put(9,0){\begin{picture}(6,6)

\put(5,2){\line(-1,3){1}}
\put(5,2){\circle*{.07}}
\put(4,5){\circle*{.07}}
\put(4.7,3.2){\makebox(.5,.5)[b]{$A_1^-$}}
\put(4.7,1.3){\makebox(.5,.5)[t]{$-\Bar{b}_1$}}
\put(3.7,5.2){\makebox(.5,.5)[b]{$-\Bar{a}_1$}}
\put(4.1,2.4){\makebox(.5,.5)[t]{$C_1^-$}}
\put(3.6,4.1){\makebox(.5,.5)[t]{$B_1^-$}}

\put(3,2){\line(-1,3){1}}
\put(3,2){\circle*{.07}}
\put(2,5){\circle*{.07}}
\put(1.8,3.2){\makebox(.5,.5)[b]{$A_2^-$}}
\put(2.7,1.3){\makebox(.5,.5)[t]{$-\Bar{b}_2$}}
\put(1.7,5.2){\makebox(.5,.5)[b]{$-\Bar{a}_2$}}
\put(2.9,2.4){\makebox(.5,.5)[t]{$C_2^-$}}
\put(2.3,4.1){\makebox(.5,.5)[t]{$B_2^-$}}

\put(4.5,3.5){\line(-1,0){2}}
\put(3.3,3.6){\makebox(.5,.5)[b]{$E^-$}}
\put(3.3,2.9){\makebox(.5,.5)[t]{$D^-$}}

\put(2,0){\makebox(3,.5)[t]{figure of $\mathbb{P}_{p}^-$}}
\end{picture} }

\end{picture}
\item[A modified pair of parallel cuts.]
\index{cut!modified pair of parallel $\sim$s}
This is a variation of the \Em{pairs of parallel cuts},
where the two sheets at the horizontal cut of $\mathbb{P}_{p}^+$ are
glued by some map $p\mapsto p$ plus some real constant, and
similarly for the two sheets at the horizontal cut of $\mathbb{P}_{p}^-$.
In this case the Riemann surface $\Spa{Y}$ has two connected components:

\setlength{\unitlength}{1cm}
\begin{picture}(15,6)
\put(0,0){\begin{picture}(6,6)

\put(1,2){\line(1,3){1}}
\put(1,2){\circle*{.07}}
\put(2,5){\circle*{.07}}
\put(.8,3.2){\makebox(.5,.5)[b]{$A_1^+$}}
\put(.8,1.3){\makebox(.5,.5)[t]{$b_1$}}
\put(1.8,5.2){\makebox(.5,.5)[b]{$a_1$}}
\put(1.4,2.4){\makebox(.5,.5)[t]{$C_1^+$}}
\put(1.9,4.1){\makebox(.5,.5)[t]{$B_1^+$}}

\put(3,2){\line(1,3){1}}
\put(3,2){\circle*{.07}}
\put(4,5){\circle*{.07}}
\put(3.7,3.2){\makebox(.5,.5)[b]{$A_2^+$}}
\put(2.8,1.3){\makebox(.5,.5)[t]{$b_2$}}
\put(3.8,5.2){\makebox(.5,.5)[b]{$a_2$}}
\put(2.6,2.4){\makebox(.5,.5)[t]{$C_2^+$}}
\put(3.2,4.1){\makebox(.5,.5)[t]{$B_2^+$}}

\put(1.5,3.5){\line(1,0){2}}
\put(2.2,3.6){\makebox(.5,.5)[b]{$D^+_1$}}
\put(2.2,2.9){\makebox(.5,.5)[t]{$D^+_2$}}

\put(1,0){\makebox(3,.5)[t]{figure of $\mathbb{P}_{p}^+$}} 
\end{picture} }

\put(9,0){\begin{picture}(6,6)

\put(5,2){\line(-1,3){1}}
\put(5,2){\circle*{.07}}
\put(4,5){\circle*{.07}}
\put(4.7,3.2){\makebox(.5,.5)[b]{$A_1^-$}}
\put(4.7,1.3){\makebox(.5,.5)[t]{$-\Bar{b}_1$}}
\put(3.7,5.2){\makebox(.5,.5)[b]{$-\Bar{a}_1$}}
\put(4.1,2.4){\makebox(.5,.5)[t]{$C_1^-$}}
\put(3.6,4.1){\makebox(.5,.5)[t]{$B_1^-$}}

\put(3,2){\line(-1,3){1}}
\put(3,2){\circle*{.07}}
\put(2,5){\circle*{.07}}
\put(1.8,3.2){\makebox(.5,.5)[b]{$A_2^-$}}
\put(2.7,1.3){\makebox(.5,.5)[t]{$-\Bar{b}_2$}}
\put(1.7,5.2){\makebox(.5,.5)[b]{$-\Bar{a}_2$}}
\put(2.9,2.4){\makebox(.5,.5)[t]{$C_2^-$}}
\put(2.3,4.1){\makebox(.5,.5)[t]{$B_2^-$}}

\put(4.5,3.5){\line(-1,0){2}}
\put(3.3,3.6){\makebox(.5,.5)[b]{$D^-_1$}}
\put(3.3,2.9){\makebox(.5,.5)[t]{$D^-_2$}}

\put(2,0){\makebox(3,.5)[t]{figure of $\mathbb{P}_{p}^-$}}
\end{picture} }

\end{picture}

In the special case that the constant of the gluing map is zero,
we could remove the two horizontal cuts. However, we will not do so,
because otherwise $k$ will not extend to a single--valued function
on $\mathbb{P}_{p}^+$.
\end{description}

\begin{Lemma} \label{gluing rule}
\index{gluing rule}
\index{neighbourhood!$\mathbb{P}_{p}^{\pm}$}
Let the data $(\Spa{Y},\infty^-,\infty^+,p)$ be given by any data
$(\Spa{Y},\infty^-,\infty^+,k)\in\Bar{\moduli}_{g,\eta}$
with $p=k_1$ (or any real component of $k$).
Then for all branches of $p$ in some neighbourhood of $\infty^+$
the Riemann surface $\Spa{Y}$ may be realized by two copies of
$\mathbb{P}^1$, which are glued along a combination of
finitely many \Em{horizontal cuts} and finitely many
\Em{pairs of parallel cuts}.
The involution $\eta$ is given by the map $p\mapsto-\Bar{p}$
interchanging $\mathbb{P}_{p}^+$ and $\mathbb{P}_{p}^-$.
Moreover, the \Em{first integral} of the data
$(\Spa{Y},\infty^-,\infty^+,k)$ is
equal to the Dirichlet integral\index{Dirichlet integral}
\index{integral!Dirichlet $\sim$}
\index{integral!first $\sim$ $\willmore$}
of the positive form
$$2\pi\sqrt{-1}\vol(\torus)
\left(dk_1+\sqrt{-1}k_2\right)\wedge
\left(\overline{dk_1+\sqrt{-1}k_2}\right)$$
over $\mathbb{P}_{p}^+$. In particular, the
\Em{first integral} is always some non--negative real number.
The \Em{first integral} vanishes, if and only if the data
$(\Spa{Y},\infty^-,\infty^+,k)$ are isomorphic to the following data:
$\Spa{Y}$ is biholomorphic to the disjoint union of
$\mathbb{P}_{k_1}^+\cup\mathbb{P}_{k_1}^-$ with the
natural marked points $\infty^-$ and $\infty^+$.
The component $k_1$  is the natural meromorphic function on both copies
of $\mathbb{P}_{k_1}$. The component $k_2$ is on $\mathbb{P}_{k_1}^+$
equal to $\sqrt{-1}k_1$ and on $\mathbb{P}_{k_1}^-$ equal to
$-\sqrt{-1}k_1$. Finally, $\eta$ is the involution induced by the map
$k\mapsto -\Bar{k}$, which interchanges the components
$\mathbb{P}_{k_1}^+$ and $\mathbb{P}_{k_1}^-$.
\end{Lemma}

\begin{proof}{\bf 1.} First we claim that the convex compact sets
  excluded from the image of $\Set{U}_{p}^+$,
  mentioned in \Em{Neighbourhood}~(i), are the
  convex hulls of the values of $p$ at finitely many zeroes of $dp$.
  For all $\Set{U}_{p}$ obeying
  conditions \Em{Neighbourhood}~(i)--(ii), the connected components of
  the boundary of $\Set{U}_{p}$ are mapped by $p$
  on convex Jordan curves around
  the excluded domains in $\mathbb{C}\subset\mathbb{P}_{p}^+$.
  Since $p$ is locally biholomorphic if and only if $dp$ has no zero,
  we may enlarge $\Set{U}_{p}$ and
  therefore shrink these Jordan curves,
  until they form the boundary of the convex hull of
  finitely many values of $p$ at zeroes of $dp$.
  In case where the image of $\Set{U}_{p}^+$ is the complement of
  finitely many  disjoint \Em{horizontal cuts},
  the Riemann surface $\Spa{Y}$ is biholomorphic to the union of
  $\mathbb{P}_{p}^+$ and $\mathbb{P}_{p}^-$
  glued in the prescribed way.

  \noindent
  {\bf 2.} Now we claim that no straight line
  with constant imaginary part of $p$,
  which does not contain any zero of $dp$,
  can connect $\Set{U}_{p}^+$ and $\Set{U}_{p}^-$.
  Let $s_0\in\mathbb{R}$ be the imaginary part of $p$ at
  one zero of $dp$, and consider the family of all those
  straight lines, whose constant imaginary part belongs to a small
  neighbourhood of $s_0$, and whose open ends, where the real
  part of $p$ goes to $-\infty$, belong either to $\Set{U}_{p}^+$ or
  $\Set{U}_{p}^-$. Some of the zeroes of $dp$ have the property,
  that if we move  around them and return to the part
  with very small real parts of $p$ we may connect
  $\Set{U}_{p}^+$ with $\Set{U}_{p}^-$. Now all zeroes of
  $dp$ with this property, and which belong to one connected component
  of the complement of $\Set{U}_{p}^+\cup \Set{U}_{p}^-$,
  have to be permuted by the involution $\eta$. Hence the number of such
  zeroes has to be even, and therefore these straight lines cannot
  connect $\Set{U}_{p}^+$ and $\Set{U}_{p}^-$.

  \noindent
  {\bf 3.} Let us call a straight line with constant imaginary part of $p$
  periodic, if there exists some period $\Pi\in\mathbb{R}$, such that
  the preimages of $p+\Pi$ coincide with the preimages of $p$.
  If the imaginary part of $p$ varies in some interval,
  such that the corresponding periodic lines contain
  no zero of $dp$, then the corresponding periods are constant,
  since $dp$ is closed.
  Obviously, all straight lines with constant imaginary part of $p$
  are either periodic, or have open ends intersecting
  one of the two neighbourhoods  $\Set{U}_{p}^+$ or $\Set{U}_{p}^-$
  of the two poles of $dp$ at $\infty^{\pm}$. Hence for generic real $s$
  (i.\ e.\ for $s$ in the complement of the finite set
  of values of the imaginary part of $p$ at all zeroes of $dp$)
  there exist exactly two non--periodic straight lines with constant
  imaginary part of $p$ equals $s$, one having non--empty intersection
  with $\Set{U}_{p}^+$ and another one having non--empty intersection
  with $\Set{U}_{p}^-$. In particular, the Riemann surfaces decomposes
  into three types of regions. On one hand \Em{regions of periodic lines}
  \index{region of periodic lines}
  and on the other hand two regions of non--periodic lines intersecting
  either $\Set{U}^+$ or $\Set{U}^-$. Moreover, the common boundaries of
  different regions with periodic lines, are periodic lines, which
  contain a zero of $dp$. Furthermore, the common boundaries of
  regions with periodic lines with one of the two latter regions
  are of the same form. The common boundaries of the latter
  two regions, however, are \Em{horizontal cuts}.

  \noindent
  {\bf 4.} Now we claim that for all real $s$,
  the part of the complement of $\Set{U}_{p}^+\cup \Set{U}_{p}^-$,
  on which the imaginary part of $p$ is equal to $s$
  is a union of straight lines, whose sum of lengths is equal to the
  sum of those segments of the straight lines with imaginary part
  equals $s$, which are excluded from $\Set{U}_{p}^+$ and
  $\Set{U}_{p}^-$. Obviously this is true for very large
  (and very small) values of $s$. On every component of a
  \Em{region of periodic lines} the imaginary part of $p$ takes the
  maximum and the minimum either on a common boundary
  with $\Set{U}^+$ or with $\Set{U}^-$. Hence the claim follows
  from the description of the boundaries
  of the three types of regions described in part 3.
  Moreover, the same arguments show that the complement of all
  \Em{regions of periodic lines} is the disjoint union
  $\Set{U}^+\cup\Set{U}^-$.

  \noindent
  {\bf 5.} We have decomposed the Riemann surface $\Spa{Y}$
  into three regions:
  the neighbourhoods $\Set{U}_{p}^+$ and $\Set{U}_{p}^-$
  (or their closures) and the \Em{regions of periodic lines}.
  We have already seen that for any value of the imaginary part of $p$,
  the corresponding lines together form two copies of lines of
  $\mathbb{P}^1$, with constant imaginary part. Consequently we will
  now see that these parts may be glued in such a way, that
  they form the two copies $\mathbb{P}_{p}^+$ and $\mathbb{P}_{p}^-$.
  In doing so we have some freedom and no canonical choice. But, as we
  already have seen all connected components of the boundary of the
  closures of $\Set{U}_{p}^+$ and $\Set{U}_{p}^-$,
  which do not form a \Em{horizontal cut},
  are glued with some \Em{regions of periodic lines},
  and for each single \Em{region of periodic lines}
  (with only one period $\Pi$) the corresponding gluing rule
  is of the from $p\mapsto p+$ some real constant.
  This together with the condition that
  for all imaginary values of $p$ the sum of the
  lengths of the corresponding periodic lines is exactly equal to
  the excluded parts of those lines with open ends intersecting
  $\Set{U}_{p}^+$ or $\Set{U}_{p}^-$
  is the essential content of the statement
  that $\Spa{Y}$ may be realized by a combination of several
  \Em{horizontal cuts} and several \Em{pair of parallel cuts}.
  In order to make this more precise,
  let us describe these combinations in some more detail.

  \noindent
  {\bf 6.} If we try two realize $\Spa{Y}$ in general
  by a combination of several \Em{horizontal cuts} and several
  \Em{pairs of parallel cuts} we have to envisage the possibility that
  different \Em{pairs of parallel cuts} penetrate each other. In
  the following two figures we show two equivalent such
  possibilities, where all pairs of parallel slits, which differ by
  some real period, are glued:

\setlength{\unitlength}{1cm}
\begin{picture}(15,6)

\put(0,0){\begin{picture}(6,6)
\put(2,5){\line(1,-4){1}}
\put(4,5){\line(1,-4){1}}

\put(1,1){\line(1,2){1.33}}
\put(2,1){\line(1,2){.66}}

\put(5,5){\line(-1,-2){.66}}
\put(6,5){\line(-1,-2){1.33}}

\put(1,0){\makebox(3,.5)[t]{figure of $\mathbb{P}_{p}^+$}} 
\end{picture} }

\put(9,0){\begin{picture}(6,6)
\put(2,5){\line(0,-1){3}}
\put(4,5){\line(0,-1){1}}

\put(1,1){\line(1,1){4}}
\put(2,1){\line(1,1){4}}

\put(3,1){\line(0,1){1}}
\put(5,1){\line(0,1){3}}

\put(1,0){\makebox(3,.5)[t]{figure of $\mathbb{P}_{p}^+$}} 
\end{picture} }

\end{picture}

In such cases we introduce additional horizontal cuts as shown in the
following figures:

\setlength{\unitlength}{1cm}
\begin{picture}(15,6)

\put(0,0){\begin{picture}(6,6)
\put(2,5){\line(0,-1){2}}
\put(4,5){\line(0,-1){2}}

\put(3,1){\line(0,1){2}}
\put(5,1){\line(0,1){2}}

\put(1,1){\line(1,2){1}}
\put(2,1){\line(1,2){1}}

\put(5,5){\line(-1,-2){1}}
\put(6,5){\line(-1,-2){1}}

\put(2,3){\line(1,0){3}}
\put(2.5,3.1){\makebox(.5,.5)[b]{$A^+_1$}}
\put(4.5,3.1){\makebox(.5,.5)[b]{$B^+_1$}}
\put(2,2.4){\makebox(.5,.5)[b]{$B^+_2$}}
\put(3.5,2.4){\makebox(.5,.5)[b]{$A^+_2$}}

\put(1,0){\makebox(3,.5)[t]{figure of $\mathbb{P}_{p}^+$}} 
\end{picture} }

\put(9,0){\begin{picture}(6,6)
\put(4,5){\line(0,-1){2}}
\put(2,5){\line(0,-1){2}}

\put(3,1){\line(0,1){2}}
\put(1,1){\line(0,1){2}}

\put(5,1){\line(-1,2){1}}
\put(4,1){\line(-1,2){1}}

\put(1,5){\line(1,-2){1}}
\put(0,5){\line(1,-2){1}}

\put(1,3){\line(1,0){3}}
\put(2.5,3.1){\makebox(.5,.5)[b]{$A^-_1$}}
\put(1,3.1){\makebox(.5,.5)[b]{$B^-_1$}}
\put(3.5,2.4){\makebox(.5,.5)[b]{$B^-_2$}}
\put(1.5,2.4){\makebox(.5,.5)[b]{$A^-_2$}}

\put(1,0){\makebox(3,.5)[t]{figure of $\mathbb{P}_{p}^-$}} 
\end{picture} }

\end{picture}

Furthermore, we permit the situation,
where that the parallel cuts of \Em{pair of parallel cuts}
are two lines with some corners, and where one is
a shift of the other by some real period. A \Em{horizontal cut}
and a \Em{pair of parallel cuts} may also penetrate each other:
Either the horizontal cut belongs to the
\Em{regions of periodic lines}, or to the part
$\Set{U}_{p}^+\cup \Set{U}_{p}^-$
where these lines are not periodic.
In the following figures, we show two
\Em{pairs of parallel cuts} with corners
and two \Em{horizontal cuts},
one inside of the part of periodic lines and one outside.

\setlength{\unitlength}{1cm}
\begin{picture}(15,6)

\put(0,0){\begin{picture}(6,6)

\put(1,1){\line(1,1){1}}
\put(3,1){\line(1,1){1}}

\put(5,1){\line(-1,1){1}}
\put(6,1){\line(-1,1){1}}

\put(2,2){\line(1,3){1}}
\put(4,2){\line(1,3){1}}
\put(5,2){\line(1,3){1}}

\put(1,3){\line(1,0){1.33}}
\put(6,3){\line(-1,0){.66}}

\put(3,4){\line(-1,0){.33}}
\put(4,4){\line(1,0){.66}}

\put(1,0){\makebox(3,.5)[t]{figure of $\mathbb{P}_{p}^+$}} 
\end{picture} }

\put(9,0){\begin{picture}(6,6)

\put(5,1){\line(-1,1){1}}
\put(3,1){\line(-1,1){1}}

\put(1,1){\line(1,1){1}}
\put(0,1){\line(1,1){1}}

\put(4,2){\line(-1,3){1}}
\put(2,2){\line(-1,3){1}}
\put(1,2){\line(-1,3){1}}

\put(5,3){\line(-1,0){1.33}}
\put(0,3){\line(1,0){.66}}

\put(3,4){\line(1,0){.33}}
\put(2,4){\line(-1,0){.66}}

\put(1,0){\makebox(3,.5)[t]{figure of $\mathbb{P}_{p}^-$}} 
\end{picture} }

\end{picture}

As a result of the preceding considerations all Riemann surfaces $\Spa{Y}$
may be realized as a combination of several \Em{horizontal cuts} with
several \Em{pairs of parallel cuts} and some additional
horizontal cuts connecting pairs of parallel cuts.

\noindent
{\bf 7.} Let us assume now that we have given some data
$(\Spa{Y},\infty^-,\infty^+,k)\in\Bar{\moduli}_{g,\eta}$, and consider
the corresponding realization of $\Spa{Y}$ as two planes
$\mathbb{P}_{p}^{\pm}$ glued in the way described above.
The positive two-form
$$2\pi\sqrt{-1}\vol(\torus)
\left(dk_1+\sqrt{-1}dk_2\right)\wedge
\left(\overline{dk_1+\sqrt{-1}dk_2}\right)$$
over the part $\mathbb{P}_{p}^+$ is equal to
$$4\pi\vol(\torus)
\left(d\Re(k_1)\wedge d\Im(k_1)+d\Re(k_2)\wedge d\Im(k_2)
+d\Re(k_1)\wedge d\Re(k_2)+d\Im(k_1)\wedge d\Im(k_2)\right).$$
On $\mathbb{P}_{p}^+$ the function $k$ is single--valued, hence
this integral is equal to the integral of the one--form
$$4\pi\vol(\torus)
\left(\Re(k_1) d\Re(k_2)-\Im(k_2) d\Im(k_1)-
\Im(k_1) d\Re(k_1)-\Im(k_2) d\Re(k_2)\right)$$
over the boundary of $\mathbb{P}_{p}^+$. On the other hand, due to
Lemma~\ref{Willmore functional} the \Em{first integral} is equal to
the integral of the one--form
$$4\pi\vol(\torus)
\left( \Re(k_1) d\Re(k_2)+\Im(k_2) d\Im(k_1)
+\sqrt{-1} \Im(k_1) d\Re(k_2) -\sqrt{-1}\Im(k_2)
d\Re(k_1)\right)$$
over the boundary of $\mathbb{P}_{p}^+$. Hence it would suffice to
prove that the integral of the four one--forms 
$\Im(k_1) d\Re(k_1)$, $\Im(k_2) d\Re(k_2)$,
$\Im(k_2) d\Im(k_1)$,
$\Im(k_1) d\Re(k_2)$ and $\Im(k_2) d\Re(k_1)$
over the boundary of $\mathbb{P}_{p}^+$ all vanish.
Let us first prove this statement for the two basic cases of a 
\Em{horizontal cut} or a \Em{pair of parallel cuts}. In the case of a
\Em{horizontal cut} the boundary of $\mathbb{P}_{p}^+$ is a cycle of
$\Spa{Y}$, which is invariant under the involution $\eta$. The
transformation of the function $k$ under this involution
implies that in this case the integrals of the one--forms
$\Im(k_1) d\Re(k_1)$, $\Im(k_2) d\Re(k_2)$,
$\Im(k_1) d\Re(k_2)$ and
$\Im(k_2) d\Re(k_1)$  over the boundary of $\mathbb{P}_{p}^+$
vanish. Since on a horizontal cut the form $d\Im(k_1)$ vanishes, the
integral $\Im(k_2)d\Im(k_1)$ over the boundary of $\mathbb{P}_{p}^+$
also vanishes. If we consider a \Em{pair of parallel cuts} we shall
first remark that all forms
$\Im(k_1) d\Re(k_1)$, $\Im(k_2) d\Re(k_2)$,
$\Im(k_2) d\Im(k_1)$,
$\Im(k_1) d\Re(k_2)$ and $\Im(k_2) d\Re(k_1)$
are invariant under a translation of $(k_1,k_2)$ by any element of
$\mathbb{R}^2$. This implies that all the integrals
along the parallel cuts vanish. Moreover,
this argument also extends to a horizontal cut connecting
two parallel cuts, which cuts two parts of $\mathbb{P}_{p}^+$
into different parts. If the horizontal cut, which connects both
parallel cuts, then it is a boundary of $\mathbb{P}_{p}^+$ and
$\mathbb{P}_{p}^-$, it is the sum of two circles of $\Spa{Y}$,
which are interchanged by $\eta$. The same arguments as in the
case of a \Em{horizontal cut} show that also in this case all the
integrals of the forms
$\Im(k_1) d\Re(k_1)$, $\Im(k_2) d\Re(k_2)$,
$\Im(k_2) d\Im(k_1)$,
$\Im(k_1) d\Re(k_2)$ and $\Im(k_2) d\Re(k_1)$
along the horizontal cut connecting the two parallel cuts vanish.
These arguments obviously extend to the combination of several
\Em{horizontal cuts} with several \Em{pair of parallel cuts} of the
general case. Moreover, these arguments also show, that the
\Em{first integral} is equal to the integral of the positive form
$2\pi\sqrt{-1}\vol(\torus)
(dk_1+\sqrt{-1}dk_2)\wedge(\overline{dk_1+\sqrt{-1}dk_2})$
over $\Set{U}_{p}^+$ plus the integral of the positive form
$2\pi\sqrt{-1}\vol(\torus)
(dk_1\wedge d\Bar{k}_1+dk_2\wedge d\Bar{k_2})$ over the
\Em{subregion of periodic lines} of $\mathbb{P}_{p}^+$.
The second summand is larger or equal to
$4\pi\vol(\torus)$ times the volume of this
\Em{subregion of periodic lines} of $\mathbb{P}_{p}^+$.
In particular, the \Em{first integral} is larger or equal to
the volume of this part.

\noindent
{\bf 8.} We conclude that the \Em{first integral} is a non--negative real
number. It is zero, if and only if $dk_2/dk_1$ is identically equal to
$\sqrt{-1}$ on $\mathbb{P}_{k_1}^+$. In this case $\Spa{Y}$ has two connected
components, on one $dk_2/dk_1$ is equal to $\sqrt{-1}$ and on the other
equal to $-\sqrt{-1}$. Since all periods of $dk_1$ and $dk_2$ are real,
this implies that they have to be zero in this case. Hence $k_1$ and
$k_2$ may be chosen to be single--valued meromorphic functions with one
pole on each connected component of $\Spa{Y}$. Therefore, both connected
components have to be biholomorphic to $\mathbb{P}^1$. This proves the
Lemma.
\end{proof}

For all data $(\Spa{Y},\infty^-,\infty^+,\eta,p)$ let
$D(\Spa{Y},\infty^-,\infty^+,\eta,p)$ denote the maximal diameters of the
compact convex subsets mentioned in \Em{Neighbourhood}~(i), which
form the complement of the image of $\Set{U}_{p}^+$.
Since the diameter is invariant under translations,
it depends only on the form $dp$.

\begin{Lemma} \label{Willmore bound}
For all $g\in\mathbb{N}_0$ and $\willmore>0$ there exists some $C>0$,
such that for all
$(\Spa{Y},\infty^-,\infty^+,\eta,k)\in
\Bar{\moduli}_{g,\lattice,\eta,\willmore}$
the diameters $D(\Spa{Y},\infty^-,\infty^+,\eta,\xx{p})$ and
$D(\Spa{Y},\infty^-,\infty^+,\eta,\yy{p})$ are not larger than $C$.
\end{Lemma}

\begin{proof} Let us assume on the contrary that there exists some
  sequence of data in $\Bar{\moduli}_{g,\lattice,\eta,\willmore}$,
  such that either the diameter of the corresponding data
  $(\Spa{Y},\infty^-,\infty^+,\eta,\xx{p})$ or the diameter of the data 
  $(\Spa{Y},\infty^-,\infty^+,\eta,\yy{p})$ tends to infinity. By passing to
  some subsequence we may assume that for all elements of the
  subsequence the second diameter is larger than the former diameter
  and that the second sequence of diameters is monotone increasing. We
  \Em{renormalize} the functions $k$ and obtain some sequence of data
  in $\Bar{\moduli}_{g,\eta}$, whose corresponding data 
  $(\Spa{Y},\infty^-,\infty^+,\eta,\yy{p})$ have constant diameter equal
  to $1$, and whose \Em{first integral} converges to zero. Some of
  the excluded domains of $\Set{U}_{\yy{p}}^+$ corresponding to the
  sequence of data $(\Spa{Y},\infty^-,\infty^+,\eta,\yy{p})$
  tend to infinity. Now we use the invariance under the shifts
  described in Lemma~\ref{covariant transformation} to transform the
  sequence into another sequence in
  $\Bar{\moduli}_{g,\lattice,\eta,\willmore}$
  of data, such that the \Em{renormalized} data have excluded domains of
  $\Set{U}_{\yy{p}}^+$ corresponding to the \Em{renormalized}
  $(\Spa{Y},\infty^-,\infty^+,\eta,\yy{p})$ inside of some finite subset of
  $\mathbb{C}$ with diameter equal to $1$.

  In the proof of Lemma~\ref{gluing rule} we showed that
  the \Em{first integral} is larger or equal to
  $4\pi\vol(\torus)$ times the
  volume of the \Em{regions of periodic lines} of $\mathbb{P}_{p}^+$.
  Since the \Em{first integral} is invariant under conformal linear
  transformations of $\lattice$, in general the \Em{first integral} is
  larger or equal to
  $4\pi\vol(\torus)/g(\gamma,\gamma)$ times
  the volume of the \Em{regions of periodic lines} of
  $\mathbb{P}_{p}^+$ with respect to the component $p=g(\gamma,k)$.
  Since the periods of $\xx{p}$ and $\yy{p}$ are integers, we
  conclude that the \Em{regions of periodic lines} of
  $\mathbb{P}_{\xx{p}}^+$ and $\mathbb{P}_{\yy{p}}^+$ of the
  \Em{renormalized} sequence converge to finitely many
  horizontal lines, respectively. Now we claim that on
  all compact subsets in the complement of these finitely many
  horizontal lines of $\mathbb{P}_{\xx{p}}$ the function
  $k_1+\sqrt{-1}k_2$ converges uniformly to zero. In fact, due to
  Cauchy's formula \cite[Chapter~IV. \S6.1]{Co1} the value of this
  function at any element $\xx{p}'$ in the image of $\Set{U}_{\xx{p}}^+$
  under the function $\xx{p}$ is equal to
  $$\frac{1}{2\pi\sqrt{-1}}\int\limits_{\partial \Set{U}_{\xx{p}}^+}
  (k_1+\sqrt{-1}k_2)\frac{d\xx{p}}{\xx{p}-\xx{p}'}=
  -\frac{1}{2\pi\sqrt{-1}}\int\limits_{\partial \Set{U}_{\xx{p}}^+}
  \ln(1-\xx{p}/\xx{p}')(dk_1+\sqrt{-1}dk_2).$$
  We remark that this formula is only valid, if $\xx{p}'$ is
  not contained inside of any of the cycles corresponding to the
  component of the boundary of $\Set{U}_{\xx{p}}^+$
  (in other words the formula is valid,
  if there exists a path not crossing the boundary of
  $\Set{U}_{\xx{p}}^+$, which connects $\xx{p}'$ with
  $\infty^+$). This formula yields the branch of the function
  $k_1+\sqrt{-1}k_2$, which vanishes at $\infty^+$.
  In particular, if we use the
  left hand side to calculate the contribution of some component of
  the boundary of $\Set{U}_{\xx{p}}^+$,
  this contribution does not depend
  on some constant, which is added to the function $k_1+\sqrt{-1}k_2$
  nearby this boundary component. In fact, due to Lemma~\ref{gluing rule}
  the function $k_1+\sqrt{-1}k_2$ extends to some
  single--valued function on $\Set{U}_{\xx{p}}^+$.
  On the other hand we may estimate
  $|(dk_1+\sqrt{-1}dk_2)/d\xx{p}|$ at every point
  $\xx{p}'$ in the image of $\Set{U}_{\xx{p}}^+$
  under $\xx{p}$, whose distance to the boundary of
  $\Set{U}_{\xx{p}}^+$ is larger than some
  $\varepsilon>0$ by some constant (depending only on the lattice) times
  $\willmore/\varepsilon$. We conclude that the function
  $k_1+\sqrt{-1}k_2$ converges to zero uniformly on compact subsets of
  the complement of the horizontal straight lines.
 
  We conclude that all accumulation points of zeroes of
  $d\yy{p}$ cannot belong to the interior of
  $\mathbb{P}_{\xx{p}}^+$ corresponding to the limit. Hence also the
  maximal difference of the real parts of the domains excluded from
  $U_{\yy{p}}^+$ corresponding to the sequences
  $(\Spa{Y},\infty^-,\infty^+,\eta,\yy{p})$ converges to zero, which
  contradicts the assumptions. Therefore, the diameters of the
  elements of $\Bar{\moduli}_{g,\lattice,\eta,\willmore}$ are
  bounded.
\end{proof}

\begin{Lemma} \label{bounded diameter}
For all fixed $g\in\mathbb{N}_0$, $C>0$ and $\willmore>0$ and any sequence
of elements $(\Spa{Y},\infty^-,\infty^+,k)\in
\Bar{\moduli}_{g,\lattice,\eta,\willmore}$ whose diameters
of the corresponding sequences of members
$D(\Spa{Y},\infty^-,\infty^+,\eta,\xx{p})$ and
$D(\Spa{Y},\infty^-,\infty^+,\eta,\yy{p})$ are not larger than $C$,
there exists some subsequence and some choice
of the corresponding cuts of
$\mathbb{P}_{\xx{p}}^+\setminus \{\infty^+\}$ described
in Lemma~\ref{gluing rule}, which converges in the sense of
\cite[\S15.4 Definition~4.1]{Co2} to some choice of
$\mathbb{P}_{\xx{p}}^+\setminus \{\infty^+\} $ of
some data $(\Spa{Y},\infty^-,\infty^+,\eta,\xx{p})$. Moreover, the function
$\yy{p}(\xx{p})$ of this subsequence converges uniformly on
compact subsets to the corresponding functions of the limit. Finally, the
\Em{first integral} of the limit is not larger than the lim inf of the
subsequence of \Em{first integrals}.
\end{Lemma}

\begin{proof} Let us choose some constant $C>0$ and consider all
data $(\Spa{Y},\infty^-,\infty^+,\eta,k)\in
\Bar{\moduli}_{g,\lattice,\eta,\willmore}$ with the
property $D(\Spa{Y},\infty^-,\infty^+,\eta,\xx{p})\leq C$. 
Let us first prove that any sequence of such data has a subsequence
and some choice of cuts described in Lemma~\ref{gluing rule} with the
property, that the corresponding domains
$\mathbb{P}_{\xx{p}}^+\setminus \{\infty^+\}$ converge in
the sense of \cite[\S15.4 Definition~4.1]{Co2} to some choice of cuts
corresponding to some data in this space. Obviously we may choose some
subsequence such that the images of the sets
$\Set{U}_{\xx{p}}^+\setminus \{\infty^+\}$ under some choice of the
branches of the function $\xx{p}$ converge to some open domain in
$\mathbb{C}$. Moreover, by passing to some subsequence we may achieve
that the imaginary values of $\xx{p}$ at all zeroes of $d\xx{p}$
converge and moreover that if $s$ is not the limit of some sequence
of imaginary parts of zeroes of $d\xx{p}$, that all periods of
periodic straight lines with imaginary part
equal to $s$ converge to some finite set of periods. By passing again
to some subsequence we may achieve in addition
that the corresponding sequence of images of $\Set{U}_{\xx{p}}^+$
under $\xx{p}$  converge to the corresponding image of some data
in this set. From the considerations
in the proof of Lemma~\ref{gluing rule} it follows that we may chose
the cuts of these subsequence in such a way, that the corresponding
domains $\mathbb{P}_{\xx{p}}^+\setminus \{\infty^+\}$ converge in
the sense of \cite[\S15.4 Definition~4.1]{Co2} to some choice of cuts
corresponding to some data in this space.

Now we claim that there exists some constant $C'>0$ such that the
mapping $\xx{p}\mapsto\yy{p}$ defines some biholomorphic mapping
from the subset $\left\{\xx{p}\in\mathbb{C}\mid
B(\xx{p},C')\subset \xx{p}[\Set{U}_{\xx{p}}^+]\right\}$
into some domain of $\mathbb{C}$.
This follows from the formula for the function
$k_1+\sqrt{-1}k_2$ used in the proof of Lemma~\ref{Willmore bound}.
We conclude that this biholomorphic mapping
$\xx{p}\mapsto\yy{p}$ maps some complement in
$\mathbb{C}$ of finitely many compact domains into a complement of
finitely many compact domains. More generally, the
\Em{complex Fermi curve} is therefore the union of finitely many
compact pieces and two open ends, which are mapped by
$\xx{p}$ and $\yy{p}$ onto the complement in
$\mathbb{C}$ of finitely many disjoint compact subsets.
If we consider the \Em{complex Fermi curve} as a
finite--sheeted covering with respect to the mapping $\yy{p}$,
then we may also apply
Lemma~\ref{gluing rule} to this covering. We conclude that the
four functions
$\exp(2\pi\sqrt{-1}\xx{p})$, $\exp(-2\pi\sqrt{-1}\xx{p})$,
$\exp(2\pi\sqrt{-1}\yy{p})$ and $\exp(-2\pi\sqrt{-1}\yy{p})$
are bounded by some constant. Moreover, this argument carries over to
some subsequence of any sequence in
$\Bar{\moduli}_{g,\lattice,\eta,\willmore}$
with the properties described in the Lemma.
With the help of Montel's theorem \cite[Chapter~VII \S2.9]{Co1} we
conclude that there exists some subsequence, such that the functions
$\yy{p}(\xx{p})$ converge uniformly on compact subsets of
$\mathbb{P}_{\xx{p}}^+\setminus \{\infty^+\}$
to some function, which has to be equal to the corresponding function
of the limit of the data.

The single--valued function $dk_2/dk_1$ has no poles on the
interior of $\mathbb{P}_{\xx{p}}^+$ and is equal to
$\sqrt{-1}$ at $\infty^+$. We conclude that for all data in
$\Bar{\moduli}_{g,\eta}$ the function
$(dk_2/dk_1-\sqrt{-1})$ converges to zero, if the absolute value of
$\xx{p}\in\mathbb{P}_{\xx{p}}^+$ converges to infinity. Since
$dk_1+\sqrt{-1}dk_2= \sqrt{-1}((dk_2/dk_1-\sqrt{-1})dk_1$ we conclude
that for all data in $\Bar{\moduli}_{g,\eta}$ and all $\varepsilon>0$
there exists some compact subset $K$ of some choice of
$\mathbb{P}_{\xx{p}}^+\setminus\{\infty^+\}$
corresponding to this data, such that the integral of the form
$$2\pi\sqrt{-1}\vol(\torus)
(dk_1+\sqrt{-1}dk_2)\wedge(\overline{dk_1+\sqrt{-1}dk_2})$$
over $K$ is an element of the interval
$[\willmore-\varepsilon,\willmore]$, where $\willmore$ denotes the
\Em{first integral} of the data. Since $\yy{p}(\xx{p})$ converges
uniformly on compact subsets of $\mathbb{P}_{\xx{p}}^+$, we conclude
that the \Em{first integral} of the limit is smaller than the sum of
any accumulation point of the sequence of values of the
\Em{first integral} plus $\varepsilon$. Since
this is true for all $\varepsilon$, the \Em{first integral} of the limit
has to be smaller or equal to the lim inf of the sequence of values
of the \Em{first integral}.
\end{proof}

It is possible to prove a sharper statement: For a given sequence
of data some of the excluded domains of corresponding images of
$\Set{U}_{p}^+$ under $p$, will tend to infinity.
To all these regions we may associate a part of
the \Em{first integral}, by integrating the corresponding one--form
around these excluded domains. In the limit the excluded domains,
which tend to infinity are removed, and therefore the
\Em{first integral} of the limit is the limit of those parts of the
\Em{first integrals}, which correspond to excluded domains inside
of bounded domains. In general, at least if the excluded domains
have enough distance from each other, all parts of the
\Em{first integral} are positive. Therefore, the limit of the
\Em{first integral} is smaller than all accumulation points of the
sequence of values of the \Em{first integrals} of any convergent
series of data.

\noindent
{\it Proof of Theorem~\ref{compactification}.} Due to
Lemma~\ref{compact metric} the space of all
closed subsets of $\mathbb{C}^2$ is a separable metrizable
space. Hence $\Bar{\moduli}_{g,\lattice,\willmore}$ is compact
if and only if any sequence has a convergent subsequence. Due to the
Lemma~\ref{bounded diameter} and Lemma~\ref{Willmore bound}
any sequence of $\Bar{\moduli}_{g,\lattice,\willmore}$ has a
subsequence, which converges in sense described above. It is
easy to see that that this implies  that the corresponding subsets
$\fermi(\Spa{Y},\infty^-,\infty^+,k)\subset\mathbb{C}^2$ converges.
Moreover, due to Lemma~\ref{bounded diameter}
the \Em{first integral} of the limit is not larger than
$\willmore$.\qed

\subsection{Limits of \Em{complex Fermi curves}}
\label{subsection limits}

In this section we shall investigate the limits of the sequences of
resolvents corresponding to weakly convergent
sequences of complex potentials $\left(U_n\right)_{n\in\mathbb{N}}$
in $\banach{2}(\torus)$ with bounded $\banach{2}$--norm.

Due to the Banach--Alaoglu theorem \cite[Theorem~IV.21]{RS1} and the
Riesz Representation theorem \cite[Chapter~13 Section~5]{Ro2}
any sequence $\left(U_n\right)_{n\in\mathbb{N}}$ in
$\banach{2}(\torus)$ with bounded $\banach{2}$--norm has a subsequence,
with the property that the corresponding
sequence of measures $U_n\Bar{U}_nd^2x$
converge weakly to a measure on $\mathbb{R}^2/\lattice$. With the help
of Lemma~\ref{weakly continuous resolvent} we conclude that the
corresponding sequence of resolvents converge, if the limit
of the measure does not contain point measures of mass
larger or equal to the constant $S_p^{-2}$.
In fact, if the weak limit of the measures
does not contain point measures of mass larger or equal to $S_p^{-2}$,
we may cover $\torus$ by open sets, whose measures with
respect to the limit of the measures is smaller than $S_p^{-2}$.
Due to the compactness of $\torus$ this open covering has a finite
subcovering. For any finite open covering, the function on
$\torus$, which associates to each $x$ the radius of the
maximal open disc around $x$, which is entirely contained in one member
of the covering, is continuous. We exclude the trivial case,
where one member of the covering contains the whole torus and
therefore all discs. So this function is the maximum of the
distances of the corresponding point to all complements of the
members of the covering. Therefore, there exists a small
$\varepsilon>0$, such that all discs with radius $2\varepsilon$
are contained in one member of the finite subcovering.
Obviously, for any member of the subcovering
there exists a continuous $[0,1]$--valued function,
whose support is contained in this member of the subcovering,
and which is equal to $1$ on those discs $B(x,\varepsilon)$,
whose extensions $B(x,2\varepsilon)$ are contained
in this member of the subcovering.
Since the sequence of measures converges weakly,
the integrals of these functions with respect to
the measures $U_n(x)\Bar{U}_n(x)d^2x$ corresponding to the sequence
are also smaller than $S_p^{-2}$,
with the exception of finitely many elements of the sequence.
This shows that with the exception of
finitely many elements of the sequence,
the norms $\|U_n\|_{2,\varepsilon}$ are smaller than some $C_p<S_p^{-1}$.

\begin{Lemma}\label{bounded point measures}
If a weak limit of the sequence of finite Baire measures
$U_n\Bar{U}_nd^2x$ on $\mathbb{R}/\lattice$ does not contain point
measures with mass larger or equal to $S_p^{-2}$,
then there exists a $C_p<S_p^{-1}$, an $\varepsilon>0$
and a subsequence of the bounded sequence $U_n$ in $\banach{2}(\torus)$
whose norms $\|\cdot\|_{2,\varepsilon}$ are smaller than $C_p$.
\qed
\end{Lemma}

Since the $\banach{2}$--norm of the sequence is
bounded, there exist only finitely many points
$z_1,\ldots, z_L$ in $\mathbb{C}/\lattice_\mathbb{C}$,
whose mass with respect to the weak limit of the measures
$U_n(x)\Bar{U}_n(x)d^2x$ is larger or equal to $S_p^{-2}$.
Here we used the complex coordinates $z$ instead of
the real coordinates $x$.
Moreover, since we consider only local potentials
(i.\ e.\ multiplication by functions) the limits of the
corresponding sequence of resolvents should yield perturbations of the
resolvents of the Dirac operator corresponding to the weak limit $U$.
Furthermore, we expect that the support of these perturbations is
contained in this finite set $\{z_1,\ldots,z_L\}$.
In the first subsection we characterize the natural candidates
for these perturbations. In the last subsection we will show
that this is indeed the case.
The proof is intricate and will proceed through several steps.

\subsubsection{Finite rank perturbations}
\label{subsubsection finite rank perturbations}

It will turn out that in general the limits of the resolvents
$\Op{R}(U_n,\Bar{U}_n,k,\lambda)$ depend non--continuously on
$\lambda$. Since we are interested in the limits of the
\Em{complex Fermi curves} we shall concentrate on the limits of
$\Op{R}(U_n,\Bar{U}_n,k,0)$. As a preparation we first derive an
explicit formula for the integral kernels of $\Op{R}(0,0,k,0)$.

These integral kernels may be expressed in terms of a
Theta function corresponding to the elliptic curve
$\mathbb{C}/\lattice_\mathbb{C}$
together with a canonical homology bases. The choice of a
\De{Fundamental domain}~\ref{fundamental domain} determines these
data, therefore we decorate the corresponding Theta function with an
index $\Delta$. Let $\theta_{\Delta}$ denote the
following Theta function:
\index{Theta function $\theta_{\Delta}$}
$$\theta_{\Delta}(z)=\frac{\xx{\gamma}_1+\sqrt{-1}\xx{\gamma}_2}
{\theta'_1(0,
\frac{\yy{\gamma}_1+\sqrt{-1}\yy{\gamma}_2}
     {\xx{\gamma}_1+\sqrt{-1}\xx{\gamma}_2})}
\theta_1\left(\frac{z}{\xx{\gamma}_1+\sqrt{-1}\xx{\gamma}_2},
\frac{\yy{\gamma}_1+\sqrt{-1}\yy{\gamma}_2}
     {\xx{\gamma}_1+\sqrt{-1}\xx{\gamma}_2}\right).$$
Here $\theta_1(z,\tau)$ is the first one--dimensional Theta--function
with modulus $\tau$ \cite[13.19]{MOT}. Since the
orientation of the basis $\xx{\gamma}$ and $\yy{\gamma}$ is
positive, the imaginary part of the modulus
$\frac{\yy{\gamma}_1+\sqrt{-1}\yy{\gamma}_2}
      {\xx{\gamma}_1+\sqrt{-1}\xx{\gamma}_2}$ is
positive and so the one--dimensional Theta--functions are well defined. The
properties of these functions implies that the functions
$\theta_{\Delta}$ are uniquely determined by the following properties:
\begin{description}
\item[Quasi--periodicity:] Under the shifts of the elements of $\lattice$
  considered as element in $\mathbb{C}$ this function transform as
  follows:
\begin{eqnarray*}
\theta_{\Delta}(z+\xx{\gamma}_1+\sqrt{-1}\xx{\gamma}_2)
&=&-\theta_{\Delta}(z)\\
\theta_{\Delta}(z+\yy{\gamma}_1+\sqrt{-1}\yy{\gamma}_2)
&=&-\exp\left(-\pi\sqrt{-1}\frac{2z+\yy{\gamma}_1+\sqrt{-1}\yy{\gamma}_2}
{\xx{\gamma}_1+\sqrt{-1}\xx{\gamma}_2}\right)\theta_{\Delta}(z).
\end{eqnarray*}
\item[Holomorphicity:] The function $z\mapsto\theta_{\Delta}(z)$ is
  a holomorphic entire function.
\item[Anti--symmetry:] Under the map $z\mapsto-z$ this function is
  anti--symmetric.
\item[Zeroes:] The zeroes of this function are exactly the elements
  of $\lattice_\mathbb{C}$. More precisely, at
  the origin it has an expansion of the form
  $\theta_{\Delta}(z)=z+\text{O}(z^3)$.
\end{description}
Now Dolbeault's Lemma 
\cite[Chapter~I Section~D 2.~Lemma]{GuRo} implies that the inverse of
the operators $\Bar{\partial}_{[k]}$ and $\partial_{[k]}$ have the
following integral kernels:
\begin{eqnarray*}
\Func{K}_1(z-z',k)\frac{d\Bar{z}'\wedge dz'}{2\pi\sqrt{-1}} &=&
\exp\left(2\pi\sqrt{-1}g(\xx{\gamma},k)
\frac{z-z'}{\xx{\gamma}_1+\sqrt{-1}\xx{\gamma}_2}\right)
\frac{\theta_{\Delta}(z-z'+z^+(k))}
{\theta_{\Delta}(z-z')\theta_{\Delta}(z^+(k))}
\frac{d\Bar{z}'\wedge dz'}{2\pi\sqrt{-1}}\\
\Func{K}_2(z-z',k)\frac{d\Bar{z}'\wedge dz'}{2\pi\sqrt{-1}} &=&
\exp\left(2\pi\sqrt{-1}g(\xx{\gamma},k)
\frac{\Bar{z}-\Bar{z}'}{\xx{\gamma}_1-\sqrt{-1}\xx{\gamma}_2}\right)
\frac{\Bar{\theta}_{\Delta}(z-z'+\Bar{z}^-(k))}
{\Bar{\theta}_{\Delta}(z-z')\Bar{\theta}_{\Delta}(\Bar{z}^-(k))}
\frac{d\Bar{z}'\wedge dz'}{2\pi\sqrt{-1}},
\end{eqnarray*}
where $z^+(k)=\vol(\torus)(\sqrt{-1}k_1-k_2)$
and $z^-(k)=\vol(\torus)(\sqrt{-1}k_1+k_2)$.
In fact, this choice implies the relations
\begin{eqnarray*}
\exp\left(2\pi\sqrt{-1}g(\xx{\gamma},k)
\frac{\yy{\gamma}_1+\sqrt{-1}\yy{\gamma}_2}
     {\xx{\gamma}_1+\sqrt{-1}\xx{\gamma}_2}\right)
\exp\left(-\pi\sqrt{-1}
\frac{2z^+(k)}{\xx{\gamma}_1+\sqrt{-1}\xx{\gamma}_2}\right)
&=&\exp(2\pi\sqrt{-1}g(\yy{\gamma},k))\\
\exp\left(2\pi\sqrt{-1}g(\xx{\gamma},k)
\frac{\yy{\gamma}_1-\sqrt{-1}\yy{\gamma}_2}
     {\xx{\gamma}_1-\sqrt{-1}\xx{\gamma}_2}\right)
\exp\left(\pi\sqrt{-1}
\frac{2z^-(k)}{\xx{\gamma}_1-\sqrt{-1}\xx{\gamma}_2}\right)
&=&\exp(2\pi\sqrt{-1}g(\yy{\gamma},k)).
\end{eqnarray*}
Hence both integral kernels define with respect to $z$
sections of the line bundle corresponding to $[k]$ and with respect to
$z'$ sections of the line bundle corresponding to $[-k]$.
The first integral kernel is well defined if and only
$z^+(k)$ is not an element of $\lattice_\mathbb{C}$ and the
second integral kernel is well defined if and
only if $\Bar{z}^-(k)$ is not an element of $\lattice_\mathbb{C}$.
The generators of the dual lattice are equal
to $\xx{\kappa}=\Op{J}\yy{\gamma}/\vol
(\torus)$ and $\yy{\kappa}=\Op{J}\xx{\gamma}/
\vol(\torus)$.
Therefore, the functions $z^{\pm}(k)$ obey the relations
\begin{align*}
z^+(k+\xx{\kappa})&=
z^+(k)+\yy{\gamma}_1+\sqrt{-1}\yy{\gamma}_2 &
z^+(k+\yy{\kappa})&=
z^+(k)+\xx{\gamma}_1+\sqrt{-1}\xx{\gamma}_2\\
z^-(k+\xx{\kappa})&=
z^-(k)-\overline{\yy{\gamma}_1+\sqrt{-1}\yy{\gamma}_2} &
z^-(k+\yy{\kappa})&=
z^-(k)-\overline{\xx{\gamma}_1+\sqrt{-1}\xx{\gamma}_2}.
\end{align*}
Hence the first condition is equivalent to the condition
that the line bundle corresponding to $[k]$ has a nontrivial
holomorphic section and $\Bar{\partial}_{[k]}$ has a kernel and the
second condition is equivalent to the condition that
this line bundle has a anti--holomorphic section and $\partial_{[k]}$ has
a kernel. More generally, these integrals kernels considered as
functions of $k$ are invariant under shifts of the dual lattice, and
therefore depend only on $[k]\in\mathbb{C}^2/\lattice\dual$.
We collect these calculations in the following

\begin{Lemma}\label{resolvent integral kernel}
\index{integral!kernel!of $\Op{R}(0,0,k,0)$}
The resolvent $\Op{R}(0,0,k,0)$ has the following integral kernel:
$$\begin{pmatrix}
0 & \Func{K}_1(z-z',k)\\
-\Func{K}_2(z-z',k) & 0
\end{pmatrix}\frac{d\Bar{z}'\wedge dz'}{2\pi\sqrt{-1}}=
\begin{pmatrix}
0 & \Func{K}_1(z-z',k)\\
-\Func{K}_2(z-z',k) & 0
\end{pmatrix}\frac{d^2x'}{\pi}.$$\qed
\end{Lemma}

The following properties characterize the natural candidates
of those perturbations, which occur as limits of sequences of
resolvents.

\begin{description}
\index{finite!rank perturbation!condition $\sim$ (i)--(iii)}
\index{condition!finite rank perturbation (i)--(iii)}
\item[Finite rank Perturbation (i)]
The resolvents fit together to a meromorphic family of operators
depending on $k\in\mathbb{C}^2$, which differs from the family of
resolvents $k\mapsto\Op{R}(V,W,k,0)$ by a family of finite rank operators.
\item[Finite rank Perturbation (ii)]
The application of the Dirac operators $\Op{D}(V,W,k)$ on all elements
of the range of the finite rank operators in (i)
are tempered distributions,
whose support is contained in a finite set $\{z_1,\ldots,z_L\}$.
The same holds for the corresponding transposed operators.
The support of these distributions is called the support of
the \Em{Finite rank Perturbation}.
\item[Finite rank Perturbation (iii)]
The corresponding perturbations of the Dirac operators
$\Op{D}(V,W,k)$ does not depend on $k$.
\end{description}
Since any  tempered distribution,
which is located in finitely many points, is a sum of derivatives of
$\delta$--functions \cite[Theorem~V.11]{RS1}
condition \Em{Finite rank Perturbation}~(ii) implies
that the difference of the resolvents minus
$\Op{R}(V,W,k,0)$ is an operator with integral kernel
$$\sum\limits_{m,n}\psi_m(z,k)\Mat{S}(k)_{m,n}\phi_n(z',k)d^2x',$$
where $\Mat{S}(k)$ is a $L\times L$--matrix--valued function
depending on $[k]\in\mathbb{C}^2/\lattice\dual$,
and $\psi_1,\ldots,\psi_K$ and $\phi_1,\ldots,\phi_K$ have the
property that the applications of
$\Op{D}(V,W,k)$ and $\Op{D}^{t}(V,W,k)$ yields derivatives of
$\delta$--functions, respectively. Therefore, these functions are
values of the derivatives of the integral kernel at the support of the
\Em{Finite rank Perturbation}. In the following ten steps
we will make this more precise. The first three steps
can be considered as an application of Krein's formula \cite{AK}.
However, these perturbations violate the general rule,
that the domain of the perturbation should contain the
domain of the unperturbed operator
(\cite[Chapter four]{Ka} and \cite[Chapter 1]{AK}).

\noindent
{\bf 1.} Let us first characterize all \Em{Finite rank Perturbations}
of $\Op{R}(0,0,k,0)$, whose support is just one point $z_l$.
All these perturbations have integral kernels of the form
$$\Psi(z,k)\Mat{S}(k)\Phi(z',k)d^2x',$$
where $\Mat{S}(k)$ is a $M\times M$ matrix depending on
$k\in\mathbb{C}^2$, whose entries are
$2\times 2$--matrices, and $\Psi(z,k)$ and $\Phi(z,k)$
are the following row--vectors and column--vectors with $M$ entries
indexed by $m=0,\ldots,M-1$, respectively:
\begin{align*}
\Psi_m(z,k)&=\frac{(-1)^m}{\pi m!}\left(\begin{smallmatrix}
0 & \Func{K}_1^{(m)}(z-z_l,k)\\
-\Func{K}_2^{(m)}(z-z_l,k) & 0
\end{smallmatrix}\right) &
\Phi_m(z,k)&=\frac{1}{\pi m!}\left(\begin{smallmatrix}
0 & \Func{K}_1^{(m)}(z_l-z,k)\\
-\Func{K}_2^{(m)}(z_l-z,k) & 0
\end{smallmatrix}\right),
\end{align*}
where $\Func{K}_1^{(m)}(z,k)$ and $\Func{K}_2^{(m)}(z,k)$
denotes the $m$--th holomorphic and anti--holomorphic derivative
with respect to $z$  of the functions $\Func{K}_1$ and $\Func{K}_2$,
respectively. This ansatz should be understood as the
$2\times 2$ matrix, which is obtained by the matrix multiplication
of the transposed vector $\Psi$ with $L$ entries times the
$L\times L$ matrix $\Mat{S}(k)$ times the vector $\Phi$
with again $L$ entries, all whose entries are $2\times 2$ matrices.
Formally it would be also possible to consider
anti--holomorphic derivatives of $\Func{K}_1$ or holomorphic
derivatives of $\Func{K}_2$. But these perturbations
have integral kernels, which are supported only at the points
$z_1,\ldots,z_L$. We shall not consider such highly singular
perturbations, since we are only interested in perturbations of the
resolvents, which can be obtained as the limit
$\varepsilon\downarrow 0$ of perturbations, whose integral kernel are
supported in the complement of
$B(z_1,\varepsilon)\cup\ldots\cup B(z_L,\varepsilon)$.

Formally the corresponding perturbations $\Op{P}$ of
$\Op{D}(0,0,k,0)$ should have integral kernels of the form
\begin{align*}
&\Breve{\Psi}(z)\Mat{P}\Breve{\Phi}(z')d^2x'
&\text{ with}\\
\Breve{\Psi}_m&=\frac{(-1)^m}{m!}\left(\begin{smallmatrix}
\Bar{\partial}^m\delta(z-z_l) & 0\\
0 & \partial^m\delta(z-z_l)
\end{smallmatrix}\right) & 
\Breve{\Phi}_m&=\frac{(-1)^m}{m!}\left(\begin{smallmatrix}
\partial^m\delta(z-z_l) & 0\\
0 & \Bar{\partial}^m\delta(z-z_l)
\end{smallmatrix}\right).\end{align*}
Again $\Breve{\Psi}$ and $\Breve{\Phi}$ are row--vectors and
column--vectors, respectively,
with the given entries indexed by $m=0,\ldots,M-1$.
They are related to $\Psi$ and $\Phi$:
\begin{align*}
\Psi(z,k)=&\Op{R}(0,0,k,0)\Breve{\Psi}(z)&
\Phi^{t}(z,k)=&\Op{R}^{t}(0,0,k,0)\Breve{\Phi}^{t}(z).
\end{align*}
Since the resolvent of the perturbed operator is equal to
$$\Op{R}(0,0,k,0)\comp
\left(\unity-\Op{P}\comp\Op{R}(0,0,k,0)\right)^{-1},$$
we should calculate the inverse operator on the right hand side.

\begin{Lemma}\label{inverse operator}
Let $\Spa{F}$ be a Banach subspace of the Banach space
$\Spa{E}$ and $\Op{A}$ a linear operator on $\Spa{E}$,
whose range is contained in $\Spa{F}$. Then the operator
$\unity_{\Spa{E}}-\Op{A}$ on $\Spa{E}$ has an inverse
if and only if the operator
$\unity_{\Spa{F}}-\Op{A}_{\Spa{F}}$ on $\Spa{F}$
has an inverse, where $\Op{A}_{\Spa{F}}$ denotes the
restriction of $\Op{A}$ to $\Spa{F}$ considered as an
operator on $\Spa{F}$. Finally, if the inverse of the latter exists,
then the inverse of the former is equal to
$$\left(\unity_{\Spa{E}}-\Op{A}\right)^{-1}=
\unity_{\Spa{E}}+
\left(\unity_{\Spa{F}}-\Op{A}_{\Spa{F}}\right)^{-1}\comp
\Op{A}.$$
\end{Lemma}

In order to simplify notation we do not distinguish between operators
with the  same range, but which are mappings to different Banach spaces. 
Therefore, $\Op{A}$ denotes either an operator on $\Spa{E}$ or
an operator from $\Spa{E}$ to $\Spa{F}$. Furthermore, the
operators $\unity_{\Spa{F}}$ and $\Op{A}_{\Spa{F}}$
denote either operators on $\Spa{F}$ or operators from
$\Spa{F}$ to $\Spa{E}$.

\begin{proof}
Since the range of $\Op{A}$ is contained in $\Spa{F}$, the operator
$\unity_{\Spa{E}}-\Op{A}$ leaves the subspace $\Spa{F}$ invariant.
If the inverse of
$\unity_{\Spa{E}}-\Op{A}$ exists, then it also leaves invariant
the subspace $\Spa{F}$,
and the restriction to the subspace $\Spa{F}$ of the inverse is
an inverse of the operator $\unity_{\Spa{F}}-\Op{A}_{\Spa{F}}$.
Conversely, if the inverse of the latter operator exists,
then the operator $\unity_{\Spa{E}}+
\left(\unity_{\Spa{F}}-\Op{A}_{\Spa{F}}\right)^{-1}\comp
\Op{A}$ is a right inverse of $\unity_{\Spa{E}}-\Op{A}$:
$$\left(\unity_{\Spa{E}}-\Op{A}\right)\comp
\left(\unity_{\Spa{E}}+
\left(\unity_{\Spa{F}}-\Op{A}_{\Spa{F}}\right)^{-1}\comp
\Op{A}\right)=\unity_{\Spa{E}}-\Op{A}+
\left(\unity_{\Spa{F}}-\Op{A}_{\Spa{F}}\right)\comp
\left(\unity_{\Spa{F}}-\Op{A}_{\Spa{F}}\right)^{-1}\comp
\Op{A} =\unity_{\Spa{E}}.$$
Moreover, since the range of $\Op{A}$ is contained in $\Spa{F}$,
the operator
$\Op{A}\comp\left(\unity_{\Spa{E}}-\Op{A}\right)$
on $\Spa{E}$ is equal to the operator
$\left(\unity_{\Spa{F}}-\Op{A}_{\Spa{F}}\right)\comp\Op{A}$.
Therefore, the operator
$\unity_{\Spa{E}}+
\left(\unity_{\Spa{F}}-\Op{A}_{\Spa{F}}\right)^{-1}\comp
\Op{A}$ is also a left inverse of $\unity_{\Spa{E}}-\Op{A}$:
$$\left(\unity_{\Spa{E}}+
\left(\unity_{\Spa{F}}-\Op{A}_{\Spa{F}}\right)^{-1}\comp
\Op{A}\right)\comp\left(\unity_{\Spa{E}}-\Op{A}\right)
=\unity_{\Spa{E}}-\Op{A}+
\left(\unity_{\Spa{F}}-\Op{A}_{\Spa{F}}\right)^{-1}\comp
\left(\unity_{\Spa{F}}-\Op{A}_{\Spa{F}}\right)\comp
\Op{A}=\unity_{\Spa{E}}.$$
\end{proof}

In order to apply this Lemma the integrals
$\int\partial^n\delta(z-z_l)\Func{K}_1^{(m)}(z-z_l,k)d^2x$ and
$\int\Bar{\partial}^n\delta(z-z_l)\Func{K}_2^{(m)}(z-z_l,k)d^2x$
have to be evaluated. In doing so we use in some neighbourhood of
$z_l$ polar coordinates $(r,\varphi)$ with
$z-z_l=r\exp\left(2\pi\sqrt{-1}\varphi\right)$ and make a partial
Fourier transformation with respect to $\varphi$: any function $f$ may be
decomposed into an infinite sum
$$f(r,\varphi)=\sum\limits_{n\in\mathbb{N}}
\fourier{f}(n,r)\exp\left(2n\pi\sqrt{-1}\varphi\right).$$
The integrals are regularized by the rule
$$\int f(z)\delta(z-z_l)d^2x=\fourier{f}(0,0)$$
and the corresponding compatible rules for the holomorphic and
anti--holomorphic derivatives of the $\delta$--function.
With these rules we obtain the formulas
$$\int\limits_{\torus}\Phi(z,k)\Breve{\Psi}(z)d^2x=
\int\limits_{\torus}\Breve{\Phi}(z)\Psi(z,k)d^2x
=\Mat{R}(k).$$
Here $\Mat{R}(k)$ denotes the $M\times M$
matrices, whose entries indexed by $m,n=0,\ldots,M-1$
are equal to
$$\Mat{R}(k)_{m,n}=\binom{m+n}{m}\frac{(-1)^{n}}{\pi}\left(\begin{smallmatrix}
0 & \Func{K}_{1,m+n}(k)\\
-\Func{K}_{2,m+n}(k) & 0
\end{smallmatrix}\right).$$
Here $\Func{K}_{1,n}$ and $\Func{K}_{2,n}$ denote
the Laurent series at $z=0$ of the two functions $\Func{K}_1$
and $\Func{K}_2$ in the integral kernels of $\Op{R}(0,0,k,0)$
(compare with Lemma~\ref{resolvent integral kernel}):
\begin{align*}
\Func{K}_1(z,k)&=z^{-1}+\sum\limits_{n=0}^{\infty}
\Func{K}_{1,n}(k)z^n&
\Func{K}_2(z,k)&=\Bar{z}^{-1}+\sum\limits_{n=0}^{\infty}
\Func{K}_{2,n}(k)\Bar{z}^n.
\end{align*}
Hence, due to Lemma~\ref{inverse operator}, the inverse of the operator
$\unity-\Op{P}\comp\Op{R}(0,0,k,0)$ is equal to $\unity$
plus the operator with integral kernel
$$\Breve{\Psi}(z)\Mat{S}(k)\Phi(z',k)d^2x',$$
where $\Mat{S}(k)$ is the composition of $\Mat{P}$ with the
inverse of the restriction of $\unity-\Mat{P}\comp\Mat{R}(k)$
to the range of $\Mat{P}$. If $\Spa{E}$ and $\Spa{F}$ denote
the kernel and the range of $\Mat{P}$, respectively,
then let $\Tilde{\Mat{P}}$ and
$\Tilde{\Mat{R}}(k)$ denote the operators
\begin{align*}
\Tilde{\Mat{P}}:&\mathbb{C}^2\otimes\mathbb{C}^M/\Spa{E}
\xrightarrow{\Mat{P}}\Spa{F}, &
\Tilde{\Mat{R}}(k):&\Spa{F}\hookrightarrow
\mathbb{C}^2\otimes\mathbb{C}^M\xrightarrow{\Mat{R}(k)}
\mathbb{C}^2\otimes\mathbb{C}^M\twoheadrightarrow
\mathbb{C}^2\otimes\mathbb{C}^M/\Spa{E}.
\end{align*}
Thus $\Mat{S}(k)$ is given by
$$\Mat{S}(k):
\mathbb{C}^2\otimes\mathbb{C}^M\twoheadrightarrow
\mathbb{C}^2\otimes\mathbb{C}^M/\Spa{E}
\xrightarrow{\left(\Tilde{\Mat{P}}^{-1}-\Tilde{\Mat{R}}(k)\right)^{-1}}
\Spa{F}\hookrightarrow\mathbb{C}^2\otimes\mathbb{C}^M.$$
We remark that by definition of $\Spa{E}$ and $\Spa{F}$
the operator $\Tilde{\Mat{P}}$ is invertible.
It is quite easy to see that all these functions $\Mat{S}(k)$ obey
for all $k,k'\in\mathbb{C}^2$ the relation
$$\Mat{S}(k)-\Mat{S}(k')=\Mat{S}(k)\comp
\left(\Mat{R}(k)-\Mat{R}(k')\right)\comp\Mat{S}(k').$$
Conversely, any solution of these equations is of the form
$$\Mat{S}(k):
\mathbb{C}^2\otimes\mathbb{C}^M\twoheadrightarrow
\mathbb{C}^2\otimes\mathbb{C}^M/\Spa{E}
\xrightarrow{\left(\Tilde{\Mat{Q}}-\Tilde{\Mat{R}}(k)\right)^{-1}}
\Spa{F}\hookrightarrow\mathbb{C}^2\otimes\mathbb{C}^M,$$
where $\Spa{F}$ and $\Spa{E}$ are subspaces of
$\mathbb{C}^2\otimes\mathbb{C}^M$ such that the co--dimension of
$\Spa{E}$ is equal to the dimension of $\Spa{F}$.
Moreover, $\Tilde{\Mat{Q}}$ is any (not necessarily invertible)
operator from $\Spa{F}$ into $\mathbb{C}^2\otimes\mathbb{C}^M/\Spa{E}$.

\noindent
{\bf 2.} If the support of the \Em{Finite rank Perturbation}
contains several points, then the index set
(analogous to $\{0,\ldots,M-1\}$)
is determined by an integral divisor
$D=\sum\limits_{l=1}^{L}M_lz_l$.
The analogous ansatz for the integral kernel of a
\Em{Finite rank Perturbation} is
$$\Psi_{D}(z,k)\Mat{S}_{D}(k)\Phi_{D}(z',k)d^2x',$$
where $\Psi_{D}$ and $\Phi_{D}$ are the analogous row--vectors and
column--vectors with $\deg(D)$ entries indexed by
$m_l\in\bigcup\limits_{l=1}^{L}\{0,\ldots,M_l-1\}$ and
$\Mat{S}_{D}$ is a $\deg(D)\times\deg(D)$--matrix,
all whose entries are $2\times 2$--matrices.
Let $\Mat{R}_{D}(k)$ denote the
$\deg(D)\times\deg(D)$--matrix, whose entries indexed by
$(m_i,n_j)\in\bigcup\limits_{i=1}^{L}\{0,\ldots,M_i-1\}\times
\bigcup\limits_{j=1}^{L}\{0,\ldots,M_j-1\}$ are given by
$$\Mat{R}_{m_i,n_j}=\begin{cases}
\displaystyle{\binom{m_i+n_i}{m_i}\frac{(-1)^{n_i}}{\pi}}
\begin{pmatrix}
0 & \Func{K}_{1,m_i+n_i}(k)\\
-\Func{K}_{2,m_i+n_i}(k) & 0
\end{pmatrix} & \text{if }i=j\\
\displaystyle{\frac{(-1)^{n_j}}{\pi m_i!n_j!}}\begin{pmatrix}
0 & \Func{K}_1^{(m_i+n_j)}(z_i-z_j,k)\\
-\Func{K}_2^{(m_i+m_j)}(z_i-z_j,k) & 0
\end{pmatrix} & \text{if }i\neq j.
\end{cases}$$
Therefore, the \Em{Finite rank Perturbations} of $\Op{R}(0,0,k,0)$
are defined as the operators with integral kernels
$$\Psi_{D}(z,k)\Mat{S}_{D}(k)\Phi_{D}(z',k)d^2x'.$$
Here $D$ is any integral divisor on $\mathbb{C}/\lattice_\mathbb{C}$
and $\Mat{S}_{D}(k)$ satisfies for all $k,k'\in\mathbb{C}^2$
the relation
$$\Mat{S}_{D}(k)-\Mat{S}_{D}(k')=\Mat{S}_{D}(k)\comp
\left(\Mat{R}_{D}(k)-\Mat{R}_{D}(k')\right)\comp\Mat{S}_{D}(k').$$

\noindent
{\bf 3.} The perturbed free resolvent
$\Op{R}(0,0,k,0)+\Op{R}(0,0,k,0)\comp\Op{S}(k)\comp\Op{R}(0,0,k,0)$
induces the following perturbation of the resolvent $\Op{R}(V,W,k,0)$:
\begin{align*}
&\Op{R}(0,0,k,0)\comp
\left(\unity+\Op{S}(k)\comp\Op{R}(0,0,k,0)\right)\comp
\left(\unity-\left(\begin{smallmatrix}
V & 0\\
0 & W
\end{smallmatrix}\right)\comp\Op{R}(0,0,k,0)\comp
\left(\unity+\Op{S}(k)\comp\Op{R}(0,0,k,0)\right)\right)^{-1}\\
&=\Op{R}(0,0,k,0)\comp
\left(\left(\unity+\Op{S}(k)\comp\Op{R}(0,0,k,0)\right)^{-1}-
\left(\begin{smallmatrix}
V & 0\\
0 & W
\end{smallmatrix}\right)\comp\Op{R}(0,0,k,0)\right)^{-1}\\
&=\Op{R}(0,0,k,0)\comp
\left(\unity-
\left(\unity+\Op{S}(k)\comp\Op{R}(0,0,k,0)\right)^{-1}\comp
\Op{S}(k)\comp\Op{R}(0,0,k,0)-
\left(\begin{smallmatrix}
V & 0\\
0 & W
\end{smallmatrix}\right)\comp\Op{R}(0,0,k,0)\right)^{-1}\\
&=\Op{R}(V,W,k,0)\comp
\left(\unity-
\left(\unity+\Op{S}(k)\comp\Op{R}(0,0,k,0)\right)^{-1}\comp
\Op{S}(k)\comp\Op{R}(V,W,k,0)\right)^{-1}.
\end{align*}
This implies that the corresponding perturbation is equal to
\begin{multline*}
\Op{R}(V,W,k,0)\comp\Op{S}(V,W,k)\comp\Op{R}(V,W,k,0)
\text{ with }\Op{S}(V,W,k)=\\
\begin{aligned}
=&\left(\unity-
\left(\unity+\Op{S}(k)\comp\Op{R}(0,0,k,0)\right)^{-1}\comp
\Op{S}(k)\comp\Op{R}(V,W,k,0)\right)^{-1}\comp
\left(\unity+\Op{S}(k)\comp\Op{R}(0,0,k,0)\right)^{-1}\comp
\Op{S}(k)&&\\
=&\left(\unity+\Op{S}(k)\comp\Op{R}(0,0,k,0)-
\Op{S}(k)\comp\Op{R}(V,W,k,0)\right)^{-1}\comp
\Op{S}(k)&&\\
=&\Op{S}(k)\comp\left(\unity+\Op{R}(0,0,k,0)\comp\Op{S}(k)-
\Op{R}(V,W,k,0)\comp\Op{S}(k)\right)^{-1}.&&
\end{aligned}\end{multline*}
In the notation of the first two steps $\Op{S}(k)$ has the integral kernel
$\Breve{\Psi}_{D}(z)\Mat{S}_{D}(k)\Breve{\Phi}_{D}(z')d^2x'$.
With the regularization rules of the first step we have
$\Mat{R}_{D}(k)=\left\langle\left\langle\Breve{\Phi}_{D},
\Op{R}_{D}(0,0,k,0)\Breve{\Psi}_{D}\right\rangle\right\rangle.$
For pairs of smooth potentials we may analogously define
$\Mat{R}_{D}(V,W,k)=
\left\langle\left\langle\Breve{\Phi}_{D},
\Op{R}_{D}(V,W,k,0)\Breve{\Psi}_{D}\right\rangle\right\rangle$.
Then the perturbation $\Op{S}(V,W,k)$ has the integral kernel
$\Breve{\Psi}_{D}(z)\Mat{S}_{D}(V,W,k)\Breve{\Phi}_{D}(z')d^2x'$
with 
\begin{eqnarray*}
\Mat{S}_{D}(V,W,k)&=&
\left(\unity+\Mat{S}_{D}(k)\comp
(\Mat{R}(k_{D})-\Mat{R}_{D}(V,W,k))\right)^{-1}\comp\Mat{S}_{D}(k)\\
&=&\Mat{S}_{D}(k)\comp\left(\unity+(\Mat{R}(k_{D})-\Mat{R}_{D}(V,W,k))
\comp\Mat{S}_{D}(k)\right)^{-1}.
\end{eqnarray*}
Moreover, due to the relation
$\Mat{S}_{D}(k)-\Mat{S}_{D}(k')=\Mat{S}_{D}(k)\comp
\left(\Mat{R}_{D}(k)-\Mat{R}_{D}(k')\right)\comp\Mat{S}_{D}(k')$,
the matrices $\Mat{S}_{D}(V,W,k)$ obey the relation
$$\Mat{S}_{D}(V,W,k)-\Mat{S}_{D}(V,W,k')=\Mat{S}_{D}(k)\comp
\left(\Mat{R}_{D}(V,W,k)-\Mat{R}_{D}(V,W,k')\right)\comp\Mat{S}_{D}(V,W,k').$$
For pairs of non--smooth potentials the matrices $\Mat{R}_{D}(V,W,k)$ are
not defined.

\noindent
{\bf 4.} However, a reformulation of these perturbations makes sense
for all pairs of potentials in $\banach{2}$. In fact, the perturbation
$\Op{R}(V,W,k,0)\comp\Op{S}(V,W,k)\comp\Op{R}(V,W,k,0)$
has the integral kernel
$$\Psi_{D}(z,k)\Mat{S}(V,W,k)\Phi_{D}(k,z')d^2x'$$
in terms of the row--vectors and column--vectors 
\begin{align*}
\Psi_{D}(z,k)&=\Op{R}(V,W,k,0)\Breve{\Psi}_{D}(z)\text{ and}&
\Phi_{D}^{t}(z,k)&=\Op{R}^{t}(V,W,k,0)\Breve{\Phi}_{D}^{t}(z).
\end{align*}
Furthermore, due to Lemma~\ref{branchpoints},
the forms $\Phi_{D}(z,k)\left(\begin{smallmatrix}
0 & d\Bar{z}\\
dz & 0
\end{smallmatrix}\right)\Psi_{D}(z,k')$ are closed on
$\Delta\setminus\{z_1,\ldots,z_L\}$,
and for all small $\varepsilon$ the integrals over the boundary of
$\Set{S}_{\varepsilon}=\bigcup\limits_{l=1}^{L}B(z_l,\varepsilon)$
are equal to
$$\int\limits_{\partial\Set{S}_{\varepsilon}}
\Phi_{D}(z,k)\begin{pmatrix}
0 & d\Bar{z}\\
dz & 0
\end{pmatrix}\Psi_{D}(z,k')=2\sqrt{-1}\left(
\Mat{R}_{D}(V,W,k)-\Mat{R}_{D}(V,W,k')\right).$$
Hence we may renormalize $\Psi_{D}(z,k)$ and $\Phi_{D}(z,k)$
in such a way, that they have finite $\banach{q}$--norm on the
complement of $\Set{S}_{\varepsilon}$.
This leads to the following definition:

\newtheorem{Finite rank Perturbations}[Lemma]{Finite rank Perturbations}
\index{Dirac operator $\Op{D}$!perturbation of the $\sim$}
\index{perturbation!$\rightarrow$ finite rank $\sim$}
\index{finite!rank perturbation}
\index{integral!kernel!of a finite rank perturbation}
\begin{Finite rank Perturbations}\label{finite rank perturbations}
of $\Op{R}(V,W,k,0)$ are defined as perturbations
by the operators with integral kernels
$$\Psi_{D}(z,k)\Mat{S}_{D}(k)\Phi_{D}(z',k)d^2x'.$$
Here $D=\sum\limits_{l=1}^{L}M_lz_l$ is any integral divisor on
$\mathbb{C}/\lattice_\mathbb{C}$, $\Psi_{D}$ and $\Phi_{D}$ are
finite row--vectors and column--vectors, whose entries belong on
$\Delta\setminus\{z_1,\ldots,z_L\}$
to the kernel and co--kernel of $\Op{D}(V,W,k)$, respectively.
Moreover, for all $l=1,\ldots,L$ these vectors obey
on a small neighbourhood of $z_l$ the equations
\begin{align*}
\begin{pmatrix}
V\left(\Bar{z}-\Bar{z}_l\right)^{M_l} & 
\partial\left(\Bar{z}-\Bar{z}_l\right)^{M_l}\\
-\Bar{\partial}\left(z-z_l\right)^{M_l} &
W\left(z-z_l\right)^{M_l}\end{pmatrix}
\Psi_{D}(z,k)&=0&
\begin{pmatrix}
V\left(z-z_l\right)^{M_l} & 
\Bar{\partial}\left(z-z_l\right)^{M_l}\\
-\partial\left(\Bar{z}-\Bar{z}_l\right)^{M_l} &
W\left(\Bar{z}-\Bar{z}_l\right)^{M_l}\end{pmatrix}
\Phi_{D}^{t}(z,k)&=0.\end{align*}
Furthermore, the representatives of these entries
in the quotients of these kernels modulo
the kernels and co--kernels of $\left(\begin{smallmatrix}
V & \partial\\
-\Bar{\partial} & W
\end{smallmatrix}\right)$ (compare with Lemma~\ref{quotients of kernels}),
do not depend on $k$
(i.\ e.\ locally there exist $\Tilde{\Psi}_{D}(z,k)$ and
$\Tilde{\Phi}_{D}(z,k)$ in the kernels and co--kernel of
$\left(\begin{smallmatrix}
V & \partial\\
-\Bar{\partial} & W
\end{smallmatrix}\right)$,
such that the differences $\Psi_{D}(z,k)-\Tilde{\Psi}_{D}(z,k)$ and
$\Phi_{D}(z,k)-\Tilde{\Phi}_{D}(z,k)$ do not depend on $k$).
Finally, for all $k,k'$ the matrices $\Mat{S}_{D}$ obey the equation
$$\Mat{S}_{D}(k)-\Mat{S}_{D}(k')=
\frac{1}{2\sqrt{-1}}\int\limits_{\partial\Set{S}_{\varepsilon}}
\Mat{S}_{D}(k)\Phi_{D}(z,k)\begin{pmatrix}
0 & d\Bar{z}\\
dz & 0
\end{pmatrix}\Psi_{D}(z,k')\Mat{S}_{D}(k').$$
\end{Finite rank Perturbations}

\begin{Remark}\label{kernels of Dirac operators}
The kernels of the operators $\left(\begin{smallmatrix}
V\left(\Bar{z}-\Bar{z}_l\right)^{M_l} & 
\partial\left(\Bar{z}-\Bar{z}_l\right)^{M_l}\\
-\Bar{\partial}\left(z-z_l\right)^{M_l} &
W\left(z-z_l\right)^{M_l}
\end{smallmatrix}\right)$ and $\left(\begin{smallmatrix}
V\left(z-z_l\right)^{M_l} & 
\Bar{\partial}\left(z-z_l\right)^{M_l}\\
-\partial\left(\Bar{z}-\Bar{z}_l\right)^{M_l} &
W\left(\Bar{z}-\Bar{z}_l\right)^{M_l}
\end{smallmatrix}\right)$ are those spinors $\psi$ and $\phi$,
such that $\left(\begin{smallmatrix}
z-z_l & 0\\
0 & \Bar{z}-\Bar{z}_l
\end{smallmatrix}\right)^{M_l}\psi$ and 
$\left(\begin{smallmatrix}
\Bar{z}-\Bar{z}_l & 0\\
0 & z-z_l
\end{smallmatrix}\right)^{M_l}\phi$ belongs to the kernel of
$\begin{pmatrix}
V\left(\frac{\Bar{z}-\Bar{z}_l}{z-z_l}\right)^{M_l} & \partial\\
-\Bar{\partial} & W\left(\frac{z-z_l}{\Bar{z}-\Bar{z}_l}\right)^{M_l}
\end{pmatrix}$ and $\begin{pmatrix}
V\left(\frac{z-z_l}{\Bar{z}-\Bar{z}_l}\right)^{M_l} & \Bar{\partial}\\
-\partial & W\left(\frac{\Bar{z}-\Bar{z}_l}{z-z_l}\right)^{M_l}
\end{pmatrix}$, respectively
(compare with Section~\ref{subsubsection local behaviour}).
\end{Remark}

\noindent
{\bf 5.} We shall generalize Cauchy's integral formula
to these elements in the kernel of Dirac operators with potentials
$V,W\in\banach{2}$.
If $\Func{K}_{\mathbb{R}^2}(V,W,z,z')\frac{d^2x'}{\pi}$
denotes the integral kernel of the resolvent
$\Op{R}_{\mathbb{R}^2}(V,W,0)$
on the bounded open set $\Set{O}\subset\mathbb{C}$, then we have
\begin{align*}
\begin{pmatrix}
V & \partial\\
-\Bar{\partial} & W
\end{pmatrix}\Func{K}_{\mathbb{R}^2}(V,W,z,z')&=
\pi\delta(z-z')\unity &
\begin{pmatrix}
V & \Bar{\partial}\\
-\partial & W
\end{pmatrix}\Func{K}_{\mathbb{R}^2}^{t}(V,W,z',z)&=
\pi\delta(z-z')\unity.
\end{align*}
Here the Dirac operator and his transposed acts on the integral kernel
as a function depending on $z$ for fixed $z'$.
The arguments of the proof of Lemma~\ref{branchpoints}
imply the

\newtheorem{Generalized Cauchy}[Lemma]{Generalized Cauchy's integral formula}
\index{generalized!Cauchy's integral formula}
\begin{Generalized Cauchy}\label{generalized cauchy}
All elements $\psi$ and $\phi$ in the kernel and co--kernel of the
Dirac operator $\left(\begin{smallmatrix}
V & \partial\\
-\Bar{\partial} & W
\end{smallmatrix}\right)$ on a small open set $\Set{O}$
obey the formula
\begin{equation*}\begin{split}
\psi(z')&=\frac{1}{2\pi\sqrt{-1}}\oint\Func{K}_{\mathbb{R}^2}(V,W,z',z)
\begin{pmatrix}
0 & d\Bar{z}\\
dz & 0
\end{pmatrix}\psi(z)\\
\phi(z')&=\frac{-1}{2\pi\sqrt{-1}}\oint\phi(z)\begin{pmatrix}
0 & d\Bar{z}\\
dz & 0
\end{pmatrix}\Func{K}_{\mathbb{R}^2}(V,W,z,z'),
\end{split}\end{equation*}
as long as the integration path surrounds $z'$
one times in the anti--clockwise--order, respectively.
\end{Generalized Cauchy}

\begin{Remark}\label{smeared closed path}
At a first look it is not clear, whether the integral along the closed
path is well defined. However, since on the complement of $\{z'\}$ the
corresponding one--forms are closed, we may extend the integration
over the closed path to an integration over a cylinder around $z'$.
This extended closed integral can be considered as the evaluation of
the action of the resolvent on a $\banach{q}$--spinor.
Moreover, since the spinors $\psi$ and $\phi$
are in general not continuous, this equality should be
understood as an equality of measurable functions.
\end{Remark}

\noindent
{\bf 6.} Now we may extend the classification of these
\De{Finite rank Perturbations}~\ref{finite rank perturbations}
of the free resolvents in the first two steps
to pairs of $\banach{2}$--potentials.
For all integral divisors $D=\sum\limits_{l=1}^{L}M_lz_l$
we introduce the sheaves analogous to $\Sh{O}_{D}$.
On the complement of the support of $D$ the corresponding sections
belong to the kernels and co--kernels of $\left(\begin{smallmatrix}
V & \partial\\
-\Bar{\partial} & W
\end{smallmatrix}\right)$. Moreover, for all $l=1,\ldots,L$ these
sections $\psi$ and $\phi$ obey in a small neighbourhood of $z_l$
\begin{align*}
\begin{pmatrix}
V\left(\Bar{z}-\Bar{z}_l\right)^{M_l} & 
\partial\left(\Bar{z}-\Bar{z}_l\right)^{M_l}\\
-\Bar{\partial}\left(z-z_l\right)^{M_l} &
W\left(z-z_l\right)^{M_l}\end{pmatrix}
\psi&=0&
\begin{pmatrix}
V\left(z-z_l\right)^{M_l} & 
\Bar{\partial}\left(z-z_l\right)^{M_l}\\
-\partial\left(\Bar{z}-\Bar{z}_l\right)^{M_l} &
W\left(\Bar{z}-\Bar{z}_l\right)^{M_l}\end{pmatrix}
\phi^{t}&=0.
\end{align*}
We chose a basis of the quotient of these sheaves
modulo the kernel and co--kernel of 
$\left(\begin{smallmatrix}
V & \partial\\
-\Bar{\partial} & W
\end{smallmatrix}\right)$ (compare with Lemma~\ref{quotients of kernels}),
respectively. These basis should be represented by local elements
$\psi_{m_l}$ and $\phi_{m_l}$ of these sheaves
nearby the singular points $z_1,\ldots,z_L$.
We arrange the former basis to a row--vector $\Psi_{D}(z)$
and the latter basis to e column--vector $\Phi_{D}(z)$
with $\deg(D)$ entries, respectively.
Moreover, we assume that the corresponding matrix of residues vanishes
$$\res\limits_{D}\Phi_{D}(z)\begin{pmatrix}
0 & d\Bar{z}\\
dz & 0
\end{pmatrix}
\Psi_{D}(z)=0.$$
If $\Func{K}(V,W,k,z,z')\frac{d^2x'}{\pi}$
denotes the integral kernel of the operator $\Op{R}(V,W,k,0)$,
then the row--vectors and columns vectors
\begin{equation*}\begin{split}
\Psi_{D}(z',k)&=\res\limits_{D}\Func{K}(V,W,k,z',z)\left(\begin{smallmatrix}
0 & d\Bar{z}\\
dz & 0
\end{smallmatrix}\right)\Psi_{D}(z)\\
\Phi_{D}(z',k)&=-\res\limits_{D}\Phi_{D}(z)\left(\begin{smallmatrix}
0 & d\Bar{z}\\
dz & 0
\end{smallmatrix}\right)\Func{K}(V,W,k,z,z')
\end{split}\end{equation*}
fulfill the conditions of
\De{Finite rank Perturbations}~\ref{finite rank perturbations}.
We remark that these functions are defined only for $z'$
outside of the cycles around the support of $D$,
which are used to calculate the residue, and for $k\not\in\fermi(V,W)$.
Due to the
\De{Generalized Cauchy's integral formula}~\ref{generalized cauchy}
the differences $\Psi_{D}(z,k)-\Psi_{D}(z)$ and
$\Phi_{D}(z,k)-\Phi_{D}(z)$ belong on a small neighbourhood
of the support of $D$ to the kernel and co--kernel of
$\left(\begin{smallmatrix}
V & \partial\\
-\Bar{\partial} & W
\end{smallmatrix}\right)$, respectively.
Since the product $\Phi_{D}(z,k)\left(\begin{smallmatrix}
0 & d\Bar{z}\\
dz & 0
\end{smallmatrix}\right)\Psi_{D}(z,k)$ is periodic the matrix of residues
$\res\limits_{D}\Phi_{D}(z,k)\left(\begin{smallmatrix}
0 & d\Bar{z}\\
dz & 0
\end{smallmatrix}\right)\Psi_{D}(z,k)$ vanishes.
Hence the following matrices of residues coincide:
$$\Mat{R}_{D}(V,W,k)=-\pi\res\limits_{D}\Phi_{D}(z)\left(\begin{smallmatrix}
0 & d\Bar{z}\\
dz & 0
\end{smallmatrix}\right)\Psi_{D}(z,k)=
\pi\res\limits_{D}\Phi_{D}(z,k)\left(\begin{smallmatrix}
0 & d\Bar{z}\\
dz & 0
\end{smallmatrix}\right)\Psi_{D}(z).$$
In general we have
$$\res\limits_{D}\Phi_{D}(z,k)\left(\begin{smallmatrix}
0 & d\Bar{z}\\
dz & 0
\end{smallmatrix}\right)\Psi_{D}(z,k')=
\Mat{R}_{D}(V,W,k)-\Mat{R}_{D}(V,W,k'),$$
and the classification of
\De{Finite rank Perturbations}~\ref{finite rank perturbations}
of $\Op{R}(0,0,k,0)$
in the first two steps carries over to a classification of
\De{Finite rank Perturbations}~\ref{finite rank perturbations}
of $\Op{R}(V,W,k,0)$. Consequently, all
\De{Finite rank Perturbations}~\ref{finite rank perturbations}
of $\Op{R}(V,W,k,0)$ have an integral kernel of the form
$\Psi_{D}(z,k)\Mat{S}_{D}(k)\Phi_{D}(z',k)d^2x'$ with
$$\Mat{S}_{D}(k)-\Mat{S}_{D}(k')=\Mat{S}_{D}(k)\comp
(\Mat{R}_{D}(V,W,k)-\Mat{R}_{D}(V,W,k')\comp\Mat{S}_{D}(k').$$

\noindent
{\bf 7.} The foregoing considerations lead to another description of
these \De{Finite rank Perturbations}~\ref{finite rank perturbations}.
The corresponding Dirac operators yield holomorphic structures
(in the sense of `quaternionic function theory' \cite{PP,BFLPP,FLPP})
with singularities at the points $\{z_1,\ldots,z_L\}$.
Let us first explain these singularities for one singular point $z_l$.
We have seen in the first step that in these cases the matrices
$\Mat{S}$ are of the form
$$\Mat{S}(k):
\mathbb{C}^2\otimes\mathbb{C}^M\twoheadrightarrow
\mathbb{C}^2\otimes\mathbb{C}^M/\Spa{E}
\xrightarrow{\left(\Tilde{\Mat{Q}}-\Tilde{\Mat{R}}(k)\right)^{-1}}
\Spa{F}\hookrightarrow\mathbb{C}^2\otimes\mathbb{C}^M.$$
Here $\Spa{F}$ and $\Spa{E}$ are subspaces of
$\mathbb{C}^2\otimes\mathbb{C}^M$ such that the co--dimension of
$\Spa{E}$ is equal to the dimension of $\Spa{F}$.
Moreover, $\Tilde{\Mat{Q}}$ is any operator from
$\Spa{F}$ into $\mathbb{C}^2\otimes\mathbb{C}^M/\Spa{E}$.
For smooth potentials we determined formally in the first step the
corresponding perturbations of the Dirac operators. Surprisingly, we
may define the kernels of these perturbed Dirac operators for all
$\banach{2}$--potentials. In fact, by choice of the row--vectors $\Psi(z)$
and column--vectors $\Phi(z)$ in the sixth step these data
$(\Spa{E},\Spa{F},\Tilde{\Mat{Q}})$ are in one--to--one correspondence
to the deformed kernels $\Spa{H}_{\text{\scriptsize\rm sing},z_l}$
with the following properties:
\begin{description}
\index{condition!singularity (i)--(ii)}
\index{singularity!condition $\sim$ (i)--(ii)}
\item[Singularity (i)]
The space $\Spa{H}_{\text{\scriptsize\rm sing},z_l}$
is a subspace of the kernel of
$\left(\begin{smallmatrix}
V(\Bar{z}-\Bar{z}_l)^{M_l} & \partial(\Bar{z}-\Bar{z}_l)^{M_l}\\
-\Bar{\partial}(z-z_l)^{M_l} & W(z-z_l)^{M_l}
\end{smallmatrix}\right)$ of co--dimension $2M_l$, i.\ e.\ the same
co--dimension as the kernel of $\left(\begin{smallmatrix}
V & \partial\\
-\Bar{\partial} & W
\end {smallmatrix}\right)$ (compare with Lemma~\ref{quotients of kernels}).
\item[Singularity (ii)] The kernel of
$\left(\begin{smallmatrix}
V(\Bar{z}-\Bar{z}_l)^{-M_l} & \partial(\Bar{z}-\Bar{z}_l)^{-M_l}\\
-\Bar{\partial}(z-z_l)^{-M_l} & W(z-z_l)^{-M_l}
\end{smallmatrix}\right)$
is a subspace of the space
$\Spa{H}_{\text{\scriptsize\rm sing},z_l}$ of co--dimension $2M_l$,
i.\ e.\ the same co--dimension as this kernel considered
as a subspace of the kernel of
$\left(\begin{smallmatrix}
V & \partial\\
-\Bar{\partial} & W
\end {smallmatrix}\right)$ (compare with Lemma~\ref{quotients of kernels}).
\end{description}
More precisely, the entries of $\Psi(z)$
fit together to a basis of the quotient of the kernel of
$\left(\begin{smallmatrix}
V(\Bar{z}-\Bar{z}_l)^{M_l} & \partial(\Bar{z}-\Bar{z}_l)^{M_l}\\
-\Bar{\partial}(z-z_l)^{M_l} & W(z-z_l)^{M_l}
\end{smallmatrix}\right)$ modulo the kernel of
$\left(\begin{smallmatrix}
V & \partial\\
-\Bar{\partial} & W
\end{smallmatrix}\right)$. On the other hand the entries of $\Phi(z)$
are dual with respect to the pairing given by the residue
to a basis of the quotient of the kernel of
$\left(\begin{smallmatrix}
V & \partial\\
-\Bar{\partial} & W
\end{smallmatrix}\right)$ modulo the kernel of
$\left(\begin{smallmatrix}
V(\Bar{z}-\Bar{z}_l)^{-M_l} & \partial(\Bar{z}-\Bar{z}_l)^{-M_l}\\
-\Bar{\partial}(z-z_l)^{-M_l} & W(z-z_l)^{-M_l}
\end{smallmatrix}\right)$. Due to the assumption that the matrix of
residues of $\Phi(z)\left(\begin{smallmatrix}
0 & d\Bar{z}\\
dz & 0
\end{smallmatrix}\right)\Psi(z)$ vanishes, these two basis
fit together to a basis of the quotient of the kernel of
$\left(\begin{smallmatrix}
V(\Bar{z}-\Bar{z}_l)^{M_l} & \partial(\Bar{z}-\Bar{z}_l)^{M_l}\\
-\Bar{\partial}(z-z_l)^{M_l} & W(z-z_l)^{M_l}
\end{smallmatrix}\right)$ modulo the kernel of
$\left(\begin{smallmatrix}
V(\Bar{z}-\Bar{z}_l)^{-M_l} & \partial(\Bar{z}-\Bar{z}_l)^{-M_l}\\
-\Bar{\partial}(z-z_l)^{-M_l} & W(z-z_l)^{-M_l}
\end{smallmatrix}\right)$.
Therefore the spaces $\Spa{E}$ and $\Spa{F}$ are subspaces of this
quotient, and the the corresponding space
$\Spa{H}_{\text{\scriptsize\rm sing},z_l}$ is the sum of the
the kernel of
$\left(\begin{smallmatrix}
V(\Bar{z}-\Bar{z}_l)^{-M_l} & \partial(\Bar{z}-\Bar{z}_l)^{-M_l}\\
-\Bar{\partial}(z-z_l)^{-M_l} & W(z-z_l)^{-M_l}
\end{smallmatrix}\right)$ plus $\Spa{E}$ plus the graph of
$\Tilde{\Mat{Q}}$ considered as an operator from $\Spa{F}$ into the
quotient of the kernel of $\left(\begin{smallmatrix}
V & \partial\\
-\Bar{\partial} & W
\end{smallmatrix}\right)$ modulo the sum of
the kernel of  $\left(\begin{smallmatrix}
V(\Bar{z}-\Bar{z}_l)^{-M_l} & \partial(\Bar{z}-\Bar{z}_l)^{-M_l}\\
-\Bar{\partial}(z-z_l)^{-M_l} & W(z-z_l)^{-M_l}
\end{smallmatrix}\right)$ plus $\Spa{E}$.

In analogy to `quaternionic functions theory' the Dirac operators
are considered as holomorphic structures. Consequently, the
\De{Finite rank Perturbations}~\ref{finite rank perturbations} give
rise to

\newtheorem{Singularities of the structure}[Lemma]{Singularities of the
  holomorphic structure}
\index{singularity!of the holomorphic structure}
\begin{Singularities of the structure}\label{singularities of the structure}
The kernel of a perturbed Dirac operators is locally in a
neighbourhood $\Set{O}$ of a singular point $z_1$ a subspace
$\Spa{H}_{\text{\scriptsize\rm sing},z_l}$ of
\begin{align*}
\left\{\psi\mid\left(\begin{smallmatrix}
z-z_l & 0\\
0 & \Bar{z}-\Bar{z}_l
\end{smallmatrix}\right)^{M_l}\psi\in
\sobolev{1,p}_{\text{\scriptsize\rm loc}}(\Set{O})\times
\sobolev{1,p}_{\text{\scriptsize\rm loc}}(\Set{O})\right\}
&\text{ with $1<p<2$}&
\text{(compare with Section~\ref{subsubsection local behaviour}),}&
\end{align*}
which fulfills conditions \Em{Singularity}~(i)--(ii).
\end{Singularities of the structure}

Obviously, for a given pair of potentials and a given $M_l$
the set of those spaces $\Spa{H}_{\text{\scriptsize\rm sing},z_l}$,
which obey condition \Em{Singularity}~(i)--(ii),
can be considered as a finite--dimensional Grassmannian
\cite[Chapter~1 \S5.]{GrHa}. 
The resolvents of the
\De{Finite rank Perturbations}~\ref{finite rank perturbations}
maps all spinors, whose support is disjoint from a small neighbourhood
$\Set{O}$ of $z_l$, into the corresponding deformed kernel
$\Spa{H}_{\text{\scriptsize\rm sing},z_l}$.
Conversely, the corresponding
\De{Finite rank Perturbations}~\ref{finite rank perturbations}
are uniquely defined by this property.

\begin{Remark}\label{existence of perturbations}
From the very beginning it is not clear, whether for all these choices of
spaces $\Spa{H}_{\text{\scriptsize\rm sing},z_l}$
(or equivalently of data $(\Spa{E},\Spa{F},\Tilde{\Mat{Q}})$)
the corresponding matrices $\Mat{S}(k)$ are indeed meromorphic.
More precisely, the finiteness of $\Mat{S}(k)$ for at least one
value has to be guaranteed. We show this in the following step
for all perturbations, which respect the involution $\eta$.
The corresponding spaces $\Spa{H}_{\text{\scriptsize\rm sing},z_l}$
are invariant under the (projective) involution
$\psi\mapsto\Op{J}\Bar{\psi}$.
\end{Remark}

\begin{Remark}\label{local contribution to the willmore functional}
The results of Section~\ref{subsubsection compactified isospectral}
imply that these singularities contribute to the Willmore functional.
More precisely, if $\fermi$ is the \Em{complex Fermi curve} of a
\De{Finite rank Perturbation}~\ref{finite rank perturbations} of
$\Op{R}(U,\Bar{U},k,0)$, then the Willmore functional is equal to
$$\willmore(\fermi)=4\|U\|_2^2+\sum\limits_{l=1}^{L}
\willmore_{\text{\scriptsize\rm sing},z_l}.$$
In the correspondence of spaces $\Spa{H}_{\text{\scriptsize\rm sing},z_l}$
and data $(\Spa{E},\Spa{F},\Tilde{\Mat{Q}})$
we realized $\Spa{F}$ as a space of `meromorphic' spinors
and $\Spa{E}$ as a space of `holomorphic' spinors.
Let $0<n_1+1<\ldots<n_J+1$ denote the sequence of pole orders of
the elements of $\Spa{F}$ and $0\leq m_1<\ldots<n_j$ denote the sequence
of \De{Orders of zeroes}~\ref{order of zeroes} of the elements of $\Spa{E}$.
(In Section~\ref{subsubsection family} and \ref{subsubsection general case}
we will see that these sequences of integers
$0\leq n_1<\ldots<n_J$ and $0\leq m_1<\ldots<m_J$ are
the \De{Orders of zeroes}~\ref{order of zeroes} at $\infty$ of the
spinors in the kernel and co--kernel of a blown up
Dirac operator on $\mathbb{P}^1$.)
The corresponding sets $\Set{N}$ introduced in
Section~\ref{subsubsection compactified isospectral} are given by
$$\Set{N}=\{-n_J-1,\ldots,-n_1-1\}\cup\mathbb{N}_0
\setminus\{m_1,\ldots,m_J\}.$$
The corresponding local contribution to the
Willmore functional is equal to
$$\willmore_{\text{\scriptsize\rm sing},z_l}=
4\pi\left(J+n_1+\ldots+n_J+m_1+\ldots+m_j\right).$$
For all $d\in\mathbb{N}$ let $\Sh{L}_{(d-1)\infty}$
denote the sheaf of spinors in the kernel
of the corresponding Dirac operator on $\mathbb{P}^1$,
which may have poles of order $d-1$ at $\infty$,
i.\ e.\ the sections obey $\text{\rm ord}_{\infty}(\psi)\geq 1-d$
(compare with Section~\ref{subsubsection local behaviour}).
The corresponding quaternionic holomorphic line bundle has degree $d$. 
Due to extension of S\'{e}rre Duality and the Riemann--Roch theorem
to these quaternionic holomorphic line bundles
(compare with \cite[Theorem~2.2]{FLPP} and
Section~\ref{subsubsection local behaviour})
for large $d$ the quaternionic dimension of
$\dim H^0\left(\mathbb{P}^1,\Sh{L}_{(d-1)\infty}\right)$ is equal to $d+1$.
The corresponding \De{Orders of zeroes}~\ref{order of zeroes}
at $\infty$ are equal to the numbers
$\{n\in\mathbb{N}_0\mid d-n-1\in\Set{N}\}$.
The corresponding numbers, which in \cite[Definition~4.2]{FLPP}
are called $\text{\rm ord}_{\infty}\left(\Sh{L}_{(d-1)\infty}\right)$,
are for large $d$ equal to
$$\text{\rm ord}_{\infty}\left(\Sh{L}_{(d-1)\infty}\right)=
\sum\limits_{\{n\in\mathbb{N}_0\mid d-n-1\in\Set{N}\}}n-\frac{1}{2}d(d-1)=
\sum\limits_{\{n\in\mathbb{N}\mid -n\in\Set{N}\}}n+
\sum\limits_{\{n\in\mathbb{N}_0\mid n\not\in\Set{N}\}}n=
\frac{\willmore_{\text{\scriptsize\rm sing},z_l}}{4\pi}.$$
Hence the local contribution to the Willmore functional is exactly
equal to the lower bound of the Pl\"ucker formula \cite[Theorem~4.7]{FLPP}
for the quaternionic holomorphic line bundle corresponding to
$\Sh{L}_{(d-1)\infty}$ for large $d$.
\end{Remark}

\begin{Remark}\label{locality}
The local structure of these perturbations yields a further condition
for the \De{Finite rank Perturbations}~\ref{finite rank perturbations}.
The data $(\Spa{E}_{D},\Spa{F}_{D},\Tilde{\Mat{Q}}_{D})$ have
diagonal block form with respect to the decomposition into $L$ blocks
corresponding to the singular points $\{z_1,\ldots,z_l\}$.
Moreover, if the pair of potentials is invariant under one or several
involutions introduced in Section~\ref{subsection reductions},
then the corresponding perturbations are covariant under the
corresponding involutions (compare with Remark~\ref{additional restrictions}).
\end{Remark}

\begin{Remark}\label{generalization}
The local structure suggests the following
generalization of Theorem~\ref{limits of resolvents}:
A subsequence of the sequence of kernels $\Spa{H}_n$
of any sequence of Dirac operators $\left(\begin{smallmatrix}
U_n & \partial\\
-\Bar{\partial} & \Bar{U}_n
\end{smallmatrix}\right)$ with bounded potentials
$U_n\in\banach{2}(\Set{O})$ on a bounded open domain
$\Set{O}\subset\mathbb{C}$ converges (in an appropriate sense)
to a kernel with
\De{Singularities of the holomorphic structure}~\ref{singularities of
  the structure}.
From this point of view the
\De{Quaternionic Singularity condition}~\ref{quaternionic condition}
states that the there should exist
a quaternionic holomorphic line bundle
whose index of specialty is at least equal to one.
Therefore, the extension to immersion of higher genus
Riemann surfaces depends on an investigation of the
subsets of the Picard groups of special divisors.
For genus one, these are exactly the \Em{complex Fermi curves}.
For higher genus, the corresponding Bloch theory give rise to
higher--dimensional analogs of \Em{complex Fermi curves}.

The spinors in the kernel of a given Dirac operator
$\left(\begin{smallmatrix}
V & \partial\\
-\Bar{\partial} & W
\end {smallmatrix}\right)$ are uniquely determined by the restrictions
to the boundary $\partial B(z_l,\varepsilon)$ of a ball $B(z_l,\varepsilon)$.
Moreover, the spaces $\Spa{H}_{\text{\scriptsize\rm sing},z_l}$
may be described by the corresponding subspaces of the Banach spaces
$\banach{q}(\partial B(z_l,\varepsilon)\times
\banach{q}(\partial B(z_l,\varepsilon)$.
The description of the latter subspaces as graphs of operators
yields a suitable notion of convergence
(compare with the description of the Grassmannian of a Hilbert space
in \cite[Chapter~7.]{PrSe}).
Obviously, this generalization extends the compactness
property to immersions of arbitrary compact Riemann surfaces
into $\mathbb{R}^3$ and $\mathbb{R}^4$.
\end{Remark}

\noindent
{\bf 8.} In this step we shall extend the asymptotic analysis of
Theorem~\ref{asymptotic analysis 1} to those
\De{Finite rank Perturbations}~\ref{finite rank perturbations},
which respect the involution $\eta$.
Due to the proof of Lemma~\ref{covariant transformation}
the resolvents $\Op{R}(U,\Bar{U},k,0)$ obey the relations
$$\Op{R}(U,\Bar{U},k+k^+_{\kappa},0)=\Op{R}(U,\Bar{U},k+k^-_{\kappa},0)=
\begin{pmatrix}
\psi_{-k^+_{\kappa}} & 0\\
0 & \psi_{-k^-_{\kappa}}
\end{pmatrix}\comp
\Op{R}(\psi_{-\kappa}U,\psi_{\kappa}\Bar{U},k,0)\comp
\begin{pmatrix}
\psi_{k^-_{\kappa}} & 0\\
0 & \psi_{k^+_{\kappa}}
\end{pmatrix}.$$
This implies that these transformations acts on the subspaces
$\Spa{H}_{\text{\scriptsize\rm sing},z_l}$ as the transformation
$$\Spa{H}_{\text{\scriptsize\rm sing},z_l}\mapsto
\begin{pmatrix}
\psi_{-k^+_{\kappa}} & 0\\
0 & \psi_{-k^-_{\kappa}}
\end{pmatrix}\Spa{H}_{\text{\scriptsize\rm sing},z_l}
=\left\{\begin{pmatrix}
\psi_{-k^+_{\kappa}} & 0\\
0 & \psi_{-k^-_{\kappa}}
\end{pmatrix}\psi\mid\psi\in\Spa{H}_{\text{\scriptsize\rm sing},z_l}
\right\}.$$
If $\Spa{H}_{\text{\scriptsize\rm sing},z_l}$ obeys conditions
\Em{Singularity}~(i)--(ii) and is invariant
under the (projective) involution $\psi\mapsto\Op{J}\Bar{\psi}$,
then in the limit $\|\kappa\|\rightarrow\infty$
the transformed spaces
$\begin{pmatrix}
\psi_{-k^+_{\kappa}} & 0\\
0 & \psi_{-k^-_{\kappa}}
\end{pmatrix}\Spa{H}_{\text{\scriptsize\rm sing},z_l}$ converge
to the kernel of the free Dirac operator. In fact, due to
Theorem~\ref{asymptotic analysis 1} the limit corresponds to a
\De{Singularity of the holomorphic structure}~\ref{singularities of
  the structure}
corresponding to the free Dirac operator.
Moreover the limits are fixed points under the former transformations.
Therefore they are invariant under
$\left(\begin{smallmatrix}
z-z_l & 0\\
0 & \unity
\end{smallmatrix}\right)$ and $\left(\begin{smallmatrix}
\unity & 0\\
0 & \Bar{z}-\Bar{z}_l
\end{smallmatrix}\right)$.
Obviously, the kernel of the free Dirac operator is the only subspace 
$\Spa{H}_{\text{\scriptsize\rm sing},z_l}$,
which fulfills conditions \Em{Singularity}~(i)--(ii)
(for the free Dirac operator), and which is invariant
with respect to the action of these two operators and with respect
to the (projective) involution $\psi\mapsto\Op{J}\Bar{\psi}$.
Now Theorem~\ref {asymptotic analysis 1} carries over to
the \Em{complex Fermi curves} of
\De{Finite rank Perturbations}~\ref{finite rank perturbations}
of $\Op{R}(U,\Bar{U},k,0)$ (with $U\in\banach{2}(\torus)$)
corresponding to spaces
$\Spa{H}_{\text{\scriptsize\rm sing},z_l}$,
which are invariant with respect to the (projective) involution
$\psi\mapsto\Op{J}\Bar{\psi}$. More precisely, in the estimates of the
corresponding eigenfunctions the norms $\|\cdot\|_q$ has to be
replaced by the norms of the corresponding
restrictions to the relative complement of
$\Set{S}_{\varepsilon}=\bigcup\limits_{l=1}^{L}B(z_l,\varepsilon)$.
In particular, the matrices $\Mat{S}_{D}(k)$ corresponding to all these  
\De{Finite rank Perturbations}~\ref{finite rank perturbations}
are indeed meromorphic functions on $\mathbb{C}^2$
(compare with Remark~\ref{existence of perturbations}),
whose polar sets are the corresponding \Em{complex Fermi curves}.

\noindent
{\bf 9.} The methods of the fourth step yield
another explanation of the relation in the definition
\De{Finite rank Perturbations}~\ref{finite rank perturbations}.
Let $\psi$ and $\phi$ be spinors with support
in $\Set{S}_{\varepsilon}$.
The proof of Lemma~\ref{branchpoints} implies the relation
\begin{multline*}
\left\langle\left\langle\Op{R}^{t}(V,W,k,0)\phi,
\psi\right\rangle\right\rangle-
\left\langle\left\langle\phi,
\Op{R}(V,W,k',0)\psi\right\rangle\right\rangle=\\
\frac{1}{2\sqrt{-1}}\int\limits_{\partial\Set{S}_{\varepsilon}}
\left\langle\left\langle\Op{R}^{t}(V,W,k,0)\phi,
\begin{pmatrix}
0 & d\Bar{z}\\
dz & 0
\end{pmatrix}\Op{R}(V,W,k',0)\psi\right\rangle\right\rangle.
\end{multline*}
The analog of this relation for the perturbed resolvents
imply the relation in the definition
\De{Finite rank Perturbations}~\ref{finite rank perturbations}.
In fact, with the help of the row--vectors $\Psi_{D}(z)$ and the
column--vectors $\Phi_{D}(z)$ introduced in the sixth step,
the entries of $\Psi_{D}$ and $\Phi_{D}$ can be distinguished by suitable
$\psi$ and $\phi$. Therefore the analog to the former relation
implies the latter.

\noindent
{\bf 10.} In general the limits of the resolvents
$\Op{R}(U_n,\Bar{U}_n,k,\lambda)$
depend non--continuously on $\lambda$, and do not obey the
\Em{first resolvent formula} \cite[Theorem~VI.5]{RS1}.
But as we already saw in
Section~\ref{subsection spectral projections},
some modifications of the Dirac operators still reflect
properties of \Em{complex Fermi curves}. In fact, since
the whole spectrum  of the family of operators
$\left(\begin{smallmatrix}
0 & \pm{\unity}\\
{\unity} & 0
\end{smallmatrix}\right)\comp\Op{D}(V,W,k)$
is determined by the \Em{complex Fermi curves}, the limits of the
corresponding resolvents might satisfy the corresponding
\Em{first resolvent formula}. More precisely, due to the relation
\begin{multline*}
\begin{pmatrix}
0 & k'_2-k_2+\sqrt{-1}(k'_1-k_1)\\
k'_2-k_2-\sqrt{-1}(k'_1-k_1) & 0
\end{pmatrix}\comp\triv{\Op{D}}(V,W,k')=\\
\pi g(k'-k,k'-k)\unity+\begin{pmatrix}
0 & k'_2-k_2+\sqrt{-1}(k'_1-k_1)\\
k'_2-k_2-\sqrt{-1}(k'_1-k_1) & 0
\end{pmatrix}\comp\triv{\Op{D}}(V,W,k)
\end{multline*}
the corresponding resolvents obey the
\cite[Theorem~VI.5]{RS1}

\newtheorem{Resolvent Formula}[Lemma]{Resolvent Formula}
\index{resolvent!$\sim$ formula}
\begin{Resolvent Formula}\label{resolvent formula}
\begin{multline*}
\triv{\Op{R}}(V,W,k,0)-\triv{\Op{R}}(V,W,k',0)=\\
\pi\triv{\Op{R}}(V,W,k,0)\comp\begin{pmatrix}
0 & k_2-k'_2+\sqrt{-1}(k_1-k'_1)\\
k_2-k'_2-\sqrt{-1}(k_1-k'_1) & 0
\end{pmatrix}\comp\triv{\Op{R}}(V,W,k',0).
\end{multline*}
\end{Resolvent Formula}

At a first look condition \De{Finite rank Perturbations}~(iii)
and the corresponding relation in the definition 
\De{Finite rank Perturbations}~\ref{finite rank perturbations}
seems to be equivalent to the analog of the
\De{Resolvent Formula}~\ref{resolvent formula}
for these 
\De{Finite rank Perturbations}~\ref{finite rank perturbations}.
However, they do not obey the
\Em{Resolvent Formula}~\ref{resolvent formula}.
More precisely, if the support of $\triv{\psi}$ and $\triv{\phi}$
are disjoint with small balls $B(z_1,\varepsilon),\ldots,B(z_L,\varepsilon)$
around the support of the \De{Finite rank Perturbations},
then the action of $\triv{\Op{D}}(V,W,k)$ and
$\triv{\Op{D}}^{t}(V,W,k)$ on
$\triv{\Op{R}}(V,W,k,0)\triv{\psi}$
and $\triv{\Op{R}}^{t}(V,W,k,0)\triv{\phi}$ yields
functions, which vanish on these balls, respectively.
We decompose the integrals appearing in the matrix element
\begin{multline*}
\left\langle\left\langle\triv{\phi},\triv{\Op{R}}^{t}(V,W,k',0)
\comp\left(\begin{smallmatrix}
0 & k'_2-k_2+\sqrt{-1}(k'_1-k_1)\\
k'_2-k_2-\sqrt{-1}(k'_1-k_1) & 0
\end{smallmatrix}\right)\comp
\triv{\Op{R}}(V,W,k,0)\triv{\psi}\right\rangle\right\rangle\\
=2\pi\sqrt{-1}\left\langle\left\langle
\triv{\Op{R}}^{t}(V,W,k',0)\triv{\phi},
\triv{\Op{R}}(V,W,k,0)\triv{\psi}
\right\rangle\right\rangle_{(k'_2-k_2,k_1-k'_1)}
\end{multline*}
of the right hand side of the
\De{Resolvent Formula}~\ref{resolvent formula}
into an integral over the complements of the small balls
$B(z_1,\varepsilon)\cup\ldots\cup B(z_L,\varepsilon)$
and an integral over these balls. Then,
due to Lemma~\ref{general branchpoints},
the latter integral is equal to an integral
of a one--form over the boundaries of these balls.
An easy calculation shows that
for all such $\triv{\psi}$ and $\triv{\phi}$ and all
\De{Finite rank Perturbations}~\ref{finite rank perturbations}
the corresponding matrix elements of the left hand side of the
\De{Resolvent Formula}~\ref{resolvent formula} coincides with
the corresponding matrix elements of the right hand side, if we
replace the integral over the balls
$B(z_1,\varepsilon),\ldots,B(z_l,\varepsilon)$ by the integral of the
corresponding one--form over the boundaries of these balls
(described in Lemma~\ref{general branchpoints}).
Furthermore, this condition does not result
in any restriction on the functions $\Mat{S}(k)$.
However, in the limit $\varepsilon\downarrow 0$ these boundary terms
do not tend to zero for all $k$.

We shall see that the limits of resolvents are such
\De{Finite rank Perturbations}~\ref{finite rank perturbations}.
In a first step we investigate the limits $\parameter{t}\rightarrow\infty$
of the resolvents of a family
$\left(U_{\parameter{t}}\right)_{\parameter{t}\in [1,\infty)}$
of potentials of the form
$U_{\parameter{t}}(z)=\parameter{t}U_{\parameter{t}=1}(\parameter{t}z)$,
where $U_{\parameter{t}=1}$ has support in $B(0,\varepsilon)$.
Therefore, the potentials $U_{\parameter{t}}$ have support in
$B(0,\parameter{t}^{-1}\varepsilon)$ and the $\banach{2}$--norm
does not depend on $\parameter{t}$.
Moreover, for all $U_0\in \banach{2}(\torus)$ the measures
$(U_0+U_{\parameter{t}})(\Bar{U}_0+\Bar{U}_{\parameter{t}})d^2x$
converge in the limit $\parameter{t}\rightarrow\infty$
to the measure $U_0\Bar{U}_0d^2x$ plus the point
measure at $z=0$ with mass equal to the square of the $\banach{2}$--norm of
$U_{\parameter{t}=1}$.
If the $\banach{2}$--norm of $U_{\parameter{t}=1}$ is smaller
than the constant $S_p^{-1}$ and $k$ does not belong to the
\Em{complex Fermi curve} $\fermi(U_0,\Bar{U}_0)$,
then Lemma~\ref{weakly continuous resolvent} implies that the family
of resolvents
$\Op{R}(U_0+U_{\parameter{t}},\Bar{U}_0+\Bar{U}_{\parameter{t}},k,0)$
converges to the resolvent $\Op{R}(U_0,\Bar{U}_0,k,0)$.
In general we shall see that this limit depends on the inverse of the
operators $\unity-\left(\begin{smallmatrix}
U_{\parameter{t}=1} & \partial\\
-\Bar{\partial} & \Bar{U}_{\parameter{t}=1}
\end{smallmatrix}\right)\comp\Op{R}_{\mathbb{R}^2}(0,0,0)$, where
$\Op{R}_{\mathbb{R}^2}(0,0,0)$ denotes the resolvent of the free
Dirac operator on $\mathbb{R}^2$. In a next step
the spectral theory of Dirac operators on $\mathbb{P}^1$ is used
in order to investigate these inverse operators.

\subsubsection{Spectral theory of Dirac operators on the sphere}
\label{subsubsection spectral theory on the sphere}

For this purpose we use the covering
$\mathbb{P}^1=\mathbb{C}\cup\mathbb{P}^1\setminus\{0\}$ and
the affine coordinates $z$ and $1/z$ on the members of this open
covering. The spin bundle on $\mathbb{P}^1$ may be described as the
trivial $\mathbb{C}^2$ bundles on both members of this covering
together with the transition function
$\left(\begin{smallmatrix}
z & 0\\
0 & -\Bar{z}
\end{smallmatrix}\right)$, which transforms the sections on
$\mathbb{C}$ into sections on $\mathbb{P}^1\setminus\{0\}$.
With respect to these coordinates the free Dirac operator
\cite[Section~3.4]{Fr1} is equal to
\begin{align*}
&(1+z\Bar{z})\begin{pmatrix}
0 & \partial_{z}\\
-\Bar{\partial}_{z} & 0
\end{pmatrix}
\text{ on $\mathbb{C}$} &\text{and equal to}&\\
(1+z\Bar{z})\begin{pmatrix}
z & 0\\
0 & -\Bar{z}
\end{pmatrix}\comp
\begin{pmatrix}
0 & \partial_{z}\\
-\Bar{\partial}_{z} & 0
\end{pmatrix}\comp
\begin{pmatrix}
z & 0\\
0 & -\Bar{z}
\end{pmatrix}^{-1}&=
(1+1/(z\Bar{z}))\begin{pmatrix}
0 & \partial_{1/z}\\
-\Bar{\partial}_{1/z} & 0
\end{pmatrix}&\text{on $\mathbb{P}^1\setminus\{0\}$}.&
\end{align*}
We conclude that if $(1+z\Bar{z})U_{z}$ denotes the potential on
$\mathbb{C}$, the corresponding potential on
$\mathbb{P}^1\setminus\{0\}$ is equal to
$(1+1/(z\Bar{z}))U_{1/z}$, where $U_{1/z}$ is given by
$$U_{1/z}(1/z)=\frac{1+z\Bar{z}}{1+1/(z\Bar{z})}
U_{z}(z)=z\Bar{z}U_{z}(z).$$
The corresponding Dirac operators are given by
\begin{align*}
(1+z\Bar{z})\begin{pmatrix}
U_{z} & \partial_{z}\\
-\Bar{\partial}_{z} & \Bar{U}_{z}
\end{pmatrix}
&\text{ on $\mathbb{C}$ and}
&(1+1/(z\Bar{z}))\begin{pmatrix}
U_{1/z} & \partial_{1/z}\\
-\Bar{\partial}_{1/z} & \Bar{U}_{1/z}
\end{pmatrix} & \text{ on $\mathbb{P}^1\setminus\{0\}$.}
\end{align*}
In particular, the integral over the density
$$U_{z}\Bar{U}_{z}\frac{d\Bar{z}\wedge dz}{2\sqrt{-1}}=
U_{1/z}\Bar{U}_{1/z}\frac{d(1/\Bar{z})\wedge d(1/z)}{2\sqrt{-1}}$$
is equal to the $\banach{2}$--norm of $U_{z}$ on $\mathbb{C}$
(or of $U_{1/z}$ on $\mathbb{P}^1\setminus\{0\}$).
Therefore, all potentials $U\in\banach{2}(\mathbb{R}^2)$ may be considered as
potentials $U_{z}$ of Dirac operators on $\mathbb{P}^1$.
Due to Dolbeault's Lemma
\cite[Chapter~I Section~D 2.~Lemma]{GuRo}
the integral kernel of the corresponding free resolvent is equal to
\begin{align*}
&\begin{pmatrix}
0 & (z-z')^{-1}\\
(\Bar{z}'-\Bar{z})^{-1} & 0
\end{pmatrix}\frac{d\Bar{z}'\wedge dz'}{2\pi\sqrt{-1}(1+z'\Bar{z}')}
\text{ on $\mathbb{C}$} &\text{and equal to}&\\
&\begin{pmatrix}
z & 0\\
0 & -\Bar{z}
\end{pmatrix}\comp
\begin{pmatrix}
0 & (z-z')^{-1}\\
(\Bar{z}'-\Bar{z})^{-1} & 0
\end{pmatrix}\comp
\begin{pmatrix}
z' & 0\\
0 & -\Bar{z}'
\end{pmatrix}^{-1}\frac{d\Bar{z}'\wedge dz'}{2\pi\sqrt{-1}(1+z'\Bar{z}')}&&\\
&=\begin{pmatrix}
0 & (1/z-1/z')^{-1}\\
(1/\Bar{z}'-1/\Bar{z})^{-1} & 0
\end{pmatrix}\frac{d(1/\Bar{z}')\wedge d(1/z')}
                  {2\pi\sqrt{-1}(1+1/(z'\Bar{z}'))}
&\text{ on $\mathbb{P}^1\setminus\{0\}$}.&
\end{align*}
We define $\banach{p}$--spinors as
$\mathbb{C}^2$--valued functions on $\mathbb{C}$,
whose components have bounded norm
$$\|f\|=
\left(\int\limits_{\mathbb{R}^2}|f(z)|^p(1+z\Bar{z})^{(\frac{p}{2}-2)}
d^2x\right)^{1/p}.$$
Due to the transformation property of spinors this norm for spinors on
$\mathbb{C}$ is equal to the norm of the corresponding spinors on
$\mathbb{P}^1\setminus\{0\}$  with respect to the coordinates $1/z$.
Since the multiplication with $(1+z\Bar{z})^{\frac{1}{2}-\frac{2}{p}}$ is an
isometry from the space of $\banach{p}$--spinors onto
$\banach{p}(\mathbb{R}^2)\times\banach{p}(\mathbb{R}^2)$, the free resolvent,
$\Op{R}_{\mathbb{P}^1}(0,0,0)$
considered as an operator from
$\banach{p}(\mathbb{R}^2)\times\banach{p}(\mathbb{R}^2)$
($\simeq \banach{p}$--spinors) into 
$\banach{q}(\mathbb{R}^2)\times\banach{q}(\mathbb{R}^2)$
($\simeq \banach{q}$--spinors), has the integral kernel
$$\begin{pmatrix}
0 & (z-z')^{-1}\\
(\Bar{z}'-\Bar{z})^{-1} & 0
\end{pmatrix}(1+z\Bar{z})^{\frac{1}{2}-\frac{2}{q}}
(1+z'\Bar{z}')^{\frac{2}{p}-\frac{3}{2}}
\frac{d^2x'}{\pi}.$$
If $\frac{1}{q}=\frac{1}{p}-\frac{1}{2}$,
then the potential $U$ represents the potential
$(1+z\Bar{z})U$ on $\mathbb{P}^1$ considered as an operator from
$\banach{q}(\mathbb{R}^2)\times\banach{q}(\mathbb{R}^2)$
($\simeq \banach{q}$--spinors) into
$\banach{p}(\mathbb{R}^2)\times\banach{p}(\mathbb{R}^2)$
($\simeq \banach{p}$--spinors).
Therefore, in this case the resolvent of the Dirac operator
on $\mathbb{P}^1$ with potential
$(1+z\Bar{z})\left(\begin{smallmatrix}
U & 0\\
0 & \Bar{U}
\end{smallmatrix}\right)$ considered as an operator from
$\banach{p}(\mathbb{R}^2)\times\banach{p}(\mathbb{R}^2)$
($\simeq \banach{p}$--spinors) into 
$\banach{q}(\mathbb{R}^2)\times\banach{q}(\mathbb{R}^2)$
($\simeq \banach{q}$--spinors) is equal to
$$\Op{R}_{\mathbb{P}^1}(U,\Bar{U},\lambda)=
\Op{R}_{\mathbb{P}^1}(0,0,\lambda)\comp\left(\unity-\begin{pmatrix}
U & 0\\
0 & \Bar{U}
\end{pmatrix}\comp\Op{R}_{\mathbb{P}^1}(0,0,\lambda)\right)^{-1}.$$
Fortunately, if $p=\frac{4}{3}$ and $q=4$, the resolvent
$\Op{R}_{\mathbb{P}^1}$ considered as an operator from
$\banach{\frac{4}{3}}(\mathbb{R}^2)\times
\banach{\frac{4}{3}}(\mathbb{R}^2)$
into $\banach{4}(\mathbb{R}^2)\times
\banach{4}(\mathbb{R}^2)$ is equal to
the resolvent $\Op{R}_{\mathbb{R}^2}(0,0,0)$
of the free Dirac operator on $\mathbb{R}^2$.

In order to carry over the spectral theory developed in
Section~\ref{subsection resolvent} we shall estimate the free
resolvent $\Op{R}_{\mathbb{P}^1}(0,0,\sqrt{-1}\lambda)$ in the
limit $\lambda\rightarrow\infty$. Due to the invariance of the
Dirac operator under the action of $SU(2,\mathbb{C})$ on
$\mathbb{P}^1$, it suffices to estimate the Green's function
$\left(\begin{smallmatrix}
a_{\mathbb{P}^1,\lambda}(x) & b_{\mathbb{P}^1,\lambda}(x)\\
c_{\mathbb{P}^1,\lambda}(x) & d_{\mathbb{P}^1,\lambda}(x)
\end{smallmatrix}\right)$, which satisfies the equation
$$\begin{pmatrix}
\sqrt{-1}\lambda & -(1+z\Bar{z})\partial\\
(1+z\Bar{z})\Bar{\partial} & \sqrt{-1}\lambda
\end{pmatrix}\begin{pmatrix}
a_{\mathbb{P}^1,\lambda}(x) & b_{\mathbb{P}^1,\lambda}(x)\\
c_{\mathbb{P}^1,\lambda}(x) & d_{\mathbb{P}^1,\lambda}(x)
\end{pmatrix}=\begin{pmatrix}
\delta(x) & 0\\
0 & \delta(x)
\end{pmatrix}.$$
Here $\delta(x)$ denotes the two--dimensional $\delta$--function on
$\mathbb{R}^2$, i.\ e.\ $\delta(x)d^2x$ denotes the point measure
at $x=0$ of mass $1$. Due to the relation
\begin{multline*}\begin{pmatrix}
\sqrt{-1}\lambda & -(1+z\Bar{z})\partial\\
(1+z\Bar{z})\Bar{\partial} & \sqrt{-1}\lambda
\end{pmatrix}\comp\begin{pmatrix}
\sqrt{-1}\lambda & (1+z\Bar{z})\partial\\
-(1+z\Bar{z})\Bar{\partial} & \sqrt{-1}\lambda
\end{pmatrix}\\=\begin{pmatrix}
-\lambda^2+(1+z\Bar{z})^2\partial\Bar{\partial}+
(1+z\Bar{z})\Bar{z}\Bar{\partial} & 0\\
0 & -\lambda^2+(1+z\Bar{z})^2\Bar{\partial}\partial+
(1+z\Bar{z})z\partial
\end{pmatrix}\end{multline*}
it suffices to calculate the Green's functions of the
two operators in the diagonal entries of the right hand side. These
Green's functions depend only on $r=\sqrt{z\Bar{z}}=\sqrt{x_1^2+x_2^2}$.
In fact, on functions depending only on $r$, both operators in the
diagonal entries of the right hand side act as the operator
$$-\lambda^2+\left(\frac{1+r^2}{2}\right)^2
\left(\frac{d^2}{dr^2}+\frac{1}{r}\frac{d}{dr}\right)+
\frac{1+r^2}{2}r\frac{d}{dr}.$$
The substitution $r=\tan(\varphi/2)$ transforms this operator into
$$-\lambda^2+\frac{d^2}{d\varphi^2}+
\frac{1}{\sin(\varphi)}\frac{d}{d\varphi}.$$
The two dimensional $\delta$--function is equal to the one--dimensional
$\delta$--function with respect to the measure $\sin(\varphi)d\varphi$.
Therefore, we shall calculate the unique Green's function
\index{Green's function!$\sim$ $\Func{G}_{\mathbb{P}^1,\lambda}$}
$\Func{G}_{\mathbb{P}^1,\lambda}$ depending on
$\varphi\in\mathbb{R}/(2\pi\mathbb{Z})$, which satisfies
$$-\lambda^2\sin(\varphi)\Func{G}_{\mathbb{P}^1,\lambda}(\varphi)+
\sin(\varphi)\Func{G}_{\mathbb{P}^1,\lambda}''(\varphi)+
\Func{G}_{\mathbb{P}^1,\lambda}'(\varphi)=\delta(\varphi),$$
where $\delta$ denotes the one--dimensional $\delta$--function on
$\varphi\in\mathbb{R}/(2\pi\mathbb{Z})$ with respect to the measure
$d\varphi$. Obviously this Green's function
$\Func{G}_{\mathbb{P}^1,\lambda}$
is odd and has a Fourier series of the form
$$\Func{G}_{\mathbb{P}^1,\lambda}(\varphi)=
\sum\limits_{n\in\mathbb{N}} a_n\sin(n\varphi).$$
Since the Fourier series of the $\delta$--functions is equal to
$$\delta(\varphi)=\frac{1}{2\pi}+\frac{1}{\pi}
\sum\limits_{n\in\mathbb{N}} \cos(n\varphi),$$
the coefficients $a_n$ satisfiy the recursion relation
$$a_{n+1}=\frac{2na_{n}+(\lambda^2+(n-1)^2)a_{n-1}-2/\pi}
               {\lambda^2+(n+1)^2}
\text{ for }n\in\mathbb{N}
\text{ with }a_1=-\frac{1}{\pi(\lambda^2+1)}\text{ and }a_0=0.$$

\begin{Lemma}\label{free sphere resolvent}
If $\lambda\in\mathbb{R}$, then the functions
$(1-\cos(\varphi))\Func{G}_{\mathbb{P}^1,\lambda}(\varphi)$
and $(1-\cos(\varphi))\Func{G}_{\mathbb{P}^1,\lambda}'(\varphi)$
are continuous functions depending on
$\varphi\in\mathbb{R}/(2\pi\mathbb{Z})$, which
converge uniformly to zero in the limit $\lambda\rightarrow\infty$.
\end{Lemma}

\begin{proof}
The functions $(1-\cos(\varphi))\Func{G}_{\mathbb{P}^1,\lambda}(\varphi)$
and $(1-\cos(\varphi))\Func{G}_{\mathbb{P}^1,\lambda}'(\varphi)$
have the Fourier series
\begin{eqnarray*}
(1-\cos(\varphi))\Func{G}_{\mathbb{P}^1,\lambda}(\varphi)&=&
\sum\limits_{n\in\mathbb{N}}
\frac{2a_n-a_{n-1}-a_{n+1}}{2}\sin(n\varphi)\text{ and}\\
(1-\cos(\varphi))\Func{G}_{\mathbb{P}^1,\lambda}'(\varphi)&=&
-\frac{a_1}{2}+\sum\limits_{n\in\mathbb{N}}
\frac{2na_n-(n-1)a_{n-1}-(n+1)a_{n+1}}{2}\cos(n\varphi),
\end{eqnarray*}
respectively. It suffices to show that the sum of
$|2a_n-a_{n-1}-a_{n+1}|$ and the sum of $|2na_n-(n-1)a_{n-1}-(n+1)a_{n+1}|$ 
exist and converge in the limit $\lambda\rightarrow\infty$ to zero.
We shall estimate the absolute values of the coefficients $a_n$.

\noindent
{\bf 1.} First we claim that for all $1<\alpha<2$ there exists a constant
$C_{\alpha}>0$, such that all absolute values $|a_n|$ are bounded by
$C_{\alpha}\frac{n^\alpha}{\lambda^2+n^2}.$ In fact, due to the
recursion relation it suffices to show the inequalities
\begin{multline*}
C_{\alpha}\frac{2n^{(1+\alpha)}}{\lambda^2+n^2}+
C_{\alpha}(n-1)^{\alpha}+\frac{2}{\pi}\leq C_{\alpha}(n+1)^{\alpha}
\;\forall n\in\mathbb{N}\setminus\{1\},\\
|a_1|\leq C_{\alpha}\frac{1}{\lambda^2+1}\text{ and }
|a_2|\leq C_{\alpha}\frac{2^{\alpha}}{\lambda^2+4}.
\end{multline*}
Since $a_2=-\frac{2}{\pi(\lambda^2+4)}$, the latter inequalities are
equivalent to $\frac{1}{\pi}\leq C_{\alpha}$ and
$\frac{2}{\pi}\leq C_{\alpha}2^{\alpha}$.
Obviously we have for all $1<\alpha<2$ and all $n\in\mathbb{N}$
\begin{multline*}
(n-1)^{\alpha}\leq
\left(1-\frac{2}{n+1}\right)^{\alpha}(n+1)^{\alpha}\leq
\left(1-\frac{2}{n+1}\right)
\left(1-\frac{2(\alpha-1)}{n+1}\right)(n+1)^{\alpha}\leq\\
(n+1)^{\alpha}-2\alpha(n+1)^{\alpha-1}+
4(\alpha-1)(n+1)^{\alpha-2}.\end{multline*}
If we insert this inequality in the expression on the left hand side
of the first inequality above, we obtain
\begin{multline*}
C_{\alpha}\frac{2n^{(1+\alpha)}}{\lambda^2+n^2}+
C_{\alpha}(n-1)^{\alpha}+\frac{2}{\pi}\leq\\
\leq 2C_{\alpha}(n+1)^{\alpha-1}+
C_{\alpha}\left((n+1)^{\alpha}-2\alpha(n+1)^{\alpha-1}+
4(\alpha-1)(n+1)^{\alpha-2}\right)+\frac{2}{\pi}.
\end{multline*}
Therefore, it suffices to establish the inequalities
\begin{align*}
\frac{1}{\pi}&\leq C_{\alpha}(\alpha-1)(n+1)^{\alpha-1}
\left(1-\frac{2}{n+1}\right)
\;\forall n\in\mathbb{N}\setminus\{1\},&
\frac{1}{\pi}\leq C_{\alpha}\text{ and }&
\frac{2}{\pi}\leq C_{\alpha}2^{\alpha}.
\end{align*}
Obviously there exist such $C_{\alpha}<\infty$ for all $1<\alpha<2$.

\noindent
{\bf 2.} In a second step we claim that for all $1<\alpha<2$
there exists a constant $C_{\alpha}'$,
such that the absolute values $|a_{n+1}-a_{n-1}|$ are bounded by
$C_{\alpha}'\frac{n^{\alpha+1}}
                 {(\lambda^2+(n+1)^2)(\lambda^2+(n-1)^2)}$.
In fact, the recursion relation implies
$$|a_{n+1}-a_{n-1}|\leq\frac{2n}{\lambda^2+(n+1)^2}
(|a_{n}|+2|a_{n-1}|).$$
Hence the second claim follows from what has already been proven,
with $C_{\alpha}'=6C_{\alpha}$.

\noindent
{\bf 3.} In a third step we prove that for all $1<\alpha<2$ and
$|\lambda|\geq\sqrt{2}$
(this assumption on $\lambda$ is not essential and can be removed)
there exists a constant $C_{\alpha}''$
such that the absolute values
$|2a_n-a_{n+1}-a_{n-1}|$ are bounded by
$C_{\alpha}''\frac{n^{\alpha}}
                 {(\lambda^2+(n+1)^2)(\lambda^2+(n-1)^2)}$.
To see this set $b_n=2a_n-a_{n+1}-a_{n-1}$ for $n\in\mathbb{N}$.
If we subtract from two times
the recursion relation above the two relations,
which are obtained by reducing and enlarging the indices by one,
then we obtain the following recursion relation:
\begin{multline*}
b_{n+1}=
\frac{-2nb_{n}+(\lambda^2+(n-2)^2)b_{n-1}+4(a_{n+1}-a_{n-1})}
     {\lambda^2+(n+2)^2}\;\forall n\in\mathbb{N}\setminus\{1\}\\
\text{with }
b_1=-\frac{6}{\pi(\lambda^2+1)(\lambda^2+4)}\text{ and }
b_2=-\frac{4(\lambda^2-5)}{\pi(\lambda^2+1)(\lambda^2+4)(\lambda^2+9)}.
\end{multline*}
Therefore, it suffices to show for all $n\in\mathbb{N}\setminus\{1\}$
$$C_{\alpha}''\frac{2n^{(\alpha+1)}}
                   {(\lambda^2+(n-1)^2)(\lambda^2+(n+1)^2)}+
C_{\alpha}''\frac{(n-1)^{\alpha}}{\lambda^2+n^2}+4|a_{n+1}-a_{n-1}|\leq
C_{\alpha}''\frac{(n+1)^{\alpha}}{\lambda^2+n^2}$$
together with the two inequalities $\frac{6}{\pi}\leq C_{\alpha}''$
and $\frac{4|\lambda^2-5|}{\pi(\lambda^2+4)}\leq C_{\alpha}''2^{\alpha}$.
If $\lambda^2\geq2$, then we have
$$\frac{(\lambda^2+n^2)n^2}{(\lambda^2+(n-1)^2)(\lambda^2+(n+1)^2)}\leq
1-\frac{\lambda^4+\lambda^2(n^2+2)-2n^2+1}
       {(\lambda^2+(n-1)^2)(\lambda^2+(n+1)^2)}\leq 1.$$
With the help of the estimates in the first and the second step
we arrive at the estimates
\begin{multline*}
2C_{\alpha}'n^{\alpha-1}\leq 2C_{\alpha}'(n+1)^{\alpha-1}
\leq C_{\alpha}''(\alpha-1)(n+1)^{\alpha-1}
\left(1-\frac{2}{n+1}\right)
\;\forall n\in\mathbb{N}\setminus\{1\},\\
\frac{6}{\pi}\leq C_{\alpha}''\text{ and }
\frac{4|\lambda^2-5|}{\pi(\lambda^2+4)}\leq C_{\alpha}''2^{\alpha}.
\end{multline*}
Again it is obvious that there exist such $C_{\alpha}''$.

Due to the proven third claim the sum $\sum\limits_{n\in\mathbb{N}}|b_n|$
exists and converges in the limit $\lambda\rightarrow\infty$ to zero.
Due to the proven second and the third claim the sum
$$\sum\limits_{n\in\mathbb{N}}|2na_n-(n+1)a_{n+1}-(n-1)a_{n-1}|\leq
\sum\limits_{n\in\mathbb{N}}n|b_n|+
\sum\limits_{n\in\mathbb{N}}|a_{n+1}-a_{n-1}|$$
exists and converges in the limit $\lambda\rightarrow\infty$ to zero.
\end{proof}

The Green's functions of the Dirac operator on
$\mathbb{P}^1$ may be expressed in terms of
$\Func{G}_{\mathbb{P}^1,\lambda}$ as
$$\begin{pmatrix}
a_{\mathbb{P}^1,\lambda}(x) & b_{\mathbb{P}^1,\lambda}(x)\\
c_{\mathbb{P}^1,\lambda}(x) & d_{\mathbb{P}^1,\lambda}(x)
\end{pmatrix}=\begin{pmatrix}
\sqrt{-1}\lambda
\Func{G}_{\mathbb{P}^1,\lambda}(2\arctan(\sqrt{z\Bar{z}})) &
(1+z\Bar{z})\partial
\Func{G}_{\mathbb{P}^1,\lambda}(2\arctan(\sqrt{z\Bar{z}}))\\
-(1+z\Bar{z})\Bar{\partial}
\Func{G}_{\mathbb{P}^1,\lambda}(2\arctan(\sqrt{z\Bar{z}})) &
\sqrt{-1}\lambda
\Func{G}_{\mathbb{P}^1,\lambda}(2\arctan(\sqrt{z\Bar{z}}))
\end{pmatrix}.$$
Hence, due to the invariance under the action of $SU(2,\mathbb{C})$ on
$\mathbb{P}^1$, this lemma implies that outside of the diagonal
the integral kernel of the resolvent $\left(\begin{smallmatrix}
\sqrt{-1}\lambda & -(1+z\Bar{z})\partial\\
(1+z\Bar{z})\Bar{\partial} & \sqrt{-1}\lambda
\end{smallmatrix}\right)^{-1}$ is a continuous
function, which converges in the limit $\lambda\rightarrow\infty$
on the complement of arbitrary small discs
(with respect to the invariant metric of $\mathbb{P}^1$)
uniformly to zero. Therefore, the arguments of Theorem~\ref{meromorph}
and Lemma~\ref{weakly continuous resolvent} carry over to the
resolvents of Dirac operators on $\mathbb{P}^1$.

\begin{Corollary}\label{weakly continuous sphere resolvents}
The resolvents define meromorphic maps $(U,\lambda)\mapsto
\Op{R}_{\mathbb{P}^1}(U,\Bar{U},\lambda)$ from the Hilbert space
$\banach{2}(\mathbb{R}^2)\times\mathbb{C}$ into the compact
operators on the $\banach{2}$--spinors.
On all subsets of $\banach{2}(\mathbb{R}^2)$,
where for some $C_p<S_p^{-1}$ and $\varepsilon>0$
the $\banach{2}$--norms of the restrictions
of the potentials to all $\varepsilon$--balls of
$\mathbb{C}\subset\mathbb{P}^1$
(with respect to the natural metric of $\mathbb{P}^1$)
are not larger than $C_p$,
this map is weakly continuous.
Finally, for $1<p<2$ the operators
$\unity-\left(\begin{smallmatrix}
U & 0\\
0 & \Bar{U}
\end{smallmatrix}\right)\comp
\Op{R}_{\mathbb{P}^1}(0,0,\lambda)$ on
$\banach{p}(\mathbb{R}^2)\times\banach{p}(\mathbb{R}^2)$
are Fredholm operators of index zero,
and invertible if and only if $\lambda$ is no eigenvalue
of the corresponding Dirac operator on $\mathbb{P}^1$.
\end{Corollary}

\begin{proof}
It remains to prove the last statement. Due to the first resolvent formula
\cite[Theorem~VI.5]{RS1} for all vectors $\psi$
and all $\lambda$ and $\lambda'$ in the complement of the spectrum
the corresponding resolvents $\Op{R}_{\lambda}$ and
$\Op{R}_{\lambda'}$ satisfy the equation
$$\Op{R}_{\lambda'}\left(\psi+
(\lambda'-\lambda)\Op{R}_{\lambda}\psi\right)=\Op{R}_{\lambda}\psi.$$
We conclude that the range of the resolvent $\Op{R}_{\lambda}$
does not depend on $\lambda$ as long as $\lambda$ does not belong to
the spectrum.  Since the
Riesz--Schauder Theorem is also true for compact operators on
Banach spaces \cite[Chapter~X.5.]{Yo}, the
Dirac operators on $\mathbb{P}^1$ may be considered as
closed operators on the space of $\banach{p}$--spinors with compact resolvents.
The inverse mapping theorem \cite[Theorem~III.11]{RS1} together with
the formula
$$\Op{R}_{\mathbb{P}^1}(U,\Bar{U},\lambda)=
\Op{R}_{\mathbb{P}^1}(0,0,\lambda)\comp\left(\unity-\begin{pmatrix}
U & 0\\
0 & \Bar{U}
\end{pmatrix}\comp\Op{R}_{\mathbb{P}^1}(0,0,\lambda)\right)^{-1}$$
for the resolvents of these Dirac operators acting on the
$\banach{p}$--spinors shows that the operator $\unity-\left(\begin{smallmatrix}
U & 0\\
0 & \Bar{U}
\end{smallmatrix}\right)\comp
\Op{R}_{\mathbb{P}^1}(0,0,\lambda)$
on the $\banach{p}$ spinors is invertible if and only if $\lambda$ does not
belong to the spectrum of the corresponding
Dirac operator. This follows also from the following variant
of the first resolvent formula \cite[Theorem~VI.5]{RS1}:
$$(\lambda-\lambda')\Op{R}_{\lambda'}
=\Op{R}_{\lambda}\comp
\left(\frac{\unity}{\lambda-\lambda'}-\Op{R}_{\lambda}\right)^{-1}
=\left(\frac{\unity}{\lambda-\lambda'}-\Op{R}_{\lambda}\right)^{-1}
\comp\Op{R}_{\lambda},$$
since, due to the arguments of the proof of
Lemma~\ref{weakly continuous resolvent} for all
$U\in\banach{2}(\mathbb{R}^2)$ and large $\lambda$ the operators
$\unity-\left(\begin{smallmatrix}
U & 0\\
0 & \Bar{U}
\end{smallmatrix}\right)\comp
\Op{R}_{\mathbb{P}^1}(0,0,\sqrt{-1}\lambda)$ on
$\banach{p}(\mathbb{R}^2)\times\banach{p}(\mathbb{R}^2)$
are invertible. Since the range of the spectral projections
$\Breve{\Op{P}}$ are contained in the domain, the operators
$$\Breve{\Op{Q}}=\Op{R}_{\mathbb{P}^1}^{-1}(0,0,\lambda)\comp
\Breve{\Op{P}}\comp\Op{R}_{\mathbb{P}^1}(0,0,\lambda)$$
are well defined projections, whenever $\lambda$ does not belong to
the spectrum of the free Dirac operator. In general the operators
$\unity-\left(\begin{smallmatrix}
U & 0\\
0 & \Bar{U}
\end{smallmatrix}\right)\comp
\Op{R}_{\mathbb{P}^1}(0,0,\lambda)$
are invertible operators from the range of the projection
$\unity-\Breve{\Op{Q}}$ onto the range of the projection
$\unity-\Breve{\Op{P}}$. Since the spectral projections of the
Dirac operators on $\mathbb{P}^1$ are finite--dimensional, these
operators are Fredholm operators of index zero, and the kernel and the
co--kernel coincides with the corresponding eigenspaces and transposed
eigenspaces, respectively.
\end{proof}

\begin{Remark}\label{fredholm on the plane}
In particular, for $p=\frac{4}{3}$ and $q=4$
and for all $U\in \banach{2}(\mathbb{R}^2)$ the operators
$\unity-\left(\begin{smallmatrix}
U & 0\\
0 & \Bar{U}
\end{smallmatrix}\right)\comp
\Op{R}_{\mathbb{R}^2}(0,0,0)$ on
$\banach{\frac{4}{3}}(\mathbb{R}^2)\times\banach{\frac{4}{3}}(\mathbb{R}^2)$
are also Fredholm operators of index zero, and invertible,
if the corresponding Dirac operator on $\mathbb{P}^1$
has a trivial kernel.
The relation of these kernels to the spectral theory
of the Dirac operators on $\mathbb{R}^2$ is not so obvious.
\end{Remark}

\subsubsection{An illustrative family of potentials}
\label{subsubsection family}

Now we are prepared to treat the limits $\parameter{t}\rightarrow\infty$
of the resolvents $\Op{R}(U_{\parameter{t}},\Bar{U}_{\parameter{t}},k,0)$
of the family $\left(U_{\parameter{t}}\right)_{\parameter{t}\in [1,\infty)}$
of potentials of the form
$U_{\parameter{t}}(z)=\parameter{t}U_{\parameter{t}=1}(\parameter{t}z)$,
where $U_{\parameter{t}=1}$ has support in a ball $B(0,\varepsilon)$
contained in $\Delta$.

\begin{Lemma} \label{nontrivial limit 1}
If the Dirac operator
$(1+z\Bar{z})\left(\begin{smallmatrix}
U_{\parameter{t}=1} & \partial\\
-\Bar{\partial} & \Bar{U}_{\parameter{t}=1}
\end{smallmatrix}\right)$ on $\mathbb{P}^1$
does not have a kernel, then the resolvents
$\Op{R}(U_{\parameter{t}},\Bar{U}_{\parameter{t}},k,0)$
considered as operators from
$\banach{p}(\Delta)\times\banach{p}(\Delta)$ into
$\banach{q'}(\Delta)\times\banach{q'}(\Delta)$
with $1<p\leq 2$ and $2\leq q'<q=\frac{2p}{2-p}$ 
converges in the limit $\parameter{t}\rightarrow \infty$
to $\Op{R}(0,0,k,0)$.
In general, for all $\varepsilon>0$ the restrictions
of these resolvents to the complement of
$B(0,\varepsilon)$ considered as operators from
$\banach{p}(\Delta\setminus B(0,\varepsilon))\times
\banach{p}(\Delta\setminus B(0,\varepsilon))$ into
$\banach{\frac{2p}{2-p}}(\Delta\setminus B(0,\varepsilon))\times
\banach{\frac{2p}{2-p}}(\Delta\setminus B(0,\varepsilon))$ with $1<p<2$
converges in the limit $\parameter{t}\rightarrow\infty$ to a
\Em{Finite rank Perturbation}~\ref{finite rank perturbations}
of $\Op{R}(0,0,k,0)$.
Moreover, the rank of the perturbation is bounded
by the dimension of the kernel of the Dirac operator
$(1+z\Bar{z})\left(\begin{smallmatrix}
U_{\parameter{t}=1} & \partial\\
-\Bar{\partial} & \Bar{U}_{\parameter{t}=1}
\end{smallmatrix}\right)$ on $\mathbb{P}^1$,
and the support of this perturbation is contained in $0\in\Delta$.
\end{Lemma}

\begin{Remark}\label{restriction}
  If $\Set{B}$ is a subset of $\Set{A}$,
  then the functions on $\Set{B}$ may be considered as functions on $\Set{A}$,
  which vanish on the relative complement of $\Set{B}$. Vice versa, the
  restriction to $\Set{B}$ of any function on $\Set{A}$ yields a function on
  $\Set{B}$. For a given operator on a space of functions on $\Set{A}$, we call
  the operator, which associates to any function on $\Set{B}$ the
  restriction to $\Set{B}$ of the action of the operator on this function
  considered as a function on $\Set{A}$, the restriction of the operator to
  $\Set{B}$ or to the space of functions on $\Set{B}$.
\end{Remark}

\begin{proof}
We proof this within three steps. In the first step we investigate the
family of resolvents transformed under the scaling transformations
$z\rightarrow\parameter{t}z$. The limit $\parameter{t}\rightarrow\infty$
of these transformations is a kind of blowing up.
Consequently, it yields a Dirac operator on $\mathbb{P}^1$.
In the second step we show that
if this blown up Dirac operator on $\mathbb{P}^1$ has a trivial kernel,
then the family of resolvents $\Op{R}(U_{\parameter{t}},\Bar{U},k,0)$
converge to the free resolvents $\Op{R}(0,0,k,0)$.
Finally, in the third step we show that a kernel
of the Dirac operator blown up to $\mathbb{P}^1$ yields a
\De{Finite rank Perturbation}~\ref{finite rank perturbations}
of the free resolvents.

\noindent
{\bf 1.\ The family of transformed resolvents.}
For all $\parameter{t}\in\mathbb{C}$ the transformations
$$\Op{I}_{p,\parameter{t}}:f\mapsto\Op{I}_{p,\parameter{t}}f
\text{ with } (\Op{I}_{p,\parameter{t}}f)(z)=
|\parameter{t}|^{2/p}f(\parameter{t}z)$$
define isometries from $\banach{p}(\mathbb{R}^2)$ onto
$\banach{p}(\mathbb{R}^2)$
and from $\banach{p}(B(0,\varepsilon))$ onto
$\banach{p}(B(0,\parameter{t}^{-1}\varepsilon))$.
Due to the explicit form of the integral kernel of
$\Op{R}(0,0,k,0)$ in Lemma~\ref{resolvent integral kernel}
this resolvent is the convolution with a meromorphic and a
anti--meromorphic function, respectively. The Laurent series of these
functions around $z=0$ converge on small discs. Therefore, for all
$k\notin\fermi(0,0)$ the restrictions of the operators
$$\begin{pmatrix}
\Op{I}_{\frac{2p}{2-p},\parameter{t}} & 0\\
0 & \Op{I}_{\frac{2p}{2-p},\parameter{t}}
\end{pmatrix}^{-1}\comp\Op{R}(0,0,k,0)\comp
\begin{pmatrix}
\Op{I}_{p,\parameter{t}} & 0\\
0 & \Op{I}_{p,\parameter{t}}
\end{pmatrix}$$
to $\banach{p}(B(0,\varepsilon))\times\banach{p}(B(0,\varepsilon))$ define
a real analytic family of operators depending on the parameter
$\parameter{t}\in [1,\infty)$.
Furthermore, since in the limit $\parameter{t}\rightarrow\infty$ only the
pole of the convolution function at $z=0$ contributes, this family
converges in this limit to the restriction of the operator
$\Op{R}_{\mathbb{R}^2}(0,0,0)$ to
$\banach{p}(B(0,\varepsilon))\times\banach{p}(B(0,\varepsilon))$.
Hence, due to Lemma~\ref{resolvent integral kernel},
this family extends to a holomorphic family
$$\frac{\parameter{t}}{|\parameter{t}|}\begin{pmatrix}
\Op{I}_{\frac{2p}{2-p},\parameter{t}} & 0\\
0 & \Op{I}_{\frac{2p}{2-p},\Bar{\parameter{t}}}
\end{pmatrix}^{-1}\comp\Op{R}(0,0,k,0)\comp
\begin{pmatrix}
\Op{I}_{p,\Bar{\parameter{t}}} & 0\\
0 & \Op{I}_{p,\parameter{t}}
\end{pmatrix}$$
depending on $\{\parameter{t}\in\mathbb{P}^1\mid |\parameter{t}|\geq 1\}$.
More precisely, if
\begin{align*}
\Func{K}_1(z,k)&=z^{-1}+\sum\limits_{n=0}^{\infty}
\Func{K}_{1,n}z^n\text{ and}&
\Func{K}_2(z,k)&=\Bar{z}^{-1}+\sum\limits_{n=0}^{\infty}
\Func{K}_{2,n}\Bar{z}^n
\end{align*}
denotes the Laurent series of the two functions $\Func{K}_1$ and
$\Func{K}_2$ in the integral kernel of $\Op{R}(0,0,k,0)$
(compare with Lemma~\ref{resolvent integral kernel}),
then the Taylor series at $\parameter{t}=\infty$
of the integral kernel of this family is given by
$$\begin{pmatrix}
0 & (z-z')^{-1}\\
(\Bar{z}'-\Bar{z})^{-1} & 0\\
\end{pmatrix}\frac{d^2x'}{\pi}+
\sum\limits_{n\in\mathbb{N}}\parameter{t}^{-n}\begin{pmatrix}
0 & \Func{K}_{1,n-1}(z-z')^{n-1}\\
-\Func{K}_{2,n-1}(\Bar{z}-\Bar{z}')^{n-1}
\end{pmatrix}\frac{d^2x'}{\pi}.$$
Since the convolution with a polynomial with respect to $z$ is a
finite rank operator, with the exception of the zero--th
Taylor coefficients of this family at $\parameter{t}=\infty$
(which is equal to the restriction
of $\Op{R}_{\mathbb{R}^2}(0,0,0)$ to $B(0,\varepsilon)$)
all Taylor coefficients at $\parameter{t}=\infty$ are
finite rank operators from
$\banach{p}(B(0,\varepsilon))\times\banach{p}(B(0,\varepsilon))$
to $\banach{\frac{2p}{2-p}}(B(0,\varepsilon)\times
\banach{\frac{2p}{2-p}}(B(0,\varepsilon)$.

Therefore, we have also a real analytic family of operators
$$\begin{pmatrix}
\Op{I}_{p,\parameter{t}} & 0\\
0 & \Op{I}_{p,\parameter{t}}
\end{pmatrix}^{-1}\comp\left(\unity-\begin{pmatrix}
U_{\parameter{t}} & 0\\
0 & \Bar{U}_{\parameter{t}}
\end{pmatrix}\comp\Op{R}(0,0,k,0)\right)\comp\begin{pmatrix}
\Op{I}_{p,\parameter{t}} & 0\\
0 & \Op{I}_{p,\parameter{t}}
\end{pmatrix}$$
on $\banach{p}(B(0,\varepsilon))\times\banach{p}(B(0,\varepsilon))$
depending on $\parameter{t}\in [1,\infty)$,
which extends to a holomorphic family
\begin{multline*}
\unity-\frac{\parameter{t}}{|\parameter{t}|}\begin{pmatrix}
U_{\parameter{t}=1} & 0\\
0 & \Bar{U}_{\parameter{t}=1}
\end{pmatrix}\comp
\begin{pmatrix}
\Op{I}_{\frac{2p}{2-p},\parameter{t}} & 0\\
0 & \Op{I}_{\frac{2p}{2-p},\Bar{\parameter{t}}}
\end{pmatrix}^{-1}\comp\Op{R}(0,0,k,0)\comp
\begin{pmatrix}
\Op{I}_{p,\Bar{\parameter{t}}} & 0\\
0 & \Op{I}_{p,\parameter{t}}
\end{pmatrix}=\\
\begin{pmatrix}
\Op{I}_{p,\Bar{\parameter{t}}} & 0\\
0 & \Op{I}_{p,\parameter{t}}
\end{pmatrix}^{-1}\comp\left(\unity-\begin{pmatrix}
U_{\parameter{t}} & 0\\
0 & \Bar{U}_{\Bar{\parameter{t}}}
\end{pmatrix}\comp\Op{R}(0,0,k,0)\right)\comp
\begin{pmatrix}
\Op{I}_{p,\Bar{\parameter{t}}} & 0\\
0 & \Op{I}_{p,\parameter{t}}
\end{pmatrix}
\end{multline*}
depending on $\{\parameter{t}\in\mathbb{P}^1\mid |\parameter{t}|\geq 1\}$,
if we extend the family $U_{\parameter{t}}$ by the formula
\begin{align*}
U_{\parameter{t}}&=\Op{I}_{p,\Bar{\parameter{t}}}\comp
\frac{\parameter{t}}{|\parameter{t}|}U_{\parameter{t}=1}\comp
\Op{I}_{\frac{2p}{2-p},\parameter{t}}^{-1}&
\Bar{U}_{\Bar{\parameter{t}}}&=
\Op{I}_{p,\parameter{t}}\comp
\frac{\parameter{t}}{|\parameter{t}|}\Bar{U}_{\parameter{t}=1}\comp
\Op{I}_{\frac{2p}{2-p},\Bar{\parameter{t}}}^{-1}.
\end{align*}
This definition does not depend on $p\in(1,2)$.
If $\parameter{t}$ is not real, then these operators are not
multiplication operators and therefore non--local, but for real
$\parameter{t}$ this definition coincides with the definition of
$U_{\parameter{t}}$ above. At $\parameter{t}=\infty$ the value of this
family is equal to the restriction of the operator
$\unity-\left(\begin{smallmatrix}
U_{\parameter{t}=1} & 0\\
0 & \Bar{U}_{\parameter{t}=1}
\end{smallmatrix}\right)\comp\Op{R}_{\mathbb{R}^2}(0,0,0)$ to
$\banach{p}(B(0,\varepsilon))\times\banach{p}(B(0,\varepsilon))$.

\noindent
{\bf 2.\ A trivial kernel of the blown up Dirac operator.}
Due to Corollary~\ref{weakly continuous sphere resolvents}
the foregoing operator has an inverse, if and only if the Dirac operator
$(1+z\Bar{z})\left(\begin{smallmatrix}
U_{\parameter{t}=1} & \partial\\
-\Bar{\partial} & \Bar{U}_{\parameter{t}=1}
\end{smallmatrix}\right)$ on $\mathbb{P}^1$
does not have a kernel. At a first sight this is only true for
$p=\frac{4}{3}$ and $q=4$. But since the function $(1+z\Bar{z})$ and
its inverse are bounded on the domain of $U_{\parameter{t}=1}$, this is
true for all $1<p<2$. We conclude that in this case and for all
sufficiently large real $\parameter{t}$ and $k\notin\fermi(0,0)$
the restrictions of the operators
$\unity-\left(\begin{smallmatrix}
U_{\parameter{t}} & 0\\
0 & \Bar{U}_{\parameter{t}}
\end{smallmatrix}\right)\comp\Op{R}(0,0,k,0)$
to $\banach{p}(B(0,\parameter{t}^{-1}\varepsilon))\times
\banach{p}(B(0,\parameter{t}^{-1}\varepsilon))$
have an inverse, which is bounded uniformly with respect to
$\parameter{t}$.
Now we apply Lemma~\ref{inverse operator} and conclude
that under the same assumptions the resolvents
$\Op{R}(U_{\parameter{t}},\Bar{U}_{\parameter{t}},k,0)$
converge in the limit $\parameter{t}\rightarrow\infty$
to the corresponding limit of the operators
$$\Op{R}(0,0,k,0)+\Op{R}(0,0,k,0)\comp
\left(\unity-\begin{pmatrix}
U_{\parameter{t}} & 0\\
0 & \Bar{U}_{\parameter{t}}
\end{pmatrix}\comp\Op{R}_{\mathbb{R}^2}(0,0,0)\right)^{-1}
\comp\begin{pmatrix}
U_{\parameter{t}} & 0\\
0 & \Bar{U}_{\parameter{t}}
\end{pmatrix}\comp\Op{R}(0,0,k,0).$$
If in this formula we view the resolvent on the left hand side as
a bounded operator from $\banach{p'}(\Delta)\times\banach{p'}(\Delta)$
into $\banach{q'}(\Delta)\times\banach{q'}(\Delta)$ with $1<p'<2$ and
$p\leq q'\leq 2p'/(2-p')$ and $k\notin\fermi(0,0)$
(compare with Lemma~\ref{free resolvent}),
then the inverse operator may be considered as an operator on
$\banach{p'}(B(0,\parameter{t}^{-1}\varepsilon))\times
\banach{p'}(B(0,\parameter{t}^{-1}\varepsilon))$
instead of
$\banach{p}(B(0,\parameter{t}^{-1}\varepsilon))\times
\banach{p}(B(0,\parameter{t}^{-1}\varepsilon))$.
If $p'<p$, then the restriction is an operator from
$\banach{p}(\Delta)$ into
$\banach{p'}(B(0,\parameter{t}^{-1}\varepsilon))$,
which due to H\"older's inequality \cite[Theorem~III.1~(c)]{RS1}
is bounded by
$(\pi\varepsilon^2)^{(p-p')/(pp')}|\parameter{t}|^{2(p'-p)/(pp')}$.
Hence, if the Dirac operator
$\left(\begin{smallmatrix}
U_{\parameter{t}=1} & \partial\\
-\Bar{\partial} & \Bar{U}_{\parameter{t}=1}
\end{smallmatrix}\right)$
on $\mathbb{P}^1$ has no kernel,
then in the limit $t\rightarrow\infty$ the resolvents
$\Op{R}(U_{\parameter{t}},\Bar{U}_{\parameter{t}},k,0)$ converge as
operators from $\banach{p}(\Delta)\times\banach{p}(\Delta)$ into
$\banach{q'}(\Delta)\times\banach{q'}(\Delta)$
with $1<p\leq 2$ and $2\geq q'<q=\frac{2p}{2-p}$
for all $k\notin\fermi(0,0)$ in the limit $t\rightarrow\infty$
to the operator $\Op{R}(0,0,k,0)$.

\noindent
{\bf 3. A non--trivial kernel of the blown up Dirac operator.}
In the third step we extend our arguments to the case,
where the Dirac operator
$(1+z\Bar{z})\left(\begin{smallmatrix}
U_{\parameter{t}=1} & \partial\\
-\Bar{\partial} & \Bar{U}_{\parameter{t}=1}
\end{smallmatrix}\right)$
on $\mathbb{P}^1$ has a non--trivial kernel.
First we shall calculate the singular Laurent coefficients at
$\parameter{t}=\infty$ of the family of inverse operators
$$\begin{pmatrix}
\Op{I}_{p,\Bar{\parameter{t}}} & 0\\
0 & \Op{I}_{p,\parameter{t}}
\end{pmatrix}^{-1}\comp\left(\unity-\begin{pmatrix}
U_{\parameter{t}} & 0\\
0 & \Bar{U}_{\Bar{\parameter{t}}}
\end{pmatrix}\comp\Op{R}(0,0,k,0)\right)^{-1}\comp\begin{pmatrix}
\Op{I}_{p,\Bar{\parameter{t}}} & 0\\
0 & \Op{I}_{p,\parameter{t}}
\end{pmatrix}$$
on $\banach{p}(B(0,\varepsilon))\times\banach{p}(B(0,\varepsilon))$
of the holomorphic family introduced above.
The original holomorphic family has a Taylor series at
$\parameter{t}=\infty$, whose zero coefficient is equal to the
restriction $\Op{A}$ of the operator
$\unity-\left(\begin{smallmatrix}
U_{\parameter{t}=1} & 0\\
0 & \Bar{U}_{\parameter{t}=1}
\end{smallmatrix}\right)\comp\Op{R}_{\mathbb{R}^2}(0,0,0)$
to $\banach{p}(B(0,\varepsilon))\times\banach{p}(B(0,\varepsilon))$,
which due to Corollary~\ref{weakly continuous sphere resolvents}
is a Fredholm operator of index zero.
Furthermore, all other Taylor coefficients are
finite rank operators, since the corresponding
Taylor coefficients of the resolvents are finite rank operators.
Finally, the difference $\Op{B}_{\parameter{t}}$ of the original
holomorphic family minus $\Op{A}$
(the value at $\parameter{t}=\infty$) are compact operators.

\begin{Lemma}\label{perturbation of inverse operators}
Let $\Op{A}$ be a Fredholm operator of index zero from the
Banach space $\Spa{E}$ into the Banach space $\Spa{G}$.
Let $\Spa{E}_{\text{\scriptsize\rm finite}}\subset\Spa{E}$
and $\Spa{F}_{\text{\scriptsize\rm finite}}\subset\Spa{G}^{\ast}$
denote the kernel and co--kernel of $\Op{A}$,
respectively. Furthermore, let
$\Spa{G}_{\text{\scriptsize\rm infinite}}$ denote the range of
$\Op{A}$ (i.\ e.\ the orthogonal complement of
$\Spa{F}_{\text{\scriptsize\rm finite}}$). If
$\Op{B}$ is another operator from $\Spa{E}$ into $\Spa{G}$
with the property, that the operator
$$\Op{B}_{\text{\scriptsize\rm finite}}:
\Spa{E}_{\text{\scriptsize\rm finite}}
\hookrightarrow\Spa{E}\xrightarrow{\Op{B}}\Spa{G}
\hookrightarrow\Spa{G}^{\ast\ast}\twoheadrightarrow
\Spa{F}_{\text{\scriptsize\rm finite}}^{\ast}$$ is invertible,
then there exist unique decompositions
\begin{align*}
\Spa{E}&=\Spa{E}_{\text{\scriptsize\rm finite}}\oplus
\Spa{E}_{\text{\scriptsize\rm infinite}}&
\Spa{G}&=\Spa{G}_{\text{\scriptsize\rm finite}}\oplus
\Spa{G}_{\text{\scriptsize\rm infinite}},
\end{align*}
where $\Spa{E}_{\text{\scriptsize\rm infinite}}$ is the orthogonal
complement of the image of
$\Spa{F}_{\text{\scriptsize\rm finite}}$ under the operator
$\Op{B}^{\ast}:\Spa{G}^{\ast}\rightarrow\Spa{E}^{\ast}$, and
$\Spa{G}_{\text{\scriptsize\rm finite}}$ is the image of
$\Spa{E}_{\text{\scriptsize\rm finite}}$
under the operator $\Op{B}$. Moreover, both operators
$\Op{A}$ and $\Op{B}$ have diagonal block form
with respect to these decompositions.
\end{Lemma}

Here the map $\Spa{G}^{\ast\ast}\twoheadrightarrow
\Spa{F}_{\text{\scriptsize\rm finite}}^{\ast}$ on the
right hand side in the definition of
$\Op{B}_{\text{\scriptsize\rm finite}}$ is the
dual map of the natural inclusion
$\Spa{F}_{\text{\scriptsize\rm finite}}\hookrightarrow
\Spa{G}^{\ast}$.

\begin{proof}
If $\Op{B}_{\text{\scriptsize\rm finite}}$ is invertible, the operator
$$\Spa{E}\xrightarrow{\Op{B}}\Spa{G}
\hookrightarrow\Spa{G}^{\ast\ast}\twoheadrightarrow
\Spa{F}_{\text{\scriptsize\rm finite}}^{\ast}
\xrightarrow{\Op{B}_{\text{\tiny\rm finite}}^{-1}}
\Spa{E}_{\text{\scriptsize\rm finite}}
\hookrightarrow\Spa{E}$$
is a projection onto $\Spa{E}_{\text{\scriptsize\rm finite}}$.
Since $\Spa{E}_{\text{\scriptsize\rm infinite}}$ is the kernel of
the map
$$\Spa{E}\xrightarrow{\Op{B}}\Spa{G}
\hookrightarrow\Spa{G}^{\ast\ast}\twoheadrightarrow
\Spa{F}_{\text{\scriptsize\rm finite}}^{\ast},$$
this space is also the kernel of the former projection.
Therefore, we have a decomposition $\Spa{E}=
\Spa{E}_{\text{\scriptsize\rm finite}}\oplus
\Spa{E}_{\text{\scriptsize\rm infinite}}$.
On the other hand, the operator
$$\Spa{G}\hookrightarrow\Spa{G}^{\ast\ast}\twoheadrightarrow
\Spa{F}_{\text{\scriptsize\rm finite}}^{\ast}
\xrightarrow{\Op{B}_{\text{\tiny\rm finite}}^{-1}}
\Spa{E}_{\text{\scriptsize\rm finite}}
\xrightarrow{\Op{B}
\mid_{\Spa{E}_{\text{\tiny\rm finite}}}}
\Spa{G}_{\text{\scriptsize\rm finite}}
\hookrightarrow\Spa{G}$$
is a projection onto $\Spa{G}_{\text{\scriptsize\rm finite}}$.
Since $\Spa{G}_{\text{\scriptsize\rm infinite}}$ is the orthogonal
complement of $\Spa{F}_{\text{\scriptsize\rm finite}}$, this
subspace is the kernel of the latter projection, and we have a
decomposition $\Spa{G}=
\Spa{G}_{\text{\scriptsize\rm finite}}\oplus
\Spa{G}_{\text{\scriptsize\rm infinite}}$.

Moreover, the operator $\Op{A}$ vanishes on
$\Spa{E}_{\text{\scriptsize\rm finite}}$ and maps
$\Spa{E}_{\text{\scriptsize\rm infinite}}$ onto
$\Spa{G}_{\text{\scriptsize\rm infinite}}$.
By definition of $\Spa{G}_{\text{\scriptsize\rm finite}}$, the
operator $\Op{B}$ maps
$\Spa{E}_{\text{\scriptsize\rm finite}}$ onto
$\Spa{G}_{\text{\scriptsize\rm finite}}$.
Finally, since $\Spa{E}_{\text{\scriptsize\rm infinite}}$
is the kernel of the map
$$\Spa{E}\xrightarrow{\Op{B}}\Spa{G}
\hookrightarrow\Spa{G}^{\ast\ast}\twoheadrightarrow
\Spa{F}_{\text{\scriptsize\rm finite}}^{\ast},$$
the operator $\Op{B}$ maps
$\Spa{E}_{\text{\scriptsize\rm infinite}}$ into
$\Spa{G}_{\text{\scriptsize\rm infinite}}$.
\end{proof}

This lemma implies that the operator $\Op{A}+\Op{B}$ is
invertible, if $\Op{B}_{\text{\scriptsize\rm finite}}$
is invertible, and the component
of $\Op{B}$ from $\Spa{E}_{\text{\scriptsize\rm infinite}}$
into $\Spa{G}_{\text{\scriptsize\rm infinite}}$ has norm smaller
than the inverse of the same component of $\Op{A}$. In fact, in
this case the direct sum of these two components are invertible.
In the application we have in mind, $\Op{A}$ is equal to the value at
$\parameter{t}=\infty$ of the family of Fredholm operators introduced
above, and $\Op{B}_{\parameter{t}}$ is equal to the difference
of this family minus the value at $\parameter{t}=\infty$.
In the limit $\parameter{t}\rightarrow\infty$ the operators
$\Op{B}_{\parameter{t}}$ become arbitrarily small.

\begin{Remark}\label{change of notation}
  In the sequel we use a notation,
  which is slightly different from the notation used
  in Lemma~\ref{perturbation of inverse operators}. In fact, the spaces,
  which are denoted by capital $\Spa{F}$ in
  Lemma~\ref{perturbation of inverse operators}, are in the
  sequel denoted by $\Op{R}_{\mathbb{R}^2}^{t}(0,0,0)\Spa{F}$. The
  latter spaces $\Spa{F}$ have the advantage, that the support of
  their elements is contained in the support of $U_{\parameter{t}=1}$.
\end{Remark}

Let $\Spa{E}_{\text{\scriptsize\rm finite}}\subset\Spa{E}=
\banach{p}(B(0,\varepsilon))\times\banach{p}(B(0,\varepsilon))$
denote the kernel of the restriction of the operator
$\unity-\left(\begin{smallmatrix}
U_{\parameter{t}=1} & 0\\
0 & \Bar{U}_{\parameter{t}=1}
\end{smallmatrix}\right)\comp\Op{R}_{\mathbb{R}^2}(0,0,0)$
to $\banach{p}(B(0,\varepsilon))\times\banach{p}(B(0,\varepsilon))$,
and $\Spa{F}_{\text{\scriptsize\rm finite}}\subset\Spa{F}=
\banach{\frac{2p}{3p-2}}(B(0,\varepsilon))\times
\banach{\frac{2p}{3p-2}}(B(0,\varepsilon))$
the co--kernel of the restriction of the operator
$\unity-\Op{R}_{\mathbb{R}^2}(0,0,0)\comp\left(\begin{smallmatrix}
U_{\parameter{t}=1} & 0\\
0 & \Bar{U}_{\parameter{t}=1}
\end{smallmatrix}\right)$ to
$\banach{\frac{2p}{2-p}}(B(0,\varepsilon))\times
\banach{\frac{2p}{2-p}}(B(0,\varepsilon))$
(i.\ e.\ the kernel of the restriction of the operator
$\unity-\left(\begin{smallmatrix}
U_{\parameter{t}=1} & 0\\
0 & \Bar{U}_{\parameter{t}=1}
\end{smallmatrix}\right)\comp\Op{R}_{\mathbb{R}^2}^{t}(0,0,0)$ to
$\banach{\frac{2p}{3p-2}}(B(0,\varepsilon))\times
\banach{\frac{2p}{3p-2}}(B(0,\varepsilon))$).
Obviously, the supports of the elements of
$\Spa{E}_{\text{\scriptsize\rm finite}}$
and $\Spa{F}_{\text{\scriptsize\rm finite}}$ are contained in
$B(0,\varepsilon)$. Therefore, these
subspaces does not depend on $\varepsilon$, as long as $B(0,\varepsilon)$
contains the support of $U_{\parameter{t}=1}$. Furthermore, the resolvent
$\Op{R}_{\mathbb{R}^2}(0,0,0)$ and the transposed resolvent map
these subspaces $\Spa{E}_{\text{\scriptsize\rm finite}}$ and
$\Spa{F}_{\text{\scriptsize\rm finite}}$ onto the kernel and
the co--kernel of the corresponding Dirac operator on $\mathbb{P}^1$,
respectively. Finally, due to the statement analogous to
(iii) of Corollary~\ref{involutions},
$\Spa{E}_{\text{\scriptsize\rm finite}}$
and $\Spa{F}_{\text{\scriptsize\rm finite}}$ are invariant under
the anti--linear (projective) involutions
$\psi\mapsto\Op{J}\Bar{\psi}$
and $\phi\mapsto\Op{J}\Bar{\phi}$, respectively.
As we have seen in Corollary~\ref{fixed points},
these involutions cannot have fixed points, and the dimensions of
$\Spa{E}_{\text{\scriptsize\rm finite}}$
and $\Spa{F}_{\text{\scriptsize\rm finite}}$ are even, say $2K$.
Now let $n_1,\ldots,n_J$ denote the unique sequence of numbers in
$\mathbb{N}_0$ with the property,\vspace{.2cm}

\noindent\hfill
\begin{minipage}[b]{16cm}
  that for all $j=1,\ldots,J$ $n_j$ is the smallest index
  such that the span of
  $\begin{pmatrix}
  \Bar{z}^{n_1}\\
  0
  \end{pmatrix},\begin{pmatrix}
  0\\
  z^{n_1}
  \end{pmatrix},\ldots,\begin{pmatrix}
  \Bar{z}^{n_j}\\
  0
  \end{pmatrix},\begin{pmatrix}
  0\\
  z^{n_j}
  \end{pmatrix}$
  does not contain any non--trivial element in the orthogonal complement
  of $\Spa{E}_{\text{\scriptsize\rm finite}}$.
\end{minipage}\vspace{.2cm}

\noindent
Furthermore, let $m_1,\ldots,m_J$ be the analogous sequence with the
property,\vspace{.2cm}

\noindent\hfill
\begin{minipage}[b]{16cm}
  that for all $j=1,\ldots,J$ $m_j$ is the smallest index
  such that the span of
  $\begin{pmatrix}
  z^{m_1}\\
  0
  \end{pmatrix},\begin{pmatrix}
  0\\
  \Bar{z}^{m_1}
  \end{pmatrix},\ldots,\begin{pmatrix}
  z^{m_j}\\
  0
  \end{pmatrix},\begin{pmatrix}
  0\\
  \Bar{z}^{m_j}
  \end{pmatrix}$
  does not contain any non--trivial element in the orthogonal complement
  of $\Spa{F}_{\text{\scriptsize\rm finite}}$.
\end{minipage}\vspace{.2cm}

\noindent
In order to ensure the existence of these sequences we remark
that for all $\psi\in\Spa{E}_{\text{\scriptsize\rm finite}}$
and $\phi\in\Spa{F}_{\text{\scriptsize\rm finite}}$
the Taylor coefficients at $z=\infty$ of the corresponding
elements $\Op{R}_{\mathbb{R}^2}(0,0,0)\psi$ and
$\Op{R}_{\mathbb{R}^2}^{t}(0,0,0)\phi$ in the kernel and the
co--kernel of the corresponding Dirac operator on $\mathbb{P}^1$
are multiples of the numbers
$\left\langle\left\langle\left(\begin{smallmatrix}
\Bar{z}^n\\
0
\end{smallmatrix}\right),\psi\right\rangle\right\rangle$,
$\left\langle\left\langle\left(\begin{smallmatrix}
0\\
z^n
\end{smallmatrix}\right),\psi\right\rangle\right\rangle$ and 
$\left\langle\left\langle\phi,\left(\begin{smallmatrix}
z^n\\
0
\end{smallmatrix}\right)\right\rangle\right\rangle$,
$\left\langle\left\langle\phi,\left(\begin{smallmatrix}
0\\
\Bar{z}^n
\end{smallmatrix}\right)\right\rangle\right\rangle$ with
$n\in\mathbb{N}_0$, respectively. Hence, due to the
\De{Strong unique continuation property}~\ref{strong unique continuation},
not all of these numbers can vanish, and the sequences $n_1,\ldots,n_J$
and $m_1,\ldots,m_J$ are well defined
(compare with Remark~\ref{local contribution to the willmore functional}).

By definition of these sequences,
the kernel $\Spa{E}_{\text{\scriptsize\rm finite}}$
and the co--kernel $\Spa{F}_{\text{\scriptsize\rm finite}}$
have unique basis $\psi_1,\ldots,\psi_{2J}$ and
$\phi_1,\ldots\phi_{2J}$ dual to the basis
$$\begin{pmatrix}
\Bar{z}^{n_1}\\
0
\end{pmatrix},\begin{pmatrix}
0\\
z^{n_1}
\end{pmatrix},\ldots,\begin{pmatrix}
\Bar{z}^{n_J}\\
0
\end{pmatrix},\begin{pmatrix}
0\\
z^{n_J}
\end{pmatrix}\text{and}
\begin{pmatrix}
z^{m_1}\\
0
\end{pmatrix},\begin{pmatrix}
0\\
\Bar{z}^{m_1}
\end{pmatrix},\ldots,\begin{pmatrix}
z^{m_J}\\
0
\end{pmatrix},\begin{pmatrix}
0\\
\Bar{z}^{m_J}
\end{pmatrix},\text{respectively.}$$
Coming back to the notation of
Section~\ref{subsubsection finite rank perturbations}
let $D$ be the divisor equal to $\left(\max\{n_J,m_J\}+1\right)$
times the point $0$ on $\mathbb{C}/\lattice_\mathbb{C}$
and $\Psi_{D}$, $\Phi_{D}$ and $\Mat{R}_{D}(k)$
the corresponding matrices (compare with the first step of
Section~\ref{subsubsection finite rank perturbations}).
Moreover, let $\Mat{N}$ and $\Mat{M}$ denote the following
$\left(\max\{n_J,m_J\}+1\right)\times J$-- and
$J\times\left(\max\{n_J,m_J\}+1\right)$--matrices,
whose entries are the following $2\times 2$--matrices indexed by
$j=1,\ldots,J$ and $n,m=0,\ldots,\max\{n_J,m_J\}$:
\begin{align*}
\Mat{N}_{n,j}&=
\begin{cases}\unity & \text{if }n=n_j\\
             0 & \text{if }n\neq n_j
\end{cases}\text{ and} &
\Mat{M}_{j,m}&=
\begin{cases}\unity & \text{if }m=m_j\\
             0 & \text{if }m\neq m_j
\end{cases}\text{, respectively.}\end{align*}

The transposed resolvent $\Op{R}_{\mathbb{R}^2}^{t}(0,0,0)$ maps
the dual of the Banach space
$\banach{\frac{2p}{2-p}}(B(0,\varepsilon))\times
\banach{\frac{2p}{2-p}}(B(0,\varepsilon))$, which is isomorphic to
$\banach{\frac{2p}{3p-2}}(B(0,\varepsilon))\times
\banach{\frac{2p}{3p-2}}(B(0,\varepsilon))$,
onto the dual of the Banach space
$\banach{p}(B(0,\varepsilon))\times\banach{p}(B(0,\varepsilon))$,
which is isomorphic to
$\banach{\frac{p}{p-1}}(B(0,\varepsilon))\times
\banach{\frac{p}{p-1}}(B(0,\varepsilon))$.
The range of the restriction of the operator
$\unity-\left(\begin{smallmatrix}
U_{\parameter{t}=1} & 0\\
0 & \Bar{U}_{\parameter{t}=1}
\end{smallmatrix}\right)\comp\Op{R}_{\mathbb{R}^2}(0,0,0)$
to $\banach{p}(B(0,\varepsilon))\times\banach{p}(B(0,\varepsilon))$
is the orthogonal complement of $\Op{R}_{\mathbb{R}^2}^{t}(0,0,0)
\Spa{F}_{\text{\scriptsize\rm finite}}$.
Let $\Spa{G}_{\text{\scriptsize\rm infinite}}\subset
\banach{\frac{2p}{2-p}}(B(0,\varepsilon))\times
\banach{\frac{2p}{2-p}}(B(0,\varepsilon))$
($=\Spa{G}$) denote this orthogonal complement.
By construction the operator $\unity-\left(\begin{smallmatrix}
U_{\parameter{t}=1} & 0\\
0 & \Bar{U}_{\parameter{t}=1}
\end{smallmatrix}\right)\comp\Op{R}_{\mathbb{R}^2}(0,0,0)$
has the kernel $\Spa{E}_{\text{\scriptsize\rm finite}}$
and the co--kernel $\Op{R}_{\mathbb{R}^2}^{t}(0,0,0)
\Spa{F}_{\text{\scriptsize\rm finite}}$.
With respect to the basis $\psi_1,\ldots,\psi_{2J}$ and the dual basis of
$\Op{R}_{\mathbb{R}^2}^{t}(0,0,0)\phi_1,\ldots,
\Op{R}_{\mathbb{R}^2}^{t}(0,0,0)\phi_{2J}$ the holomorphic family
$$\Op{B}_{\parameter{t},\text{\scriptsize\rm finite}}:
\Spa{E}_{\text{\scriptsize\rm finite}}\hookrightarrow\Spa{E}
\xrightarrow{\Op{B}_{\parameter{t}}}\Spa{G}\hookrightarrow
\Spa{G}^{\ast\ast}\twoheadrightarrow
\left(\Op{R}_{\mathbb{R}^2}^{t}(0,0,0)
\Spa{F}_{\text{\scriptsize\rm finite}}\right)^{\ast}$$
takes the matrix form:
$$\frac{-1}{\parameter{t}}\begin{pmatrix}
\parameter{t}^{-m_1-n_1}\Mat{R}_{m_1,n_1} &\hdots &
\parameter{t}^{-m_1-n_J}\Mat{R}_{m_1,n_J}\\
\vdots &\hdots &\vdots\\
\parameter{t}^{-m_J-n_1}\Mat{R}_{m_J,n_1} &\hdots &
\parameter{t}^{-m_J-n_J}\Mat{R}_{m_J,n_J}.
\end{pmatrix}$$
Here $\Mat{R}_{m,n}$ denotes the enrties of the matrix $\Mat{R}_{D}(k)$.
In the sequel we will assume that this matrix is invertible for
$\parameter{t}\neq\infty$. This is equivalent to the
statement that the determinant of the $J\times J$--matrix
$\Mat{M}\comp\Mat{R}_{D}(k)\comp\Mat{N}$ does not vanish.
In the eighth step of Section~\ref{subsubsection finite rank perturbations}
we showed that the result of Section~\ref{subsection fermi curve}
carry over to these \Em{complex Fermi curves}.

If $n_J=J-1$ ($\iff n_j=j-1\;\forall j=1,\ldots,J$),
then the projection onto
$\Spa{E}_{\text{\scriptsize\rm finite}}$
$$\Spa{E}\xrightarrow{\Op{B}_{\parameter{t}}}\Spa{G}
\hookrightarrow\Spa{G}^{\ast\ast}\twoheadrightarrow
\left(\Op{R}_{\mathbb{R}^2}^{t}(0,0,0)
\Spa{F}_{\text{\scriptsize\rm finite}}\right)^{\ast}
\xrightarrow{\Op{B}_{\parameter{t},\text{\tiny\rm finite}}^{-1}}
\Spa{E}_{\text{\scriptsize\rm finite}}
\hookrightarrow\Spa{E}$$ 
has no singularity at $\parameter{t}=\infty$.
Unfortunately for all other values of $n_J$ it has a singularity.
In fact, if $\psi\in\Spa{E}$ has the property that
\begin{align*}
\left\langle\left\langle
\left(\begin{smallmatrix}
\Bar{z}^n\\
0
\end{smallmatrix}\right),\psi\right\rangle\right\rangle&=
\begin{cases}0 &\text{for }n=0,\ldots,N-1\\
a &\text{ for }n=N
\end{cases}&
\left\langle\left\langle
\left(\begin{smallmatrix}
0\\
z^n
\end{smallmatrix}\right),\psi\right\rangle\right\rangle&=
\begin{cases}0 &\text{for }n=0,\ldots,N-1\\
b &\text{ for }n=N,
\end{cases}
\end{align*}
with  $N\in\mathbb{N}_0$ and $(a,b)\in\mathbb{P}^1$,
then the image of $\psi$ under
$$\Spa{E}\xrightarrow{\Op{B}_{\parameter{t}}}\Spa{G}
\hookrightarrow\Spa{G}^{\ast\ast}\twoheadrightarrow
\left(\Op{R}_{\mathbb{R}^2}^{t}(0,0,0)
\Spa{F}_{\text{\scriptsize\rm finite}}\right)^{\ast}$$
has with respect to the dual basis of
$\Op{R}_{\mathbb{R}^2}^{t}(0,0,0)\phi_1,\ldots,
\Op{R}_{\mathbb{R}^2}^{t}(0,0,0)\phi_{2J}$
and for large $\parameter{t}$ the form
$$\left(\frac{1}{\pi\parameter{t}}+\text{\bf{O}}(\parameter{t}^{-2})\right)
\begin{pmatrix}
\binom{m_1+N}{m_1}\frac{(-1)^{N}}{\parameter{t}^{m_1+N}}
\left(\begin{smallmatrix}
-b\Func{K}_{1,m_1+N}\\
a\Func{K}_{2,m_1+N}
\end{smallmatrix}\right)\\
\vdots\\
\binom{m_J+N}{m_J}\frac{(-1)^{N}}{\parameter{t}^{m_J+N}}
\left(\begin{smallmatrix}
-b\Func{K}_{1,m_J+N}\\
a\Func{K}_{2,m_J+N}
\end{smallmatrix}\right)\end{pmatrix}.$$
Obviously the inverse of
$\Op{B}_{\parameter{t},\text{\scriptsize\rm finite}}$
applied on this vector can have a singularity at $\parameter{t}=\infty$,
only if $N\in\{0,\ldots,n_J\}\setminus\{n_1,\ldots,n_J\}$.
Analogously, if $m_J=J-1$ ($\iff m_j=j-1\;\forall j=1,\ldots,J$),
then the projection of $\Spa{F}$ onto
$\Spa{F}_{\text{\scriptsize\rm finite}}$ has no singularity at
$\parameter{t}=\infty$.

We conclude that the resolvents of the family
$(U_{\parameter{t}},\Bar{U}_{\Bar{\parameter{t}}})$ differ from the free
resolvent by
\begin{multline*}
\Op{R}(U_{\parameter{t}},\Bar{U}_{\Bar{\parameter{t}}},k,0)-
\Op{R}(0,0,k,0)=\\
\begin{aligned}
=&\Op{R}(0,0,k,0)\comp\begin{pmatrix}
\Op{I}_{p,\Bar{\parameter{t}}} & 0\\
0 & \Op{I}_{p,\parameter{t}}
\end{pmatrix}\comp
\left(\Op{A}+\Op{B}_{\parameter{t}}\right)^{-1}\comp\begin{pmatrix}
\Op{I}_{p,\Bar{\parameter{t}}} & 0\\
0 & \Op{I}_{p,\parameter{t}}
\end{pmatrix}^{-1}\comp\begin{pmatrix}
U_{\parameter{t}} & 0\\
0 & \Bar{U}_{\Bar{\parameter{t}}}
\end{pmatrix}\comp\Op{R}(0,0,k,0)&&\\
=&\Op{R}(0,0,k,0)\comp\begin{pmatrix}
\Op{I}_{p,\Bar{\parameter{t}}} & 0\\
0 & \Op{I}_{p,\parameter{t}}
\end{pmatrix}\comp
\left(\Op{A}+\Op{B}_{\parameter{t}}\right)^{-1}\comp
\left(\begin{smallmatrix}
U_{\parameter{t}=1} & 0\\
0 & \Bar{U}_{\parameter{t}=1}
\end{smallmatrix}\right)\comp\begin{pmatrix}
\Op{I}_{\frac{2p}{2-p},\parameter{t}} & 0\\
0 & \Op{I}_{\frac{2p}{2-p},\Bar{\parameter{t}}}
\end{pmatrix}^{-1}\comp\Op{R}(0,0,k,0).&&
\end{aligned} \end{multline*}
We shall calculate these differences in the limit
$\parameter{t}\rightarrow\infty$ with the arguments used in the first
two steps combined with Lemma~\ref{perturbation of inverse operators}.
By definition of the sequence $n_1,\ldots,n_J$ and due to
Lemma~\ref{resolvent integral kernel} we have
\begin{eqnarray*}
\Op{R}(0,0,k,0)\comp\begin{pmatrix}
\Op{I}_{p,\parameter{t}} & 0\\
0 & \Op{I}_{p,\parameter{t}}
\end{pmatrix}\psi_{2j-1}&=&
\frac{|\parameter{t}|^{2/p}(-1)^{n_j}}{\pi n_j!\parameter{t}^{n_j+2}}
\left(\begin{pmatrix}
0\\
-\Func{K}_2^{(n_j)}
\end{pmatrix}+\text{\bf{O}}(\parameter{t}^{-1})\right)\text{ and}\\
\Op{R}(0,0,k,0)\comp\begin{pmatrix}
\Op{I}_{p,\parameter{t}} & 0\\
0 & \Op{I}_{p,\parameter{t}}
\end{pmatrix}\psi_{2j}&=&
\frac{|\parameter{t}|^{2/p}(-1)^{n_j}}{\pi n_j!\parameter{t}^{n_j+2}}
\left(\begin{pmatrix}
\Func{K}_1^{(n_j)}\\
0
\end{pmatrix}+\text{\bf{O}}(\parameter{t}^{-1})\right)
\end{eqnarray*}
on the complement of the support of $U_{\parameter{t}}$
and for all real $\parameter{t}$ and $j=1,\ldots,J$
Here $\Func{K}_1^{(n_j)}$ and $\Func{K}_2^{(n_j)}$ denotes
the holomorphic and anti--holomorphic $n_j$--th derivatives of the
functions $\Func{K}_1$ and $\Func{K}_2$, respectively.
Moreover, the transpose of
$\Op{I}_{\frac{2p}{2-p},\parameter{t}}^{-1}$ is
$\Op{I}_{\frac{2p}{3p-2},\parameter{t}}$. Therefore we also have
\begin{eqnarray*}
\Op{R}^{t}(0,0,k,0)\comp\begin{pmatrix}
\Op{I}_{\frac{2p}{3p-2},\parameter{t}} & 0\\
0 & \Op{I}_{\frac{2p}{3p-2},\parameter{t}}
\end{pmatrix}\phi_{2k-1}&=&
\frac{|\parameter{t}|^{(1-2/p)}}{\pi m_j!\parameter{t}^{m_j}}
\left(\begin{pmatrix}
0\\
-\Func{K}_1^{(m_j)}
\end{pmatrix}+\text{\bf{O}}(\parameter{t}^{-1})\right)\text{ and}\\
\Op{R}^{t}(0,0,k,0)\comp\begin{pmatrix}
\Op{I}_{\frac{2p}{3p-2},\parameter{t}} & 0\\
0 & \Op{I}_{\frac{2p}{3p-2},\parameter{t}}
\end{pmatrix}\phi_{2k}&=&
\frac{|\parameter{t}|^{(1-2/p)}}{\pi m_j!\parameter{t}^{m_j}}
\left(\begin{pmatrix}
\Func{K}_2^{(m_j)}\\
0
\end{pmatrix}+\text{\bf{O}}(\parameter{t}^{-1})\right),
\end{eqnarray*}
on the complement of the support of $U_{\parameter{t}}$
and for all real $\parameter{t}$ and for all $j=1,\ldots,J$.
If we endow $\Spa{E}/\Spa{E}_{\text{\scriptsize\rm finite}}$
with the topology of the dual Banach space
of the orthogonal complement of 
$\Spa{E}_{\text{\scriptsize\rm finite}}$ in the dual of the
Banach space
$\banach{p}(B(0,\varepsilon))\times\banach{p}(B(0,\varepsilon))$,
which is isomorphic to $\banach{\frac{p}{p-1}}(B(0,\varepsilon))\times
\banach{\frac{p}{p-1}}(B(0,\varepsilon))$,
then due to Corollary~\ref{weakly continuous sphere resolvents}
the operator $\Op{A}$ induces an invertible operator from
$\Spa{E}/\Spa{E}_{\text{\scriptsize\rm finite}}$ onto
$\Spa{G}_{\text{\scriptsize\rm infinite}}$.
The identity of $\Spa{E}$ minus the projection of $\Spa{E}$ onto
the subspace $\Spa{E}_{\text{\scriptsize\rm finite}}$, which was
introduced in Lemma~\ref{perturbation of inverse operators}, defines
an operator from
$\Spa{E}/\Spa{E}_{\text{\scriptsize\rm finite}}$ onto the
corresponding complementary subspace
$\Spa{E}_{\text{\scriptsize\rm infinite}}$.
We have seen above that this operator is not bounded in the limit
$\parameter{t}\rightarrow\infty$, if $n_J\neq J-1$. Nevertheless, this
divergence is overcompensated by the convergence of the corresponding
limits of $\Op{R}(0,0,k,0)\comp\left(\begin{smallmatrix}
\Op{I}_{p,\parameter{t}} & 0\\
0 & \Op{I}_{p,\parameter{t}}
\end{smallmatrix}\right)\psi_j$. Hence the composition of the
projection onto $\Spa{E}_{\text{\scriptsize\rm infinite}}$,
considered as an operator from
$\Spa{E}/\Spa{E}_{\text{\scriptsize\rm finite}}$ into $\Spa{E}$,
with the operator $\Op{R}(0,0,k,0)\comp\left(\begin{smallmatrix}
\Op{I}_{p,\parameter{t}} & 0\\
0 & \Op{I}_{p,\parameter{t}}
\end{smallmatrix}\right)$ converges to zero
in the limit $\parameter{t}\rightarrow\infty$.
On the other hand, the operator $\left(\begin{smallmatrix}
U_{\parameter{t}=1} & 0\\
0 & \Bar{U}_{\parameter{t}=1}
\end{smallmatrix}\right)$ maps
$\Spa{G}_{\text{\scriptsize\rm infinite}}$ onto the dual space of
$\Spa{F}/\Spa{F}_{\text{\scriptsize\rm finite}}$,
since $\left(\begin{smallmatrix}
U_{\parameter{t}=1} & 0\\
0 & \Bar{U}_{\parameter{t}=1}
\end{smallmatrix}\right)\comp\Op{R}_{\mathbb{R}^2}^{t}(0,0,0)$
acts on $\Spa{F}_{\text{\scriptsize\rm finite}}$ as the identity.
The analogous projection of $\Spa{F}$ onto the complementary
subspace $\Spa{F}_{\text{\scriptsize\rm infinite}}$ defines an
operator from
$\Spa{F}/\Spa{F}_{\text{\scriptsize\rm infinite}}$ onto
$\Spa{F}$, which is not bounded in the limit
$\parameter{t}\rightarrow\infty$, if $m_j\neq J-1$. Analogously this
divergence is overcompensated by the convergence of the corresponding
limits $\Op{R}^{t}(0,0,k,0)\comp\left(\begin{smallmatrix}
\Op{I}_{\frac{2p}{3p-2},\parameter{t}} & 0\\
0 & \Op{I}_{\frac{2p}{3p-2},\parameter{t}}
\end{smallmatrix}\right)\phi_j$. To sum up, in all cases
regardless of whether $n_j$ or $m_j$ is equal to $J-1$ or not,
the arguments of the first two steps combined with the application of
Lemma~\ref{perturbation of inverse operators} yields
that only the part
$$\Spa{G}\hookrightarrow\Spa{G}^{\ast\ast}\twoheadrightarrow
\left(\Op{R}_{\mathbb{R}^2}^{t}(0,0,0)\Spa{F}\right)^{\ast}
\xrightarrow{\Op{B}_{\parameter{t},\text{\tiny\rm finite}}^{-1}}
\Spa{E}_{\text{\scriptsize\rm finite}}\hookrightarrow\Spa{E}$$
of $(\Op{A}+\Op{B}_{\parameter{t}})^{-1}$ contributes
in the limit $\parameter{t}\rightarrow\infty$ to
$\Op{R}(U_{\parameter{t}},\Bar{U}_{\parameter{t}},k,0)$. Therefore,
and with the notation of
Section~\ref{subsubsection finite rank perturbations},
the operator $\lim\limits_{\parameter{t}\rightarrow\infty}
\Op{R}(U_{\parameter{t}},\Bar{U}_{\parameter{t}},k,0)-\Op{R}(0,0,k,0)$
has the integral kernel
\begin{align*}
\Psi_{D}(z)\Mat{S}_{D}(k)\Phi_{D}(z')d^2x'&\text{ with}&
\Mat{S}_{D}(k)=-\Mat{N}\comp
\left(\Mat{M}\comp\Mat{R}_{D}(k)\comp\Mat{N}\right)^{-1}
\comp\Mat{M}.\end{align*}
Obviously this is an integral kernel of a
\De{Finite rank Perturbation}~\ref{finite rank perturbations}.

More precisely, we shall show that for all $\varepsilon>0$,
the restrictions (compare with Remark~\ref{restriction})
of the resolvents $\Op{R}(U_{\parameter{t}},\Bar{U}_{\parameter{t}},k,0)$
to the complements of $B(0,\varepsilon)$,
considered as operators from
$\banach{p}(\Delta\setminus B(0,\varepsilon))\times
\banach{p}(\Delta\setminus B(0,\varepsilon))$ to 
$\banach{\frac{2p}{2-p}}(\Delta\setminus B(0,\varepsilon))\times
\banach{\frac{2p}{2-p}}(\Delta\setminus B(0,\varepsilon))$
with $1<p< 2$,
converges in the limit $\parameter{t}\rightarrow\infty$
to the corresponding restriction of the operator above.
In fact, due to Lemma~\ref{resolvent integral kernel} the restriction
of the resolvent $\Op{R}(0,0,k,0)$ considered as an operator from
$\banach{p}(B(0,\parameter{t}^{-1}\varepsilon))\times
\banach{p}(B(0,\parameter{t}^{-1}\varepsilon))$ into
$\banach{\infty}(\Delta\setminus B(0,\varepsilon))\times
\banach{\infty}(\Delta\setminus B(0,\varepsilon))$
is for all $k\notin\fermi(0,0)$ and large $\parameter{t}$ bounded.
Therefore, the arguments of the first two steps carry over and
prove the claim, provided the Dirac operator $\left(\begin{smallmatrix}
U_{\parameter{t}=1} & \partial\\
-\Bar{\partial} & \Bar{U}_{\parameter{t}=1}
\end{smallmatrix}\right)$ on $\mathbb{P}^1$ has no kernel.
Consequently the arguments of the third step imply the claim for
general potentials $U_{\parameter{t}=1}\in \banach{2}(B(0,\varepsilon))$.
\end{proof}

\subsubsection{The order of zeroes of spinors}
\label{subsubsection local behaviour}

For any potential $U\in \banach{2}(\mathbb{R}^2)$
the restriction $U_{\parameter{t}}$ of
$z\mapsto \parameter{t}U(\parameter{t}z)$ to any small ball
$B(0,\varepsilon)$ yields an analogous family of potentials.
In order to carry over Lemma~\ref{nontrivial limit 1}
to these families, we shall investigate the limits
$\lim\limits_{\parameter{t}\rightarrow\infty}\left(\begin{smallmatrix}
\Op{I}_{p,\parameter{t}} & 0\\
0 & \Op{I}_{p,\parameter{t}}
\end{smallmatrix}\right)\psi$ of spinors $\psi$
in the kernel of a Dirac operator on $\mathbb{P}^1$
with $\banach{2}$--potentials.
Obviously these limits depend only on the restriction of $\psi$ to
arbitrary small neighbourhoods of $z=\infty$. Moreover, if we use the
coordinates $z^{-1}$ instead of $z$ this limit is the same as
$\lim\limits_{\parameter{t}\rightarrow 0}\left(\begin{smallmatrix}
\Op{I}_{p,\parameter{t}} & 0\\
0 & \Op{I}_{p,\parameter{t}}
\end{smallmatrix}\right)\psi$
of a spinor $\psi$ in the kernel of $\left(\begin{smallmatrix}
V & \partial\\
-\Bar{\partial} & W
\end {smallmatrix}\right)$ on some neighbourhood
$0\ni\Set{O}\subset\mathbb{C}$.

In the following discussion concerning these limits we make use
of the \Em{Lorentz spaces} \index{Lorentz spaces} $\banach{p,q}$.
These rearrangement invariant Banach spaces are an extension of the
family of the usual Banach spaces $\banach{p}$ indexed by an additional
parameter $1\leq q\leq\infty$ for $1<p<\infty$. For $p=1$ or
$p=\infty$ we consider only the \Em{Lorentz spaces} $\banach{p,\infty}$,
which in these cases are isomorphic to $\banach{p}$ 
(\cite[Chapter~V. \S3.]{SW},\cite[Chapter~4 Section~4.]{BS} 
and \cite[Chapter~1. Section~8.]{Zi}).
We recall some properties of these Banach spaces:
\begin{description}
\item[(i)] For $1<p\leq\infty$ the \Em{Lorentz spaces} $\banach{p,p}$
  coincide with the usual $\banach{p}$--spaces. Moreover, the
  \Em{Lorentz space} $\banach{1,\infty}$ coincides with the usual Banach
  space $\banach{1}$.
\item[(ii)] On a finite measure space the \Em{Lorentz space}
  $\banach{p,q}$ is contained in $\banach{p',q'}$
  either if $p>p'$ or if $p=p'$ and $q\leq q'$.
\end{description}
In \cite{O} H\"older's inequality and Young's inequality are generalized
to these \Em{Lorentz spaces}
(\cite[Chapter~4 Section~7.]{BS} and \cite[Chapter~2. Section~10.]{Zi}):

\newtheorem{Generalized Hoelder}[Lemma]{Generalized H\"older's inequality}
\index{generalized!H\"older's inequality}
\begin{Generalized Hoelder}\label{generalized hoelder}
Either for $1/p_1+1/p_2=1/p_3<1$ and $1/q_1+1/q_2\geq 1/q_3$ or for
$1/p_1+1/p_2=1$, $1/q_1+1/q_2\geq 1$ and $(p_3,q_3)=(1,\infty)$
there exists some constant $C>0$ with
$$\|fg\|_{(p_3,q_3)}\leq C\|f\|_{(p_1,q_1)}\|g\|_{(p_2,q_2)}.$$
\end{Generalized Hoelder}
\newtheorem{Generalized Young}[Lemma]{Generalized Young's inequality}
\index{generalized!Young's inequality}
\begin{Generalized Young}\label{generalized young}
Either for $1/p_1+1/p_2-1=1/p_3>0$ and $1/q_1+1/q_2\geq 1/q_3$ or for
$1/p_1+1/p_2=1$, $1/q_1+1/q_2\geq 1$ and $(p_3,q_3)=(\infty,\infty)$
there exists some constant $C>0$ with
$$\|f\ast g\|_{(p_3,q_3)}\leq C\|f\|_{(p_1,q_1)}\|g\|_{(p_2,q_2)}.$$
\end{Generalized Young}

Therefore, the resolvents of the Dirac operators
on $\torus$ or $\mathbb{P}^1$ or $\mathbb{R}^2$
are bounded  operators
from the $\banach{1}$--spinors into the $\banach{2,\infty}$--spinors,
from the $\banach{2,1}$--spinors into the continuous spinors,
from the $\banach{p}$--spinors into the $\banach{q,p}$--spinors,
and finally from the $\banach{p,q}$--spinors into the $\banach{q}$--spinors,
with $1<p<2$ and $q=2p/(2-p)$.
Moreover, the Sobolev constant $S_p$
(compare with Lemma~\ref{free resolvent} and Remark~\ref{sobolev 1})
is equal to the corresponding norm $\|f\|_{2,\infty}$
times the corresponding constant of the
\De{Generalized Young's inequality}~\ref{generalized young}.

Obviously the restriction of $z\mapsto\psi(\parameter{t}z)$ to any ball
$B(0,\varepsilon)$ depends for small $\parameter{t}$ only on the restriction
of $\psi$ to an arbitrary small neighbourhood of $0$. Let
$\Set{O}\subset\mathbb{C}$ be a bounded open neighbourhood of $0$
and let $V,W\in \banach{2}(\Set{O})$ be two potentials on this domain.
If $\psi$ is any element in the kernel of the corresponding Dirac operator
$\left(\begin{smallmatrix}
V & \partial\\
-\Bar{\partial} & W
\end{smallmatrix}\right)$ on this domain, then
$\left(\unity-\Op{R}_{\mathbb{R}^2}(0,0,0)\comp
\left(\begin{smallmatrix}
V & 0\\
0 & W
\end{smallmatrix}\right)\right)\psi$ defines on this domain an element
of the kernel of the free Dirac operator $\left(\begin{smallmatrix}
0 & \partial\\
-\Bar{\partial} & 0
\end{smallmatrix}\right)$. Vice versa, if we diminish the
neighbourhood $\Set{O}\ni 0$ and consequently the $\banach{2}$--norms of
$V$ and $W$, such that the sum
$$\sum\limits_{n=0}^{\infty}
\left(\Op{R}_{\mathbb{R}^2}(0,0,0)\comp\left(\begin{smallmatrix}
V & 0\\
0 & W
\end{smallmatrix}\right)\right)^n=\left(\unity-
\Op{R}_{\mathbb{R}^2}(0,0,0)\comp\left(\begin{smallmatrix}
V & 0\\
0 & W
\end{smallmatrix}\right)\right)^{-1}$$
converges to a bounded operator on
$\banach{q}(\Set{O})\times\banach{q}(\Set{O})$
with $2<q<\infty$,
then this operator transforms all elements in the kernel of the free
Dirac operator on $\Set{O}$ into the kernel of the Dirac operator
$\left(\begin{smallmatrix}
V & \partial\\
-\Bar{\partial} & W
\end{smallmatrix}\right)$ on this domain. Hence the latter kernel is
isomorphic to the former kernel. Moreover, the former kernel contains
all functions $\psi$, whose components are holomorphic and
anti--holomorphic functions on $\Set{O}$, respectively.
This shows that any spinor
$\psi\in\sobolev{1,p}(\Set{O})\times\sobolev{1,p}(\Set{O})$
with $1<p<2$ in the kernel of the Dirac operator
$\left(\begin{smallmatrix}
V & \partial\\
-\Bar{\partial} & W
\end{smallmatrix}\right)$ on $\Set{O}$ with potentials
$V,W\in \banach{2}(\Set{O})$ are on small neighbourhood equal to
$\left(\unity-\Op{R}_{\mathbb{R}^2}(0,0,0)\comp
\left(\begin{smallmatrix}
V & 0\\
0 & W
\end{smallmatrix}\right)\right)^{-1}$
applied on an element in the kernel of the free Dirac operator
on $\Set{O}$. Due to Weyl's Lemma
\cite[\S18.4. Lemma~4.10]{Co2} all elements in the kernel of
the free Dirac operator are continuous functions.
This implies that the spinor belongs to
$$\psi\in\bigcap\limits_{1<p<2}
\sobolev{1,p}_{\text{\scriptsize\rm loc}}(\Set{O})\times
\sobolev{1,p}_{\text{\scriptsize\rm loc}}(\Set{O})
\subset\bigcap\limits_{2<q<\infty}
\banach{q}_{\text{\scriptsize\rm loc}}(\Set{O})
\times\banach{q}_{\text{\scriptsize\rm loc}}(\Set{O}).$$
Moreover, for $1<p<2$ and $\frac{1}{p}=\frac{1}{2}+\frac{1}{q}$
the $\banach{q}(\Set{O}\times\banach{q}(\Set{O})$--norms
and the $\sobolev{1,p}(\Set{O})\times\sobolev{1,p}(\Set{O})$--norms
of a spinor in the kernel of a Dirac operator with potentials
$V,W\in\banach{2}(\Set{O})$ are equivalent.
Since the kernel of the free Dirac operator is a direct sum of
Bergman spaces,
i.\ e.\ the closed subspaces of $\banach{q}(B(0,\varepsilon))$
of holomorphic and anti--holomorphic functions, respectively
(\cite[Chapter~1.]{HKZ} and \cite[Chapter~4.]{Zh}),
all kernels of Dirac operators with potentials
$V,W\in\banach{2}(\Set{O})$ are closed subspaces of
$\banach{q}(\Set{O})\times\banach{q}(\Set{O})$ ($2<q<\infty$).
In the sequel these kernels are considered as subspaces of
these Banach spaces.

If a spinor
$\psi\in\sobolev{1,p}(\Set{O})\times\sobolev{1,p}(\Set{O})$
with $1<p<2$ belongs to the kernel of the Dirac operator
$\left(\begin{smallmatrix}
V & \partial\\
-\Bar{\partial} & W
\end{smallmatrix}\right)$ on $\Set{O}$ with potentials
$V,W\in \banach{2}(\Set{O})$, then obviously the spinor
$\left(\begin{smallmatrix}
z & 0\\
0 & \Bar{z}
\end{smallmatrix}\right)\psi$ belongs to the kernel
of the Dirac operator $\left(\begin{smallmatrix}
V\Bar{z}/z & \partial\\
-\Bar{\partial} & Wz/\Bar{z}
\end{smallmatrix}\right)$.

\newtheorem{Order of zeroes}[Lemma]{Order of zeroes}
\index{order of zeroes $\text{\rm ord}_{z}(\psi)$}
\begin{Order of zeroes}\label{order of zeroes}
The order of a zero of $\psi$ at $z=0$ is defined
as the largest integer $m$, such that
$\left(\begin{smallmatrix}
z & 0\\
0 & \Bar{z}
\end{smallmatrix}\right)^{-m}\psi\in
\sobolev{1,p}(\Set{O})\times\sobolev{1,p}(\Set{O})$
with $1<p<2$ belongs to the kernel of the
Dirac operator $\left(\begin{smallmatrix}
V(z/\Bar{z})^m & \partial\\
-\Bar{\partial} & W(\Bar{z}/z)^{m}
\end{smallmatrix}\right)$.
Due to the
\De{Strong unique continuation property}~\ref{strong unique continuation}
this number is finite and denoted by $\text{\rm ord}_{0}(\psi)$.
\end{Order of zeroes}

If the potentials $V$ and $W$ belong to $\banach{2,1}(\Set{O})$,
then the eigen--spinors of the corresponding Dirac operators are
continuous spinors. We claim that in this case $\psi$ vanishes at
$z=0$ if and only if $\Tilde{\psi}=\left(\begin{smallmatrix}
z & 0\\
0 & \Bar{z}
\end{smallmatrix}\right)^{-1}\psi$ belongs to the kernel of
the Dirac operator $\left(\begin{smallmatrix}
\Tilde{V} & \partial\\
-\Bar{\partial} & \Tilde{W}
\end{smallmatrix}\right)=\left(\begin{smallmatrix}
Vz/\Bar{z} & \partial\\
-\Bar{\partial} & W\Bar{z}/z
\end{smallmatrix}\right)$. In fact, since the resolvent
$\Op{R}_{\mathbb{R}^2}(0,0,0)$ has the integral kernel
$$\begin{pmatrix}
0 & (z-z')^{-1}\\
(\Bar{z}'-\Bar{z})^{-1} & 0
\end{pmatrix}\frac{d\Bar{z}'\wedge dz'}{2\pi\sqrt{-1}},$$
the operator
$\left(\begin{smallmatrix}
z & 0\\
0 & \Bar{z}
\end{smallmatrix}\right)\comp\Op{R}_{\mathbb{R}^2}(0,0,0)\comp
\left(\begin{smallmatrix}
\Bar{z} & 0\\
0 & z
\end{smallmatrix}\right)^{-1}-\Op{R}_{\mathbb{R}^2}(0,0,0)$
has the integral kernel
$\left(\begin{smallmatrix}
0 & 1/z'\\
-1/\Bar{z}' & 0
\end{smallmatrix}\right)\frac{d\Bar{z}'\wedge dz'}{2\pi\sqrt{-1}}$.
Consequently, the operator
$\unity-\Op{R}_{\mathbb{R}^2}(0,0,0)\comp\left(\begin{smallmatrix}
\Tilde{V} & 0\\
0 & \Tilde{W}
\end{smallmatrix}\right)$
maps $\Tilde{\psi}$ into the kernel of the
free Dirac operator on $\Set{O}$, if $\psi$ vanishes at $z=0$.
Vice versa, if $\Tilde{\psi}$ belongs to the kernel
of the Dirac operator $\left(\begin{smallmatrix}
\Tilde{V} & \partial\\
-\Bar{\partial} & \Tilde{W}
\end{smallmatrix}\right)$,
then $\Tilde{\psi}$ is continuous, since the potentials
$\Tilde{V}$ and $\Tilde{W}$ belong to $\banach{2,1}(\Set{O})$.
Therefore, the spinor $\psi=\left(\begin{smallmatrix}
z & 0\\
0 & \Bar{z}
\end{smallmatrix}\right)\Tilde{\psi}$ vanishes at $z=0$ and obviously
belongs to the kernel of the Dirac operator $\left(\begin{smallmatrix}
V & \partial\\
-\Bar{\partial} & W
\end{smallmatrix}\right)$ on $\Set{O}$.
This proves

\begin{Lemma}\label{order of zeroes 1}
\index{order of zeroes $\text{\rm ord}_{z}(\psi)$}
Let $V$ and $W$ be potentials in
$\banach{2,1}(\Set{O})$ on a bounded open neighbourhood
of $0\in\Set{O}\subset\mathbb{C}$ and $\psi$
a non--trivial spinor in the kernel of $\left(\begin{smallmatrix}
V & \partial\\
-\Bar{\partial} & W
\end{smallmatrix}\right)$ on $\Set{O}$.
Then there exists a unique $\mathbb{N}_0\ni m=\text{\rm ord}_{0}(\psi)$
such that
$\left(\begin{smallmatrix}
z & 0\\
0 & \Bar{z}
\end{smallmatrix}\right)^{-m}\psi$
is a continuous spinor in the kernel of $\left(\begin{smallmatrix}
V(z/\Bar{z})^m & \partial\\
-\Bar{\partial} & W(\Bar{z}/z)^m
\end{smallmatrix}\right)$, which does not vanish
in a small neighbourhood of $0\in\Set{O}$.
In particular, all zeroes are isolated.\qed
\end{Lemma}

If the potentials $V$ and $W$ belong to $\banach{2}(\Set{O})$,
then the space of all spinors, which have at $0\in\Set{O}$ a zero of
positive order, is a subspace of the kernel of the Dirac operator
$\left(\begin{smallmatrix}
V & \partial\\
-\Bar{\partial} & W
\end{smallmatrix}\right)$. Due to our considerations above,
this subspace is for small open neighbourhoods $\Set{O}\ni 0$
the image of the kernel of the free Dirac operator on $\Set{O}$
under the operator
$$\left(\begin{smallmatrix}
z & 0\\
0 & \Bar{z}
\end{smallmatrix}\right)\comp\sum\limits_{n=0}^{\infty}
\left(\Op{R}_{\mathbb{R}^2}(0,0,0)\comp\left(\begin{smallmatrix}
Vz/\Bar{z} & 0\\
0 & W\Bar{z}/z
\end{smallmatrix}\right)\right)^n.$$
The composition of this operator with the operator
$\unity-\Op{R}_{\mathbb{R}^2}(0,0,0)\comp\left(\begin{smallmatrix}
V & 0\\
0 & W
\end{smallmatrix}\right)$ is an operator on
$\banach{q}(\Set{O})\times\banach{q}(\Set{O})$.
For small $\Set{O}$ this composition belongs to small
neighbourhoods of $\left(\begin{smallmatrix}
z & 0\\
0 & \Bar{z}
\end{smallmatrix}\right)$ in the space of operators on
$\banach{q}(\Set{O})\times\banach{q}(\Set{O})$.
Therefore, for all spinors
$\psi\in\sobolev{1,p}(\Set{O})\times\sobolev{1,p}(\Set{O})$
in the kernel of $\left(\begin{smallmatrix}
V & \partial\\
-\Bar{\partial} & W
\end{smallmatrix}\right)$ there exists an unique
$\psi'\in\sobolev{1,p}(\Set{O})\times\sobolev{1,p}(\Set{O})$
in the same kernel, which has a zero of positive order at $0$,
such that $\left(\unity-
\Op{R}_{\mathbb{R}^2}(0,0,0)\comp\left(\begin{smallmatrix}
V & 0\\
0 & W
\end{smallmatrix}\right)\right)(\psi+\psi')$
is a constant spinor in the kernel of the free Dirac operator.
In particular, the quotient of this kernel modulo the subspace
of all spinors, which have at $0$ a zero of positive order,
is a complex two--dimensional space
(for potentials of the form $(U,\Bar{U})$
a quaternionic one--dimensional space).

\begin{Lemma}\label{quotients of kernels}
The quotient spaces of the kernel of the Dirac operator
$\left(\begin{smallmatrix}
V & \partial\\
-\Bar{\partial} & W
\end{smallmatrix}\right)$ with potentials $V,W\in\banach{2}(\Set{O})$
on a bounded open neighbourhood of $0\in\Set{O}\subset\mathbb{C}$
divided by the subspace of all spinors,
which have at $0$ a zero of at least $n$--th order,
is a complex $2n$--dimensional space
(for potentials of the form $(U,\Bar{U})$
a quaternionic $n$--dimensional space).\qed
\end{Lemma}

Actually these arguments show that if a spinor
in the kernel of $\left(\begin{smallmatrix}
V & \partial\\
-\Bar{\partial} & W
\end{smallmatrix}\right)$ has no zero at $0$, then it has no zero in a
neighbourhood of $0$. By definition of
\De{Order of zeroes}~\ref{order of zeroes}
this implies

\begin{Corollary}\label{isolated zeroes}
The zeroes of any element of the kernel of the Dirac operator
$\left(\begin{smallmatrix}
V & \partial\\
-\Bar{\partial} & W
\end{smallmatrix}\right)$ with potentials $V,W\in\banach{2}(\Set{O})$
on a bounded open set $\Set{O}\subset\mathbb{C}$
are isolated.\qed
\end{Corollary}

From the topological point of view the continuity of a cross section
is indispensable in order to determine the topological invariants of
the corresponding bundles. But the spinors in the kernels
of Dirac operators with potentials $V,W\in\banach{2}$ are
in general not continuous.
We shall overcome this problem with the help of sheaf theory.
In fact, the elements of these kernels of a Dirac operators with
potentials $V,W\in\banach{2}$ define sheaves,
which are very similar to sheaves of holomorphic functions.
In fact, due to the
\De{Strong unique continuation property}~\ref{strong unique continuation},
all local germs have at most one global continuation.
Furthermore, due to Lemma~\ref{quotients of kernels}
the corresponding local sheaves have the same properties
as the local sheaves of holomorphic functions.
Finally, due to Lemma~\ref{branchpoints}
the residue can be extended to a pairing of `meromorphic' elements in
the kernel of Dirac operator with potentials $V,W\in\banach{2}$.
In fact, if locally near
$z_l\in\Set{O}\subset\mathbb{C}$ $\psi$ and $\phi$ obey
\begin{align*}
\begin{pmatrix}
V\left(\Bar{z}-\Bar{z}_l\right)^{M_l} & 
\partial\left(\Bar{z}-\Bar{z}_l\right)^{M_l}\\
-\Bar{\partial}\left(z-z_l\right)^{M_l} &
W\left(z-z_l\right)^{M_l}\end{pmatrix}
\psi&=0&
\begin{pmatrix}
V\left(z-z_l\right)^{M_l} & 
\Bar{\partial}\left(z-z_l\right)^{M_l}\\
\partial\left(\Bar{z}-\Bar{z}_l\right)^{M_l} &
W\left(\Bar{z}-\Bar{z}_l\right)^{M_l}\end{pmatrix}
\phi^{t}&=0,
\end{align*}
respectively, then the residue of $\psi$ and $\phi$ at $z=z_l$ is defined as
$$\res\limits_{z=z_l}\left(\phi\begin{pmatrix}
0 & d\Bar{z}\\
dz & 0
\end{pmatrix}\psi\right)=\frac{1}{2\pi\sqrt{-1}}
\int\limits_{\partial B(0,\varepsilon)}
\phi\begin{pmatrix}
0 & d\Bar{z}\\
dz & 0
\end{pmatrix}\psi.$$
Due to Remark~\ref{sobolev embedding}
and the Sobolev embedding theorem \cite[5.4 Theorem]{Ad}
the restriction of $\psi$ and $\phi$ to one--dimensional submanifolds
belong to all $\banach{q}$ with $q<\infty$.
Therefore, the integral on the right hand side is well defined.
Moreover, due to Lemma~\ref{branchpoints} this residue does not depend
on $\varepsilon$ or more generally on the path around $z=0$.
As a consequence the sheaf--theoretic proofs
of the Riemann--Roch theorem \cite[Theorem~16.9]{Fo}
and S\'{e}rre Duality \cite[\S17]{Fo} carry over
to these sheaves of elements in the kernel of
Dirac operators with potentials $V,W\in\banach{2}$.

In order to establish the relation between Dirac operators with
potentials and quaternionic holomorphic structures
\cite[Definition~2.1]{FLPP} on quaternionic line bundles, we identify
the quaternions with all complex $2\times 2$--matrices of the form
$\left(\begin{smallmatrix}
a & b\\
-\Bar{b} & \Bar{a}
\end{smallmatrix}\right)$.
If we consider a $\mathbb{C}^2$--valued function
$\psi=\left(\begin{smallmatrix}
\psi_1\\
\psi_2
\end{smallmatrix}\right)$ on an open set $\Set{O}\subset\mathbb{C}$
as a quaternionic--valued function
$\left(\begin{smallmatrix}
\psi_1 & -\Bar{\psi}_2\\
\psi_2 & \Bar{\psi}_1
\end{smallmatrix}\right)$, then $-\Op{J}$ times the Dirac operator
$$-\Op{J}\comp\begin{pmatrix}
U & \partial\\
-\Bar{\partial} & \Bar{U}
\end{pmatrix}=\begin{pmatrix}
\Bar{\partial} & -\Bar{U}\\
U & \partial
\end{pmatrix}$$
yields a quaternionic holomorphic structure
on the trivial quaternionic line bundle on $\Set{O}$
endowed with the complex structure of multiplication on the left
with complex numbers $\mathbb{C}\subset\mathbb{Q}$.
In particular, the action of $\sqrt{-1}$ is given by left--
multiplication with $\left(\begin{smallmatrix}
\sqrt{-1} & 0\\
0 & -\sqrt{-1}
\end{smallmatrix}\right)$. In fact, we have
$$\begin{pmatrix}
\Bar{\partial}\psi_1 & -\Bar{\partial}\Bar{\psi}_2\\
\partial\psi_2 & \partial\Bar{\psi}_1
\end{pmatrix}+\begin{pmatrix}
0 & -\Bar{U}\\
U & 0
\end{pmatrix}\begin{pmatrix}
\psi_1 & -\Bar{\psi}_2\\
\psi_2 & \Bar{\psi}_1
\end{pmatrix}=\begin{pmatrix}
\Bar{\partial}\psi_1-\Bar{U}\psi_2 & -\overline{\partial\psi_2+U\psi_1}\\
\partial\psi_2+U\psi_1 & \overline{\Bar{\partial}\psi_1-\Bar{U}\psi_2}
\end{pmatrix}.$$
Hence the corresponding Hopf field is equal to $Q=-\Bar{U}d\Bar{z}$.
Lemma~\ref{order of zeroes 1} now implies, that if $U$ belongs to
$\banach{2,1}(\Set{O})$, then any holomorphic section of
the corresponding quaternionic holomorphic line bundle defines a
non--vanishing section of another
quaternionic holomorphic line bundle,
whose degree differs from the degree of the former line bundle
by the number of zeroes of the section of the former line bundle.
Moreover, the Hopf field $\Tilde{Q}$ of the latter line bundle belongs
locally also to $\banach{2,1}(\Set{O})$ and the Willmore energy
\cite[Definition~2.3]{FLPP} of the latter is the same
as the Willmore energy of the former.
If on $\Set{O}$ $\psi$ is a non--vanishing spinor of the kernel of
$\left(\begin{smallmatrix}
U & \partial\\
-\Bar{\partial} & \Bar{U}
\end{smallmatrix}\right)$, then the corresponding
quaternionic flat connection may be decomposed into
a holomorphic structure and an anti--holomorphic structure
\cite[(19)]{FLPP}. The holomorphic structure is equal to
the holomorphic structure, which is induced by the composition
of the Dirac operator with $-\Op{J}$. Furthermore, the
anti--holomorphic structure has the form
\begin{equation*}\begin{split}\begin{pmatrix}
\partial\psi_1 & -\partial\Bar{\psi}_2\\
\Bar{\partial}\psi_2 & \Bar{\partial}\Bar{\psi}_1
\end{pmatrix}+\begin{pmatrix}
B & A\\
-\Bar{A} & \Bar{B}
\end{pmatrix}\begin{pmatrix}
\psi_1 & -\Bar{\psi}_2\\
\psi_2 & \Bar{\psi}_1
\end{pmatrix}=0\\
\begin{aligned}
\text{with}&&&
\begin{aligned}
B & =
-\frac{\Bar{\psi}_1\partial\psi_1+\psi_2\partial\Bar{\psi}_2}
      {\psi_1\Bar{\psi_1}+\psi_2\Bar{\psi}_2} &
A & =
-\frac{\Bar{\psi}_2\partial\psi_1-\psi_1\partial\Bar{\psi}_2}
      {\psi_1\Bar{\psi_1}+\psi_2\Bar{\psi}_2}.
\end{aligned}\end{aligned}\end{split}\end{equation*}
Due to the decomposition \cite[(20)]{FLPP} the operator $\partial+B$
defines an anti--holomorphic structure on the corresponding
complex line bundle, and $A$ yields the part of the corresponding
anti--holomorphic structure on the quaternionic line bundle, which
anti--commutes with the action of $\sqrt{-1}$ (i.\ e.\ the left
multiplication with $\sqrt{-1}\in\mathbb{C}\subset\mathbb{Q}$).
If the Hopf field belongs locally to $\banach{2,1}(\Set{O})$, then the
eigen--spinors are continuous and belong to the Sobolev space
$\sobolev{1,2}(\Set{O})\times\sobolev{1,2}(\Set{O})$.
Hence for non--vanishing $\psi$ the potentials $A$ and $B$ belong to
$\banach{2}(\Set{O})$.

\begin{Lemma}\label{local dual potentials}
Let $U\in\banach{2}(\Set{O})$ be a potential with small
$\banach{2}$--norm, such that the spinor
\begin{align*}
\psi&=\left(\unity-\Op{R}_{\mathbb{R}^2}(0,0,0)\comp\left(\begin{smallmatrix}
U & 0\\
0 & \Bar{U}
\end{smallmatrix}\right)\right)^{-1}\left(\begin{smallmatrix}
a\\
b
\end{smallmatrix}\right)&
\text{with $(a,b)\in\mathbb{P}^1$}&
\end{align*}
is a well defined non--vanishing spinor on $\Set{O}$.
Then the parts $A$ and $B$ of the corresponding
anti--holomorphic structure, which anti--commute and commute
with the action of $\sqrt{-1}$,
have finite $\banach{2}_{\text{\scriptsize\rm loc}}$--
and $\banach{2,\infty}_{\text{\scriptsize\rm loc}}$--norms, respectively.
Finally, the derivative $\Bar{\partial}B$ of $B$ defines a real signed
measure on $\mathbb{P}^1$, which contains no point measures.
\end{Lemma}

\begin{proof}
We just showed that for non--vanishing spinors $\psi$ in the kernel of
the Dirac operator $\left(\begin{smallmatrix}
U & \partial\\
-\Bar{\partial} & \Bar{U}
\end{smallmatrix}\right)$ on $\Set{O}$ with a potential
$U\in\banach{2,1}(\Set{O})$ the potentials $A$ and $B$ of the corresponding
anti--holomorphic structure belong to $\banach{2}(\Set{O})$.
The zero curvature equation of the connection may be written as
$$\begin{pmatrix}
\Bar{\partial}B & \Bar{\partial}A\\
-\partial\Bar{A} & \partial\Bar{B}
\end{pmatrix}-\begin{pmatrix}
B & 0\\
0 & \Bar{B}\end{pmatrix}\comp\begin{pmatrix}
0 & -\Bar{U}\\
U & 0
\end{pmatrix}-\begin{pmatrix}
0 & A\\
-\Bar{A} & 0
\end{pmatrix}\comp\begin{pmatrix}
B & A\\
-\Bar{A} & \Bar{B}
\end{pmatrix}=\begin{pmatrix}
0 & -\partial\Bar{U}\\
\Bar{\partial}U & 0
\end{pmatrix}-\begin{pmatrix}
0 & -\Bar{U}\\
U & 0
\end{pmatrix}^2.$$
The diagonal and off--diagonal terms of this equation describes
the vanishing of the parts of the curvature of the flat connection,
which commutes and anti--commutes with the action of $\sqrt{-1}$
\cite[proof of Theorem~4.4]{FLPP}:
\begin{align*}
\Bar{\partial}B&=
U\Bar{U}-A\Bar{A} &
\Bar{\partial}A+\partial\Bar{U}&=\Bar{B}A-B\Bar{U}.
\end{align*}
More invariantly these equations may be written as
\begin{align*}
d\left(Bdz\right)&=
Ud\Bar{z}\wedge\overline{Ud\Bar{z}}+
Adz\wedge\overline{Adz} &
d\left(Adz-\Bar{U}d\Bar{z}\right)&=
\overline{Bdz}\wedge Adz+
\Bar{U}d\Bar{z}\wedge Bdz.
\end{align*}
Hence for $U\in\banach{2,1}(\Set{O})$ the statement follows.

We shall extend this proof with a limiting argument
to small potentials $U\in\banach{2}(\Set{O})$.
In fact, for any small $U\in\banach{2}(\Set{O})$ we choose a
sequence of potentials in $\banach{2,1}(\Set{O})$ with limit $U$.
We extend all potentials to a slightly larger open domain $\Set{O}'$
which contains the closure $\Bar{\Set{O}}$, so that they vanish on the
relative complement of $\Set{O}$ in $\Set{O}'$.
Obviously, the corresponding sequence of spinors $\psi_n$ defined
in the lemma extend to $\Set{O}'$.
By definition these spinors are smooth on
$\Set{O}'\setminus\Bar{\Set{O}}$.
Furthermore, the sequence of integrals of the corresponding
one--forms $B_ndz$ along a closed path in
$\Set{O}'\setminus\Bar{\Set{O}}$ around $\Set{O}$ converges.
Since the sequence of measures $U_n\Bar{U}_nd^2x$ converges,
this implies that the sequence $A_n$ is a bounded sequence in
$\banach{2}(\Set{O})$. 
Due to the Banach--Alaoglu theorem \cite[Theorem~IV.21]{RS1},
this sequence $A_n$ has a weakly convergent subsequence
with limit $A$. Also the sequence of real signed measures
$(U_n\Bar{U}_n-A_n\Bar{A}_n)d^2x$ on $\Set{O}$
has a weakly convergent subsequence.
Finally, due to the equations $\Bar{\partial}B=U\Bar{U}-A\Bar{A}$,
the sequence of functions $B_n$ is bounded
in the \Em{Lorentz space} $\banach{2,\infty}(\Set{O})$.
Due to \cite[Chapter~2 Theorem~2.7. and Chapter~4 Corollary~4.8.]{BS}
this \Em{Lorentz space} is the dual space
of the corresponding \Em{Lorentz space} $\banach{2,1}$.
The sequence $B_n$ has also a weakly convergent subsequence with limit $B$
and $\Bar{\partial}B$ considered as a measure is equal to the limit of
the measures $\left(U_n\Bar{U}_n-A_n\Bar{A}_n\right)d^2x$.
Since the sequence of spinors $\psi_n$
converges in the Banach space of $\banach{q}$--spinors,
and since the sequences $A_n$ and $B_n$ both converge weakly,
the limit $\psi$ is anti--holomorphic with respect to
the anti--holomorphic structure defined by the limits $A$ and $B$.

It remains to prove that the function
$\Bar{\partial}B=-\Bar{\partial}\partial\ln
\left(\psi_1\Bar{\psi}_1+\psi_2\Bar{\psi}_2\right)$
considered as a measure contains no point measures.
If this finite Baire measure contains a point measure at $z=z'$
of negative mass smaller or equal to $-n\pi$,
then, due to the following Lemma~\ref{zygmund estimate},
the spinor $\Tilde{\psi}=\left(\begin{smallmatrix}
z-z' & 0\\
0 & \Bar{z}-\Bar{z}'
\end{smallmatrix}\right)^{-n}\psi$ belongs to
$\bigcap\limits_{q<\infty}\banach{q}_{\text{\scriptsize\rm loc}}(\Set{O})
\times\banach{q}_{\text{\scriptsize\rm loc}}(\Set{O})$.
This implies that $\psi$ has a zero of order $n$ at $z'$.
Hence, due to our assumptions, the masses of all point measures are
larger than $-\pi$. Again the following Lemma~\ref{zygmund estimate}
implies that $\psi/(\psi_1\Bar{\psi}_1+\psi_2\Bar{\psi}_2)$ is a
$\banach{2}_{\text{\scriptsize\rm loc}}$--spinor in the kernel of
$\left(\begin{smallmatrix}
\partial & A\\
-\Bar{A} & \Bar{\partial}
\end{smallmatrix}\right)$. Since these kernels are contained in
$\bigcap\limits_{q<\infty}\banach{q}_{\text{\scriptsize\rm loc}}(\Set{O})
\times\banach{q}_{\text{\scriptsize\rm loc}}(\Set{O})$,
Lemma~\ref{zygmund estimate} implies
that this measure contains no point measures.
\end{proof}

\begin{Lemma}\label{zygmund estimate}
If $\Set{O}$ denotes a bounded open subset of $\mathbb{C}$,
then for all finite signed Baire measures $d\mu$ on $\Set{O}$
\cite[Chapter~13 Section~5]{Ro2}
there exists a function $h$ in the Zygmund space
$\banach{}_{\exp}(\Set{O})$
\cite[Chapter~4 Section~6.]{BS} such that
$-\Bar{\partial}\partial h=d\mu$ (in the sense of distributions).
Moreover, if for a suitable $\varepsilon>0$ all $\varepsilon$--balls of $\Set{O}$
have measure smaller than $\pi/q$
with respect to the positive part $d\mu^+$ of the
Hahn decomposition of the finite signed Baire measure
$d\mu$ on $\Set{O}$ \cite[Chapter~11 Section~5]{Ro2},
then the exponentials $\exp(h)$
of all $h\in \banach{}_{\exp}(\Set{O})$ with
$-\Bar{\partial}\partial h=d\mu$ belong to
$\banach{q}_{\text{\scriptsize\rm loc}}(\Set{O})$.
Conversely, if the positive part $d\mu^+$ contains a point measure
with mass $\pi/q$, then the corresponding functions
$h\in \banach{}_{\exp}(\Set{O})$ with
$-\Bar{\partial}\partial h=d\mu$ does not belong to
$\banach{q}_{\text{\scriptsize\rm loc}}(\Set{O})$.
\end{Lemma}

\begin{proof}
Due to Dolbeault's Lemma \cite[Chapter~I Section~D 2.~Lemma]{GuRo}
the convolution with the function $-\frac{2}{\pi}\ln|z|$
defines a right inverse of the operator $-\Bar{\partial}\partial$.
Now we claim that the restriction of this convolution operator
defines a bounded operator from $\banach{1}(\Set{O})$ into the
Zygmund space $\banach{}_{\exp}(\Set{O})$.
Since the domain $\Set{O}$ is bounded, the claim
is equivalent to the analogous statement about the restriction
to $\Set{O}$  of the convolution with the non--negative function
$$f(z)=\begin{cases}-\frac{2}{\pi}\ln|z| & \text{if }|z|<1\\
                    0 & \text{if }1\leq |z|
       \end{cases}.$$
Associated to this function $f$ is its distribution function $\mu_f$
and its non--increasing rearrangement $f^{\ast}$
(\cite[Chapter~II \S3. Chapter~V \S3.]{SW},
\cite[Chapter~2 Section~1.]{BS} and \cite[Chapter~1. Section~8.]{Zi}):
\begin{align*}
\mu_f(s)&=\pi\exp\left(-\pi s\right) &
f^{\ast}(t)&=\begin{cases}-\frac{\ln\left(t/\pi\right)}{\pi} &
                        \text{if }0\leq t\leq\pi\\
                       0 & \text{if } \pi\leq t
          \end{cases}.
\end{align*}
If $g\in\banach{1}(\Set{O})$,
then $g^{\ast\ast}(t)=\frac{1}{t}\int_{0}^{t}g^{\ast}(s)ds$
is bounded by $\|g\|_1/t$, since the $\banach{1,\infty}$--norm
$\|g\|_{(1,\infty)}=\sup\{tg^{\ast\ast}(t)\mid t>0\}=
\int_0^{\infty}g^{\ast}(t)dt$
coincides with the $\banach{1}$--norm \cite[Chapter~V (3.9)]{SW}.
Therefore, \cite[(1.8.14) and (1.8.15)]{Zi} in the proof of
\cite[1.8.8.~Lemma]{Zi} (borrowed from \cite[Lemma~1.5.]{O})
implies that the non--increasing rearrangement $h^{\ast}(t)$
of the convolution $h=f\ast g$ is bounded by
\begin{eqnarray*}
h^{\ast}(t)\leq h^{\ast\ast}(t)\leq h_2^{\ast\ast}(t)+h_1^{\ast\ast}(t)&\leq &
g^{\ast\ast}(t)\int\limits_{f^{\ast}(t)}^{\infty}\mu_{f}(s)ds-
\int\limits_{t}^{\infty}sg^{\ast\ast}(s)df^{\ast}(s)\\
&\leq &\frac{\|g\|_1}{t}\exp\left(-\pi f^{\ast}(t)\right)-
\|g\|_1\int\limits_{t}^{\infty}df^{\ast}(s)\\
&\leq &\frac{\|g\|_1}{\pi}+\|g\|_1f^{\ast}(t).
\end{eqnarray*}
Since by definition the non--increasing rearrangement $h^{\ast}(t)$
vanishes for all arguments,
which are larger than the Lebegues measure of $\Set{O}$,
we conclude the validity of the following estimate:
$$\int\limits_{0}^{\infty}
\left(\exp\left(q h^{\ast}(t)\right)-1\right)dt\leq
\int\limits_{0}^{|\Set{O}|}
\exp\left(q h^{\ast}(t)\right)dt\leq
|\Set{O}|\exp\left(q\|g\|_1/\pi\right)\int\limits_{0}^{|\Set{O}|}
\left(t/\pi\right)^{-\|g\|_1q/\pi}dt,$$
with an obvious modification when $\pi<|\Set{O}|$.
Due to a standard argument \cite[Chapter~2 Exercise~3.]{BS}
the finiteness of this integral is equivalent to the statement that
$\exp|h|$ belongs to $\banach{q}(\Set{O})$.
To sum up, the exponential $\exp(h)$ of the convolution $h=f\ast g$
belongs to $\banach{q}(\Set{O})$, if $q<\frac{\pi}{\|g\|_1}$.
This proves the claim.
In particular, for all $g\in \banach{1}(\Set{O})$ there exists
an element $h\in \banach{}_{\exp}(\Set{O})$ with
$-\Bar{\partial}\partial h=g$.

Due to \cite[Chapter~4 Theorem~6.5]{BS}
$\banach{}_{\exp}(\Set{O})$ is the dual space of
the Zygmund space $\banach{}\log\banach{}(\Set{O})$. Hence we shall improve the
previous estimate and show that the convolution with
$-\frac{2}{\pi}\ln|z|$ defines a bounded operator from
$\banach{}\log\banach{}(\Set{O})\subset \banach{1}(\Set{O})$ into
$C(\Set{O})\subset \banach{}_{\exp}(\Set{O})$.
By definition of the norm \cite[Chapter~4 Definition~6.3.]{BS}
$$\|g\|_{\banach{}\log\banach{}}=-\frac{1}{|\Set{O}|}
\int\limits_{0}^{|\Set{O}|}g^{\ast}(t)\ln(t/|\Set{O}|)dt=
\int\limits_{0}^{|\Set{O}|}g^{\ast\ast}(t)dt$$
we may improve the previous estimate to
\cite[(1.8.14) and (1.8.15)]{Zi}
$$h^{\ast\ast}(t)\leq
g^{\ast\ast}(t)\int\limits_{f^{\ast}(t)}^{\infty}\mu_{f}(s)ds-
\int\limits_{t}^{\infty}sg^{\ast\ast}(s)df^{\ast}(s)\leq
\frac{1}{\pi}\int\limits_{0}^{t}g^{\ast}(s)ds+
\int\limits_{t}^{\pi}g^{\ast\ast}(s)ds\leq \|g\|_{\banach{}\log\banach{}}.$$
This implies that in this case $h^{\ast\ast}(t)$ is bounded,
and consequently $h\in \banach{\infty}(\Set{O})$.
Furthermore, since the function $\ln|z|$ is continuous for $z\neq 0$,
the convolution with $-\frac{2}{\pi}\ln|z|$ is a bounded operator from
$\banach{}\log\banach{}(\Set{O})$ into the Banach space $C(\Set{O})$.
Finally, the dual of this operator yields a bounded operator from
the Banach space of finite signed Baire measures on $\Set{O}$
\cite[Chapter~13 Section~5 25.~Riesz Representation Theorem]{Ro2}
into $\banach{}_{\exp}(\Set{O})$.
More precisely, if the measure of $\Set{O}$ with respect to a
finite positive measure $d\mu$ is smaller than $\pi/q$, then the
exponential $\exp(h)$ of the corresponding function $h=f\ast d\mu$
belongs to $\banach{q}(\Set{O})$.

Due to Weyl's Lemma \cite[\S18.4. Lemma~4.10]{Co2} the difference of
two arbitrary functions $h_1$ and $h_2$
with $-\Bar{\partial}\partial h_1=-\Bar{\partial}\partial h_2=d\mu$
is analytic. Therefore, it suffices to show the second and third
statement of the lemma for the convolution of $-\frac{2}{\pi}\ln|z|$
with $d\mu$. Due to the boundedness of $\Set{O}$
we may neglect that part of this convolution,
where the former function is negative. Therefore, we may neglect the
negative part of $d\mu$ in order to bound the exponential $\exp(h)$.
The decomposition of the convolution into
an $\varepsilon$--near and an $\varepsilon$--distant part
analogous to the decomposition in the proof of
Lemma~\ref{weakly continuous resolvent} completes the proof.
\end{proof}

\begin{Proposition}\label{order of zeroes 2}
\index{order of zeroes $\text{\rm ord}_{z}(\psi)$}
Let $U$ be a potential in $\banach{2}(\Set{O})$ on a bounded open
neighbourhood of $0\in\Set{O}\subset\mathbb{C}$
and $\psi$ a non--trivial spinor in the kernel
of the Dirac operator $\left(\begin{smallmatrix}
U & \partial\\
-\Bar{\partial} & \Bar{U}
\end{smallmatrix}\right)$ on $\Set{O}$.
For any sequence of parameters $\parameter{t}_n\in(0,1]$ with limit $0$
the sequence $\psi_{\parameter{t}_n}$ of the family
$$\psi_{\parameter{t}}=\frac{\left(\begin{smallmatrix}
\Op{I}_{q,\parameter{t}} & 0\\
0 & \Op{I}_{q,\parameter{t}}
\end{smallmatrix}\right)\psi}{\left\|\left(\begin{smallmatrix}
\Op{I}_{q,\parameter{t}} & 0\\
0 & \Op{I}_{q,\parameter{t}}
\end{smallmatrix}\right)\psi\right\|_{\banach{q}(B(0,1))\times
\banach{q}(B(0,1))}}\text{ (with $2<q<\infty$)}$$
has a subsequence, whose restriction to any
$\banach{q}(B(0,\varepsilon))\times\banach{q}(B(0,\varepsilon))$
with $0<\varepsilon<\infty$ converges to $\left(\begin{smallmatrix}
z & 0\\
0 & \Bar{z}
\end{smallmatrix}\right)^{\text{\scriptsize\rm ord}_0(\psi)}
\left(\begin{smallmatrix}
a\\
b
\end{smallmatrix}\right)$
with $(a,b)\in\mathbb{P}^1$.
\end{Proposition}

\begin{proof}
Let us first prove the statement for the spinors considered in
Lemma~\ref{local dual potentials}.
These spinors belong locally to the kernels of the operators
\begin{align*}
\begin{pmatrix}
\partial & 0\\
0 & \Bar{\partial}
\end{pmatrix} &+\begin{pmatrix}
B & A\\
-\Bar{A} & \Bar{B}
\end{pmatrix} &\text{and }
\begin{pmatrix}
\Bar{\partial} & 0\\
0 & \partial
\end{pmatrix}&+\begin{pmatrix}
0 & -\Bar{U}\\
U & 0
\end{pmatrix}.
\end{align*}
Due to the \De{Generalized H\"older's inequality}~\ref{generalized hoelder}
the operators $\left(\begin{smallmatrix}
B & A\\
-\Bar{A} & \Bar{B}
\end{smallmatrix}\right)$ and $\left(\begin{smallmatrix}
0 & -\Bar{U}\\
U & 0
\end{smallmatrix}\right)$ define bounded operators from
$\banach{q}_{\text{\scriptsize\rm loc}}\times
\banach{q}_{\text{\scriptsize\rm loc}}$
into $\banach{p,q}_{\text{\scriptsize\rm loc}}\times
\banach{p,q}_{\text{\scriptsize\rm loc}}$,
with $2<q<\infty$ and $1/p=1/q+1/2$. 
Furthermore, if we restrict these operators to arbitrary small
balls $B(0,\varepsilon)$, then we may achieve
that the norms of these operators become arbitrarily small.
On the other hand,
due to the \De{Generalized Young's inequality}~\ref{generalized young},
the restrictions of the resolvents
$-\left(\begin{smallmatrix}
\partial & 0\\
0 & \Bar{\partial}
\end{smallmatrix}\right)^{-1}=
\Op{J}\comp\Op{R}_{\mathbb{R}^2}(0,0,0)$ and
$-\left(\begin{smallmatrix}
\Bar{\partial} & 0\\
0 & \partial
\end{smallmatrix}\right)^{-1}=
\Op{R}_{\mathbb{R}^2}(0,0,0)\comp\Op{J}$
define bounded operators from
$\banach{p,q}(B(0,\varepsilon))\times\banach{p,q}(B(0,\varepsilon))$ into
$\banach{q}(B(0,\varepsilon))\times\banach{q}(B(0,\varepsilon))$.
Therefore, the discussion of spinors in the kernel of Dirac operators
previous to of Lemma~\ref{order of zeroes 1}
carries over and show that the spinors
\begin{align*}
\left(\unity-\Op{J}\comp\Op{R}_{\mathbb{R}^2}(0,0,0,)
\comp\begin{pmatrix}
B & A\\
-\Bar{A} & \Bar{B}
\end{pmatrix}\right)&\psi&
\left(\unity-\Op{R}_{\mathbb{R}^2}(0,0,0)\comp\Op{J}
\comp\begin{pmatrix}
0 & -\Bar{U}\\
U & 0
\end{pmatrix}\right)&\psi
\end{align*}
belong to the kernel of $\left(\begin{smallmatrix}
\partial & 0\\
0 & \Bar{\partial}
\end{smallmatrix}\right)$ and $\left(\begin{smallmatrix}
\Bar{\partial} & 0\\
0 & \partial
\end{smallmatrix}\right)$, respectively.
The components of spinors in these kernels belong to the
Bergman spaces, i.\ e.\ the closed subspaces of $\banach{q}(B(0,\varepsilon))$
of holomorphic and anti--holomorphic functions, respectively
(\cite[Chapter~1.]{HKZ} and \cite[Chapter~4.]{Zh}).
When $\varepsilon$ tends to zero,
then these spinors converge to $\psi$.
This implies that the components of $\psi$ converge
to the corresponding Bergman projections of $\psi$
(\cite[Chapter~1.2.]{HKZ} and \cite[\S4.3]{Zh})
and the analogous projections of $\psi$
onto the anti--holomorphic functions, respectively.
These Bergman projections define bounded operators
onto the closed subspaces of $\banach{q}(B(0,\varepsilon))$
of holomorphic functions
(\cite[Theorem~1.10]{HKZ} and \cite[4.2.3.~Theorem]{Zh}).
Since the Bergman projections map the anti--holomorphic functions
onto the constant functions
(and analogously the projections onto the anti--holomorphic functions
maps the holomorphic functions onto the constant functions),
$\psi$ converges for $\varepsilon\downarrow 0$ to a constant spinor.
The compactness of the space of all constant spinors
$\left(\begin{smallmatrix}
a\\
b
\end{smallmatrix}\right)$ with
$\left\|\left(\begin{smallmatrix}
a\\
b
\end{smallmatrix}\right)
\right\|_{\banach{q}(B(0,1))\times\banach{q}(B(0,1))}^q
=\pi\left(|a|^q+|b|^q\right)=1$ proves the statement for those spinors,
which were considered in Lemma~\ref{local dual potentials}.

If the norms of the restriction of the potentials $U$, $A$ and $B$
to $B(0,\parameter{t})$ are very small compared with
$\varepsilon^{2/q'-2/q}$ with $0<\varepsilon<1$ and $q'<q<\infty$,
then these spinors obey the estimate
$$\varepsilon^{2/q'}\|\psi\|_{\banach{q}(B(0,\parameter{t}))\times 
            \banach{q}(B(0,\parameter{t}))}\leq
\|\psi\|_{\banach{q}(B(0,\parameter{t}\varepsilon))\times
            \banach{q}(B(0,\parameter{t}\varepsilon))}$$
uniformly with respect to $U$, $A$ and $B$ and for $\varepsilon$
bounded from below.
In fact, this follows from the arguments above, since the constant
spinors obey this estimate for $q'=q$.
Reiterated application implies that
these spinors fulfill for all $q'<q$
in the limit $\varepsilon\downarrow 0$ the estimate
$\text{\bf{O}}\left(
  \varepsilon^{2/q'}\right)\leq
  \left\|\psi\right\|_{\banach{q}(B(0,\varepsilon))\times
  \banach{q}(B(0,\varepsilon))}$.
On the other hand, since all spinors $\psi$ in the kernel of
$\left(\begin{smallmatrix}
U & \partial\\
-\Bar{\partial} & \Bar{U}
\end{smallmatrix}\right)$ belong to
$\bigcap\limits_{2<q<\infty}
\banach{q}_{\text{\scriptsize\rm loc}}(\Set{O})\times
\banach{q}_{\text{\scriptsize\rm loc}}(\Set{O})$,
H\"older's inequality \cite[Theorem~III.1~(c)]{RS1}
implies for all $q<q''$ the estimate
$\left\|\psi\right\|_{\banach{q}(B(0,\varepsilon))\times
  \banach{q}(B(0,\varepsilon))}\leq\text{\bf{O}}\left(
\varepsilon^{\text{\scriptsize\rm ord}_0(\psi)+2/q''}\right)$.
Now the transformation
$\psi\mapsto\Tilde{\psi}=\left(\begin{smallmatrix}
z & 0\\
0 & \Bar{z}
\end{smallmatrix}\right)^{-\text{\scriptsize\rm ord}_0(\psi)}\psi$
together with Lemma~\ref{quotients of kernels} shows
that the statement for general $\psi$ follows from the statement for
those $\psi$ considered in Lemma~\ref{local dual potentials}.
\end{proof}

Summing up, we have the following equivalent definitions of the
\Em{order of a zero} of $\psi$:
\begin{description}
\index{order of zeroes $\text{\rm ord}_{z}(\psi)$}
\item[(i)] $\text{\rm ord}_0(\psi)$ is equal to the largest
  $m\in\mathbb{N}_0$, such that $\left(\begin{smallmatrix}
  z & 0\\
  0 & \Bar{z}
  \end{smallmatrix}\right)^{-m}\psi$ belongs to the kernel of
  $\left(\begin{smallmatrix}
  \Bar{\partial} & -\Bar{U}(\Bar{z}/z)^m\\
  U(z/\Bar{z})^m & \partial
  \end{smallmatrix}\right)$
  ($\iff$ \De{Order of zeroes}~\ref{order of zeroes}).
\item[(ii)] For all $2\leq q'<q<q''\leq\infty$ we have in the limit
  $\varepsilon\downarrow 0$
  $$\text{\bf{O}}\left(
  \varepsilon^{\text{\scriptsize\rm ord}_0(\psi)+2/q'}\right)\leq
  \left\|\psi\right\|_{\banach{q}(B(0,\varepsilon))\times
  \banach{q}(B(0,\varepsilon))}
  \leq\text{\bf{O}}\left(
  \varepsilon^{\text{\scriptsize\rm ord}_0(\psi)+2/q''}\right).$$
\item[(iii)] $\text{\rm ord}_0(\psi)$ is equal
  to the order of the zero of the limit constructed in
  Proposition~\ref{order of zeroes 2}
\end{description}

We do not know, whether in general these limits of
$\psi_{\parameter{t}_n}$ coincide for all sequences $\parameter{t}_n$
converging to zero.

\subsubsection{The general case}\label{subsubsection general case}

Now we are prepared to treat the general case. We shall construct
a convergent subsequence of a sequence of resolvents
$\Op{R}(U_n,\Bar{U}_n,k,0)$ corresponding to a bounded
sequence of complex potentials $U_n\in \banach{2}(\torus)$.
Firstly, due to the Banach--Alaoglu theorem
\cite[Theorem~IV.21]{RS1}, we may assume that the sequence of
potentials is weakly convergent, and that the sequence of
non--negative $\banach{1}$--functions $U_n\Bar{U}_n$ considered as
finite Baire measures on $\torus$ \cite[Chapter~13]{Ro2}
is also weakly convergent.
If the latter limit does not contain point measures
with mass larger or equal to $S_p^{-2}$
then due to Lemma~\ref{weakly continuous resolvent}
the corresponding resolvents converge.

\begin{Theorem}\label{limits of resolvents}
Any bounded sequence $U_n$ of complex potentials in
$\banach{2}(\mathbb{C}/\lattice_\mathbb{C})$
has a subsequence, whose resolvents $\Op{R}(U_n,\Bar{U}_n,k,0)$
converge weakly to a
\De{Finite rank Perturbation}~\ref{finite rank perturbations}
of the resolvent $\Op{R}(U,\Bar{U},k,0)$
of a weak limit of the sequence.
More precisely, the support of the
\De{Finite rank Perturbation}~\ref{finite rank perturbations}
is contained in the finite set of singular points $\{z_1,\ldots,z_L\}$,
which have measure larger or equal to $\pi$
with respect to a weak limit of the sequence of measures
$U_n\Bar{U}_nd^2x$.
Moreover, for all $\varepsilon>0$ the restrictions of the subsequence
of resolvents to the complements of
$\Set{S}_{D,\varepsilon}=B(z_1,\varepsilon)\cup\ldots\cup B(z_l,\varepsilon)$
converge when considered as operators from
$\banach{p}(\Delta\setminus\Set{S}_{D,\varepsilon})\times
\banach{p}(\Delta\setminus\Set{S}_{D,\varepsilon})$ to
$\banach{q'}(\Delta\setminus\Set{S}_{D,\varepsilon})\times
\banach{q'}(\Delta\setminus\Set{S}_{D,\varepsilon})$
with $1<p\leq 2$ and $2\leq q'<q=\frac{2p}{2-p}$.
\end{Theorem}

\begin{proof}
We prove this theorem in four steps.
Since we shall prove the limits of the sequences of resolvents to be
\De{Finite rank Perturbations}~\ref{finite rank perturbations},
whose support is contained in the finite set of singular points,
we shall decompose the sequence of potentials into a sum
of sequences of potentials with disjoint support.
More precisely, the limit of one sequence,
which we call regular, yields the ordinary weak limit.
The limits of the other sequences, which we call singular,
and whose domains converge to the singular points, yields the
\De{Finite rank Perturbations}~\ref{finite rank perturbations}.
Hence in the first step we set up this decomposition.
In the second step we treat the cases of single decompositions,
in which for each singular point we have one sequence of
singular potentials.
In the third step we treat the more complicated cases,
in which for some of the sequences of singular potentials
the procedure of decomposition has to be iterated.

Since we consider only local potentials
(i.\ e.\ multiplication operators)
the corresponding limits are also local in the sense of
Remark~\ref{locality}. Therefore we treat only the case of one
singular point (i.\ e.\ $L=1$). The interested reader should fill in
the straight forward modifications of our arguments for the case of
several singular points. The crucial issue is the replacement of the
Taylor coefficients of the family of resolvents in the beginning of
the proof of Lemma~\ref{nontrivial limit 1} by the Taylor coefficients
of a \De{Finite rank Perturbation}~\ref{finite rank perturbations},
whose support is disjoint from $z=0$.

\noindent
{\bf 1.\ The decomposition.}
We shall decompose the sequence of potentials $U_n$ into a sum
$$U_n=U_{\text{\scriptsize\rm reg},n}+
\sum\limits_{l=1}^L U_{\text{\scriptsize\rm sing},n,l}$$
of potentials with disjoint support.
Here $U_{\text{\scriptsize\rm sing},n,1},\ldots,
U_{\text{\scriptsize\rm sing},n,L}$
are the restrictions of $U_n$ to small disjoint balls
$B(z_1,\varepsilon_{n,l}),\ldots,B(z_l,\varepsilon_{n,L})$,
whose radii $\varepsilon_{n,l}$ tend to zero.
Consequently, $U_{\text{\scriptsize\rm reg},n}$
are the restrictions of $U_n$ to the complements
of the union of these balls.
More precisely, we assume
\begin{description}
\index{condition!decomposition (i)--(ii)}
\index{decomposition!condition $\sim$ (i)--(ii)}
\item[Decomposition (i)] For all $l=1,\ldots,L$ the weak limit
  of the sequence of finite Baire measures
  $U_{\text{\scriptsize\rm sing},n,l}
  \Bar{U}_{\text{\scriptsize\rm sing},n,l}d^2x$
  \cite[Chapter~13]{Ro2}
  is equal to the point measures of the weak limit of
  $U_n\Bar{U}_nd^2x$ at $z_l$.
  Consequently, the weak limit of the measures
  $U_{\text{\scriptsize\rm reg},n}\Bar{U}_{\text{\scriptsize\rm reg},n}d^2x$
  is equal to the weak limit of the measures $U_n\Bar{U}_nd^2x$
  minus the corresponding point measures at $z_1\ldots,z_L$.
\end{description}
Obviously there are many possible choices of the sequences
$\varepsilon_{n,l}$ with this property
(e.\ g.\ for a unique choice of $\varepsilon_{n,l}$ the square of
the  $\banach{2}$--norm of $U_{\text{\scriptsize\rm sing},n,l}$ is equal
to the mass of the point measure at $z_l$ of the weak limit of
$U_n\Bar{U}_nd^2x$). Locally we may consider the
potentials $U_{\text{\scriptsize\rm sing},n,l}$ as potentials on
$\mathbb{P}^1$. We want to treat each of these $L$ sequences
of potentials analogously to a subsequence
$U_{\parameter{t}_n}$ of an
\Em{illustrative family of potentials}~\ref{subsubsection family}
with $\parameter{t}_n\rightarrow\infty$.
In general we should use arbitrary M\"obius transformations
instead of the scaling transformations
$z\mapsto\parameter{t}z$.
The action of the M\"obius group $SL(2,\mathbb{C})/\mathbb{Z}_2$
on $\mathbb{P}^1$
$$SL(2,\mathbb{C})\ni\begin{pmatrix}
a & b\\
c & d
\end{pmatrix}:\mathbb{P}^1\rightarrow\mathbb{P}^1,
z\mapsto\frac{az+b}{cz+d}$$
induces a unitary representation $\pi$
of the M\"obius group on $\banach{2}(\mathbb{C})$
\cite[Chapter~II. \S4.]{Kn}
$$\pi\left(\begin{pmatrix}
a & b\\
c & d
\end{pmatrix}\right):\banach{2}(\mathbb{C})\rightarrow
\banach{2}(\mathbb{C}),\text{ with }
\pi\left(\begin{pmatrix}
a & b\\
c & d
\end{pmatrix}\right)U(z)=
\frac{1}{|a-cz|^2}
U\left(\frac{dz-b}{a-cz}\right).$$
In doing so we transform each of these $L$ sequences of potentials
(considered as potentials on $\mathbb{P}^1)$ by a sequence of
M\"obius transformations $g_{n,l}$ in such a way that the
transformed sequence of potentials has the following property:
\begin{description}
\index{condition!decomposition (i)--(ii)}
\index{decomposition!condition $\sim$ (i)--(ii)}
\item[Decomposition (ii)] There exists some $\varepsilon>0$, such that
  the $\banach{2}$--norm of the restrictions of the transformed potentials
  $\pi\left(g_{n,l}\right)U_{\text{\scriptsize\rm sing},n,l}$
  to all $\varepsilon$--balls
  (with respect to the usual metric of $\mathbb{P}^1$) is bounded by
  the constant $C_p<S_p^{-1}$.
\end{description}
Such decompositions do not always exist. But we shall see now that,
if all points of $\mathbb{P}^1$ have measure smaller than $2S_p^{-2}$
with respect to the weak limit of the finite Baire measures
$U_n\Bar{U}_nd^2x$, then these decompositions indeed exist.
The free Dirac operator on $\mathbb{P}^1$ is invariant
under the compact subgroup
$SU(2,\mathbb{C})/\mathbb{Z}_2$ of the M\"obius group
($\simeq SL(2,\mathbb{C})/\mathbb{Z}_2$) as well as
the usual metric on $\mathbb{P}^1$.
Therefore, due to the global Iwasawa decomposition
\cite[Chapter~VI Theorem~5.1]{He},
it suffices to consider in place of the whole M\"obius group only the
semidirect product of the scaling transformations
($z\mapsto\exp(t)z$ with $t\in\mathbb{R}$)
with the translations ($z\mapsto z+z_0$ with $z_0\in\mathbb{C}$).
In the sequel we assume that all $g_{n,l}$
belong to this subgroup of the M\"obius group.

\begin{Lemma}\label{optimal Moebius}
If the square of the $\banach{2}$--norm of a complex potential $U$ on
$\mathbb{P}^1$ is smaller than $2S_p^{-2}$,
then there exists a constant $C_p<S_p^{-1}$ and a M\"obius transformation
such that the $\banach{2}$--norm of the restrictions
of the transformed potentials to all balls of radius $\pi/6$
is not larger than $C_p$.
\end{Lemma}

\begin{proof}
Let $r_{\max}(U)$ be the maximum of the set
$$\left\{r\mid
\text{the $\banach{2}$--norms of the restrictions of $U$
      to all balls of radius $r$ are not larger than $C_p$}\right\}.$$
Since the $\banach{2}$--norm of the restriction of $U$ to a ball
depends continuously on the center and the diameter of the ball,
this set has indeed a maximum.
Moreover, there exist balls with radius $r_{\max}(U)$,
on which the restricted potential has $\banach{2}$--norm equal to $C_p$.

We claim that there exists a M\"obius transformation $h$,
such that $r_{\max}(\pi(h)U)$ is the supremum of the set of all
$r_{\max}(\pi(g)U)$, where $g$ runs through the M\"obius group.
Let $g_n$ be a maximizing sequence of this set, i.\ e.\ 
the limit of the sequence $r_{\max}(\pi(g_n)U)$
is equal to the supremum of the former set. Since $r_{\max}(\pi(g)U)$ is
equal to $r_{\max}(U)$, if $g$ belongs to the subgroup
$SU(2,\mathbb{C})/\mathbb{Z}_2$ of isometries of the M\"obius group,
and due to the global Iwasawa decomposition
\cite[Chapter~VI Theorem~5.1]{He},
the sequence $g_n$ may be chosen in the semidirect product of the
scaling transformations $z\mapsto\exp(t)z$
with the translations $z\mapsto z+z_0$.
It is quite easy to see, that if the values of $t$ and $z_0$
corresponding to a sequence $g_n$ of such M\"obius transformations
are not bounded, then there exist arbitrary small balls,
on which the $\banach{2}$--norms of the restrictions of $\pi(g_n)U$
have subsequences converging to the $\banach{2}$--norm of $U$.
Hence if the $\banach{2}$--norm of $U$ is larger than $C_p$,
then the maximizing sequence of M\"obius transformations
can be chosen to be bounded and therefore has a convergent subsequence.
In this case the continuity implies the claim.
If the $\banach{2}$--norm of $U$ is not larger than $C_p$,
then $r_{\max}(\pi(g)U)$ does not depend on $g$ and the claim is obvious.

If the $\banach{2}$--norm of $U$ is smaller than $\sqrt{2}C_p$,
then all $r_{\max}(U)$--balls,
on which the restriction of $U$ has $\banach{2}$--norm equal to $C_p$,
have pairwise non--empty intersection.
In particular, all of them have non--empty intersection
with one of these balls.
Consequently, if $r_{\max}(U)$ is smaller than $\pi/6$,
then these $r_{\max}(U)$--balls are contained in one hemisphere.
In this case there exists a M\"obius transformation $g$,
which enlarges $r_{\max}(U)$ (i.\ e.\ $r_{\max}(U)<r_{\max}(\pi(g)U)$).
We conclude that there exist a M\"obius transformation $g$,
such that $r_{\max}(\pi(g)U)$ is smaller than $\pi/6$.
\end{proof}

The upper bound $2S_p^{-2}$ is sharp,
because for a sequence of $\banach{2}$--potentials on $\mathbb{P}^1$,
whose square of the absolute values considered as a sequence
of finite Baire measures converges weakly to the sum
of two point measures of mass $S_p^{-2}$ at opposite points,
the corresponding sequence of maxima of $r_{\max}(\pi(h)\cdot)$
converges to zero.
But the lower bound $\pi/6$ is of course not optimal.

If the $\banach{2}$--norms of the potentials
$U_{\text{\scriptsize\rm sing},n,l}$ are smaller than $\sqrt{2}C_p$,
then this lemma ensures the existence of M\"obius transformations
$g_{n,l}$ with the property \Em{Decomposition~(ii)}.
In general we showed the existence of a sequence of
M\"obius transformations $h_{n,l}$, which maximizes
$r_{\max}(\pi(g)U_{\text{\scriptsize\rm sing},n,l})$.
If for one $l=1,\ldots,L$ the correspondong sequences
$\pi\left(h_{n,l}\right)U_{\text{\scriptsize\rm sing},n,l}$
do not obey condition \Em{Decomposition~(ii)},
then we apply this procedure of decomposition to
the corresponding sequence of potentials
$\pi\left(h_{n,l}\right)U_{\text{\scriptsize\rm sing},n,l}$
on $\mathbb{P}^1$.
Consequently, we decompose the sequence
$U_{\text{\scriptsize\rm sing},n,l}$ into a finite sum
of potentials with disjoint support, such that the corresponding
potentials $\pi\left(g_{n,l}\right)U_{\text{\scriptsize\rm sing},n,l}$
on $\mathbb{P}^1$ obey
the analogous conditions \Em{Decomposition}~(i).
Due to Lemma~\ref{optimal Moebius} after finitely many iterations
of this procedure of decomposing the potentials into finite
sums of potentials with disjoint support, we arrive at a decomposition
$$U_n=U_{\text{\scriptsize\rm reg},n}+
\sum\limits_{l=1}^{L'} U_{\text{\scriptsize\rm sing},n,l}$$
of potentials with disjoint support.
More precisely, the potentials
$U_{\text{\scriptsize\rm sing},n,1},\ldots,
U_{\text{\scriptsize\rm sing},n,L'}$ are restrictions of $U_n$
either to small balls or to the relative complements
of finitely many small balls inside of small balls.
In particular, the domains of these potentials are
excluded either from the domain of $U_{\text{\scriptsize\rm reg},n}$,
or from the domain of another $U_{\text{\scriptsize\rm sing},n,l}$.
The former potentials obey condition \Em{Decomposition}~(i)
and the latter obey condition
\begin{description}
\index{condition!decomposition (i)--(ii)}
\index{decomposition!condition $\sim$ (i)--(ii)}
\item[Decomposition (i')] If the domain of
  $U_{\text{\scriptsize\rm sing},n,l'}$ is excluded from the domain
  of $U_{\text{\scriptsize\rm sing},n,l}$,
  then the weak limit of the sequence of finite Baire measures
  $\pi\left(g_{n,l}\right)U_{\text{\scriptsize\rm sing},n,l'}
  \pi\left(g_{n,l}\right)\Bar{U}_{\text{\scriptsize\rm sing},n,l'}d^2x$
  on $\mathbb{C}\subset\mathbb{P}^1$ converges weakly
  to the point measure of the weak limit of the sequence
  $\pi\left(g_{n,l}\right)U_{\text{\scriptsize\rm sing},n,l}
  \pi\left(g_{n,l}\right)\Bar{U}_{\text{\scriptsize\rm sing},n,l}d^2x$
  at some point of $\mathbb{C}$, whose measure with respect to the
  latter limit is not smaller than $S_p^{-2}$.
\end{description}
All these potentials $U_{\text{\scriptsize\rm sing},n,1},\ldots,
U_{\text{\scriptsize\rm sing},n,L'}$ obey condition
\Em{Decomposition}~(ii).
We remark that if the weak limit of the finite Baire measures
$\pi\left(h_{n,l}\right)U_{\text{\scriptsize\rm sing},n,l}
\pi\left(h_{n,l}\right)\Bar{U}_{\text{\scriptsize\rm sing},n,l}d^2x$
on $\mathbb{C}\subset\mathbb{P}^1$,
where $h_{n,l}$ denotes the sequence of M\"obius transformations
minimizing $R\left(\pi(g)U_{\text{\scriptsize\rm sing},n,l}\right)$,
contains at $z=\infty$ a point measure,
whose mass is not smaller than $S_p^{-2}$,
then we decompose from the sequence
$\pi\left(h_{n,l}\right)U_{\text{\scriptsize\rm sing},n,l}$
potentials, whose domains are the complement of a large ball in the
domains of these potentials.
In these cases the domains of the analog to the regular sequence
of the decomposition are excluded from the domains
of the analog to the singular sequence,
whose $\banach{2}$--norm accumulates at $z=\infty$.
Since the M\"obius transformations corresponding to the former are
faster divergent then the M\"obius transformations of the latter,
the latter should be considered as less singular than the former.
Therefore, also in this case the domains of the more singular sequences
are excluded from the domains of the less singular sequences.

By construction and due to Lemma~\ref{weakly continuous resolvent} the
sequence of resolvents
$\Op{R}(U_{\text{\scriptsize\rm reg},n},
\Bar{U}_{\text{\scriptsize\rm reg},n},k,0)$ converges in the norm
limit to the resolvent $\Op{R}(U,\Bar{U},k,0)$ of the weak limit
$U$ of the squence of potentials. We shall include the
perturbations of the sequences
$U_{\text{\scriptsize\rm sing},n,1},\ldots,
U_{\text{\scriptsize\rm sing},n,L'}$.
In doing so we start with the most singular of these sequences
(i.\ e.\ those $U_{\text{\scriptsize\rm sing},n,l}$, whose domains
are small balls) and include subsequently those sequences,
from whose domains the domains of the former sequences are excluded.
In the last step we include the sequences
$U_{\text{\scriptsize\rm reg},n}$.

\noindent
{\bf 2.\ The case of single decompositions.}
Here we assume that the decomposition obeys conditions
\Em{Decomposition}~(i)--(ii)
(i.\ e.\ the domains of all $U_{\text{\scriptsize\rm sing},n,l}$
are small balls) and $L=1$.
Due to Corollary~\ref{weakly continuous sphere resolvents}
and condition \Em{Decomposition}~(ii)
a subsequence of the sequence of resolvents
$\Op{R}_{\mathbb{P}^1}\left(
\pi\left(g_{n,1}\right)U_{\text{\scriptsize\rm sing},n,1},
\pi\left(g_{n,1}\right)\Bar{U}_{\text{\scriptsize\rm sing},n,1},
0\right)$
converges to the resolvent
$\Op{R}_{\mathbb{P}^1}\left(
U_{\text{\scriptsize\rm sing},1},
\Bar{U}_{\text{\scriptsize\rm sing},1},0\right)$
of the weak limit $U_{\text{\scriptsize\rm sing},1}$ of
$\pi\left(g_{n,1}\right)U_{\text{\scriptsize\rm sing},n,1}$.
Therefore, the arguments in Lemma~\ref{nontrivial limit 1} concerning
the case that the Dirac operator $(1+z\Bar{z})\left(\begin{smallmatrix}
U_{\parameter{t}=1} & \partial\\
-\Bar{\partial} & \Bar{U}_{\parameter{t}=1}
\end{smallmatrix}\right)$ on $\mathbb{P}^1$
has a trivial kernel, carry over to the case that the
Dirac operator
$(1+z\Bar{z})\left(\begin{smallmatrix}
U_{\text{\tiny\rm sing},1} & \partial\\
-\Bar{\partial} & \Bar{U}_{\text{\tiny\rm sing},1}
\end{smallmatrix}\right)$ on $\mathbb{P}^1$
has a trivial kernel.

In case that the latter Dirac operator
has a non--trivial kernel, we follow again the lines of the proof
of Lemma~\ref{nontrivial limit 1} concerning the analogous case.
In the notation of the proof of Lemma~\ref{nontrivial limit 1}
we use the formula (compare with Remark~\ref{fredholm on the plane})
\begin{multline*}
\Op{R}\left(U_{\text{\scriptsize\rm sing},n,1},
            \Bar{U}_{\text{\scriptsize\rm sing},n,1},k,0\right)-
\Op{R}\left(0,0,k,0\right)=
\Op{R}\left(0,0,k,0\right)\comp\\
\comp\begin{pmatrix}
\Op{I}_{\frac{4}{3},n} & 0\\
0 & \Op{I}_{\frac{4}{3},n}
\end{pmatrix}\comp
\left(\Op{A}_n+\Op{B}_n\right)^{-1}\comp\begin{pmatrix}
\pi(g_{n,1})U_{\text{\scriptsize\rm sing},n,1} & 0\\
0 & \pi(g_{n,1})\Bar{U}_{\text{\scriptsize\rm sing},n,1}
\end{pmatrix}\comp\begin{pmatrix}
\Op{I}_{4,n} & 0\\
0 & \Op{I}_{4,n}
\end{pmatrix}^{-1}\comp
\Op{R}\left(0,0,k,0\right).
\end{multline*}
Here $\Op{I}_{p,n}$ is the analog sequence of isometries
induced by the sequence of M\"obius transformations $g_{n,1}$.
Moreover, $\Spa{E}_{\text{\scriptsize\rm finite}}\subset
\banach{\frac{4}{3}}(\mathbb{R}^2)\times\banach{\frac{4}{3}}(\mathbb{R}^2)$
is the kernel of the operator
$\unity-\left(\begin{smallmatrix}
U_{\text{\tiny\rm sing},1} & 0\\
0 & \Bar{U}_{\text{\tiny\rm sing},1}
\end{smallmatrix}\right)\comp\Op{R}_{\mathbb{R}^2}(0,0,0)$ and
$\Spa{F}_{\text{\scriptsize\rm finite}}\subset
\banach{\frac{4}{3}}(\mathbb{R}^2)\times\banach{\frac{4}{3}}(\mathbb{R}^2)$
is the co--kernel of the operator
$\unity-\Op{R}_{\mathbb{R}^2}(0,0,0)\comp\left(\begin{smallmatrix}
U_{\text{\tiny\rm sing},1} & 0\\
0 & \Bar{U}_{\text{\tiny\rm sing},1}
\end{smallmatrix}\right)$ (compare with Remark~\ref{fredholm on the plane}).
Furthermore, the two sequences of operators $\Op{A}_n$ and $\Op{B}_n$
are the corresponding decompositions introduced in
Lemma~\ref{perturbation of inverse operators}
of the restriction of the sequence of operators
$$\begin{pmatrix}
\Op{I}_{\frac{4}{3},n} & 0\\
0 & \Op{I}_{\frac{4}{3},n}
\end{pmatrix}^{-1}\comp\left(\unity-\begin{pmatrix}
U_{\text{\scriptsize\rm sing},n,1} & 0\\
0 & \Bar{U}_{\text{\scriptsize\rm sing},n,1}
\end{pmatrix}\comp\Op{R}(0,0,k,0)\right)\comp\begin{pmatrix}
\Op{I}_{\frac{4}{3},n} & 0\\
0 & \Op{I}_{\frac{4}{3},n}
\end{pmatrix}$$
to the domains of $\pi(g_{n,1})U_{\text{\scriptsize\rm sing},n,1}$.
More precsiely, let $B(z_1,\varepsilon_{n,1})$ denote the sequence of
domains of $U_{\text{\scriptsize\rm sing},n,1}$.
Then these restrictions to the M\"obius transformed domains
$g_{n,1}B(z_1,\varepsilon_{n,1})\subset\mathbb{R}^2$
are considered as operators on
$\banach{\frac{4}{3}}(\mathbb{R}^2)\times\banach{\frac{4}{3}}(\mathbb{R}^2)$.
Also the subspaces, which are denoted in
Lemma~\ref{perturbation of inverse operators} by capital $\Spa{F}$,
are here again denoted by $\Op{R}_{\mathbb{R}^2}^{t}(0,0,0)\Spa{F}$
(compare with Remark~\ref{change of notation}).
By definition of $\Spa{E}_{\text{\scriptsize\rm finite}}$ and
$\Spa{F}_{\text{\scriptsize\rm finite}}$ the spaces
$\Op{R}_{\mathbb{R}^2}(0,0,0)\Spa{E}_{\text{\scriptsize\rm finite}}$
and $\Op{R}_{\mathbb{R}^2}^{t}(0,0,0)\Spa{F}_{\text{\scriptsize\rm finite}}$
are the kernel and co--kernel of the Diarc operator 
$\left(\begin{smallmatrix}
U_{\text{\tiny\rm sing},1} & \partial\\
-\Bar{\partial} & \Bar{U}_{\text{\tiny\rm sing},1}
\end{smallmatrix}\right)$, respectively.
So the numbers $n_1,\ldots,n_J$ and $m_1,\ldots,m_J$ here are defined as
those numbers, which occur as \De{Order of zeroes}~\ref{order of zeroes}
of the elements of the latter two spaces, respectively. In the situation
of Lemma~\ref{nontrivial limit 1} this definition gives the same numbers
as the definition therein.

With the help of Lemma~\ref{inverse operator}
we may now again determine the limits of the resolvents
$\Op{R}\left(U_{\text{\scriptsize\rm sing},n,l},
\Bar{U}_{\text{\scriptsize\rm sing},n,l},k,0\right)$.
Here the family $\Op{B}_{\parameter{t}}$
in the proof of Lemma~\ref{nontrivial limit 1}
is replaced by the sequence of finite rank operators $\Op{B}_n$.
We remark that the M\"obius transformations $g_{n,1}$
belong to the semidirect product of the scaling transformations
($z\mapsto\exp(t)z$ with $t\in\mathbb{R}$)
with the translations ($z\mapsto z+z_0$ with $z_0\in\mathbb{C}$).
Since both resolvents $\Op{R}(0,0,k,\lambda)$ and
$\Op{R}_{\mathbb{R}^2}(0,0,\lambda)$
are invariant under the translations,
it suffices to take the dilatations into account,
and that was done in Lemma~\ref{nontrivial limit 1}.
However, in addition to the arguments in the proof of that Lemma
we encounter two complications. Firstly, the domains
$g_{n,1}B(z_1,\varepsilon_{n,1})$ depend on $n$.
For all elements $\psi\in\Spa{E}_{\text{\scriptsize\rm finite}}$
and elements $\phi\in\Spa{F}_{\text{\scriptsize\rm finite}}$ let
\begin{eqnarray*}
\psi_{n,\text{\scriptsize\rm ren}}&=&
\frac{\left.\psi\right|_{g_{n,1}B(z_1,\varepsilon_{n,1})}}{
\left\|\Op{R}_{\mathbb{R}^2}(0,0,0)\psi\right\|_{
  \banach{4}\left(\mathbb{R}^2\setminus
  g_{n,1}B(z_1,\varepsilon_{n,1})\right)\times
  \banach{4}\left(\mathbb{R}^2\setminus
  g_{n,1}B(z_1,\varepsilon_{n,1})\right)}}\\
\phi_{n,\text{\scriptsize\rm ren}}&=&
\frac{\left.\phi\right|_{g_{n,1}B(z_1,\varepsilon_{n,1})}}{
\left\|\Op{R}_{\mathbb{R}^2}^{t}(0,0,0)\phi\right\|_{
  \banach{4}\left(\mathbb{R}^2\setminus
  g_{n,1}B(z_1,\varepsilon_{n,1})\right)\times
  \banach{4}\left(\mathbb{R}^2\setminus
  g_{n,1}B(z_1,\varepsilon_{n,1})\right)}},
\end{eqnarray*}
denote two sequences of restricted and renormalized spinors.
In order to carry over the proof of Lemma~\ref{nontrivial limit 1}
to the present situation we replace the pairings of the elements
$\psi\in\Spa{E}_{\text{\scriptsize\rm finite}}$ with the functions
$\left(\begin{smallmatrix}
\Bar{z}^m\\
0
\end{smallmatrix}\right)$ and $\left(\begin{smallmatrix}
0\\
z^m
\end{smallmatrix}\right)$ and of the elements
$\phi\in\Spa{E}_{\text{\scriptsize\rm finite}}$ with the functions
$\left(\begin{smallmatrix}
z^m\\
0
\end{smallmatrix}\right)$ and $\left(\begin{smallmatrix}
0\\
\Bar{z}^m
\end{smallmatrix}\right)$.
So we use the renormalized pairings
\begin{align*}
\lim\limits_{n\rightarrow\infty}
\left\langle\left\langle\psi_{n,\text{\scriptsize\rm ren}},
                        \left(\begin{smallmatrix}
\Bar{z}^m\\
0
\end{smallmatrix}\right)\right\rangle\right\rangle, &
\lim\limits_{n\rightarrow\infty}
\left\langle\left\langle\psi_{n,\text{\scriptsize\rm ren}},
                        \left(\begin{smallmatrix}
0\\
z^m
\end{smallmatrix}\right)\right\rangle\right\rangle &\text{and }
\lim\limits_{n\rightarrow\infty}
\left\langle\left\langle\phi_{n,\text{\scriptsize\rm ren}},
                        \left(\begin{smallmatrix}
z^m\\
0
\end{smallmatrix}\right)\right\rangle\right\rangle, &
\lim\limits_{n\rightarrow\infty}
\left\langle\left\langle\phi_{n,\text{\scriptsize\rm ren}},
                        \left(\begin{smallmatrix}
0\\
\Bar{z}^m
\end{smallmatrix}\right)\right\rangle\right\rangle.
\end{align*} 
We remark that due to Proposition~\ref{order of zeroes 2} the limits
\begin{align*}
\lim\limits_{n\rightarrow\infty}&
\Op{R}_{\mathbb{R}^2}(0,0,0)\left(\begin{smallmatrix}
\Op{I}_{\frac{4}{3},n} & 0\\
0 & \Op{I}_{\frac{4}{3},n}
\end{smallmatrix}\right)\psi_{n,\text{\scriptsize\rm ren}}
&\text{and }
\lim\limits_{n\rightarrow\infty}&
\Op{R}_{\mathbb{R}^2}^{t}(0,0,0)\left(\begin{smallmatrix}
\Op{I}_{\frac{4}{3},n} & 0\\
0 & \Op{I}_{\frac{4}{3},n}
\end{smallmatrix}\right)\phi_{n,\text{\scriptsize\rm ren}}
\end{align*}
exist and are non--trivial for non--trivial $\psi$ and $\phi$.

Secondly, we have to take the differences of
$\pi\left(g_{n,l}\right)U_{\text{\scriptsize\rm sing},n,l}$ from
$U_{\text{\scriptsize\rm sing},l}$ into account. Therefore in the
present situation the sequence $\Op{B}_n$ differs from the
corresponding sequence in Lemma~\ref{nontrivial limit 1} by a
sequence of finite rank operators.
The corresponding data $(\Spa{E},\Spa{F},\Tilde{\Mat{Q}})$
belong to a finite--dimensional Grassmannian
(compare with the seventh step of
Section~\ref{subsubsection finite rank perturbations}).
These Grassmannians are compact
\cite[Chapter~1 \S5.]{GrHa}.
Consequently those subsequences of resolvents,
which correspond to convergent sequences, converge to
\De{rank one Perturbations}~\ref{finite rank perturbations}.

Now we shall include into our considerations the sequence
$U_{\text{\scriptsize\rm reg},n}$.
For this purpose we use the tools of the third step of
Section~\ref{subsubsection finite rank perturbations}.
More precisely, we replace in the corresponding formulas
the operators $\Op{S}(k)$ by the sequence of operators
$$\Op{S}_n(k)=\begin{pmatrix}
\Op{I}_{\frac{4}{3},n} & 0\\
0 & \Op{I}_{\frac{4}{3},n}
\end{pmatrix}\comp
\left(\Op{A}_n+\Op{B}_n\right)^{-1}\comp\begin{pmatrix}
\pi(g_{n,1})U_{\text{\scriptsize\rm sing},n,1} & 0\\
0 & \pi(g_{n,1})\Bar{U}_{\text{\scriptsize\rm sing},n,1}
\end{pmatrix}\comp\begin{pmatrix}
\Op{I}_{4,n} & 0\\
0 & \Op{I}_{4,n}
\end{pmatrix}^{-1}.$$
Then the corresponding perturbations of
$\Op{R}\left(U_{\text{\scriptsize\rm reg},n},
                \Bar{U}_{\text{\scriptsize\rm reg},n},k,0\right)$
are equal to
\begin{multline*}
\Op{R}\left(U_{\text{\scriptsize\rm reg},n},
                \Bar{U}_{\text{\scriptsize\rm reg},n},k,0\right)\comp
\Op{S}_n\left(U_{\text{\scriptsize\rm reg},n},
              \Bar{U}_{\text{\scriptsize\rm reg},n},k\right)\comp
\Op{R}\left(U_{\text{\scriptsize\rm reg},n},
                \Bar{U}_{\text{\scriptsize\rm reg},n},k,0\right)
\text{, with}\\
\Op{S}_n\left(U_{\text{\scriptsize\rm reg},n},
              \Bar{U}_{\text{\scriptsize\rm reg},n},k\right)=
\left(\unity+\Op{S}_n(k)\comp\Op{R}(0,0,k,0)-
\Op{S}_n(k)\comp\Op{R}\left(U_{\text{\scriptsize\rm reg},n},
                          \Bar{U}_{\text{\scriptsize\rm reg},n},k,0\right)
\right)^{-1}\comp\Op{S}_n(k).
\end{multline*}
Moreover, in the present situation for all
$\psi\in\Spa{E}_{\text{\scriptsize\rm finite}}$
and $\phi\in\Spa{F}_{\text{\scriptsize\rm finite}}$
we redefine the sequences
\begin{eqnarray*}
\psi_{n,\text{\scriptsize\rm ren}}&=&
\frac{\left.\psi\right|_{g_{n,1}B(z_1,\varepsilon_{n,1})}}{
\left\|\Op{R}_{\mathbb{R}^2}\left(
\left.U_{\text{\scriptsize\rm reg},n}\right|_\Set{O},
\left.\Bar{U}_{\text{\scriptsize\rm reg},n}\right|_\Set{O},0\right)
\comp\left(\begin{smallmatrix}
\Op{I}_{\frac{4}{3},n} & 0\\
0 & \Op{I}_{\frac{4}{3},n}
\end{smallmatrix}\right)\psi\right\|_{
  \banach{4}\left(\Set{O}\setminus B(z_1,\varepsilon)\right)\times
  \banach{4}\left(\Set{O}\setminus B(z_1,\varepsilon)\right)}}\\
\phi_{n,\text{\scriptsize\rm ren}}&=&
\frac{\left.\phi\right|_{g_{n,1}B(z_1,\varepsilon_{n,1})}}{
\left\|\Op{R}_{\mathbb{R}^2}^{t}\left(
\left.U_{\text{\scriptsize\rm reg},n}\right|_\Set{O},
\left.\Bar{U}_{\text{\scriptsize\rm reg},n}\right|_\Set{O},0\right)
\comp\left(\begin{smallmatrix}
\Op{I}_{\frac{4}{3},n} & 0\\
0 & \Op{I}_{\frac{4}{3},n}
\end{smallmatrix}\right)\psi\right\|_{
  \banach{4}\left(\Set{O}\setminus B(z_1,\varepsilon)\right)\times
  \banach{4}\left(\Set{O}\setminus B(z_1,\varepsilon)\right)}}.
\end{eqnarray*}
Here $\Set{O}$ is a small open neighbourhood
(in order to ensure the existence of the resolvent),
which contains a ball $B(z_1,\varepsilon)$,
and $\left.U_{\text{\scriptsize\rm reg},n}\right|_\Set{O}$
are the restrictions of the  sequence of potentials
$U_{\text{\scriptsize\rm reg},n}$ to $\Set{O}$.
We claim that this normalization essentially does not depend on
$\Set{O}$ and $\varepsilon$. More precisely, for all non--trivial
$\psi$ and $\phi$ and $k\not\in\fermi(U,\Bar{U})$
and all appropriate choices of $\Set{O}$ and $\varepsilon$ the limits
\begin{eqnarray*}
\psi(k)&=&\lim\limits_{n\rightarrow\infty}
\Op{R}\left(U_{\text{\scriptsize\rm reg},n},
            \Bar{U}_{\text{\scriptsize\rm reg},n},k,0\right)
\comp\left(\begin{smallmatrix}
\Op{I}_{\frac{4}{3},n} & 0\\
0 & \Op{I}_{\frac{4}{3},n}
\end{smallmatrix}\right)\psi_{n,\text{\scriptsize\rm ren}}\\
\phi(k)&=&\lim\limits_{n\rightarrow\infty}
\Op{R}^{t}\left(U_{\text{\scriptsize\rm reg},n},
                \Bar{U}_{\text{\scriptsize\rm reg},n},k,0\right)
\comp\left(\begin{smallmatrix}
\Op{I}_{\frac{4}{3},n} & 0\\
0 & \Op{I}_{\frac{4}{3},n}
\end{smallmatrix}\right)\phi_{n,\text{\scriptsize\rm ren}}
\end{eqnarray*}
converge on $\Delta$ to non--trivial functions $\psi(k)$ and $\phi(k)$,
which belong on $\Delta\setminus\{z_1\}$ to the kernel and co--kernel
of $\Op{D}(U,\Bar{U},k)$, respectively,
and fulfill on a neighbourhood of $z_1$ the equations
\begin{align*}
\begin{pmatrix}
U\left(\Bar{z}-\Bar{z}_1\right)^{M_1} & 
\partial\left(\Bar{z}-\Bar{z}_1\right)^{M_1}\\
-\Bar{\partial}\left(z-z_1\right)^{M_1} &
\Bar{U}\left(z-z_1\right)^{M_1}\end{pmatrix}
\psi(k)&=0&\text{and }
\begin{pmatrix}
U\left(z-z_1\right)^{M_1} & 
\Bar{\partial}\left(z-z_1\right)^{M_1}\\
\partial\left(\Bar{z}-\Bar{z}_1\right)^{M_1} &
\Bar{U}\left(\Bar{z}-\Bar{z}_1\right)^{M_1}\end{pmatrix}
\phi(k)&=0.
\end{align*}
In fact, the sequence of resolvents
$\Op{R}\left(U_{\text{\scriptsize\rm reg},n},
            \Bar{U}_{\text{\scriptsize\rm reg},n},k,0\right)$
converges as compact operators from
$\banach{\frac{4}{3}}(\Delta)\times\banach{\frac{4}{3}}(\Delta)$
into $\banach{q}(\Delta)\times\banach{q}(\Delta)$ with $q<4$.
If we apply Proposition~\ref{order of zeroes 2}
to the integral kernels of these operators
with fixed $z$ and $z'$ in the complement of $B(z_1,\varepsilon)$,
then the claim follows from the considerations concerning the
analogous limits with trivial regular potentials
$U_{\text{\scriptsize\rm reg},n}$.
These limits yield the entries of the row--vector $\Psi_{D}(z,k)$
and column--vector $\Phi_{D}(z,k)$
(compare with \De{Finite rank Perturbations}~\ref{finite rank perturbations}).

Now we show that in the present situation again the sequence
of operators $\Op{A}_n$ does not contribute to the limit of
$\Op{R}(U_n,\Bar{U}_n,k,0)$.
If the numbers $n_J$ and $m_J$ are equal to $J$, then the sequence of
projections onto $\Spa{E}_{\text{\scriptsize\rm finite}}$ and
$\Spa{F}_{\text{\scriptsize\rm finite}}$ are again bounded.
Due to Lemma~\ref{weakly continuous resolvent}
we may replace the potentials $U_{\text{\scriptsize\rm reg},n}$ by the
restrictions $U'_{\text{\scriptsize\rm reg},n}$
of $U_n$ to the complement of $B(z_1,\varepsilon'_{n,1})$
with a sequence of larger radii
$0<\varepsilon_{n,1}<\varepsilon'_{n,1}$
converging to zero, such that
$\varepsilon_{n,1}/\varepsilon'_{n,1}$ converges to zero.
By inserting in the formula for the perturbations of
$\Op{R}\left(U'_{\text{\scriptsize\rm reg},n},
                \Bar{U}'_{\text{\scriptsize\rm reg},n},k,0\right)$
the equations
\begin{eqnarray*}
\Op{R}\left(U'_{\text{\scriptsize\rm reg},n},
            \Bar{U}'_{\text{\scriptsize\rm reg},n},k,0\right)&=&
\left(\unity-\Op{R}(0,0,k,0)\comp\begin{pmatrix}
U'_{\text{\scriptsize\rm reg},n} & 0\\
0 & \Bar{U}'_{\text{\scriptsize\rm reg},n}
\end{pmatrix}\right)^{-1}\comp\Op{R}(0,0,k,0)\\
&=&\Op{R}(0,0,k,0)\comp\left(\unity-\begin{pmatrix}
U'_{\text{\scriptsize\rm reg},n} & 0\\
0 & \Bar{U}'_{\text{\scriptsize\rm reg},n}
\end{pmatrix}\comp\Op{R}(0,0,k,0)\right)^{-1}
\end{eqnarray*}
we may achieve that the free resolvent $\Op{R}(0,0,k,0)$ occurs
at all places on the left hand sides and on the right hand sides
of the operators $\Op{S}_n(k)$. Hence it suffices to estimate the
operator $\Op{R}(0,0,k,0)$ as an operator from
$\banach{\frac{4}{3}}(\Delta\setminus B(0,\varepsilon'))\times
\banach{\frac{4}{3}}(\Delta\setminus B(0,\varepsilon'))$ into
$\banach{4}(B(0,\varepsilon))\times\banach{4}(B(0,\varepsilon))$.
The fourth power of the norm of this operator is smaller than
the integral
\begin{eqnarray*}
\int\limits_{B(0,\varepsilon)}
\int\limits_{\Delta\setminus B(0,\varepsilon')}
\frac{1}{|z-z'|^4}\frac{d^2x'd^2x}{\pi} &\leq &
\int\limits_{B(0,\varepsilon)}
\int\limits_{\mathbb{R}^2\setminus B(0,\varepsilon')}
\frac{1}{|z-z'|^4}\frac{d^2x'd^2x}{\pi} \\
&\leq &
\int\limits_{0}^{\varepsilon}
\int\limits_{\varepsilon'}^{\infty}\int\limits_{0}^{\pi}
\frac{4d\theta RdR rdr}{(R^2+r^2+2Rr\cos(\theta))^2}\\
&\leq &
\int\limits_{0}^{\varepsilon}
\int\limits_{\varepsilon'/r}^{\infty}\int\limits_{0}^{\pi}
\frac{4d\theta (R/r)d(R/r)}{(1+(R/r)^2+2(R/r)\cos(\theta))^2}\frac{dr}{r}\\
&\leq &
2\pi\int\limits_{0}^{\varepsilon}
\left(\frac{1}{((\varepsilon'/r)^2-1)^2}+\frac{1}{(\varepsilon'/r)^2-1}\right)
\frac{dr}{r}\\
&\leq &\pi\int\limits_{(\varepsilon'/\varepsilon)^2}^{\infty}
\frac{(1+x-1)dx}{x(x-1)^2}\leq\frac{\pi}{(\varepsilon'/\varepsilon)^2-1}.
\end{eqnarray*}
Here we used the equalities
\begin{multline*}
\int\limits_{b_0}^{\infty}\int\limits_{[0,\pi]}
\frac{4d\theta bdb}{(1+b^2-2b\cos(\theta))^2}=
\int\limits_{b_0}^{\infty}\int\limits_{[0,\pi]}
\frac{d\theta bdb}{b^2(\frac{1+b^2}{2b}-\cos(\theta))^2}=
\int\limits_{b_0}^{\infty}\frac{\pi\frac{1+b^2}{2b}bdb}
        {b^2\left(\left(\frac{1+b^2}{2b}\right)^2-1\right)^{3/2}}\\
=2\pi\int\limits_{b_0}^{\infty}\frac{(1+b^2)2bdb}{(b^2-1)^3}=
2\pi\int\limits_{b_0}^{\infty}
\left(\frac{2}{(b^2-1)^3}+\frac{1}{(b^2-1)^2}\right)db^2=
2\pi\left(\frac{1}{(b_0^2-1)^2}+\frac{1}{b_0^2-1}\right).
\end{multline*}
In fact, due to \cite[Chapter~4 Section~5.3]{Ah} the integral
$\int_0^\pi\frac{d\theta}{a-\cos(\theta)}$ is for $a>1$
equal to $\pi\left(a^2-1\right)^{-1/2}$. Furthermore, the
derivatives with respect to $a$ yields
with $a=\frac{1+b^2}{2b}$ and $b>1$ this formula.
Since $\varepsilon_{n,1}/\varepsilon'_{n,1}$ converges to zero,
the operators $\Op{A}_n$ do not contribute to the limit.
Moreover, due to the classification of
\De{Finite rank Perturbations}~\ref{finite rank perturbations}
in the sixth step of
Section~\ref{subsubsection finite rank perturbations}
the space of \De{rank one Perturbations}~\ref{finite rank perturbations}
of $\Op{R}(U,\Bar{U},k,0)$ is also compact.
We remark that due to the considerations of the ninth step in
Section~\ref{subsubsection finite rank perturbations},
the limits obey the equation in the definition of
\De{rank one Perturbations}~\ref{finite rank perturbations}.
Hence for $n_J=n_J=J$ the resolvents $\Op{R}(U_n,\Bar{U}_n,k,0)$
have subsequences, which converge to 
\De{Finite rank Perturbations}~\ref{finite rank perturbations}
of $\Op{R}(U,\Bar{U},k,0)$.
The arguments in Lemma~\ref{nontrivial limit 1}
concerning the general case of values of the sequences
$n_1,\ldots,n_J$ and $m_1,\ldots,m_J$ carry over,
and we have proven that all sequences, which allow a decomposition
obeying conditions \Em{Decomposition}~(i)--(ii) have subsequences,
whose corresponding resolvents converge to
\De{Finite rank Perturbations}~\ref{finite rank perturbations}
of $\Op{R}(U,\Bar{U},k,0)$.

\noindent
{\bf 3.\ The case of iterated decompositions.}
Finally, we extend our arguments to the case that the decomposition
obeys condition \Em{Decomposition}~(ii').
In this case for all singular points $z_l\in\{z_1,\ldots,z_L\}$
the limits of the sequence of resolvents corresponding
to the most singular sequences of potentials at $z_1$
may contribute to analogous \De{Finite rank Perturbations}
of the resolvents of Dirac operators on $\mathbb{P}^1$.
More generally, the \De{Finite rank Perturbations}
of the sequence of resolvents corresponding
to those sequences of potentials, whose domains are excluded of one
particular sequence of potentials, may contribute to analogous
\De{Finite rank Perturbations}
of the sequence of resolvents corresponding
to the particular sequence of potentials.
More precisely, an easy calculation shows
that the resolvents of \De{Finite rank Perturbations} of the usual
Dirac operators do not exist for different values of $\lambda$.
In case of resolvents on $\torus$
we overcame this problem by using the family of resolvents
$\Op{R}(U,\Bar{U},k,0)$.
As we already mentioned, the corresponding operators
$\triv{\Op{R}}(U,\Bar{U},k,0)$ may indeed be considered as a
family of resolvents, if we multiply the corresponding operators with
off--diagonal matrices (e.\ g.\ with $-\Op{J}$).
On $\mathbb{P}^1$ we may use the same trick and multiply the usual
Dirac operator on the left with the matrix $-\Op{J}$.
The corresponding resolvents of the transformed operators are
calculated by the usual method of perturbing the operator
$(1+z\Bar{z})\left(\begin{smallmatrix}
U & \partial\\
-\Bar{\partial} & \Bar{U}
\end{smallmatrix}\right)$ by the the potential
$(1+z\Bar{z})\left(\begin{smallmatrix}
0 & \lambda/(1+z\Bar{z})\\
-\lambda/(1+z\Bar{z}) & 0
\end{smallmatrix}\right)$.
Since $\lambda/(1+z\Bar{z})$ belongs to $\banach{2}(\mathbb{C})$,
the methods of Section~\ref{subsection resolvent}
can be used in order to develop the spectral theory
of this perturbed Dirac operators.
Therefore, the product of the inverse of these operators
multiplied on the right with $\Op{J}$ yields the resolvents of the
Dirac operators on $\mathbb{P}^1$
multiplied on the left with $-\Op{J}$.
With this trick the kernel of the \De{Finite rank Perturbations}
of Dirac operators on $\mathbb{P}^1$ analogous to the
\De{Finite rank Perturbations}~\ref{finite rank perturbations}
are well defined.
The whole discussion of the local behaviour of non--trivial elements
in the kernel of Dirac operators on $\mathbb{P}^1$
in Section~\ref{subsubsection local behaviour}
extends obviously to the elements in the kernels of
\De{Finite rank Perturbations} of Dirac operators on $\mathbb{P}^1$.
Whenever we consider \De{Finite rank Perturbations}
of Dirac operators on $\mathbb{P}^1$ in the following paragraphs,
then these operators should be considered as implicitly defined by
the corresponding \De{Finite rank Perturbations}
of these resolvents of Dirac--like operators.

\begin{Remark}
The generalization of Theorem~\ref{limits of resolvents}
as suggested in Remark~\ref{generalization}
has an obvious extension to these \De{Finite rank Perturbations}
of Dirac operators on $\mathbb{P}^1$.
\end{Remark}

Let us consider the situation, where at one singular point
(i.\ e.\ $L=1$) two sequences (i.\ e.\ $L'=2$) of potentials
$U_{\text{\scriptsize\rm sing},n,1}$ and
$U_{\text{\scriptsize\rm sing},n,2}$ accumulate in such a way,
that the transformed potentials
$\pi\left(g_{n,1}\right)U_{\text{\scriptsize\rm sing},n,2}$
accumulates at one point of $\mathbb{C}\subset\mathbb{P}^1$.
Here $g_{n,1}$ and $g_{n,2}$ denotes the corresponding sequences
of M\"obius transformations, such that conditions
\Em{Decomposition}~(i)--(ii) and (i')--(ii) are fulfilled,
respectively.
In a first step we determine the limit of the transformed
of the resolvents
$\Op{R}\left(U_{\text{\scriptsize\rm sing},n,1}+
U_{\text{\scriptsize\rm sing},n,2},
\Bar{U}_{\text{\scriptsize\rm sing},n,1}+
\Bar{U}_{\text{\scriptsize\rm sing},n,2},k,0\right)$
under the M\"obius transformations $g_{n,1}$.
In doing so we use a variant of the line of arguments concerning the
case of a decomposition of potentials obeying condition
\Em{Decomposition}~(i)--(ii) with one singular point (i.\ e.\ $L=L'=1$).
In fact, firstly we determine the sequence of resolvents
$\Op{R}\left(U_{\text{\scriptsize\rm sing},n,2},
\Bar{U}_{\text{\scriptsize\rm sing},n,2},k,0\right)$
transformed under the M\"obius transformations $g_{n,1}$.
Since the sequence of resolvents on corresponding to the
Dirac operators on $\mathbb{P}^1$ with potentials
$(1+z\Bar{z})\pi\left(g_{n,2}\right)U_{\text{\scriptsize\rm sing},n,2}$
converge, in the calculation of the limit of the former sequence
we have to transform the sequences of resolvents $\Op{R}(0,0,k,0)$
under the M\"obius transformations $g_{n,1}^{-1}g_{n,2}$.
Therefore, the corresponding sequence of scaling parameters
$\parameter{t}_n$ should be inserted in the Taylor series of
Lemma~\ref{nontrivial limit 1}.
Consequently, the arguments showing the convergence of a subsequence
in case conditions \Em{Decompostion}~(i)--(ii) are fulfilled,
carry over and show firstly that the sequence of resolvents
$\Op{R}\left(U_{\text{\scriptsize\rm sing},n,2},
\Bar{U}_{\text{\scriptsize\rm sing},n,2},k,0\right)$
transformed under the M\"obius transformations $g_{n,1}$
converges to an analogous \De{Finite rank Perturbation} of the
resolvents $\Op{R}_{\mathbb{P}^1}(0,0,0)$,
and secondly that the sequence of resolvents
$\Op{R}\left(U_{\text{\scriptsize\rm sing},n,1}+
U_{\text{\scriptsize\rm sing},n,2},
\Bar{U}_{\text{\scriptsize\rm sing},n,1}+
\Bar{U}_{\text{\scriptsize\rm sing},n,2},k,0\right)$
transformed under the M\"obius transformations $g_{n,1}$ converges
to an analogous \Em{Finite rank Perturbation} of
$\Op{R}_{\mathbb{P}^1}\left(
\pi\left(g_{n,1}\right)U_{\text{\scriptsize\rm sing},n,1},
\pi\left(g_{n,1}\right)\Bar{U}_{\text{\scriptsize\rm sing},n,1},0\right)$.
In case that these \De{Finite rank Perturbations} are not defined,
since the corresponding operators have non--trivial kernels,
we determine the whole family of resolvents parametrized by $\lambda$
of the Dirac--like operators mentioned above.

In a second step the sequence of resolvents
$$\Op{R}\left(U_{\text{\scriptsize\rm reg},n}+
U_{\text{\scriptsize\rm sing},n,1}+
U_{\text{\scriptsize\rm sing},n,2},
\Bar{U}_{\text{\scriptsize\rm reg},n}+
\Bar{U}_{\text{\scriptsize\rm sing},n,1}+
\Bar{U}_{\text{\scriptsize\rm sing},n,2},k,0\right)$$
is determined analogously to the case
of a sequence obeying conditions \Em{Decomposition}~(i)--(ii).
However, the
\Em{Finite rank perturbations}~\ref{finite rank perturbations}
of $\Op{R}(U,\Bar{U},k,0)$ are not induced by
the kernels and co--kernels of that Dirac operator on $\mathbb{P}^1$,
whose potentials are the weak limit of the sequence of potentials
$(1+z\Bar{z})\pi\left(g_{n,1}\right)U_{\text{\scriptsize\rm sing},n,1}$.
Rather they are induced by the kernels and co--kernels of that
\Em{Finite rank perturbation} of the former
Dirac operator on $\mathbb{P}^1$,
which is the limit of the sequence of Dirac operators
with potentials $(1+z\Bar{z})\pi\left(g_{n,1}\right)\left(
U_{\text{\scriptsize\rm sing},n,1}+U_{\text{\scriptsize\rm sing},n,1}
\right)$.

These arguments extend to the general case.
In fact, we start with the most singular
sequences of potentials. Inductively all kernels of the already
determined \De{Finite rank Perturbations} on the domains excluded
from one sequence of potentials give rise to
\De{Finite rank Perturbations} of the corresponding sequence of
Dirac operators on $\mathbb{P}^1$ and finally on
$\torus$.
\end{proof}

\begin{Remark}\label{extension of the symplectic form}
Lemma~\ref{residue} extends to these limits of sequences of
potentials. More precisely, the restriction of the operator
$\triv{\Op{A}}_f$ (defined in the proof of this lemma)
to the complement of the union of small balls around
the singular points $\Set{S}_{\varepsilon}$
(compare with Remark~\ref{restriction})
is well defined for all
\De{Finite rank Perturbations}~\ref{finite rank perturbations}.
Since the Dirac operator is a local operator, the commutator in (ii)
yields again the restriction of the corresponding variations
$(\var V_f,\var W_f)$ to the complement of $\Set{S}_{\varepsilon}$.
Moreover, for functions, whose support is contained in this complement,
the statements (iii) and (iv) are still valid.
Since the $\banach{2}$--functions, whose support are contained in the
complement of some $\Set{S}_{\varepsilon}$ with $\varepsilon>0$, are
dense in the corresponding Hilber space,
Proposition~\ref{global meromorphic function}
extends as well to these
\De{Finite rank Perturbations}~\ref{finite rank perturbations}.
\end{Remark}

\subsection{The compactified moduli spaces}
\label{subsection compactified moduli}

Due to the previous Section~\ref{subsection limits} all sequences
in $\moduli_{\lattice,\eta,\willmore}$ have subsequences,
which converge to a \Em{complex Fermi curve} of a
\De{Finite rank Perturbation}~\ref{finite rank perturbations}.
Since these moduli spaces $\moduli_{\lattice,\eta,\willmore}$
are subspaces of a separable compact metrizable space
(compare with Lemma~\ref{compact metric}),
the closures $\Bar{\moduli}_{\lattice,\eta,\willmore}$
\index{moduli space! $\Bar{\moduli}_{\lattice,\eta},
                      \Bar{\moduli}_{\lattice,\eta,\willmore},
                      \Bar{\moduli}_{\lattice,\eta,\sigma},
                      \Bar{\moduli}_{\lattice,\eta,\sigma,\willmore}$|(}
\index{moduli space!compactification of the $\sim$}
\index{compactification!of the moduli space}
are the sets of all accumulation points of these spaces.
Analogously $\Bar{\moduli}_{\lattice,\eta,\sigma,\willmore}$ denotes
the subspace of $\Bar{\moduli}_{\lattice,\eta,\willmore}$,
whose \Em{complex Fermi curves} are invariant under the involution
$\sigma$, and $\Bar{\moduli}_{\lattice,\eta}$ and
$\Bar{\moduli}_{\lattice,\eta,\sigma}$ denote the unions
\begin{align*}
\Bar{\moduli}_{\lattice,\eta}&=\bigcup\limits_{\willmore>0}
\Bar{\moduli}_{\lattice,\eta,\willmore} &
\Bar{\moduli}_{\lattice,\eta,\sigma}&=\bigcup\limits_{\willmore>0}
\Bar{\moduli}_{\lattice,\eta,\sigma,\willmore}.
\end{align*}
\index{moduli space! $\Bar{\moduli}_{\lattice,\eta},
                      \Bar{\moduli}_{\lattice,\eta,\willmore},
                      \Bar{\moduli}_{\lattice,\eta,\sigma},
                      \Bar{\moduli}_{\lattice,\eta,\sigma,\willmore}$|)}
In this section we shall investigate these compactified moduli spaces
$\Bar{\moduli}_{\lattice,\eta,\willmore}$.
Due to the
\De{Strong unique continuation property}~\ref{strong unique continuation}
we may replace in Lemma~\ref{continuity of fermi curves}
the \Em{compact open topology} of the resolvents $\Op{R}(V,W,k,0)$
by the \Em{compact open topology} of the restrictions
of the resolvents to the complement of $\Set{S}_{\varepsilon}$
(compare with Remark~\ref{restriction} and Theorem~\ref{limits of resolvents}).
Therfore Theorem~\ref{limits of resolvents} implies that
all elements of $\Bar{\moduli}_{\lattice,\eta}$
are \Em{complex Fermi curves} of
\De{Finite rank Perturbations}~\ref{finite rank perturbations}.
In particular, we shall see that the elements of
$\Bar{\moduli}_{\lattice,\eta,\willmore}$ may be realized by gluing
two planes $\mathbb{C}_{\xx{p}}^{\pm}$
(or $\mathbb{C}_{\yy{p}}^{\pm}$) along a combination of possibly
infinitely many \Em{horizontal cuts} and finitely many
\Em{pairs of parallel cuts} as described in Lemma~\ref{gluing rule}. 

\begin{Proposition}\label{asymptotic analysis 2}
\index{asymptotic analysis}
The elements $\fermi\in\Bar{\moduli}_{\lattice,\eta}$
may be realized by gluing two planes
$\mathbb{C}^{\pm}_{\xx{p}}$ along
at most infinitely many \Em{horizontal cuts} and
finitely many \Em{pairs of parallel cuts} as decsribed in
Lemma~\ref{gluing rule}. More precisely, for all
$\fermi\in\Bar{\moduli}_{\lattice,\eta}$
and all small $\varepsilon>0$
(the excluded discs $B(g(\xx{\gamma},k_{\kappa^{\pm}}),\varepsilon)$
should be disjoint) there exists a $\delta>0$, such that
$\mathbb{C}^{\pm}_{\xx{p}}$ does not contain any cuts in
\begin{eqnarray*}
\left\{\xx{p}\in\mathbb{C}\mid |\xx{p}|>1/\delta\right\}
\bigcap\limits_{\kappa\in\lattice\dual}
\left\{\xx{p}\in\mathbb{C}\mid
\left|\xx{p}-g(\xx{\gamma},k_{\kappa^{\pm}})\right|>\varepsilon\right\}&=&\\
=\left\{\xx{p}\in\mathbb{C}\mid |\xx{p}|>1/\delta\right\}
\setminus\bigcup\limits_{\kappa\in\lattice\dual\setminus\{0\}}
B(g(\xx{\gamma},k^{\pm}_{\kappa}),\varepsilon)&,&
\end{eqnarray*}
and at most one \Em{horizontal cut} in each of the excluded discs
$B(g(\xx{\gamma},k_{\kappa^{\pm}}),\varepsilon)$ indexed
by
$\kappa\in\lattice\dual_{\delta}=
\left\{\kappa\in\lattice\dual\mid
|g(\xx{\gamma},k^+_{\kappa})|>1/\delta\right\}$.
\index{index set!of horizontal cuts $\lattice\dual_{\delta}$}
Finally, the disc $B(0,1/\delta)$ contains finitely many
\Em{horizontal cuts} and \Em{pairs of parallel cuts}. The total
number of these cuts is bounded
by the number of $\kappa\in\lattice\dual$,
whose discs $B(g(\xx{\gamma},k_{\kappa^{\pm}}),\varepsilon)$
are not disjoint from $B(0,1/\delta)$.
\end{Proposition}

In spite of being in general different from the similar
index set in Theorem~\ref{asymptotic analysis 1} we denote the
index set $\lattice\dual_{\delta}$ by the same symbol. Since
$|g(\xx{\gamma},k^+_{\kappa})|=\frac{\|\xx{\gamma}\|}{2}\|\kappa\|=
\frac{\|\xx{\gamma}\|}{\sqrt{2}}\|k^+_{\kappa}\|=
\frac{\|\xx{\gamma}\|}{\sqrt{2}}\|k^-_{\kappa}\|$
the transformation
$\delta\mapsto\frac{\|\xx{\gamma}\|}{\sqrt{2}}\delta$
transforms the latter index set into the former.
Moreover, if the cuts are not disjoint, then the total
number of these cuts is given by the half of the degree of the
zero divisor of $d\xx{p}$ restricted to the corresponding set of
the \Em{complex Fermi curve}.

\begin{proof}
In a first step we improve the arguments of
Theorem~\ref{asymptotic analysis 1} and show that all
\Em{complex Fermi curves} in $\moduli_{\lattice,\eta,\willmore}$
have a representation of the form just described.
Due to Theorem~\ref{asymptotic analysis 1} near infinity
the \Em{complex Fermi curve} looks like two open planes,
which are connected by infinitely many small handles.
We shall show that these handles
may be represented by \Em{horizontal cuts}.
Due to arguments used in the proof of Theorem~\ref{asymptotic analysis 1}
these handles are small perturbations of the ordinary double point
of the free \Em{complex Fermi curve}
at $k=0$. For $k=0$ the kernel of $\triv{\Op{D}}(0,0,k)$ is the
two-dimensional space of the constant functions
($\sim\mathbb{C}\times\mathbb{C}\subset
\banach{p}(\torus)\times\banach{p}(\torus)$).
Therefore, we decompose the resolvent
$\triv{\Op{R}}(0,0,k,\lambda)$
into a sum of the reduced resolvent plus the singular matrix
$\left(\begin{smallmatrix}
\lambda & -\pi(k_2+\sqrt{-1}k_1)\\
-\pi(k_2-\sqrt{-1}k_1) & \lambda
\end{smallmatrix}\right)^{-1}$ acting on the constant functions:
$$\triv{\Op{R}}(0,0,k,\lambda)=
\triv{\Op{R}}_{\text{\scriptsize\rm red}}(0,0,k,\lambda)+
\triv{\Op{S}}(0,0,k,\lambda).$$
For small $k$ the reduced resolvent is holomorphic and bounded
and the singular part is a rank--two operator
with a pole at $\lambda^2=\pi^2g(k,k)$.
Consequently the resolvent $\triv{\Op{R}}(V,W,k,\lambda)$
has the form
\begin{multline*}
\triv{\Op{R}}(V,W,k,\lambda)=
\triv{\Op{R}}_{\text{\scriptsize\rm red}}(0,0,k,\lambda)
\comp\left(\unity-\left(\begin{smallmatrix}
V & 0\\
0 & W
\end{smallmatrix}\right)\comp
\triv{\Op{R}}_{\text{\scriptsize\rm red}}(0,0,k,\lambda)
\right)^{-1}+\\
+\left(\unity-
\triv{\Op{R}}_{\text{\scriptsize\rm red}}(0,0,k,\lambda)
\comp\left(\begin{smallmatrix}
V & 0\\
0 & W
\end{smallmatrix}\right)\right)^{-1}\comp
\triv{\Op{S}}(V,W,k,\lambda)\comp
\left(\unity-\left(\begin{smallmatrix}
V & 0\\
0 & W
\end{smallmatrix}\right)\comp
\triv{\Op{R}}_{\text{\scriptsize\rm red}}(0,0,k,\lambda)
\right)^{-1}.
\end{multline*}
Here $\triv{\Op{S}}(V,W,k,\lambda)$
is the following $2\times 2$--matrix
(acting on the constant functions):
\begin{multline*}
\triv{\Op{S}}(V,W,k,\lambda)=\\
\begin{aligned}
=&\triv{\Op{S}}(0,0,k,\lambda)\comp
\left(\unity-\left(\unity-\left(\begin{smallmatrix}
V & 0\\
0 & W
\end{smallmatrix}\right)\comp
\triv{\Op{R}}_{\text{\scriptsize\rm red}}(0,0,k,\lambda)
\right)^{-1}\comp\left(\begin{smallmatrix}
V & 0\\
0 & W
\end{smallmatrix}\right)\comp
\triv{\Op{S}}(0,0,k,\lambda)\right)^{-1}&&\\
=&\left(\unity-\triv{\Op{S}}(0,0,k,\lambda)
\comp\left(\begin{smallmatrix}
V & 0\\
0 & W
\end{smallmatrix}\right)\comp\left(\unity-
\triv{\Op{R}}_{\text{\scriptsize\rm red}}(0,0,k,\lambda)
\comp\left(\begin{smallmatrix}
V & 0\\
0 & W
\end{smallmatrix}\right)\right)^{-1}\right)^{-1}\comp
\triv{\Op{S}}(0,0,k,\lambda).&&
\end{aligned}\end{multline*}
The arguments in Section~\ref{subsection resolvent} carry over and
show that for small $V,W\in \banach{2}(\torus)$
the operators $\unity-\left(\begin{smallmatrix}
V & 0\\
0 & W
\end{smallmatrix}\right)\comp
\triv{\Op{R}}_{\text{\scriptsize\rm red}}(0,0,k,\lambda)$
and $\unity-
\triv{\Op{R}}_{\text{\scriptsize\rm red}}(0,0,k,\lambda)
\comp\left(\begin{smallmatrix}
V & 0\\
0 & W
\end{smallmatrix}\right)$ are invertible. More precisely,
the restriction of the operator
$\triv{\Op{R}}_{\text{\scriptsize\rm red}}(0,0,k,\lambda)\comp
\left(\unity-\left(\begin{smallmatrix}
V & 0\\
0 & W
\end{smallmatrix}\right)\comp
\triv{\Op{R}}_{\text{\scriptsize\rm red}}(0,0,k,\lambda)
\right)^{-1}$ to the orthogonal complement of the constant functions
is the resolvent of the restriction of $\triv{\Op{D}}(V,W,k)$ to
this orthogonal complement.

\begin{Remark}\label{orthogonal complement}
Dirac operators with potentials $\triv{\Op{D}}(V,W,k)$
do not preserve this orthogonal complement
in contrast to the the corresponding free Dirac operators
$\triv{\Op{D}}(0,0,k)$. But with the help of the natural inclusion
of the orthogonal complement of the constant functions into
$\banach{p}(\torus)\times\banach{p}(\torus)$
and the natural projection of this Banach space onto the orthogonal
complement of the constant functions the restriction to this
complement is well defined.
\end{Remark}

By the first resolvent formula \cite[Theorem~VI.5]{RS1}
we may express the resolvent $\Op{R}_{\lambda'}$ with argument
$\lambda'$ in terms of the resolvent $\Op{R}_{\lambda}$
with argument $\lambda$ and the resolvent with argument
$\frac{1}{\lambda-\lambda'}$ of the resolvent with argument $\lambda$:
$$(\lambda-\lambda')\Op{R}_{\lambda'}
=\Op{R}_{\lambda}\comp
\left(\frac{\unity}{\lambda-\lambda'}-\Op{R}_{\lambda}\right)^{-1}
=\left(\frac{\unity}{\lambda-\lambda'}-\Op{R}_{\lambda}\right)^{-1}
\comp\Op{R}_{\lambda}.$$
Now the operators $\triv{\Op{D}}(0,0,k)\comp
\triv{\Op{R}}_{\text{\scriptsize\rm red}}(0,0,k,\lambda)$
are invertible on
$\banach{p}(\torus)\times\banach{p}(\torus)$
with $1<p<2$, and the operators
$\triv{\Op{R}}_{\text{\scriptsize\rm red}}(0,0,k,\lambda)
\comp\triv{\Op{D}}(0,0,k)$ are invertible on
$\banach{q}(\torus)\times\banach{q}(\torus)$
with $2<q<\infty$. Therefore the arguments of
Section~\ref{subsection resolvent} show that the operators
$\unity-\left(\begin{smallmatrix}
V & 0\\
0 & W
\end{smallmatrix}\right)\comp
\triv{\Op{R}}_{\text{\scriptsize\rm red}}(0,0,k,\lambda)$ on
$\banach{p}(\torus)\times\banach{p}(\torus)$
and the operators $\unity-
\triv{\Op{R}}_{\text{\scriptsize\rm red}}(0,0,k,\lambda)
\comp\left(\begin{smallmatrix}
V & 0\\
0 & W
\end{smallmatrix}\right)$ on
$\banach{q}(\torus)\times\banach{q}(\torus)$
are invertible, if and only if $\lambda$ does not belong
to the spectrum of the restriction of $\triv{\Op{D}}(V,W,k)$
to the orthogonal complement of the constant functions.
Furthermore, the proof of Lemma~\ref{weakly continuous resolvent}
shows that the inverse of the former operators considered as operators
from $\banach{p}(\torus)\times\banach{p}(\torus)$
to $\banach{p'}(\torus)\times\banach{p'}(\torus)$
with $1<p'<p<2$ and the inverse of the latter operators considered as
operators from
$\banach{q}(\torus)\times\banach{q}(\torus)$
to $\banach{q'}(\torus)\times\banach{q'}(\torus)$
with $2<q'<q<\infty$
depend weakly continuously on the potentials on the
subsets described in Lemma~\ref{weakly continuous resolvent}.
Hence for all potentials
$V,W\in \banach{2}(\torus)$ and all small open
neighbourhoods $\Set{O}\subset\mathbb{C}^2$ of $0$
and all $\varepsilon>0$ there exists a $\delta>0$
such that for all $\kappa\in\lattice\dual_{\delta}$
and all $k\in\Set{O}$ the restriction of the operators
\begin{eqnarray*}
&&\left(\unity-\begin{pmatrix}
\psi_{-\kappa}V & 0\\
0 & \psi_{\kappa}W
\end{pmatrix}\comp
\triv{\Op{R}}_{\text{\scriptsize\rm red}}(0,0,k,0)
\right)^{\-1}\comp\begin{pmatrix}
\psi_{-\kappa}V & 0\\
0 & \psi_{\kappa}W
\end{pmatrix}\\
&=&\begin{pmatrix}
\psi_{-\kappa}V & 0\\
0 & \psi_{\kappa}W
\end{pmatrix}\comp\left(\unity-
\triv{\Op{R}}_{\text{\scriptsize\rm red}}(0,0,k,0)
\comp\begin{pmatrix}
\psi_{-\kappa}V & 0\\
0 & \psi_{\kappa}W
\end{pmatrix}\right)^{-1}
\end{eqnarray*}
to the constant functions defines a $2\times 2$--matrix
$\Mat{A}(\psi_{-\kappa}V,\psi_{\kappa}W,k)$,
which is bounded by $\varepsilon$.
More precisely, these $2\times 2$--matrices
$\Mat{A}(\psi_{-\kappa}V,\psi_{\kappa}W,k)$
are defined as the matrix elements of the given operators
with respect to the natural orthonormal basis
of the constant functions
$\simeq \mathbb{C}\times\mathbb{C}\subset
\banach{2}(\torus)\times\banach{2}(\torus)$.
Due to Lemma~\ref{inverse operator},
for all $\kappa\in\lattice\dual_{\delta}$
and all $k\in\Set{O}$ the spectrum of the operator
$\Op{D}(\psi_{-\kappa}v,\psi_{\kappa}W,k)$
is given by the polar set of the inverse of the $2\times 2$--matrices
$\unity+\left(\begin{smallmatrix}
0 & \pi(k_2+\sqrt{-1}k_1)\\
\pi(k_2-\sqrt{-1}k_1) & 0
\end{smallmatrix}\right)^{-1}\comp
\Mat{A}(\psi_{-\kappa}V,\psi_{\kappa}W,k)$.
Obviously this polar set is equal to the zero set of the determinant
$$\det\left(\left(\begin{smallmatrix}
0 & \pi(k_2+\sqrt{-1}k_1)\\
\pi(k_2-\sqrt{-1}k_1) & 0
\end{smallmatrix}\right)+
\Mat{A}(\psi_{-\kappa}V,\psi_{\kappa}W,k)\right).$$

Since the operator
$\Op{F}\comp
\triv{\Op{R}}_{\text{\scriptsize\rm red}}(0,0,k,0)
\comp\Op{F}^{-1}$
acts on $\ell_{\frac{p}{p-1}}(\lattice\dual)\times
\ell_{\frac{p}{p-1}}(\lattice\dual)$ as the multiplication with the
following sequence of matrices indexed by $\kappa\in\lattice\dual$:
$$\begin{cases}0 &\text{if }\kappa=0,\\
\displaystyle\frac{-1}{\pi}\begin{pmatrix}
0 & (k_2+\kappa_2-\sqrt{-1}(k_1+\kappa_1))^{-1}\\
(k_2+\kappa_2+\sqrt{-1}(k_1+\kappa_1))^{-1} & 0
\end{pmatrix} & \text{if }\kappa\neq 0,
\end{cases}$$
the partial derivatives of this operator with respect
to $k_1$ and $k_2$ acts as the multiplication with the partial
derivatives of this sequence of matrices, respectively.
The arguments in the proof of Lemma~\ref{weakly continuous resolvent}
concerning the operators
$\Op{F}\comp\Op{R}(0,0,k,\sqrt{-1}\lambda)\comp\Op{F}^{-1}$
imply that the partial derivatives of the operator--valued function
$k\mapsto\triv{\Op{R}}_{\text{\scriptsize\rm red}}(0,0,k,0)$
define for all $1<p<q<\infty$ and all small $k$ bounded operators from
$\banach{p}(\torus)\times\banach{p}(\torus)$ into 
$\banach{q}(\torus)\times\banach{q}(\torus)$.
Therefore, the arguments concerning the convergence of the matrices
$\Mat{A}(\psi_{-\kappa}V,\psi_{\kappa}W,k)$ in the limit
$g(\kappa,\kappa)\rightarrow\infty$ also imply
that the partial derivatives of the function
$k\mapsto\Mat{A}(\psi_{-\kappa}V,\psi_{\kappa}W,k)$
converge uniformly with respect to $k\in\Set{O}$ to zero in the
limit $g(\kappa,\kappa)\rightarrow\infty$.

Consequently for all $V,W\in \banach{2}(\torus)$ there exists
a $\delta>0$ and a sequence $k_{\kappa}$ of elements of $\Set{O}$
indexed by $\lattice\dual_{\delta}$
such that the values of the matrix--valued functions
$k\mapsto\left(\begin{smallmatrix}
0 & \pi(k_2+\sqrt{-1}k_1)\\
\pi(k_2-\sqrt{-1}k_1) & 0
\end{smallmatrix}\right)+\Mat{A}(\psi_{-\kappa}V,\psi_{\kappa}W,k)$
at $k=k_{\kappa}$ are diagonal matrices.
In fact, the bounds on the partial derivatives of the function
$k\mapsto\Mat{A}(\psi_{-\kappa}V,\psi_{\kappa}W,k)$ imply that
the mapping, which maps $k$ to the unique $k'\in\Set{O}$
such that
$\left(\begin{smallmatrix}
0 & \pi(k_2'+\sqrt{-1}k_1')\\
\pi(k_2'-\sqrt{-1}k_1') & 0
\end{smallmatrix}\right)+\Mat{A}(\psi_{-\kappa}V,\psi_{\kappa}W,k)$
is a diagonal $2\times 2$--matrix, is a contraction.
Therefore, due to the contraction mapping principle
\cite[Theorem~V.18]{RS1}, this mapping has a unique fixed point,
which yields the desired $k_{\kappa}$.
Let $\fourier{\mathrm{v}}(V,W,\kappa)$ and
$\fourier{\mathrm{w}}(V,W,-\kappa)$ denote the
diagonal elements of these matrices:
$$\Mat{A}(\psi_{-\kappa}V,\psi_{\kappa}W,k_{\kappa})=
\left(\begin{smallmatrix}
\fourier{\mathrm{v}}(V,W,\kappa) & 0\\
0 & \fourier{\mathrm{w}}(V,W,-\kappa)
\end{smallmatrix}\right).$$
We conclude that the intersections
$\fermi(\psi_{-\kappa}V,\psi_{\kappa}W)\cap\Set{O}$
are small perturbations of the subvarieties
$\left\{k\in\Set{O}\mid g(k-k_{\kappa},k-k_{\kappa})=
\fourier{\mathrm{v}}(V,W,\kappa)\fourier{\mathrm{w}}(V,W,-\kappa)\right\}$.
In fact, the bounds of the partial derivatives of the function
$k\mapsto\Mat{A}(\psi_{-\kappa}V,\psi_{\kappa}W,k)$
imply for all $\kappa\in\lattice\dual_{\delta}$ the estimates
$$\left\|\Mat{A}(\psi_{-\kappa}V,\psi_{\kappa}W,k)-
\left(\begin{smallmatrix}
\fourier{\mathrm{v}}(V,W,\kappa) & 0\\
0 & \fourier{\mathrm{w}}(V,W,-\kappa)
\end{smallmatrix}\right)\right\|\leq\varepsilon\left\|k-k_{\kappa}\right\|$$
uniformly for all $k\in\Set{O}$.

If the mapping $k\mapsto\|\Mat{A}(\psi_{-\kappa}V,\psi_{\kappa}W,k)\|$
and the corresponding partial derivatives of first order have on
the small ball $\Set{O}=B(k_{\kappa},\varepsilon')$ an upper bound
$\varepsilon$, which is small compared with $\varepsilon'$,
then this zero set defines a two--sheeted covering over a small ball
$\xx{p}\in B(g(\xx{\gamma},k_\kappa),\varepsilon'')$ with
$0<\varepsilon\ll\varepsilon''\ll\varepsilon'$.
In particular, on $B(g(\xx{\gamma},k_\kappa),\varepsilon'')$
this two--sheeted covering has either an ordinary double point
connecting the two sheets and no other zero of $d\xx{p}$,
or two simple branch point connecting the two sheets
and no other zero of $d\xx{p}$.
Again Lemma~\ref{covariant transformation} shows that this yields an
asymptotic description of the \Em{complex Fermi curves}

In case of a pair of potentials of the form $(U,\Bar{U})$,
the map $\eta$ transforms the matrices
$\Mat{A}(\psi_{-\kappa}U,\psi_{\kappa}\Bar{U},k)$ like
$\Bar{\Mat{A}}(\psi_{-\kappa}U,\psi_{\kappa}\Bar{U},k)=
\Op{J}\comp\Mat{A}(\psi_{-\kappa}U,\psi_{\kappa}\Bar{U},-\Bar{k})
\comp\Op{J}^{-1}$ (compare with Section~\ref{subsection reductions}).
Consequently the corresponding sequences $k_{\kappa}$
are fixed points of $\eta$,
and the sequences $\fourier{\mathrm{w}}(U,\Bar{U},-\kappa)$ are the complex
conjugate of the sequences $\fourier{\mathrm{v}}(U,\Bar{U},\kappa)$,
which we denote by $\fourier{\mathrm{u}}(U,\kappa)$.
The invariance of a \Em{complex Fermi curve} of potentials
of the form $(U,\Bar{U})$ under the involution $\eta$ implies
that in case of vanishing $\fourier{\mathrm{u}}(U,\kappa)$
this part of the \Em{complex Fermi curve}
gives rise to an ordinary double point and no cut
and in case of non--vanishing $\fourier{\mathrm{u}}(U,\kappa)$
to one \Em{horizontal cut}. On the other hand,
Theorem~\ref{asymptotic analysis 1} implies that the sets
$\Set{V}^{\pm}_{\varepsilon,\delta}$ are mapped by $\xx{p}$
biholomorphically onto open subsets of $\mathbb{C}^{\pm}_ {\xx{p}}$.
Finally, the continuity of the zeroes of $d\xx{p}$ under deformations
of $U\in \banach{2}(\torus)$ implies that the arithmetic genus
of $\fermi$ is preserved. More precisely,
if we count the singularities as zeroes of $d\xx{p}$, whose order is
equal to twice the
\De{Local contribution to the arithmetic genus}~\ref{local contribution},
then the total order of zeroes of $d\xx{p}$ is locally preserved and
equal to twice the number of elements of $\lattice\dual$.
Therefore, $\Set{V}^{\pm}_{\varepsilon,\delta}$
has no zeroes of $d\xx{p}$,
the small handles excluded from these sets
have two zeroes of $d\xx{p}$
(either as an ordinary double point or as two simple branch points)
and the number of zeroes of $d\xx{p}$ on the compact part
is twice the number of elements of
$\lattice\dual\setminus\lattice\dual_{\delta}$.
This proves the theorem for all elements
$\fermi\in\moduli_{\lattice,\eta,\willmore}$.

Due to Theorem~\ref{limits of resolvents} the accumulation points of
sequences in $\Bar{\moduli}_{\lattice,\eta,\willmore}$
are \Em{complex Fermi curves} of
\De{Finite rank Perturbations}~\ref{finite rank perturbations}
of Dirac operators $\Op{D}(U,\Bar{U},k)$
with complex potentials in
$\left\{U\in \banach{2}(\torus)\mid 4\|U\|_2^2\leq\willmore\right\}$.
In the eighth step of
Section~\ref{subsubsection finite rank perturbations}
we extended the asymptotic analysis of
Theorem~\ref{asymptotic analysis 1} to these
\De{Finite rank Perturbations}~\ref{finite rank perturbations}.
Due to the
{Strong unique continuation property}~\ref{strong unique continuation}
we may use the restrictions (compare with Remark~\ref{restriction})
of the corresponding resolvents to the complement of small balls
around the singular points and the corresponding norms.
Therefore, also the \Em{complex Fermi curves} of
\De{Finite rank Perturbations}~\ref{finite rank perturbations}
have a representation of the prescribed form.
\end{proof}

Due to this theorem all \Em{complex Fermi curves} with
dividing real parts have finite genus.
Moreover, all $\fermi\in\Bar{\moduli}_{\lattice,\eta,\sigma,\willmore}$
with dividing real parts have a small open neighbourhood $\Set{U}$
with respect to the subspace topology of the closed subsets of
$\overline{\mathbb{C}^2}$ (compare with Lemma~\ref{compact metric})
such that the \Em{geometric genera} of all \Em{complex Fermi curves}
in $\Set{U}$ with dividing real parts are bounded.
In fact, with the exception of finitely many indexes $\kappa$
the corresponding handles of a \Em{complex Fermi curve}
with dividing real part have ordinary double points
of the form $(y,\eta(y))$. Furthermore,
a small bound on the deformations of the functions
$k_1\pm\sqrt{-1}k_2$ on the compact part guarantees that
these ordinary double points are deformed into
ordinary double points of the same form,
since a deformation into \Em{horizontal cuts}
contradicts the dividing property of the real part.
On the other hand, the subspace of \Em{complex Fermi curves} with
dividing real parts is closed, since a small deformation of a cut
connecting the two planes $\mathbb{C}^{\pm}_{\xx{p}}$
disjoint from the real part yields a cut with the same properties.
Hence the compactness of $\Bar{\moduli}_{\lattice,\eta,\willmore}$ implies

\begin{Corollary}\label{dividing genus bound}
For all $\willmore>0$ there exists a natural number
$g_{\max}\in\mathbb{N}$
depending only on the lattice $\lattice$ and $\willmore$
such that all \Em{complex Fermi curves}
$\fermi\in\Bar{\moduli}_{\lattice,\eta,\sigma,\willmore}$
with dividing real parts are contained in
$\Bar{\moduli}_{g_{\max},\lattice,\eta,\sigma,\willmore}$.
\qed
\end{Corollary}

\begin{Remark}\label{not uniform bound}
These numbers $g_{\max}$ are not bounded
uniformly for all lattices and become arbitrary large,
whenever $\vol(\torus)\min\{g(\kappa,\kappa)\mid
\kappa\in\lattice\dual\setminus\{0\}\}$ becomes arbitrary small.
In fact, for all $\kappa\in\lattice\dual\setminus\{0\}$ with
$\pi^2\vol(\torus)g(\kappa,\kappa)\leq\willmore$
the \Em{complex Fermi curves} of the constant potentials have
one real double point.
\end{Remark}

Proposition~\ref{asymptotic analysis 2}
suggests the following characterization of the elements of
$\Bar{\moduli}_{\lattice,\eta,\willmore}$.
Instead of the data $(\Spa{Y},\infty^-,\infty^+,k)$ we consider all data
$(\Spa{Y},k)$ of complex spaces $\Spa{Y}$ together with a
holomorphic multi--valued function $k$ on $\Spa{Y}$,
which satisfies conditions \Em{Quasi--momenta}~(ii)--(iii)
of Section~\ref{subsection complex Fermi curves of finite genus}
and the following conditions \Em{Quasi--momenta}~(i') and
\Em{Quasi--momenta}~(vi).
\begin{description}
\index{condition!quasi--momenta (i') and (vi)|(}
\index{quasi--momenta!condition $\sim$ (i') and (vi)|(}
\item[Quasi--momenta (i')]
  For both components $\xx{p}$ and $\yy{p}$ of
  $k$ the Riemann surface $\Spa{Y}$ may be realized by gluing
  two planes $\mathbb{C}_{\xx{p}}^{\pm}$
  (or $\mathbb{C}_{\yy{p}}^{\pm}$) along
  a countable  combination of \Em{horizontal cuts} and
  \Em{pairs of parallel cuts} as described in
  Lemma~\ref{gluing rule}. Moreover, there exists some branch of the
  function $k$ with the property that the function $k_1+\sqrt{-1}k_2$
  on $\mathbb{C}_{\xx{p}}^+$ converges to zero,
  whenever the function $\xx{p}$ tends to infinity.
\item[Quasi--momenta (vi)]
  The integral of the two-form
  $$2\pi\sqrt{-1}\vol(\torus)
  \left(dk_1+\sqrt{-1}dk_2\right)\wedge
  \left(\overline{dk_1+\sqrt{-1}dk_2}\right)$$
  over $\mathbb{C}_{\xx{p}}^+$ is less or equal to $\willmore$.
\index{condition!quasi--momenta (i') and (vi)|)}
\index{quasi--momenta!condition $\sim$ (i') and (vi)|)}
\end{description}
Due to Proposition~\ref{asymptotic analysis 2} the elements of
$\Bar{\moduli}_{\lattice,\eta}$ obey these conditions.
We conjecture that these conditions characterize the elements of
$\Bar{\moduli}_{\lattice,\eta,\willmore}$ for all $\willmore>0$.

Let us remark, that the mapping
$$(V,W)\mapsto\left(
\fourier{\mathrm{v}}(V,W,\kappa),\fourier{\mathrm{w}}(V,W,\kappa)
\right)_{\kappa\in\lattice\dual_{\delta}}$$
is a nonlinear perturbation of the Fourier transform.
For small $\|V\|$ and $\|W\|$ the sequences
$\fourier{\mathrm{v}}(V,W,\kappa)$
and $\fourier{\mathrm{w}}(V,W,\kappa)$ are defined for all
$\kappa\in\lattice\dual$. 

\begin{Remark}~\label{overlapping bound}
In order to define all entries of these sequences,
we have to ensure that the handles connecting the parts near
$k^+_{\kappa}$ with the part near $k^-_{\kappa}$ of the
\Em{complex Fermi curve} do not overlap.
One can show that an upper bound on $4\|V\|_2^2$ and $4\|W\|_2^2$,
which prevents such an overlapping essentially is given by
$\pi^2\vol(\torus)\min\left\{g(\kappa,\kappa)\mid
\kappa\in\lattice\dual\setminus\{0\}\right\}.$
Moreover, there does not exist a non--vanishing upper bound,
which prevents this overlapping uniformly for all lattices $\lattice$
(compare with Remark~\ref{not uniform bound}).
\end{Remark}

For the zero potentials the functions $\Mat{A}(0,0,k)$
vanishes identically.
Furthermore, the derivative of the mapping
$(V,W)\mapsto\left(
\fourier{\mathrm{v}}(V,W,\kappa),\fourier{\mathrm{w}}(V,W,\kappa)
\right)_{\kappa\in\lattice\dual}$ at the zero potentials is equal to
the Fourier transform.
An easy calculation shows that the \Em{complex Fermi curves}
of pairs of potentials of the form
$(u\psi_{\kappa},\Bar{u}\psi_{-\kappa})$
with $u\in\mathbb{C}$ are equal to
$$\fermi(u\psi_{\kappa},\Bar{u}\psi_{-\kappa})=
\bigcup\limits_{\kappa'\in\lattice\dual}
\left\{k\in\mathbb{C}^2\mid 
\pi^2g(k+k^+_{\kappa}+\kappa',k+k^+_{\kappa}+\kappa')=u\Bar{u}\right\}.$$
Moreover, with the help of the arguments of the proof of
Lemma~\ref{singularity} it is easy to see that if
$4u\Bar{u}\leq \min\{g(\kappa,\kappa)\mid
\kappa\in\lattice\dual\setminus\{0\}\}$
($\iff\willmore\leq\pi^2\vol(\torus)
\min\{g(\kappa,\kappa)\mid
\kappa\in\lattice\dual\setminus\{0\}\}$)
these pairs of potentials $(u\psi_{\kappa},\Bar{u}\psi_{-\kappa})$
are all potentials of the form $(U,\Bar{U})$,
whose \Em{complex Fermi curves} have \Em{geometric genus} zero.
This leads to a local description of the mapping
\begin{align*}
\left\{U\in \banach{2}(\torus)\mid 
4\|U\|_2^2\leq\willmore\right\}
&\rightarrow\moduli_{\lattice,\eta,\willmore},
&U&\mapsto\fermi(U,\Bar{U}).
\end{align*}
We shall see that for small $\willmore$ this mapping
is a nonlinear perturbation of the mapping
\begin{align*}
\banach{2}(\torus)&\rightarrow
\ell^+_1(\lattice\dual)=
\left\{\parameter{t}\in
\ell_1(\lattice\dual)\mid\parameter{t}(\kappa)\in\mathbb{R}^+_0
\;\forall\kappa\in\lattice\dual\right\},
&U&\mapsto\parameter{t},\text{ with }
\parameter{t}(\kappa)=|\fourier{U}(\kappa)|^2.
\end{align*}
Here $\fourier{U}$ denotes the Fourier transform of $U$.
Moreover, together with the symplectic form, which is equal to the
imaginary part of the Hermitian form of $\banach{2}(\torus)$,
this mapping defines an infinite--dimensional
completely integrable system. This is a general feature of
integrable systems related to Lax operators.
The mappings from the potentials to the spectral curves
(i\ e.\ a complete sets of integrals)
of these integrable systems are nonlinear perturbations
of the mapping from the potentials to the
absolute values of their Fourier coefficients,
which together with the corresponding symplectic form
represents simple examples of
infinite--dimensional integrable systems.
The preimages of the latter map are compact tori,
whose dimensions are equal to the number of non--vanishing moduli.
The preimages of the former map are the \Em{isospectral sets},
which in general may be identified with open sets of the
corresponding Jacobian varieties
(compare with Lemma~\ref{existence of potentials}).
Therefore, these dimensions are equal to the genera of the
corresponding \Em{complex Fermi curves}.
Due to a theorem of Littlewood
(\cite{Li,FR} and \cite[Chapter~3 \S1.6.]{Kis})
almost all elements of preimages of the latter map have finite
$\banach{q}$--norm for all finite $q$.
This leads to the conjecture that all
\Em{complex Fermi curves} $\fermi(U,\Bar{U})$
are equal to the \Em{complex Fermi curve} of
some potential $U'\in\banach{q}(\torus)$, and that the
$\banach{q}$--norm of $U'$ is bounded by some constant depending on
$q$ and the maximum of all $\varepsilon>0$
such that $\|U\|_{\varepsilon,2}$ is smaller than $S_p^{-1}$
(compare with Lemma~\ref{weakly continuous resolvent}).

A precise version of this perturbation is that the moduli space
may be realized by gluing homeomorphically open subsets of
$\ell^+_1(\lattice\dual)$.
Moreover, these parameterizations of the moduli space has the property,
that the genus of the corresponding \Em{complex Fermi curves}
is equal to the number of non--vanishing moduli.
We remark that these local parameterizations of the moduli space
endows it with a finer topology than the topology introduced
in Section~\ref{subsection compactified bounded genus}.
In contrast to the former topology
with respect to which the mapping
$U\mapsto\fermi(U,\Bar{U})$ is `almost' weakly continuous,
the latter topology is the finest topology, such that this mapping is
continuous with respect to the usual $\banach{2}$--topology.
Let us call a topology of the moduli space with the property,
that each element of the moduli space has an open neighbourhood,
which is homeomorphic to an open subset of $\ell^+_1(\lattice\dual)$,
\index{structure!$\ell^+_1(\lattice\dual)$--$\sim$}
a $\ell^+_1(\lattice\dual)$--structure.

\begin{Theorem}\label{l1-structure}
The moduli space $\moduli_{\lattice}=
\left\{\fermi(V,W)\mid V,W\in\banach{2}(\torus)\right\}$
\index{moduli space!$\moduli_{\lattice}$}
is a holomorphic $\ell_1(\lattice\dual)$--manifold.
Furthermore, locally there exists a $\delta>0$,
such that the moduli with indexes
$\kappa\in\lattice\dual_{\delta}$ are unique.
Moreover, the compactified moduli space
$\Bar{\moduli}_{\lattice,\eta}$ has an
asymptotic $\ell^+_1(\lattice\dual)$--structure
(i.\ e.\ it is a closed subset of a $\ell_1(\lattice\dual)$--manifold,
and the restriction to the moduli with indexes
$\kappa\in\lattice\dual_{\delta}$ with a suitable $\delta$ yields an
open mapping into $\ell^+_1(\lattice\dual_{\delta})$).
Again these moduli are asymptotically unique.
Finally, the \Em{complex Fermi curves} of finite genus are
parameterized by elements of $\ell^+_1(\lattice\dual)$
with finitely many non--vanishing entries.
In particular, they are dense with respect to
the corresponding topologies.
\end{Theorem}

\begin{proof}
The proof of this theorem is a consequence
of the improved asymptotic analysis
of Proposition~\ref{asymptotic analysis 2}.
Again we first treat the subspace
$\moduli_{\lattice,\eta}$
and then in a second step we extend our arguments to
$\Bar{\moduli}_{\lattice,\eta}$.

\noindent
{\bf 1.} First we define the moduli $\parameter{t}(\kappa)$ for all
all \Em{complex Fermi curves} of pairs of potentials
of the form $(U,\Bar{U})$ and $\kappa\in\lattice\dual_{\delta}$,
with an appropriate $\delta$ (depending on $U$)
defined in the proof of Proposition~\ref{asymptotic analysis 2}.
$$\parameter{t}(\kappa)=
\pi\oint\limits_{\text{\scriptsize \Em{horizontal cut} of
                $\mathbb{C}^+_{\xx{p}}$ with index $\kappa$.}}
k_1dk_2.$$
Due to Lemma~\ref{gluing rule}
these \Em{horizontal cuts} are invariant under $\eta$.
Therefore, the integrals of $dk$ along these cycles vanish.
Also the integral in the definition of $\parameter{t}(\kappa)$
does not depend on the branch of $k$.
Moreover, due to the invariance under $\eta$, the well defined moduli
$\parameter{t}(\kappa)$ are real.
Since the form $k_1dk_2$ is closed, we may integrate it
along any closed cycles of an open set of the form
$\Set{V}^+_{\varepsilon,\delta}$
with an appropriate $\varepsilon>0$ and $\delta>0$ 
around the excluded disc near $k^+_{\kappa}$.
This extends the definition of the moduli to all pairs of
potentials with $V,W\in\banach(\torus)$.

Due to Lemma~\ref{Willmore functional} the contribution of
the \Em{horizontal cut} with index $\kappa$
to the \Em{first integral} is equal to
$4\vol(\torus)\parameter{t}(\kappa)$.
Consequently, in case that all these handles do not overlap
and all moduli $\parameter{t}(\kappa)$ are well defined,
the \Em{first integral} $4\int\limits_{\torus}$ is equal to
$4\vol(\torus)\sum\limits_{\kappa\in\lattice\dual}\parameter{t}(\kappa)$,
whenever this sum converges.

\noindent
{\bf 2.} Now we claim that for pairs of potentials of the form
$(U,\Bar{U})$ the well--defined $\parameter{t}(\kappa)$ are positive,
which imply that they belong to $\ell^+_1(\lattice\dual_{\delta})$.
The real part of $k$ defines a mapping from the
$\fermi(U,\Bar{U})$'s into $\mathbb{R}^2$
(which may depend on the branch of $k$).
An easy calculation shows that this mapping
is locally a homeomorphism, if and only if $d\yy{p}/d\xx{p}$ is not real.
Due to Proposition~\ref{asymptotic analysis 2}
asymptotically the function $d\yy{p}/d\xx{p}$
is real only on small disjoint cycles,
which are homolog to the \Em{horizontal cuts} indexed by
$\lattice\dual_{\delta}$.
Moreover, all these cycles are invariant under $\eta$ and on each of
these cycles the function $d\yy{p}/d\xx{p}$ takes all values of
$\mathbb{P}^1_\mathbb{R}$ exactly twice.
Since $\eta$ does not have fixed points, this implies that these
cycles are mapped by the real part of $k$ onto smooth Jordan curves
in $\torus\dual$ around some element of the form $\kappa/2$
(which are the images of fixed points of $\eta$ under the real part of $k$).
If we cut the complement of the compact part of the
\Em{complex Fermi curve} described in
Theorem~\ref{asymptotic analysis 1}
along these cycles into two pieces,
then the real part of the unique branch of $k$,
whose linear combinations $k_1\pm\sqrt{-1}k_2$ become at
$\infty^{\pm}$ arbitrary small, respectively,
maps these two components homeomorphically onto the complement
of some compact part of $\mathbb{R}^2$
and those unique discs around $\kappa/2$ and $-\kappa/2$,
whose boundary is the Jordan curve with index $\kappa$, respectively.
Finally, they are glued along these Jordan curves by the map
$k\mapsto k-\kappa$. The arguments of the proof Lemma~\ref{gluing rule}
show that $\parameter{t}(\kappa)$ is equal to $\pi$ times the volume
of the corresponding disc in $\mathbb{R}^2$,
whose boundary is the image of the cycle with index $\kappa$
under the real part of $k$, plus or minus $\pi$ times the volume
of the analogous disc,
whose boundary is the images of this cycle under the imaginary
part of $k$. Furthermore, these arguments also show that
$\parameter{t}(\kappa)$ is equal to $\pi$ times the volume
of the unique disc in $\mathbb{R}^2$, whose boundary is the image of
the \Em{horizontal cut} with index $\kappa$
under the real part of $k$.
But locally in a small tubular neighbourhood
of the unique cycle with index $\kappa$,
on which $d\yy{p}/d\xx{p}$ takes real values,
the real part of $k$ takes values in the complement of the disc,
whose boundary is the image of this cycle.
We conclude that $\parameter{t}(\kappa)$ is larger or equal to the
volume of this disc and therefore non--negative.
Finally, $\parameter{t}(\kappa)$ vanishes if and only if
the corresponding part of the \Em{complex Fermi curve} has an
ordinary double point and does not have a \Em{horizontal cut}.
To sum up, we have proven that the mapping $U\mapsto\parameter{t}$
maps $\banach{2}(\torus)$ into $\ell^+_1(\lattice\dual_{\delta})$
(with some suitable $\delta$ depending on $U$), and
$4\vol(\torus)\|\parameter{t}\|_1\leq 4\|U\|_2^2$.

These arguments show that image under the imaginary part of $k$
of those unique cycles, on which $d\yy{p}/d\xx{p}$ takes real values,
have negative `volume'. Since the imaginary part of $k$ is preserved
under $\eta$, the imaginary part of $k$ maps these cycles not
injectively into $\mathbb{R}^2$. In fact, for the
constant potentials the unique non--trivial cycle is real and the
imaginary part of $k$ maps this cycle onto the origin.
In general the preimages contain at least two points
of the \Em{complex Fermi curves},
which are interchanged by $\eta$. Therefore, the images of these cycles
under the imaginary part of $k$ look like smooth curves

\noindent
\begin{minipage}[t]{11cm}
  \begin{align*}
  S^1&\rightarrow\mathbb{R}^2&
  \exp(2\pi\sqrt{-1}\varphi)&\mapsto(x,y),
  \end{align*}
  with the property that the quotient $dy/dx=(dy/d\varphi)/(dx/d\varphi)$
  takes each value in $\mathbb{P}^1_\mathbb{R}$ exactly once.
  Such curves cannot be immersed
  (i.\ e.\ they must have common zeroes of $dy/d\varphi$ and $dx/d\varphi$),
  and the exterior parts of these curves are locally convex in the
  complement of the non--immersed points. This yields another explanation
  why the corresponding contribution to the \Em{first integral} is
  negative. The figure on the right gives the image of a non--trivial example.
  \end{minipage}\hfill
\setlength{\unitlength}{1mm}
\begin{picture}(45,-40)
\qbezier(0,0)(22.5,-13)(45,0)
\qbezier(0,0)(22.5,-13)(22.5,-39)
\qbezier(45,0)(22.5,-13)(22.5,-39)
\end{picture}

\noindent
{\bf 3.} Now we shall show that the mapping
$(V,W)\mapsto\parameter{t}$
is an open mapping from $\banach{2}(\torus)\times\banach{2}(\torus)$ into
$\ell_1(\lattice\dual_{\delta})$
with suitable $\delta$ depending on $V$ and $W$.
For this purpose we shall first prove that the mapping
$(V,W)\mapsto(\fourier{\mathrm{v}}(V,W,\cdot),\fourier{\mathrm{w}}(V,W,\cdot))$
is an open mapping from $\banach{2}(\torus)\times\banach{2}(\torus)$ into
$\ell_2(\lattice\dual_{\delta})\times\ell_2(\lattice\dual_{\delta}\times)$
with suitable $\delta$ depending on $V$ and $W$.

\begin{Lemma}\label{continuous derivative}
The derivatives of the mapping $(V,W)\mapsto
\left(\fourier{\mathrm{v}}(V,W,\cdot),\fourier{\mathrm{w}}(V,W,\cdot)\right)$
are bounded operators from
$\banach{2}(\torus)\times\banach{2}(\torus)$
into $\ell_2(\lattice\dual_{\delta})\times\ell_2(\lattice\dual_{\delta})$
(with suitable $\delta$ depending on $(V,W)$)
\end{Lemma}

\begin{proof}
First we investigate for the same $k\in\Set{O}$ as in the proof of
Proposition~\ref{asymptotic analysis 2} in three steps
the derivative of the mapping
\begin{align*}
(V,W)&\mapsto\Mat{A}_{V,W,k}&\text{with }
\Mat{A}_{V,W,k}(x)=\sum\limits_{\kappa\in\lattice\dual_{\delta}}
\psi_{\kappa}(x)
\Mat{A}(\psi_{-\kappa}V,\psi_{\kappa}W,k)\;\forall x\in\torus
\end{align*}
In doing so we we have to estimate the Fourier transform
$\Mat{A}_{V,W,k}$ of the sequence of matrices
$\Mat{A}(\psi_{\kappa}V,\psi_{\kappa}W,k)$ indexed by
$\kappa\in\lattice\dual_{\delta}$.
For this purpose we make use of the Bochner $q$--integrable functions
from $\torus$ into some Banach space $\Spa{E}$.
They are denoted by $\banach{q}(\torus)\Tilde{\otimes}_{\Delta_q}\Spa{E}$
\cite[Section~7.]{DF}.
More precisely, we shall consider sequences of elements of
$\banach{q}(\torus)$ indexed by $\kappa\in\lattice\dual$
(or $\lattice\dual_{\delta}$) as elements of
$\banach{q}(\torus)\Tilde{\otimes}_{\Delta_q}\ell_{\cdot}(\lattice\dual)$.
Consequently, the Fourier transform $\unity\otimes\Op{F}$ maps such
tensor products into the tensor products
$\banach{q}(\torus)\Tilde{\otimes}_{\Delta_q}\banach{\cdot}(\torus)$.
The derivative
$\frac{\var\Mat{A}(\psi_{-\kappa}V,\psi_{\kappa}W,k)}{(\var V,\var W)}$
is the restriction of the following operator to the constant functions:
\begin{multline*}
\left(\unity-
\begin{pmatrix}
\psi_{-\kappa}V & 0\\
0 & \psi_{\kappa}W
\end{pmatrix}\comp
\triv{\Op{R}}_{\text{\scriptsize\rm red}}(0,0,k,0)
\right)^{-1}\comp\begin{pmatrix}
\psi_{-\kappa}\var V & 0\\
0 & \psi_{\kappa}\var W
\end{pmatrix}\comp\\
\comp\left(\unity-
\triv{\Op{R}}_{\text{\scriptsize\rm red}}(0,0,k,0)
\comp\begin{pmatrix}
\psi_{-\kappa}V & 0\\
0 & \psi_{\kappa}W
\end{pmatrix}\right)^{-1}.
\end{multline*}

\noindent
{\bf 1.} In a first step we show that for all $k\in\Set{O}$
and for all $1<p<2$ and $1/p+1/q=1$ both entries of the functions
$$(x,x')\mapsto\sum\limits_{\kappa\in\lattice\dual}
\Tilde{\Op{R}}_{\text{\scriptsize\rm red}}(0,0,k,0)\comp\begin{pmatrix}
\psi_{-\kappa}V & 0\\
0 & \psi_{\kappa}W
\end{pmatrix}\begin{pmatrix}
\psi_0\\
\psi_0
\end{pmatrix}(x)\psi_{\kappa}(x')$$
belong to the Banach space
$\banach{q}(\torus)\Tilde{\otimes}_{\Delta_q}\banach{p}(\torus)$.
We view the elements of this Banach space as measurable functions
$f$ on $\torus\times\torus$ with norm
$$\|f\|_{q,p}=\left(\int\limits_{\torus}\left(\int\limits_{\torus}
|f(x,x')|^pd^2x'\right)^{\frac{q}{p}}d^2x\right)^{1/q}.$$
The off--diagonal entries of the resolvents
$\Tilde{\Op{R}}_{\text{\scriptsize\rm red}}(0,0,k,0)$
are convolutions with functions in the \Em{Lorentz spaces}
$\banach{2,\infty}(\torus)$
(compare with Lemma~\ref{resolvent integral kernel})
and the diagonal entries vanish. Moreover, the sum
$\sum\limits_{\kappa\in\lattice\dual} \psi_{\kappa}(x)$ is equal to
$\vol^{1/2}(\torus)$ times the Dirac's $\delta$--function on $\torus$.
Therefore, we have to show that for $f\in\banach{2,\infty}(\torus)$ and
$g\in\banach{2}(\torus)$ the function $(x,x')\mapsto f(x-x')g(x')$
belongs to $\banach{q}(\torus)\Tilde{\otimes}_{\Delta_q}\banach{p}(\torus)$.
The distribution function $\mu_{|f|^p}$
and the non--increasing rearrangement $\left(|f|^p\right)^{\ast}$ of
the function $|f|^p$ may be expressed in terms of the
distribution function $\mu_f$ and the non--increasing rearrangement
$f^{\ast}$ of $f$:
\begin{align*}
\mu_{|f|^p}(\lambda)&=\mu_f(\lambda^p)&
\left(|f|^p\right)^{\ast}(t)&=\left(f^{\ast}(t)\right)^p.
\end{align*}
Hence the functions $|f|^p$ and $|g|^p$ belong to the \Em{Lorentz space}
$\banach{\frac{2}{p},\infty}(\torus)$ and $\banach{\frac{2}{p}}(\torus)$
and the corresponding norms are equal to
\begin{align*}
\left\| |f|^p\right\|_{\frac{2}{p},\infty}&=\left\|f\right\|_{2,\infty}^p
\text{ and}
& \left\| |g|^p\right\|_{\frac{2}{p}}&=\left\|g\right\|_2^p
\text{, respectively.}\end{align*}
Due to the \De{Generalized Young's inequality}~\ref{generalized young}
the $\banach{\frac{q}{p}}$--norm of the convolution of $|f|^p$ with
$|g|^p$ is bounded by some constant times
$\left\| |f|^p\right\|_{\frac{2}{p},\infty}\cdot
\left\| |g|^p\right\|_{\frac{2}{p}}$.
We remark that the relation $1/q+1/p=1$ implies $p/2+p/2=p/q+1$.
To sum up, the
$\banach{q}(\torus)\Tilde{\otimes}_{\Delta_q}\banach{p}(\torus)$--norm
of $(x,x')\mapsto f(x-x')g(x')$
is bounded by the corresponding constant of the
\De{Generalized Young's inequality}~\ref{generalized young}
times $\|f\|_{2,\infty}\cdot\|g\|_2$.

\noindent
{\bf 2.} Secondly we show that for the same $p,q$ and
the same $k\in\Set{O}$ as before
all entries of the Fourier transform of the sequence of operators
$$\Tilde{\Op{R}}_{\text{\scriptsize\rm red}}(0,0,k,0)\comp\begin{pmatrix}
\psi_{-\kappa}V & 0\\
0 & \psi_{\kappa}W
\end{pmatrix}\text{ and }
\left(\unity-
\Tilde{\Op{R}}_{\text{\scriptsize\rm red}}(0,0,k,0)\comp\begin{pmatrix}
\psi_{-\kappa}V & 0\\
0 & \psi_{\kappa}W
\end{pmatrix}\right)^{-1}$$
are bounded operators on
$\banach{q}(\torus)\Tilde{\otimes}_{\Delta_q}\banach{p}(\torus)$.
Due to the Hahn Decomposition Theorem \cite[Chapter~11 Section~5]{Ro2}
all $\banach{p}$--functions may be decomposed into a finite complex
linear combination of positive $\banach{p}$--functions.
Therefore, due to \cite[Section~7.3. Theorem]{DF} the multiplication
with $V$ and $W$ define bounded operators from
$\banach{q}(\torus)\Tilde{\otimes}_{\Delta_q}\banach{p}(\torus)$ to
$\banach{\frac{2q}{q+2}}(\torus)\Tilde{\otimes}_{\Delta_\frac{2q}{q+2}}
\banach{p}(\torus)$. Multiplication by
$\psi_{\pm\kappa}$ induces on these tensor products the operators
$f(x,x')\mapsto g(x,x')$ with $g(x,x')=f(x,x'\pm x)$,
which are obviously isometries. Finally, the entries of
$\Tilde{\Op{R}}_{\text{\scriptsize\rm red}}(0,0,k,0)$ are
convolution operators with $\banach{2,\infty}$--functions and therefore,
due to the Hahn Decomposition Theorem \cite[Chapter~11 Section~5]{Ro2},
regular in the sense of \cite[Section~7.3.]{DF},
which again due to \cite[Section~7.3. Theorem]{DF}
induce bounded operators from
$\banach{\frac{2q}{q+2}}(\torus)\Tilde{\otimes}_{\Delta_\frac{2q}{q+2}}
\banach{p}(\torus)$ into
$\banach{q}(\torus)\Tilde{\otimes}_{\Delta_q}\banach{p}(\torus)$.
This shows the claim for the first operator.

In the proof of Proposition~\ref{asymptotic analysis 2}
we mentioned that the restrictions of the sequence of operators
$\left(\unity-\triv{\Op{R}}_{\text{\scriptsize\rm red}}(0,0,k,\lambda)
\comp\left(\begin{smallmatrix}
\psi_{-\kappa}V & 0\\
0 & \psi_{\kappa}W
\end{smallmatrix}\right)\right)^{-1}\comp
\triv{\Op{R}}_{\text{\scriptsize\rm red}}(0,0,k,\lambda)$
to the orthogonal complement of the constant functions
are the resolvents of the restrictions of
$\triv{\Op{D}}(\psi_{-\kappa}V,\psi_{\kappa}W,k)$ to
this orthogonal complement
(compare with Remark~\ref{orthogonal complement}). This generalizes to
the Fourier transform of these sequences of operators indexed by
$\kappa\in\lattice\dual_{\delta}$ acting on the $\banach{p}$--function
on $\torus$ with values in $\banach{q}(\torus)$.
Hence the restriction of the Fourier transform of the former sequence
to the orthogonal complement of
$1\otimes\banach{p}(\torus)\subset
\banach{q}(\torus)\Tilde{\otimes}_{\Delta_q}\banach{p}(\torus)$
is the resolvent of the restriction of the Fourier transform
of the latter sequence.
Now the methods of Section~\ref{subsection resolvent} carry over to
these resolvents and imply that for all $V,W\in\banach{2}(\torus)$
there exists a $\delta>0$, such that all entries
of the Fourier transform of the sequence of operators
\begin{align*}
\left(\unity-\triv{\Op{R}}_{\text{\scriptsize\rm red}}(0,0,k,0)
\comp\begin{pmatrix}
\psi_{-\kappa}V & 0\\
0 & \psi_{\kappa}W
\end{pmatrix}\right)^{-1}&\text{ and}&
\left(\unity-\triv{\Op{R}}_{\text{\scriptsize\rm red}}^{t}(0,0,k,0)
\comp\begin{pmatrix}
\psi_{-\kappa}V & 0\\
0 & \psi_{\kappa}W
\end{pmatrix}\right)^{-1}&
\end{align*}
are bounded operator on
$\banach{q}(\torus)\Tilde{\otimes}_{\Delta_q}\banach{p}(\torus)$.

\noindent
{\bf 3.} In the third step we prove that
for the same $k\in\Set{O}$ as before
the $\banach{2}$--norms of the diagonal entries
and for all $1<p<2$ the $\banach{p}$--norms
of the off-- diagonal entries of $\var\Mat{A}_{V,W,k}$
may be estimated in terms of $\|\var V\|_2$ and $\|\var W\|_2$.
These entries of $\var\Mat{A}_{V,W,k}$ are integrals of
pointwise products of the sequences $\psi_{-\kappa}\var V$ and
$\psi_{\kappa}\var W$ with the entries of the sequence of functions
\begin{eqnarray*}
\left(\unity-\triv{\Op{R}}_{\text{\scriptsize\rm red}}^{t}(0,0,k,0)
\comp\begin{pmatrix}
\psi_{-\kappa}V & 0\\
0 & \psi_{\kappa}W
\end{pmatrix}\right)^{-1}\begin{pmatrix}
\psi_0\\
\psi_0
\end{pmatrix}=&\begin{pmatrix}
\psi_0\\
\psi_0
\end{pmatrix}\\
+\left(\unity-\triv{\Op{R}}_{\text{\scriptsize\rm red}}^{t}(0,0,k,0)
\comp\begin{pmatrix}
\psi_{-\kappa}V & 0\\
0 & \psi_{\kappa}W
\end{pmatrix}\right)^{-1}\comp
\triv{\Op{R}}_{\text{\scriptsize\rm red}}^{t}(0,0,k,0)
\comp\begin{pmatrix}
\psi_{-\kappa}V & 0\\
0 & \psi_{\kappa}W
\end{pmatrix}&\begin{pmatrix}
\psi_0\\
\psi_0
\end{pmatrix}
\end{eqnarray*}
and the entries of the sequence of functions
\begin{eqnarray*}
\left(\unity-\triv{\Op{R}}_{\text{\scriptsize\rm red}}(0,0,k,0)
\comp\begin{pmatrix}
\psi_{-\kappa}V & 0\\
0 & \psi_{\kappa}W
\end{pmatrix}\right)^{-1}\begin{pmatrix}
\psi_0\\
\psi_0
\end{pmatrix}=&\begin{pmatrix}
\psi_0\\
\psi_0
\end{pmatrix}\\
+\left(\unity-\triv{\Op{R}}_{\text{\scriptsize\rm red}}(0,0,k,0)
\comp\begin{pmatrix}
\psi_{-\kappa}V & 0\\
0 & \psi_{\kappa}W
\end{pmatrix}\right)^{-1}\comp
\triv{\Op{R}}_{\text{\scriptsize\rm red}}(0,0,k,0)
\comp\begin{pmatrix}
\psi_{-\kappa}V & 0\\
0 & \psi_{\kappa}W
\end{pmatrix}&\begin{pmatrix}
\psi_0\\
\psi_0
\end{pmatrix},
\end{eqnarray*}
respectively.
Due to the first two steps the Fourier transforms of the second terms
of the latter two sequences of functions belong to the tensor products
$\banach{q}(\torus)\Tilde{\otimes}_{\Delta_q}\banach{p}(\torus)$.
Since the Fourier transform of pointwise multiplication
is the convolution \cite[Theorem~IX.3]{RS2},
we shall estimate the mapping
\begin{multline*}
\banach{q_1}(\torus)\Tilde{\otimes}_{\Delta_{q_1}}\banach{p_1}(\torus)\times 
\banach{q_2}(\torus)\Tilde{\otimes}_{\Delta_{q_2}}\banach{p_2}(\torus)
\rightarrow
\banach{q_3}(\torus)\Tilde{\otimes}_{\Delta_{q_3}}\banach{p_3}(\torus),\\
(f,g)\mapsto h\text{ with }
h(x,x')=\int\limits_{\torus}f(x,x'-x'')g(x,x'')d^2x''.
\end{multline*}
Due to Young's inequality \cite[Section~IX.4]{RS2}
for $1/p_1+1/p_2=1/p_3+1$ the convolution is a bounded bilinear mapping
from the values of the $\banach{q_1}$--function $f$
and the values of the $\banach{q_2}$--function $g$
into the values of the $\banach{q_3}$--function $h$.
Therefore, H\"older's inequality \cite[Theorem~III.1~(c)]{RS1}
yields the estimate
$$\left\|h\right\|_{q_3,p_3}\leq
\left\|f\right\|_{q_1,p_1}\left\|g\right\|_{q_2,p_2}$$
for $1/q_1+1/q_2=1/q_3$.
We remark that the relations $1/p_1+1/p_2=1/p_3+1$ and
$1/q_1+1/q_2=1/q_3$ are compatible with the relations
$1/p_1+1/q_1=1/p_2+1/q_2=1/p_3+1/q_3=1$.
Hence for all $1<p<\infty$ and $1/p+1/q=1$ the Fourier transform
of the pointwise multiplication of the entries of
$$\left(\unity-\triv{\Op{R}}_{\text{\scriptsize\rm red}}^{t}(0,0,k,0)
\comp\begin{pmatrix}
\psi_{-\kappa}V & 0\\
0 & \psi_{\kappa}W
\end{pmatrix}\right)^{-1}\comp
\triv{\Op{R}}_{\text{\scriptsize\rm red}}^{t}(0,0,k,0)
\comp\begin{pmatrix}
\psi_{-\kappa}V & 0\\
0 & \psi_{\kappa}W
\end{pmatrix}\begin{pmatrix}
\psi_0\\
\psi_0
\end{pmatrix}$$
with the entries of
$$\left(\unity-\triv{\Op{R}}_{\text{\scriptsize\rm red}}(0,0,k,0)
\comp\begin{pmatrix}
\psi_{-\kappa}V & 0\\
0 & \psi_{\kappa}W
\end{pmatrix}\right)^{-1}\comp
\triv{\Op{R}}_{\text{\scriptsize\rm red}}(0,0,k,0)
\comp\begin{pmatrix}
\psi_{-\kappa}V & 0\\
0 & \psi_{\kappa}W
\end{pmatrix}\begin{pmatrix}
\psi_0\\
\psi_0
\end{pmatrix}$$
belong to $\banach{q}(\torus)\Tilde{\otimes}_{\Delta_q}
\banach{p}(\torus)$.
If $\var V,\var W\in\banach{p}(\torus)$,
then the integrals of the products with $\var V$ and $\var W$
are bounded regular (in the sense of \cite[Section~7.3.]{DF}) maps
from $\banach{q}(\torus)$ into $\mathbb{C}$.
Since the Fourier transform of the pointwise multiplication with the
sequences $\psi_{\pm\kappa}$ indexed by $\kappa\in\lattice\dual$
is the bounded operator $f\mapsto g$ with $g(x,x')=f(x,x'\pm x)$ on
$\banach{q}(\torus)\Tilde{\otimes}_{\Delta_q}\banach{p}(\torus)$,
we finally obtain, again due to \cite[Section~7.3. Theorem]{DF},
that the operator, which maps $(\var V,\var W)$ to the
Fourier transform of the restrictions of the sequence of operators
\begin{multline*}
\begin{pmatrix}
\psi_{-\kappa}V & 0\\
0 & \psi_{\kappa}W
\end{pmatrix}\comp\triv{\Op{R}}_{\text{\scriptsize\rm red}}(0,0,k,0)
\comp
\left(\unity-\begin{pmatrix}
\psi_{-\kappa}V & 0\\
0 & \psi_{\kappa}W
\end{pmatrix}\comp\triv{\Op{R}}_{\text{\scriptsize\rm red}}(0,0,k,0)
\right)^{-1}\comp\\
\comp\begin{pmatrix}
\psi_{-\kappa}\var V & 0\\
0 & \psi_{\kappa}\var W
\end{pmatrix}\comp\\
\comp\left(\unity-\triv{\Op{R}}_{\text{\scriptsize\rm red}}(0,0,k,0)
\comp\begin{pmatrix}
\psi_{-\kappa}V & 0\\
0 & \psi_{\kappa}W
\end{pmatrix}\right)^{-1}\comp
\triv{\Op{R}}_{\text{\scriptsize\rm red}}(0,0,k,0)\comp\begin{pmatrix}
\psi_{-\kappa}V & 0\\
0 & \psi_{\kappa}W
\end{pmatrix}
\end{multline*}
to the constant functions, is a bounded operator from
$\banach{p}(\torus)\times\banach{p}(\torus)$ into the
$\banach{p}$--functions on $\torus$
with values in the $2\times 2$--matrices.

It remains to estimate integrals over pointwise products
of the sequences $\psi_{-\kappa}V$ and $\psi_{\kappa}W$
with the entries of the sequence of functions
\begin{align*}
\left(\begin{smallmatrix}
\psi_{-\kappa}\var V & 0\\
0 & \psi_{\kappa}\var W
\end{smallmatrix}\right)\comp
\triv{\Op{R}}_{\text{\scriptsize\rm red}}^{t}(0,0,k,0)
\comp\left(\begin{smallmatrix}
\psi_{-\kappa}\var V\\
\psi_{\kappa}\var W
\end{smallmatrix}\right)&\text{ and}&
\left(\begin{smallmatrix}
\psi_{-\kappa}\var V & 0\\
0 & \psi_{\kappa}\var W
\end{smallmatrix}\right)\comp
\triv{\Op{R}}_{\text{\scriptsize\rm red}}(0,0,k,0)
\comp\left(\begin{smallmatrix}
\psi_{-\kappa}\var V\\
\psi_{\kappa}\var W
\end{smallmatrix}\right)
\end{align*}
and the entries of the sequence of functions
\begin{align*}
\left(\unity-\triv{\Op{R}}_{\text{\scriptsize\rm red}}^{t}(0,0,k,0)
\comp\left(\begin{smallmatrix}
\psi_{-\kappa}V & 0\\
0 & \psi_{\kappa}W
\end{smallmatrix}\right)\right)^{-1}\left(\begin{smallmatrix}
\psi_0\\
\psi_0
\end{smallmatrix}\right)&\text{and}&
\left(\unity-\triv{\Op{R}}_{\text{\scriptsize\rm red}}(0,0,k,0)
\comp\left(\begin{smallmatrix}
\psi_{-\kappa}V & 0\\
0 & \psi_{\kappa}W
\end{smallmatrix}\right)\right)^{-1}\left(\begin{smallmatrix}
\psi_0\\
\psi_0
\end{smallmatrix}\right),&
\end{align*}
respectively. With similar arguments
(some $V,W$ are replaced  by $\var V,\var W$)
we decompose the corresponding operators
again into sums of two operators, whose second terms yield
bounded operators from $\banach{2}(\torus)\times\banach{2}(\torus)$
into the $\banach{2}$--functions on $\torus$
with values in the $2\times 2$--matrices.

The images of the remaining terms are the Fourier transforms
of the restrictions of the sequences
\begin{eqnarray*}
\begin{pmatrix}
\psi_{-\kappa}\var V & 0\\
0 & \psi_{\kappa}\var W
\end{pmatrix}\comp\triv{\Op{R}}_{\text{\scriptsize\rm red}}(0,0,k,0)
\comp\begin{pmatrix}
\psi_{-\kappa}V & 0\\
0 & \psi_{\kappa}W
\end{pmatrix}\\
\begin{pmatrix}
\psi_{-\kappa} V & 0\\
0 & \psi_{\kappa} W
\end{pmatrix}\comp\triv{\Op{R}}_{\text{\scriptsize\rm red}}(0,0,k,0)
\comp\begin{pmatrix}
\psi_{-\kappa}\var V & 0\\
0 & \psi_{\kappa}\var W
\end{pmatrix}
\end{eqnarray*}
to the constant functions, respectively.
These correspond to off--diagonal matrices and induce
bounded operators from
$\banach{p}(\torus)\times\banach{p}(\torus)$ into the
$\banach{p}$--functions on $\torus$
with values in the $2\times 2$--off--diagonal matrices for all $1<p<2$.

\noindent
{\bf 4.} By definition of
$\fourier{\mathrm{v}}(V,W,\cdot),\fourier{\mathrm{w}}(V,W,\cdot)$
the variation of the off--diagonal entries of
$\Mat{A}_{V,W,k}$ yield variations of the sequences $k_{\kappa}$
and contribute indirectly through the partial derivatives
of the diagonal parts of $\Mat{A}_{V,W,k}$
with respect to $k$ to the variations of
\begin{align*}
\mathrm{v}&=\sum\limits_{\kappa\in\lattice\dual_{\delta}}
\fourier{\mathrm{v}}(V,W,\kappa) &
\mathrm{w}&=\sum\limits_{\kappa\in\lattice\dual_{\delta}}
\fourier{\mathrm{w}}(V,W,\kappa)
\end{align*}
The considerations concerning these partial derivatives in the proof
of Proposition~\ref{asymptotic analysis 2} combined with the methods
of the first three steps imply that for suitable $k\in\Set{O}$
all $\banach{q}$--norms with $1<q<\infty$ of these partial derivatives
are bounded (in terms of $\|V\|_2$ and $\|W\|_2$).
Since the Fourier transform of pointwise multiplication
is the convolution \cite[Theorem~IX.3]{RS2},
the variations of the off--diagonal terms of $\Mat{A}_{V,W,k}$
contribute to the variations of $\mathrm{v},\mathrm{w}$ over the
convolutions with the partial derivatives of the diagonal parts of
$\Mat{A}_{V,W,k}$ with respect to $k$. Due to Young's inequality
\cite[Section~IX.4]{RS2}, for all $q<\infty$
the $\banach{q}$--norms of these contributions are finite
and the lemma is proven.
\end{proof}

\noindent
{\it Continuation of the proof of Theorem~\ref{l1-structure}.}
The proof of this lemma implies that for all
$V,W\in\banach{2}(\torus)$ there exists a $\delta>0$ such that
the mapping
\begin{align*}
(V,W)&\mapsto(\mathrm{v},\mathrm{w})&\text{with }
\begin{pmatrix}
\mathrm{v}\\
\mathrm{w}
\end{pmatrix}&=\sum\limits_{\kappa\in\lattice\dual_{\delta}}
\begin{pmatrix}
\psi_{\kappa}\fourier{\mathrm{v}}(V,W,\kappa)\\
\psi_{\kappa}\fourier{\mathrm{w}}(V,W,\kappa)
\end{pmatrix}
\end{align*}
is locally a holomorphic mapping from an open neighbourhood of
$(V,W)\in\banach{2}(\torus)\times\banach{2}(\torus)$ onto an open
neighbourhood of $(\mathrm{v},\mathrm{w})\in
\banach{2}(\torus)\times\banach{2}(\torus)$.
In fact, in order to estimate
the $\banach{2}$--norms of $\mathrm{v}$ and $\mathrm{w}$ we integrate
the derivative along a path to the zero potentials with
$\mathrm{v}=\mathrm{w}=0$.
Furthermore, for small $\delta$ locally the restrictions
of the derivatives of this mapping to the inverse Fourier transform of
$\ell_2(\lattice\dual_{\delta})\times\ell_2(\lattice\dual_{\delta})$
are bounded invertible operators.
Hence the implicit function theorem
\cite[Supplementary material~V.5 Theorem~S.11]{RS1}
implies that this mapping is open.
In particular, the \Em{finite type potentials} are dense in
$\banach{2}(\torus)$. We remark that for potentials of the form
$(U,\Bar{U})$ the positivity of the moduli combined with
the asymptotic analysis of Proposition~\ref{asymptotic analysis 2}
implies that the $\banach{2}$--norms of
$\mathrm{u}=\sum\limits_{\kappa\in\lattice\dual_{\delta}}
\psi_{\kappa}\fourier{\mathrm{u}}(U,\kappa)$ are bounded
in terms of $\|U\|_2$. Moreover, the methods of
Lemma~\ref{continuous derivative} can be used in order to estimate
directly the $\banach{2}$--norms of $\mathrm{v}$ and $\mathrm{w}$
in terms of $\|V\|_2$ and $\|W\|_2$.

Due to the estimate
$\left\|\Mat{A}(\psi_{-\kappa}V,\psi_{\kappa}W,k)-
\left(\begin{smallmatrix}
\fourier{\mathrm{v}}(V,W,\kappa) & 0\\
0 & \fourier{\mathrm{w}}(V,W,-\kappa)
\end{smallmatrix}\right)\right\|\leq\varepsilon\left\|k-k_{\kappa}\right\|$
in the proof of Proposition~\ref{asymptotic analysis 2}
these arguments extend to the moduli $\parameter{t}$ and show
that the mapping $(V,W)\mapsto\parameter{t}$ is a holomorphic
open mapping from $\banach{2}(\torus)\times\banach{2}(\torus)$
into $\ell_1(\lattice\dual_{\delta})$ with suitable $\delta$
depending on $V$ and $W$.

\noindent
{\bf 4.} Now we shall define locally all moduli
$\parameter{t}(\kappa)$. We shall complete the cycles of the handles
index by $\kappa\in\lattice\dual_{\delta}$ to a set of generators
of a Lagrangian subgroup of the first homology group of the
\Em{complex Fermi curves} with respect to the intersection form.
In addition we assume that the integrals of $dk$ along these cycles
vanish. Moreover, in neighbourhoods of \Em{complex Fermi curves}
with singularities, we assume that a small neighbourhood of all
deformed singularities contain as many cycles, as the
\De{Local contribution to the arithmetic genus}~\ref{local contribution}.
The cycles of the handles described in
Theorem~\ref{asymptotic analysis 1} and
Proposition~\ref{asymptotic analysis 2} have these properties.
Due to Proposition~\ref{asymptotic analysis 2}
the elements of this basis may be indexed by $\kappa\in\lattice\dual$.
For \Em{complex Fermi curves} of the form $\fermi(U,\Bar{U})$
we complete the \Em{horizontal cuts} indexed by
$\kappa\in\lattice\dual_{\delta}$ to a set of generators
of those cycles in the first homology group,
which are invariant under the corresponding action of $\eta$.
Consequently we define all moduli as
$$\parameter{t}(\kappa)=
\pi\oint\limits_{\text{\scriptsize cycle with index $\kappa$.}}
k_1dk_2.$$
Since the integral of $dk$ along these cycles vanish,
the corresponding moduli do not depend on the choice of the branch.
Moreover, the moduli of \Em{complex Fermi curves} of the form
$\fermi(U,\Bar{U})$ are real.
Due to Proposition~\ref{smooth moduli}, these moduli define
locally coordinates of the manifold $\moduli_{g,\lattice}$.

We shall now show that locally these moduli completely determine the
\Em{complex Fermi curves};
i.\ e.\ if two potentials have locally the same moduli,
then they have the same \Em{complex Fermi curve}.
Let $\Set{U}$ be a small open subset of
$\banach{2}(\torus)\times\banach{2}(\torus)$,
on which a complete choice of moduli is defined.
Two pairs of potentials in $\Set{U}$ whose moduli coincide,
are limits of two convergent sequences in $\Set{U}$,
whose \Em{complex Fermi curves} have finite \Em{geometric genus},
and whose sequences of moduli coincide. Therefore, due to the
continuity of the map $(V,W)\mapsto\fermi(V,W)$, it suffices to show
that the moduli completely determine the \Em{complex Fermi curves}
of finite \Em{geometric genus}.
This follows from Proposition~\ref{smooth moduli}.

To sum up, these moduli define locally holomorphic coordinates of
$\moduli_{\lattice}$, which take values in $\ell_1(\lattice\dual)$.
In particular, $\moduli_{\lattice}$ is a
$\ell_1(\lattice\dual)$--manifold.

\noindent
{\bf 5.} We conclude that all \Em{complex Fermi curves} in
$\moduli_{\lattice}$,
whose normalization is invariant under $\eta$ without fixed points,
belongs to
$\Bar{\moduli}_{\lattice,\eta,\willmore}$ with
$\willmore=4\int\limits_{\torus}V(x)W(x)d^2x$.
This follows from Lemma~\ref{existence of potentials} and the fact
that the potentials, whose \Em{complex Fermi curves} have finite
\Em{geometric genus},
are dense in $\banach{2}(\torus)\times\banach{2}(\torus)$.
Consequently the moduli space
$\moduli_{\lattice,\eta}$
is a closed real subspace of $\moduli_{\lattice}$
(with respect to the $\ell_1(\lattice\dual)$--topology).
Moreover, the moduli with indexes in $\lattice\dual_{\delta}$
are non--negative. Therefore, this moduli space has
an asymptotic $\ell^+_1(\lattice\dual_{\delta})$--structure.

We may extend this local description of the moduli space to
neighbourhoods of \Em{complex Fermi curves} of
\De{Finite rank Perturbations}~\ref{finite rank perturbations}
of a limit of sequences of potentials in $\banach{2}(\torus)$.
In fact, to these sequences we may add potentials, whose support is
contained in the complement of the unions $\Set{S}_{D,\varepsilon}$
of small balls around the singular points in $\torus$
(compare with Theorem~\ref{limits of resolvents} and
Remark~\ref{extension of the symplectic form}).
Since the union
$$\bigcup\limits_{\varepsilon>0}\left\{\var U\in\banach{2}(\torus)\mid
\text{\rm support}(\var U)\cap\Set{S}_{D,\varepsilon}=\emptyset\right\}$$
is dense in $\banach{2}(\torus)$ all elements of
$\Bar{\moduli}_{\lattice,\eta}$
are contained in an open set, which is mapped by the restricted moduli
onto an open set of $\ell^+_1(\lattice\dual_{\delta})$.
\end{proof}

In particular, we proved

\begin{Corollary}\label{eta invariant}
All \Em{complex Fermi curves} $\fermi(V,W)$ with potentials
$V,W\in \banach{2}(\torus)$,
whose normalizations are invariant under the involution $\eta$
but without fixed points, belong to
$\Bar{\moduli}_{\lattice,\eta,\willmore}$ with
$\willmore=4\int\limits_{\torus}V(x)W(x)d^2x$.
\qed
\end{Corollary}

\begin{Remark}~\label{additional restrictions}
The boundedness of the \Em{first integral} of the
\Em{complex Fermi curves} yields additional restrictions of
\De{Finite rank Perturbations}~\ref{finite rank perturbations}
of $\Op{R}(U,\Bar{U},k,0)$.
In fact, it is quite easy to see that the \Em{complex Fermi curves} of
the \De{Finite rank Perturbations}~\ref{finite rank perturbations}
of the free resolvent have \Em{disconnected normalization},
if the \Em{first integral} is bounded.
However, there exists such \Em{rank one Perturbations},
whose normalizations are connected.
The corresponding \Em{complex Fermi curves} might correspond to
potentials $U\in\banach{2,\infty}(\torus)$.
\end{Remark}

Due to the estimates established in the proof of
Proposition~\ref{asymptotic analysis 2},
the derivatives $\var\parameter{t}(\kappa)/\var U$
of the moduli are locally analytic vector fields on
$\banach{2}(\torus)$. For all compact Riemann surfaces $\Spa{Y}$ the
long exact cohomology sequence
$$0\rightarrow \mathbb{Z}\rightarrow \mathbb{C}\xrightarrow{\exp}
\mathbb{C}^{\ast}\rightarrow H^1(\Spa{Y},\mathbb{Z})\rightarrow
H^1(\Spa{Y},\Sh{O})\xrightarrow{\exp}H^1(\Spa{Y},\Sh{O}^{\ast})\rightarrow
H^2(\Spa{Y},\mathbb{Z})\rightarrow 0$$
corresponding to the exact sequence
$$0\rightarrow\mathbb{Z}\hookrightarrow\Sh{O}
\xrightarrow{\exp}\Sh{O}^{\ast}\rightarrow 1$$
of sheaves shows that the first cohomology group $H^1(\Spa{Y},\mathbb{Z})$
considered as a subspace of the Lie algebra of the Picard group is
exactly the kernel of the exponential mapping of the Picard group.
By definition of the moduli $\parameter{t}(\kappa)$ and the one--forms
$\Omega_{V,W}$ (compare with Lemma~\ref{regular form})
the variations $\var\parameter{t}(\kappa)$
with respect to $(\var V,\var W)$ is equal to the integral
$\Omega_{V,W}(\var V,\var W)$ over the unique cycle
of the handle with index $\kappa$.
For \Em{finite type potentials} Lemma~\ref{residue} implies therefore
that the Hamiltonian vector fields of the moduli
$\parameter{t}(\kappa)$ induce periodic flows with period $1$.
Since the \Em{finite type potentials} are dense, the same is true for
all potentials. Since in the limit $g(\kappa,\kappa)\rightarrow\infty$
the norm of these vector fields converges to zero,
for all $U\in\banach{2}(\torus)$ and all $\varepsilon>0$ there exists
some $\delta>0$ such that the flows corresponding
to the Hamiltonian vector fields of all moduli $\parameter{t}(\kappa)$
with $\kappa\in\lattice\dual_{\delta}$ maps $U$ for all times into
$B(U,\varepsilon)$. In particular, for all $U\in\banach{2}(\torus)$
and all small $\varepsilon>0$ the isospectral set of $U$ contains an
infinite--dimensional torus in $B(U,\varepsilon)$,
which is homeomorphic to
$\left(\mathbb{R}/\mathbb{Z}\right)^{\left\{\kappa\in
\lattice\dual_{\delta}\mid\parameter{t}(\kappa)>0\right\}}$.
Also the corresponding action is induced by the Hamiltonian vector fields
of the moduli with indexes $\kappa\in\lattice\dual_{\delta}$.

\section{The \De{Generalized Willmore functional}}
\label{section generalized}
\subsection{The \De{Weak Singularity condition}}
\label{subsection weak singularity}

Not all real potentials $U$ correspond to some immersion of a
torus into $\mathbb{R}^3$. In this section we want to find some
characteristic properties of
all \Em{complex Fermi curves}, which correspond to immersions. All
\Em{complex Fermi curves}, which correspond to immersions have,
roughly speaking, some higher order singularity:

\begin{Lemma} \label{fourth singularity}
Assume that the Dirac operator 
$\left(\begin{smallmatrix}
U & \partial \\
-\Bar{\partial} & U
\end{smallmatrix}\right)$
with real potential $U$ has an eigenfunction
$\chi(x)=\left(\begin{smallmatrix}
\chi_1(x)\\
\chi_2(x)
\end{smallmatrix}\right)$
with eigenvalue zero, which obeys $\chi(x+\gamma)=\pm\chi(x)$ for all
$x\in \mathbb{R}^2$ and all periods $\gamma\in\lattice$ and the
\De{Equivalent form of the Periodicity condition}~\ref{equivalent form}.
  Then there exists some
  $\kappa\in\lattice\dual$, whose equivalence class $[\kappa/2]$
  belongs to the \Em{complex Fermi curve}
  $\fermi(U,U)/\lattice\dual$.
  Moreover, the preimage of this point in the normalization of
  $\fermi(U,U)/\lattice\dual$ contains  either at least
  four elements or two elements, which are zeroes of both
  differentials $d\xx{p}$ and $d\yy{p}$.
\end{Lemma}

\begin{proof} The transformation of the eigenfunction under
  shifts of the periods implies that there exists some element
  $\kappa$ of the dual lattice, such that $[\kappa/2]$ belongs to
  $\fermi(U,U)/\lattice\dual$. Due to (ii) of
  Corollary~\ref{fixed points} the
  preimage of this element $[\kappa/2]$ in the normalization of
  $\fermi(U,U)/\lattice\dual$ has an even number of
  elements. Now let us assume that the preimage has two elements $y$
  and $y'$. 
  Due to Corollary~\ref{fixed points} $y'$ is equal to
  $\rho\left(\sigma(y)\right)=\sigma\left(\rho(y)\right)$. Thus
  the differentials $dk_1$ and $dk_2$ have zeroes either at both
  points $y$ and $y'$ or at none of them. Let us assume on the
  contrary to the statement of the lemma that there
  exists some period $\xx{\gamma}$, such that the differential
  $dg(\xx{\gamma},k)$ is not equal to zero at $y$ and $y'$.
  Thus the projection
  $\Breve{\Op{P}}_{\xx{\gamma}}([\kappa/2])$ defined in (iii) of
  Lemma~\ref{projection 2} has rank two. Due to (v) of
  Lemma~\ref{projection 2} the eigenfunction $\chi$, which fulfills the
  conditions of the theorem, is invariant under
  $\Breve{\Op{P}}_{\xx{\gamma}}([\kappa/2])$.
  Hence, this eigenfunction $\chi$
  and $\Op{J}\Bar{\chi}$ span the eigenspace of this
  projection. The normalization of $\fermi(U,U)/\lattice\dual$ is a
  smooth Riemann surface. Hence the pullback of the eigenfunctions
  $\psi$ and $\phi$ to the normalization can be chosen to be locally
  holomorphic without any zeroes. Since one of the involutions
  $\sigma$ and $\rho$
  permutes $y$ and $y'$ and the other involution leaves
  these points invariant, the span of $\phi(y)$ and $\phi(y')$ is
  equal to the span of $\Op{J}\psi(y)$ and $\Op{J}\psi(y')$
  and to the span of $\Bar{\psi}(y)$ and $\Bar{\psi}(y')$. Thus the
  \De{Equivalent form of the Periodicity condition}~\ref{equivalent form}
  implies that $\Op{P}_{\xx{\gamma}}([\kappa/2])\psi(y)$
  and $\Op{P}_{\xx{\gamma}}([\kappa/2])\psi(y')$
  are equal to zero.
  This is a contradiction to (vi) of Lemma~\ref{projection 2}.
  Thus both differentials $dk_1$ and $dk_2$ are zero
  at these two elements $y$ and $y'$.
\end{proof}

For any singularity of the \Em{complex Fermi curve} we may
distinguish between multiple points and cusps:
\begin{description}
\item[Multiple point:]
\index{multiple point}
\index{singularity!$\rightarrow$ multiple point}
   A point of the normalization of the
  \Em{complex Fermi curve} in the preimage of some singularity is
  called a \Em{multiple point}, if there exists some $\gamma\in\lattice$,
  whose differential $dg(\gamma,k)$ does not vanish at this point.
\item[Cusp:]
\index{cusp}
\index{singularity!$\rightarrow$ cusp}
  A point of the normalization of the
  \Em{complex Fermi curve} in the preimage of some singularity is
  called \Em{cusp}, if it is a zero of both differentials $dk_1$ and
  $dk_2$.
\end{description}
Now we may paraphrase the condition of Lemma~\ref{fourth singularity}.

\begin{Weak Singularity condition}\label{weak singularity condition}
\index{singularity!condition!weak $\sim$}
\index{condition!weak singularity $\sim$}
  (See also \ref{weak singularity condition 1}.) There exists some
  $\kappa\in\lattice\dual$, such that $[\kappa/2]$ belongs to the
  \Em{complex Fermi curve}. Moreover, the preimage
  of this point in the normalization of $\fermi/\lattice\dual$
  contains either at least four \Em{multiple points}, or at least two
  \Em{cusps}.
\end{Weak Singularity condition}

If we decompose the preimage of $[\kappa/2]$ in the normalization of
$\fermi(U,U)/\lattice\dual$ into orbits of the two involutions $\sigma$
and $\rho$ three possible orbits can occur.
\begin{description}\index{orbit}
\item[Orbit of type 1:]\index{orbit!of type~1}
  Two fixed points of $\rho$, which are permuted by
  $\sigma$.
\item[Orbit of type 2:]\index{orbit!of type~2}
  Two fixed points of $\sigma$, which are permuted by
  $\rho$.
\item[Orbit of type 3:]\index{orbit!of type~3}
  Four points and both involutions $\sigma$ and
  $\rho$ have no fixed point.
\end{description}
In particular, the
\De{Weak Singularity condition}~\ref{weak singularity condition}
permits the following situations:
\begin{description}
\item[Multiple points O1/O1:]
\index{multiple point}
\index{orbit!of type~1}
  Two \Em{Orbits of type~1} containing
  \Em{multiple points} (8).
\item[Multiple points O1/O2:]
\index{multiple point}
\index{orbit!of type~1}
\index{orbit!of type~2}
  One \Em{Orbit of type~1} containing
  \Em{multiple points} and one \Em{Orbit of type~2} containing
  \Em{multiple points} (7).
\item[Multiple points O2/O2:]
\index{multiple point}
\index{orbit!of type~2}
  Two \Em{Orbits of type~2} containing
  \Em{multiple points} (6).
\item[Multiple points O3:]
\index{multiple point}
\index{orbit!of type~3}
  One \Em{Orbit of type~3} containing
  \Em{multiple points} (8).
\item[Cusps O1:]
\index{cusp}
\index{orbit!of type~1}
  One \Em{Orbit of type~1} containing \Em{cusps} (8).
\item[Cusps O2:]
\index{cusp}
\index{orbit!of type~2}
  One \Em{Orbit of type~2} containing \Em{cusps} (23).
\item[Cusps O3:]
\index{cusp}
\index{orbit!of type~3}
  One \Em{Orbit of type~3} containing \Em{cusps} (37).
\end{description}
These situations differ in the possible orders of the
singularity. The numbers in the round brackets are the
lowest \De{Local contribution to the arithmetic genus}
\index{arithmetic genus!local contribution to the $\sim$|(}
\index{local contribution!to the $\rightarrow$ arithmetic genus|(}
of the \Em{complex Fermi curves}
The last two cases yield very high contributions and they are
exceptional situations.

\newtheorem{Local contribution}[Lemma]{Local contribution to the
  arithmetic genus}
\begin{Local contribution}\label{local contribution}
If $\Sh{O}_{\Spa{Y},\text{\scriptsize\rm normal}}$ is the
normalization sheaf (i.\ e.\ the direct image of the
structure sheaf of the normalization under the normalization map
\cite[Chapter~6. \S4.]{GrRe}) of a
compact pure one--dimensional complex space $\Spa{Y}$, then the
long exact cohomology sequence \cite[\S15.]{Fo} corresponding
to the exact sequence of sheaves
$$0\rightarrow\Sh{O}_\Spa{Y}\rightarrow
\Sh{O}_{\Spa{Y},\text{\scriptsize\rm normal}}\rightarrow
\Sh{O}_{\Spa{Y},\text{\scriptsize\rm normal}}/\Sh{O}_\Spa{Y}
\rightarrow 0$$
shows that the arithmetic genus of $\Spa{Y}$, which is equal to
$\dim H^1\left(\Spa{Y},\Sh{O}_\Spa{Y}\right)$, is equal to the
\Em{geometric genus}, which is equal to
$\dim H^1\left(\Spa{Y},
\Sh{O}_{\Spa{Y},\text{\scriptsize\rm normal}}\right)$, plus
$\dim H^0\left(\Spa{Y},
\Sh{O}_{\Spa{Y},\text{\scriptsize\rm normal}}/\Sh{O}_\Spa{Y}\right)$.
The support of this sheaf is the singular locus of $\Spa{Y}$
\cite[Chapter~6. \S5.3.]{GrRe},
and the space of global sections is equal to the direct sum of all
non--trivial stalks. The dimension of these stalks is the
local contribution to the arithmetic genus
at the corresponding singularities
\cite[Chapter~IV \S1.4. and \S2.7.]{Se}.
\index{arithmetic genus!local contribution to the $\sim$|)}
\index{local contribution!to the $\rightarrow$ arithmetic genus|)}
\end{Local contribution}

The lemma shows that the \Em{complex Fermi curve} of some potential, which
corresponds to some immersion, obeys some condition. In
Section~\ref{subsection singularity} we will see that this
condition on the \Em{complex Fermi curve} of some real potential implies
neither that this potential is a
\De{Weierstra{\ss} potential}~\ref{Weierstrass potentials} nor
that the \Em{complex Fermi curve} is a
\De{Weierstra{\ss} curve}~\ref{Weierstrass curves}.
In fact, in Section~\ref{subsection singularity} we shall formulate
the \De{Singularity condition}~\ref{singularity condition}, which is
equivalent to the condition that the \Em{complex Fermi curve} is a
\De{Weierstra{\ss} curve}~\ref{Weierstrass curves}. Moreover, we will
describe the subset of the corresponding \Em{isospectral sets}, which
contain the \De{Weierstra{\ss} potentials}~\ref{Weierstrass potentials}.
For \Em{Complex Fermi curves}, which are invariant under the
involution $\sigma$, however,
the \De{Weak Singularity condition}~\ref{weak singularity condition},
is equivalent to the
\De{Quaternionic Singularity condition}~\ref{quaternionic condition}
combined with the condition on the singularity
being invariant under $\sigma$. In particular, the
\De{Weak Singularity condition}~\ref{weak singularity condition}
is stronger than the
\De{Quaternionic Singularity condition}~\ref{quaternionic condition}.
In some sense the
\De{Weak Singularity condition}~\ref{weak singularity condition}
seems to be the most convenient combination of the advantages of
immersion into $\mathbb{R}^4$ (corresponding to the
\De{Quaternionic Singularity condition}~\ref{quaternionic condition}),
which yield stronger conditions on relative minimizers,
and the advantages of immersion into $\mathbb{R}^3$,
whose \Em{complex Fermi curves} have an additional involution $\sigma$.
In fact, in contrast to the
\De{Singularity condition}~\ref{singularity condition}
the \De{Weak Singularity condition}~\ref{weak singularity condition}
allows to deform \Em{Orbits of type~1} and \Em{Orbits of type~3}
into \Em{Orbits of type~2}. This deformation is essential for the
proof that the relative minimizers have dividing real parts.
Moreover, in Section~\ref{subsection absolute minimizer}
we shall see that this property is almost as restrictive
as the \Em{real--$\sigma$--hyperellipticity},
which characterizes minimal tori in $S^3$ \cite{Hi}.

In Section~\ref{subsection compactified moduli} we
have seen that the moduli space has a compactification containing
\Em{complex Fermi curves} of
\Em{Finite rank perturbations}~\ref{finite rank perturbations}.
Hence it is natural to extend
the Willmore functional not only to all
\De{Generalized Weierstra{\ss} potentials}~\ref{generalized potentials},
but also to all
\De{Generalized Weierstra{\ss} curves}.

\newtheorem{Generalized Weierstrass curves}[Lemma]{Generalized
  Weierstra{\ss} curves}
\index{Weierstra{\ss}!generalized $\sim$ curve}
\index{curve!generalized Weierstra{\ss} $\sim$}
\index{generalized!Weierstra{\ss}!curve}
\begin{Generalized Weierstrass curves}\label{generalized curves}
\Em{Complex Fermi curves} in
$\Bar{\moduli}_{\lattice,\eta,\sigma}$, which obey the
\De{Weak Singularity condition}~\ref{weak singularity condition}
are called \De{Generalized Weierstra{\ss} curves}.
\end{Generalized Weierstrass curves}

Due to Lemma~\ref{compact metric}, for all $\willmore>0$ the spaces
$\Bar{\moduli}_{\lattice,\eta,\sigma,\willmore}$ are compact with
respect to the subspace topology of the space of compact subsets of
$\overline{\mathbb{C}^2}$. Since the subset consisting of
\De{Generalized Weierstra{\ss} curves}~\ref{generalized curves}
is closed the restriction of the \Em{first integral} to this subspace
is lower semi--continuous. This restriction is called

\newtheorem{Generalized Willmore functional}[Lemma]{Generalized
  Willmore functional}
\index{Willmore!functional!generalized $\sim$}
\index{generalized!Willmore functional}
\begin{Generalized Willmore functional}\label{generalized functional}
\end{Generalized Willmore functional}

In Theorem~\ref{l1-structure} we have seen that the compactified
moduli space has an $\ell^+_1(\lattice\dual)$--structure.
Consequently, this compactified moduli space is very similar to a
manifold and some specific regular one--forms correspond to
tangent vectors.

\subsection{Relative minimizers of the \De{Generalized Willmore functional}}
\label{subsection relative minimizers}
\index{relative minimizer|(}
\index{minimizer!relative $\sim$|(}

In this section we investigate relative minimizers of the
\De{Generalized Willmore functional}~\ref{generalized functional},
i.\ e.\ minimizers of the restrictions to open sets \cite{Str}.
The main result is that all relative minimizers have
\begin{description}
\index{real!part!dividing $\sim$}
\index{dividing!real part}
\item[Dividing real parts.] A \Em{complex Fermi curve}
  $\fermi\in\Bar{\moduli}_{\lattice,\eta,\sigma}$ has a dividing real part,
  if the relative complement of the real part
  (i.\ e.\ the set of fixed points of $\rho$)
  in the normalization of $\fermi/\lattice\dual$
  has two connected components \cite[11.6]{BCR}
\end{description}

\begin{Remark}
Due to Proposition~\ref{asymptotic analysis 2} 
(compare with Corollary~\ref{dividing genus bound}),
all \Em{complex Fermi curves} with \Em{dividing real parts}
have finite \Em{geometric genus}.
In Corollary~\ref{finite type} we shall show with the arguments of
Proposition~\ref{global meromorphic function} that all
\De{Constrained Willmore tori}~\ref{constrained Willmore tori}
have finite \Em{geometric genus}.
Unfortunately these arguments do not apply to the
\De{Generalized Willmore functional}~\ref{generalized functional}.
More precisely, in order to apply
Proposition~\ref{global meromorphic function} to the
\De{Generalized Willmore functional}~\ref{generalized functional},
we have to prove that all variations of
\De{Generalized Weierstra{\ss} potentials}~\ref{generalized potentials}
may be characterized by the finite--order vanishing of the corresponding
one--form at the preimage in the normalization of the singularity
described in the
\De{Weak Singularity condition}~\ref{weak singularity condition}
analogous to Proposition~\ref{variation}. But the subvariety of
\De{Generalized Weierstra{\ss} potentials}~\ref{generalized potentials}
may have singularities, whose \Em{complex Fermi curves} may be of
infinite \Em{geometric genus}.
More precisely, Lemma~\ref{residue} implies that the corresponding 
\Em{complex Fermi curves} have global meromorphic functions, but
they may take the same values at $\infty^{\pm}$.
Therefore, the corresponding \Em{complex Fermi curves} can have
infinite \Em{geometric genus}. Varieties with singularities have
different notions of tangent spaces \cite[Chapter~7]{Wh}.
In particular, differentiable paths contained in the variety
correspond to elements of the tangent cone
\cite[Chapter~7 Section~2 and Theorem~3C]{Wh}.
Hence in order to carry over the methods of Corollary~\ref{finite type}
it would suffice to characterize the elements of the tangent cone
by a finite--order vanishing of the corresponding
one--forms at the preimage in the normalization
of the singularity mentioned in the 
\De{Weak Singularity condition}~\ref{weak singularity condition}.
We did not succeed in proving this. For this reason the
characterization of relative minimizers of the
\De{Generalized Willmore functional}~\ref{generalized functional}
in this Section has to apply to
\De{Generalized Weierstra{\ss} curves}~\ref{generalized curves}
of infinite \Em{geometric genus}. For this purpose we proved
Theorem~\ref{l1-structure}, which will turn out to be the essential
tool of the extension of the deformation theory devolved in
Section~\ref{subsection deformation} to \Em{complex Fermi curves}
of infinite \Em{geometric genus}.
\end{Remark}

Of special importance is the subclass of
\begin{description}
\index{real!--$\sigma$--hyperelliptic}
\item[Real--$\sigma$--hyperelliptic] \Em{complex Fermi curves}.
  In this case the normalization of the quotient of
  $\fermi/\lattice\dual$ modulo $\sigma$ is biholomorphic
  to $\mathbb{P}^1$ with two points removed.
  Moreover, on this quotient the involutions $\eta$ and $\rho$
  induce the unique involution of $\mathbb{P}^1$ with one real cycle.
\end{description}

\begin{Remark}\label{real hyperelliptic}
A real hyperelliptic Riemann surface is determined by
a collection of branch points on $\mathbb{P}^1$, which are
invariant under some anti--holomorphic involution of
$\mathbb{P}^1$. Since $\mathbb{P}^1$ has two equivalence
classes of anti--holomorphic involutions, these real hyperelliptic
Riemann surfaces also splits into two classes. In the first case
the anti--holomorphic involution on $\mathbb{P}^1$
has one real cycle, and in the second case the real part is empty.
\end{Remark}

Due to Corollary~\ref{regular fermi} the normalization is not
necessarily connected. If the real part shrinks to one or several
singular points we obtain as a another subclass
the \Em{complex Fermi curves} with
\begin{description}
\index{normalization!disconnected}
\index{disconnected!normalization}
\item[Disconnected normalization.]
\end{description}
Since the fixed points of $\sigma$ are the
branch points of the natural two--sheeted covering of the
\Em{complex Fermi curve} over the quotient modulo $\sigma$,
due the Riemann--Hurwitz formula \cite[17.14]{Fo}
their number has to be even. Therefore, each connected component of
a disconnected normalization has besides the fixed point at
$\infty^{\pm}$ in addition an odd number of fixed points of $\sigma$.
Since these fixed points of $\sigma$ may be deformed
into real cycles, all \Em{complex Fermi curves} in
$\Bar{\moduli}_{\lattice,\eta,\sigma}$
with disconnected normalization are indeed limits of
\Em{complex Fermi curves} with \Em{dividing real parts}.

\begin{Remark}\label{real part}
Both anti--holomorphic involutions $\rho$ and $\eta$ induce on
the quotient modulo $\sigma$ the same anti--holomorphic involution,
whose fixed points are the real part of the quotient modulo $\sigma$.
Since $\eta$ does not have any fixed points on the normalization,
the preimage of the real part of the normalization
of the quotient modulo $\sigma$ in normalization is the real part
of the normalization. Therefore a \Em{complex Fermi curve} has a
\Em{dividing real part} if and only if the normalization of the
quotient modulo the involution $\sigma$ has a
\Em{dividing real part}.
\end{Remark}

\begin{Theorem} \label{relative minimizers}
The relative minimizers of the
\De{Generalized Willmore functional}~\ref{generalized functional}
belong to one of the following classes:
\begin{description}
\item[Minimizers with dividing real parts:]
  The normalization of the quotient of
  $\fermi/\lattice\dual$ modulo the involution $\sigma$ is a
  real hyperelliptic Riemann surface, whose hyperelliptic involution does
  not have any real fixed points.
  More precisely, the real part is a two--sheeted covering of
  the real cycle of $\mathbb{P}^1$,
  which is the real part of the quotient modulo the hyperelliptic
  involution.
  Moreover, the preimage in the normalization of the singularity
  described in the
  \De{Weak Singularity condition}~\ref{weak singularity condition}
  contains \Em{Multiple points O1/O1}, and no other orbit.
  Finally, the function $d\yy{p}/d\xx{p}$ takes different values on
  the two \Em{Orbits of type~1}, and the exterior derivative of this
  function does not vanish on the two \Em{Orbits of type~1}.
\item[Real--$\sigma$--hyperelliptic minimizers:]
  The minimizers are \Em{real--$\sigma$--hyperelliptic} and
  the pre\-im\-age in the normalization of the singularity
  described in the
  \De{Weak Singularity condition}~\ref{weak singularity condition}
  contains either \Em{Multiple points O1/O1},
  or \Em{Multiple points O1/O2}, or \Em{Cusps O1}, and no other orbit.
\item[Minimizers with disconnected normalization:]
  The minimizers have \Em{disconnected normalizations} and
  the preimage in the normalization of the singularity described in the
  \De{Weak Singularity condition}~\ref{weak singularity condition}
  contains either \Em{Multiple points O2/O2}, or
  \Em{Multiple points O3}, or \Em{Cusps O2} or \Em{Cusps O3}.
\end{description}
\end{Theorem}

\begin{proof} In a first step we show that all relative minimizers in
$\bigcup\limits_{g\in\mathbb{N}}
\Bar{\moduli}_{g,\lattice,\eta,\sigma}$ belong to one of the three
cases described in the Theorem. In a second step we extend our
arguments to all relative minimizers in
$\Bar{\moduli}_{\lattice,\eta,\sigma}$.

\noindent
{\bf 1. Local minimizers of finite \Em{geometric genus}.}
Let us describe the general construction of local
flows on the space $\bigcup\limits_{g\in\mathbb{N}}
\Bar{\moduli}_{g,\lattice,\eta,\sigma}$.
Due to Proposition~\ref{smooth moduli} the
complex moduli space $\moduli_{g,\lattice}$
of algebraic curves with arithmetic genus equal to $g$
is a $(g+1)$--dimensional complex manifold. Moreover, at each
element $\Spa{Y}$ of this manifold the tangent space is isomorphic to the
space of regular forms of the algebraic curve $\Tilde{\Spa{Y}}$, which is
obtained from $\Spa{Y}$ by identifying the
two marked points $\infty^-$ and $\infty^+$ to an ordinary double point.
The \De{Generalized Weierstra{\ss} curves}~\ref{generalized curves}
are elements of the moduli space
$\Bar{\moduli}_{g,\lattice,\eta,\sigma}$, and may be
characterized as those elements of the complex moduli space, which
are invariant under the involutions $\eta:k\mapsto-\Bar{k}$ and
$\sigma:k\mapsto -k$, and on whose normalizations the first
involution does not have fixed points. Due to
Lemma~\ref{transformations} these involutions induce involutions on
the moduli space\index{involution!of the moduli space}
$\moduli_{g,\lattice}$, and those \Em{complex Fermi curves},
which are invariant under the former
involutions, are the fixed points of the latter involutions.
For all \Em{complex Fermi curves} in
$\Bar{\moduli}_{g,\lattice,\eta,\sigma}$ we may find some open
neighbourhood of the corresponding element of
$\moduli_{g,\lattice}$, which is invariant under these
involutions of the complex moduli space. Due to the proof of
Proposition~\ref{smooth moduli} the regular forms of $\Tilde{\Spa{Y}}$
form a holomorphic vector bundle on the complex moduli space, which
is isomorphic to the tangent bundle. Hence any regular form $\omega$
fulfilling the equations $\eta^{\ast}\omega=\Bar{\omega}$ and
$\sigma^{\ast}\omega=\omega$ on $\Spa{Y}$ may be extended to a holomorphic
vector field on some neighbourhood of $\Spa{Y}$ in the
complex moduli space, which fulfills the corresponding
extensions of these equations to this
invariant neighbourhood. The flows corresponding to
these vector fields leave invariant the fixed point set
of the two involutions on the moduli space.
With this method we construct local flows
on the fixed point set of both involutions of the
complex moduli space. In a second step we investigate the
fixed points of the involutions $\eta$ on the corresponding
invariant \Em{complex Fermi curves}. It may happen that in spite of
$\eta$ on $\Spa{Y}$ having no fixed point, this is not true for some
neighbourhood of $\Spa{Y}$ in the fixed point set of both involutions on
the complex moduli space. We will see that
for some one--dimensional families $\fermi_{\parameter{t}}$
of \Em{complex Fermi curves} being invariant under $\eta$,
this involution has fixed points on the normalizations of
$\fermi_{\parameter{t}}$ for $\parameter{t}<0$ and no fixed points
for $\parameter{t}\geq 0$.
In these situations we call those regular forms
fulfilling $\eta^{\ast}\omega=\Bar{\omega}$ and $\sigma^{\ast}\omega=\omega$
\Em{admissible},
\index{form!admissible $\sim$}
\index{admissible form}
whose flows corresponding to positive times do not
leave the subset of \Em{complex Fermi curves} with no
fixed points of $\eta$ on the corresponding normalizations.
The set of these \Em{admissible} regular forms considered as elements
of the tangent space of the moduli space
(Proposition~\ref{smooth moduli} and Theorem~\ref{l1-structure})
is called real tangent semi cone along the subvariety of
\De{Generalized Weierstra{\ss} curves}~\ref{generalized curves}
\cite[Chapter~7 Section~2 and Theorem~3C]{Wh}.

Let us apply these considerations to those singularities of $\Spa{Y}$,
which are different from the singularity described in the
\De{Weak Singularity condition}~\ref{weak singularity condition}.
The former singularities may be removed without leaving the subset of
\De{Generalized Weierstra{\ss} curves}~\ref{generalized curves}.
We shall distinguish between the following cases:
\begin{description}
\item[Ordinary double points of the form $(y,\eta(y))$.]
Without loss of generality we may assume that both components of
$k$ are zero at these two points, and that the differential
$d\xx{p}$ does not vanish at these double points.
The Taylor series of the function $\yy{p}(\xx{p})$ at these
two points has to look like
$\yy{p}=a_1\xx{p}+a_2\xx{p}^2+\ldots$ and
$\yy{p}=\Bar{a}_1\xx{p}-\Bar{a}_2\xx{p}^2+\ldots$, respectively.
Therefore, the \Em{complex Fermi curve} is locally given by the
equation $R(\xx{p},\yy{p})=0$. More precisely, $R(\xx{p},\yy{p})=$
$$(\yy{p})^2-\yy{p}
\left((a_1+\Bar{a}_1)\xx{p}+(a_2-\Bar{a}_2)\xx{p}^2+\ldots\right)
+a_1\Bar{a}_1\xx{p}^2
+(a_2\Bar{a}_1-a_1\Bar{a}_2)\xx{p}^3+\ldots=0.$$
Let us assume in addition that the imaginary part of $a_1$ does not
vanish. Otherwise the function $d\yy{p}/d\xx{p}$ takes the same
values at both elements of the double point, and we have some
higher oder singularity. The form
$$\omega=\frac{1+\text{\bf{O}}(\xx{p})+\text{\bf{O}}(\yy{p})}
              {\partial R(\xx{p},\yy{p})/
                      \partial \yy{p}}d\xx{p}$$
is regular and the deformation of any extension of this regular form
to a holomorphic vector field of the moduli space is locally
given by the equation
$$R(\xx{p},\yy{p})+t(1+\text{\bf{O}}(\xx{p})
+\text{\bf{O}}(\yy{p}))+\text{\bf{O}}(t^2)=0.$$
The vector field may be chosen to be invariant under the involutions
$\sigma$ and $\rho$. In this case the whole family corresponding to
real $t$ is invariant under these involutions. Fixed points of $\eta$
correspond in our coordinates to elements, whose coordinates
$\xx{p}$ and $\yy{p}$ are imaginary. In general the real
part of $k$ has to take values in $\lattice\dual/2$ at these
double points of the form $(y,\eta(y))$. For simplicity
we assumed that $k$ vanishes at the double point.
Now it is easy to see that the involution $\eta$ does
not have fixed points on the part of the family corresponding to
$t\leq 0$. We conclude that any \Em{admissible} form may have poles
of first--order at ordinary double points of the form $(y,\eta(y))$,
but the residue has to be purely imaginary and the sign of the
imaginary part has to be the same as the sign of the imaginary part of
$d\yy{p}/d\xx{p}$ at these points.
\item[Higher oder singularities of the form $(y,\eta(y)$.]
There are two possibilities:
either the function $d\yy{p}/d\xx{p}$ takes the same value at
these two points, or both differentials $d\xx{p}$ and $d\yy{p}$
vanish at these two points. Let us first consider the case that for
some choice of generators of $\lattice$ the differential $d\xx{p}$ does
not vanish at both points of $(y,\eta(y))$. In this case the
first Taylor coefficient $a_1$ introduced above is real. If we add to
$\yy{p}(\xx{p})$ some polynomial of $\xx{p}$,
which takes for imaginary--valued $\xx{p}$ imaginary values,
we may achieve that the two Taylor series are given by
$\yy{p}=a_l\xx{p}^l+a_{l+1}\xx{p}^{l+1}+\ldots$ and
$\yy{p}=-\Bar{a}_l(-\xx{p})^l-a_{l+1}(-\xx{p})^{l+1}+\ldots$,
where $l$ is the lowest index with $a_l+(-1)^l\Bar{a}_l\neq 0$.
Now it is easy to see that the singular
part of of an \Em{admissible} regular form at this singularity has to
fulfill some positivity condition. The same is true in the other cases
of higher singularities.
\item[Singularities, which are disjoint from their image under $\eta$.]
In these cases any regular form $\omega$ fulfilling
$\sigma^{\ast}\omega=\omega$ and $\eta^{\ast}\omega=\Bar{\omega}$ is
\Em{admissible}.
\end{description}

If a
\De{Generalized Weierstra{\ss} curve}~\ref{generalized curves}
admits a local flow decreasing the \Em{first integral}, it cannot
correspond to a relative minimizer of the
\De{Generalized Willmore functional}~\ref{generalized functional}.
Therefore, we shall try to characterize all those
situations, where such flows do not exist. The main task is to
determine all vector fields, whose flows preserves the
\Em{Weak singularity condition}~\ref{weak singularity condition}.
Again we use the local description of some one--sheeted coverings of
the \Em{complex Fermi curves} by an equation
$R(\xx{p},\yy{p})=0$. In order to simplify notations, we will
always assume that the function $k$ (and therefore also the
components $\xx{p}$ and $\yy{p}$) is zero at the singularity
described in the
\De{Weak singularity condition}~\ref{weak singularity condition}.
Since the elements of the form $[\kappa/2]$ are fixed points of
the two involutions $\eta$ and $\sigma$, this simplification is
compatible with these two involutions.
Moreover, in some situations it will be more convenient to choose
the second component $\yy{p}$ to be a complex component
$\yy{p}=g(\yy{\mu},k)$ of $k$, like at the end of
Section~\ref{subsection spectral projections}. For the construction
of the complex manifold of deformations of
\Em{complex Fermi curves} in Proposition~\ref{smooth moduli}
we covered $\Spa{Y}$ by some open smooth set
and small neighbourhoods of the singularities. In order to apply
this construction to deformations of \Em{complex Fermi curves}
fulfilling the
\De{Weak singularity condition}~\ref{weak singularity condition}
we modify the representation of $\Spa{Y}$ in a small neighbourhood of this
singularity. The function $\xx{p}$ depends on the choice of
the generators of the lattice. Obviously there exists some choice,
such that $d\xx{p}$ has no zero at all \Em{multiple points} in the
preimage of $[\kappa/2]$ described in the
\De{Weak singularity condition}~\ref{weak singularity condition}.
At each \Em{cusp} of this preimage, the differentials of
almost all components of $k$ have an zero of some order
$d-1\in\mathbb{N}$. The differential of a specific
complex linear combination of $k_1$ and $k_2$ has a zero of
higher--order $m-1$. Moreover, the natural numbers $d<m$ have to be
co--prime, because otherwise the preimage of the normalization of
this point contains several elements. Obviously we may also assume
that for all \Em{cusps} in the preimage of $[\kappa/2]$ the
differential $d\xx{p}$ has a zero of minimal--order, compared with
all other components of $k$. This implies that the function
$\yy{p}/\xx{p}$ extends to a holomorphic function on some
neighbourhood of the preimage of $[\kappa/2]$. Now we describe $\Spa{Y}$
in some small neighbourhood of this preimage, by some equation of
the form $R_{\sigma}(\xx{p}^2,\yy{p}/\xx{p})=0$. Again we may
choose this function to be a polynomial with respect to
$\yy{p}/\xx{p}$, whose degree is equal to the
local number of sheets of the covering map $\xx{p}$, and whose
highest coefficient is equal to $1$. Due to this choice
the other coefficients are unique holomorphic functions
depending on $\xx{p}$. If the small neighbourhoods are invariant
under $\sigma$, this equation depends only on $\xx{p}^2$. In fact, in
this case both functions $\yy{p}/\xx{p}$ and $\xx{p}^2$ are
invariant under $\sigma$, and this equation yields a local description
of $\Spa{Y}_{\sigma}$. If we consider $\yy{p}$ as a function of $\xx{p}$
and the moduli (as we did in Proposition~\ref{smooth moduli}) the
regular form $\omega$ describing a local deformation is equal to
$$\omega=\frac{\partial \yy{p}}{\partial t}d\xx{p}=
1/2
\frac{\partial (\yy{p}/\xx{p})}{\partial t}d(\xx{p}^2)=
-\frac{\partial R_{\sigma}(\xx{p}^2,\yy{p}/\xx{p},t)/\partial t}
      {2\partial R_{\sigma}(\xx{p}^2,\yy{p}/\xx{p},t)/
        \partial (\yy{p}/\xx{p})}
d(\xx{p}^2).$$
Since the regular forms $\omega$ of $\Tilde{\Spa{Y}}$ fulfilling the
equation $\sigma^{\ast}\omega=\omega$ are the pullbacks of the
regular forms of the quotient $\Tilde{\Spa{Y}}_{\sigma}$
of $\Tilde{\Spa{Y}}$ modulo $\sigma$ under the
canonical map $\Tilde{\Spa{Y}}\rightarrow\Tilde{\Spa{Y}}_{\sigma}$,
we conclude with the help of Lemma~\ref{dualizing sheaf}
that if the form $\omega$
is locally regular with respect to the description of
$\Tilde{\Spa{Y}}_{\sigma}$ by the equation
$R_{\sigma}(\xx{p}^2,\yy{p}/\xx{p})=0$, then it describes a
deformation of \Em{complex Fermi curves}, on which the function
$\yy{p}/\xx{p}$ extends locally near the preimage of
$[\kappa/2]$ to a holomorphic function. In particular, all these
\Em{complex Fermi curves} obey the
\Em{Weak singularity condition}~\ref{weak singularity condition}.
Hence we may modify the construction of the complex manifold of
deformations in Proposition~\ref{smooth moduli} and obtain a
complex submanifold of deformations of \Em{complex Fermi curves},
which are invariant under the involution $\sigma$, and obey the
\De{Weak singularity condition}~\ref{weak singularity condition}.

Therefore, due to the considerations above, all
\Em{admissible} regular forms of this modified description of
$\Tilde{\Spa{Y}}_{\sigma}$ correspond to deformations of
\De{Generalized Weierstra{\ss} curves}~\ref{generalized curves}.
As a consequence of this modification of the description of
$\Tilde{\Spa{Y}}_{\sigma}$, the \Em{complex Fermi curves} have to be
described in some neighbourhood of $[\kappa/2]$ as
complete intersections in $\mathbb{C}^2$, which are invariant under
$\sigma$. Therefore, over \Em{Orbits of type~1}
the involution $\eta$ has fixed points on $\Tilde{\Spa{Y}}$. But over
\Em{Orbits of type~2} and \Em{Orbits of type~3} we may choose $\Spa{Y}$ in
such a way that $\eta$ does not have fixed points.
In particular, in the situations \Em{Multiple points O2/O2},
\Em{Multiple points O3}, \Em{Cusps O2} and \Em{Cusps O3} all
regular forms of this choice of $\Tilde{\Spa{Y}}_{\sigma}$ correspond to
deformations of
\De{Generalized Weierstra{\ss} curves}~\ref{generalized curves}.
If the corresponding $\Spa{Y}_{\sigma}$ is not disconnected
(compare with Remark~\ref{disconnected moduli}), then there exists
such deformations with decreasing
\De{Generalized Willmore functional}~\ref{generalized functional}.

Before we consider the remaining situations let us explain the
structure of these deformations.
Since the pullbacks of all local holomorphic forms of the
normalization of the quotient modulo $\sigma$ under the canonical
two--sheeted covering from the normalization onto the
normalization of the quotient modulo $\sigma$ have a zero at all
fixed points of $\sigma$, the deformations corresponding to such
forms preserve fixed points of $\sigma$ and therefore
\Em{Orbits of type~2}.
An easy calculation shows that in general these four situations are
deformed into two (and sometimes even more) \Em{Orbits of type~2}.
In fact, a double point of the
form $(y,\sigma(y))$ should be considered as the limit of two
colliding branch points (which are fixed points of $\sigma$)
of this two--sheeted covering. At these singularities the functions
$\xx{p}$ and $\yy{p}$ and their first derivatives take the same
value. Therefore, the
\De{Local contributions to the arithmetic genus}~\ref{local contribution}
this singularity is at least equal to $2$. The deformation into
a double point containing two fixed points of $\sigma$ changes the
\Em{geometric genus} of the \Em{complex Fermi curve}, but
does not change the \Em{geometric genus} of the quotient
modulo $\sigma$. Consequently, the local contribution of these
double points containing two fixed points of $\sigma$ to the
arithmetic genus is at least equal to $1$.

In the remaining cases \Em{Multiple points O1/O1},
\Em{Multiple points O1/O2} and \Em{Cusps O1} it suffices to
determine the subspace of \Em{admissible} regular forms of this
modified quotient $\Tilde{\Spa{Y}}_{\sigma}$ modulo $\sigma$. We observe
that in these situations the normalization is connected. In fact,
if the normalization is disconnected, then the involution $\rho$
interchanges both infinities $\infty^{\pm}$ and therefore also both
connected components. Therefore, in this case the real part is empty,
which contradicts the existence of an \Em{Orbit of type~1}.
In the following discussion of these situations
we shall derive sufficient conditions on the real regular forms
of $\Tilde{\Spa{Y}}_{\sigma}$ to be \Em{admissible}.
We will see that if a regular form takes special values
at the preimage of $[\kappa/2]$ in the real part
of the normalization, then it is \Em{admissible}.

More precisely, let us use the description of the
\Em{complex Fermi curves} by equations of the form $R(\xx{p},\yy{p})=0$.
Here $R$ may be chosen to be a polynomial with respect to $\yy{p}$.
If its degree is equal to the number of those sheets
of the covering map $\xx{p}$,
which contain one of the elements of $\Spa{Y}$ under consideration,
and if the highest coefficient is equal to $1$, then $R$ is unique.
Consequently the other coefficients are holomorphic functions
depending on $\xx{p}$.
Due to Lemma~\ref{dualizing sheaf} the regular forms are locally of
the form
$$\frac{f(\xx{p},\yy{p})}
       {\partial R(\xx{p},\yy{p})/\partial \yy{p}}d\xx{p}.$$
Here again $f$ may be chosen
to be a polynomial with respect to $\yy{p}$,
whose degree is smaller than the degree of $R$.
Due to the definition of the regular form
corresponding to a deformation, $-f$ is equal to the derivative of $R$
with respect to the deformation parameter of the corresponding family.
The Taylor coefficients of $R$ define locally holomorphic functions
on the moduli space $\moduli_{g,\lattice}$
constructed in Proposition~\ref{smooth moduli}.
We will obtain some conditions on the derivatives
of these Taylor coefficients with respect to the moduli,
which are equal to the corresponding Taylor coefficients of $f$.
Then we apply the following Lemma,
which is a consequence of S\'{e}rre duality
(\cite[\S18.2]{Fo} and \cite[Chapter~IV \S3.10.]{Se}).

\begin{Lemma}\label{regular forms}
Let $\Tilde{\Spa{Y}}_{\sigma}$ be the quotient modulo $\sigma$ corresponding
to some $\fermi(\Spa{Y},\infty^-,\infty^+,k)\in
\Bar{\moduli}_{g,\lattice,\eta,\sigma}$, and let $y_1,\ldots,y_l$ be
finitely many points of $\Tilde{\Spa{Y}}_{\sigma}$ in the complement
of the double point $(\infty^-,\infty^+)$. The values of all
real regular forms of $\Tilde{\Spa{Y}}_{\sigma}$, whose residue at
$\infty^+$ is equal to $-\sqrt{-1}$ is the orthogonal complement
of all meromorphic functions of $\Tilde{\Spa{Y}}_{\sigma}$,
which are holomorphic on
$\Tilde{\Spa{Y}}_{\sigma}\setminus\{(\infty^-,\infty^+),y_1,\ldots,y_l\}$
and bounded on some small neighbourhood of $(\infty^-,\infty^+)$ with
respect to the pairing given by the total residue.\qed
\end{Lemma}

In particular, if $\Spa{Y}$ corresponds to a
relative minimizer, then it must not contain an \Em{admissible}
regular form, whose flow decreases the
\De{Generalized Willmore functional}~\ref{generalized functional}.
Hence, due to the lemma, it has a meromorphic function with poles
contained in the preimage of $[\kappa]/2]$ described in the
\De{Weak Singularity condition}~\ref{weak singularity condition}
under the mapping $k$.
In Proposition~\ref{global meromorphic function} and
Theorem~\ref{constrained function} we construct this function
with the help of some special tools of the integrable systems having
Lax operators.

\begin{Corollary}\label{implicit function}
The exterior derivative of a linear combination
of Taylor coefficients of $R$ at a singularity
corresponding to a regular form, which has poles at all points
of the preimage in the normalization of this singularity
cannot vanish.
\end{Corollary}

\begin{proof}
The only relation on the singular parts of meromorphic differentials
on smooth Riemann surfaces $\Spa{Y}$ is that the sum of the residues vanishes.
Due to Lemma~\ref{dualizing sheaf} this relation
is even locally satisfied for all regular forms of $\Spa{Y}$,
considered as meromorphic differentials on
the normalization of $\Spa{Y}$. Hence the exterior derivative
of such a Taylor coefficient considered as a holomorphic function
on the moduli space cannot vanish.
\end{proof}

This corollary implies that in the remaining situations
\Em{Multiple points O1/O2}, \Em{Multiple points O1/O1} and
\Em{Cusps O1} the preimage of $[\kappa/2]$
contains no other \Em{Orbit of type~1}. In fact, otherwise
we choose a one--sheeted covering $\Spa{Y}$ of the
\Em{complex Fermi curve} with no fixed points of $\eta$,
and two points, which contain all \Em{Orbits of type~1}.
Now we may apply the corollary to all coefficients
of order less than $2$ in the situation of
\Em{Multiple points O1/O1/O2},
and to all coefficients of order less than $3$
in the situations \Em{Multiple points O1/O1/O1},
\Em{Cusps O1/Multiple points O1} and \Em{Cusps O1/O1}.

Let us now consider these situations.
\begin{description}
\item[\Em{Multiple points O1/O2}.]
\index{multiple point}
\index{orbit!of type~1}
\index{orbit!of type~2}
Due to our general consideration made above we have only to consider
the \Em{Orbit of type~1}. The Taylor expansions
of the function $\yy{p}$ in terms of $\xx{p}$ at the two points
of the \Em{Orbit of type~1} has to be of the form
$\yy{p}=a_1\xx{p}+a_2\xx{p}^2+a_3\xx{p}^3+\ldots$ and
$\yy{p}=a_1\xx{p}-a_2\xx{p}^2+a_3\xx{p}^3-\ldots$,
respectively. Let $2m$ be the smallest even index,
whose coefficient $a_{2m}$ does not vanish. After adding to the
function $\yy{p}(\xx{p})$ some odd function with respect to
$\xx{p}$, we may
achieve that all coefficient with odd index vanish. We conclude
that the quotient modulo $\sigma$ may be described locally by some
equation of the form
$(\yy{p}/\xx{p})^2=a_{2m}^2(\xx{p}^2)^{2m-1}
+2a_{2m}a_{2m+2}(\xx{p}^2)^{2m}+\ldots$.
If we choose a regular form $\omega$, which is locally given by
$$\omega=(2\yy{p}/\xx{p})^{-1}(1+\text{\bf{O}}(\xx{p}^2))d\xx{p}^2=
\pm a_{2m}^{-1}(\xx{p}^2)^{1-m}(1+\text{\bf{O}}(\xx{p}^2))d\xx{p},$$
and extend this form to
some holomorphic section of the tangent bundle of the
moduli space on some neighbourhood of $\Spa{Y}$, then the
corresponding family of \Em{complex Fermi curves} may be described
by the equation
$$(\yy{p}/\xx{p})^2=a_{2m}^2(\xx{p}^2)^{2m-1}+
2a_{2m}a_{2m+2}(\xx{p}^2)^{2m}
+\ldots+t(1+\text{\bf{O}}(\xx{p}^2))+\text{\bf{O}}(t^2).$$
This shows that for small non--zero $t$ these
\Em{complex Fermi curves} have double points of the form $(y,y')$,
where both points are fixed points of $\sigma$. The corresponding form
has a pole of order $m-1$ on the normalization of $\Spa{Y}_{\sigma}$ 
over the \Em{Orbit of type~1}.

If we include the involution $\eta$ into our considerations, we see that
on the normalizations of the corresponding \Em{complex Fermi curves}
the involution $\eta$ may have fixed points. More precisely, all the
coefficients $a_1,a_2,\ldots$ are real, if the corresponding
\Em{multiple points} belong to the real part. Then the last formula
shows that at least for small $t$ on the part of the
deformation family corresponding to $t<0$
the involution $\eta$ does not have fixed points. The corresponding
form on the normalization of the quotient modulo $\sigma$ has a pole
of order $l-1$. With the help of Lemma~\ref{regular forms}
we conclude that for a relative minimizer
$m$ is equal to $1$, and the \Em{complex Fermi curve}
is \Em{real--$\sigma$--hyperelliptic}.
\item[\Em{Cusps O1}.]
\index{cusp}
\index{orbit!of type~1}
We consider first a double point of the form $(y,\sigma(y))$
containing two \Em{cusps} identified
to one higher order singularity (the
\De{Local contributions to the arithmetic genus}~\ref{local contribution}
of these singularities is at least equal to $8$).
In this case for almost all lattice vectors $\gamma$ the function
$g(\gamma,k)$ has a zero of order $d$, and for some element
$\mu\in \mathbb{C}^2$ the function $g(\mu,k)$
has a zero of order $m>d$. These two numbers has to be co--prime.
Otherwise the normalization would contain several points over these
\Em{cusps}. Moreover, due to Corollary~\ref{implicit function}
for a relative minimizer $d$ has to be equal to $2$.
In fact, otherwise we choose a one--sheeted covering $\Spa{Y}$, on which
$\eta$ does not have fixed points. An application of this corollary
to all coefficients of $R(\xx{p},\yy{p})$ of order less than $3$
yields a deformation with decreasing \Em{first integral}.
In order to simplify notation we may assume that
$\xx{p}$ has a zero of order $2$ and $\yy{p}$ has a zero of
order $m$ at the points $(y,\sigma(y))$ of the normalization.
This implies that the normalization is locally nearby both points of
$(y,\sigma(y))$ a two--sheeted covering with respect to $\xx{p}$.
They obey locally some equation
$\yy{p}^2=a_1\xx{p}^m+$ higher--order terms,
and $\yy{p}^2=-a_1\xx{p}^m+$ higher--order terms, respectively.
This implies
$$R_{\sigma}(\xx{p}^2,\yy{p}/\xx{p})=
\left(\yy{p}/\xx{p}\right)^4-a_1^2\left(\xx{p}^2\right)^{m-2}+
\text{ higher--order terms }=0.$$
If we choose a regular form $\omega$, which is locally given by
$$\omega=\frac{1+\text{\bf{O}}(\xx{p}^2)+\text{\bf{O}}(\yy{p}/\xx{p})}
              {\partial R_{\sigma}(\xx{p}^2,\yy{p}/\xx{p})/
                  \partial (\yy{p}/\xx{p})}d(\xx{p}^2),$$
the flow corresponding to an appropriate extension of this form to a
holomorphic vector field of the moduli space, will deform the
equation above locally to an equation of the form
$$R_{\sigma}(\xx{p}^2,\yy{p}/\xx{p})
+t(1+\text{\bf{O}}(\xx{p}^2)+\text{\bf{O}}(\yy{p}/\xx{p}))
+\text{\bf{O}}(t^2)=0.$$
We conclude that the multiple point $(y,\sigma(y))$ is deformed into
four fixed points of $\sigma$. An easy calculation shows that this
regular form $\omega$ may have a pole at the multiple points of order
$3(m-3)$. This order is zero, if and only if $m=3$.

If we include $\eta$ into our considerations, then we see again that
only on half of the deformation family the involution $\eta$ does
not have fixed points. Therefore, Lemma~\ref{regular forms}
implies again that $m$ is equal to $3$ and that
the \Em{complex Fermi curve} is \Em{real-$\sigma$--hyperelliptic}.
\item[\Em{Multiple points O1/O1}.]
\index{multiple point}
\index{orbit!of type~1}
If we choose $\Tilde{\Spa{Y}}_{\sigma}$ in such a way that it contains two
different singular points over the two \Em{Orbits of type~1}, we may
apply the arguments used in \Em{Multiple points O1/O2}
concerning the \Em{Orbit of type~1}. This implies that one of the
following cases occurs:
\begin{description}
\item[Lowest order singularity.] The $m$'s corresponding to both
  \Em{Orbits of type~1} are equal to $1$ and the quotient
  $\Spa{Y}_{\sigma}$ of $\Spa{Y}$ modulo $\sigma$ is hyperelliptic.
  Moreover, if $f$ denotes the real function on $\Spa{Y}_{\sigma}$,
  which is invariant under the hyperelliptic involution,
  then the values of the function
  $d\left(d\yy{p}/d\xx{p}\right)/df$ at the two
  \Em{Orbits of type~1} have the same sign.
\item[Higher order singularity.] One of the $m$'s is equal to $1$ and $\Spa{Y}$ is
  \Em{real--$\sigma$--hyperelliptic}.
\end{description}
In order to complete the proof
it suffices to show that the case of the
\Em{Lowest singularity} satisfies the conditions of
\Em{Minimizers with dividing real parts}.

We claim that, if the function $d\yy{p}/d\xx{p}$ takes
the same values at both \Em{Orbits of type~1},
then $\Spa{Y}$ cannot correspond to a local minimizer.
In this case we choose $\Tilde{\Spa{Y}}_{\sigma}$
to have one singular point over both \Em{Orbits of type~1}.
The analogous calculation of
$R_{\sigma}(\xx{p}^2,\yy{p}/\xx{p})$ and the
corresponding deformations of lowest--order show
that in this case there exists \Em{admissible} regular forms
of the modification of $\Tilde{\Spa{Y}}_{\sigma}$
with prescribed singularity at the two \Em{Orbits of type~1}.
Hence there exists deformations decreasing the
\De{Generalized Willmore functional}~\ref{generalized functional}.
We remark that in case of the \Em{Lowest singularity} the
\De{Local contributions to the arithmetic genus}~\ref{local contribution}
is not smaller than $8$. If $d\yy{p}/d\xx{p}$ takes the same value
at both \Em{Orbits of type~1}, then the
\De{Local contributions to the arithmetic genus}~\ref{local contribution}
is not smaller than $12$. In Lemma~\ref{type 1 and 2}~(ii)
we will see that in this case the \Em{isospectral set} contains only
\De{Weierstra{\ss} potentials}~\ref{Weierstrass potentials}.

Now we show that in case of the \Em{Lowest singularity}
the hyperelliptic involution of $\Spa{Y}_{\sigma}$ does not have
real fixed points. Let us assume on the contrary that the
hyperelliptic involution of the quotient modulo $\sigma$ has
real fixed points. This
implies that on all connected components of the real part of the
quotient modulo $\sigma$, the form $df$ has exactly two zeroes, and
furthermore $f$ maps each connected component of the real part onto
some real interval and besides the two endpoints of the interval (which
have to be fixed points of the hyperelliptic involution and therefore
zeroes of $df$) the preimage of each element of this interval
contains exactly two points, at which the form $df$ has different signs
(with respect to any orientation of the real part). This implies that
one real cycle of the quotient modulo $\sigma$ contains the two points
corresponding to the two \Em{Orbits of type~1}, and that the form
$d(d\yy{p}/d\xx{p})$ has different signs at these two points
(with respect to any orientation of this real cycle).
The preimage of this real cycle under the natural map from the
\Em{complex Fermi curve} onto the quotient
modulo $\sigma$ consists of either one or two real cycles.

\noindent
\begin{minipage}[b]{10cm}
Let us now consider
the images of these real cycles under the map $k$.
These images are one or two closed curves in $\torus\dual$, which are
invariant under the involution $k\mapsto -k$. Moreover, if we identify
each pair of points, which are interchanged by this involution, then
we obtain one connected component. Finally, nearby the four points
corresponding to the two \Em{Orbits of type~1} the images has to look
like four paths, which all intersect in one point in two different
tangent lines. More precisely, to each path there exists another path,
which touches the path in this point. In the neighbouring figure
we indicate theses images nearby the two \Em{Orbits of type~1}
and denote the open ends in the anti--clockwise--order
by $A$, $B$, $C$, $D$, $E$, $F$, $G$ and $H$. 
\end{minipage}\hfill
{\setlength{\unitlength}{25pt}
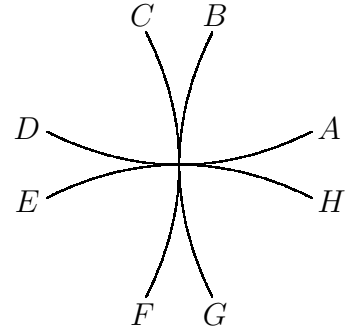
\begin{picture}(6,6)
\put(1,1){\begin{picture}(4,4)
\qbezier(0,1.5)(2,2.5)(4,1.5)
\qbezier(0,2.5)(2,1.5)(4,2.5)
\qbezier(2.5,4)(1.5,2)(2.5,0)
\qbezier(1.5,4)(2.5,2)(1.5,0)
\end{picture} }
\put(5.1,3.3){\makebox(.5,.5)[l]{$A$}}
\put(3.3,5.1){\makebox(.5,.5)[b]{$B$}}
\put(2.2,5.1){\makebox(.5,.5)[b]{$C$}}
\put(.4,3.3){\makebox(.5,.5)[r]{$D$}}
\put(.4,2.2){\makebox(.5,.5)[r]{$E$}}
\put(2.2,.4){\makebox(.5,.5)[t]{$F$}}
\put(3.3,.4){\makebox(.5,.5)[t]{$G$}}
\put(5.1,2.2){\makebox(.5,.5)[l]{$H$}}
\end{picture}}

\noindent
We conclude that there are exactly four possibilities to connect the
open ends. We denote these possibilities by combining those capitals
into common brackets, which correspond to connected open ends. In fact, an open
end of any path has to be connected to an open end of another path,
which does not have the same tangent. Also, due to the condition that
$d(d\yy{p}/d\xx{p})$ has different signs at the two
\Em{Orbits of type~1},
each open end has to be connected with some other open end, which in
our schematic picture may be obtained from the first open end by some
rotation through $\pm 90$ degrees. Therefore, we have the
following possibilities:
\begin{description}
\item[1.] $(A,C),(E,G),(B,D),(F,H)$
\item[2.] $(A,C),(E,G),(B,H),(F,D)$
\item[3.] $(A,G),(E,C),(B,D),(F,H)$
\item[4.] $(A,G),(E,C),(B,H),(F,D)$
\end{description}
We should remark that the integral of the form $dk$ along a way
connecting these open ends does not have to vanish, but due to the
invariance under the involution $k\mapsto -k$, the integral has to be
minus the integral of the corresponding image under this
involution. Let $0$ denote the image of the two \Em{Orbits of type~1}
under the map $k$. We conclude that in the first two cases the
sequence of open ends denoted by $\ldots-0-A-C-0-E-G-0-\ldots$ and in
the third and in the fourth case the sequence of open ends denoted by
$\ldots-0-A-G-0-E-C-0-\ldots$ describes a closed curve in
$\mathbb{R}^2$, which is defined up to translation by some
element of $\lattice\dual$. On the other hand, in the first and the
third case the sequence of open ends denoted by
$\ldots-0-B-D-0-F-H-0-\ldots$ and in the seconds and in the fourth
case the sequence of open ends denoted by
$\ldots-0-B-H-0-F-D-0-\ldots$ also describes a closed curve in
$\mathbb{R}^2$, which again is defined up to translation by some
element of $\lattice\dual$. Moreover, since at all
real double points the function $f$ has to take the same value,
these two closed curves may have selfintersections, but they must not
intersect each other and the corresponding shifted curves outside of
the points corresponding to the two \Em{Orbits of type~1}.
In fact, if these points are identified with some other point of the
corresponding curves, then the preimage of the corresponding element
$[\kappa/2]$ under the normalization map contains
at least three \Em{Orbits of type~1}.
However, this case was excluded in the above considerations. If we
chose any orientation of these cycles, they represent some chains of
the singular homology of $\mathbb{R}^2$. Since the first homology
groups of $\mathbb{R}^2$ are trivial, the intersection numbers of
these two cycles have to be zero
\index{intersection number}
(\cite[\S73]{ST} and \cite[Chapter~VII \S4]{Do}).
An easy calculation of the local contributions
to the intersection numbers
\index{intersection number!local contribution to the $\sim$}
\index{local contribution!to the intersection number}
of the two cycles shows,
that those, which are not zero, have the same sign
(depending on the choice of the orientation).
Moreover, for some shifted representatives a local contribution is
not zero. Hence these cycles has to intersect in some other points,
which is a contradiction.
Therefore, the hyperelliptic involution of the quotient
modulo $\sigma$ does not have real fixed points, and the real part
divides the quotient into two connected components.
Moreover, the real part of the quotient modulo $\sigma$
has one connected component,
if the \Em{geometric genus} is even,
and two connected components, if the \Em{geometric genus} is odd.
Finally, the real part of the \Em{complex Fermi curve}
divides the normalization into two connected components.
\end{description}

\noindent
{\bf 2. Relative minimizers of infinite \Em{geometric genus}.}
If $\fermi(V,W)$ is a \Em{complex Fermi curve} of potentials
$V,W\in\banach{2}(\torus)$,
which is invariant under $\eta$ without fixed points,
then due to Corollary~\ref{eta invariant} all forms of the form
$\Omega_{V,W}(\var V,\var W)$,
whose complex conjugate are equal to the pullback under $\eta$,
correspond to infinitesimal deformations of \Em{complex Fermi curves}
in $\Bar{\moduli}_{\lattice,\eta}$.
More precisely, the arguments of the proof of
Theorem~\ref{l1-structure} extend to all \Em{complex Fermi curves},
which are
\begin{description}
\item[(i)] invariant under $\eta$ without fixed points,
\item[(ii)] and have a finite value of the \Em{first integral}.
\end{description}
Consequently, the correspondence between regular one--forms and
one--dimensional families of deformations of
\Em{complex Fermi curves} obeying the
\De{Weak Singularity condition}~\ref{weak singularity condition}
extends from $\bigcup\limits_{g\in\mathbb{N}}
\Bar{\moduli}_{g,\lattice,\eta,\sigma}$ to
$\bigcup\limits_{\willmore>0}
\Bar{\moduli}_{\lattice,\eta,\sigma,\willmore}$.
Therefore, Proposition~\ref{global meromorphic function} shows
that the same conclusion holds in the general case
(compare with Remark~\ref{extension of the symplectic form}).
In particular, all relative minimizers
of the \De{Generalized Willmore functional}~\ref{generalized functional},
whose Willmore functional is smaller than $\willmore$,
belong to $\Bar{\moduli}_{g_{\max},\lattice,\eta,\sigma,\willmore}$,
with the corresponding integer $g_{\max}$ introduced in
Corollary~\ref{dividing genus bound}.
\end{proof}

At the end of this section we want to determine
all possible real double points of a relative minimizer.
In the situations \Em{Multiple points O1/O2}, \Em{Cusps O1}
and a \Em{Higher singularity} of \Em{Multiple points O1/O1}
the function $f$ on $\Spa{Y}_{\sigma}$ of degree one
has to take the same values at the double points.
Therefore, they are double points of the form $(y,\sigma(y))$.

In the situation of two \Em{Lowest singularities} of
\Em{Multiple points O1/O1} we know that in the notation
introduced above we have the following possibilities to connect the
open ends nearby the two \Em{Orbits of type~1}:

\noindent
\begin{minipage}[b]{8cm}
\begin{description}
\item[1.] $(A,B),(E,F),(C,D),(G,H)$
\item[2.] $(A,B),(E,F),(C,H),(G,D)$
\item[3.] $(A,F),(E,B),(C,D),(G,H)$
\item[4.] $(A,F),(E,B),(C,H),(G,D)$
\item[5.] $(A,D),(E,H),(B,C),(F,G)$
\item[6.] $(A,D),(E,H),(B,G),(F,C)$
\item[7.] $(A,H),(E,D),(B,C),(F,G)$
\item[8.] $(A,H),(E,D),(B,G),(F,C)$
\end{description}
\end{minipage}
{\setlength{\unitlength}{25pt}
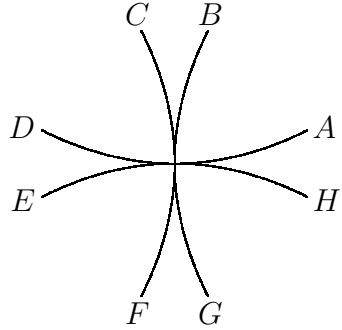
\begin{picture}(6,6)
\put(1,1){\begin{picture}(4,4)
\qbezier(0,1.5)(2,2.5)(4,1.5)
\qbezier(0,2.5)(2,1.5)(4,2.5)
\qbezier(2.5,4)(1.5,2)(2.5,0)
\qbezier(1.5,4)(2.5,2)(1.5,0)
\end{picture} }
\put(5.1,3.3){\makebox(.5,.5)[l]{$A$}}
\put(3.3,5.1){\makebox(.5,.5)[b]{$B$}}
\put(2.2,5.1){\makebox(.5,.5)[b]{$C$}}
\put(.4,3.3){\makebox(.5,.5)[r]{$D$}}
\put(.4,2.2){\makebox(.5,.5)[r]{$E$}}
\put(2.2,.4){\makebox(.5,.5)[t]{$F$}}
\put(3.3,.4){\makebox(.5,.5)[t]{$G$}}
\put(5.1,2.2){\makebox(.5,.5)[l]{$H$}}
\end{picture}}

\noindent
Again we use the cycles corresponding to the first two pairs
and the last two pairs (for example in the first case
the cycles denoted by $\ldots-0-A-B-0-E-F-0-\ldots$
and $\ldots-0-C-D-0-G-H-0\ldots$).
Again they correspond to closed cycles in $\mathbb{R}^2$
and their intersection numbers has to vanish.
A direct calculation shows that the local contribution
to the intersection numbers at the two \Em{Orbits of type~2}
cancels. Therefore, the number of points
with non--zero contribution to these intersection numbers
has to be even. Furthermore, due to the invariance under $\sigma$,
these points have to occur in pairs with equal contribution to the
intersection numbers. Finally, due to Lemma~\ref{regular forms},
the function $f$ has to take the same value at
any two points of the normalization, which belong to one singularity.
Hence the pairs of such double points with non--zero contribution to
the intersection number occur in different pairs,
where $f$ takes different values,
or four double points, where $f$ takes only one value.
Now we claim that the first case is impossible.
In fact, if we choose two of these pairs with opposite local
contribution to the intersection number,
we may use these double points in order to connect
the cycles nearby the two \Em{Orbits of type~1}
in such a way, that the differential $df$ has different sign
at the two \Em{Orbits of type~1}
with respect to any orientation of these cycles.
The considerations concerning this case in the proof of
Theorem~\ref{relative minimizers} implies that then
we have a double point containing two real points
of the normalization, where $f$ takes different values.
This contradicts Lemma~\ref{regular forms}.
In the second case we have again two \Em{Orbits of type~1}
over some other element of $\lattice\dual/2\lattice\dual$.
In this case the same considerations apply to this singularity.

We conclude that we may deform all these \Em{Orbits of type~1}
by some arbitrary small deformation increasing
the \Em{first integral} into \Em{Orbits of type~2}.
After this deformation the real part will not
have any double point.
We claim that the function $k_1+\sqrt{-1}k_2$ maps
the connected component $\Spa{Y}^{+}$ of the complement of the real part
in the normalization, which contains $\infty^+$,
onto the relative complement in a finite--sheeted covering of the torus
$\mathbb{C}/\lattice\dual_\mathbb{C}$
\index{lattice!dual $\sim$ $\lattice\dual_\mathbb{C}$}
of finitely many disjoint sets of this torus.
Here $\lattice\dual_\mathbb{C}=
\{\kappa_1+\sqrt{-1}\kappa_2\in\mathbb{C}\mid\kappa\in\lattice\dual\}$
denotes the complex version of the dual lattice.
We shall prove the analogous stronger statement for the corresponding
quotients modulo the involution $\sigma$.
The Weierstra{\ss} function $\wp$ corresponding to the complex torus
$\mathbb{C}/\lattice\dual_\mathbb{C}$ (see e.g. \cite[13.12]{MOT})
represents this complex torus as a two--sheeted covering
over $\mathbb{P}^1$,
whose hyperelliptic involution is induced by the
involution $\sigma:k\mapsto -k$.
The four elements of $\lattice\dual/2\lattice\dual$ are the four
fixed points of this involution.
Consequently the composition of the function
$k_1+\sqrt{-1}k_2$ with $\wp$ yields a function on the quotient
$\Spa{Y}^+_{\sigma}$ of $\Spa{Y}^+$ modulo the involution $\sigma$.
Nearby the two \Em{Orbits of type~2}
(the deformed \Em{Orbits of type~1})
this function maps $\Spa{Y}^+_{\sigma}$ onto the complement of
finitely many (either one or two) disjoint regions of $\mathbb{P}^1$.
Since the real part has no double points,
and therefore does not intersect itself,
the same is true for the whole of $\Spa{Y}^+_{\sigma}$.
Consequently we obtain

\begin{Corollary}\label{deformed minima}
All relative minima may be deformed by some arbitrary small deformation
increasing the \Em{first integral} into a
\De{Generalized Weierstra{\ss} curve}~\ref{generalized curves}
obeying the following conditions:
\begin{description}
\item[(i)] The real part divides the normalization
  into two connected components
  $\Spa{Y}^-_{\text{\scriptsize\rm normal}}$ and
  $\Spa{Y}^+_{\text{\scriptsize\rm normal}}$
  containing $\infty^-$ and $\infty^+$, respectively.
\item[(ii)] The function $k_1+\sqrt{-1}k_2$ maps
  $\Spa{Y}^+_{\text{\scriptsize\rm normal}}$
  into the relative complement in a finite--sheeted covering of
  $\mathbb{C}/\lattice\dual_\mathbb{C}$
  of finitely many disjoint closed subsets
  of this complex torus.
\item[(iii)] All \Em{multiple points}, which connect the two
  components $\Spa{Y}^+_{\text{\scriptsize\rm normal}}$ and
  $\Spa{Y}^-_{\text{\scriptsize\rm normal}}$ are of the form $(y,\eta(y))$,
  where the imaginary part of the function $d\yy{p}/d\xx{p}$
  has the same sign as at the corresponding points $\infty^{\pm}$,
  respectively.
\item[(iv)] For some $[\kappa/2]\in\lattice\dual/2\lattice\dual$
  the normalization of the \Em{complex Fermi curve} contains either two
  \Em{Orbits of type~2} or one \Em{Orbit of type~3} over this
  element.\qed
\end{description}
\index{relative minimizer|)}
\index{minimizer!relative $\sim$|)}
\end{Corollary}

Corollary~\ref{dividing genus bound} suggests,
that for large $g$ all relative minimizers in
$\Bar{\moduli}_{g,\lattice,\eta,\willmore}$ fulfill the conditions of
Theorem~\ref{relative minimizers}.
A proof of such a statement could make it possible to prove the
Willmore conjecture without making use of the results of
Sections~\ref{subsection limits}--\ref{subsection compactified moduli}.
At an earlier state of this project, we thought to have a proof of
such a statement. However, we overlooked the possibility of
\Em{real cusps}
(compare with Proposition~\ref{modified smooth moduli} and 
Remark~\ref{examples of real cusps}).
We did not succeed to overcome this problem.
Afterwards we succeed to prove the compactness of the moduli spaces
without any bound on the \Em{geometric genus} of the
\Em{complex Fermi curves}.

\subsection{Absolute minimizers of the
  \Em{Generalized Willmore functional}}
\label{subsection absolute minimizer}

In Section~\ref{subsection relative minimizers} we obtained
a classification of relative minimizers of the
\De{Generalized Willmore functional}~\ref{generalized functional}.
In this section we will determine
all relative minimizers of the 
\De{Generalized Willmore functional}~\ref{generalized functional},
whose \Em{first integral} is not larger than $8\pi$. The
\De{Weak Singularity condition}~\ref{weak singularity condition}
is, roughly speaking, a condition on
four points of the \Em{complex Fermi curve}.
As a preparation we shall first determine
the minimum of the restriction of the \Em{first integral}
to the space of all \Em{complex Fermi curves},
which  contain a fixed element $k'\in\mathbb{R}^2$.

\begin{Lemma} \label{unique real point minimizers}
For all elements $[k']\in\torus\dual$
there exists a unique absolute minimizer
of the restriction of the \Em{first integral}
to the space of all \Em{complex Fermi curves} in
$\bigcup\limits_{\willmore>0}\Bar{\moduli}_{\lattice,\eta,\sigma,\willmore}$,
which contains $[k']$.
These minimizers fit together to a continuous map
\index{minimizer!$\fermi_{\min}(\cdot)$}
\index{family!$\fermi_{\min}(\cdot)$}
$$\fermi_{\min}:\torus\dual\rightarrow
\Bar{\moduli}_{3,\lattice,\eta,\sigma,4\pi},
[k]\mapsto\fermi_{\min}([k])$$
into the set of \Em{real--$\sigma$--hyperelliptic}
\Em{complex Fermi curves} in
$\Bar{\moduli}_{3,\lattice,\eta,\sigma,4\pi}$.
The restriction of the \Em{first integral}
to the image of this mapping yields a homeomorphism onto $[0,4\pi]$.
Therefore, we may parameterize this family as well by $\willmore$:
$$\fermi_{\min}:[0,4\pi]\rightarrow
\Bar{\moduli}_{3,\lattice,\eta,\sigma,4\pi},
\willmore\mapsto\fermi_{\min}(\willmore).$$
\end{Lemma}

\begin{proof}
The space of all $\fermi\in\Bar{\moduli}_{\lattice,\eta,\sigma}$,
which contain an element $[k']\in\torus\dual$ is a subvariety, which
is locally given by the zero set of at most two analytic functions.
Due to the compactness this subvariety contains a minimizer
of the \Em{first integral}.

The preimage of $[k']$ (and $[-k']$) in the normalization of
$\fermi$ contains either fixed points of $\rho$ (i.\ e.\ real points)
or pairs, which are interchanged by $\rho$.
If in the latter case the
\De{Structure sheaf}~\ref{structure sheaf}
contains functions, which take different values at a pair,
then a \Em{B\"acklund transformation}
with the corresponding values of the eigenfunctions $\psi$
transforms the potential into a potential,
whose \De{Structure sheaf}~\ref{structure sheaf} contains only
functions taking the same value at these pairs.
We remark that the formulas of Lemma~\ref{Baecklund transformation}
also apply to
\Em{Finite rank perturbations}~\ref{finite rank perturbations}.
Nearby this transformed potential the subvariety
described in the lemma is the zero set of one analytic function.
If the preimage of $[k']$ contains real points,
then the subvariety described in the Lemma
is the zero set of one analytic function.
Therefore, we may assume that nearby the minimizer the subvariety
described in the lemma is equal to the zero set of one analytic
function.

Then Proposition~\ref{global meromorphic function} shows that
that if $\fermi\in\moduli_{\lattice,\eta,\sigma}$
is a relative minimum of the restriction
of the \Em{first integral} to this subvariety,
then this \Em{complex Fermi curve} is $\sigma$--hyperelliptic
(i.\ e.\ the quotient modulo the involution $\sigma$
is isomorphic to $\mathbb{P}^1$ with two points removed).
We remark that these arguments lead to the same conclusion,
if the derivative of that single analytic function vanishes,
whose zero set is the subvariety described in the Lemma.
Moreover, in case that the preimage of $[k']$
in the normalization contains pairs,
which are interchanged by $\rho$,
then the corresponding function constructed in
Proposition~\ref{global meromorphic function} would have no poles,
and the corresponding \Em{complex Fermi curve}
has a disconnected normalization. (We shall see that this is the case
for $[k']=0$ and those $[k']$, which are the unique real elements of
the unique \Em{complex Fermi curves} with disconnected normalization,
whose \Em{first integral} is equal to $4\pi$).
If the preimage of $[k']$ in the normalization contains
a fixed point of $\rho$, then the \Em{complex Fermi curve} is
\Em{real $\sigma$--hyperelliptic}.
To sum up, the minima in $\moduli_{\lattice,\eta,\sigma}$
have dividing real parts.

These conclusions are also true for minima in
$\bigcup\limits_{\willmore>0}\Bar{\moduli}_{\lattice,\eta,\sigma,\willmore}$.
In fact the construction of $(\var V_f,\var W_f)$ in
Lemma~\ref{residue} extends to the corresponding variations
on the complement of the union $\Set{S}_{D,\varepsilon}$ of the union of
arbitrary small balls around the singular set described in
Theorem~\ref{limits of resolvents}
(compare with Remark~\ref{extension of the symplectic form}).
Since the union over $\varepsilon$ of all subspaces of
$\banach{2}(\torus)\times\banach{2}(\torus)$ of elements, which vanish
on $\Set{S}_{D,\varepsilon}$ is dense, the arguments given above
extends to all minima of the \Em{first integral}
restricted to the subvariety described in the Lemma.
Therefore, due to Corollary~\ref{dividing genus bound},
these minima have finite genus.

We shall see in Lemma~\ref{disconnected} that the \Em{first integral}
of any \Em{complex Fermi curve} with disconnected normalization is a
multiple of $4\pi$. The unique \Em{complex Fermi curve},
whose \Em{first integral} vanishes
(i.\ e.\ the free \Em{complex Fermi curve} introduced in
Section~\ref{subsection fermi curve})
contains only the real element $[k']=0$.
Obviously, the corresponding \Em{complex Fermi curves} depend
continuously on $[k']$. Therefore, for small non--zero $k'$
these \Em{real--$\sigma$--hyperelliptic} minimizers have genus zero.
These \Em{complex Fermi curves} are given by the equation
$g(k,k)=g(k',k')$. More precisely, if $g(k',k')$
is smaller than $1/4$ times the square of the length of
the shortest non--zero dual lattice vector,
this \Em{complex Fermi curve} is the minimizer and all
double points of this \Em{complex Fermi curve}
are of the form $(y,\eta(y))$.
If $g(k',k')$ is equal to $1/4$ times
the square of the length of the shortest
non--zero dual lattice element,
then the \Em{complex Fermi curve} contains
\Em{Orbits of type~1} over a finite number of elements of
$\lattice\dual/2\lattice\dual$.
If the length of $k'$ becomes larger,
the \Em{complex Fermi curve} will be deformed
into a \Em{complex Fermi curve}, with as much
\Em{Orbits of type~2}, as the former \Em{complex Fermi curve} has
\Em{Orbits of type~1}. We will see that in general there is only one
\Em{Orbit of type~2}. In this case the real part has two connected
components, and the integral of $dk$ over any component is
non--zero. The map $k_1+\sqrt{-1}k_2$ maps the closure of $\Spa{Y}^+$
onto some part of the elliptic curve
$\mathbb{C}/\lattice\dual_\mathbb{C}$.
Moreover, this map has no branch point,
and therefore may be considered as a subset of this elliptic curve.
Obviously, the image of the closure of
$\Spa{Y}^+$ in the elliptic curve covers all points $[k_1]$, whose
\Em{complex Fermi curve} $\fermi_{\min}(k_1)$ belongs to the
family of the deformation of
$\fermi_{\min}(k')$ into $\fermi_{\min}(0)$. If we move $k'$
into the complement of the
image of the closure of $\Spa{Y}^+$ under the map $k_1+\sqrt{-1}k_2$, 
then this subset will become larger
until the boundaries meet each other. In this case the corresponding
\Em{complex Fermi curve} has again some \Em{Orbits of type~1}
over some elements of $\lattice\dual/2\lattice\dual$. 
In fact, since it is defined as a relative minimizer,
all double points of the whole family have to be either
of the form $(y,\eta(y))$ or the function introduced in
Proposition~\ref{global meromorphic function}
has to take the same values at these points,
which implies that they are of the form $(y,\sigma(y))$.
Since both points are invariant under $\rho$
they fit together to an \Em{Orbit of type~1} over some element
of $\lattice\dual/2\lattice\dual$. If we enlarge the
length of $g(k',k')$ further, these \Em{Orbits of type~1} will again
be deformed into \Em{Orbits of type~2}. The maps of the closures of
the corresponding $\Spa{Y}^+$ will still define some one--sheeted covering of
the corresponding elliptic curve. But the boundary will look like
one or two circle, which are homologous to zero.
If we move $k'$ inside the part,
which corresponds to the complement of the image of $\Spa{Y}^+$ under the map
$k_1+\sqrt{-1}k_2$ in the elliptic curve, the \Em{first integral} will
be enlarged and any \Em{Orbit of type~1} of the boundary will be
deformed into an \Em{Orbit of type~2}. Finally, the whole elliptic
curve is covered, and the \Em{first integral} is equal to $4\pi$. This
shows that for all $k'$ the \Em{first integral} is not larger than
$4\pi$. Moreover, by construction this family is periodic
with respect to the dual lattice and therefore defines
some continuous map from
$\torus\dual$ into the set of $\sigma$--hyperelliptic
\Em{complex Fermi curves}, which contains the preimage of this map.
We conclude that for all $\willmore\in [0,4\pi]$ there exists one
\Em{real--$\sigma$--hyperelliptic} \Em{complex Fermi curve},
all whose double points are of the form $(y,\eta(y))$.
Obviously $\fermi_{\min}(k_1)$ coincide with $\fermi_{\min}(k_2)$,
if and only if $k_2$ belongs to the real part of
$\fermi_{\min}(k_1)$ which is equivalent that $k_1$ belongs to
the real part of $\fermi_{\min}(k_2)$. 

We now claim that this family contains all
\Em{real--$\sigma$--hyperelliptic} \Em{complex Fermi curves},
whose \Em{first integral} is not larger than $4\pi$,
and which are relative minimizers of the restrictions of the
\Em{first integral} to the \Em{complex Fermi curves}
containing any element of the real part
of the given \Em{complex Fermi curve}.
In fact, all these elements have to be
\Em{real--$\sigma$--hyperelliptic}, and they
can have only double points of the form $(y,\eta(y))$. Moreover, we
may deform any such \Em{complex Fermi curve} within this class into
the unique \Em{complex Fermi curve}, whose \Em{first integral}
is zero. Hence this family has to meet the other family at some point.
But the space of \Em{real--$\sigma$--hyperelliptic}
\Em{complex Fermi curves} is locally a one--dimensional space
with some multiple points,
where it has \Em{Orbits of type~1}. Any such \Em{Orbit of type~1}
may be either deformed into double points of the form
$(y,\eta(y))$ not belonging to the real part, or by changing the
\Em{first integral} in the other direction either into
\Em{Orbits of type~2} or into real double points not of the form
$(y,\eta(y))$. In particular, the space of
\Em{real--$\sigma$--hyperelliptic} \Em{complex Fermi curves},
all whose double points are of the form $(y,\eta(y)$, is locally a
one--dimensional manifold. Also the \Em{first integral} defines an
injective function to $\mathbb{R}^+$. This proves the claim.

We conclude that this family contains all
\Em{real--$\sigma$--hyperelliptic} \Em{complex Fermi curves},
whose \Em{first integral} is smaller than $4\pi$,
and which contain only double points of
the form $(y,\eta(y))$.
\end{proof}

In particular, if the absolute minimizers of the restriction of the
\Em{first integral} to the space of all \Em{complex Fermi curves}, which
contain an element of $\lattice\dual/2\lattice\dual$,
are of the following form:

\begin{Corollary} \label{simplified minimizers}
For all elements $[\kappa/2]\in\lattice\dual/2\lattice\dual$ there exists
a unique \Em{complex Fermi curve} $\fermi_{\min}(\kappa/2)$,
which minimizes the restriction of the \Em{first integral}
to the space of all \Em{complex Fermi curves},
which contain the element $[\kappa/2]$ ($\fermi_{\min}(0)$
is the unique \Em{complex Fermi curve}
with vanishing \Em{first integral}). All these
minimizers are $\sigma$--hyperelliptic and the \Em{geometric genus}
is at most equal to two. More precisely, the genus is equal to the
number of \Em{Orbits of type~2}, which are contained in the
normalization of $\Spa{Y}$. Each \Em{Orbit of type~2} is mapped onto
one non--zero element of $\lattice\dual/2\lattice\dual$, and over
all non--zero elements of $\lattice\dual/2\lattice\dual$ there
exists at most one \Em{Orbit of type~2}. Finally, the
\Em{first integral} of $\fermi_{\min}(\kappa/2)$ is
larger or equal to the \Em{first integral} of
$\fermi_{\min}(\kappa'/2)$, if and only if $[\kappa'/2]$
is contained in $\fermi_{\min}(\kappa/2)$. \qed
\end{Corollary}

\begin{Remark}\label{order of minima}
In the proof of Lemma~\ref{unique real point minimizers} we have
actually shown, that the \Em{first integral} of
$\fermi_{\min}(\kappa/2)$ is
larger than the \Em{first integral} of $\fermi_{\min}(\kappa'/2)$,
if and only if $\fermi_{\min}(\kappa/2)$ contains an
\Em{Orbit of type~2} over $[\kappa'/2]$. Moreover, the
\Em{first integrals} are equal, if and
only if both \Em{complex Fermi curves} coincide and contain
\Em{Orbits of type~1} over both elements.
\end{Remark}

The proof of Lemma~\ref{unique real point minimizers}
shows that there is either one element in
$\lattice\dual/2\lattice\dual$
or two different elements of
$\left(\mathbb{R}^2\setminus\lattice\dual/2\right)/\lattice\dual$, 
whose corresponding minimum of the \Em{first integral} is equal to $4\pi$.
Moreover, the proof of this lemma shows
that the real part of this \Em{complex Fermi curve} is empty.
In fact, these points are the limits of the real parts of
the family of zero energy limits with lower values of the
\Em{first integral} as described in this Lemma.
Since the real part of any \Em{real--$\sigma$--hyperelliptic}
\Em{complex Fermi curve}
divides the normalization into two connected components,
we conclude that this \Em{complex Fermi curve} has 
disconnected normalization. Let us consider
these \Em{complex Fermi curves}
with disconnected normalization in some more detail.
Due to Theorem~\ref{asymptotic analysis 1}
they have finite \Em{geometric genus}.

\begin{Lemma}\label{disconnected}
\index{normalization!disconnected}
\index{disconnected!normalization}
If $\fermi\in\bigcup\limits_{g\in\mathbb{N}}
\Bar{\moduli}_{g,\lattice,\eta}$ has a disconnected normalization,
then the component
$\Bar{\fermi}_{\text{\scriptsize\rm normal}}^+$
of the normalization of $\Bar{\fermi}$,
which contains $\infty^+$, is a finite--sheeted covering of
the unique complex torus,
which is conformally equivalent to $\torus$.
Moreover, the \Em{first integral} of $\fermi$
is $4\pi$ times the number of sheets of this covering.
Furthermore, for all conformal classes there exists a
real two--dimensional family of such \Em{complex Fermi curves}
with disconnected normalization in
$\bigcup\limits_{g\in\mathbb{N}}
\Bar{\moduli}_{g,\lattice,\eta,4\pi}$,
and only one such \Em{complex Fermi curve} in
$\bigcup\limits_{g\in\mathbb{N}}
\Bar{\moduli}_{g,\lattice,\eta,\sigma,4\pi}$.
\end{Lemma}

\begin{proof}
If $\fermi\in\bigcup\limits_{g\in\mathbb{N}}
\Bar{\moduli}_{g,\lattice,\eta}$ has disconnected normalization,
then the functions $k_1\pm\sqrt{-1}k_2$ extend to
holomorphic multi--valued functions on the connected component
$\Spa{Y}_{\text{\scriptsize\rm normal}}^{\pm}$
of the normalization of $\Spa{Y}$,
which contain $\infty^{\pm}$, respectively.
More precisely, they define holomorphic functions
into the complex torus
$\mathbb{C}/\lattice\dual_\mathbb{C}$.
Due to the Open Mapping Theorem \cite[Chapter~IV \S7. 7.5]{Co1}
these mappings are either constant or open.
Due to Lemma~\ref{gluing rule} in the first case the
\Em{complex Fermi curve} is the unique \Em{complex Fermi curve},
whose \Em{first integral} vanishes. Due to the compactness of
$\Spa{Y}_{\text{\scriptsize\rm normal}}^{\pm}$
in the second case these mappings are finitely--sheeted coverings.
Moreover, due to Lemma~\ref{gluing rule},
the \Em{first integral} is equal to
$4\pi$ times the number of sheets of
these coverings. With the help of the formula of
Lemma~\ref{Willmore functional} it can be shown that
in these cases the \Em{first integral} is equal to $4\pi$
times the intersection number of those cycles in
$H_1(\Spa{Y}_{\text{\scriptsize\rm normal}}^+,\mathbb{Z})$,
whose intersection numbers with any cycle in
$H_1(\Spa{Y}^+_{\text{\scriptsize\rm normal}},\mathbb{Z})$
is equal to the integral of $d\xx{p}$ and $d\yy{p}$ over the cycle,
respectively. Due to Poincar\'{e} duality this intersection number
is again equal to the number of sheets of these coverings.

If the \Em{first integral} is equal to $4\pi$, the normalization
$\Bar{\fermi}^+_{\text{\scriptsize\rm normal}}$ has to be
biholomorphic to the complex torus
$\mathbb{C}/\lattice\dual_\mathbb{C}$
with one marked point $0$ corresponding to $\infty^{\pm}$.
Furthermore, the one--forms
$d\xx{p}$ and $d\yy{p}$ on the complex torus
are uniquely determined by their singular part at $\infty^+$
(compare with condition \Em{Quasi--momenta}~(i) in
Section~\ref{subsection complex Fermi curves of finite genus})
and the condition that for all $\kappa\in\lattice\dual$
the integral along the corresponding cycle of
$\mathbb{C}/\lattice\dual_\mathbb{C}$ is equal to
$g(\xx{\gamma},\kappa)$ and $g(\yy{\gamma},\kappa)$, respectively.
We conclude that the function $k$ on
$\Bar{\fermi}^{\pm}_{\text{\scriptsize\rm normal}}$
are determined up to some additive constants $k^+$ and $k^-$.
Since on $\Bar{\fermi}^{\pm}_{\text{\scriptsize\rm normal}}$
the function $k_1\pm\sqrt{-1}k_2$ maps $\infty^{\pm}$ onto $0$,
and the involution $k\mapsto -\Bar{k}$ interchanges
$\Bar{\fermi}^{\pm}_{\text{\scriptsize\rm normal}}$, 
there exists only a real two--dimensional family of such
\Em{complex Fermi curves}
(compare with Proposition~\ref{smooth moduli} and
Remark~\ref{disconnected moduli}).
The invariance under the involution
$k\mapsto -k$ determines the function $k$ uniquely
up to shifts by dual lattice elements.
\end{proof}

We shall construct the unique \Em{complex Fermi curve} in
$\Bar{\moduli}_{3,\eta,\sigma,4\pi}$
with disconnected normalization explicitly.
Obviously it is enough to determine the meromorphic multi--valued
functions $\xx{p}$ and $\yy{p}$ on the component
$\Bar{\fermi}^+_{\text{\scriptsize\rm normal}}$.
In the foregoing lemma we have
shown that this Riemann surface is biholomorphic to the elliptic curve
$\mathbb{C}/\lattice\dual_\mathbb{C}$.
We realize this elliptic curve as the quotient
$\mathbb{C}/\left(2\omega\mathbb{Z}+2\omega'\mathbb{Z}\right)$,
where $2\omega=\xx{\kappa}_1+\sqrt{-1}\xx{\kappa}_2$ and
$2\omega'=\yy{\kappa}_1+\sqrt{-1}\yy{\kappa}_2$. After
some translation we may identify the point $\infty^+$ with the point
$z=0$. Therefore, the involution $\sigma$ is induced by the map
$z\mapsto -z$. The functions $\xx{p}$ and $\yy{p}$ satisfy
$\xx{p}(z+2\omega)=\xx{p}(z)+1,
\xx{p}(z+2\omega')=\xx{p}(z)$ and
$\yy{p}(z+2\omega)=\yy{p}(z),
\yy{p}(z+2\omega')=\yy{p}(z)+1$.
Since $\xx{p}$ and $\yy{p}$ have poles of first--order at
$z=0$ this uniquely determines these functions, and they may be
expressed in terms of the Weierstra{\ss} functions
(see e.g. \cite[13.12]{MOT}):
\begin{align*}
\xx{p}(z)&=\sqrt{-1}\frac{\omega'\zeta(z)-\eta' z}{\pi}&
\yy{p}(z)&=-\sqrt{-1}\frac{\omega\zeta(z)-\eta z}{\pi}.
\end{align*}
This implies that $(\xx{p},\yy{p})$ is at $\omega$
equal to $(-1/2,0)$, at $\omega'$ equal to $(0,-1/2)$
and at $\omega+\omega'$ equal to $(-1/2,-1/2)$.
Together with $0$ these are all fixed points of $\sigma$,
therefore the normalization of this \Em{complex Fermi curve}
has over all non--zero elements of
$\lattice\dual/2\lattice\dual$ exactly one
\Em{Orbit of type~2} and no other \Em{Orbit of type~2}.
This follows also from the proof of
Lemma~\ref{unique real point minimizers},
since the unique \Em{complex Fermi curve}
with disconnected normalization
whose \Em{first integral} is equal to $4\pi$,
is the same \Em{complex Fermi curve} as the unique element
of the family described in Lemma~\ref{unique real point minimizers}
whose \Em{first integral} is equal to $4\pi$. Moreover,
since $d\yy{p}-\tau d\xx{p}$ is equal to $-dz/\omega$
the integral of the form $dz$ from any element of
$\lattice\dual/2\lattice\dual$ to any other element
has to be an element of $\omega\mathbb{Z}+\omega'\mathbb{Z}$.
This shows that these \Em{Orbits of type~2}
are all elements of $\lattice\dual/2\lattice\dual$
of this \Em{complex Fermi curve}.
Therefore, there does not exists a 
\De{Generalized Weierstra{\ss} curves}~\ref{generalized curves}
with disconnected normalization, whose
\Em{first integral} is smaller than $8\pi$.
However, we may use this \Em{complex Fermi curve} in order to construct
\De{Generalized Weierstra{\ss} curves}~\ref{generalized curves}
with disconnected normalization, whose \Em{first integral}
is equal to $8\pi$.

In Lemma~\ref{gluing rule} we cut the \Em{complex Fermi curves}
into two parts,
on both of which the function $k$ is single--valued.
Each part contains one of the two marked points $\infty^{\pm}$,
and the condition
that $k_1\pm\sqrt{-1}k_2$ should vanish at $\infty^{\pm}$
uniquely determines the branch of $k$ on both parts.
The differences between these functions along the cuts
are elements of $\lattice\dual$.
Also it might happen, that these elements generate
a sublattice. Let us denote this sublattice by
\index{lattice!effective dual $\sim$
       $\lattice\dual_{\text{\scriptsize\rm effective}}$}
\index{effective!dual lattice
       $\lattice\dual_{\text{\scriptsize\rm effective}}$}
$\lattice\dual_{\text{\scriptsize\rm effective}}\subset
\lattice\dual$ and the corresponding dual lattice by
\index{lattice!effective $\sim$ $\lattice_{\text{\scriptsize\rm effective}}$}
\index{effective!lattice $\lattice_{\text{\scriptsize\rm effective}}$}
$\lattice_{\text{\scriptsize\rm effective}}\supset
\lattice$. More precisely, we define the effective period lattice
$\lattice_{\text{\scriptsize\rm effective}}$ as the dual lattice of
  $\lattice\dual_{\text{\scriptsize\rm effective}}$:
$$\lattice_{\text{\scriptsize\rm effective}}=
\left\{\gamma\in\mathbb{R}^2\mid
g(\gamma,\kappa)\in\mathbb{Z}
\;\forall\kappa\in\lattice\dual_{\text{\scriptsize\rm effective}}
\right\}.$$
If $\lattice\dual_{\text{\scriptsize\rm effective}}$ is
contained in a one--dimensional subspace of $\mathbb{R}^2$,
there exists a lattice vector $\gamma\in\lattice$, which is orthogonal
to the whole lattice
$\lattice\dual_{\text{\scriptsize\rm effective}}$.
This implies that $g(\gamma,k)$ extends to a
single--valued meromorphic function on
the normalization of the \Em{complex Fermi curve},
and the corresponding flow acts trivially on all
line bundles of the normalization
of this \Em{complex Fermi curves}.
In this case $\lattice_{\text{\scriptsize\rm effective}}$
is by definition no longer any lattice.
However, the intersection of this group with the
one--dimensional subspace generated by
$\lattice\dual_{\text{\scriptsize\rm effective}}$,
is a one--dimensional lattice. 
In the most degenerated case
$\lattice\dual_{\text{\scriptsize\rm effective}}$
is equal to $\{0\}$,
and $\lattice_{\text{\scriptsize\rm effective}}$
is equal to $\mathbb{R}^2$.
This implies that both components of $k$ extend to
single--valued meromorphic functions
on the normalization of the \Em{complex Fermi curve},
and both flows act trivially on the line bundles of
the normalization of the \Em{complex Fermi curve}.
In this case the \Em{geometric genus} is zero
and the corresponding potentials of the form $(U,\Bar{U})$ have
only one non--vanishing Fourier coefficient
(compare with the motivation of Theorem~\ref{l1-structure}).
In particular, the corresponding real potentials are constant.
More generally, if the normalization of a given
\Em{complex Fermi curve} is connected,
then due to Lemma~\ref{existence of potentials},
there exists a complex potential $U$,
which is invariant under translations of all elements of
$\lattice_{\text{\scriptsize\rm effective}}$,
and whose \Em{complex Fermi curve} is equal to the given one.

\begin{Remark}\label{effective lattice}
The integral of $dk$ over cycles defines a map from
$H_1(\Spa{Y}_{\text{\scriptsize\rm normal}}, \mathbb{Z})$
into $\lattice\dual$. The image of this map is contained in the
sublattice $\lattice\dual_{\text{\scriptsize\rm effective}}
\subset\lattice\dual$. But due to condition \Em{Quasi--momenta}~(i)
in Section~\ref{subsection complex Fermi curves of finite genus}
the image might be a proper sublattice of this sublattice.
However, in these cases the transformation of
Lemma~\ref{covariant transformation} transforms the
\Em{complex Fermi curve} into other isomorphic
\Em{complex Fermi curve} without changing the image of this map,
and such that this image coincides with the sublattice
corresponding to the transformed \Em{complex Fermi curves}.
\end{Remark}

Now let us assume that $\fermi\in
\bigcup\limits_{\willmore>0}\Bar{\moduli}_{\lattice,\eta,\willmore}$
is a \Em{complex Fermi curve}, whose effective dual lattice
$\lattice\dual_{\text{\scriptsize\rm effective}}$
is a proper sublattice of $\lattice\dual$.
Obviously, there exists a proper sublattice
$\Breve{\lattice}\dual$ with the following properties:
\begin{description}
\index{lattice!dual sub--$\sim$ $\Breve{\lattice}\dual$}
\item[Dual sublattice (i)] The proper sublattice
  $\Breve{\lattice}\dual\subset\lattice\dual$ contains
  $\lattice\dual_{\text{\scriptsize\rm effective}}$.
\item[Dual sublattice (ii)] The order of
  $\lattice\dual/\Breve{\lattice}\dual$ is finite.
\end{description}
The dual lattice denoted by $\Breve{\lattice}$ contains $\lattice$,
and the order of $\Breve{\lattice}/\lattice$ is the same as
the order of the quotient in \Em{Dual sublattice}~(ii).
For all \Em{complex Fermi curves} $\Breve{\fermi}\in
\bigcup\limits_{g\in\mathbb{N}}
\Bar{\moduli}_{g,\Breve{\lattice},\eta}$,
the finite union of the subvarieties
$\Breve{\fermi}\subset\mathbb{C}^2$
shifted by representatives contained in all elements of
$\lattice\dual/\Breve{\lattice}\dual$
yields a \Em{complex Fermi curve}
$\fermi\in
\bigcup\limits_{g\in\mathbb{N}}
\Bar{\moduli}_{g,\lattice,\eta}$,
whose effective dual lattice is contained in $\Breve{\lattice}\dual$.
Moreover, this relation provides a one--to--one correspondence
between \Em{complex Fermi curves}
with lattice $\Breve{\lattice}$
and \Em{complex Fermi curves} with lattice $\lattice$,
whose effective dual lattice is contained in $\Breve{\lattice}\dual$.
The quotient $\Breve{\fermi}/\Breve{\lattice}\dual$
is a one--sheeted covering over the quotient
$\fermi/\lattice\dual$.
In particular, the normalizations of these quotients are isomorphic.
The \Em{first integral} of the \Em{complex Fermi curve}
$\fermi$ with lattice $\lattice$
is $\left|\Breve{\lattice}/\lattice\right|$ times
the \Em{first integral} of the \Em{complex Fermi curve}
$\Breve{\fermi}$ with lattice $\Breve{\lattice}$.

Furthermore, if the \Em{complex Fermi curve}
with lattice $\Breve{\lattice}$
is a \De{Generalized Weierstra{\ss} curve}~\ref{generalized curves},
then the \Em{complex Fermi curve}
with lattice $\lattice$ is also a
\De{Generalized Weierstra{\ss} curve}~\ref{generalized curves}.
However, there are some situations, where the
\Em{complex Fermi curve} with lattice $\lattice$
is a \De{Generalized Weierstra{\ss} curve}~\ref{generalized curves},
whereas the \Em{complex Fermi curve}
with lattice $\Breve{\lattice}$ is not a
\De{Generalized Weierstra{\ss} curve}~\ref{generalized curves}.
Of special importance is the following situation:
Let $\fermi$ be a \Em{complex Fermi curve}
with lattice $\Breve{\lattice}$,
containing two different elements of
$\Breve{\lattice}\dual/2\Breve{\lattice}\dual$,
whose difference is mapped under the canonical map
$$\Breve{\lattice}\dual/2\Breve{\lattice}\dual\rightarrow
\lattice\dual/2\lattice\dual$$
onto zero. Then the \Em{complex Fermi curve}
with lattice $\lattice$ is a
\De{Generalized Weierstra{\ss} curve}~\ref{generalized curves}
in contrast to the \Em{complex Fermi curve}
with lattice $\Breve{\lattice}$
not necessarily being a
\De{Generalized Weierstra{\ss} curve}~\ref{generalized curves}.
We will see that this is true for the minimizers of the
\De{Generalized Willmore functional}~\ref{generalized functional}.

For any lattice $\lattice\dual$ there exist exactly three
sublattices $\Breve{\lattice}\dual$,
with the property that the quotient lattice
$\lattice\dual/\Breve{\lattice}\dual$
is isomorphic to $\mathbb{Z}_2$.
In fact, all these lattices have to contain the lattice
$2\lattice\dual$, and exactly one non--zero element of
$\lattice\dual/2\lattice\dual$ is contained
in the image of the natural map
$\Breve{\lattice}\dual/2\Breve{\lattice}\dual\rightarrow
\lattice\dual/2\lattice\dual$.
Hence we may parameterize these sublattices
by the non--zero elements of this group.
For all non--zero $[\kappa/2]\in\lattice\dual/2\lattice\dual$
let $\lattice\dual_{[\kappa/2]}$
denote the unique sublattice with the following properties:
\begin{description}
\index{lattice!dual $\sim$ $\lattice\dual_{[\kappa/2]}$}
\item[$\lattice\dual_{[\kappa/2]}$ (i)]
  The group $\lattice\dual/\lattice\dual_{[\kappa/2]}$
  is isomorphic to $\mathbb{Z}_2$.
\item[$\lattice\dual_{[\kappa/2]}$ (ii)]
  The image of the natural map
  $\lattice\dual_{[\kappa/2]}/2\lattice\dual_{[\kappa/2]}\rightarrow
  \lattice\dual/2\lattice\dual$ is equal to $\{0,[\kappa/2]\}$.
\end{description}
The above considerations have the following consequence:
If $\Breve{\fermi}$ is the unique \Em{complex Fermi curve},
with dual lattice $\lattice\dual_{[\kappa/2]}$
with disconnected normalization,
and whose \Em{first integral} is equal to $4\pi$,
then the normalization of the corresponding
\Em{complex Fermi curve} $\fermi$
with dual lattice $\lattice\dual$
has two \Em{Orbits of type~2} over the element $[\kappa/2]$
and one \Em{Orbit of type~2} over the zero element of
$\lattice\dual/2\lattice\dual$.
One can prove that these are all elements of the form
$\lattice\dual/2\lattice\dual$.
Moreover, the \Em{first integral} is equal to $8\pi$.
In particular, we have three
\De{Generalized Weierstra{\ss} curves}~\ref{generalized curves},
whose \Em{first integrals} are equal to $8\pi$.
We will now show that these are all
\De{Generalized Weierstra{\ss} curves}~\ref{generalized curves}
with disconnected normalization,
whose \Em{first integrals} are equal to $8\pi$.

\begin{Lemma} \label{admissible disconnected}
The set of
\De{Generalized Weierstra{\ss} curves}~\ref{generalized curves}
with disconnected normalizations,
whose \Em{first integrals} are equal to $8\pi$,
consists of exactly three elements. More precisely, for each non--zero element
$[\kappa/2]\in\lattice\dual/2\lattice\dual$
there exists exactly one such \Em{complex Fermi curve} fulfilling the
\De{Weak Singularity condition}~\ref{weak singularity condition}
with this element $[\kappa/2]$ and no
such \Em{complex Fermi curve} fulfilling the
\De{Weak Singularity condition}~\ref{weak singularity condition}
with the zero element of
$\lattice\dual/2\lattice\dual$.
The effective dual lattice of these
\De{Generalized Weierstra{\ss} curves}~\ref{generalized curves}
is equal to $\lattice\dual_{[\kappa/2]}$.
\end{Lemma}

\begin{proof}
Given a
\De{Generalized Weierstra{\ss} curve}~\ref{generalized curves}
$\fermi$ with disconnected normalization,
whose \Em{first integral} is equal to $8\pi$,
let $\Bar{\fermi}^{\pm}_{\text{\scriptsize\rm normal}}$
denote the components of the normalization,
which contain $\infty^{\pm}$.
Due to the properties of the restriction of the
multi--valued function $k$ to
$\Bar{\fermi}^+_{\text{\scriptsize\rm normal}}$,
this compact Riemann surface is the spectral curve
of an elliptic solution of the Kadomtsev--Petviashvili Equation.
More precisely, on this compact Riemann surface there exists a
\Em{Baker--Akhiezer} function corresponding to these solutions,
and the corresponding elliptic curve is the unique elliptic curve,
which is conformally equivalent to the flat torus $\torus$.
In \cite[Theorem~4.]{Kr2} these spectral curves are constructed
as $N$--sheeted coverings of the
corresponding elliptic curves.
Due to the general construction of the algebraic-geometric solutions
\cite{Kr1,DKN} the space, on which these solutions are defined,
is mapped linearly into the Jacobian variety of the
corresponding spectral curve. Hence for elliptic solutions
the elliptic curve is mapped linearly into
the corresponding Jacobian variety. The dual map
\cite[Chapter~2 \S4]{LB} is a map
from the Jacobian variety onto the elliptic curve.
The composition of the canonical map from the spectral curve into
the Jacobian variety with this map makes the spectral curve of
an elliptic solutions in a canonical way
to a finite--sheeted covering of the corresponding elliptic curve.
In the present situation the function $k_1+\sqrt{-1}k_2$
yields this canonical covering
(compare with Lemma~\ref{disconnected}).
Due to Lemma~\ref{disconnected} the
\Em{first integrals} of the corresponding
\Em{complex Fermi curves} are equal to
$4N\pi$. Hence in the present situation $N$ is equal to $2$.
Due to \cite[Lemma~2.]{Kr2} the genus of
$\Bar{\fermi}^+_{\text{\scriptsize\rm normal}}$
is therefore not larger than $2$
(in \cite[Lemma~2.]{Kr2} the \Em{arithmetic genus}
of the spectral curve is shown to be equal to $N$).

Since the \Em{complex Fermi curve} is assumed to obey the
\De{Weak Singularity condition}~\ref{weak singularity condition},
the connected component
$\Bar{\fermi}^+_{\text{\scriptsize\rm normal}}$
has a singularity.
Due to Proposition~\ref{smooth moduli} there exist families of
deformations of these \Em{complex Fermi curves} within the set of
\Em{complex Fermi curves} with disconnected normalization,
whose \Em{first integral} is equal to $8\pi$,
such that this singularity is removed.
The components of the normalizations of the these families,
which contain $\infty^+$, correspond to these
elliptic solutions of the Kadomtsev--Petviashvili Equation.
We conclude that the genus of
$\Bar{\fermi}^+_{\text{\scriptsize\rm normal}}$
is at most equal to $1$. On the other hand, the \Em{first integral}
of the unique \Em{complex Fermi curve},
whose normalization has two connected components of genus zero,
is zero.
Therefore, $\Bar{\fermi}^+_{\text{\scriptsize\rm normal}}$
has genus one and
due to the Riemann--Hurwitz formula \cite[17.14]{Fo}
it is an unbranched two--sheeted covering of the
unique elliptic curve,
which is conformally equivalent to $\mathbb{R}/\lattice$.
This implies that the effective dual lattice has to be a
proper sublattice and that
$\lattice\dual/\lattice\dual_{\text{\scriptsize\rm effective}}$
is isomorphic to $\mathbb{Z}_2$.
We conclude that the \Em{complex Fermi curve} is one of the three
\De{Generalized Weierstra{\ss} curves}~\ref{generalized curves}
constructed above. 
\end{proof}

With the help of Lemma~\ref{unique real point minimizers}
we may construct a whole family of
\De{Generalized Weierstra{\ss} curves}~\ref{generalized curves},
whose \Em{first integrals} are not larger than $8\pi$.
In this lemma we constructed a family of \Em{complex Fermi curves}
parameterized by the corresponding values of the \Em{first integral}
$\willmore\in[0,4\pi]$.
For all non--zero elements 
$[\kappa/2]\in\lattice\dual/2\lattice\dual$ let
$$[0,8\pi]\rightarrow \Bar{\moduli}_{3,\lattice,\eta,\sigma,8\pi},
\willmore\mapsto \fermi_{\min,[\kappa/2]}(\willmore/2)$$
\index{family!$\fermi_{\min,[\kappa/2]}(\cdot)$}
denote this family of \Em{complex Fermi curves}
with lattice $\lattice_{[\kappa/2]}\subset\lattice$
considered as a family of \Em{complex Fermi curves}
with lattice $\lattice$.
Since the \Em{first integrals} of these \Em{complex Fermi curves},
considered as a \Em{complex Fermi curve}
with lattice $\lattice$,
are twice the \Em{first integrals} of these
\Em{complex Fermi curves}, considered as
\Em{complex Fermi curves}
with lattice $\lattice_{[\kappa/2]}$,
the parameter $\willmore$ is equal to the \Em{first integral}
of the corresponding \Em{complex Fermi curves}.
From the proof of Lemma~\ref{unique real point minimizers}
we conclude that there exists an unique $\willmore_{[\kappa/2]}$,
such that the \Em{complex Fermi curve}
$\fermi_{\min,[\kappa/2]}(\willmore/2)$ considered as a
\Em{complex Fermi curve}
with lattice $\lattice_{[\kappa/2]}$
contains both elements in the preimage of $[\kappa/2]$
under the natural map
$\lattice\dual_{[\kappa/2]}/2\lattice\dual_{[\kappa/2]}\rightarrow
\lattice\dual/2\lattice\dual$, if and only if $\willmore\geq
\willmore_{[\kappa/2]}$. For all non--zero elements
$[\kappa/2]\in\lattice\dual/2\lattice\dual$ this yields a family of
\Em{complex Fermi curves} fulfilling the
\De{Weak Singularity condition}~\ref{weak singularity condition}
with the element $[\kappa/2]$ parameterized by the corresponding
\Em{first integrals} $\willmore\in[\willmore_{[\kappa/2]},8\pi]$.
In the proof of Lemma~\ref{unique real point minimizers} we described the
singularities of the family
$\fermi_{\min,[\kappa/2]}(\willmore)$
with $\willmore\in[0,4\pi]$,
considered as a family of \Em{complex Fermi curves}
with lattice $\lattice_{[\kappa/2]}$.
But besides these singularities this family
also contains singularities,
which are not singularities of the same family
considered as a family of \Em{complex Fermi curves}
with lattice $\lattice_{[\kappa/2]}$.
More precisely, we have seen above that the subvarieties in
$\mathbb{C}^2/\lattice\dual$ of this
family considered as a family of \Em{complex Fermi curves}
with lattice $\lattice$
are one--sheeted coverings of the subvarieties in
$\mathbb{C}^2/\lattice\dual_{[\kappa/2]}$ of this
family considered as a family of \Em{complex Fermi curves}
with lattice $\lattice_{[\kappa/2]}$.
Those elements of the corresponding \Em{complex Fermi curves}
with lattice $\lattice$,
where this map is not biholomorphic, are exactly the
singularities of these \Em{complex Fermi curves},
which are not singularities of the corresponding
\Em{complex Fermi curves}
with lattice $\lattice_{[\kappa/2]}$.

\begin{Lemma} \label{admissible family}
For all $\willmore\in[\willmore_{[\kappa/2]},8\pi]$ the
\Em{complex Fermi curve} $\fermi_{\min,[\kappa/2]}(\willmore/2)$
has besides the singularity of the
\De{Weak Singularity condition}~\ref{weak singularity condition}
and the singularities of this \Em{complex Fermi curve}
considered as a \Em{complex Fermi curve}
with lattice $\lattice_{ [\kappa/2]}$
only ordinary double points of the form $(y,\eta(y))$,
where the function $dk_2/dk_1$ has
positive imaginary part at those points of these double points,
which belong to
$\Bar{\fermi}^+_{\text{\scriptsize\rm normal}}$,
and no other singularities.
\end{Lemma}

\begin{proof}
Let $\fermi(t)$ be a continuous family of \Em{complex Fermi curves}
parameterized by $\parameter{t}\in(-\varepsilon,\varepsilon)$,
which for $\parameter{t}<0$ has a continuous family of
double points of the form $(y(\parameter{t}),\eta(y(\parameter{t}))$
with the property that the imaginary part of
$dk_2/dk_1$ at $y(\parameter{t})$ is positive.
Moreover, let us assume in addition
that for $\parameter{t}=0$ the imaginary part of this function at
$y(0)$ is zero,
and that in some neighbourhood of $(y(0),\eta(y(0)))$ the
\De{Local contributions to the arithmetic genus}~\ref{local contribution}
does not depend on $\parameter{t}$.
Then at least for small $t<0$ there exists another continuous
family of double points $(y'(\parameter{t}),\eta(y'(\parameter{t}))$,
with the property that the imaginary part of $dk_2/dk_1$
at $y'(\parameter{t})$ is negative,
and the limit $\parameter{t}\rightarrow 0$ of $y(\parameter{t})$ and
$y'(\parameter{t})$ coincide. In fact,
all elements $y$, at which the real part of the function $k$ takes
some value in $\lattice\dual/2\lattice\dual$, belong to
double points of the form $(y,\eta(y))$. The two form
$d\Re(k_1)\wedge d\Re(k_2)$ vanishes on the normalization,
if the function $dk_2/dk_1$ is real. Besides these zeroes all
other zeroes of $d\Re(k_1)\wedge d\Re(k_2)$
on the normalization are \Em{multiple points}. The former have
real co--dimension one and the latter real co--dimension two.
Hence the image under the map $\Re(k)$ in
$\mathbb{R}^2$ is locally a finite--sheeted covering. Furthermore,
the lines, where the function $dk_2/dk_1$ takes real values,
are boundaries of these coverings,
where two different sheets with opposite orientation meet each other.

Now let us consider some continuous family
$\fermi(\parameter{t})$ of \Em{complex Fermi curves}, whose real
part divides the normalization into two connected components.
Moreover, the \Em{complex Fermi curve} $\fermi(0)$
should have the property, that all
singularities are double
points of the form $(y,\eta(y))$, where $dk_2/dk_1$ has positive
imaginary part at the point of $(y,\eta(y))$ contained in the
component of $\infty^+$. Then this property is conserved for all
$0\leq \parameter{t}$ as long as the real part of
$\fermi(\parameter{t})$ does not contain an
element of $\lattice\dual/2\lattice\dual$. In fact, due to
the assumption these double points can loose this property only if
they meet the boundary between
$\Bar{\fermi}^+_{\text{\scriptsize\rm normal}}$ and
$\Bar{\fermi}^-_{\text{\scriptsize\rm normal}}$,
which is the real part.

We distinguished between singularities of these \Em{complex Fermi curves}
considered as \Em{complex Fermi curves}
with lattice $\lattice_{[\kappa/2]}$
and singularities,
which are not singularities of these \Em{complex Fermi curves}
considered as \Em{complex Fermi curves}
with lattice $\lattice_{[\kappa/2]}$.
Ordinary double points of the form $(y,\eta(y))$
belong to the first class,
if at these points the real part of the function $k$
takes values in the image of the canonical map
$\lattice\dual_{[\kappa/2]}/2\lattice\dual_{[\kappa/2]}
\rightarrow\lattice\dual/2\lattice\dual$, and otherwise to the
second class. From the above considerations it
follows, that it suffices to show that for one \Em{complex Fermi curve} of
this family the statement of the lemma is true and that the real part
of this family does not contain an element in
$\lattice\dual/2\lattice\dual$, which does not belong to the image
of the canonical map
$\lattice\dual_{[\kappa/2]}/2\lattice\dual_{[\kappa/2]}
\rightarrow\lattice\dual/2\lattice\dual$.

We first consider the rectangular lattice
$\lattice=\mathbb{Z}^2\subset\mathbb{R}^2$, and the family
corresponding to the non--zero element
$[(1/2,1/2)]\in\lattice\dual/2\lattice\dual$.
The corresponding lattice $\lattice_{[(1/2,1/2)]}$
is generated by the elements $(1/2,1/2)$ and
$(1/2,-1/2)$, and the corresponding dual lattice
$\lattice\dual_{[(1/2,1/2)]}$ by the elements $(1,1)$ and
$(1,-1)$. It is quite easy to see (and we shall do this at the end of
this section), that in this case $\willmore_{[(1/2,1/2)]}$ is equal
to $2\pi^2$, and that the corresponding \Em{complex Fermi curve} may be
described by the equation $g(k,k)=1/2$. An easy calculation shows that
besides the singularity of the
\De{Weak Singularity condition}~\ref{weak singularity condition} all
singularities of this \Em{complex Fermi curve} are
ordinary double points of the form $(y,\eta(y))$,
where the imaginary part of the function
$dk_2/dk_1$ is positive at those elements of these double points,
which belong to
$\Bar{\fermi}^+_{\text{\scriptsize\rm normal}}$.
The real part of this
\Em{complex Fermi curve} is a circle of radius $1/\sqrt{2}$. From
Lemma~\ref{unique real point minimizers} we conclude that all real parts
of the family $\fermi_{\min,[(1/2,1/2)]}(\willmore/2)$ with
$\willmore\in(2\pi^2,8\pi]$ cover the complement in the
$\torus\dual_{[(1/2,1/2)]}$ of the closed disc
around the origin with radius $1/\sqrt{2}$. Obviously, this complement
contains only one element of
$\lattice\dual/2{\lattice}\dual_{[(1/2,1/2)]}$, which belongs to
the kernel of the canonical map
$\lattice\dual/2{\lattice}\dual_{[(1/2,1/2)]}\rightarrow
\lattice\dual/2{\lattice}\dual$.
This shows that for the rectangular class and the family corresponding
to the element $[(1/2,1/2)]\in\lattice\dual/2\lattice\dual$ the
statement of the lemma is true.

Now we consider all these families as one family parameterized by
the \Em{first integral} and the conformal classes.
It is quite easy to see,
that the groups $\lattice\dual/2\lattice\dual$ build a
non--trivial bundle over the moduli space of conformal
classes (compare with $\moduli_1$ at the end of this section).
More precisely, the non--zero elements of this group
fit to a connected three--sheeted covering over this
moduli space, and the zero element to a one--sheeted
covering. This follows from the explicit calculation of the
conformal equivalence classes of the lattices $\lattice_{[\kappa/2]}$
in dependence of the conformal equivalence classes of $\lattice$
given at the end of this section. We conclude
that this family is connected. With the above considerations
it suffices to show hat the real part of all elements of this
family does not contain any element of
$\lattice\dual/2\lattice\dual$, which is not contained in the
image of the canonical map
$\lattice\dual/2{\lattice}\dual_{[(1/2,1/2)]}\rightarrow
\lattice\dual/2{\lattice}\dual$. In order to prove this we
restrict this family again to one conformal class and some
non--zero element $[\kappa/2]\in\lattice\dual/2\lattice\dual$.
Due to the proof of Lemma~\ref{unique real point minimizers}
the real parts of all these families
$\fermi_{\min,[\kappa/2]}(\willmore)$, with $\willmore\in[0,4\pi]$
meet every element of $\torus\dual$ exactly
twice. In particular, each of the two elements of
$\lattice\dual/2\lattice\dual$, which are not contained in the
image of the canonical map
$\lattice\dual/2{\lattice}\dual_{[\kappa/2]}\rightarrow
\lattice\dual/2{\lattice}\dual$, is met twice.
However, the singularity of the
\De{Weak Singularity condition}~\ref{weak singularity condition}
always includes four double points,
which are not double points of any
\De{Generalized Weierstra{\ss} curve}~\ref{generalized curves}
of this family considered as a \Em{complex Fermi curve}
with lattice $\lattice_{[\kappa/2]}$.
We conclude that the real parts of the family
$\fermi_{\min,[\kappa/2]}(\willmore/2)$ meet all those elements
of $\lattice\dual/2\lattice\dual$, which are not contained
in the image of the canonical map
$\lattice\dual/2{\lattice}\dual_{[\kappa/2]}\rightarrow
\lattice\dual/2{\lattice}\dual$ for values of $\willmore$
in $[0,\willmore_{[\kappa/2]})$, and not for values in
$[\willmore_{[\kappa/2]},8\pi]$. This completes the proof.
\end{proof}

Now we can classify all relative minimizers of the
\De{Generalized Willmore functional}~\ref{generalized functional},
whose \Em{first integrals} are not larger than $8\pi$.

\begin{Proposition} \label{small relative minimizers}
\index{minimizer!of the generalized Willmore functional}
The restriction of the \Em{first integral} to the space of all
\Em{complex Fermi curves}, which obey the
\De{Weak Singularity condition}~\ref{weak singularity condition} for some
element $[\kappa/2]\in\lattice\dual/2\lattice\dual$ has exactly
one relative minimizer, whose \Em{first integral} is not larger than $8\pi$,
if $[\kappa/2]$ is not zero, and otherwise no such relative minimizer.
The \Em{complex Fermi curves} of these relative minimizers are
$\fermi_{\min,[\kappa/2]}(\willmore_{[\kappa/2]}/2)$.
In particular, the corresponding values of the
\Em{first integral} are equal to $\willmore_{[\kappa/2]}$.
\end{Proposition}

As a preparation for the proof of this Proposition
we first give the following 

\begin{Lemma} \label{positive form}
Let $\Spa{Y}$ be a compact  connected Riemann surface with some
anti--holomorphic involution $\rho$, such that the real part
$\Spa{Y}_\mathbb{R}$ (= the set of fixed points of $\rho$)
divides $\Spa{Y}$ into
two connected components $\Spa{Y}^+$ and $\Spa{Y}^-$ with
$\partial \Spa{Y}^+=\partial \Spa{Y}^-=\Spa{Y}_\mathbb{R}$.
Let $y^+\in \Spa{Y}^+$ and $y^-\in \Spa{Y}^-$ be two different points,
which are interchanged by the involution $\rho$. Then there exists a
meromorphic real (i.\ e.\ $\rho^{\ast}\omega=\Bar{\omega}$) differential
form $\omega$, which has poles of first--order at $y^+$ and $y^-$ and
no other poles, and is positive on the real part, with respect to
the orientation induced from the orientation of $\Spa{Y}^+$ on the boundary
$\Spa{Y}_\mathbb{R}=\partial \Spa{Y}^+$.
\end{Lemma}

\begin{proof} Let $\Omega_{y^++y^-}$ be the space of
  real meromorphic forms with a pole of at most first--order at $y^+$
  and $y^-$ and with no other poles. We shall first prove that for
  any choice of finitely many points $y_1,\ldots,y_n$ of
  $\Spa{Y}_\mathbb{R}$ there exists some
  $\omega\in\Omega_{y^++y^-}$, which
  is positive at these points with respect to the given orientation.
  The values of all elements of $\Omega_{y^++y^-}$ at these points in
  terms of local real positive holomorphic
  forms $dz_1,\ldots,dz_n$ form a linear subsapce $V$ of $\mathbb{R}^n$.
  If this subspace contains no element $v$,
  all whose components are positive,
  then the orthogonal complement
  (with respect to the Euclidean scalar product)
  of this subspace contains some non zero element
  $c\in\mathbb{R}^n$, all whose components are non--negative.
  In fact, due to the
  Separating hyperplane theorem \cite[Theorem~V.4]{RS1}
  there exists a linear functional on $\mathbb{R}^n$,
  whose values on the open convex cone of all elements,
  all whose components are positive,
  are larger or equal to the values on $V$.
  This functional may be written as the Euclidean scalar product
  with some element $c\in\mathbb{R}^n$,
  which has the desired properties.
  This element $c$ in the orthogonal complement of $V$
  defines some Mittag-Leffler distribution,
  and due to \cite[Theorem~18.2]{Fo}
  this Mittag-Leffler distribution has a global solution $f$.
  Obviously $\rho^{\ast}\Bar{f}$ solves
  the same Mittag-Leffler distribution.
  So let us assume that $f$ is real (i.\ e.\ $f=\rho^{\ast}\Bar{f}$).
  Since all entries of $c$ are non--negative
  the restriction of the function $f$
  to any connected component of $\Spa{Y}_\mathbb{R}$
  takes all real values with the same multiplicity.
  But $\Omega_{y^++y^-}$ contains besides the real holomorphic forms
  also meromorphic forms with poles of at most first--order
  at $y^+$ and $y^-$.
  The same arguments as in the proof of
  \cite[Theorem~18.1 and 18.2]{Fo} show,
  that in this case $f$ has a zero at $y^+$ and $y^-$,
  which is impossible.
  Hence there exists some element $\omega\in\Omega_{y^++y^-}$,
  which is positive at $y_1,\ldots,y_n$.

  Now we choose any norm on $\Omega_{y^++y^-}$ (they are all
  equivalent, since $\Omega_{y^++y^-}$ is a finite--dimensional
  vector space).
  The subspace of $\Omega_{y^++y^-}$ consisting of forms,
  whose norm is equal to $1$ is compact.
  All subsets of this compact space of forms,
  which are non--negative at finitely many points of $\Spa{Y}_\mathbb{R}$,
  are closed. Moreover, due to our first claim
  the intersection of finitely many of such subsets is non--empty.
  Hence the intersection of all these subsets is also non--empty.
  Hence there exists some form
  $\omega\in\Omega_{y^++y^-}\setminus\{0\}$,
  which is non--negative on $\Spa{Y}_\mathbb{R}$.
  This form can have only finitely many zeroes.
  If we add some small multiple of an element of $\Omega_{y^++y^-}$,
  which is positive at all the zeroes of $\omega$
  contained in $\Spa{Y}_\mathbb{R}$,
  we obtain an element of $\Omega_{y^++y^-}$,
  which is positive on $\Spa{Y}_\mathbb{R}$.
\end{proof}

\noindent
{\it Proof of Proposition~\ref{small relative minimizers}.}
Let us first prove that we may deform all relative minimizers of the
restriction of the \Em{first integral} to the space of
\De{Generalized Weierstra{\ss} curves}~\ref{generalized curves} into a
\De{Generalized Weierstra{\ss} curve}~\ref{generalized curves}, whose
normalization has two connected components. Let us consider the class
of \Em{complex Fermi curves}, which obey
conditions~(i)--(iv) of Corollary~\ref{deformed minima}.
Since the complement of the image of the map described in (ii) are
disjoint subsets, the total area of this complement is smaller than
the volume of $\torus\dual$. Therefore,
due to Lemma~\ref{gluing rule}, on these \Em{complex Fermi curves} the
\Em{first integral} takes values in the interval $(4(d-1)\pi,4d\pi]$,
where $d$ denotes the number of sheets in (ii).
Since we are interested in \Em{complex Fermi curves},
whose \Em{first integral} belongs to $(4\pi,8\pi]$,
in our case this number of sheets is equal to two.
Due to Corollary~\ref{deformed minima} we may deform all
relative minimizers of the restriction of the
\Em{first integral} to the space of
\De{Generalized Weierstra{\ss} curves}~\ref{generalized curves}
by an arbitrary small deformation into a
\Em{complex Fermi curve} within this class.

Now we claim that we may deform any element of this class
within this class into a \Em{complex Fermi curve},
with disconnected normalization.
Due to Proposition~\ref{smooth moduli} and
Lemma~\ref{positive form} we may deform the \Em{complex Fermi curves}
within this class in such a way that the disjoint closed subsets
mentioned in (ii) become smaller and smaller,
if the dividing real parts of the \Em{complex Fermi curves}
do not have \Em{cusps}. Due to Proposition~\ref{modified smooth moduli}
and the remarks following this proposition this condition
is necessary in order to preserve condition~(i).
Hence we may shrink these pieces mentioned in (ii)
until the boundary of these pieces intersect each other.
In this case the real part of the deformed \Em{complex Fermi curve}
has some multiple points. Either they are
double points of the form $(y,\sigma(y))$,
or they are pairs of double points
of the form $(y,y')$ and $\sigma(y),\sigma(y'))$.
In the second case the deformations of
Proposition~\ref{smooth moduli} allow
to split these pieces into several pieces and shrink them further.
In the first case the deformations of \Em{Orbits of type~1}
into \Em{Orbits of type~2} constructed in the proof of
Proposition~\ref{smooth moduli} concerning \Em{Multiple points O1/O2}
and \Em{Multiple points O1/O1} have the same effect.
Due to Corollary~\ref{dividing genus bound} the space of all
\Em{complex Fermi curves} with dividing real part, whose
\Em{first integral} is not larger than $8\pi$ is compact
and contained in $\Bar{\moduli}_{g_{\max},\lattice,\eta,\sigma,8\pi}$.
Furthermore, the proof of this lemma actually shows
that at all zeroes of $d\xx{p}$ and $d\yy{p}$
the imaginary part of $k$ is bounded.
Hence the formula for the \Em{first integral}
in Lemma~\ref{gluing rule} shows
that the \Em{first integral} is continuous on this space,
in contrast to the lower-continuity
on $\Bar{\moduli}_{g,\lattice,\eta}$.
Hence we may shrink the complement of the image of the map
described in (ii) until it is empty.
The normalization of the corresponding
\Em{complex Fermi curve} has two connected components.
Since we have started with a \Em{complex Fermi curve},
whose \Em{first integral} is smaller than $8\pi$ the final
\Em{complex Fermi curve} of this deformation is one of the three
\De{Generalized Weierstra{\ss} curves}~\ref{generalized curves}
described in Lemma~\ref{admissible disconnected}.
Obviously properties (i)--(ii) are preserved under these deformations.

It remains to show that also conditions~(iii)--(iv) are preserved,
and that the \Em{complex Fermi curves} within this class
do not have real cusps. The considerations in the proof of
Lemma~\ref{admissible family} concerning ordinary double points
of the form $(y,\eta(y))$ imply that condition~(iii)
is preserved under these deformations.
Due to deformations of \Em{Multiple points O2/O2} and
\Em{Multiple points O3} constructed in the proof of
Theorem~\ref{relative minimizers} we may always deform
by some arbitrary small deformation an
\Em{Orbit of type~3} into two \Em{Orbits of type~2},
and this deformation may be chosen within this class.
Due to Lemma~\ref{positive form} this deformation may be chosen
in such a way, that these pieces never contain the element
corresponding to $[\kappa/2]$.
Hence the two \Em{Orbits of type~2}
in the normalizations of this family
can only be deformed into an \Em{Orbit of type~3} over $[\kappa/2]$.
In this case again the deformations of \Em{Multiple points O2/O2}
and \Em{Multiple points O3} constructed in the proof of
Theorem~\ref{relative minimizers} show that the deformation can be
continued by deforming this \Em{Orbit of type~3}
back into two \Em{Orbits of type~2}.
Therefore, condition~(iv) is preserved.

Finally, a real \Em{cusp} contradicts conditions~(i)--(iii) with
$\willmore\in(4\pi,8\pi]$.
In fact, if both differentials $d\xx{p}$ and $d\yy{p}$
have zeroes of at least second order,
then at least three sheets of the function $k_1+\sqrt{-1}k_2$
have a branch point at the real \Em{cusp}.
Obviously this contradicts condition~(ii) for two--sheeted coverings
($\iff\willmore\in(4\pi,8\pi]$).
Furthermore, all other real \Em{cusps} are locally of the form
$\yy{p}^2=\xx{p}^{2m+1}$.
All deformations of such \Em{cusps} contain either
no singularities nearby the real \Em{cusp},
or double points of the form $(y,\rho(y))$,
or ordinary real double points,
or real \Em{cusps} of the same form of less order
(corresponding to smaller $m$).
In the first case due to condition~(i) all zeroes of $d\xx{p}$
nearby the deformed real \Em{cusp} have to be real,
and condition~(ii) is violated.
The second case contradicts condition~(iii),
and the third case condition~(ii).
Finally, the fourth case is excluded by induction.
This shows that we may deform all relative minima within the class of
\Em{complex Fermi curves} obeying conditions~(i)--(iv) into a
\De{Generalized Weierstra{\ss} curve}~\ref{generalized curves}
with disconnected normalization.

In a second step we claim that, if the number of sheets in
condition (ii) is chosen to be equal to two, (in this case the
\Em{first integral} takes values in the interval $(4\pi,8\pi]$),
these classes are equal to the families
$\fermi_{\min,[\kappa/2]}(\willmore/2)$,
with $\willmore\in(W_{[\kappa/2]},8\pi]$. Due to the first claim it
suffices to show that this family is the only possible deformation of
the
\De{Generalized Weierstra{\ss} curves}~\ref{generalized curves}
described in Lemma~\ref{admissible disconnected} within this class.
Due to Lemma~\ref{admissible family} all elements of these family
may not be deformed by a small deformation within the class of
\Em{complex Fermi curves} with dividing real parts
into a \Em{complex Fermi curve},
whose effective dual lattice is not contained in
$\lattice\dual_{[\kappa/2]}$. Lemma~\ref{unique real point minimizers}
implies that all double points of these \Em{complex Fermi curves}
considered as \Em{complex Fermi curves}
with lattice $\lattice_{[\kappa/2]}$
have only two types of singularities:
\begin{description}
\item[Ordinary double points.] Ordinary double points
  of the form $(y,\eta(y))$,
  where the imaginary part of the function $dk_2/dk_1$
  takes positive values at those points of these double points,
  which belong to $\Bar{\fermi}^+_{\text{\scriptsize\rm normal}}$.
\item[Real double points.] Double points of the form $(y,\eta(y))$,
  which belong to the real part.
  In this case the function $dk_2/dk_1$
  takes the same real value at both elements of these double points.
\end{description}
More precisely, the \Em{real double points} occur only for those
\Em{complex Fermi curves} of this family, which are described in
Corollary~\ref{simplified minimizers}. Moreover, due to the proof of
Lemma~\ref{unique real point minimizers} the only possible deformations
of these \Em{complex Fermi curves} within the class of
\Em{complex Fermi curves}
obeying conditions~(i)--(ii) is the family described in
Lemma~\ref{unique real point minimizers}. This proves the second
claim.

We conclude that the \Em{complex Fermi curves}
$\fermi_{\min,[\kappa/2]}(\willmore_{[\kappa/2]}/2)$
are the unique relative minimizers
of the restriction of the \Em{first integral} to the subset of all
\Em{complex Fermi curves} fulfilling the
\De{Weak Singularity condition}~\ref{weak singularity condition}
with the element $[\kappa/2]\in\lattice\dual/2\lattice\dual$.
\qed

\begin{Remark}\label{examples of real cusps}
For special conformal classes
there exist \Em{real--$\sigma$--hyperelliptic} 
\De{Generalized Weierstra{\ss} curves}~\ref{generalized curves}
with \Em{dividing real parts} and \Em{real cups},
whose \Em{first integral} belong to
$(8\pi,8\pi+\varepsilon)$. In fact, those conformal classes,
for which the two unique real elements of $\fermi_{\min}(4\pi)$
belong to a small neighbourhood of an element in
$\lattice\dual/2\lattice\dual$, have a \Em{complex Fermi curve}
with \Em{dividing real part}, two \Em{Orbits of type~2}
over two different non--trivial elements of
$\lattice\dual/2\lattice\dual$ and a \Em{real cusp}.
The corresponding \Em{first integral}
belongs to $(4\pi,4\pi+\varepsilon)$.
This can be seen by continuation of the family
$\fermi_{\min}(\willmore)$ constructed in
Lemma~\ref{unique real point minimizers} along another path.

\noindent
\begin{minipage}[t]{15cm}
At the unique element of the family, where the real part has
one \Em{Orbit of type~1} and two \Em{Orbits of type~2} over all
non--trivial elements of $\lattice\dual/2\lattice\dual$
we increase the \Em{first integral} along the unique path of
\Em{real--$\sigma$--hyperelliptic} \Em{complex Fermi curves},
which have real \Em{multiple points}
(the deformed \Em{Orbits of type~1}). Since the real part of these
\Em{complex Fermi curves} is a self--intersecting circle in
$\mathbb{R}^2$, this family has to end at a
\Em{real--$\sigma$--hyperelliptic} \Em{complex Fermi curve} with two
\Em{real cusps}. This is indicated in the neighbouring figure.
\end{minipage}\hfill
\setlength{\unitlength}{.5cm}
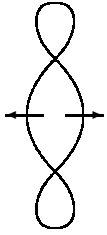
\begin{picture}(3,-6)
\qbezier(1.5,-1.5)(1,-1)(1,-.5)
\qbezier(1,-.5)(1,0)(1.5,0)
\qbezier(1.5,0)(2,0)(2,-.5)
\qbezier(2,-.5)(2,-1)(1.5,-1.5)
\qbezier(1.5,-4.5)(1,-5)(1,-5.5)
\qbezier(1,-5.5)(1,-6)(1.5,-6)
\qbezier(1.5,-6)(2,-6)(2,-5.5)
\qbezier(2,-5.5)(2,-5)(1.5,-4.5)
\qbezier(1.5,-1.5)(0,-3)(1.5,-4.5)
\qbezier(1.5,-1.5)(3,-3)(1.5,-4.5)
\put(1.2,-3){\vector(-1,0){1}}
\put(1.8,-3){\vector(1,0){1}}
\end{picture}
\end{Remark}

\begin{Remark}\label{trivial spin structure}
It is not difficult to see, that there exists a
\De{Generalized Weierstra{\ss} curve}~\ref{generalized curves}
with disconnected normalization, which corresponds to the trivial
spin structure, and whose \Em{first integral} is equal to $12\pi$.
Therefore, it seems to be possible to determine with similar methods
the minimum of the \Em{first integral} on the space of all
\De{Generalized Weierstra{\ss} curves}~\ref{generalized curves}
corresponding to the trivial spin structure.
\end{Remark}

Due to this proposition the absolute minimum of the
\De{Generalized Willmore functional}~\ref{generalized functional}
is equal to the minimum of $\willmore_{[\kappa/2]}$,
where $[\kappa/2]$ runs over the non--zero element of
$\lattice\dual/2\lattice\dual$.
In order to determine these minima
in dependence of the conformal classes,
we have to combine the mappings
$\lattice\mapsto\lattice_{[\kappa/2]}$ and the
\Em{complex Fermi curves} $\fermi_{\min}(\kappa/2)$
described in Corollary~\ref{simplified minimizers}.
More precisely, for any \Em{complex Fermi curve}
$\Breve{\fermi}_{\min}(\Breve{\kappa}/2)$
described in Corollary~\ref{simplified minimizers}
with lattice $\Breve{\lattice}$,
and which contains another non--zero
$[\Breve{\kappa}'/2]\in\Breve{\lattice}\dual/
2\Breve{\lattice}\dual$,
there exists a unique sublattice $\lattice\subset\Breve{\lattice}$
and some unique non--zero
$[\kappa/2]\in\lattice\dual/2\lattice\dual$,
which contains $[\Breve{\kappa}/2]$ and $[\Breve{\kappa}'/2]$
in the preimage under the canonical map
$\Breve{\lattice}\dual/2\Breve{\lattice}\dual
\rightarrow\lattice\dual/2\lattice\dual$.
Despite some degenerate situations, this other element of
$\lattice\dual/2\lattice\dual$ belongs to an \Em{Orbit of type~2}.
Consequently, for genus one this other element is unique,
and for genus two there are two choices.
The lattice $\Breve{\lattice}$ is equal to $\lattice_{[\kappa/2]}$.
Moreover, the \Em{complex Fermi curve}
$\Breve{\fermi}_{\min}(\Breve{\kappa}/2)$
with lattice $\lattice_{[\kappa/2]}$,
considered as a \Em{complex Fermi curve} with lattice $\lattice$,
minimizes the restriction of the \Em{first integral}
to the space of \Em{complex Fermi curves} fulfilling the
\De{Weak Singularity condition}~\ref{weak singularity condition}
with the element $[\kappa/2]$.
Furthermore, since the order of
$\Breve{\lattice}/\lattice$ is two,
the \Em{first integral} of this minimizer
is twice the \Em{first integral} of
$\Breve{\fermi}_{\min}([\kappa/2)$.
Due to Corollary~\ref{simplified minimizers} all
$[\Breve{\kappa}/2]\in
\Breve{\lattice}\dual/2\Breve{\lattice}\dual$,
whose \Em{complex Fermi curve}
$\Breve{\fermi}_{\min}(\Breve{\kappa}/2)$ have genus larger than zero,
fulfill this condition.
Moreover, if the genus is zero, but $\fermi_{\min}(\Breve{\kappa}/2)$
has two \Em{Orbits of type~1} over two different elements of
$\Breve{\lattice}\dual/2\Breve{\lattice}\dual$,
this \Em{complex Fermi curve} also fulfills this condition.
But from Corollary~\ref{simplified minimizers} we conclude
that if the genus of $\fermi_{\min}(\Breve{\kappa}/2)$
is larger than one, there exists another
$[\Breve{\kappa}'/2]\in
\Breve{\lattice}\dual/2\Breve{\lattice}\dual$,
which fulfills this condition,
and whose \Em{first integral} is smaller.
In order to determine the absolute minimizer of the
\De{Generalized Willmore functional}~\ref{generalized functional}
for all lattices,
it suffices therefore to calculate all \Em{complex Fermi curves}
$\fermi_{\min}(\Breve{\kappa}/2)$ of genus at most one.

In order to do this we parameterize the space of all lattices
$\lattice\subset \mathbb{R}^2$. Any linear conformal transformation
$$\Mat{A}:\mathbb{R}^2\rightarrow\mathbb{R}^2,x\mapsto \Mat{A}x,$$
where $\Mat{A}$ is a real $2\times 2$ matrix with
$\Mat{A}^{t}\Mat{A}=\det^2(\Mat{A})\unity$,
induces a transformation of the lattice
$\lattice$ into $\Mat{A}\lattice=\left\{\Mat{A}\gamma\mid
\gamma\in\lattice\right\}$, the dual lattice
$\lattice\dual$ into $\det^{-2}(\Mat{A})\Mat{A}^{t}\lattice\dual=
\left\{\det^{-2}(\Mat{A})\Mat{A}^{t}\kappa
\mid\kappa\in\lattice\dual\right\}$
and the \Em{complex Fermi curve} $\fermi$ into
$\det^{_2}(\Mat{A})\Mat{A}^{t}\fermi=
\{[\det^{-2}(\Mat{A})\Mat{A}^{t}k]\in\mathbb{C}^2/
\det^{-2}(\Mat{A})\Mat{A}^{t}\lattice\dual\mid
[k]\in \fermi\}$. Since the
\De{Weak Singularity condition}~\ref{weak singularity condition}
and the \Em{first integral} is invariant under these transformations,
we should parameterize the conformal equivalence classes
\index{conformal!class!of the lattice $\lattice$}
\index{lattice!conformal class of the $\sim$ $\lattice$}
of these lattices.
It is convenient to use generators of the lattice $\lattice$
in order to parameterize this space.
If $\xx{\gamma}$ and $\yy{\gamma}$ are some generators
of the lattice $\lattice$ with positive orientation,
there exists a unique conformal transformation,
which transforms $\xx{\gamma}$ into $(1,0)$.
The other generator $\yy{\gamma}$ is parameterized by
$\tau=\yy{\gamma}_1+\sqrt{-1}\yy{\gamma}_2$.
Since they have positive orientation,
we conclude that $\Im(\tau)>0$.
Hence the upper half plane
$\mathbb{C}^+=\{\tau\in\mathbb{C}\mid\Im(\tau)>0\}$
parameterizes the space of all conformal equivalence classes
of lattices up to the choice of two generators
with positive orientation.
The modular group $SL(2,\mathbb{Z})$ describes the transformation
of generators with positive orientation into each other.
Hence the space of conformal equivalence classes of lattices
is equal to the space of orbits of the action of the modular group
on the upper--half plane.
There exist many fundamental regions of this action,
but the following choice is the most common:
\index{conformal!class!moduli space $\moduli_1$ of $\sim$es}
\index{moduli space!$\moduli_1$ of conformal classes}
$$\moduli_1=\left\{\tau\in\mathbb{C}^+\mid
-1/2\leq \Im(\tau)\leq 1/2\text{ and } 1\leq |\tau|\right\}.$$
More precisely, the modular group is generated by the elements
$S=\left(\begin{smallmatrix}
0 & -1\\
1 & 0\end{smallmatrix}\right)$
and
$T=\left(\begin{smallmatrix}
1 & 1\\
0 & 1\end{smallmatrix}\right)$.
The corresponding action on the upper--half plane
is given by $S:\tau\mapsto-1/\tau$ and
$T:\tau\mapsto\tau+1$.
The first transformation induces an identification
of that part of the boundary of $\moduli_1$,
which belongs to the unit circle, with itself,
and the second induces an identification of that part of the boundary,
which belongs to the line $\Re(\tau)=-1/2$, with the part of the boundary,
which belongs to the line $\Re(\tau)=1/2$.
The corresponding dual lattice is generated by the elements
$\xx{\kappa}_1+\sqrt{-1}\xx{\kappa}_2=-\sqrt{-1}\tau/\Im(\tau)$
and $\yy{\kappa}_1+\sqrt{-1}\yy{\kappa}_2=\sqrt{-1}/\Im(\tau)$.
In the sequel we we will denote the
natural choice of generators corresponding to these moduli $\tau$
by $\xx{\gamma}_{\tau} =(1,0)$ and
$\yy{\gamma}_{\tau}=(\Re(\tau),\Im(\tau))$.
Analogously the corresponding generators of the dual lattice
are denoted by $\xx{\kappa}_{\tau}=(1,-\Re(\tau)/\Im(\tau))$ and
$\yy{\kappa}_{\tau}=(0,1/\Im(\tau))$.

\begin{description}
\item[All $\fermi_{\min}({\kappa/2})$ of genus zero:]
If $\Spa{Y}$ is a \Em{complex Fermi curve} of genus zero,
both components of $k$ extend to single--valued
meromorphic functions with first--order poles
at $\infty^-$ and $\infty^+$.
If $\Spa{Y}$ is invariant under the involution $\sigma$,
the function $g(k,k)$, extends to a global holomorphic function
and therefore has to be constant.
Since $\Spa{Y}$ is invariant under the involution
$\eta:k\mapsto -\Bar{k}$ without fixed points,
this constant has to be a non--negative real number.
The \Em{first integral} $\willmore$ is equal to
$4\pi^2\vol(\torus)$ times this constant,
and the \Em{complex Fermi curve} is given by the
equation
$g(k,k)=\willmore/(4\pi^2\vol(\torus)$.
Strictly speaking, the \Em{complex Fermi curve}
is the union of the solutions of this equation shifted by all
$\kappa\in\lattice\dual$.
We conclude from Corollary~\ref{simplified minimizers}
that we may obtain all $\fermi_{\min}(\kappa/2)$,
by finding a non--zero dual lattice vector
$\kappa\in\lattice\dual$ with the smallest length, and the
corresponding \Em{complex Fermi curve} is given by the equation
$g(k,k)=g(\kappa,\kappa)/4.$
Obviously, for all lattices corresponding to $\tau\in\moduli_1$
the second generator $\yy{\kappa}_{\tau}$ has minimal length
among all non--zero dual lattice vectors.
However, if $|\tau|=1$, then the first generator $\xx{\kappa}_{\tau}$
has the same length. We conclude that if $|\tau|=1$,
then both \Em{complex Fermi curves}
$\fermi_{\min}(\xx{\kappa}_{\tau}/2)$ and
$\fermi_{\min}(\yy{\kappa}_{\tau}/2)$
coincide and are given by the equation $g(k,k)=1/(4\Im^2(\tau))$.
In all other cases only $\fermi_{\min}(\yy{\kappa}_{\tau}/2)$
is given by this equation.
Furthermore, these are all \Em{complex Fermi curves} of genus zero
of the form $\fermi_{\min}(\kappa/2)$.
The corresponding values of the \Em{first integral}
are equal to $\pi^2/\Im(\tau)$.
This implies that there exists some $\fermi_{\min}(\kappa/2)$
of genus zero, which contains some other element of
$\lattice\dual/2\lattice\dual$, if and only if
$|\tau|=1$ and $\kappa$ is either
$\xx{\kappa}_{\tau}$ or $\yy{\kappa}_{\tau}$.
In this case the \Em{first integral} is equal to $\pi^2/\Im(\tau)$.
In particular, the minimum of the \Em{first integral}
of these cases is attained for $\tau=(0,1)$ and
the corresponding value of the \Em{first integral}
is equal to $\pi^2$.
Therefore, the lowest critical value of the
\De{Generalized Willmore functional}~\ref{generalized functional}
is at most $2\pi^2$.
Let us now see that this is in fact the absolute minimum.
\item[All $\fermi_{\min}({\kappa/2})$ of genus one:]
All these \Em{complex Fermi curves} are
\Em{real $\sigma$--hyperelliptic} and have exactly one
\Em{Orbit of type~2} over the element
$[\yy{\kappa}_{\tau}/2]\in\lattice\dual/2\lattice\dual$.
Moreover, $\infty^{\pm}$ are two other fixed points
of the involution $\sigma$ and these four fixed points are
all fixed points of the involution $\sigma$,
since due to the Hurwitz formula (see e.g. \cite[17.14]{Fo})
the involution $\sigma$ has four fixed points.
Moreover, there has to exist an \Em{Orbit of type~1}
over some other element of $\lattice\dual/2\lattice\dual$.
Therefore, let us first determine the moduli space of all
elliptic Riemann surfaces with two commuting involutions,
one of which is holomorphic ($\sigma)$
and the other is anti--holomorphic ($\rho$),
whose composition has no fixed points,
and besides the \Em{Orbits of type~2} over
$[\yy{\kappa}_{\tau}/2]$ and the points $\infty^{\pm}$
one \Em{Orbit of type~1}.
We describe any such Riemann surface as the quotient
of the complex plane $\{z\mid z\in \mathbb{C}\}$
by some two--dimensional lattice.
Since all translations induce some biholomorphic map of the
Riemann surface into itself, we may chose the origin
to correspond to the point $\infty^+$.
Hence the involution $\sigma$ is induced by the map $z\mapsto -z$.
Since the Riemann surface has \Em{dividing real part},
one generator of $H_1(\Spa{Y},\mathbb{Z})$ may be chosen
to be invariant under $\rho$ and the other to be anti--invariant.
These two generators of $H_1(\Spa{Y},\mathbb{Z})$ correspond
to generators of the lattice.
By multiplication by some non--zero complex number
we may always achieve that the invariant generator is purely imaginary,
the anti--invariant generator is equal to $1$,
and $\rho$ is induced by the map $z\mapsto -\Bar{z}+$ some constant.
Since $\rho$ and $\sigma$ commute,
twice this constant has to be a period,
and the condition that $\rho\comp\sigma$ has no fixed points
implies that this constant has to be $1/2$ times the real generator.
Hence the lattice is generated by a real positive period
usually denoted by $2\omega$ (which may be chosen to be $1$)
and an imaginary period usually denoted by $2\omega'$,
which may be chosen to have positive imaginary part.
The two involutions are induced by the maps $\sigma:z\mapsto-z$
and $\rho:z\mapsto -\Bar{z}+\omega$.
This implies that $\infty^-$ corresponds to the point $\omega$,
and the \Em{Orbit of type~2} over $[\yy{\kappa}_{\tau}/2]$
corresponds to the points $\omega'$ and $\omega+\omega'$.
The real part corresponds to the lines $\Re(z)=\pm\omega/2$.
Therefore, the \Em{Orbit of type~1} corresponds to two points
$z_1\in\omega/2+\sqrt{-1}\mathbb{R}/(2\omega'\mathbb{Z})$ and
$z_2=-z_1$. We conclude that the moduli space of such
Riemann surfaces is given by the quotient of the two generators
$\omega'/\omega\in\sqrt{-1}\mathbb{R}^+$ and some element $z_1$ of
$\omega/2+\sqrt{-1}\mathbb{R}/(2\omega'\mathbb{Z})$.
This is a real two--dimensional manifold.

The integral of $dk$ along all cycles of $H_1(\Spa{Y},\mathbb{Z})$,
which are anti--invariant under $\rho$ has to vanish.
Therefore the lattice $\lattice\dual_{\text{\scriptsize\rm effective}}$
is at most one--dimensional. On the other hand, since $\Spa{Y}$ has an
\Em{Orbit of type~2} over $[\kappa'/2]$,
there exist two generators $\xx{\kappa}$ and $\yy{\kappa}$
of $\lattice\dual$, such that
$\lattice\dual_{\text{\scriptsize\rm effective}}
=\Check{N}\yy{\kappa}\mathbb{Z}$, with some odd number $\Check{N}$.
Apriori these two generators are not necessarily the same as the generators
$\xx{\kappa}_{\tau}$ and $\yy{\kappa}_{\tau}$,
but we will see that they indeed are.
If $\xx{\gamma}$ and $\yy{\gamma}$ are the
dual generators of $\lattice$, then the function
$\xx{p}=g(\xx{\gamma},k)$ extends to a single--valued
meromorphic function with poles only at $0$ and $\omega$.
The properties $\sigma^{\ast}k=-k$ and $\rho^{\ast}\Bar{k}=k$ implies
that $\xx{p}$ may be expressed by the Weierstra{\ss} functions
(see e.g. \cite[13.12]{MOT})
$$\xx{p}=\alpha\left(\zeta(z)-\zeta(z+\omega)+\eta\right)=
\frac{-\alpha}{2} \frac{\wp'(z)}{\wp(z)-e_1},$$
with some non--zero real number $\alpha$.
Since $z_1$ and $-z_1$ correspond to an \Em{Orbit of type~1}
over some non--zero element of $\lattice\dual/2\lattice\dual$
different from $[\yy{\kappa}/2]$,
we conclude that $\alpha$ is given by
$$\alpha=\frac{\Hat{N}}{2\left(\zeta(z_1)-\zeta(z_1+\omega)+\eta\right)}=
-\Hat{N}\frac{\wp(z_1)-e_1}{\wp'(z_1)},$$
with some odd number $\Hat{N}$.
The conditions on the function $\yy{p}$ imply
that it may be expressed by the Weierstra{\ss} functions:
$$\yy{p}=\frac{\Check{N}}{ 2\pi\sqrt{-1}}\left(2\eta z
-\omega\left(\zeta(z)+\zeta(z+\omega)-\eta\right)
\right)+\beta p,$$
with some real $\beta$.
Again the condition that $z_1$ and $-z_1$
correspond to an \Em{Orbit of type~1}
over some non--zero element of  $\lattice\dual/2\lattice\dual$
different from $[\yy{\kappa}/2]$ implies that $\beta$ is given by
$$\beta=-\frac{\Check{N}}{\Hat{N}\pi\sqrt{-1}}\left(2\eta z_1
-\omega\left(\zeta(z_1)+\zeta(z_1+\omega)-\eta\right)
\right)+N,$$
with some integer $N$.
At the point $\infty^+$ the function $\yy{p}/\xx{p}$
takes the value $\tau$.
Hence the parameters $\tau$ of the conformal equivalence classes
of the lattice are given by
\begin{eqnarray*}
\Re{\tau}&=&\beta=-\frac{\Check{N}}{\Hat{N}\pi\sqrt{-1}}\left(2\eta z_1
-\omega\left(\zeta(z_1)+\zeta(z_1+\omega)-\eta\right)
\right)+N\\
\Im{\tau}&=&-\frac{\Check{N}\omega}{2\Hat{N}\pi}
\frac{\wp'(z_1)}{\wp(z_1)-e_1}
\end{eqnarray*}
Finally, due to Lemma~\ref{Willmore functional}
the \Em{first integral} is equal to the residue of the form
$8\pi^2\sqrt{-1}$ $\vol(\torus)k_2dk_1=
8\pi^2\sqrt{-1}\yy{p}d\xx{p}$ at $\infty^+$
and therefore equal to
$$\willmore=8\pi\alpha \Check{N}\left(\eta +e_1\omega\right)=
-8\pi\Hat{N}\Check{N}\frac{\wp(z_1)-e_1}{\wp'(z_1)}
\left(\eta +e_1\omega\right).$$
The function $\frac{\wp'(z)}{\wp(z)-e_1}$ is real on the lines
$\Im(z)=0,\Im(\omega')$ and on the lines $\Re(z)=\pm\omega/2$.
In particular, $d\xx{p}$ has four zeroes at the points
$\pm\omega/2$ and $\pm\omega/2+\omega'$.
Let us calculate the values of this function at these four points.
The square of this function is a rational function of $\wp$
$$\left(\frac{\wp'(z)}{\wp(z)-e_1}\right)^2=
4\frac{(\wp(z)-e_2)(\wp(z)-e_3)}{\wp(z)-e_1}.$$
Therefore, at the zeroes of $dp$ $\wp$ takes the values
$(\wp-e_1)^2=2e_1^2+e_2e_3$ and  more precisely
$\wp(\pm\omega/2)=e_1+\sqrt{2e_1^2+e_2e_3}$ and
$\wp(\pm\omega/2+\omega')=e_1-\sqrt{2e_1^2+e_2e_3}$.
Also the function $-\frac{\wp'(z)}{\wp(z)-e_1}$ takes at the points
$\pm\omega/2$ the values
$\pm\sqrt{\frac{
3e_1\sqrt{2e_1^2+e_2e_3}+(4e_1^2+e_2e_3)}{\sqrt{2e_1^2+e_2e_3}}}$
and at the points $\pm\omega/2+\omega'$ the values
$\pm\sqrt{\frac{
3e_1\sqrt{2e_1^2+e_2e_3}-(4e_1^2+e_2e_3)}{\sqrt{2e_1^2+e_2e_3}}}$.
If $\omega'/\omega\in\sqrt{-1}\mathbb{R}^+$,
then the numbers $e_1,e_2$ and $e_3$
satisfy the inequalities $e_3\leq e_2\leq e_1$.
Hence the restriction of the function $-\frac{\wp'(z)}{\wp(z)-e_1}$ to
$\Re(z)=\omega/2$ is positive and takes the maximum at $\omega/2$ and
the minimum at $\omega/2+\omega'$. Also the restriction of this
function to $\Re(z)=-\omega/2$ is negative and takes the minimum at
$-\omega/2$ and the maximum at $-\omega/2+\omega'$.
By definition of $\fermi_{\min}(\kappa/2)$
they have to be minimizers of the \Em{first integral}.
This implies that $\Hat{N}=1$ and $\Check{N}=1$.
On the other hand the function $x$ is equal to $N$ if $z_1=\omega/2$
and equal to $1+N$ if $z_1=\omega/2+\omega'$.

With the help of the limits of the Weierstra{\ss} functions
(see e.g. \cite[13.15]{MOT}) it is easy to see
that the limit $\omega'/\omega\rightarrow\infty\sqrt{-1}$ is
exactly the genus zero case on the intersection of the circle
$|\tau|=1$ with the upper half-plane.
Hence the Riemann surface is given by the equation
$g(k,k)=1/(4\Im^2(\tau))$ and the \Em{first integral} is equal to
$\willmore=\pi^2/\Im(\tau)$.
In the other limit $\omega'/\omega\rightarrow 0\sqrt{-1}$
the Riemann surface tends to the singular curve
of two copies of $\mathbb{P}^1$,
which are connected by two ordinary double points.
This is also the limit $y\rightarrow\infty$ of
the unique \Em{complex Fermi curves},
whose normalizations have two connected components,
and whose \Em{first integrals} are equal to $4\pi$.
In particular, in the limit $y\rightarrow\infty$
all the Riemann surfaces $\fermi_{\min}(\xx{\kappa}/2)$
tend to this unique singular curve,
and the \Em{first integral} converges to $4\pi$.
We conclude that with the choice
$\Hat{N}=1,\Check{N}=1,N=0$ the parameter space
$(\omega'/\omega,z_1)\in\sqrt{-1}\mathbb{R}^+\times
\left(\omega/2+\omega'[-1/2,1/2]\right)$ is mapped onto the domain
$\left\{\tau\in\mathbb{C}^+\mid\Re(\tau)\in[-1,1]\text{ and }\tau|>1\right\}$.
Not all of these \Em{complex Fermi curves}
are Riemann surfaces of the form
$\fermi_{\min}(\xx{\kappa}_{\tau}/2)$. But all those elements,
which are mapped into the fundamental domain $\moduli_1$,
are of this form. These are all \Em{complex Fermi curves}
of the form $\fermi_{\min}(\kappa/2)$ of genus one.
\end{description}

Let us now calculate the absolute minimum of the \Em{first integral}
of this family. For any value of the parameter $\omega'/\omega$, the
minimum of the \Em{first integral} is attained for $z_1=\omega/2$ and
the value of the \Em{first integral} at this point is equal to
$$\willmore=8\pi(\eta+e_1\omega)\sqrt{\frac{\sqrt{2e_1^2+e_2e_3}}{
3e_1\sqrt{2e_1^2+e_2e_3}+(4e_1^2+e_2e_3)}}.$$
This expression does not depend on the choice of $\omega$,
it depends only on $\omega'/\omega\in\sqrt{-1}\mathbb{R}^+$.
In the following lemma we will prove that this function
has no critical point for $\omega'/\omega\in\mathbb{R}^+\sqrt{-1}$.
As we have shown before, in the limit
$\omega'/\omega\rightarrow\infty\sqrt{-1}$
this function converges to $\pi^2$,
and in the limit $\omega'/\omega\downarrow 0\sqrt{-1}$ to $4\pi$.
Since this function is analytic, this implies
that it is strictly monotone decreasing and
that the infimum is attained in the first limit,
which is the genus zero case described above.
This implies that the restriction of the \Em{first integral}
to the space of
\De{Generalized Weierstra{\ss} curves}~\ref{generalized curves}
is not smaller than $2\pi^2$, and that the minimum is attained
by some unique \Em{complex Fermi curve},
which corresponds to the \Em{Clifford torus}, as we will see below.

\begin{Lemma} \label{no critical point}
The \Em{first integral} of the family of \Em{complex Fermi curves}
$\fermi_{\min}(\xx{\kappa}_{\tau}/2)$ of genus one described above,
has no critical point as a function depending on
$\{\tau\in\mathbb{C}^+\mid\Re(\tau)\in[-1,1]\text{ and }\tau|>1\}$.
\end{Lemma}

\begin{proof}
We use an extension of the deformations constructed in
Proposition~\ref{smooth moduli}. In fact,
if we consider $\yy{p}$ as a function depending on $\xx{p}$
and the moduli parameter, the form
$\omega=\frac{\partial \yy{p}}{\partial t}d\xx{p}$
corresponding to some change of the conformal class $(x,y)$
has poles of third--order at $\infty^{\pm}$. More precisely,
the singular part at $\infty^{\pm}$ is equal to
$$\left(\frac{\partial \Re(\tau)}{\partial t}
\pm\sqrt{-1}\frac{\partial \Im(\tau)}{\partial t}\right)\xx{p}d\xx{p}
\pm 8\pi^2\sqrt{-1}
\frac{\partial \willmore}{\partial t}\xx{p}^{-1}d\xx{p}.$$
We conclude that if we consider all
meromorphic differentials on the normalization
of the quotient modulo $\sigma$,
which have poles of second--order at $\infty^{\pm}$,
the analogous construction of Proposition~\ref{smooth moduli}
describes deformations of the conformal class.
Consequently, a relative minimizer of this family cannot have a real
(i.\ e.\ $\rho^{\ast}\omega=\eta^{\ast}\omega=\Bar{\omega}$)
meromorphic one--form on the normalization
of the quotient modulo $\sigma$,
which has at $\infty^{\pm}$ a second--order pole and
a non--vanishing residue and no other poles,
and which vanishes at the real point corresponding to the
\Em{Orbit of type~1}.
Since all \Em{complex Fermi curves} of this family are
\Em{real $\sigma$--hyperelliptic}
the normalizations of the quotient modulo $\sigma$
is isomorphic to $\mathbb{P}^1$ with one real cycle.
Obviously there always exists a meromorphic one--form
with the desired properties.
Hence this family does not have a
relative minimizer of the \Em{first integral}.
\end{proof}

\begin{Remark}\label{analogous deformations}
It is quite easy to describe the corresponding deformations
and the analogous deformations of the family of genus two
by ordinary differential equations in terms of the coefficients
of some polynomials, whose zero set describes the
normalization of the \Em{complex Fermi curves} of this family,
and the coefficients of some polynomials,
which describe the corresponding meromorphic one--forms.
This yields a numerical method to calculate the corresponding
\Em{complex Fermi curves} and their values of the
\Em{first integrals} (compare with \cite[Section~2.2]{GS1}).
\end{Remark}

\begin{Remark}\label{minimal tori}
The analogous calculation for \Em{real--$\sigma$--hyperelliptic}
\De{Generalized Weierstra{\ss} curves}~\ref{generalized curves}
with two \Em{Orbits of type~1} shows, that for Willmore tori
(i.\ e.\ critical points with respect to deformations changing
the conforam class)
the normalization of the quotient modulo $\sigma$
is isomorphic to $\mathbb{P}^1$,
with the two marked points $z=0$ and $z=\infty$,
the anti--holomorphic involution $z\mapsto 1/\Bar{z}$,
and the two real points $z=\pm 1$
corresponding to the two \Em{Orbits of type~1}.
The corresponding \Em{complex Fermi curves} are exactly the
\De{Weierstra{\ss} curves}~\ref{Weierstrass curves}
of minimal tori in $S^3$ described in \cite{Hi}.
\end{Remark}

\index{minimizer!of the generalized Willmore functional|(}
In conclusion, we want to describe the mappings
$\lattice\mapsto\lattice_{[\kappa/2]}$
for the three non--zero elements $[\kappa/2]$ of
$\lattice\dual/2\lattice\dual$ on $\moduli_1$
and the corresponding
\De{Generalized Weierstra{\ss} curves}~\ref{generalized curves}.
In the following list we state explicitly the
conformal equivalence class of $\lattice_{[\kappa/2]}$ in dependence of
the conformal equivalence class of $\lattice$,
considered as a mapping from domains in $\moduli_1$
into domains in $\moduli_1$.
Moreover, we state the corresponding elements
in the preimage of $[\kappa/2]$ under the canonical map
$\lattice\dual_{[\kappa/2]}/2\lattice\dual_{[\kappa/2]}
\rightarrow\lattice\dual/2\lattice\dual$.
Furthermore, we state the corresponding \Em{complex Fermi curves}
with lattice $\lattice_{[\kappa/2]}$
introduced in Corollary~\ref{simplified minimizers},
to which the
\De{Generalized Weierstra{\ss} curve}~\ref{generalized curves}.
$\fermi_{\min,[\kappa/2]}(\willmore_{[\kappa/2]}/2)$
corresponds.
Finally, we give the corresponding
\Em{geometric genus}.
Due to our calculations above the \Em{geometric genus} of
$\fermi_{\min}(\yy{\kappa}_{\tau}/2)$ is zero,
the \Em{geometric genus} of
$\fermi_{\min}(\xx{\kappa}_{\tau}/2)$ is one
if $|\tau|>1$, and finally the \Em{geometric genus} of
$\fermi_{\min}(\xx{\kappa}_{\tau}/2+\yy{\kappa}_{\tau}/2)$
is two if $\Re(\tau)\neq -1/2$ and $\Re(\tau)\neq 1/2$.

The last two cases describe the neighbourhood of the
\Em{Clifford torus}. For all rectangular conformal classes
(i.\ e.\ $\Re(\tau)=0$) the minimizers corresponding to
$[\xx{\kappa}_{\tau}/2+\yy{\kappa}_{\tau}/2]$ have two
\Em{Orbits of type~1}, all other minimizers have
one \Em{Orbit of type~1} and one \Em{Orbit of type~2}.
We shall see in the next section,
that only this part of the family corresponds to
\De{Weierstra{\ss} curves}~\ref{Weierstrass curves}.
But due to the generalization of
the Weierstra{\ss} representation to immersions into
$\mathbb{R}^4$, the whole family corresponds to
\Em{quaternionic Weierstra{\ss}
  curves}\index{Weierstra{\ss}!quaternionic $\sim$
  curve}\index{curve!quaternionic Weierstra{\ss}
  $\sim$}\index{quaternionic!Weierstra{\ss} curve},
i.\ e.\ to immersions into $\mathbb{R}^4$
(compare with Lemma~\ref{existence of potentials} and \cite{PP}).
\index{minimizer!of the generalized Willmore functional|)}

We conjecture that for all conformal classes our lower bound for the
Willmore functional on tori in $\mathbb{R}^3$ in dependence of the
conformal class describes the absolute minimum of the
restriction of the Willmore functional to tori in $\mathbb{R}^4$
of the given conformal class.\vfill

\begin{enumerate}
\item $[\xx{\kappa}_{\tau}/2]$:

\noindent
\begin{minipage}[t]{10cm}
\begin{enumerate}
\item For $|\tau|\geq 2$:
\begin{itemize}
\item $\tau\mapsto\tau'=\tau/2$.
\item The preimage of $[\xx{\kappa}_{\tau}/2]$ contains
$[\xx{\kappa}_{\tau'}/2]$ and
$[\xx{\kappa}_{\tau'}/2+\yy{\kappa}_{\tau'}/2]$.
\item $\fermi_{\min,[\xx{\kappa}_{\tau}/2]}
\left(\willmore_{[\xx{\kappa}_{\tau}/2]}/2\right)$
corresponds to $\fermi_{\min}
\left({\xx{\kappa}_{\tau'}/2+\yy{\kappa}_{\tau'}/2}\right)$.
\item They have \Em{geometric genus} two.
\end{itemize}
\item For $|\tau|\leq 2$, $|\tau+2|\geq 2$ and
  $|\tau-2|\geq 2$:
\begin{itemize}
\item $\tau\mapsto\tau'=-2/\tau$.
\item The preimage of $[\xx{\kappa}_{\tau}/2]$ contains
$[\yy{\kappa}_{\tau'}/2]$ and
$[\xx{\kappa}_{\tau'}/2+\yy{\kappa}_{\tau'}/2]$.
\item $\fermi_{\min,[\xx{\kappa}_{\tau}/2]}
\left(\willmore_{[\xx{\kappa}_{\tau}/2]}/2\right)$
corresponds to $\fermi_{\min}
\left({\xx{\kappa}_{\tau'}/2+\yy{\kappa}_{\tau'}/2}\right)$.
\item If $|\tau+2|\neq 2$ and $|\tau-2|\neq 2$,
then they have \Em{geometric genus} two, otherwise \Em{geometric genus} one.
\end{itemize}
\item For $|\tau+2|\leq 2$:
\begin{itemize}
\item $\tau\mapsto\tau'=-2/\tau-1$.
\item The preimage of $[\xx{\kappa}_{\tau}/2]$ contains
$[\xx{\kappa}_{\tau'}/2]$ and $[\yy{\kappa}_{\tau'}/2]$.
\item $\fermi_{\min,[\xx{\kappa}_{\tau}/2]}
\left(\willmore_{[\xx{\kappa}_{\tau}/2]}/2\right)$
corresponds to $\fermi_{\min}
\left({\xx{\kappa}_{\tau'}/2}\right)$.
\item They have \Em{geometric genus} one.
\end{itemize}
\item For $|\tau-2|\leq 2$:
\begin{itemize}
\item $\tau\mapsto\tau'=-2/\tau+1$.
\item The preimage of $[\xx{\kappa}_{\tau}/2]$ contains
$[\xx{\kappa}_{\tau'}/2]$ and $[\yy{\kappa}_{\tau'}/2]$.
\item $\fermi_{\min,[\xx{\kappa}_{\tau}/2]}
\left(\willmore_{[\xx{\kappa}_{\tau}/2]}/2\right)$
corresponds to $\fermi_{\min}
\left({\xx{\kappa}_{\tau'}/2}\right)$.
\item They have \Em{geometric genus} one.
\end{itemize}
\end{enumerate}
\end{minipage}\hfill
{\setlength{\unitlength}{3cm}
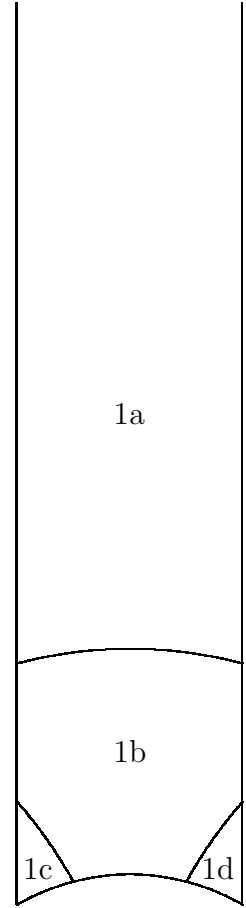
\begin{picture}(1,-6)
\put(0,-6){\begin{picture}(1,5)
\put(-.5,.866){\line(0,1){4}}
\put(.5,.866){\line(0,1){4}}
\qbezier(-.5,.866)(-.382,.934)(-.25,.968)
\qbezier(-.25,.968)(0,1.033)(.25,.968)
\qbezier(.25,.968)(.382,.934)(.5,.866)
\put(-.07,3){1a}
\qbezier(-.5,1.936)(0,2.066)(.5,1.936)
\put(-.07,1.5){1b}
\qbezier(-.5,1.323)(-.356,1.159)(-.25,.968)
\put(-.47,.98){1c}
\qbezier(.5,1.323)(.356,1.159)(.25,.968)
\put(.32,.98){1d}
\end{picture}}
\end{picture}}
\vfill

\item $[\yy{\kappa}_{\tau}/2]$:

\noindent
\begin{minipage}[t]{10cm}
\begin{enumerate}
\item For $-1/4\leq \Re(\tau)\leq 1/4$:
\begin{itemize}
\item $\tau\mapsto\tau'=2\tau$.
\item The preimage of $[\yy{\kappa}_{\tau}/2]$ contains
$[\yy{\kappa}_{\tau'}/2]$ and
$[\xx{\kappa}_{\tau'}/2+\yy{\kappa}_{\tau'}/2]$.
\item $\fermi_{\min,[\yy{\kappa}_{\tau}/2]}
\left(\willmore_{[\yy{\kappa}_{\tau}/2]}/2\right)$
corresponds to $\fermi_{\min}
\left({\xx{\kappa}_{\tau'}/2+\yy{\kappa}_{\tau'}/2}\right)$.
\item If $\Re(\tau)\neq\pm 1/4$,
then they have \Em{geometric genus} two, otherwise \Em{geometric genus} one.
\end{itemize}
\item For $\Re(\tau)\leq -1/4$:
\begin{itemize}
\item $\tau\mapsto\tau'=2\tau+1$.
\item The preimage of $[\yy{\kappa}_{\tau}/2]$ contains
$[\xx{\kappa}_{\tau'}/2]$ and $[\yy{\kappa}_{\tau'}/2]$.
\item $\fermi_{\min,[\yy{\kappa}_{\tau}/2]}
\left(\willmore_{[\yy{\kappa}_{\tau}/2]}/2\right)$
corresponds to $\fermi_{\min}
\left(\xx{\kappa}_{\tau'}/2\right)$.
\item They have \Em{geometric genus} one.
\end{itemize}
\item For $1/4\leq \Re(\tau)$:
\begin{itemize}
\item $\tau\mapsto\tau'=2\tau-1$.
\item The preimage of $[\yy{\kappa}_{\tau}/2]$ contains
$[\xx{\kappa}_{\tau'}/2]$ and $[\yy{\kappa}_{\tau'}/2]$.
\item $\fermi_{\min,[\yy{\kappa}_{\tau}/2]}
\left(\willmore_{[\yy{\kappa}_{\tau}/2]}/2\right)$
corresponds to $\fermi_{\min}
\left(\xx{\kappa}_{\tau'}/2\right)$.
\item They have \Em{geometric genus} one.
\end{itemize}
\end{enumerate}
\end{minipage}\hfill
{\setlength{\unitlength}{3cm}
\begin{picture}(1,-5)
\put(0,-5){\begin{picture}(1,5)
\put(-.5,.866){\line(0,1){4}}
\put(.5,.866){\line(0,1){4}}
\qbezier(-.5,.866)(-.382,.934)(-.25,.968)
\qbezier(-.25,.968)(0,1.033)(.25,.968)
\qbezier(.25,.968)(.382,.934)(.5,.866)
\put(-.07,2.5){2a}
\put(-.25,.968){\line(0,1){3.89}}
\put(.25,.968){\line(0,1){3.89}}
\put(-.45,2.5){2b}
\put(.30,2.5){2c}
\end{picture}}
\end{picture}}\vfill

\noindent
\item $[\xx{\kappa}_{\tau}/2+\yy{\kappa}_{\tau}/2]$:
\begin{enumerate}
\item For $2\leq|\tau+1|$, $2\leq|\tau-1|$ and $\Re(\tau)\leq 0$:
\begin{itemize}
\item $\tau\mapsto\tau'=(\tau+1)/2$.
\item The preimage of $[\xx{\kappa}_{\tau}/2+\yy{\kappa}_{\tau}/2]$
contains $[\yy{\kappa}_{\tau'}/2]$ and
$[\xx{\kappa}_{\tau'}/2+\yy{\kappa}_{\tau'}/2]$.
\item $\fermi_{\min,[\xx{\kappa}_{\tau}/2+\yy{\kappa}_{\tau}/2]}
\left(\willmore_{[\xx{\kappa}_{\tau}/2+\yy{\kappa}_{\tau}/2]}/2\right)$
corresponds to $\fermi_{\min}
\left(\xx{\kappa}_{\tau'}/2+\yy{\kappa}_{\tau'}/2\right)$.
\item If $\Re(\tau)\neq 0$, then they have \Em{geometric genus} two,
otherwise \Em{geometric genus} one.
\end{itemize}
\item For $2\leq|\tau+1|$, $2leq|\tau-1|$ and $0\leq\Re(\tau)$:
\begin{itemize}
\item $\tau\mapsto\tau'=(\tau-1)/2$.
\item The preimage of $[\xx{\kappa}_{\tau}/2+\yy{\kappa}_{\tau}/2]$
contains $[\yy{\kappa}_{\tau'}/2]$ and
$[\xx{\kappa}_{\tau'}/2+\yy{\kappa}_{\tau'}/2]$.
\item $\fermi_{\min,[\xx{\kappa}_{\tau}/2+\yy{\kappa}_{\tau}/2]}
\left(\willmore_{[\xx{\kappa}_{\tau}/2+\yy{\kappa}_{\tau}/2]}/2\right)$
corresponds to $\fermi_{\min}
\left(\xx{\kappa}_{\tau'}/2+\yy{\kappa}_{\tau'}/2\right)$.
\item If $\Re(\tau)\neq 0$, then they have \Em{geometric genus} two,
otherwise \Em{geometric genus} one.
\end{itemize}

\noindent
\begin{minipage}[t]{9cm}
\item For $|\tau+1|\leq 2$ and $2\leq |\tau-1|$:
\begin{itemize}
\item $\tau\mapsto\tau'=-2/(\tau+1)$.
\item The preimage of $[\xx{\kappa}_{\tau}/2+\yy{\kappa}_{\tau}/2]$
contains $[\xx{\kappa}_{\tau'}/2]$ and
$[\xx{\kappa}_{\tau'}/2+\yy{\kappa}_{\tau'}/2]$.
\item $\fermi_{\min,[\xx{\kappa}_{\tau}/2+\yy{\kappa}_{\tau}/2]}
\left(\willmore_{[\xx{\kappa}_{\tau}/2+\yy{\kappa}_{\tau}/2]}/2\right)$
corresponds to $\fermi_{\min}
\left(\xx{\kappa}_{\tau'}/2+\yy{\kappa}_{\tau'}/2\right)$.
\item If $|\tau-1|\neq 2$, then they have \Em{geometric genus} two,
otherwise \Em{geometric genus} one.
\end{itemize}
\item For $2\leq|\tau+1|$ and $|\tau-1|\leq 2$:
\begin{itemize}
\item $\tau\mapsto\tau'=-2/(\tau-1)$.
\item The preimage of $[\xx{\kappa}_{\tau}/2+\yy{\kappa}_{\tau}/2]$
contains $[\xx{\kappa}_{\tau'}/2]$ and
$[\xx{\kappa}_{\tau'}/2+\yy{\kappa}_{\tau'}/2]$.
\item $\fermi_{\min,[\xx{\kappa}_{\tau}/2+\yy{\kappa}_{\tau}/2]}
\left(\willmore_{[\xx{\kappa}_{\tau}/2+\yy{\kappa}_{\tau}/2]}/2\right)$
corresponds to $\fermi_{\min}
\left(\xx{\kappa}_{\tau'}/2+\yy{\kappa}_{\tau'}/2\right)$.
\item If $|\tau+1|\neq 2$, then they have \Em{geometric genus} two,
otherwise \Em{geometric genus} one.
\end{itemize}
\item For $|\tau+1|\leq 2$, $|\tau-1|\leq 2$ and $\Re(\tau)\leq 0$:
\begin{itemize}
\item $\tau\mapsto\tau'=(\tau-1)/(\tau+1)$.
\item The preimage of $[\xx{\kappa}_{\tau}/2+\yy{\kappa}_{\tau}/2]$
contains $[\xx{\kappa}_{\tau'}/2]$ and $[\yy{\kappa}_{\tau'}/2]$.
\item $\fermi_{\min,[\xx{\kappa}_{\tau}/2+\yy{\kappa}_{\tau}/2]}
\left(\willmore_{[\xx{\kappa}_{\tau}/2+\yy{\kappa}_{\tau}/2]}/2\right)$
corresponds to $\fermi_{\min}
\left(\xx{\kappa}_{\tau'}/2\right)$.
\item If $\Re(\tau)\neq 0$, then they have \Em{geometric genus} one,
otherwise \Em{geometric genus} zero.
\end{itemize}
\end{minipage}\hfill
{\setlength{\unitlength}{3cm}
\begin{picture}(1,-5)
\put(0,-5){\begin{picture}(1,5)
\put(-.5,.866){\line(0,1){4}}
\put(.5,.866){\line(0,1){4}}
\qbezier(-.5,.866)(-.382,.934)(-.25,.968)
\qbezier(-.25,.968)(0,1.033)(.25,.968)
\qbezier(.25,.968)(.382,.934)(.5,.866)
\put(0,1){\line(0,1){3.83}}
\put(-.33,3){3a}
\put(.17,3){3b}
\qbezier(-.5,1.936)(-.236,1.868)(0,1.732)
\qbezier(.5,1.936)(.236,1.868)(0,1.732)
\qbezier(-.5,1.323)(-.283,1.568)(0,1.732)
\qbezier(.5,1.323)(.283,1.568)(0,1.732)
\put(-.41,1.65){3c}
\put(.25,1.65){3d}
\put(-.3,1.25){3e}
\put(.14,1.25){3f}
\end{picture}}
\end{picture}}
\item For $|\tau+1|\leq 2$, $|\tau-1|\leq 2$ and $0\leq\Re(\tau)$:
\begin{itemize}
\item $\tau\mapsto\tau'=-(\tau+1)/(\tau-1)$.
\item The preimage of $[\xx{\kappa}_{\tau}/2+\yy{\kappa}_{\tau}/2]$
contains $[\xx{\kappa}_{\tau'}/2]$ and $[\yy{\kappa}_{\tau'}/2]$.
\item $\fermi_{\min,[\xx{\kappa}_{\tau}/2+\yy{\kappa}_{\tau}/2]}
\left(\willmore_{[\xx{\kappa}_{\tau}/2+\yy{\kappa}_{\tau}/2]}/2\right)$
corresponds to $\fermi_{\min}
\left(\xx{\kappa}_{\tau'}/2\right)$.
\item If $\Re(\tau)\neq 0$, then they have \Em{geometric genus} one,
otherwise \Em{geometric genus} zero.
\end{itemize}
\end{enumerate}
\end{enumerate}

\section{\De{Weierstra{\ss} potentials}}\label{section weierstrass}
\subsection{The \De{Singularity condition}}\label{subsection singularity}

In Lemma~\ref{fourth singularity} we proved
that the \Em{complex Fermi curve} of some potential,
which corresponds to some immersion has to fulfill the
\De{Weak Singularity condition}~\ref{weak singularity condition}.
But this was only a necessary
condition on the \Em{complex Fermi curve} to be a
\De{Weierstra{\ss} curve}~\ref{Weierstrass curves}.
In this section we want to investigate two problems.
\begin{description}
\item[First problem:] Find a necessary and sufficient condition on
  the \Em{complex Fermi curve} to be a
  \De{Weierstra{\ss} curve}~\ref{Weierstrass curves}.
\item[Second problem:] Given a
  \De{Weierstra{\ss} curve}~\ref{Weierstrass curves}, find a
  necessary and sufficient condition on a potential, which corresponds
  to the given \Em{complex Fermi curve}, to be a
  \De{Weierstra{\ss} potential}~\ref{Weierstrass potentials}.
  Section~\ref{subsection spectral projections} suggests a
  decomposition of all potentials corresponding to a given
  \Em{complex Fermi curve} into potentials,
  which corresponds to some
  one--sheeted coverings of \Em{complex Fermi curve}.
  In particular, we are interested in
  the question, which of these subclasses of potentials contains no
  \De{Weierstra{\ss} potentials}~\ref{Weierstrass potentials}, which
  of these contains some
  \De{Weierstra{\ss} potential}~\ref{Weierstrass potentials}, and
  which of these subclasses are contained in the set of
  \De{Weierstra{\ss} potentials}~\ref{Weierstrass potentials}.
\end{description}
If $U$ is a potential, whose \Em{complex Fermi curve} fulfills the
\Em{ Weak Singularity condition}~\ref{weak singularity condition},
then at least one of the following cases takes place:
\begin{enumerate}
\item \label{cusp} The preimage of $[\kappa/2]$ contains at least
    two \Em{cusps}. (Due to Corollary~\ref{fixed points}~(iii)
    the number of such elements has to be even.)
\item \label{double point} The preimage of $[\kappa/2]$
  contains two different elements $y$ and $y'$,
  such that the values of the eigenfunction at these two points
  are linearly dependent.
  (Due to Corollary~\ref{fixed points}~(ii), in this case $\eta(y)$
  and $\eta(y')$ are two other points with this property.)
\item \label{non-singular} The preimage of $[\kappa/2]$ contains only
  \Em{multiple points} and the values of the eigenfunction
  at all elements of the preimage are pairwise linearly independent.
  We may assume that $d\xx{p}$ has no zero in the preimage of
  $[\kappa/2]$. 
  \begin{enumerate}
  \item \label{type 3} The preimage contains an
    \Em{Orbit of type~3}.
  \item \label{3 type 1 and 2} The preimage contains alltogether at
    least three \Em{Orbits of type~1} or \Em{type 2}.
  \item \label{2 type 1} The preimage contains exactly two
    \Em{Orbits of type~1} and no other orbit.
  \item \label{2 type 2} The preimage contains exactly two
    \Em{Orbits of type~2} and no other orbit.
  \item \label{2 type 1 and 2} The preimage contains exactly one
    \Em{Orbit of type~1}, one \Em{Orbit of type~2}, and no
    \Em{Orbit of type~3}.
  \end{enumerate}
\end{enumerate}
In cases~\ref{cusp} and case~\ref{double point}
the preimage contains a singularity of the
\De{Structure sheaf}~\ref{structure sheaf}
in contrast to case~\ref{non-singular}.
Both former cases can take place simultaneously.
Also case~\ref{type 3} and case~\ref{3 type 1 and 2} can take place
simultaneously.
In case~\ref{cusp} $y$ is a branch point
with respect to the coverings over $p\in\mathbb{C}$ and
$p'\in\mathbb{C}$ described in the proof of
Lemma~\ref{fourth singularity}.
Hence due to Lemma~\ref{branchpoints} the function $\psi(y)$ obeys the
\De{Equivalent form of the Periodicity condition}~\ref{equivalent form}.
Analogously Lemma~\ref{branchpoints} shows that in
case~\ref{double point} the  value of the eigenfunction at $y$, and in
case~\ref{type 3} the values of the eigenfunction at all four points
of the \Em{Orbit of type~3} fulfill the
\De{Equivalent form of the Periodicity condition}~\ref{equivalent form}.

\begin{Lemma} \label{type 1 and 2}
\begin{description}
\item[(i)] In case~\ref{3 type 1 and 2} there always exists some
  eigenfunction, which fulfills the
  \De{Equivalent form of the Periodicity condition}~\ref{equivalent form}.
\item[(ii)] In case~\ref{2 type 1}--\ref{2 type 1 and 2} there
  exists some eigenfunction, which fulfills the
  \De{Equivalent form of the Periodicity condition}~\ref{equivalent form},
  if and only if the functions
  $d\yy{p}/d\xx{p}$ take the same value at two different elements
  of the preimage of $[\kappa/2]$, which belong to different orbits.
\end{description}
\end{Lemma}

\begin{proof} In case~\ref{non-singular} there always exists two
  generators $\xx{\gamma}$ and $\yy{\gamma}$ of the
  period lattice $\lattice$,
  whose differentials $d\xx{p}=dg(\xx{\gamma},k)$ and
  $d\yy{p}=dg(\yy{\gamma},k)$
  have no zeroes in the preimage of $[\kappa/2]$. Thus, due to
  Lemma~\ref{projection 2}, the ranges of the projections
  $\Breve{\Op{P}}_{\xx{\gamma}}([\kappa/2)$ and
  $\Breve{\Op{P}}_{\yy{\gamma}}([\kappa/2)$
  are spanned by the values of the eigenfunction at the elements of the
  preimage of $[\kappa/2]$. If $y_1$ and $y_2=\eta(y_1)$
  are two elements of an \Em{Orbit of type~1} or \Em{type 2}, then
  let $\psi_1$ and $\psi_2$ denote the corresponding values of the
  eigenfunction. Due to Corollary~\ref{fixed points} we may assume that
  $\psi_2$ is equal to $-\Op{J}\Bar{\psi}_1$, and therefore $\psi_1$ is
  equal to $\Op{J}\Bar{\psi}_2$. The bilinear forms
  $\langle\langle \Op{J}\psi_i,\psi_j\rangle\rangle_{\xx{\gamma}}$ and
  $\langle\langle \Op{J}\psi_i,\psi_j\rangle\rangle_{\yy{\gamma}}$
  form two $2\times2$ matrices. These bilinear
  forms are symmetric and transforms under the anti--unitary map
  $\chi\mapsto \Op{J}\Bar{\chi}$ like
  $$\langle\langle \Op{J}\Bar{\chi},
  \Op{J}\Bar{\xi}\rangle\rangle_{\xx{\gamma}}
  = -\overline{\langle\langle\chi,\xi\rangle\rangle}_{\xx{\gamma}}
  \text{ and }
  \langle\langle \Op{J}\Bar{\chi},
  \Op{J}\Bar{\xi}\rangle\rangle_{\yy{\gamma}}
  = -\overline{\langle\langle\chi,\xi\rangle\rangle}_{\yy{\gamma}}.$$
  Hence these matrices have the form
  $\left(\begin{smallmatrix}
  \alpha & \beta\\
  \beta & -\Bar{\alpha}
  \end{smallmatrix}\right)$ and
  $\left(\begin{smallmatrix}
  \alpha' & \beta'\\
  \beta' & -\Bar{\alpha'}
  \end{smallmatrix}\right)$, with complex numbers
  $\alpha,\beta,\alpha'$ and $\beta'$, respectively. Moreover, the
  matrices, whose entries are given as
  $\langle\langle\Bar{\psi}_i,\psi_j\rangle\rangle_{\xx{\gamma}}$ and
  $\langle\langle\Bar{\psi}_i,\psi_j\rangle\rangle_{\yy{\gamma}}$,
  are of the form
  $\left(\begin{smallmatrix}
  -\beta & \Bar{\alpha}\\
  \alpha & \beta
  \end{smallmatrix}\right)$ and
  $\left(\begin{smallmatrix}
  -\beta' & \Bar{\alpha}'\\
  \alpha' & \beta'
  \end{smallmatrix}\right)$, respectively. Since these bilinear forms
  are hermitian, $\beta$ and $\beta'$ have to be real. Due to
  Lemma~\ref{involutions} and Lemma~\ref{branchpoints} the numbers
  $\alpha$ and $\alpha'$ have to be zero if we consider an
  \Em{Orbit of type~1}, and the numbers $\beta$ and $\beta'$ have to
  be zero if we consider an \Em{Orbit of type~2}. Due
  to Lemma~\ref{projection 2} the other numbers are not equal to
  zero. Due to the proof of property~(iii) of Lemma~\ref{projection 2}
  the value of the function $d\yy{p}/d\xx{p}$
  at the point $y_1$ has to be
  equal to $\beta'/\beta$ if we consider an \Em{Orbit of type~1}, and
  equal to $\alpha'/\alpha$, if we consider an \Em{Orbit of type~2}.
  If in the first case we change the basis $\psi_1,\psi_2$ into
  $\frac{\psi_1+\psi_2}{\sqrt{2}},\frac{\psi_2-\psi_1}{\sqrt{2}}$,
  then the matrices corresponding to the bilinear forms
  $\langle\langle \Op{J}\cdot,\cdot\rangle\rangle_{\xx{\gamma}}$,
  $\langle\langle \Op{J}\cdot,\cdot\rangle\rangle_{\yy{\gamma}}$,
  $\langle\langle \Bar{\cdot},\cdot\rangle\rangle_{\xx{\gamma}}$ and
  $\langle\langle \Bar{\cdot},\cdot\rangle\rangle_{\yy{\gamma}}$
  transforms to 
  $\left(\begin{smallmatrix}
  \beta & 0\\
  0 & -\beta
  \end{smallmatrix}\right)$,
  $\left(\begin{smallmatrix}
  \beta' & 0\\
  0 & -\beta'
  \end{smallmatrix}\right)$,
  $\left(\begin{smallmatrix}
  0 & \beta\\
  \beta & 0
  \end{smallmatrix}\right)$ and
  $\left(\begin{smallmatrix}
  0 & \beta'\\
  \beta' & 0
  \end{smallmatrix}\right)$.
  As well for \Em{Orbits of type~1} as for \Em{Orbits of type~2}
  after some renormalization these matrices may be brought into the form
  $$\begin{pmatrix}
  1 & 0\\
  0 & -1
  \end{pmatrix}, \;
  \begin{pmatrix}
  \frac{d\yy{p}}{d\xx{p}}(y_1) & 0\\
  0 & -\overline{\frac{d\yy{p}}{d\xx{p}}}(y_1)
  \end{pmatrix}, \;
  \begin{pmatrix}
  0 & 1\\
  1 & 0
  \end{pmatrix}, \text{ and }
  \begin{pmatrix}
  0 & \overline{\frac{d\yy{p}}{d\xx{p}}}(y_1)\\
  \frac{d\yy{p}}{d\xx{p}}(y_1) & 0
  \end{pmatrix}.$$ Moreover, $\overline{\frac{d\yy{p}}{d\xx{p}}}(y_1)$ is
  equal to $\frac{d\yy{p}}{d\xx{p}}(y_2)$, and both numbers are real and
  therefore equal if we consider an \Em{Orbit of type~1}.
  Now let us prove statement~(ii). For arbitrary
  non--zero complex numbers $\alpha$ and $\beta$,
  which are the values of $d\yy{p}/d\xx{p}$
  at two different elements belonging to two different orbits,
  we have to determine the solutions of the following four equations
  with complex numbers $z_1,z_2,z_3$ and $z_4$:
  \begin{align*}
  (z_1\Bar{z}_2+\Bar{z}_1z_2)+(z_3\Bar{z}_4+\Bar{z}_3z_4) &=0
  &(z_1^2-z_2^2)+(z_3^2-z_4^2) &=0\\
  (\alpha z_1\Bar{z}_2+\Bar{\alpha}\Bar{z}_1z_2) +
  (\beta z_3\Bar{z}_4+\Bar{\beta}\Bar{z}_3z_4) &=0
  & (\alpha z_1^2-\Bar{\alpha}z_2^2) +
  (\beta z_3^2-\Bar{\beta} z_4^2) &=0
  \end{align*}
  These equations imply
  \begin{eqnarray*}
  (\alpha-\Bar{\alpha})^2(z_1\Bar{z}_2)^2 &=&
  \left((\Bar{\alpha}-\beta)z_3\Bar{z}_4+
  (\Bar{\alpha}-\Bar{\beta})\Bar{z}_3z_4\right)^2 \\
  &=&-\left((\Bar{\alpha}-\beta)z_3^2-
  (\Bar{\alpha}-\Bar{\beta})z_4^2\right)
  \left((\Bar{\alpha}-\Bar{\beta})\Bar{z}_3^2-
  (\Bar{\alpha}-\beta)\Bar{z}_4^2\right).
  \end{eqnarray*}
  The difference of these two equations yields
  $$(\Bar{\alpha}-\beta)(\Bar{\alpha}-\Bar{\beta})
  (z_3\Bar{z}_3+z_4\Bar{z}_4)^2=0.$$
  Thus either $\alpha$ is equal to $\beta$ or $\alpha$ is equal to
  $\Bar{\beta}$ or the only solution is the trivial solution
  (i.\ e.\ the complex numbers $z_1,z_2,z_3$ and $z_4$ are all zero).
  If $\alpha$ is equal to $\beta$,
  then $z_3=\sqrt{-1}z_1$ and $z_4=-\sqrt{-1}z_2$ is a solution
  for arbitrary $z_1$ and $z_2$.
  If $\alpha$ is equal to $\Bar{\beta}$,
  then $z_3=z_2$ and $z_4=-z_1$ is a solution
  for arbitrary $z_1$ and $z_2$. This proves (ii).
  The proof of (i) is similar.
\end{proof}

Whether case~\ref{cusp} takes place or not depends only on the
\Em{complex Fermi curve} and not on the choice of the potential,
in contrast, whether case~\ref{double point} takes place
depends on the choice of the potential.
In the rest of this section we want to find some property
of the \Em{complex Fermi curve},
which is equivalent to the existence of some potential,
which leads to the given \Em{complex Fermi curve}
and which is the potential of some immersion.
To all subcases of case~\ref{non-singular}
there might exist an analogous subcase of case~\ref{double point}.
Thus the mains question is,
whether to a potential of
cases~\ref{2 type 1}--\ref{2 type 1 and 2}
there exists another potential,
which leads to the same \Em{complex Fermi curve},
but corresponds to the analogous subcase of case~\ref{double point}.
Fortunately there exists some method to transform
subcases of case~\ref{non-singular}
into subcases of case~\ref{double point},
the so-called \Em{B\"acklund transformation}
\index{B\"acklund transformation}
(\cite{EK},\cite[Chapter~6.]{MS} and \cite[Section~4]{LMcL}).

\begin{Lemma} \label{Baecklund}
\begin{description}
\item[(i)] To all potentials of
  cases~\ref{2 type 1} there exists another potential,
  which leads to the same \Em{complex Fermi curve}
  and corresponds to the analogous subcase of case~\ref{double point}.
\item[(ii)] If the preimage of $[\kappa/2]$ contains exactly four
  \Em{multiple points}, at which the function $d\yy{p}/d\xx{p}$ takes
  more than two values, then there exists no immersion,
  which leads to this \Em{complex Fermi curve}.
\end{description}
\end{Lemma}

\begin{proof}
If in case~(i) the values of the eigenfunction $\psi$ at the
four elements of the two \Em{Orbits of type~1} are linearly independent,
then one of the B\"acklund transformations described in
Lemma~\ref{Baecklund transformation} has the desired properties.

The condition in case~(ii) excludes
\Em{Multiple points O1/O1} and \Em{Multiple points O3}.
It remains to consider the cases \Em{Multiple points O1/O2} and
\Em{Multiple points O2/O2}. Due to Lemma~\ref{type 1 and 2}
we have to show that in these cases there do not exist
real potentials with the given \Em{complex Fermi curve},
whose eigenfunctions $\psi$ take linearly dependent values
at these four \Em{multiple points}. In case of
\Em{Multiple points O1/O2} this follows from
Corollary~\ref{fixed points}.
If in case \Em{Multiple points O2/O2} the values of the eigenfunction
$\psi$ at the four points are linearly dependent, then,
due to Corollary~\ref{fixed points} and Lemma~\ref{projection 2},
for any non--vanishing normalization of the eigenfunction $\psi$
the functions
$\langle\langle\Op{J}\sigma^{\ast}\psi,\psi\rangle\rangle_{\xx{\gamma}}$
has a zero at the two \Em{Orbits of type~2}.
Since this function is invariant under $\sigma$,
it has to have zeroes of even order at the fixed points of $\sigma$.
Consequently, the corresponding
\De{Local contribution to the arithmetic genus}~\ref{local contribution}
of the \De{Structure sheaf}~\ref{structure sheaf}
has to be larger than one,
and the function $d\yy{p}/d\xx{p}$ has to take the same value
at two points of different \Em{Orbits of type~2}.
\end{proof}

Thus the desired property of the \Em{complex Fermi curve} is the

\newtheorem{Singularity condition}[Lemma]{Singularity condition}
\begin{Singularity condition}\label{singularity condition}
\index{singularity!condition}
\index{condition!singularity $\sim$}
  There exists some
  $\kappa\in\lattice\dual$, such that $[\kappa/2]$ belongs to the
  \Em{complex Fermi curve} $\fermi(U,U)/\lattice\dual$.
  Moreover, the preimage
  of this point in the normalization of $\fermi(U,U)/\lattice\dual$
  contains either two zeroes of both differentials 
  $d\xx{p}$ and $d\yy{p}$, or at least six elements,
  or exactly four elements. Furthermore, in the last case
  the function $d\yy{p}/d\xx{p}$ takes the same value
  at two different elements, which are not interchanged by $\eta$.
\end{Singularity condition}

The discussion above is summarized in

\begin{Theorem} \label{theorem singularity condition}
For a real potential $U\in\banach{2}(\torus)$
we have the following characterizations:
\begin{description}
\item[(i)] The \Em{complex Fermi curve} $\fermi(U,U)$ is a
  \De{Weierstra{\ss} curve}~\ref{Weierstrass curves} if and only if
  the \De{Singularity condition}~\ref{singularity condition} is fulfilled.
\item[(ii)] Let $\fermi(U,U)$ fulfill the
  \De{Singularity condition}~\ref{singularity condition}.
  Then the potential $U$ is a
  \De{Weierstra{\ss} potential}~\ref{Weierstrass potentials}
  if one of the cases~\ref{cusp},\ref{double point},\ref{type 3} or
  \ref{3 type 1 and 2} takes place. Furthermore,
  in the remaining cases~\ref{2 type 1}--\ref{2 type 1 and 2}
  the same is true if and only if the functions
  $d\yy{p}/d\xx{p}$ take the same value at two different elements
  of the preimage of $[\kappa/2]$, which belong to different orbits.\qed
\end{description}
\end{Theorem}

\subsection{Variations of immersions}\label{subsection variations}

Due to the Weierstra{\ss} representation \cite{Ta1,Fr2}
the conformal class of some immersion is given by the
conformal class of the flat torus $\torus$ with the
Euclidean metric $g(\cdot,\cdot)$.
Thus for all variations of some immersion with potential $U$,
which do not change the conformal class,
there exists some variation $\var U$ of the potential $U$,
which is a sufficiently smooth function on $\torus$.
Let us now characterize some variations
$\var U$ of the potential $U$,
which induce a variation of the
immersion and do not change the conformal class.

\begin{Proposition} \label{variation}
Let $U$ be the potential of some immersion and
$\var U\in C(\torus)$ some real variation of this
potential. If for all $\xx{\mu}\in\mathbb{C}^2$ and all regular forms
$\omega$ on the associated \Em{complex Fermi curve}
the pullback of the function
$\Omega_{U,U}(\var U, \var U)/dg(\xx{\mu},k)\cdot
\omega/dg(\xx{\mu},k)$ to the
normalization of $\fermi(U,U)/\lattice\dual$ is holomorphic on some
neighbourhood of the preimage of $[\kappa/2]$ and equal to zero at
all elements of this preimage, then there exists some variation of
the immersion, which corresponds to the variation $\var U$ and
does not change the conformal class.
\end{Proposition}

In order to prove this proposition we need some preparation.

\begin{Lemma} \label{variation of projection}
Let $\xx{\mu}\in\mathbb{C}^2$ satisfy the condition that
$dg(\xx{\mu},k)$ does
not vanish identically on any connected component of the normalization
of $\fermi(V,W)$. Also let $\chi$ be an eigenfunction
corresponding to the element $k'$ of $\fermi(V,W)$.
Assume that for all regular forms $\omega$ of the
\Em{complex Fermi curve} the pullback of the function
$\Omega_{V,W}(\var V,\var W)/dg(\xx{\mu},k)\cdot
\omega/dg(\xx{\mu},k)$ to the
normalization of $\fermi(V,W)$ is holomorphic
on some neighbourhood of the
preimage of $k'$ and equal to zero at all elements of the preimage
of $k'$. Then there exists a variation $\var_{\xx{\mu}}\chi$ of $\chi$ 
associated with the variations $\var V$ and $\var W$
(this means
$$\begin{pmatrix}
\var V & 0\\
0 & \var W
\end{pmatrix}
\begin{pmatrix}
\chi_1\\
\chi_2
\end{pmatrix} +
\begin{pmatrix}
V & \partial_{[k']}\\
\Bar{\partial}_{[k']} & W
\end{pmatrix}
\begin{pmatrix}
\var_{\xx{\mu}}\chi_1\\
\var_{\xx{\mu}}\chi_2
\end{pmatrix} = 0 \text{ and }
\var_{\xx{\mu}}\chi(x+\xx{\gamma})=
\exp\left(2\pi\sqrt{-1}g(\xx{\gamma},k')\right)\var_{\xx{\mu}}\chi(x)$$
for all $x\in\mathbb{R}^2$ and all periods $\xx{\gamma}\in\lattice$),
such that $\langle\langle \phi,
\var_{\xx{\mu}}\chi\rangle\rangle_{\xx{\mu}}$ vanishes for all
eigenfunctions $\phi$ of the transposed Dirac operator, which
correspond to the element $k'$.
\end{Lemma}

\begin{proof} Arguments similar to those in the proof of
  condition~(iii) of Lemma~\ref{projection 2}
  show that for all $\xx{\mu},\yy{\mu}\in\mathbb{C}^2$ on the
  \Em{complex Fermi curve} the identity 
  $$\langle\langle\phi(k),
  \psi(k)\rangle\rangle_{\xx{\mu}} dg(\xx{\mu},k)=
  \langle\langle\phi(k),
  \psi(k)\rangle\rangle_{\yy{\mu}} dg(\yy{\mu},k)$$ holds. Thus the
  first assumption of the lemma implies that the denominator of the
  projection $\triv{\Op{P}}_{\xx{\mu}}$
  does not vanishes identically on the
  \Em{complex Fermi curve}, and therefore this projection
  may be restricted to the
  \Em{complex Fermi curve}.
  Moreover, we choose some $\yy{\mu}\in\mathbb{C}^2$,
  which is linearly independent of $\xx{\mu}$, and parameterize the
  \Em{complex Fermi curve} by $\xx{p}=g(\xx{\mu},k)$
  and $\yy{p}=g(\yy{\mu},k)$. Now let us
  consider the normalizations of the \Em{complex Fermi curves}
  $\fermi(V+t\var V,W+t\var W)$ locally as  covering spaces
  over $(\xx{p},t)\in\mathbb{C}^2$, and let us use the convention
  $\partial \xx{p}/\partial t =0$. We conclude that the local sum of the
  function $\triv{\Op{P}}_{\xx{\mu}}$ over all sheets,
  which contain the element $(k',t=0)$ (compare with Remark~\ref{local sum}),
  is locally a holomorphic function depending on
  $\xx{p}$ and $t$. Hence the partial derivative of this function with
  respect to $t$ at the point $\xx{p}=\xx{p}'=g(\xx{\mu},k'),t=0$
  is some finite rank operator denoted by
  $\var \Breve{\triv{\Op{P}}}_{\xx{\mu}}(k')$. Now we claim
  that $\var_{\xx{\mu}}\chi=\psi_{k'}\var
  \Breve{\triv{\Op{P}}}_{\xx{\mu}}(k')\psi_{-k'}\chi$
  (or equivalently $\var_{\xx{\mu}}\triv{\chi}=\var
  \Breve{\triv{\Op{P}}}_{\xx{\mu}}(k')\triv{\chi}$)
  has the required properties. Since the local sum of the function
  $\triv{\Op{P}}_{\xx{\mu}}$ is a projection valued function,
  we have the identity
  $\var \Breve{\triv{\Op{P}}}_{\xx{\mu}}(k')\comp
  \Breve{\triv{\Op{P}}}_{\xx{\mu}}(k') +
  \Breve{\triv{\Op{P}}}_{\xx{\mu}}(k')\comp\var
  \Breve{\triv{\Op{P}}}_{\xx{\mu}}(k') =
  \var \Breve{\triv{\Op{P}}}_{\xx{\mu}}(k')$, and therefore also
  $\Breve{\triv{\Op{P}}}_{\xx{\mu}}(k')\comp\var
  \Breve{\triv{\Op{P}}}_{\xx{\mu}}(k')
  \comp\Breve{\triv{\Op{P}}}_{\xx{\mu}}(k')=0$.
  This implies the last property.
  All projections $\triv{\Op{P}}_{\xx{\mu}}$
  project onto double periodic functions, which implies
  \begin{align*}
  \var_{\xx{\mu}}\chi(x+\xx{\gamma})&=
  \exp\left(2\pi\sqrt{-1}g(\xx{\gamma},k')\right)\var_{\xx{\mu}}\chi(x)
  &\text{for all }x\in\mathbb{R}^2\text{ and all }\xx{\gamma}\in\lattice.&
  \end{align*}
  Due to our above discussion of Lemma~\ref{regular form},
  the function $\Omega_{V,W}(\var V,\var W)/d\xx{p}$
  is then proportional to $\partial\yy{p}/\partial t$.
  Therefore, the assumption on the form
  $\Omega_{V,W}(\var V, \var W)$ implies that the pullback of
  $\partial\yy{p}/\partial t\cdot\triv{\Op{P}}_{\xx{\mu}}$
  to the normalization of $\fermi(V,W)$
  is holomorphic on some neighbourhood of the preimage of $k'$
  and equal to zero at all elements of the preimage of $k'$. 
  This implies the following operator relation:
  \begin{multline*} 
  \left(\sum\limits_{i=1}^{l}
  \left(\triv{\Op{D}}_{\xx{\mu},\yy{\mu}}(\xx{p}')
  -\yy{p}'\pi\unity\right)^{i-1}\comp
  \var \triv{\Op{D}}_{\xx{\mu},\yy{\mu}}\comp 
  \left(\triv{\Op{D}}_{\xx{\mu},\yy{\mu}}(\xx{p}')
  -\yy{p}'\pi\unity\right)^{l-i}\right)\comp
  \Breve{\triv{\Op{P}}}_{\xx{\mu}}(k') +\\
  \left(\triv{\Op{D}}_{\xx{\mu},\yy{\mu}}(\xx{p})
  -\yy{p}'\pi\unity\right)^{l}\comp
  \var\Breve{\triv{\Op{P}}}_{\xx{\mu}}(k')=0.
  \end{multline*}
  Here $l$ is the dimension of the range of
  $\Breve{\triv{\Op{P}}}_{\xx{\mu}}(k')$, 
  $\var\triv{\Op{D}}_{\xx{\mu},\yy{\mu}}=
  \begin{pmatrix}
  0 & \frac{\var W}{\yy{\nu}_2-\sqrt{-1}\yy{\nu}_1}\\
  \frac{\var V}{\yy{\nu}_2+\sqrt{-1}\yy{\nu}_1} & 0
  \end{pmatrix}$ denotes the variation of
  the operator $\triv{\Op{D}}_{\xx{\mu},\yy{\mu}}(\xx{p})$,
  and $\xx{p}'=g(\xx{\mu},k')$ and $\yy{p}'=g(\yy{\mu},k')$
  are the components of $k'$.
  Since $\psi_{-k'}\chi$ is a proper eigenfunction of 
  $\triv{\Op{D}}_{\xx{\mu},\yy{\mu}}(\xx{p})$
  with the eigenvalues $\yy{p}'\pi$, the
  variation $\psi_{-k'}\var_{\xx{\mu}}\chi$ satisfies
  $$\left(\triv{\Op{D}}_{\xx{\mu},\yy{\mu}}(\xx{p})
  -\yy{p}'\pi\unity\right)^{l-1} \comp
  \var\triv{\Op{D}}_{\xx{\mu},\yy{\mu}}\psi_{-k'}\chi +
  \left(\triv{\Op{D}}_{\xx{\mu},\yy{\mu}}(\xx{p})
  -\yy{p}'\pi\unity\right)^{l}
  \psi_{-k'}\var_{\xx{\mu}}\chi=0.$$
  Since $\Breve{\triv{\Op{P}}}_{\xx{\mu}}(k')$
  is the spectral projection onto
  the generalized eigenspace of
  $\triv{\Op{D}}_{\xx{\mu},\yy{\mu}}(\xx{p})$
  associated with eigenvalues $\yy{p}'\pi$,
  the restriction of the
  operator $\triv{\Op{D}}_{\xx{\mu},\yy{\mu}}(\xx{p})
  -\yy{p}'\pi\unity$ to the kernel of the
  operator $\Breve{\triv{\Op{P}}}_{\xx{\mu}}(k')$ is invertible.
  Thus we have
  $\var\triv{\Op{D}}_{\xx{\mu},\yy{\mu}}\psi_{-k'}\chi +
  \left(\triv{\Op{D}}_{\xx{\mu},\yy{\mu}}(\xx{p})
  -\yy{p}'\pi\unity\right)
  \psi_{-k'}\var_{\xx{\mu}}\chi=0$,
  which is equivalent to
  $\begin{pmatrix}
  \var V & 0\\
  0 & \var W
  \end{pmatrix}
  \begin{pmatrix}
  \chi_1\\
  \chi_2
  \end{pmatrix} +
  \begin{pmatrix}
  V & \partial_{[k']} \\
  -\Bar{\partial}_{[k']} & W
  \end{pmatrix}
  \begin{pmatrix}
  \var_{\xx{\mu}}\chi_1\\
  \var_{\xx{\mu}}\chi_2
  \end{pmatrix} = 0$.
\end{proof}

\noindent
{\it Proof of Proposition~\ref{variation}.}
Let $E([\kappa/2])$ be the
eigenspace of the Dirac operator $\Op{D}_{[\kappa/2]}$
with real potential
$U$ corresponding to the eigenvalue
$\lambda=0$. Due to Lemma~\ref{involutions} this operator is
self--adjoint, and $E([\kappa/2])$ is equal to the range of
$\Breve{\Op{P}}([\kappa/2)$.
Moreover, the anti--unitary operator $\chi\mapsto
\Op{J}\Bar{\chi}$ leaves this subspace invariant. Due to the
\De{Equivalent form of the Periodicity condition}~\ref{equivalent form}
all non--zero elements $\chi$ of $E([\kappa/2])$,
which satisfy the
following two conditions, induce an immersion:
\begin{description}
\item[(A)] for all $\xx{\mu}\in\mathbb{C}^2$ the bilinear form
  $\langle\langle \Op{J}\chi,
  \chi\rangle\rangle_{\xx{\mu}}$
  is zero.
\item[(B)] for all $\xx{\mu}\in\mathbb{C}^2$ the bilinear form
  $\langle\langle\Bar{\chi},
  \chi\rangle\rangle_{\xx{\mu}}$
  is zero.
\end{description}
If $\var U$ is some variation of the potential $U$,
which satisfies
the assumption of Proposition~\ref{variation}, then due to
Lemma~\ref{variation of projection} there exists some variation
$\var\chi$ of $\chi$ associated with the variation
$\var U$ of the potential $U$:
\begin{align*}
\begin{pmatrix}
\var U & 0\\
0 & \var U
\end{pmatrix}
\begin{pmatrix}
\chi_1\\
\chi_2
\end{pmatrix} +
\begin{pmatrix}
U & \partial_{[\kappa/2]} \\
-\Bar{\partial}_{[\kappa/2]} & U
\end{pmatrix}
\begin{pmatrix}
\var\chi_1\\
\var\chi_2
\end{pmatrix}&=0&\text{and }
\var\chi(x+\xx{\gamma})&=(-1)^{g(\xx{\gamma},\kappa)}\chi(x)
\end{align*}
for all $x\in\mathbb{R}^2$ and all periods $\xx{\gamma}\in\lattice$.
Obviously we may add to such a variations $\var\chi$ some arbitrary
element of $E([\kappa/2])$ and again obtain such a variation. 
In a second step we show that there exists some $\var\chi$ among
these variations, which satisfies
\begin{description}
\item[(A')] for all $\xx{\mu}\in\mathbb{C}^2$ the bilinear form
  $\langle\langle \Op{J}\chi,
  \var\chi\rangle\rangle_{\xx{\mu}}$
  is zero.
\item[(B')] for all $\xx{\mu}\in\mathbb{C}^2$ the bilinear form
  $\langle\langle\Bar{\chi},
  \var\chi\rangle\rangle_{\xx{\mu}}$
  is zero.
\end{description}
Again, due to Lemma~\ref{variation of projection},
for all $\xx{\mu}\in\mathbb{C}^2$ 
there exists such a variation $\var_{\xx{\mu}}\chi$,
such that $\langle\langle \Op{J}\xi,
\var_{\xx{\mu}}\chi\rangle\rangle_{\xx{\mu}}$ vanishes for all $\xi\in
E([\kappa/2])$. In particular, $\langle\langle \Op{J}\chi,
\var_{\xx{\mu}}\chi\rangle\rangle_{\xx{\mu}}$ and
$\langle\langle \Bar{\chi},
\var_{\xx{\mu}}\chi\rangle\rangle_{\xx{\mu}}$ vanish.
For linearly independent $\xx{\mu}$ and $\yy{\mu}$ the difference
$\var_{\xx{\mu}}\chi-\var_{\yy{\mu}}\chi$ can be non--zero,
but has to be an element of $E([\kappa/2])$.
Now we again choose two generators $\xx{\gamma}$ and
$\yy{\gamma}$ of the period lattice. Let
$\langle\cdot,\cdot\rangle_1$ denote the restriction of
$\langle\langle \Op{J}\cdot,
\cdot\rangle\rangle_{\xx{\gamma}}$ and
$\langle\cdot,\cdot\rangle_2$ the restriction of
$\langle\langle \Op{J}\cdot,
\cdot\rangle\rangle_{\yy{\gamma}}$ to
$E([\kappa/2])\times E([\kappa/2])$.
Both bilinear forms are symmetric and satisfy
\begin{align*}
\langle\Op{J}\Bar{\chi},\Op{J}\Bar{\xi}\rangle_1&=
-\overline{\langle\chi,\xi\rangle}_1&\text{and }
\langle \Op{J}\Bar{\chi},\Op{J}\Bar{\xi}\rangle_2&=
-\overline{\langle\chi,\xi\rangle}_2.
\end{align*}

\begin{Lemma} \label{bilinear forms}
One of the following two cases hold:
\begin{description}
\item[(a)] The four linear forms 
  $\langle\chi,\cdot\rangle_1$,
  $\langle\chi,\cdot\rangle_2$,
  $\langle \Op{J}\Bar{\chi},\cdot\rangle_1$ and
  $\langle \Op{J}\Bar{\chi},\cdot\rangle_2$ on $E([\kappa/2])$ are
  linearly independent.
\item[(b)] There exists a non--zero linear combination $\chi'$ of
  $\chi$ and $\Op{J}\Bar{\chi}$, such that the linear forms
  $\langle\chi',\cdot\rangle_1$ and
  $\langle\chi',\cdot\rangle_2$ on $E([\kappa/2])$ are
  linearly dependent.
\end{description}
\end{Lemma}

\begin{proof} If the form
  $\langle\chi,\cdot\rangle_1$ is zero, then
  \Em{case~(b)} takes place with $\chi'=\chi$.
  Now let us assume that this linear form does not vanish.
  The subspace of $E([\kappa/2])$,
  on which the two linear
  forms $\langle\chi,\cdot\rangle_1$ and 
  $\langle \Op{J}\Bar{\chi},\cdot\rangle_1$ vanish, contains
  the subspace, on which all four linear forms vanish. Furthermore,
  both subspaces are invariant under the anti--unitary operator
  $\xi\mapsto \Op{J}\Bar{\xi}$. The same arguments as in the proof of
  Corollary~\ref{fixed points}~(iii) show that
  either these two subspaces coincide,
  or the co--dimension of the second subspace in
  $E([\kappa/2])$ is equal to four, which implies that \Em{case~(a)}
  takes place. If the two subspaces of $E([\kappa/2])$ coincide, then
  there exist two complex numbers $\alpha$ and $\beta$, such that the
  linear forms fulfill the following equations:
  \begin{align*}
  \langle\chi,\cdot\rangle_2&=
  \alpha\langle\chi,\cdot\rangle_1+
  \beta\langle\Op{J}\Bar{\chi},\cdot\rangle_1&\text{and }
  \langle \Op{J}\Bar{\chi},\cdot\rangle_2&=
  -\Bar{\beta}\langle\chi,\cdot\rangle_1+
  \Bar{\alpha}\langle \Op{J}\Bar{\chi},\cdot\rangle_1.
  \end{align*}
  If $\left(\begin{smallmatrix}
  \gamma \\
  \var
  \end{smallmatrix}\right)$ is any non--zero eigenvector of the matrix
  $\left(\begin{smallmatrix}
  \alpha & -\Bar{\beta}\\
  \beta & \Bar{\alpha}
  \end{smallmatrix}\right)$, then $\chi'=\gamma\chi+\var
  \Op{J}\Bar{\chi}$ fulfills the condition of \Em{case~(b)}.
\end{proof}

\noindent
{\it Continuation of the proof of Proposition~\ref{variation}.} In
\Em{case~(a)} we may add to some $\var_{\xx{\mu}}\chi$ some element of
$E([\kappa/2])$, such that the sum fulfills \Em{conditions~(A')} and
\Em{(B')}. In \Em{case~(b)} we may assume that the linear forms
$\alpha\langle\chi',\cdot\rangle_1+
\langle\chi',\cdot\rangle_2$ and 
$\Bar{\alpha}\langle \Op{J}\Bar{\chi}',\cdot\rangle_1+
\langle \Op{J}\Bar{\chi}',\cdot\rangle_2$ vanish.
Indeed, if the form $\langle\chi',\cdot\rangle_1$ vanishes,
then we replace $\xx{\gamma}$ and $\yy{\gamma}$
(i.\ e.\ $\langle\chi',\cdot\rangle_1$ and
.$\langle\chi',\cdot\rangle_2$)
Due to Lemma~\ref{variation of projection}
the variations $\var_{\xx{\gamma}}\chi'$ obey the equations
$$\langle\chi',\var_{\xx{\gamma}}\chi'\rangle_1=0=
\langle \Op{J}\Bar{\chi'},\var_{\xx{\gamma}}\chi'\rangle_1.$$
If $\xx{\mu}$ is equal to $\alpha\xx{\gamma}+\yy{\gamma}$
and $\Bar{\xx{\mu}}$
equal to $\Bar{\alpha}\xx{\gamma}+\yy{\gamma}$, then the variations
$\var_{\xx{\mu}}\chi'$ and $\var_{\Bar{\xx{\mu}}}{\chi'}$ satisfy
$$\langle\langle \Op{J}\chi',
\var_{\xx{\mu}}\chi'\rangle\rangle_{\xx{\mu}} =0=
\langle\langle \Bar{\chi}',
\var_{\Bar{\xx{\mu}}}\chi'\rangle\rangle_{\Bar{\xx{\mu}}}.$$
Since the differences
$\var_{\xx{\gamma}}\chi'-\var_{\xx{\mu}}\chi'$ and 
$\var_{\xx{\gamma}}\chi'-\var_{\Bar{\xx{\mu}}}\chi'$ belong to
$E([\kappa/2])$, the variation $\var_{\xx{\gamma}}\chi'$ satisfies
\begin{align*}
\langle\langle \Op{J}\chi',
\var_{\xx{\gamma}}\chi'\rangle\rangle_{\xx{\mu}} &=&
\langle\langle \Op{J}\chi',
\var_{\xx{\gamma}}\chi'-\var_{\xx{\mu}}\chi'
\rangle\rangle_{\xx{\mu}} &=&
\alpha\langle\chi',
\var_{\xx{\gamma}}\chi'-\var_{\xx{\mu}}\chi'\rangle_1+
\langle\chi',
\var_{\xx{\gamma}}\chi'-\var_{\xx{\mu}}\chi'\rangle_2 &=0\\
-\langle\langle \Bar{\chi}',
\var_{\xx{\gamma}}\chi'\rangle\rangle_{\Bar{\xx{\mu}}} &=&
-\langle\langle \Bar{\chi}',
\var_{\xx{\gamma}}\chi'-\var_{\Bar{\xx{\mu}}}\chi'
\rangle\rangle_{\Bar{\xx{\mu}}} &=&
\Bar{\alpha}\langle \Op{J}\Bar{\chi'},
\var_{\xx{\gamma}}\chi'-\var_{\Bar{\xx{\mu}}}\chi'\rangle_1+
\langle \Op{J}\Bar{\chi'},
\var_{\xx{\gamma}}\chi'-\var_{\Bar{\xx{\mu}}}\chi'\rangle_2 &=0,
\end{align*}
and therefore also the relations
\begin{align*}
\langle\chi',\var_{\xx{\gamma}}\chi'\rangle_2 &=&
\alpha\langle\chi',\var_{\xx{\gamma}}\chi'\rangle_1+
\langle\chi',\var_{\xx{\gamma}}\chi'\rangle_2 &=&
\langle\langle \Op{J}\chi',
\var_{\xx{\gamma}}\chi'\rangle\rangle_{\xx{\mu}} &=0\\
\langle \Op{J}\Bar{\chi'},\var_{\xx{\gamma}}\chi'\rangle_2 &=&
\Bar{\alpha}\langle \Op{J}\Bar{\chi'},
\var_{\xx{\gamma}}\chi'\rangle_1+
\langle \Op{J}\Bar{\chi'},\var_{\xx{\gamma}}\chi'\rangle_2 &=&
-\langle\langle \Bar{\chi}',
\var_{\xx{\gamma}}\chi'\rangle\rangle_{\Bar{\xx{\mu}}} &=0.
\end{align*}
These arguments carry over to $\Op{J}\Bar{\chi}'$ and show that
there exists a variation $\var_{\xx{\gamma}}\Op{J}\Bar{\chi}'$
such that the relation
\begin{align*}
\langle\chi',\var_{\xx{\gamma}}\Op{J}\Bar{\chi}'\rangle_1 &=0 &
\langle\chi',\var_{\xx{\gamma}}\Op{J}\Bar{\chi}'\rangle_2 &=0\\
\langle\Op{J}\Bar{\chi}',
\var_{\xx{\gamma}}\Op{J}\Bar{\chi}'\rangle_1 &=0 &
\langle\Op{J}\Bar{\chi}',
\var_{\xx{\gamma}}\Op{J}\Bar{\chi}'\rangle_2 &=0
\end{align*}
hold. Since $\chi$ is a linear combination of $\chi'$ and $\Op{J}\Bar{\chi}'$
there exists a variation $\var\chi$ with the desired properties.
\qed

\subsection{Constrained Willmore tori}\label{subsection constrained}

An immersion, which is a stationary point of the Willmore functional
is called a Willmore torus. The following terminology can be found for
example in \cite[7.10]{Wi2}:

\newtheorem{Constrained Willmore tori}[Lemma]{Constrained Willmore tori}
\begin{Constrained Willmore tori}\label{constrained Willmore tori}
\index{Willmore!torus!constrained $\sim$}
Variations of immersions, which do not change the conformal class,
are called \Em{conformal} variations.
An immersion, which is a stationary point of the Willmore functional
with respect to all \Em{conformal} variations of the immersion
is called a \De{Constrained Willmore torus}.
\end{Constrained Willmore tori}

Obviously all Willmore tori are \De{Constrained Willmore tori}.
In this section we want to prove the following

\begin{Theorem} \label{constrained function}
Let $U$ be the potential of a
\De{Constrained Willmore torus}~\ref{constrained Willmore tori}.
Then the normalization of $\fermi(U,U)/\lattice\dual$
has a meromorphic function $f$ with finitely many poles
at the preimage of the element
$[\kappa/2]$ described in the
\De{Singularity condition}~\ref{singularity condition}. Moreover,
for arbitrary small $\delta'>0$ there exists some $\delta,\varepsilon>0$,
such that $f$ maps the open set $\Set{V}^+_{\varepsilon,\delta}$ into
$\{z\in \mathbb{C}\mid|z+1|<\delta'\}$ and the set
$\Set{V}^-_{\varepsilon,\delta}$  into
$\{z\in \mathbb{C}\mid|z-1|<\delta'\}$. Finally, the function $f$
is invariant under $\sigma$ and equal to $-\rho^{\ast}\Bar{f}$.
\end{Theorem}

\begin{proof} Let $U$ be the potential of a
\De{Constrained Willmore torus}~\ref{constrained Willmore tori}.
Due to Proposition~\ref{variation} there exist finitely
many meromorphic functions $f_1,\ldots,f_l$ on some neighbourhood
$\Set{U}$ of the element $[\kappa/2]$ of the \Em{complex Fermi curve}
$\fermi(U,U)/\lattice\dual$ such that all real variations $\var U$,
whose form $\Omega_{U,U}(\var U,\var U)$ multiplied with any of these
functions has zero residue on $\Set{U}$, induce some variation of the
immersion. The pullbacks of all these functions to the normalization of
$\fermi(U,U)/\lattice\dual$ are holomorphic in the complement of the
preimage of $[\kappa/2]$. Moreover, the neighbourhood
$\Set{U}$ may be chosen to be invariant under $\sigma$ and $\rho$,
and the functions $f_1,\ldots,f_l$ can be chosen to be invariant under
$\sigma$ and to satisfy $-\rho^{\ast}\Bar{f}_i=f_i$, $i=1,\ldots,l$. In
fact, since $\Omega_{U,U}(\var U,\var U)$ is invariant under $\sigma$
and satisfies
$\rho^{\ast}\left(\Bar{\Omega}_{U,U}(\var U,\var U)\right)=
\Omega_{U,U}(\var U,\var U)$, for all meromorphic functions $f$ on
$\Set{U}$ the total residue of $f\cdot\Omega_{U,U}(\var U,\var U)$ is
equal to the residue of
$\sigma^{\ast}\left(f\cdot\Omega_{U,U}(\var U,\var U)\right)=
\sigma^{\ast}\left(f\right)\cdot\Omega_{U,U}(\var U,\var U)$ and equal to
the complex conjugate of the residue of
$\rho^{\ast}\left(\Bar{f}\cdot\Bar{\Omega}_{U,U}(\var U,\var U)\right)=
\rho^{\ast}\left(\Bar{f}\right)\cdot\Omega_{U,U}(\var U,\var U)$.
Now let
$\triv{\Op{A}}_1(\cdot),\ldots,\triv{\Op{A}}_l(\cdot)$
denote the
corresponding operator--valued meromorphic functions defined in
Lemma~\ref{residue} and $\triv{\Op{A}}_{0}(\cdot)$
the constant function
defined in (iv) of Lemma~\ref{compatibility}. Due to
Lemma~\ref{compatibility} the corresponding variations are of the form
$(-\var U_0,\var U_0),\ldots,(-\var U_l,\var U_l)$,
with real variations $\var U_0,\ldots,\var U_l\in
\banach{2}_{\mathbb{R}}(\torus)$
and with $\var U_0$ being equal to $2U$.
The right hand side of the formula~(iv) in Lemma~\ref{residue}
for the variations
$(\var U,\var U)$ and $(\var U_0,\var U_0)$ is equal to
$$\frac{1}{\pi\sqrt{-1}(\kappa_1\kappa_2'-\kappa_2\kappa_1')}
  \int\limits_{\torus}
  \var U\var U_0 d^2x = 
\frac{2}{\pi\sqrt{-1}(\kappa_1\kappa_2'-\kappa_2\kappa_1')}
  \int\limits_{\torus}
  \var U U d^2x,$$
and therefore proportional to the variation of the
Willmore functional. Since $U$ is the potential of a
\De{Constrained Willmore torus}~\ref{constrained Willmore tori},
this integral must vanish, whenever all the integrals
$$\frac{1}{\pi\sqrt{-1}(\kappa_1\kappa_2'-\kappa_2\kappa_1')}
  \int\limits_{\torus}
  \var U\var U_i d^2x,\; i=1,\ldots,l$$
vanish. We conclude that there exist some real constants
$c_1,\ldots,c_l$, such that
$\var U_0+c_1\var U_1+\ldots+\var U_l$ is identically
zero. Now Lemma~\ref{residue} implies that the meromorphic function
$$\triv{\Op{A}}(\cdot)=
\triv{\Op{A}}_0(\cdot)+
c_1\triv{\Op{A}}_1(\cdot)+\ldots+\triv{\Op{A}}_l(\cdot)$$
commutes with the operator--valued function 
$\triv{\Op{D}}_{\xx{\gamma},\yy{\gamma}}(\cdot)$.
Thus these two operators can be diagonalized
simultaneously, and the eigenvalue of $\triv{\Op{A}}(\cdot)$
defines some global meromorphic function $f$
on the \Em{complex Fermi curve}:
$$f = \tr\left(\triv{\Op{P}}_{\xx{\gamma}}\comp
  \triv{\Op{A}}(p)\right),$$
where $p$ is equal to the corresponding coordinate of the variable
of $\triv{\Op{P}}_{\xx{\gamma}}$.
Property (i) of Lemma~\ref{residue} ensures
that this function is invariant under the action of the dual lattice
and therefore defines a meromorphic function on
$\fermi(U,U)/\lattice\dual$.
The estimate (iii) of Lemma~\ref{residue}
and Theorem~\ref{asymptotic analysis 1}
implies the required estimate on $f$.
Due to the definition of $f$ the pullback of the difference
$f-c_1f_1-\ldots-c_lf_l$ to the normalization of
$\fermi(U,U)/\lattice\dual$ is holomorphic on the complement of the
preimage of $[\kappa/2]$. The other required properties of $f$
follows from Lemma~\ref{compatibility}.
\end{proof}

\begin{Corollary} \label{finite type}
\index{Willmore!torus!constrained $\sim$}
All \De{Constrained Willmore tori}~\ref{constrained Willmore tori}
are of finite type.
\end{Corollary}

\begin{proof} For a given $0<\delta'<1/2$ there exist positive
  $\delta,\varepsilon>0$, such that the meromorphic function $f^2-1$
  is bounded on
  $\Set{V}^+_{\varepsilon,\delta}\cup\Set{V}^-_{\varepsilon,\delta}$
  by $\delta'<1/2$. Thus, due to the Maximum principle (see e.g.
  \cite[Chapter~4. Theorem~12.]{Ah}), this bound holds also on all but
  finitely many handles. But the set $\{z\in\mathbb{C}\mid
  |z^2-1|<1/2\}$ decomposes into two connected components. Therefore,
  the normalizations of all but finitely many handles have two
  connected components. This proves the corollary.
\end{proof}

This corollary is similar to the finite type theorem of
\cite{PS}. The main difference is that the authors consider from the
very beginning immersions with special properties, and which do not have to
be double periodic. These properties result in a global meromorphic
function on the corresponding Riemann surface.
Afterwards they study double periodic solutions, and obtain 
some restrictions on the Riemann surface. In particular, the Riemann
surface has to be of finite genus. In our approach we consider
from the very beginning only double periodic immersions. Afterwards we
investigate the consequences of the property, that the immersion is
a \De{Constrained Willmore torus} and again obtain the proof of the
existence of some global meromorphic function on the 
Riemann surface. Now, due to the
special form of Riemann surfaces corresponding to double periodic
solutions, we can also conclude that the Riemann surface has
finite genus. The advantage of our approach in the present case is
that the global meromorphic function can be of different degree and
therefore will result in different properties of the immersions. In
fact, our construction shows that the number of poles is related to
the dimension of the eigenspace corresponding to $[\kappa/2]$. The
results of \cite{GS1} show that even in the case of surfaces of
revolutions this dimension can be arbitrary large. But
Theorem~\ref{theorem singularity condition} suggests that for
relative minimizers the degree of $f$ should be equal to four. Moreover,
this global meromorphic function should map
the complement of that subset of the Riemann surface,
which consists of the two points corresponding to
infinite energy, into the complement of some finite subset of
$\mathbb{P}^1$. Only the existence of such a global
meromorphic function will imply that the potentials satisfy some
partial differential equation. Otherwise it would be quite
complicated to formulate the property,
that the \Em{complex Fermi curve} has
a global meromorphic function in terms of the potentials.
This property of the global meromorphic function
is related to the existence of umbilic points,
which are excluded in the description of
Willmore tori with the help of the sin-Gordon equation (see
\cite{Bo,FPPS,BB}).

It can happen that not only the \Em{complex Fermi curve}, but the whole
complex Bloch variety of some Dirac operator can be compactified to some
projective variety. Therefore, the question arises, whether the whole
complex Bloch varieties of Willmore tori can be compactified to projective
varieties. Here we should mention, that the \Em{complex Fermi curve} alone
does not determine the whole spectrum. Moreover, the self--adjoint
Dirac operator
$\left(\begin{smallmatrix}
U & \partial\\
-\Bar{\partial} & U
\end{smallmatrix} \right)$
has the same \Em{complex Fermi curve} as the similar operator
$\left(\begin{smallmatrix}
\partial & -U\\
U &\Bar{\partial}
\end{smallmatrix} \right)=
\left(\begin{smallmatrix}
U & \partial\\
-\Bar{\partial} & U
\end{smallmatrix} \right)\comp\Op{J}$, but the associated
complex Bloch varieties are completely different.
From the very beginning we could have used the second operator.
Indeed, the operators
$\triv{\Op{D}}_{\xx{\gamma},\yy{\gamma}}(p)$,
introduced in Section~\ref{subsection variations}
are more closely related to the last operator
than with the original Dirac operator. The second
operator has a remarkable property: the \Em{complex Fermi curves}
corresponding to different energies are all isomorphic.
In fact, if we multiply an eigenfunction with the function
$\psi_{k}$, with $k=(\lambda,0)$, the energy is shifted by
$\lambda\pi\sqrt{-1}$. Thus the whole complex Bloch variety does not
contain more information than the \Em{complex Fermi curve} and may be
compactified to a projective variety if and only if the
\Em{complex Fermi curve} can be compactified to an algebraic curve.
This is not the case for the Dirac operator.
Also the property, that the whole
complex Bloch variety $\bloch(U,U)/\lattice\dual$ can be
compactified to a projective variety, is much more restrictive.

Obviously, the\De{Periodicity condition}~\ref{periodicity condition}
makes sense for spinors in the kernel of a
\De{Finite rank Perturbation}~\ref{finite rank perturbations}.
The corresponding conformal mappings from $\torus$ into $\mathbb{R}^3$
can have finitely many poles. Moreover, the whole discussion of
Section~\ref{subsection singularity} carries over to these conformal
mappings with poles. We remark that due to the relation of the
\De{Finite rank Perturbation}~\ref{finite rank perturbations} with the
\De{Singularities of the holomorphic structure}~\ref{singularities of
  the structure}
the closed forms, which were used to define the immersions, have no
monodromy around the singular points. However, suitable inversions
transform these conformal mappings with poles into
conformal mappings without poles from $\torus$ into $\mathbb{R}^3$.
Due to the conformal invariance of the \Em{complex Fermi curves}
\cite{GS2} these inversions do not change the
\Em{complex Fermi curves}. Therefore we may extend the
Willmore functional to those \Em{complex Fermi curves} in
$\Bar{\moduli}_{\lattice\eta,\sigma}$, which obey the
\De{Singularity condition}~\ref{singularity condition}.
These \Em{complex Fermi curves}, which obey the
\De{Singularity condition}~\ref{singularity condition} are obviously
locally zeroes of holomorphic functions on the compactified moduli space
$\Bar{\moduli}_{\lattice\eta,\sigma}$.
Hence we obtain

\begin{Corollary}\label{existence of constrained minimizers}
The restrictions of the Willmore functional to an
arbitrary conformal class has a minimizer of \Em{finite type}.
In particular, the corresponding immersion is analytic. 
\qed
\end{Corollary}

Finally, we remark that this result carries over to immersion into
$\mathbb{R}^4$.
More precisely, for all $k\in\mathbb{R}^2$ the restrictions of the
\Em{first integral} to all \Em{complex Fermi curves} in
$\Bar{\moduli}_{\lattice,\eta}$,
which obey the
\De{Quaternionic Singularity condition}~\ref{quaternionic condition}
over an element, whose imaginary part is equal to $k$,
has a minimizer. Due to the results of `quaternionic function theory'
developed by F.\ Pedit and U.\ Pinkall \cite{PP,FLPP},
these functionals are the restrictions of the Willmore functional in
$\mathbb{R}^4$ to fixed conformal classes and
fixed holomorphic line bundles
(underlying the quaternionic line bundle).
We expect that similar arguments may show that these minimizers are of
\Em{finite type}.
Since the space of holomorphic bundles is compact, the restrictions of
the Willmore functional in $\mathbb{R}^4$ to arbitrary conformal
classes has also a minimizer.
Since the whole family constructed in
Section~\ref{subsection absolute minimizer} obey the
\De{Quaternionic Singularity condition}~\ref{quaternionic condition},
for non--rectangular conformal classes the corresponding minimizers
in $\mathbb{R}^4$ take lower values than the corresponding minimizers
in $\mathbb{R}^3$

\begin{theindex}

  \item accent tilde, 16
  \item admissible form, 160
  \item arithmetic genus, 39
    \subitem local contribution to the $\sim$, 156
    \subitem of the $\rightarrow$ structure sheaf, 39
  \item asymptotic analysis, 27, 139

  \indexspace

  \item B\"acklund transformation, 47--50, 196
  \item Baker--Akhiezer function ($\rightarrow$ eigenfunction), 43
    \subitem condition   $\sim$ (i)--(ii), 52--53
    \subitem divisor of the   $\sim$, 45
  \item Bloch
    \subitem complex $\sim$ variety $\bloch$, 25
    \subitem variety, 25
  \item blowing up
    \subitem of the moduli space, 74

  \indexspace

  \item Carleman inequality, 23
  \item Clifford torus, 3
  \item compactification
    \subitem of an isospectral set, 51
    \subitem of the moduli space, 80, 139
  \item condition
    \subitem Baker--Akhiezer function (i)--(ii), 52--53
    \subitem decomposition (i)--(ii), 130--132
    \subitem divisor (i)--(iii), 45
    \subitem divisor (iv), 57
    \subitem finite rank perturbation (i)--(iii), 94
    \subitem form (i)--(iv), 80--81
    \subitem neighbourhood (i)--(ii), 81--82
    \subitem periodicity $\sim$, 5
    \subitem quasi--momenta  (i)--(iv), 43--44
    \subitem quasi--momenta  (v), 68
    \subitem quasi--momenta (i') and (vi), 144
    \subitem quasi--momenta (ii'), 80
    \subitem quaternionic singularity $\sim$, 7
    \subitem singularity $\sim$, 7, 197
    \subitem singularity (i)--(ii), 102
    \subitem weak singularity $\sim$, 7, 155
  \item conformal
    \subitem class
      \subsubitem moduli space $\moduli_1$ of $\sim$es, 3, 185
      \subsubitem of an immersion, 3
      \subsubitem of the lattice $\lattice$, 185
    \subitem immersion, 4
    \subitem mapping, 8, 10
  \item connected components
    \subitem of the normalization, 28
  \item coordinate
    \subitem $\xx{q}$ and $\yy{q}$, 15
    \subitem $x$, $z$, 15
  \item curve
    \subitem complex Fermi $\sim$ $\fermi$, 6, 25
    \subitem Fermi $\sim$, 6
    \subitem generalized Weierstra{\ss} $\sim$, 157
    \subitem quaternionic Weierstra{\ss}   $\sim$, 190
    \subitem spectral $\sim$, 6
    \subitem Weierstra{\ss} $\sim$, 6
  \item cusp, 73, 155, 156, 165
    \subitem of the moduli space, 73
  \item cut
    \subitem horizontal $\sim$, 83
    \subitem modified pair of parallel $\sim$s, 84
    \subitem pair of parallel $\sim$s, 83

  \indexspace

  \item Davey--Stewartson equation, 26
  \item decomposition
    \subitem condition $\sim$ (i)--(ii), 130--132
  \item deformation
    \subitem of a complex Fermi curve, 68
    \subitem of a complex space, 68
  \item Dirac operator $\Op{D}$, 15
    \subitem eigenfunction $\psi$ of the $\sim$, 29
    \subitem generalized $\sim$, 56
    \subitem perturbation of the $\sim$, 99
    \subitem resolvent of the $\sim$ $\Op{R}$, 17
  \item Dirichlet integral, 84
  \item disconnected
    \subitem normalization, 158, 173
  \item dividing
    \subitem real part, 157
  \item divisor
    \subitem condition $\sim$ (i)--(iii), 45
    \subitem condition $\sim$ (iv), 57
    \subitem of the Baker--Akhiezer function, 45

  \indexspace

  \item effective
    \subitem dual lattice        $\lattice\dual_{\text{\scriptsize\rm effective}}$, 
		175
    \subitem lattice $\lattice_{\text{\scriptsize\rm effective}}$, 175
  \item eigenfunction ($\rightarrow$ Baker--Akhiezer function), 29
    \subitem $\psi$   of the Dirac operator, 29
    \subitem  $\phi$   of the transposed Dirac operator, 29
    \subitem normalization   of the $\sim$, 42
  \item existence
    \subitem of a minimizer, 3
    \subitem of a potential, 44

  \indexspace

  \item family
    \subitem $\fermi_{\min,[\kappa/2]}(\cdot)$, 178
    \subitem $\fermi_{\min}(\cdot)$, 170
  \item Fermi curve, 6
    \subitem complex $\sim$ $\fermi$, 6, 25
  \item finite
    \subitem arithmetic genus, 42
    \subitem geometric genus, 42
    \subitem rank perturbation, 99
      \subsubitem condition $\sim$ (i)--(iii), 94
    \subitem type potentials, 42
  \item form
    \subitem admissible $\sim$, 160
    \subitem condition $\sim$ (i)--(iv), 80--81
    \subitem regular $\sim$, 34, 36, 58
  \item Fourier transform $\Op{F}$, 21
  \item fundamental domain, 15

  \indexspace

  \item generalized
    \subitem Cauchy's integral formula, 100
    \subitem Dirac operator, 56
    \subitem H\"older's inequality, 120
    \subitem Weierstra{\ss}
      \subsubitem curve, 9, 157
      \subsubitem potential, 8, 9
    \subitem Willmore functional, 157
    \subitem Young's inequality, 121
  \item genus
    \subitem $\rightarrow$ arithmetic $\sim$, 39
    \subitem $\rightarrow$ geometric $\sim$, 39
  \item geometric genus, 39
    \subitem finite $\sim$, 42
  \item gluing rule, 84
  \item Green's function
    \subitem $\sim$ $\Func{G}_{\lambda}$, 17
    \subitem $\sim$ $\Func{G}_{\mathbb{P}^1,\lambda}$, 109

  \indexspace

  \item handle, 27
  \item Hausdorff metric, 78
  \item Hilbert bundle, 16

  \indexspace

  \item index set
    \subitem of handles $\lattice\dual_{\delta}$, 27
    \subitem of horizontal cuts $\lattice\dual_{\delta}$, 140
  \item integral
    \subitem Dirichlet $\sim$, 84
    \subitem first $\sim$ $\willmore$, 7, 79, 84
    \subitem kernel
      \subsubitem of $\Op{R}(0,0,k,0)$, 94
      \subsubitem of $\triv{\Op{R}}(0,0,0,\sqrt{-1}\lambda)$, 18
      \subsubitem of a finite rank perturbation, 99
  \item intersection number, 168
    \subitem local contribution to the $\sim$, 168
  \item involution
    \subitem $\eta$, 29
    \subitem $\rho$, 29
    \subitem $\sigma$, 29
    \subitem of the complex Bloch variety, 29
    \subitem of the moduli space, 159
  \item isospectral
    \subitem set, 12
      \subsubitem compactification of an $\sim$, 51
    \subitem transformation, 31, 59
      \subsubitem infinitesimal $\sim$, 59

  \indexspace

  \item lattice
    \subitem $\lattice$, $\lattice_\mathbb{C}$, 15
    \subitem conformal class of the $\sim$ $\lattice$, 185
    \subitem dual $\sim$ $\lattice\dual$, 15
    \subitem dual $\sim$ $\lattice\dual_\mathbb{C}$, 169
    \subitem dual $\sim$ $\lattice\dual_{[\kappa/2]}$, 176
    \subitem dual $\sim$ generator
      \subsubitem $\xx{\kappa}$ and $\yy{\kappa}$, 15
    \subitem dual sub--$\sim$ $\Breve{\lattice}\dual$, 176
    \subitem effective $\sim$ $\lattice_{\text{\scriptsize\rm effective}}$, 
		175
    \subitem effective dual $\sim$        $\lattice\dual_{\text{\scriptsize\rm effective}}$, 
		175
    \subitem generator
      \subsubitem $\xx{\gamma}$ and $\yy{\gamma}$, 15
  \item Lax operator, 6
  \item local contribution
    \subitem to the $\rightarrow$ arithmetic genus, 156
    \subitem to the intersection number, 168
    \subitem to the Willmore energy $\willmore_{\text{\scriptsize\rm sing}}$, 
		55
  \item Lorentz spaces, 120

  \indexspace

  \item M--curve, 10
  \item minimizer
    \subitem $\fermi_{\min}(\cdot)$, 170
    \subitem existence of a $\sim$, 3
    \subitem of the generalized Willmore functional, 180, 189--190
    \subitem relative $\sim$, 157--170
  \item moduli space, 8, 67
    \subitem $\moduli_1$ of conformal classes, 3, 185
    \subitem $\moduli_{\lattice,\eta},                      \moduli_{\lattice,\eta,\willmore},                      \moduli_{\lattice,\eta,\sigma},                      \moduli_{\lattice,\eta,\sigma,\willmore}$, 
		79
    \subitem $\moduli_{\lattice}$, 145
    \subitem $\moduli_{g,\lattice,\eta},                      \moduli_{g,\lattice,\eta,\willmore},                      \moduli_{g,\lattice,\eta,\sigma},                      \moduli_{g,\lattice,\eta,\sigma,\willmore}$, 
		79
    \subitem $\moduli_{g,\lattice}$, 68
    \subitem  $\Bar{\moduli}_{\lattice,\eta},                       \Bar{\moduli}_{\lattice,\eta,\willmore},                       \Bar{\moduli}_{\lattice,\eta,\sigma},                       \Bar{\moduli}_{\lattice,\eta,\sigma,\willmore}$, 
		139
    \subitem  $\Bar{\moduli}_{g,\eta}$, 80
    \subitem  $\Bar{\moduli}_{g,\lattice,\eta},                       \Bar{\moduli}_{g,\lattice,\eta,\willmore},                       \Bar{\moduli}_{g,\lattice,\eta,\sigma},                       \Bar{\moduli}_{g,\lattice,\eta,\sigma,\willmore}$, 
		80
    \subitem blowing up of the $\sim$, 74
    \subitem compactification of the $\sim$, 80, 139
    \subitem cusp of the $\sim$, 73
    \subitem tangent space of the $\sim$, 68
  \item multiple point, 155, 156, 164, 165

  \indexspace

  \item neighbourhood
    \subitem $\Set{U}^{\pm}_{\varepsilon,\delta}$, 27
    \subitem $\Set{U}_{p}^{\pm}$, 82
    \subitem $\Set{V}^{\pm}_{\varepsilon,\delta}$, 27
    \subitem $\mathbb{P}_{p}^{\pm}$, 84
    \subitem condition $\sim$ (i)--(ii), 81--82
  \item normalization
    \subitem connected components of the $\sim$, 28
    \subitem disconnected, 158, 173
    \subitem of finite genus, 42
    \subitem of the eigenfunction, 42

  \indexspace

  \item orbit, 155
    \subitem of type~1, 155, 156, 164, 165
    \subitem of type~2, 155, 156, 164
    \subitem of type~3, 155, 156
  \item order of zeroes $\text{\rm ord}_{z}(\psi)$, 122, 128, 129
  \item ordinary double point
    \subitem --s of the free \Em{complex Fermi curve}        $(k^+_{\kappa},k^-_{\kappa})$, 
		27
    \subitem at infinity $(\infty^-,\infty^+)$, 44

  \indexspace

  \item parity transformation $\parity$, 31
  \item periodicity condition, 5
    \subitem equivalent form of the $\sim$, 5
  \item perturbation
    \subitem $\rightarrow$ finite rank $\sim$, 99
  \item Picard group, 12, 59

  \indexspace

  \item quasi--momenta
    \subitem $\xx{p}$ and $\yy{p}$, 15
    \subitem $k$, 15
    \subitem condition $\sim$  (i)--(iv), 43--44
    \subitem condition $\sim$  (v), 68
    \subitem condition $\sim$ (i') and (vi), 144
    \subitem condition $\sim$ (ii'), 80
  \item quaternionic
    \subitem complex line bundle, 45
    \subitem function theory, 7--9
    \subitem line bundle, 4
    \subitem singularity condition, 7
    \subitem Weierstra{\ss} curve, 190

  \indexspace

  \item real
    \subitem --$\sigma$--hyperelliptic, 158
    \subitem part, 30
      \subsubitem dividing $\sim$, 157
  \item region of periodic lines, 85
  \item relative minimizer, 157--170
  \item residue, 59, 79
  \item resolvent
    \subitem $\sim$ formula, 105
    \subitem of the Dirac operator $\Op{R}$, 17
    \subitem weak continuity of the $\sim$, 19

  \indexspace

  \item S\'{e}rre duality, 12, 59
  \item sheaf
    \subitem dualizing $\sim$, 32
    \subitem structure $\sim$, 39
  \item singularity
    \subitem $\rightarrow$ cusp, 155
    \subitem $\rightarrow$ multiple point, 155
    \subitem condition, 7, 197
      \subsubitem quaternionic $\sim$, 7
      \subsubitem weak $\sim$, 7, 155
    \subitem condition $\sim$ (i)--(ii), 102
    \subitem of the holomorphic structure, 55, 102
  \item Sobolev
    \subitem constant $S_p$, 17
  \item spectral
    \subitem curve, 6
    \subitem projection
      \subsubitem $\Op{P}$, 32, 34
      \subsubitem $\Op{P}_{\xx{\gamma}}$, 36
      \subsubitem $\triv{\Op{P}}_{\mu}$, 40
  \item structure
    \subitem $\ell^+_1(\lattice\dual)$--$\sim$, 145
    \subitem differentiable $\sim$, 9
    \subitem sheaf, 39
      \subsubitem arithmetic genus of the $\sim$, 39
  \item symplectic form, 60

  \indexspace

  \item tangent
    \subitem space of the moduli space, 68
    \subitem space of the phase space, 12, 59
  \item Theta function $\theta_{\Delta}$, 93
  \item topology
    \subitem compact open $\sim$, 78
    \subitem finite $\sim$, 78
    \subitem of the moduli space, 78
  \item trivialization, 16

  \indexspace

  \item unique continuation property, 23--25

  \indexspace

  \item Weierstra{\ss}
    \subitem curve, 6
    \subitem generalized $\sim$ curve, 9, 157
    \subitem generalized $\sim$ potential, 8, 9
    \subitem potential, 5, 9
    \subitem quaternionic $\sim$   curve, 190
  \item Willmore
    \subitem conjecture, 3
    \subitem functional, 3
      \subsubitem generalized $\sim$, 157
    \subitem torus, 4
      \subsubitem constrained $\sim$, 4, 201, 202

\end{theindex}

\begin{thebibliography}{MMM}
\bibitem[Ad]{Ad} R.\ A.\ Adams: Sobolev spaces. Pure and Applied
  Mathematics {\bf 65}. Academic Press, New York (1975).
\bibitem[Ah]{Ah} L.\ V.\ Ahlfors: Complex analysis, 
  third edition. McGraw--Hill, Inc., Singapore (1979).
\bibitem[A-McK-M]{AMcKM} H.\ Airault, H.\ McKean, J.\ Moser: Rational
  and elliptic solutions of the Korteweg--de Vries Equation and a
  related many--body problem. Commun.\ Pure and appl.\ Math.\
  {\bf 30}, 95--148 (1977).
\bibitem[A-K]{AK} S.\ Albeverio, P.\ Kurasov: Singular perturbations
  of differential operators.
  London Mathematical Society Lecture Note Series {\bf 271}.
  Cambridge University Press, Cambridge (2000).
\bibitem[A-C-G-H]{ACGH} E.\ Arbarello, M.\ Cornalba, P.\ A.\
  Griffiths, J.\ Harris: Geometry of algebraic curves. Volume I.
  Grundlehren der mathematischen Wissenschaften {\bf 267}.
  Springer, New York (1985).
\bibitem[Ba]{Ba} M.\ V.\ Babich: Willmore surfaces, 4--particle Toda
  lattice and double coverings of hyperelliptic surfaces. Amer.\
  Math.\ Soc.\ Transl.\ {\bf 174}, 143--168 (1996).
\bibitem[B-B]{BB} M.\ V.\ Babich, A.\ I.\ Bobenko: Willmore tori with
  umbilic lines and minimal surfaces in hyperbolic space. Duke Math.\
J.\ {\bf 72}, 151--185 (1993).
\bibitem[Ba-St]{BaSt} C.\ B\u{a}nic\u{a},
  O.\ St\u{a}n\u{a}\c{s}il\u{a}: Algebraic methods in the
  global theory of complex spaces. Editura Academiei, Bucuresti, John
  Wiley \& Sons, London (1976).
\bibitem[Bat]{MOT} A\ Erdelyi (ed.), W.\ Magnus, F.\ Oberhettinger, F.\
  G.\ Tricome: Bateman Manuscript Project, Higher Transzendental
  Functions Volume II. McGraw--Hill Book Company, Inc., New York (1953).
\bibitem[B-B-E-I-M]{BBEIM} E.\ D.\ Belokos, A.\ I.\ Bobenko, V.\ Z.\
  Enol'skii, A.\ R.\ Its, V.\ B.\ Matveev: Algebro--geometric approach
  to nonlinear integrable equations. Springer Series in Nonlinear
  Dynamics, Springer, Berlin (1994).
\bibitem[B-S]{BS} C.\ Bennett, R.\ Sharpley: Interpolation of
  operators.  Pure and Applied Mathematics {\bf 129}. Academic Press,
  Orlando (1988).
\bibitem[Bo]{Bo} A.\ I.\ Bobenko: Constant mean curvature surfaces and
  integrable equations. Uspekhhi Mt.\ Nauk {\bf 46:4}, 3--42 (1991);
  English transl.\, Russian Math.\ Surveys {\bf 46:4}, 1--45 (1991).
\bibitem[B-C-R]{BCR} J.\ Bochnak, M.\ Coste, M.--F.\ Roy: Real algebraic
  geometry.  Ergebnisse der Mathematik und ihrer Grenzgebiete 3. Folge
  {\bf 36}. Springer, Berlin (1998).
\bibitem[B-F-L-P-P]{BFLPP} F.\ Burstall, D.\ Ferus K.\ Leschke,
  F.\ Pedit, U.\ Pinkall: Conformal geometry of surfaces in $S^4$ and
  quaternions.  Lecture Notes in Mathematics {\bf 1772}.
  Springer, Berlin, New York (2002).
\bibitem[Ca]{Ca} T.\ Carleman: Sur un probleme d'unicite pour les
  sytemes d'equations aux derivees partielles a doux variables
  independantes. Ark.\ Mat.\ Astron.\ Fys.\ {\bf B 26} No.\ 17, 1--9
  (1939).
\bibitem[Ch]{Ch} B.\ Chen: Geometry of submanifolds. Pure and Applied
  Mathematics {\bf 22}. Marcel Dekker, New York (1973).
\bibitem[Co-I]{Co1} J.\ B.\ Conway: Functions of one complex
  variable. Graduate Texts in Mathematics {\bf 11}. Springer, New York
  (1973).
\bibitem[Co-II]{Co2} J.\ B.\ Conway: Functions of one complex variable
  II. Graduate Texts in Mathematics {\bf 159}. Springer, New York (1995).
\bibitem[DF]{DF} A.\ Defant, K.\ Floret: Tensor norms and operator
  ideals. North--Holland Mathematical Studies
  {\bf 176}. North--Holland, Amsterdam (1993).
\bibitem[Do]{Do} A.\ Dold: Lectures on algebraic
  topology. Grundlehren der mathematischen Wissenschaften
  {\bf 200}. Springer, Berlin (1972).
\bibitem[D-K-N]{DKN} B.\ A.\ Dubrovin, I.\ M.\ Krichever, S.\ P.\
  Novikov: Integrable systems I. In: V.\ I.\ Arnold, S.\ P.\ Novikov
  (eds.) Dynamical Systems IV. Encyclopedia of Mathematical Sciences
  vol.\ 4, pp. 173--280. Springer, Berlin, Heidelberg (1990).
\bibitem[E-K]{EK} F.\ Ehlers, H.\ Kn\"orrer: An algebro--geometric
  interpretation of the B\"acklund--transformation for the
  Korteweg--de Vries equation.
  Comment.\ Math.\ Helv.\ {\bf 57}, 1--10 (1982).
\bibitem[Ei]{Ei} L.\ P.\ Eisenhart: A treatise on the differential
  geometry of curves and surfaces. Ginn and Company, Boston, (1909).
\bibitem[F-K]{FK} H.\ M.\ Farkas, I.\ Kra: Riemann surfaces.
  Graduate Texts in Mathematics {\bf 71}. Springer, New York (1980).
\bibitem[F-K-T]{FKT} J.\ Feldman, H.\ Kn\"orrer, E.\ Trubowitz:
  Riemann surfaces of infinite genus I--IV. ETH Z\"urch preprints
  (1993).
\bibitem[F-L-P-P]{FLPP} D.\ Ferus, K.\ Leschke, F.\ Pedit,
  U.\ Pinkall: Quaternionic holomorphic geometry: Pl\"ucker formula,
  Dirac eigenvalue estimates and energy estimates of harmonic 2--tori.
  Invent.\ Math.\ {\bf 146}, 507--593 (2001).
\bibitem[F-P-P-S]{FPPS} D.\ Ferus, F.\ Pedit, U.\ Pinkall,
  I.\ Sterlin: Minimal tori in $S^4$. J.\ reine angew.\ Math.\ {\bf 429},
  1--47 (1992).
\bibitem[F-R]{FR} A.\ Fig\`{a}--Talamanca, D.\ Rider: A theorem of
  Littlewood and lacunary series for compact groups. Pacific J.\ of
  Math.\ {\bf 16}, 505--514 (1966).
\bibitem[Fi]{Fi} G.\ Fischer: Complex analytic geometry. Lecture Notes
  in Mathematics {\bf 538}. Springer, Berlin, Heidelberg (1976).
\bibitem[Fo]{Fo} O.\ Forster: Lectures on Riemann surfaces. Graduate
  Texts in Mathematics {\bf 81}. Springer, New York (1981).
\bibitem[Fr-1]{Fr1} T.\ Friedrich: Dirac--Operatoren in der
  Riemannschen Geomtrie. Vieweg Verlag, Braunschweig (1997).
  engl. transl.: T.\ Friedrich: Dirac operators in Riemannian
  geometry. Graduate Studies in Mathematics {\bf 25},
  American Mathematical Society, Providence, Rhode Island (2000).
\bibitem[Fr-2]{Fr2} T.\ Friedrich: On the spinor representation of
  surfaces in Euclidean $3$--space. J.\ Geom.\ Phys.\ {\bf 28}, 143--157
  (1998). 
\bibitem[G-K-T]{GKT} D.\ Gieseker, H.\ Kn\"orrer, E.\ Trubowitz:
  The geometry of algebraic Fermi curves.
  Perspectives in Math. {\bf 14}. Academic Press, Inc., Boston (1993).
\bibitem[G-K]{GK} S.\ G.\ Gindikin, G.\ M.\ Khenkin (eds.):
  Several complex variables IV. Encyclopedia of Mathematical Sciences
  vol.\ 10. Springer, Berlin, Heidelberg (1990).
\bibitem[G-J]{GJ} J.\ Glimm, A.\ Jaffe: Quantum physics. Springer, New
  York 1987.
\bibitem[Gr]{Gr} H.\ Grauert: Der Satz von Kuranishi f\"ur kompakte
  komplexe R\"aume. Invent.\ Math.\ {\bf 25}, 107--142 (1974).
\bibitem[G-P-R]{GPR} H.\ Grauert, Th.\ Peternell, R.\ Remmert (eds.):
  Several complex variables VII. Encyclopedia of Mathematical Sciences
  vol.\ 74. Springer, Berlin, Heidelberg (1994).
\bibitem[Gr-Re]{GrRe} H.\ Grauert, R.\ Remmert: Coherent analytic
  sheaves. Grundlehren der mathematischen Wissenschaften
  {\bf 265}. Springer, Berlin, Heidelberg (1984).
\bibitem[Gr-Ha]{GrHa} P.\ Griffiths, J.\ Harris: Principles of
  algebraic geometry. Pure and applied mathematics. John Wiley \&
  Sons, New York (1978).
\bibitem[G-S-1]{GS1} P.\ Grinevich, M.\ U.\ Schmidt: Period preserving
  flows and the moduli space of periodic solutions of soliton
  equations. Physica D {\bf 87}, 73--98 (1995).
\bibitem[G-S-2]{GS2} P.\ Grinevich, M.\ U.\ Schmidt: Conformal
  invariant functionals of immersions of tori into
  $\mathbb{R}^3$. Journal of Geometry and Physics {\bf 26}, 51--78
  (1998).
\bibitem[G-H]{GH} B.\ H.\ Gross, J.\ Harris: Real algebraic
  curves. Ann.\ scient.\ \'{E}c. Norm.\ Sup.\ {\bf 14}, 157--182
  (1981).
\bibitem[Gu-Ro]{GuRo} R.\ C.\ Gunning, H.\ Rossi: Analytic functions
  of several complex variables. Prentice--Hall, Inc., Englewood Cliffs
  (1965). 
\bibitem[Har]{Har} R.\ Hartshorne: Algebraic geometry.
  Graduate Texts in Mathematics {\bf 52}. Springer, New York (1977).
\bibitem[H-K-Z]{HKZ} H.\ Hedenmalm, B.\ Korenblum, K.\ Zhu: Theory of
  Bergman spaces. Graduate Texts in Mathematics {\bf 199}. Springer,
  New York (2000).
\bibitem[He]{He} S.\ Helgason: Differential geometry, Lie groups and
  symmetric spaces.  Pure and Applied Mathematics {\bf 80}.
  Academic Press, New York (1978).
\bibitem[Hi]{Hi} N.\ Hitchin: Harmonic maps from a $2$--torus into the
  $3$--sphere. J. Diff. Geometry {\bf 31}, 627--710 (1990).
\bibitem[Ho]{Ho} N.\ R.\ Howes: Modern analysis and topology.
  Universitext. Springer, New York (1995).
\bibitem[Je]{Je} D.\ Jerison: Carleman inequalities for the Dirac and
  Laplace operators and unique continuation. Advances in Math {\bf 62},
  118--134 (1986).
\bibitem[Ka]{Ka} T.\ Kato: Perturbation theory of linear operators. 
  Grundlehren der mathematischen Wissenschaften
  {\bf 132}. Springer, Berlin (1966).
\bibitem[Ke]{Ke} J.\ L.\ Kelley: General topology. Graduate
  Texts in Mathematics {\bf 27}. Springer, New York (1955).
\bibitem[Ki-1]{Ki1} Y.\ M.\ Kim: Unique continuation theorem for the
  Dirac operator and the Laplace operator. MIT Thesis, (1989).
\bibitem[Ki-2]{Ki2} Y.\ M.\ Kim: Carleman inequalities for the Dirac
  operator and strong unique continuation. Proceedings of the AMS
  {\bf 123}, 2103--2111 (1995).
\bibitem[Kis]{Kis} S.\ V.\ Kislyakov: Classical themes of Fourier
  analysis. In V.\ P.\ Khavin, N.\ K.\ Nilol'skij (eds.): Commutative
  harmonic analysis~I. Encyclopedia of Mathematical Sciences
  vol.\ 15. Springer, Berlin, Heidelberg (1991).
\bibitem[Kn]{Kn} A.\ W.\ Knapp: Representation theory of semisimple
  groups. Princeton Mathematical Series {\bf 36}.
  Princeton University Press, Princeton (1970).
\bibitem[Ko]{Ko} K.\ Kodaira: Complex manifolds and deformations of
  complex structures. Grund\-leh\-ren der mathematischen Wissenschaften
  {\bf 283}. Springer, New York (1986).
\bibitem[Kon]{Kon} B.\ G.\ Konopelchenko: Introduction to
  multidimensional integrable equations. Plenum Press, New York (1992).
\bibitem[Kr-1]{Kr1} I.\ M.\ Krichever: Methods of algebraic geometry in
  the theory of non--linear equations. Russ.\ Math.\ Surv.\ {\bf 32} No.6,
  185--213 (1977).
\bibitem[Kr-2]{Kr2} I.\ M.\ Krichever: Elliptic solutions of the
  Kadomtsev--Petviashvilli equation and integrable systems of
  particles. Funct. Anal. Appl. {\bf 14}, 282--289 (1981).
\bibitem[Ku]{Ku} P.\ Kuchment: Flochet theory for partial differential
equations. Operator Theory Advances and Applications {\bf 60}.
Birkh\"auser, Basel 1993.
\bibitem[Kun]{Kun} E.\ Kunz: Holomorphe Differentialformen auf
  algebraischen Variet{\"a}ten mit Singularit{\"a}ten I. Manuscripta
  Math.\ {\bf 15}, 91--108 (1975).
\bibitem[La]{La} S.\ Lang: Differential and Riemannian manifolds. 
  Graduate Texts in Mathematics {\bf 160}. Springer, New York
  (1995).
\bibitem[L-B]{LB} H.\ Lange, C.\ Birkenhake: Complex abelian
  varieties. Grundlehren der mathematischen Wissenschaften {\bf 302}.
  Springer, Berlin (1992).
\bibitem[Law]{Law} H.\ B.\ Lawson: Complete minimal surfaces in
  $S^3$. Ann. of Math. {\bf 92}, 335--374 (1970).
\bibitem[L-Y]{LY} P.\ Li, S.\ T.\ Yau: A new conformal invariant and
  its applications to the Willmore conjecture and the first eigenvalue
  of compact surfaces. Invent.\ Math.\ {\bf 69}, 269--291 (1982).
\bibitem[L-McL]{LMcL} Y.\ Li, D.\ W.\ McLaughlin: Morse and Melnikov
  functions for NLS Pde's. Commun.\ Math.\ Phys.\ {\bf 162}, 175--214
  (1994). 
\bibitem[Li]{Li} J.\ E.\ Littlewood: On the mean value of power
  series. Proc.\ London Math.\ Soc.\ {\bf 25}, 528--556 (1924).
\bibitem[Ma]{Ma} N.\ Mandache: Some remarks concerning unique
  continuation for the Dirac operator. Letters Math. Phys. {\bf 31},
  85--92 (1994).
\bibitem[M-O]{MO} V.\ A.\ Marchenko, I.\ V.\ Ostrovskii: A
  characterization of the spectrum of Hill's operator. 
  Math.\ USSR Sbornik {\bf 26}, 493--554 (1975).
\bibitem[M-S]{MS} V.\ B.\ Matveev, M.\ A.\ Salle: Darboux
  transformations and solitons.  Springer Series in Nonlinear
  Dynamics, Springer, Berlin (1991).
\bibitem[Mi]{Mi} E.\ Michael: Topologies on spaces of
  subsets. Trans.\ Am.\ Math.\ Soc.\ {\bf 71}, 152--182 (1951).
\bibitem[O]{O} R.\ O'Neil: Convolution operators and $L(p,q)$ spaces.
  Duke Math.\ J.\ {\bf 30}, 129-142 (1963).
\bibitem[Pi-1]{Pi1} U.\ Pinkall: Hopf tori in $S^3$. Invent.\ Math.\
  {\bf 81}, 379--386 (1985).
\bibitem[Pi-2]{Pi2} U.\ Pinkall: Regular homotopy classes of immersed
  surfaces. Topology {\bf 24}, 421--434 (1985).
\bibitem[P-P]{PP} F.\ Pedit, U.\ Pinkall: Quaternionic analysis on
  Riemann surfaces and differential geometry. Doc.\ Math.\ J.\ DMV,
  Extra Volume ICM 1998, Vol.~II, 389--400 (1999).
\bibitem[P-S]{PS} U.\ Pinkall, I.\ Sterling: On the classification of
  constant mean curvature tori. Annals of Math. {\bf 130}, 407--451
  (1989). 
\bibitem[Pr-Se]{PrSe} A.\ Pressley, G.\ Segal: Loop goups.
  Oxford Mathematical Monographs. Oxford University Press, New York (1986).
\bibitem[R-S-I]{RS1} M.\ Reed, B.\ Simon: Methods of modern
  mathematical physics vol.\ I: Functional analysis. Academic Press,
  Inc., London (1980).
\bibitem[R-S-II]{RS2} M.\ Reed, B.\ Simon: Methods of modern
  mathematical physics vol.\ II: Fourier Analysis,
  Self--Adjointness. Academic Press, Inc., London (1975).
\bibitem[R-S-IV]{RS4} M.\ Reed, B.\ Simon: Methods of modern
  mathematical physics vol.\ IV: Analysis of operators. Academic
  Press, Inc., London (1978).
\bibitem[Ro-1]{Ro1} H.\ L.\ Royden: The Riemann--Roch Theorem.
  Comment.\ Math.\ Helvetici {\bf 34}, 37--51 (1960).
\bibitem[Ro-2]{Ro2} H.\ L.\ Royden: Real analysis. 3 rd edition.
  Prentice--Hall, Inc., Englewood Cliffs (1988).
\bibitem[Sch]{Sch} M.\ U.\ Schmidt: Integrable systems and Riemann
  surfaces of infinite genus. Memoirs of the American Math.\ Soc.\
  {\bf 581} (1996).
\bibitem[S-T]{ST} H.\ Seifert, W.\ Threlfall: A textbook of
  topology. Pure and Applied Mathematics {\bf 89}. Academic Press,
  New York (1980).
\bibitem[Se]{Se} J.--P.\ S\'{e}rre: Algeraic groups and class fields.
  Graduate Texts in Mathematics {\bf 117}. Springer, New York (1988). 
\bibitem[Si-1]{Si1} L.\ Simon: Existence of Willmore
  surfaces. Miniconference, Canberra on geometry and partial
  differential equations. Proc.\ Cent.\ Math.\ Anal., Australian
  National University {\bf 10}, 187--216 (1986).
\bibitem[Si-2]{Si2} L.\ Simon: Existence of surfaces minimizing the
  Willmore functional. Commun.\ in Anal.\ and Geom.\ {\bf 1},
  281--326 (1993). 
\bibitem[So]{So} C.\ D.\ Sogge: Fourier integrals in classical
  analysis. Cambridge tracts in mathematics {\bf 105}. Cambridge
  University Press, Cambridge (1993).
\bibitem[St]{St} E.\ M.\ Stein: Singular integrals and
  differentiability properties of functions.
  Princeton Mathematical Series {\bf 30}.
  Princeton University Press, Princeton (1970).
\bibitem[S-W]{SW} E.\ M.\ Stein, G.\ Weiss: Introduction to Fourier
  analysis on Euclidean spaces. Princeton Mathematical Series {\bf 32}.
  Princeton University Press, Princeton 1971.
\bibitem[Str]{Str} M.\ Struwe: Variational Methods, Ergebnisse der
  Mathematik und ihrer Grenzgebiete 3. Folge {\bf 34}. Springer,
  Berlin (1990).
\bibitem[Ta-1]{Ta1} I.\ A.\ Taimanov: Modified Novikov--Veselov equation
  and differential geometry of surfaces. Am.\ Math.\ Soc.\ Transl.\
  {\bf 179}, 133--151 (1997).
\bibitem[Ta-2]{Ta2} I.\ A.\ Taimanov: The Weiestra{\ss} representation
  of closed surfaces in $\mathbb{R}^3$.
  Funct.\ Anal.\ Appl. {\bf 32}, 258--267 (1998).
  and differential geometry of surfaces. Am.\ Math.\ Soc.\ Transl.\
  {\bf 179}, 133--151 (1997).
\bibitem[Tk]{Tk} V.\ Tkachenko: Spectra of non--selfadjoint Hill's
  operators and a class of Riemann surfaces. Ann.\ Math.\ {\bf 143},
  181--231 (1996).
\bibitem[Wh]{Wh} H.\ Whitney: Complex analytic varieties.
  Addisson--Wesley Publishing Company, Inc., Menlo Park, London (1972).
\bibitem[Wi-1]{Wi1} T.\ J.\ Willmore: Total curvature in Riemannian
  geometry. Ellis Horwood Limited, Chicester (1982).
\bibitem[Wi-2]{Wi2} T.\ J.\ Willmore: Riemannian Geometry. Clarendon Press,
  Oxford (1993).
\bibitem[Wo]{Wo} T.\ H.\ Wolff: Recent work on sharp estimates in
  second--order elliptic unique continuation problems. Journal of
  Geometric Analysis {\bf 3}, 621--650 (1993).
\bibitem[Yo]{Yo} K.\ Yosida: Functional analysis. 3 rd edition.
  Grundlehren der mathematischen Wissenschaften {\bf 123}.
  Springer, Berlin (1971).
\bibitem[Zh]{Zh} K.\ Zhu: Operator Theory in functions spaces.
  Pure and Applied Mathematics {\bf 139}. Marcel Dekker, Inc.,
  New York (1990).
\bibitem[Zi]{Zi} W.\ P.\ Ziemer: Weakly differentiable functions.
  Graduate Texts in Mathematics {\bf 120}. Springer, New York
  (1989).
\end{thebibliography}
\end{document}